%% file: main.tex
\documentclass{article}
\usepackage[english]{babel}
\usepackage[utf8]{inputenc}
\usepackage[T1]{fontenc}
\usepackage{graphicx}
\usepackage[a4paper,top=3cm,bottom=2cm,left=3cm,right=3cm,marginparwidth=2.0cm]{geometry}
\usepackage{amssymb, amsmath, amsthm}
\usepackage{hyperref}
\usepackage{authblk}
\usepackage{longtable}
\usepackage{multicol}
\usepackage{subcaption} 
\usepackage{booktabs}
\usepackage{threeparttable}
\usepackage[dvipsnames]{xcolor}
\usepackage{changepage}

\DeclareMathOperator{\rot}{rot}
\DeclareMathOperator{\reid}{Reid}

\title{A table of knotoids in $S^3$ up to seven crossings}

\author[1,2]{Boštjan Gabrovšek}
\author[3]{Paolo Cavicchioli}

\affil[1]{Faculty of Education, University of Ljubljana}
\affil[2]{Rudolfovo -- Science and Technology Centre Novo mesto}
\affil[3]{IMFM -- Institute of Mathematics, Physics and Mechanics, Ljubljana}

\affil[ ]{\textit{Emails:} bostjan.gabrovsek@pef.uni-lj.si, paolo.cavicchioli@gmail.com}

\date{\today}

\begin{document}

\maketitle

\begin{abstract}
We present a complete classification of spherical knotoids with up to six crossings and conjecture that our classification up to seven crossings is complete.
Our work extends the tradition of knot tabulation to the setting of knotoids introduced by Turaev. 
We describe the methods used to enumerate diagrams, simplify them, and distinguish equivalence classes using a collection of invariants including the Kauffman bracket, the Arrow polynomial, the Affine index polynomial, the Mock Alexander polynomial, and the Yamada polynomial of the  closure. 
We also investigate the chirality and rotational symmetries of these knotoids. 
Applications to protein entanglement illustrate the importance of such classifications. 
\end{abstract}

\section{Introduction}

The theory of knotoids, introduced by Turaev~\cite{Turaev2012}, provides a natural extension of classical knot theory. A knotoid is defined as an equivalence class of generic immersions of an oriented unit interval $[0, 1]$ into a surface, typically the 2-sphere $S^2$ or the plane $\mathbb{R}^2$. 

\begin{figure}[ht]
    \centering
    \begin{subfigure}[b]{0.3\textwidth}
        \includegraphics[width=\textwidth, page=1]{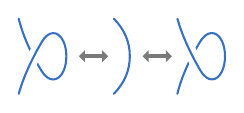}
        \caption{I}
    \end{subfigure}
    \begin{subfigure}[b]{0.3\textwidth}
        \includegraphics[width=\textwidth, page=2]{images.pdf}
        \caption{II}
    \end{subfigure}
    \begin{subfigure}[b]{0.3\textwidth}
        \includegraphics[width=\textwidth, page=3]{images.pdf}
        \caption{III}
    \end{subfigure}
    \caption{Reidemeister moves of knotoid diagrams.}
    \label{fig:reid}
\end{figure}

\begin{figure}[ht]\centering
        \includegraphics[page=4]{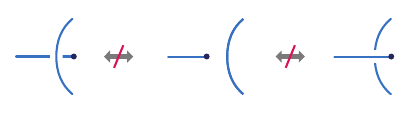}
    \caption{Forbidden moves}
    \label{fig:forbidden}
\end{figure}

More formally, if a $\Sigma$ is a surface, a knotoid diagram $D$ in $\Sigma$ is a generic immersion of the interval $[0,1]$ in the interior of $\Sigma$ whose only singularities are transversal double points endowed with over/undercrossing data. Two diagrams are considered equivalent if one can be transformed into the other by a finite sequence of Reidemeister moves I, II, and III depicted in \autoref{fig:reid}, and planar isotopy. A \emph{knotoid} is an equivalence class of such diagrams. A crucial distinction from classical knot theory is that endpoints cannot be pulled across strands. The two forbidden moves are shown in \autoref{fig:forbidden}. This restriction allows for the existence of nontrivial open-ended entanglements.


Knotoids defined in $S^2$ are called \emph{spherical knotoids}, while those in $\mathbb{R}^2$ are called \emph{planar knotoids}. There exists a natural surjective map from planar to spherical knotoids, but it is not injective: non-equivalent planar knotoids may become equivalent in $S^2$. In this paper we will work with spherical knotoids exclusively.

\begin{figure}[ht]\centering
        \includegraphics[page=5, scale=0.8]{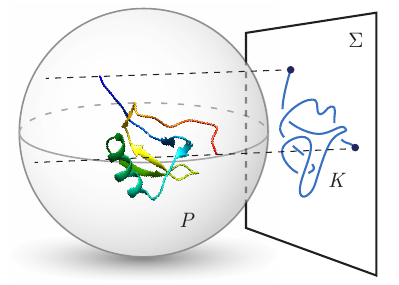}
    \caption{We obtain a knotoid representing a protein's backbone $P$ by placing the protein in a large enough sphere and projecting it to the tangent plane $\Sigma$ or to the sphere directly.}
    \label{fig:projection}
\end{figure}

Beyond their mathematical interest, knotoids have important applications in molecular biology. Proteins are open-ended polymers whose folded structures may form non-trivial entanglements \cite{Goundaroulis2017, gugumcu2022invariants, bruno2024knots}. Classical knot theory requires closing the chain into a closed loop, a process that alters the geometry and can obscure topological features. The knotoid approach avoids artificial closure and allows the study of open chains directly by studying the statistics of different projections of the protein's backbone to a surface as depicted in \autoref{fig:projection}~\cite{Dorier2018}.

This methodology has led to the development of tools such as Knoto-ID \cite{Dorier2018}, which provides topological “fingerprints” of protein chains.

The classification of knotoids follows the long tradition of knot tabulation, initiated by Tait and extended through large-scale enumerations \cite{hoste2005enumeration, hoste1998first, burton2020next}. Results up to six crossings, covering both planar and spherical knotoids, were reported in a 2019 preprint \cite{goundaroulis2019systematic}. In this paper we present an independent and complete classification of spherical knotoids up to six crossings and a conjectured classification up to seven crossings, where are unable to verify uniqueness of only 14 knottoids due to a lack of strong enough knotoid invariants.

In the case of classical knots, large-scale classification is made feasible by the geometry of knot complements. The complement of a knot in $S^3$ is a 3-manifold, and by Thurston’s theorem almost all such complements are hyperbolic. This allows one to compute canonical hyperbolic triangulations, which provide strong and efficient tools for distinguishing knots on a massive scale \cite{burton2020next}.

For knotoids, however, the complement is not a 3-manifold and as a result, we cannot appeal to hyperbolic geometry, and the classification must instead proceed by brute-force knotoid simplifications via Reidemeister moves and by comparing a range of invariants. Some of these invariants are themselves computationally demanding, such as the Yamada polynomial of a knotoid closure, which makes large-scale classification a different type of computational challenge.

\section{Knotoid invariants}

The main idea of the classification procedure can be summarized as follows:
\begin{enumerate}
    \item Generate all possible knotoid diagrams up to a given number of crossings.
    \item Compute a set of invariants for each diagram.
    \item For diagrams sharing identical invariant values, simplify them by sequences of Reidemeister moves to obtain minimal representatives and remove duplicates.
\end{enumerate}

The classification relies on the following collection of invariants, many of which are extensions of classical knot invariants.

\paragraph{Bracket polynomial.} The Bracket polynomial for knotoids was introduced by Turaev \cite{Turaev2012} by a state-sum expansion equivalent to the (Kauffman) bracket polynomial defined in \cite{kauffman1987state}. The Bracket polynomial is a special case of the Kauffman bracket skein module generalized for knotoids in \cite{gabrovvsek2023invariants,diamantis2021knotoids, gabrovsek2025bracket}.

\paragraph{Arrow polynomial.} First defined for virtual knots, the arrow polynomial was extended to classical and virtual knotoids by Gügümcü and Kauffman \cite{Gueguemcue2017a}. It is built from an oriented state expansion of the bracket polynomial. In each state, pairs of cusps cancel; surviving features contribute $\Lambda$–variables, and the resulting \emph{$\Lambda$–degree} gives a lower bound on the \emph{height} of a spherical knotoid (the minimum number of crossings a shortcut between endpoints must pass under, minimized over equivalent diagrams).

\paragraph{Affine index polynomial.} Originally defined for virtual knots, this polynomial was extended to classical and virtual knotoids by Gügümcü and Kauffman \cite{Gueguemcue2017a}. It is constructed from an integer labeling of the arcs of a flat knotoid diagram, obtained by ignoring over/under data. Its degree provides a lower bound for the height of the knotoid.

\paragraph{Mock Alexander polynomial.}  
Recently introduced in \cite{ellis2024mock}, the Mock Alexander polynomial is defined for \emph{admissible} diagrams, where the number of regions equals the number of crossings. For spherical knotoids, a diagram with $n$ crossings has $n+1$ regions, so admissibility is obtained by \emph{starring} one region, effectively removing it from the state calculations. By canonically starring the region adjacent to the tail, one obtains a well-defined invariant for spherical knotoids.

\begin{figure}[ht]\centering
        \includegraphics[page=6]{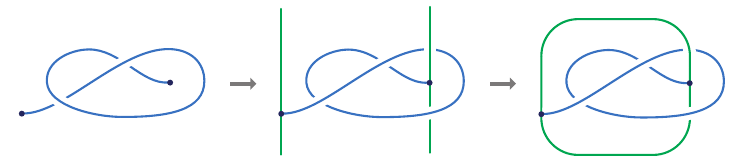}
    \caption{The closure of a knotoid.}
    \label{fig:closure}
\end{figure}

\paragraph{Yamada polynomial of the closure.}  
The Yamada polynomial is a topological invariant of spatial graphs \cite{Yamada1989}, defined by combining the skein relations of the (Kauffman) bracket polynomial for knots with those of the Tutte polynomial for graphs.  
The double-sided closure of a knotoid diagram (see \cite{Turaev1990, gugumcu2016survey, Gugumcu2018Thesis, kodokostas2019rail, gugumcu2022invariants}) produces a $\Theta$–curve \cite{dabrowski2024theta}, which can be viewed as a 3–regular spatial graph, see \autoref{fig:closure}.
For such graphs, the Yamada polynomial is a \emph{topological–vertex} invariant, in contrast to 4–valent graphs, where it is only a \emph{rigid–vertex} invariant \cite{Kauffman1989}.  
For knots, the Yamada polynomial specializes to the Dubrovnik version of the Kauffman polynomial.

\subsection{Prime knotoids and symmetries of knotoids}
\begin{figure}[ht]
    \centering
    \begin{subfigure}[b]{0.3\textwidth} 
        \includegraphics[scale=1.2, page=7]{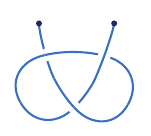}
        \caption{knot-like} \label{fig:knotlike}
    \end{subfigure}
    \begin{subfigure}[b]{0.3\textwidth}
        \includegraphics[scale=1.2, page=8]{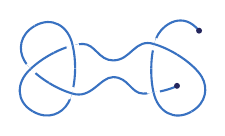}
        \caption{composite}\label{fig:composite}
    \end{subfigure}
    \begin{subfigure}[b]{0.3\textwidth}
        \includegraphics[scale=1.2, page=9]{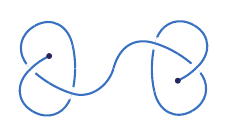}
        \caption{concatenation}\label{fig:concatenation}
    \end{subfigure}
    \caption{Structural types of knotoids.}
    \label{fig:struct}
\end{figure}

A knotoid diagram is called \emph{knot-like} if its two endpoints lie in the same region/face of the diagram \cite{Turaev2012} (\autoref{fig:knotlike}). Knot-like knotoids are in one-to-one correspondence classical knots.

A knotoid is a \emph{composite}, if it is a connected sum of a knotoid and a diagram of a classical knot (\autoref{fig:composite}).
A knotoid is a \emph{concatenation} if you can obtain it by join the head of one knotoid to the tail of another (\autoref{fig:concatenation}).

We call a knotoid \emph{prime}, if it cannot be expressed as a composite or concatenation of two non-trivial knots or knotoids, respectively.

In the tabulation of knotoids we also consider to classify knotoid symmetries. We define the following two key properties based on involutions of knotoid diagrams.

\begin{figure}[ht]
    \centering
    \begin{subfigure}[b]{0.4\textwidth}
        \includegraphics[scale=1.2, page=10]{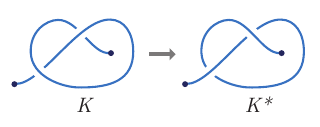}
        \caption{mirror} \label{fig:mirror}
    \end{subfigure} \hspace{4em}
    \begin{subfigure}[b]{0.4\textwidth}
        \includegraphics[scale=1.2, page=11]{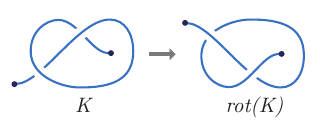}
        \caption{rotation}\label{fig:rotation}
    \end{subfigure}
    \caption{Two involutions on knotoids.}
    \label{fig:sym}
\end{figure}

\paragraph{Chirality.} The \emph{mirror image} of a knotoid diagram $K$, denoted $K^*$, is obtained by switching all over-crossings to an under-crossings and vice-versa (\autoref{fig:mirror}). A knotoid is \emph{achiral} (or \emph{amphicheiral}) if it is equivalent to its mirror image ($K = K^*$). Otherwise it is \emph{chiral} ($K \neq K^*$).

\paragraph{Rotatability.} The \emph{rotation} (or \emph{symmetry} \cite{ellis2024mock}) of a knotoid $K$, which we denote by $\rot(K)$, corresponds to a reflection of the diagram in $\mathbb{R}^2$ with respect to a vertical line, which extends to a self-homeomorphism of $S^2$. This can be visualized as a $180^\circ$ degree rotation of the diagram (\autoref{fig:rotation}). A knotoid is \emph{rotatable} if it is equivalent to its rotation ($K = \rot(K)$). Otherwise, it is \emph{non-rotatable} ($K \neq \rot(K)$).

\begin{figure}[ht]
    \centering
        \includegraphics[scale=1, page=12]{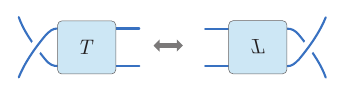}
    \caption{A flype.}
    \label{fig:flype}
\end{figure}

We also note here the \emph{flype} operation. Following \autoref{fig:flype}, a flype is a $180^\circ$ rotation of a 4-tangle $T$ inside of a diagram, where to neighboring strands meet at a crossing outside of $T$. Note that a flype does not change the equivalence class of the knotoids, but such a transformation would require a large amount of Reidemeister moves. When using brute-force approaches to determine knot-equivalences, it is often optimal to include a flype as an elementary move, as the authors have argued in \cite{gabrovvsek2017tabulation}. 

\section{Knotoid notation}

All operations (mirroring, rotation, Reidemeister moves, flypes) are programmed in Python using the \texttt{KnotPy} library~\cite{knotpy}. The full source code, results, and the Appendix are available at \cite{gabrovsek_knotoids_github}. Internally, the knotoid diagram is stored in an EM (Ewing-Millett) type notation \cite{ewing1991load}. The Appendix also provides a more commonly used DT (Dowker–Thistlethwaite) notation \cite{dowker1983classification}.

\begin{figure}[ht]
    \centering
    \begin{subfigure}{0.33\textwidth}
        \centering
        \includegraphics[width=10em,page=1]{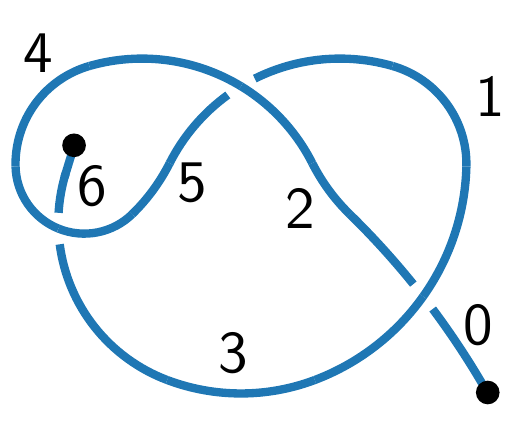}
        \caption{PD notation.}
        \label{fig:pd_example}
    \end{subfigure}
    \qquad
    \begin{subfigure}{0.33\textwidth}
        \centering
        \includegraphics[width=10em,page=2]{notation.pdf}
        \caption{EM notation.}
        \label{fig:em_example}
    \end{subfigure}
    \caption{Labeling for PD and EM knotoid notations for the knotoid $K3_2$.}
    \label{fig:notation}
\end{figure}

\paragraph{PD notation.}
Let $K$ be a knotoid diagram with $n$ arcs. 
We enumerate the arcs consecutively from $0$ to $n-1$. 
For each vertex (endpoint or crossing), we record the incident arcs in counterclockwise order, beginning with an undercrossing arc (in case of a crossing).
The PD code is then the ordered list consisting of the two endpoint entries and the entries corresponding to all crossings. 
For the diagram shown in \autoref{fig:pd_example}, the PD code is
\[ \mathrm{PD}(K3_2) =\texttt{ 
\text{``}[0],
[0,1,2,3],
[1,4,5,2],
[3,5,6,4],
[6]\text{''}}.
\]

\paragraph{EM notation.}
Let $K$ be a knotoid diagram. 
We label the endpoints and crossings consecutively by letters $A,B,C,\ldots$
At each vertex (endpoint or crossing), the incident half-edges are enumerated in clockwise order by $0,1,2,3$ (or by just $0$ in the case of an endpoint). 

For each vertex, we record the half-edges to which its half-edges are connected in order. Listing these data vertex by vertex produces the EM code. For the diagram in \autoref{fig:em_example}, the EM code is
\[\mathrm{EM}(K3_2) = \texttt{
\text{``}B0,
A0C0C3D0,
B1D3D1B2,
B3C2E0C1,
D2\text{''}}.
\]

Also note that each diagram can be uniquely expressed in a \emph{canonical EM form} (or canonical PD form), i.e.~there is a 1-1 correspondence between canonical EM notations and planar diagrams.
The canonical form is obtained by enumerating crossings and ingoing strands in all possible ways and taking the lexicographical minimal EM notation, a similar approach that was used in \cite{gabrovvsek2017tabulation}. 

\section{Methodology}
Our tabulation of prime spherical knotoids up to seven crossings is conducted through a multi-step computational procedure designed to systematically generate, simplify, and distinguish knotoid diagrams.

\begin{figure}[ht]
    \centering
        \includegraphics[scale=1, page=13]{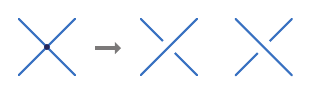}
    \caption{Replacing a 4-degree vertex with two possible crossing types.}
    \label{fig:crossing}
\end{figure}

\paragraph{Step 1: Generation of candidate diagrams.}
The initial set of candidates is generated by first creating all non-isomorphic planar graphs up to 9 vertices with connectivity 1 and minimal degree 1 using \texttt{plantri} software [75, 76], here we obtain 15,969 such graphs. Out of these planar graphs, we only take planar graphs with two degree one vertices and $n-2$ degree 4 vertices. We obtain 2479 such planar diagrams.
Next, we convert the graphs into knotoid diagrams by replacing every 4-degree vertex into a positive or negative crossing as depicted in \autoref{fig:crossing}. We also insert all rotations of the diagams, as \texttt{plantri} treats rotations as isomorphic diagrams. At this step we also remove diagram that have multiple connected components (linkoids \cite{gabrovvsek2023invariants}). We obtain 160,825 different knotoid diagrams, that is, different canonical EM codes.

\paragraph{Step 2: Simplification and filtering.}
At this step, we first  perform simple crossing reducing Reidemeister moves I and II and obtain 28,447 different knotoid diagrams. We remove 11,612 diagrams that represent connected sums (are composite) and 1163 diagrams that have bridges (are concatenations). 
Next, we perform brute-force Reidemeister moves for determining the minimal representative of a given knotoid. In more detail, we compute the whole space of Reidmeister moves, i.e. systematically perform all Reidemeister moves until no new diagrams appear in their canonical representation. 

We perform simplifications consecutively using the following parameters:
\begin{enumerate}
    \item simplification using Reidemeister moves (I, II, and III) that do not increase the number of crossings. This step reduced the number of candidates to 5863 knotoids,
    \item simplification using one Reidemeister II increasing move (reduction to 4049 knotoids),
    \item simplification using one Reidemeister II increasing move and flypes (reduction to 3115 knotoids),
    \item simplification using two  Reidemeister II increasing moves (reduction to 2600 knotoids).
\end{enumerate}

\paragraph{Step 3: Computation of invariants and partitioning}
For each of the 2600 remaining candidate diagrams, we computed a collection of polynomial invariants: the Kauffman bracket, the arrow polynomial, the Mock Alexander polynomial, the affine index polynomial, and the Yamada polynomial of the double closure \cite{kodokostas2019rail}. We note here that the Yamada polynomial is by far the most computatioaly expensive invairant to compute. We obtain 49 diagram with a unique set of invariants, which are now part of the final classification, and 772 groups of diagrams that share their invariant values.

\paragraph{Step 4: Exhaustive brute-force search}
In the next steps we systematicaly perform Reidemeister moves simultaniosly on groups of diagrams. In more detail, let $K_i$ be a knotoid diagram, we denote by $\reid_r(K_i)$ the set of all diagrams after systematicaly performing Reidemeister moves with $r \geq 0$ crossing-increasing II moves. Let now $K = \{K_1, K_2,\ldots,K_n\}$ is a set of knotoiod diagrams for which all invairants have the same value. In our exhaustive search, we compute $\reid_r(K_i)$ for every $i=1,\ldots,n$ and regard diagrams $K_i$ and $K_j$ equal if $\reid_r(K_i) \cap \reid_r(K_j) \neq \emptyset$. Note that this method is superior to just compare minimal representatives of a diagram after Reidemeister moves, as it compares the whole space of Reidemseiter between diagrams. In practice, this roughly doubles the number of Reidemeister moves performed, since the search proceeds from both knotoids simultaneously until the two exploration trees meet at a common intermediate diagram.

We perform such searches consecutively using the following parameters:

\begin{enumerate}
    \item searching the Reidemeister space with two crossing increasing II moves including flypes (reduction 752 groups),
\end{enumerate}

We now tested with three crossing increasing moves and no further reductions were being made. Thus, we conjecture at this point that diagrams were completely classified, but 752 groups of diagrams consisting only groups containing 2 or 4 knotoids could not be distinguished using invariants - they are either chiral knotoids, where invariants do not detect chirality, or, in the majority of cases, we conjecture that are non-rotatable knotoids where invariants fail to detect the non-rotatability.

We therefore continue with comparing Reidemeister spaces simultaneously, but inserting also the rotation of a knotoid into the same group. In detail, for a group $K = \{K_1, K_2,\ldots,K_n\}$, we declare two knotoids $K_i$ and $K_j$ to be equivalent if $(\reid_r(K_i) \cup \reid_r(\rot K_i)) \cap (\reid_r(K_j) \cup \reid_r(\rot K_j)) \neq \emptyset$.

We peform such searches consecutively using the following parameters:

\begin{enumerate}
    \item Searching the Reidemeister space with one crossing increasing II move including flypes (reduction 68 groups containing duplicates),
    \item Searching the Reidemeister space with two crossing increasing II moves including flypes (reduction final 28 groups containing duplicates).
\end{enumerate}
As before, we check the Reidemeister space with three crossing increasing moves and no diagrams are further reduced.

\paragraph{Step 5: Final analysis.}
For all classified knotoids, we eliminate mirror duplicates as follows. 
For each knotoid $K$, we compute the set of invariants of its mirror image $K^*$. 
If this set coincides with the invariants of $K$, then $K$ is achiral; otherwise, it is chiral. 
For each chiral pair $(K, K^*)$, we retain only one representative.

This procedure is consistent with the structure observed in the final reduction. 
Among the $28$ groups consisting of two knotoids each (in total $56$ knotoids), the groups occur in mirror pairs. 
More precisely, for every group
$K = \{K_{1}, K_{2}\},$
there exists a group
$K' = \{K'_{1}, K'_{2}\}$
such that
$
\{K_{1}^*, K_{2}^*\} = \{K'_{1}, K'_{2}\}.
$
We therefore conjecture that these $56$ knotoids form $28$ mirror pairs.

We could not find that the remaining $28$ knotoids are related by rotation, yet share all computed invariants. In analogy with classical knot theory, where mutant knots are indistinguishable by many standard invariants, we suggest that these knotoids may be regarded as \emph{mutant-type} pairs.

In the results (\autoref{tab:per_crossing} and the Appendix) we conjecture that these knotoids are non-rotatable.

\section{Results}
Our computational enumeration and classification method produced a comprehensive table of prime spherical knotoids up to 7 crossings. 

\subsection{Per-Crossing Statistics}
Table~\ref{tab:per_crossing} details the number of distinct prime spherical knotoids found at each crossing number from 0 to 7. The table provides a breakdown based on the chirality and rotatability of the knotoids.
Up to 6 crossings (see the Appendix), achiral knotoids are: $0_1$, $4_1$, $6_3$, $6_{22}$, and $6_{48}$, rotatable knotoids are $0_1$, $3_1$, $4_1$, $5_1$, $5_2$, $5_{14}$, $5_{15}$, $5_{24}$, $5_{25}$, $6_1$, $6_2$, $6_3$, $6_{31}$, $6_{34}$, $6_{81}$, and $6_{82}$.
In the whole classification, we have only three fully achiral (achiral and rotatable) knotoids, namely $0_1$, $4_1$, and $6_3$, which are all knot-like knotoids that coincide with the fully achiral knots.

\begin{table}[h!]
\centering
\begin{threeparttable}
\caption{Number of prime spherical knotoids up to 7 crossings.}
\label{tab:per_crossing}

\begin{tabular}{@{}lcccccc@{}}
\toprule
\textbf{Crossings} & \textbf{Total} & \multicolumn{2}{c}{\textbf{Chiral}} & \multicolumn{2}{c}{\textbf{Rotatable}} & \textbf{Possible} \\
\cmidrule(lr){3-4}\cmidrule(lr){5-6}
& & \textbf{Yes} & \textbf{No} & \textbf{Yes} & \textbf{No}\tnote{\textcolor{MidnightBlue}{*}} & \textbf{duplicates}\\
\midrule
0 & 1 & 0 & 1 & 1 & 0 & 0\\
1 & 0 & 0 & 0 & 0 & 0 & 0\\
2 & 1 & 1 & 0 & 0 & 1 & 0\\
3 & 2 & 2 & 0 & 1 & 1 & 0\\
4 & 8 & 7 & 1 & 1 & 7 & 0\\
5 & 25  & 25  & 0  & 6 & 19 & 0\\  
6 & 82  & 79  & 3  & 7 & 75 & 0\\  
7 & 308 & 280 & 28 & 23 & 285 & 14\\
\midrule
\textbf{Total} & \textbf{427} & \textbf{394} & \textbf{33} & \textbf{37} & \textbf{390}  & \textbf{14} \\
\bottomrule
\end{tabular}


\begin{tablenotes}
\footnotesize
\item[\textcolor{MidnightBlue}{*}] 
Conjectured.
\end{tablenotes}

\end{threeparttable}
\end{table}

\subsection{Performance of Invariants}
The successful classification of the knotoids relies on the distinguishing power of the invariants. Table~\ref{tab:invariants} summarizes the performance of each polynomial in identifying unique knotoids among the final set of candidates. The Yamada polynomial of the double-sided closure proved to be the most powerful invariant in our study, closely followed by the Arrow and Kauffman bracket polynomial.

\begin{table}[h!]
\centering
\caption{Distinguishing power of the computed invariants on the final set of 427 knotoids.}
\label{tab:invariants}
\begin{tabular}{@{}lcc@{}}
\toprule
\textbf{Invariant} & \textbf{Unique} & \textbf{Non-unique} \\ \midrule
Kauffman Bracket & 370 & 57 \\
Arrow Polynomial & 382 & 45 \\
Mock Alexander Polynomial & 342 & 85 \\
Affine Index Polynomial & 0 & 427 \\
Yamada Polynomial (of closure) & 397 & 30 \\ \bottomrule
\end{tabular}
\end{table}

\subsection{Knotoid Table Images}
\begin{figure}[ht]
\centering
\setlength{\tabcolsep}{8pt}
\renewcommand{\arraystretch}{1.2}

\begin{tabular}{cccc}
\includegraphics[page=1,width=0.22\textwidth]{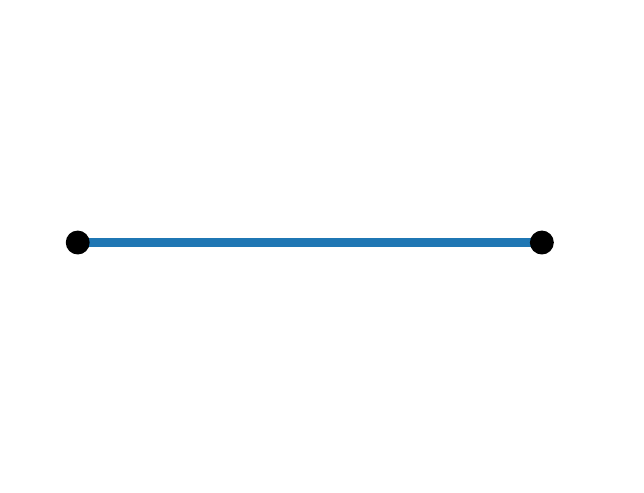} &
\includegraphics[page=2,width=0.22\textwidth]{knotoids.pdf} &
\includegraphics[page=3,width=0.22\textwidth]{knotoids.pdf} &
\includegraphics[page=4,width=0.22\textwidth]{knotoids.pdf} \\
$K0_1$ & $K2_1$ & $K3_1$ & $K3_2$ \\
\includegraphics[page=5,width=0.22\textwidth]{knotoids.pdf} &
\includegraphics[page=6,width=0.22\textwidth]{knotoids.pdf} &
\includegraphics[page=7,width=0.22\textwidth]{knotoids.pdf} &
\includegraphics[page=8,width=0.22\textwidth]{knotoids.pdf} \\
$K4_1$ & $K4_2$ & $K4_3$ & $K4_4$ \\
\includegraphics[page=9,width=0.22\textwidth]{knotoids.pdf} &
\includegraphics[page=10,width=0.22\textwidth]{knotoids.pdf} &
\includegraphics[page=11,width=0.22\textwidth]{knotoids.pdf} &
\includegraphics[page=12,width=0.22\textwidth]{knotoids.pdf} \\
$K4_5$ & $K4_6$ & $K4_7$ & $K4_8$
\end{tabular}

\caption{Knotoid diagrams up to 4 crossings.} \label{fig:4c}
\end{figure}

The complete enumeration of knotoid diagrams is extensive. In \autoref{fig:4c} we present knotoids up to 4 crossings and in \autoref{fig:5c} we present knotoids up to 5 crossings. The full table can be found in the Appendix and at the GitHub repository \cite{gabrovsek_knotoids_github}.
 
\begin{figure}[ht]
\centering
\setlength{\tabcolsep}{8pt}
\renewcommand{\arraystretch}{1.2}

\begin{tabular}{cccc}
\includegraphics[page=13,width=0.22\textwidth]{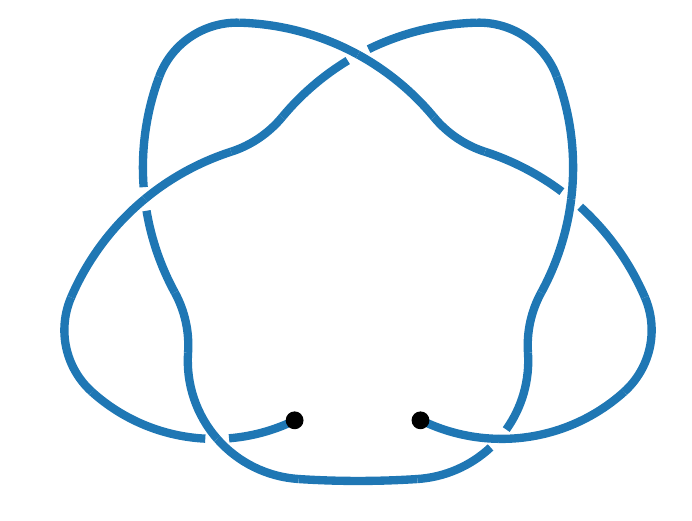} &
\includegraphics[page=14,width=0.22\textwidth]{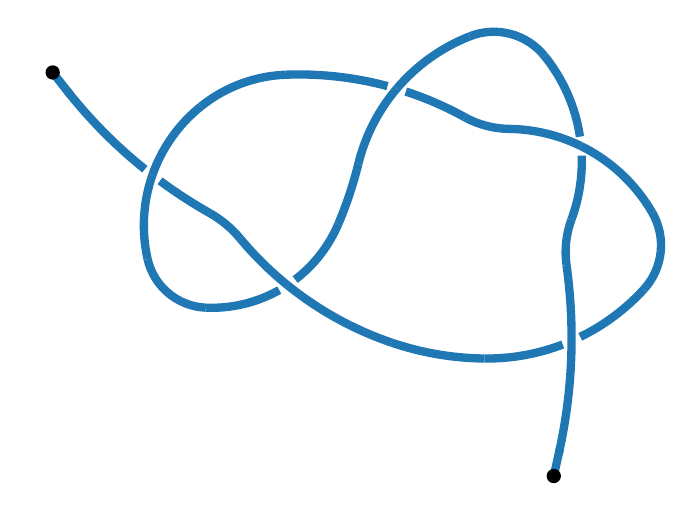} &
\includegraphics[page=15,width=0.22\textwidth]{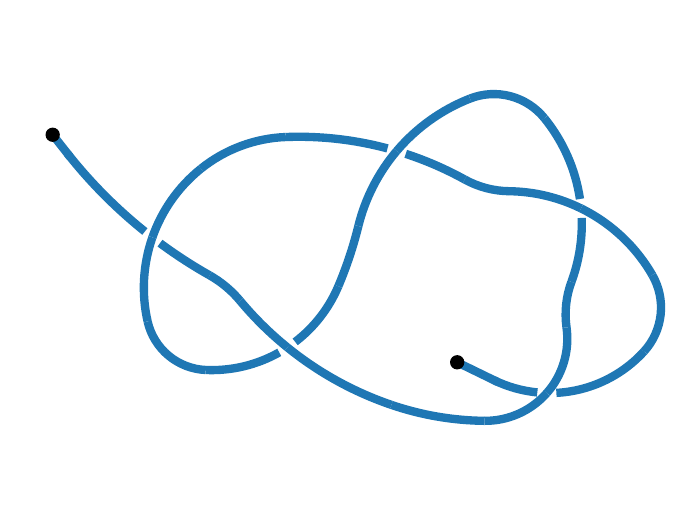} &
\includegraphics[page=16,width=0.22\textwidth]{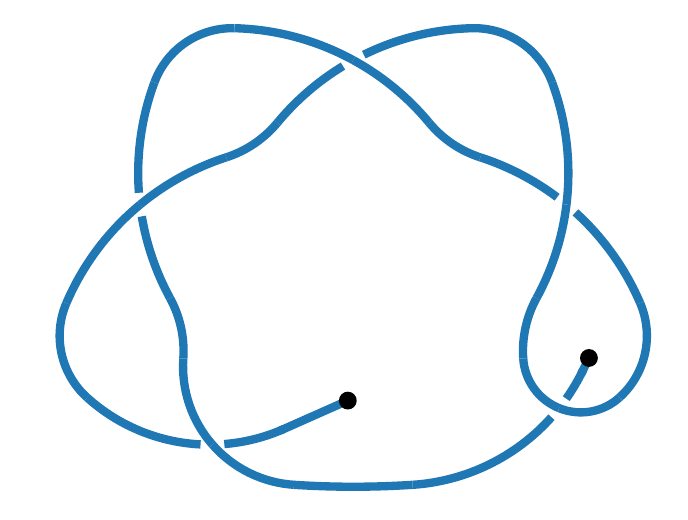} \\[-0.7em]
$K5_1$ & $K5_2$ & $K5_3$ & $K5_4$ \\
\includegraphics[page=17,width=0.22\textwidth]{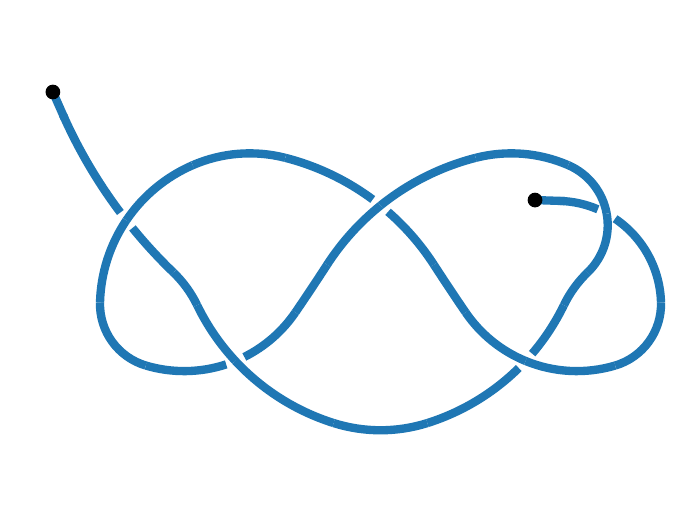} &
\includegraphics[page=18,width=0.22\textwidth]{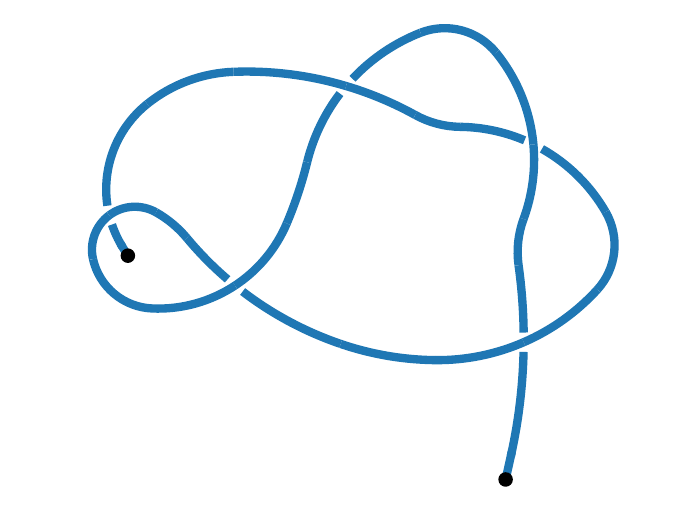} &
\includegraphics[page=19,width=0.22\textwidth]{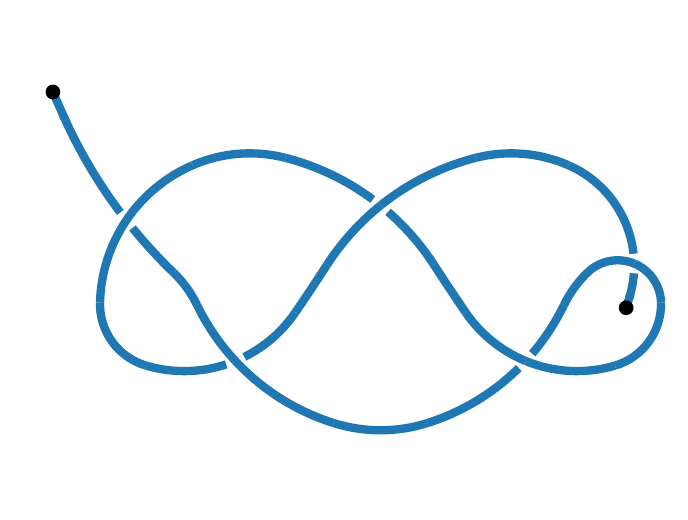} &
\includegraphics[page=20,width=0.22\textwidth]{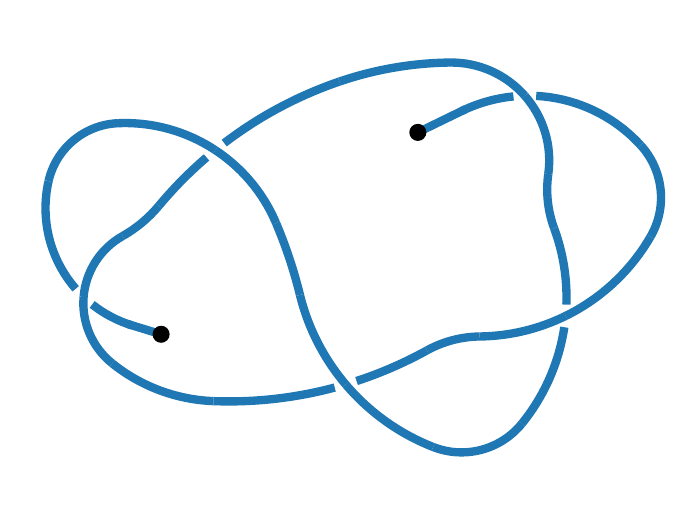} \\[-0.7em]
$K5_5$ & $K5_6$ & $K5_7$ & $K5_8$ \\
\includegraphics[page=21,width=0.22\textwidth]{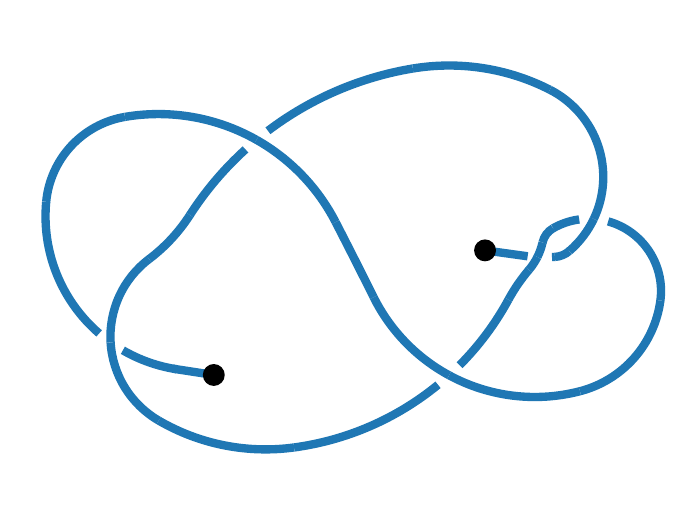} &
\includegraphics[page=22,width=0.22\textwidth]{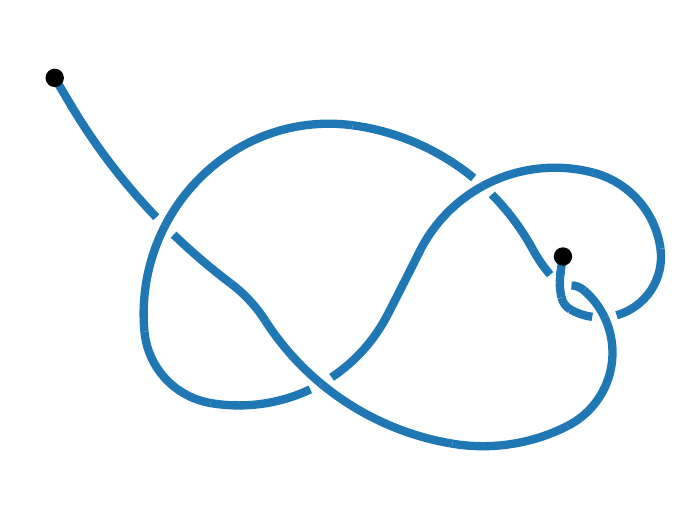} &
\includegraphics[page=23,width=0.22\textwidth]{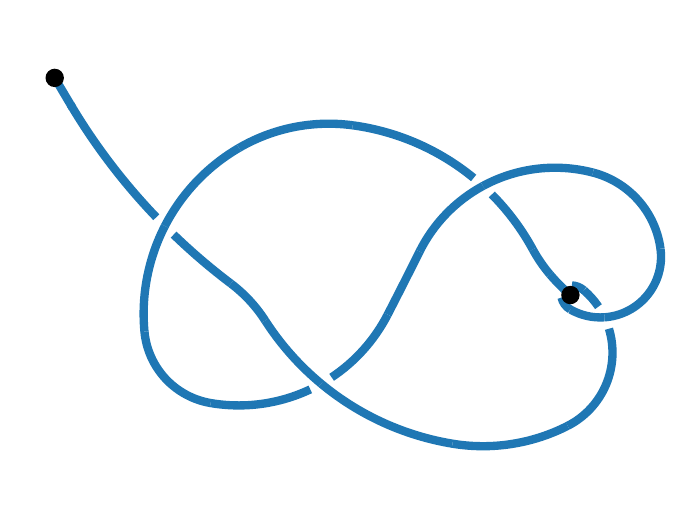} &
\includegraphics[page=24,width=0.22\textwidth]{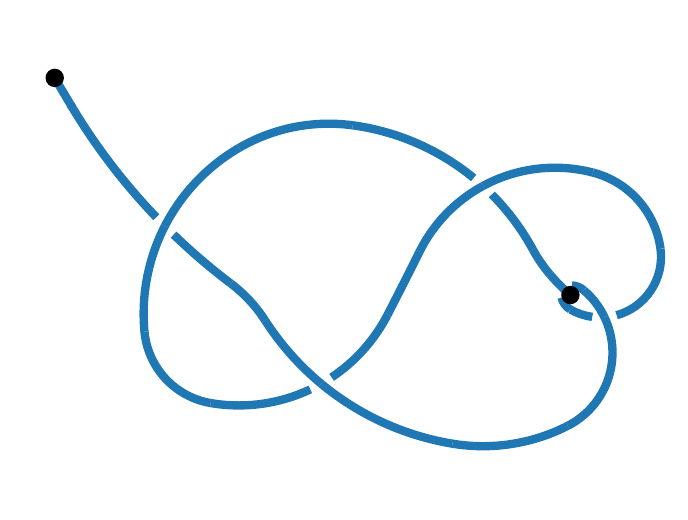} \\[-0.7em]
$K5_9$ & $K5_{10}$ & $K5_{11}$ & $K5_{12}$ \\
\includegraphics[page=25,width=0.22\textwidth]{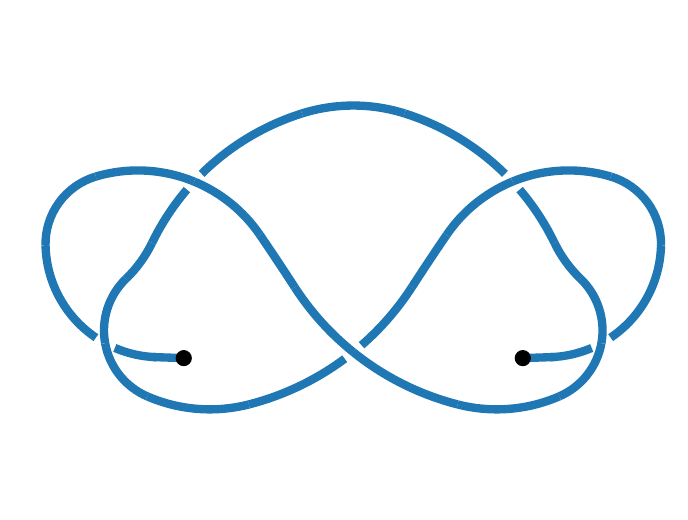} &
\includegraphics[page=26,width=0.22\textwidth]{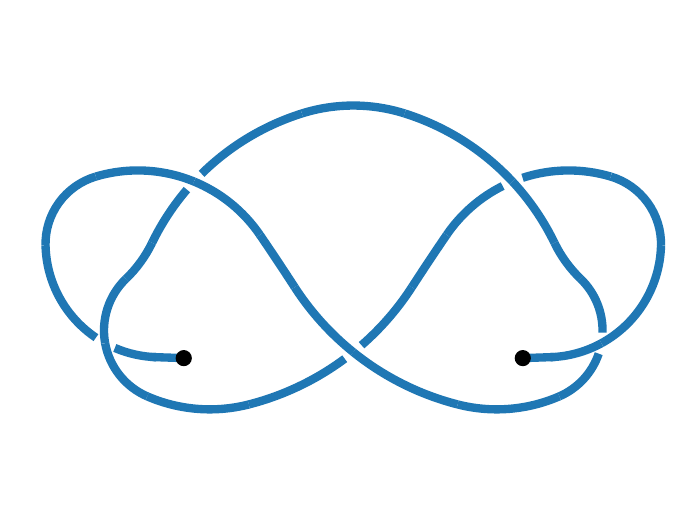} &
\includegraphics[page=27,width=0.22\textwidth]{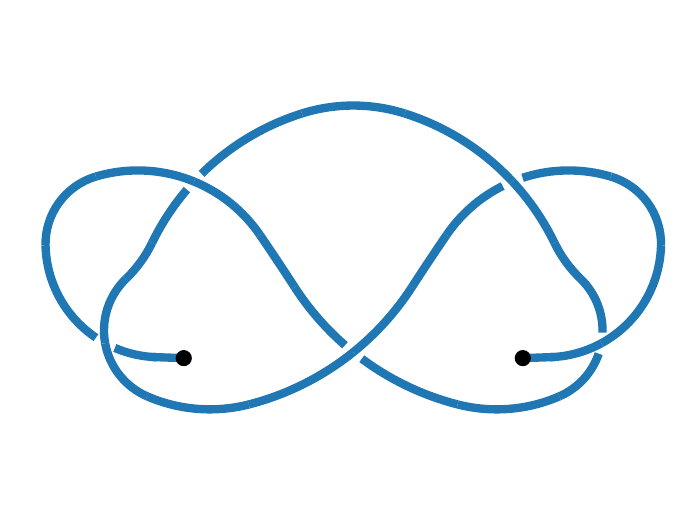} &
\includegraphics[page=28,width=0.22\textwidth]{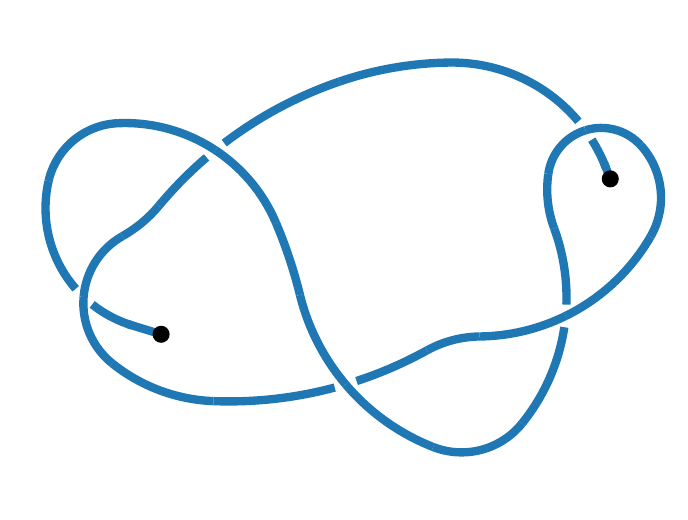} \\[-0.7em]
$K5_{13}$ & $K5_{14}$ & $K5_{15}$ & $K5_{16}$ \\
\includegraphics[page=29,width=0.22\textwidth]{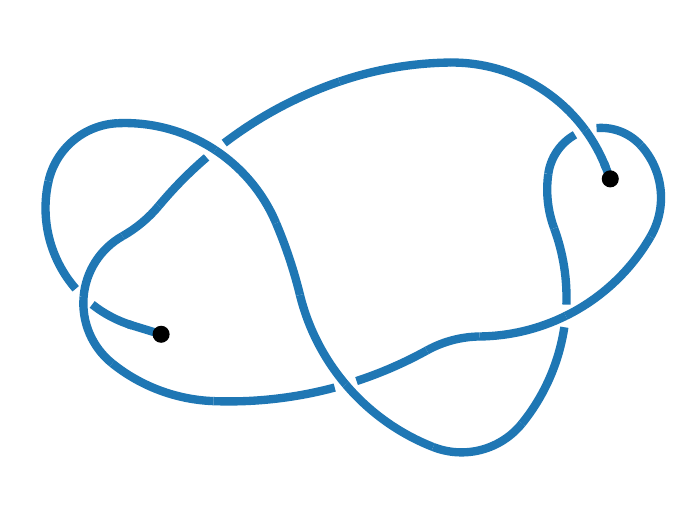} &
\includegraphics[page=30,width=0.22\textwidth]{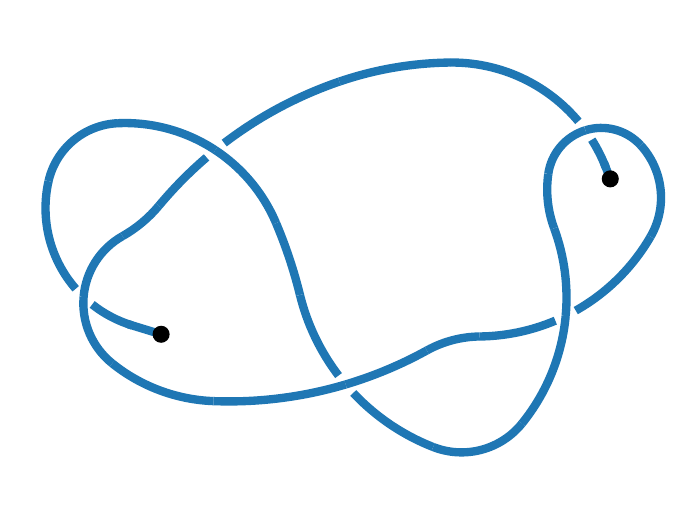} &
\includegraphics[page=31,width=0.22\textwidth]{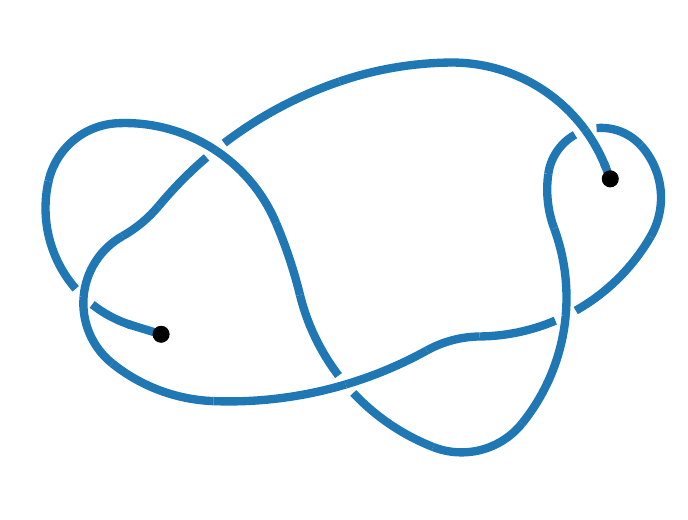} &
\includegraphics[page=32,width=0.22\textwidth]{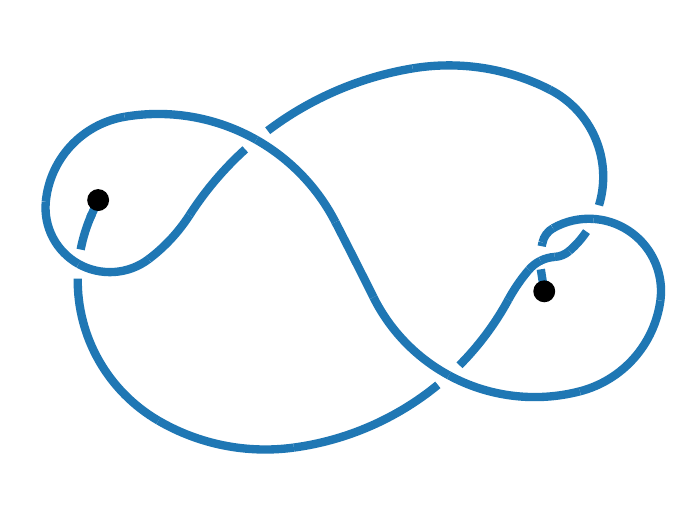} \\[-0.7em]
$K5_{17}$ & $K5_{18}$ & $K5_{19}$ & $K5_{20}$ \\
\includegraphics[page=33,width=0.22\textwidth]{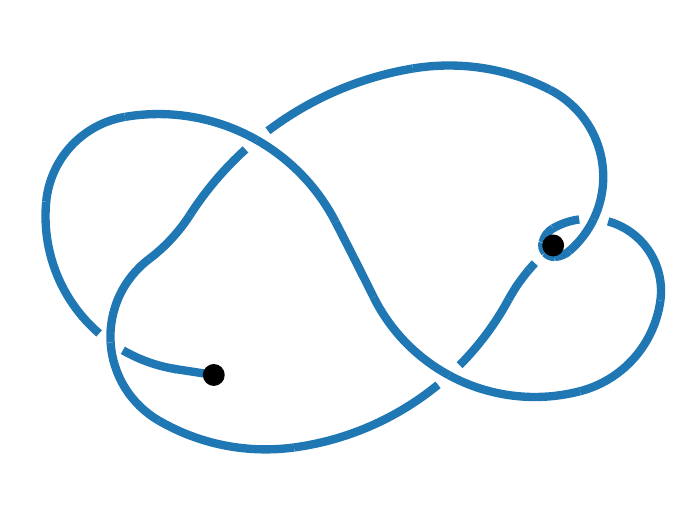} &
\includegraphics[page=34,width=0.22\textwidth]{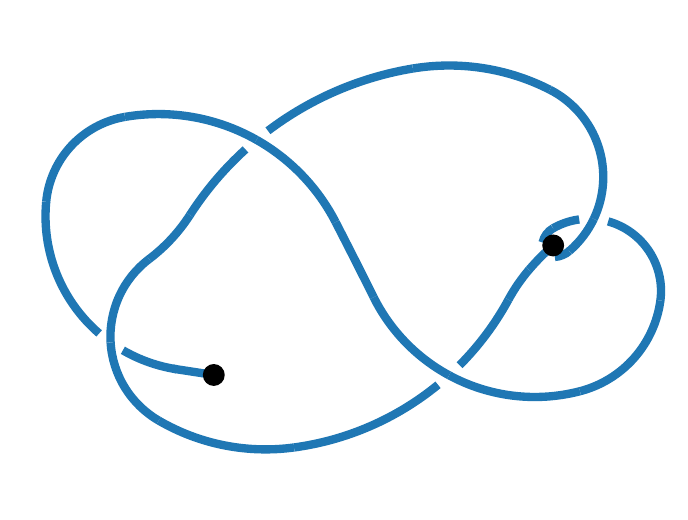} &
\includegraphics[page=35,width=0.22\textwidth]{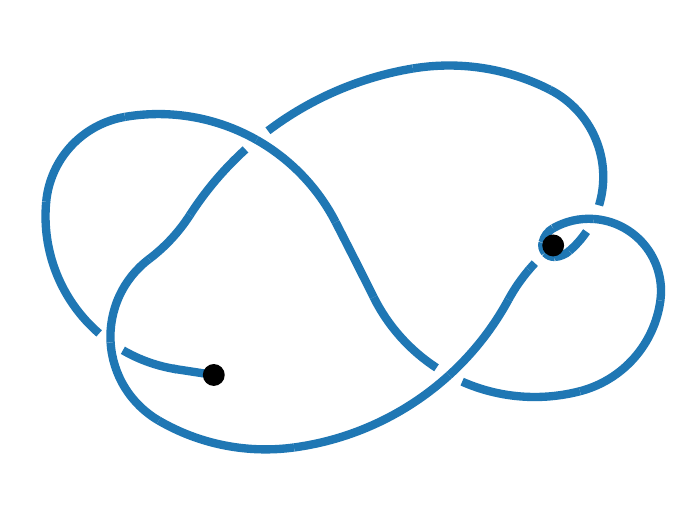} &
\includegraphics[page=36,width=0.22\textwidth]{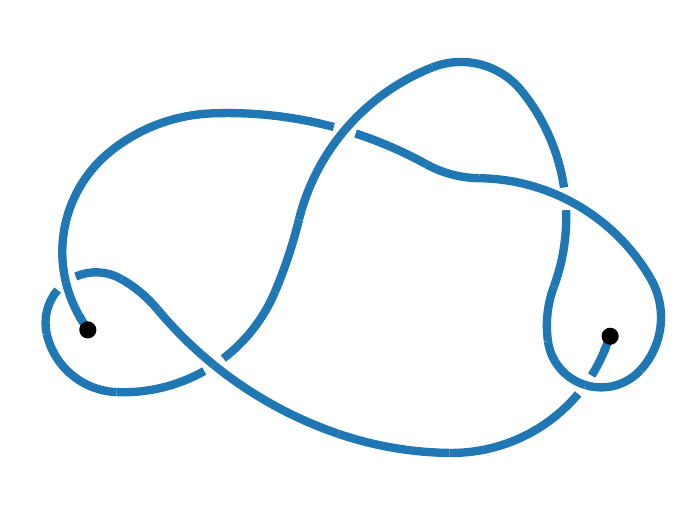} \\[-0.7em]
$K5_{21}$ & $K5_{22}$ & $K5_{23}$ & $K5_{24}$ \\
\includegraphics[page=38,width=0.22\textwidth]{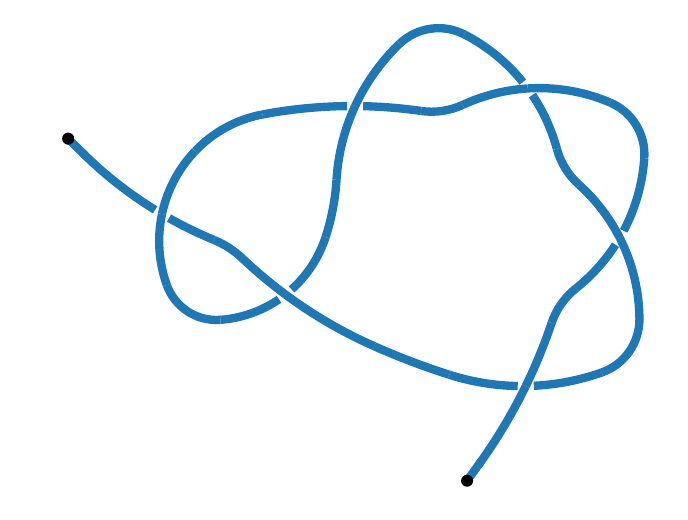} & & & \\[-0.7em]
$K5_{25}$ & & &
\end{tabular}

\caption{Knotoid diagrams with 5 crossings.} \label{fig:5c}
\end{figure}

In summary, this work establishes a classification of prime spherical knotoids up to seven crossings, identifying a total of 427 distinct types with 14 possible duplicates. This tabulation serves as a foundational resource for future research in knotoid theory and its practical applications in analyzing the topological structures of proteins.



\subsection*{Acknowledgments}
B. Gabrovšek was financially supported by the Slovenian Research and Innovation Agency grants J1-4031, N1-0278, and program P1-0292. P. Cavicchioli was financially supported by the Slovenian Research and Innovation Agency grant J1-4031 and program P1-0292. The authors declare that they have no conflict of interest.

\bibliographystyle{abbrv}
\bibliography{biblio}

\newgeometry{left=1cm,right=1cm, top=1cm, bottom=2cm}
\input{table-insert}

\end{document}

%% file: table-insert.tex
\vspace*{\fill}
\begin{center}
{\Huge Appendix}
\end{center}
\vspace*{\fill}

\section*{Number of crossings: 0}

\noindent{\color{gray!40}\rule{\textwidth}{0.4pt}}
\vspace{0.9\baselineskip}
\noindent \begin{minipage}[t]{0.25\textwidth}
\vspace{0pt}
\centering
\includegraphics[page=1,width=\linewidth]{knotoids.pdf}
\end{minipage}
\hfill
\begin{minipage}[t]{0.73\textwidth}
\vspace{0pt}
\raggedright
\textbf{Name:} {\large{$\mathbf{K0_{1}}$}} (achiral, rotatable) \\ \textbf{EM:} {\small\texttt{[0],[0]}} \\ \textbf{EM:} {\small\texttt{(B0, A0)}} \\ \textbf{Kauffman bracket:} {\scriptsize $1$} \\ \textbf{Arrow:} {\scriptsize $1$} \\ \textbf{Mock:} {\scriptsize $1$} \\ \textbf{Affine:} {\scriptsize $0$} \\ \textbf{Yamada:} {\scriptsize $-A^{4} - A^{3} - 2A^{2} - A - 1$}
\end{minipage}

\noindent{\color{gray!40}\rule{\textwidth}{0.4pt}}
\vspace{0.9\baselineskip}

\section*{Number of crossings: 2}

\noindent{\color{gray!40}\rule{\textwidth}{0.4pt}}
\vspace{0.9\baselineskip}
\noindent \begin{minipage}[t]{0.25\textwidth}
\vspace{0pt}
\centering
\includegraphics[page=2,width=\linewidth]{knotoids.pdf}
\end{minipage}
\hfill
\begin{minipage}[t]{0.73\textwidth}
\vspace{0pt}
\raggedright
\textbf{Name:} {\large{$\mathbf{K2_{1}}$}} (chiral, non-rotatable$^{*}$) \\ \textbf{EM:} {\small\texttt{[0],[0,1,2,3],[1,3,4,2],[4]}} \\ \textbf{EM:} {\small\texttt{(B0, A0C0C3C1, B1B3D0B2, C2)}} \\ \textbf{Kauffman bracket:} {\scriptsize $A^{8} + A^{6} - A^{2}$} \\ \textbf{Arrow:} {\scriptsize $A^{-4} + L_1/A^{6} - L_1/A^{10}$} \\ \textbf{Mock:} {\scriptsize $w^{2} + w - 1/w$} \\ \textbf{Affine:} {\scriptsize $t - 2 + 1/t$} \\ \textbf{Yamada:} {\scriptsize $-A^{12} - A^{11} - A^{10} - A^{9} - A^{8} - A^{6} - A^{4} + 1$}
\end{minipage}

\noindent{\color{gray!40}\rule{\textwidth}{0.4pt}}
\vspace{0.9\baselineskip}

\section*{Number of crossings: 3}

\noindent{\color{gray!40}\rule{\textwidth}{0.4pt}}
\vspace{0.9\baselineskip}
\noindent \begin{minipage}[t]{0.25\textwidth}
\vspace{0pt}
\centering
\includegraphics[page=3,width=\linewidth]{knotoids.pdf}
\end{minipage}
\hfill
\begin{minipage}[t]{0.73\textwidth}
\vspace{0pt}
\raggedright
\textbf{Name:} {\large{$\mathbf{K3_{1}}$}} (chiral, rotatable) \\ \textbf{EM:} {\small\texttt{[0],[0,1,2,3],[1,4,5,2],[3,5,4,6],[6]}} \\ \textbf{EM:} {\small\texttt{(B0, A0C0C3D0, B1D2D1B2, B3C2C1E0, D3)}} \\ \textbf{Kauffman bracket:} {\scriptsize $A^{14} + A^{6} - A^{2}$} \\ \textbf{Arrow:} {\scriptsize $A^{-4} + A^{-12} - 1/A^{16}$} \\ \textbf{Mock:} {\scriptsize $w^{2} - 1 + w^{-2}$} \\ \textbf{Affine:} {\scriptsize $0$} \\ \textbf{Yamada:} {\scriptsize $A^{13} + A^{12} + 2A^{11} + 2A^{10} + 2A^{9} + A^{8} + A^{7} - A^{6} - A^{5} - 2A^{4} - A^{3} + 1$}
\end{minipage}

\noindent{\color{gray!40}\rule{\textwidth}{0.4pt}}
\vspace{0.9\baselineskip}
\noindent \begin{minipage}[t]{0.25\textwidth}
\vspace{0pt}
\centering
\includegraphics[page=4,width=\linewidth]{knotoids.pdf}
\end{minipage}
\hfill
\begin{minipage}[t]{0.73\textwidth}
\vspace{0pt}
\raggedright
\textbf{Name:} {\large{$\mathbf{K3_{2}}$}} (chiral, non-rotatable$^{*}$) \\ \textbf{EM:} {\small\texttt{[0],[0,1,2,3],[1,4,5,2],[3,5,6,4],[6]}} \\ \textbf{EM:} {\small\texttt{(B0, A0C0C3D0, B1D3D1B2, B3C2E0C1, D2)}} \\ \textbf{Kauffman bracket:} {\scriptsize $-A^{11} + A^{7} + A^{5} - A^{3} - A$} \\ \textbf{Arrow:} {\scriptsize $A^{8} - A^{4} - A^{2}L_1 + 1 + L_1/A^{2}$} \\ \textbf{Mock:} {\scriptsize $-w^{2} - w + 2 + 1/w$} \\ \textbf{Affine:} {\scriptsize $-t + 2 - 1/t$} \\ \textbf{Yamada:} {\scriptsize $A^{15} - A^{13} - A^{11} - A^{10} - A^{8} - A^{6} - A^{3} - 1$}
\end{minipage}

\noindent{\color{gray!40}\rule{\textwidth}{0.4pt}}
\vspace{0.9\baselineskip}

\section*{Number of crossings: 4}

\noindent{\color{gray!40}\rule{\textwidth}{0.4pt}}
\vspace{0.9\baselineskip}
\noindent \begin{minipage}[t]{0.25\textwidth}
\vspace{0pt}
\centering
\includegraphics[page=5,width=\linewidth]{knotoids.pdf}
\end{minipage}
\hfill
\begin{minipage}[t]{0.73\textwidth}
\vspace{0pt}
\raggedright
\textbf{Name:} {\large{$\mathbf{K4_{1}}$}} (achiral, rotatable) \\ \textbf{EM:} {\small\texttt{[0],[0,1,2,3],[1,4,5,2],[3,5,6,7],[4,8,7,6],[8]}} \\ \textbf{EM:} {\small\texttt{(B0, A0C0C3D0, B1E0D1B2, B3C2E3E2, C1F0D3D2, E1)}} \\ \textbf{Kauffman bracket:} {\scriptsize $-A^{17} + A^{13} - A^{9} + A^{5} - A$} \\ \textbf{Arrow:} {\scriptsize $A^{8} - A^{4} + 1 - 1/A^{4} + A^{-8}$} \\ \textbf{Mock:} {\scriptsize $-w^{2} + 3 - 1/w^{2}$} \\ \textbf{Affine:} {\scriptsize $0$} \\ \textbf{Yamada:} {\scriptsize $-A^{16} + A^{12} - A^{10} - A^{9} - 2A^{8} - A^{7} - A^{6} + A^{4} - 1$}
\end{minipage}

\noindent{\color{gray!40}\rule{\textwidth}{0.4pt}}
\vspace{0.9\baselineskip}
\noindent \begin{minipage}[t]{0.25\textwidth}
\vspace{0pt}
\centering
\includegraphics[page=6,width=\linewidth]{knotoids.pdf}
\end{minipage}
\hfill
\begin{minipage}[t]{0.73\textwidth}
\vspace{0pt}
\raggedright
\textbf{Name:} {\large{$\mathbf{K4_{2}}$}} (chiral, non-rotatable$^{*}$) \\ \textbf{EM:} {\small\texttt{[0],[0,1,2,3],[1,4,5,2],[6,7,8,3],[4,8,7,5],[6]}} \\ \textbf{EM:} {\small\texttt{(B0, A0C0C3D3, B1E0E3B2, F0E2E1B3, C1D2D1C2, D0)}} \\ \textbf{Kauffman bracket:} {\scriptsize $A^{16} + A^{14} - A^{10} + A^{6} - A^{2}$} \\ \textbf{Arrow:} {\scriptsize $A^{-8} + L_1/A^{10} - L_1/A^{14} + L_1/A^{18} - L_1/A^{22}$} \\ \textbf{Mock:} {\scriptsize $w^{4} + w^{3} - w + 1/w - 1/w^{3}$} \\ \textbf{Affine:} {\scriptsize $2t - 4 + 2/t$} \\ \textbf{Yamada:} {\scriptsize $A^{17} + A^{16} + A^{15} + A^{14} + 2A^{13} + A^{12} + A^{11} + A^{10} - A^{8} - A^{6} - A^{4} - A + 1$}
\end{minipage}

\noindent{\color{gray!40}\rule{\textwidth}{0.4pt}}
\vspace{0.9\baselineskip}
\noindent \begin{minipage}[t]{0.25\textwidth}
\vspace{0pt}
\centering
\includegraphics[page=7,width=\linewidth]{knotoids.pdf}
\end{minipage}
\hfill
\begin{minipage}[t]{0.73\textwidth}
\vspace{0pt}
\raggedright
\textbf{Name:} {\large{$\mathbf{K4_{3}}$}} (chiral, non-rotatable$^{*}$) \\ \textbf{EM:} {\small\texttt{[0],[0,1,2,3],[1,4,5,2],[3,5,6,7],[7,6,8,4],[8]}} \\ \textbf{EM:} {\small\texttt{(B0, A0C0C3D0, B1E3D1B2, B3C2E1E0, D3D2F0C1, E2)}} \\ \textbf{Kauffman bracket:} {\scriptsize $A^{14} - A^{10} + 2A^{6} + A^{4} - A^{2} - 1$} \\ \textbf{Arrow:} {\scriptsize $A^{-4} - 1/A^{8} + 2/A^{12} + L_1/A^{14} - 1/A^{16} - L_1/A^{18}$} \\ \textbf{Mock:} {\scriptsize $2w^{2} + w - 2 - 1/w + w^{-2}$} \\ \textbf{Affine:} {\scriptsize $t - 2 + 1/t$} \\ \textbf{Yamada:} {\scriptsize $2A^{17} + A^{16} + 2A^{14} + A^{13} + 2A^{11} + A^{9} - A^{8} - 2A^{5} - A^{2} + 1$}
\end{minipage}

\noindent{\color{gray!40}\rule{\textwidth}{0.4pt}}
\vspace{0.9\baselineskip}
\noindent \begin{minipage}[t]{0.25\textwidth}
\vspace{0pt}
\centering
\includegraphics[page=8,width=\linewidth]{knotoids.pdf}
\end{minipage}
\hfill
\begin{minipage}[t]{0.73\textwidth}
\vspace{0pt}
\raggedright
\textbf{Name:} {\large{$\mathbf{K4_{4}}$}} (chiral, non-rotatable$^{*}$) \\ \textbf{EM:} {\small\texttt{[0],[0,1,2,3],[1,4,5,2],[3,6,7,8],[4,8,6,5],[7]}} \\ \textbf{EM:} {\small\texttt{(B0, A0C0C3D0, B1E0E3B2, B3E2F0E1, C1D3D1C2, D2)}} \\ \textbf{Kauffman bracket:} {\scriptsize $A^{15} - A^{11} - A^{9} + A^{7} + A^{5} - A^{3} - A$} \\ \textbf{Arrow:} {\scriptsize $-A^{18}L_1 + A^{14}L_1 + A^{12} - A^{10}L_1 - A^{8} + A^{6}L_1 + A^{4}$} \\ \textbf{Mock:} {\scriptsize $-2w - 1 + 2/w + 2/w^{2}$} \\ \textbf{Affine:} {\scriptsize $-2t + 4 - 2/t$} \\ \textbf{Yamada:} {\scriptsize $A^{18} - A^{16} - A^{13} + A^{12} + A^{11} + A^{9} - A^{8} - 2A^{6} - A^{5} - 2A^{3} - A^{2} - 1$}
\end{minipage}

\noindent{\color{gray!40}\rule{\textwidth}{0.4pt}}
\vspace{0.9\baselineskip}
\noindent \begin{minipage}[t]{0.25\textwidth}
\vspace{0pt}
\centering
\includegraphics[page=9,width=\linewidth]{knotoids.pdf}
\end{minipage}
\hfill
\begin{minipage}[t]{0.73\textwidth}
\vspace{0pt}
\raggedright
\textbf{Name:} {\large{$\mathbf{K4_{5}}$}} (chiral, non-rotatable$^{*}$) \\ \textbf{EM:} {\small\texttt{[0],[0,1,2,3],[1,4,5,2],[3,6,7,4],[5,7,8,6],[8]}} \\ \textbf{EM:} {\small\texttt{(B0, A0C0C3D0, B1D3E0B2, B3E3E1C1, C2D2F0D1, E2)}} \\ \textbf{Kauffman bracket:} {\scriptsize $-A^{10} + A^{8} + 2A^{6} + A^{4} - A^{2} - 1$} \\ \textbf{Arrow:} {\scriptsize $-L_1/A^{2} + A^{-4} + 2L_1/A^{6} + L_2/A^{8} - L_1/A^{10} - L_2/A^{12}$} \\ \textbf{Mock:} {\scriptsize $-w^{3} + 2w + 1 - 1/w$} \\ \textbf{Affine:} {\scriptsize $t^{2} - 2 + t^{-2}$} \\ \textbf{Yamada:} {\scriptsize $-A^{18} - 2A^{15} - A^{12} + A^{11} - A^{10} - 2A^{8} - A^{7} - A^{5} + A^{3} + 1$}
\end{minipage}

\noindent{\color{gray!40}\rule{\textwidth}{0.4pt}}
\vspace{0.9\baselineskip}
\noindent \begin{minipage}[t]{0.25\textwidth}
\vspace{0pt}
\centering
\includegraphics[page=10,width=\linewidth]{knotoids.pdf}
\end{minipage}
\hfill
\begin{minipage}[t]{0.73\textwidth}
\vspace{0pt}
\raggedright
\textbf{Name:} {\large{$\mathbf{K4_{6}}$}} (chiral, non-rotatable$^{*}$) \\ \textbf{EM:} {\small\texttt{[0],[0,1,2,3],[1,4,5,2],[3,6,7,4],[7,8,6,5],[8]}} \\ \textbf{EM:} {\small\texttt{(B0, A0C0C3D0, B1D3E3B2, B3E2E0C1, D2F0D1C2, E1)}} \\ \textbf{Kauffman bracket:} {\scriptsize $-A^{12} + A^{8} + 2A^{6} - A^{2}$} \\ \textbf{Arrow:} {\scriptsize $-A^{6}L_1 + A^{2}L_1 + L_2 + 1 - L_2/A^{4}$} \\ \textbf{Mock:} {\scriptsize $-w^{3} - w^{2} + w + 2$} \\ \textbf{Affine:} {\scriptsize $t^{2} - t - 1/t + t^{-2}$} \\ \textbf{Yamada:} {\scriptsize $-A^{16} - A^{15} + A^{14} + A^{11} - A^{10} - 2A^{8} - A^{7} - 2A^{6} - 2A^{5} + A + 1$}
\end{minipage}

\noindent{\color{gray!40}\rule{\textwidth}{0.4pt}}
\vspace{0.9\baselineskip}
\noindent \begin{minipage}[t]{0.25\textwidth}
\vspace{0pt}
\centering
\includegraphics[page=11,width=\linewidth]{knotoids.pdf}
\end{minipage}
\hfill
\begin{minipage}[t]{0.73\textwidth}
\vspace{0pt}
\raggedright
\textbf{Name:} {\large{$\mathbf{K4_{7}}$}} (chiral, non-rotatable$^{*}$) \\ \textbf{EM:} {\small\texttt{[0],[0,1,2,3],[1,4,5,2],[6,7,4,3],[5,7,8,6],[8]}} \\ \textbf{EM:} {\small\texttt{(B0, A0C0C3D3, B1D2E0B2, E3E1C1B3, C2D1F0D0, E2)}} \\ \textbf{Kauffman bracket:} {\scriptsize $A^{10} + A^{8} + A^{6} - A^{4} - A^{2}$} \\ \textbf{Arrow:} {\scriptsize $A^{-8} + L_1/A^{10} + L_2/A^{12} - L_1/A^{14} - L_2/A^{16}$} \\ \textbf{Mock:} {\scriptsize $w^{4} + w^{3} - w - 1 + w^{-2}$} \\ \textbf{Affine:} {\scriptsize $t^{2} + t - 4 + 1/t + t^{-2}$} \\ \textbf{Yamada:} {\scriptsize $-A^{18} - A^{17} - A^{16} - A^{15} - A^{14} - A^{8} - A^{6} - A^{4} + A^{2} + 1$}
\end{minipage}

\noindent{\color{gray!40}\rule{\textwidth}{0.4pt}}
\vspace{0.9\baselineskip}
\noindent \begin{minipage}[t]{0.25\textwidth}
\vspace{0pt}
\centering
\includegraphics[page=12,width=\linewidth]{knotoids.pdf}
\end{minipage}
\hfill
\begin{minipage}[t]{0.73\textwidth}
\vspace{0pt}
\raggedright
\textbf{Name:} {\large{$\mathbf{K4_{8}}$}} (chiral, non-rotatable$^{*}$) \\ \textbf{EM:} {\small\texttt{[0],[0,1,2,3],[1,4,5,2],[6,7,4,3],[7,8,6,5],[8]}} \\ \textbf{EM:} {\small\texttt{(B0, A0C0C3D3, B1D2E3B2, E2E0C1B3, D1F0D0C2, E1)}} \\ \textbf{Kauffman bracket:} {\scriptsize $-2A^{11} - A^{9} + A^{7} + 2A^{5} - A$} \\ \textbf{Arrow:} {\scriptsize $L_2/A^{4} + A^{-4} + L_1/A^{6} - L_2/A^{8} - 2L_1/A^{10} + L_1/A^{14}$} \\ \textbf{Mock:} {\scriptsize $w^{2} + w - 1 - 2/w + w^{-2} + w^{-3}$} \\ \textbf{Affine:} {\scriptsize $t^{2} - 2 + t^{-2}$} \\ \textbf{Yamada:} {\scriptsize $A^{19} + A^{18} + A^{15} - A^{14} - A^{13} + A^{12} + 2A^{9} + A^{8} + 2A^{7} + A^{5} - 2A^{3} + A^{2} - 1$}
\end{minipage}

\noindent{\color{gray!40}\rule{\textwidth}{0.4pt}}
\vspace{0.9\baselineskip}

\section*{Number of crossings: 5}

\noindent{\color{gray!40}\rule{\textwidth}{0.4pt}}
\vspace{0.9\baselineskip}
\noindent \begin{minipage}[t]{0.25\textwidth}
\vspace{0pt}
\centering
\includegraphics[page=13,width=\linewidth]{knotoids.pdf}
\end{minipage}
\hfill
\begin{minipage}[t]{0.73\textwidth}
\vspace{0pt}
\raggedright
\textbf{Name:} {\large{$\mathbf{K5_{1}}$}} (chiral, rotatable) \\ \textbf{EM:} {\small\texttt{[0],[0,1,2,3],[1,4,5,2],[3,6,7,8],[4,9,10,5],[6,10,9,7],[8]}} \\ \textbf{EM:} {\small\texttt{(B0, A0C0C3D0, B1E0E3B2, B3F0F3G0, C1F2F1C2, D1E2E1D2, D3)}} \\ \textbf{Kauffman bracket:} {\scriptsize $A^{22} + A^{14} - A^{10} + A^{6} - A^{2}$} \\ \textbf{Arrow:} {\scriptsize $A^{-8} + A^{-16} - 1/A^{20} + A^{-24} - 1/A^{28}$} \\ \textbf{Mock:} {\scriptsize $w^{4} - w^{2} + 1 - 1/w^{2} + w^{-4}$} \\ \textbf{Affine:} {\scriptsize $0$} \\ \textbf{Yamada:} {\scriptsize $A^{19} + A^{18} + 2A^{17} + 2A^{16} + 2A^{15} + A^{14} + A^{13} - A^{8} - A^{7} - 2A^{6} - A^{5} - A^{4} + A^{2} + 1$}
\end{minipage}

\noindent{\color{gray!40}\rule{\textwidth}{0.4pt}}
\vspace{0.9\baselineskip}
\noindent \begin{minipage}[t]{0.25\textwidth}
\vspace{0pt}
\centering
\includegraphics[page=14,width=\linewidth]{knotoids.pdf}
\end{minipage}
\hfill
\begin{minipage}[t]{0.73\textwidth}
\vspace{0pt}
\raggedright
\textbf{Name:} {\large{$\mathbf{K5_{2}}$}} (chiral, rotatable) \\ \textbf{EM:} {\small\texttt{[0],[0,1,2,3],[1,4,5,2],[3,5,6,7],[4,8,9,10],[10,9,7,6],[8]}} \\ \textbf{EM:} {\small\texttt{(B0, A0C0C3D0, B1E0D1B2, B3C2F3F2, C1G0F1F0, E3E2D3D2, E1)}} \\ \textbf{Kauffman bracket:} {\scriptsize $A^{20} - A^{16} + 2A^{12} - A^{8} + A^{4} - 1$} \\ \textbf{Arrow:} {\scriptsize $A^{-4} - 1/A^{8} + 2/A^{12} - 1/A^{16} + A^{-20} - 1/A^{24}$} \\ \textbf{Mock:} {\scriptsize $2w^{2} - 3 + 2/w^{2}$} \\ \textbf{Affine:} {\scriptsize $0$} \\ \textbf{Yamada:} {\scriptsize $A^{19} + A^{17} + 2A^{16} + A^{15} + 2A^{14} + 2A^{13} + A^{11} - A^{10} - A^{8} - A^{7} - A^{5} - A^{4} + 1$}
\end{minipage}

\noindent{\color{gray!40}\rule{\textwidth}{0.4pt}}
\vspace{0.9\baselineskip}
\noindent \begin{minipage}[t]{0.25\textwidth}
\vspace{0pt}
\centering
\includegraphics[page=15,width=\linewidth]{knotoids.pdf}
\end{minipage}
\hfill
\begin{minipage}[t]{0.73\textwidth}
\vspace{0pt}
\raggedright
\textbf{Name:} {\large{$\mathbf{K5_{3}}$}} (chiral, non-rotatable$^{*}$) \\ \textbf{EM:} {\small\texttt{[0],[0,1,2,3],[1,4,5,2],[3,5,6,7],[8,9,10,4],[9,8,7,6],[10]}} \\ \textbf{EM:} {\small\texttt{(B0, A0C0C3D0, B1E3D1B2, B3C2F3F2, F1F0G0C1, E1E0D3D2, E2)}} \\ \textbf{Kauffman bracket:} {\scriptsize $A^{20} - A^{16} + 2A^{12} - 2A^{8} - A^{6} + A^{4} + A^{2}$} \\ \textbf{Arrow:} {\scriptsize $A^{8} - A^{4} + 2 - 2/A^{4} - L_1/A^{6} + A^{-8} + L_1/A^{10}$} \\ \textbf{Mock:} {\scriptsize $-2w^{2} - w + 4 + 1/w - 1/w^{2}$} \\ \textbf{Affine:} {\scriptsize $-t + 2 - 1/t$} \\ \textbf{Yamada:} {\scriptsize $-A^{20} + A^{19} - 2A^{17} + A^{16} - 2A^{14} - A^{12} - A^{10} + A^{8} - 2A^{7} + A^{5} - A^{4} + A^{2} - 1$}
\end{minipage}

\noindent{\color{gray!40}\rule{\textwidth}{0.4pt}}
\vspace{0.9\baselineskip}
\noindent \begin{minipage}[t]{0.25\textwidth}
\vspace{0pt}
\centering
\includegraphics[page=16,width=\linewidth]{knotoids.pdf}
\end{minipage}
\hfill
\begin{minipage}[t]{0.73\textwidth}
\vspace{0pt}
\raggedright
\textbf{Name:} {\large{$\mathbf{K5_{4}}$}} (chiral, non-rotatable$^{*}$) \\ \textbf{EM:} {\small\texttt{[0],[0,1,2,3],[1,4,5,2],[3,6,7,8],[4,9,10,5],[6,10,9,8],[7]}} \\ \textbf{EM:} {\small\texttt{(B0, A0C0C3D0, B1E0E3B2, B3F0G0F3, C1F2F1C2, D1E2E1D3, D2)}} \\ \textbf{Kauffman bracket:} {\scriptsize $-A^{19} + A^{15} + A^{13} - A^{11} - A^{9} + A^{7} + A^{5} - A^{3} - A$} \\ \textbf{Arrow:} {\scriptsize $A^{16} - A^{12} - A^{10}L_1 + A^{8} + A^{6}L_1 - A^{4} - A^{2}L_1 + 1 + L_1/A^{2}$} \\ \textbf{Mock:} {\scriptsize $-2w^{2} - 2w + 3 + 2/w$} \\ \textbf{Affine:} {\scriptsize $-2t + 4 - 2/t$} \\ \textbf{Yamada:} {\scriptsize $A^{21} - A^{19} - A^{16} + A^{14} - A^{13} - A^{10} - A^{8} + A^{7} - A^{6} - A^{5} + A^{4} - A^{3} - A^{2} - 1$}
\end{minipage}

\noindent{\color{gray!40}\rule{\textwidth}{0.4pt}}
\vspace{0.9\baselineskip}
\noindent \begin{minipage}[t]{0.25\textwidth}
\vspace{0pt}
\centering
\includegraphics[page=17,width=\linewidth]{knotoids.pdf}
\end{minipage}
\hfill
\begin{minipage}[t]{0.73\textwidth}
\vspace{0pt}
\raggedright
\textbf{Name:} {\large{$\mathbf{K5_{5}}$}} (chiral, non-rotatable$^{*}$) \\ \textbf{EM:} {\small\texttt{[0],[0,1,2,3],[1,4,5,2],[3,5,6,7],[4,8,9,6],[10,9,8,7],[10]}} \\ \textbf{EM:} {\small\texttt{(B0, A0C0C3D0, B1E0D1B2, B3C2E3F3, C1F2F1D2, G0E2E1D3, F0)}} \\ \textbf{Kauffman bracket:} {\scriptsize $-A^{19} - A^{17} + A^{15} + 2A^{13} - A^{11} - 2A^{9} + 2A^{5} - A$} \\ \textbf{Arrow:} {\scriptsize $A^{4} + A^{2}L_1 - 1 - 2L_1/A^{2} + A^{-4} + 2L_1/A^{6} - 2L_1/A^{10} + L_1/A^{14}$} \\ \textbf{Mock:} {\scriptsize $-w^{4} - w^{3} + 2w^{2} + 3w - 3/w + w^{-3}$} \\ \textbf{Affine:} {\scriptsize $0$} \\ \textbf{Yamada:} {\scriptsize $-A^{20} + A^{18} - A^{17} - A^{16} + 2A^{15} - A^{14} - 2A^{13} + A^{12} - A^{11} - 2A^{10} - A^{8} - 2A^{6} + A^{5} + 2A^{4} - 2A^{3} + 2A^{2} + A - 2$}
\end{minipage}

\noindent{\color{gray!40}\rule{\textwidth}{0.4pt}}
\vspace{0.9\baselineskip}
\noindent \begin{minipage}[t]{0.25\textwidth}
\vspace{0pt}
\centering
\includegraphics[page=18,width=\linewidth]{knotoids.pdf}
\end{minipage}
\hfill
\begin{minipage}[t]{0.73\textwidth}
\vspace{0pt}
\raggedright
\textbf{Name:} {\large{$\mathbf{K5_{6}}$}} (chiral, non-rotatable$^{*}$) \\ \textbf{EM:} {\small\texttt{[0],[0,1,2,3],[1,4,5,2],[3,6,7,8],[4,9,6,5],[9,8,10,7],[10]}} \\ \textbf{EM:} {\small\texttt{(B0, A0C0C3D0, B1E0E3B2, B3E2F3F1, C1F0D1C2, E1D3G0D2, F2)}} \\ \textbf{Kauffman bracket:} {\scriptsize $-A^{19} - A^{17} + A^{15} + A^{13} - A^{11} - 2A^{9} + A^{7} + 2A^{5} - A$} \\ \textbf{Arrow:} {\scriptsize $L_1/A^{2} + A^{-4} - L_1/A^{6} - 1/A^{8} + L_1/A^{10} + 2/A^{12} - L_1/A^{14} - 2/A^{16} + A^{-20}$} \\ \textbf{Mock:} {\scriptsize $w^{3} + 2w^{2} - w - 2 + 1/w + 2/w^{2} - 1/w^{3} - 1/w^{4}$} \\ \textbf{Affine:} {\scriptsize $2t - 4 + 2/t$} \\ \textbf{Yamada:} {\scriptsize $A^{21} + 2A^{18} + A^{17} - A^{16} + 2A^{15} + A^{14} - A^{13} + A^{12} + A^{11} - A^{10} - A^{8} + A^{7} - A^{6} + 3A^{4} - 2A^{3} + A - 1$}
\end{minipage}

\noindent{\color{gray!40}\rule{\textwidth}{0.4pt}}
\vspace{0.9\baselineskip}
\noindent \begin{minipage}[t]{0.25\textwidth}
\vspace{0pt}
\centering
\includegraphics[page=19,width=\linewidth]{knotoids.pdf}
\end{minipage}
\hfill
\begin{minipage}[t]{0.73\textwidth}
\vspace{0pt}
\raggedright
\textbf{Name:} {\large{$\mathbf{K5_{7}}$}} (chiral, non-rotatable$^{*}$) \\ \textbf{EM:} {\small\texttt{[0],[0,1,2,3],[1,4,5,2],[3,5,6,7],[4,8,9,6],[7,9,10,8],[10]}} \\ \textbf{EM:} {\small\texttt{(B0, A0C0C3D0, B1E0D1B2, B3C2E3F0, C1F3F1D2, D3E2G0E1, F2)}} \\ \textbf{Kauffman bracket:} {\scriptsize $-A^{18} + 2A^{14} + A^{12} - 2A^{10} - A^{8} + 2A^{6} + 2A^{4} - A^{2} - 1$} \\ \textbf{Arrow:} {\scriptsize $-A^{6}L_1 + 2A^{2}L_1 + 1 - 2L_1/A^{2} - 1/A^{4} + 2L_1/A^{6} + 2/A^{8} - L_1/A^{10} - 1/A^{12}$} \\ \textbf{Mock:} {\scriptsize $-w^{3} + 3w + 2 - 3/w - 2/w^{2} + w^{-3} + w^{-4}$} \\ \textbf{Affine:} {\scriptsize $0$} \\ \textbf{Yamada:} {\scriptsize $-A^{21} + A^{20} + A^{19} - 3A^{18} + A^{17} + 2A^{16} - 2A^{15} + 2A^{14} + 3A^{13} + 2A^{11} + 2A^{9} - A^{8} - A^{7} + 3A^{6} - 2A^{5} - 2A^{4} + 2A^{3} - A^{2} - A + 1$}
\end{minipage}

\noindent{\color{gray!40}\rule{\textwidth}{0.4pt}}
\vspace{0.9\baselineskip}
\noindent \begin{minipage}[t]{0.25\textwidth}
\vspace{0pt}
\centering
\includegraphics[page=20,width=\linewidth]{knotoids.pdf}
\end{minipage}
\hfill
\begin{minipage}[t]{0.73\textwidth}
\vspace{0pt}
\raggedright
\textbf{Name:} {\large{$\mathbf{K5_{8}}$}} (chiral, non-rotatable$^{*}$) \\ \textbf{EM:} {\small\texttt{[0],[0,1,2,3],[1,4,5,2],[3,6,7,4],[8,9,10,5],[6,10,9,7],[8]}} \\ \textbf{EM:} {\small\texttt{(B0, A0C0C3D0, B1D3E3B2, B3F0F3C1, G0F2F1C2, D1E2E1D2, E0)}} \\ \textbf{Kauffman bracket:} {\scriptsize $-A^{19} + A^{15} + A^{13} - 2A^{11} - 2A^{9} + A^{7} + 2A^{5} - A$} \\ \textbf{Arrow:} {\scriptsize $A^{4} - 1 - L_1/A^{2} + 2/A^{4} + 2L_1/A^{6} - 1/A^{8} - 2L_1/A^{10} + L_1/A^{14}$} \\ \textbf{Mock:} {\scriptsize $-w^{4} - w^{3} + 2w^{2} + 2w - 1 - 2/w + w^{-2} + w^{-3}$} \\ \textbf{Affine:} {\scriptsize $-t + 2 - 1/t$} \\ \textbf{Yamada:} {\scriptsize $-A^{20} + 2A^{18} - A^{16} + 3A^{15} + A^{14} - A^{13} + 2A^{12} + A^{11} - A^{10} + A^{9} + A^{7} - 2A^{6} + 2A^{4} - 3A^{3} + A^{2} + 2A - 1$}
\end{minipage}

\noindent{\color{gray!40}\rule{\textwidth}{0.4pt}}
\vspace{0.9\baselineskip}
\noindent \begin{minipage}[t]{0.25\textwidth}
\vspace{0pt}
\centering
\includegraphics[page=21,width=\linewidth]{knotoids.pdf}
\end{minipage}
\hfill
\begin{minipage}[t]{0.73\textwidth}
\vspace{0pt}
\raggedright
\textbf{Name:} {\large{$\mathbf{K5_{9}}$}} (chiral, non-rotatable$^{*}$) \\ \textbf{EM:} {\small\texttt{[0],[0,1,2,3],[1,4,5,2],[3,6,7,4],[8,9,6,5],[9,8,10,7],[10]}} \\ \textbf{EM:} {\small\texttt{(B0, A0C0C3D0, B1D3E3B2, B3E2F3C1, F1F0D1C2, E1E0G0D2, F2)}} \\ \textbf{Kauffman bracket:} {\scriptsize $-A^{18} + 2A^{14} + A^{12} - 2A^{10} - 2A^{8} + 2A^{6} + 2A^{4} - 1$} \\ \textbf{Arrow:} {\scriptsize $-A^{12} + 2A^{8} + A^{6}L_1 - 2A^{4} - 2A^{2}L_1 + 2 + 2L_1/A^{2} - L_1/A^{6}$} \\ \textbf{Mock:} {\scriptsize $w^{4} + w^{3} - 2w^{2} - 2w + 3 + 2/w - 1/w^{2} - 1/w^{3}$} \\ \textbf{Affine:} {\scriptsize $t - 2 + 1/t$} \\ \textbf{Yamada:} {\scriptsize $A^{20} - A^{19} - 2A^{18} + 3A^{17} - 3A^{15} + 3A^{14} - 2A^{12} + A^{11} - A^{10} - 3A^{8} + A^{6} - 4A^{5} + 2A^{3} - 3A^{2} + 2$}
\end{minipage}

\noindent{\color{gray!40}\rule{\textwidth}{0.4pt}}
\vspace{0.9\baselineskip}
\noindent \begin{minipage}[t]{0.25\textwidth}
\vspace{0pt}
\centering
\includegraphics[page=22,width=\linewidth]{knotoids.pdf}
\end{minipage}
\hfill
\begin{minipage}[t]{0.73\textwidth}
\vspace{0pt}
\raggedright
\textbf{Name:} {\large{$\mathbf{K5_{10}}$}} (chiral, non-rotatable$^{*}$) \\ \textbf{EM:} {\small\texttt{[0],[0,1,2,3],[1,4,5,2],[3,5,6,7],[7,8,9,4],[6,9,8,10],[10]}} \\ \textbf{EM:} {\small\texttt{(B0, A0C0C3D0, B1E3D1B2, B3C2F0E0, D3F2F1C1, D2E2E1G0, F3)}} \\ \textbf{Kauffman bracket:} {\scriptsize $A^{16} + A^{14} - 2A^{10} + 2A^{6} + A^{4} - A^{2} - 1$} \\ \textbf{Arrow:} {\scriptsize $A^{-8} + L_1/A^{10} - 2L_1/A^{14} + 2L_1/A^{18} + A^{-20} - L_1/A^{22} - 1/A^{24}$} \\ \textbf{Mock:} {\scriptsize $w^{4} + w^{3} - 2w - 1 + 2/w + w^{-2} - 1/w^{3}$} \\ \textbf{Affine:} {\scriptsize $t - 2 + 1/t$} \\ \textbf{Yamada:} {\scriptsize $A^{21} + A^{20} + 3A^{17} + A^{16} + 3A^{14} + A^{13} - A^{12} + A^{11} - A^{10} - 2A^{8} - 3A^{5} + A^{4} + A^{3} - A^{2} + 1$}
\end{minipage}

\noindent{\color{gray!40}\rule{\textwidth}{0.4pt}}
\vspace{0.9\baselineskip}
\noindent \begin{minipage}[t]{0.25\textwidth}
\vspace{0pt}
\centering
\includegraphics[page=23,width=\linewidth]{knotoids.pdf}
\end{minipage}
\hfill
\begin{minipage}[t]{0.73\textwidth}
\vspace{0pt}
\raggedright
\textbf{Name:} {\large{$\mathbf{K5_{11}}$}} (chiral, non-rotatable$^{*}$) \\ \textbf{EM:} {\small\texttt{[0],[0,1,2,3],[1,4,5,2],[3,5,6,7],[4,7,8,9],[9,10,8,6],[10]}} \\ \textbf{EM:} {\small\texttt{(B0, A0C0C3D0, B1E0D1B2, B3C2F3E1, C1D3F2F0, E3G0E2D2, F1)}} \\ \textbf{Kauffman bracket:} {\scriptsize $A^{14} - 2A^{10} - A^{8} + 2A^{6} + 2A^{4} - 1$} \\ \textbf{Arrow:} {\scriptsize $A^{8} - A^{4}L_2 - A^{4} - A^{2}L_1 + L_2 + 1 + 2L_1/A^{2} - L_1/A^{6}$} \\ \textbf{Mock:} {\scriptsize $-w^{2} - w + 3 + 2/w - 1/w^{2} - 1/w^{3}$} \\ \textbf{Affine:} {\scriptsize $-t^{2} + 2 - 1/t^{2}$} \\ \textbf{Yamada:} {\scriptsize $A^{19} + A^{18} - A^{17} + A^{15} - 2A^{14} - 2A^{13} + A^{12} - A^{11} - A^{10} + A^{9} - A^{8} - 2A^{6} + A^{5} - 2A^{3} + 2A^{2} - 2$}
\end{minipage}

\noindent{\color{gray!40}\rule{\textwidth}{0.4pt}}
\vspace{0.9\baselineskip}
\noindent \begin{minipage}[t]{0.25\textwidth}
\vspace{0pt}
\centering
\includegraphics[page=24,width=\linewidth]{knotoids.pdf}
\end{minipage}
\hfill
\begin{minipage}[t]{0.73\textwidth}
\vspace{0pt}
\raggedright
\textbf{Name:} {\large{$\mathbf{K5_{12}}$}} (chiral, non-rotatable$^{*}$) \\ \textbf{EM:} {\small\texttt{[0],[0,1,2,3],[1,4,5,2],[3,5,6,7],[7,8,9,4],[9,10,8,6],[10]}} \\ \textbf{EM:} {\small\texttt{(B0, A0C0C3D0, B1E3D1B2, B3C2F3E0, D3F2F0C1, E2G0E1D2, F1)}} \\ \textbf{Kauffman bracket:} {\scriptsize $A^{15} - 2A^{11} - A^{9} + 2A^{7} + 2A^{5} - A^{3} - 2A$} \\ \textbf{Arrow:} {\scriptsize $-A^{18}L_1 + 2A^{14}L_1 + A^{12} - 2A^{10}L_1 - A^{8}L_2 - A^{8} + A^{6}L_1 + A^{4}L_2 + A^{4}$} \\ \textbf{Mock:} {\scriptsize $w^{3} + w^{2} - 3w - 2 + 2/w + 2/w^{2}$} \\ \textbf{Affine:} {\scriptsize $-t^{2} - t + 4 - 1/t - 1/t^{2}$} \\ \textbf{Yamada:} {\scriptsize $A^{22} - A^{20} + 2A^{19} + A^{18} - 3A^{17} + A^{16} - 3A^{14} - A^{12} + A^{11} - A^{10} + A^{9} + 2A^{8} - 3A^{7} - A^{6} + A^{5} - 2A^{4} - A^{3} + A^{2} - 1$}
\end{minipage}

\noindent{\color{gray!40}\rule{\textwidth}{0.4pt}}
\vspace{0.9\baselineskip}
\noindent \begin{minipage}[t]{0.25\textwidth}
\vspace{0pt}
\centering
\includegraphics[page=25,width=\linewidth]{knotoids.pdf}
\end{minipage}
\hfill
\begin{minipage}[t]{0.73\textwidth}
\vspace{0pt}
\raggedright
\textbf{Name:} {\large{$\mathbf{K5_{13}}$}} (chiral, non-rotatable$^{*}$) \\ \textbf{EM:} {\small\texttt{[0],[0,1,2,3],[1,4,5,2],[3,6,7,4],[5,7,8,9],[9,8,10,6],[10]}} \\ \textbf{EM:} {\small\texttt{(B0, A0C0C3D0, B1D3E0B2, B3F3E1C1, C2D2F1F0, E3E2G0D1, F2)}} \\ \textbf{Kauffman bracket:} {\scriptsize $A^{14} + 2A^{12} - A^{10} - 2A^{8} - A^{6} + A^{4} + A^{2}$} \\ \textbf{Arrow:} {\scriptsize $A^{2}L_1 + 2 - L_1/A^{2} - 2/A^{4} - L_1/A^{6} + A^{-8} + L_1/A^{10}$} \\ \textbf{Mock:} {\scriptsize $-w^{2} + 3 - 1/w^{2}$} \\ \textbf{Affine:} {\scriptsize $0$} \\ \textbf{Yamada:} {\scriptsize $-A^{20} + A^{19} + A^{18} - A^{17} + 2A^{16} + 2A^{15} + 2A^{13} - 2A^{10} - A^{9} + A^{8} - A^{7} + A^{6} + 2A^{5} + A^{2} - 1$}
\end{minipage}

\noindent{\color{gray!40}\rule{\textwidth}{0.4pt}}
\vspace{0.9\baselineskip}
\noindent \begin{minipage}[t]{0.25\textwidth}
\vspace{0pt}
\centering
\includegraphics[page=26,width=\linewidth]{knotoids.pdf}
\end{minipage}
\hfill
\begin{minipage}[t]{0.73\textwidth}
\vspace{0pt}
\raggedright
\textbf{Name:} {\large{$\mathbf{K5_{14}}$}} (chiral, rotatable) \\ \textbf{EM:} {\small\texttt{[0],[0,1,2,3],[1,4,5,2],[3,6,7,4],[7,8,9,5],[6,9,8,10],[10]}} \\ \textbf{EM:} {\small\texttt{(B0, A0C0C3D0, B1D3E3B2, B3F0E0C1, D2F2F1C2, D1E2E1G0, F3)}} \\ \textbf{Kauffman bracket:} {\scriptsize $-A^{13} + A^{9} + 2A^{7} - 2A^{3} - A$} \\ \textbf{Arrow:} {\scriptsize $A^{16} - A^{12} - 2A^{10}L_1 + 2A^{6}L_1 + A^{4}$} \\ \textbf{Mock:} {\scriptsize $-w^{2} - 2w + 1 + 2/w + w^{-2}$} \\ \textbf{Affine:} {\scriptsize $-2t + 4 - 2/t$} \\ \textbf{Yamada:} {\scriptsize $A^{17} - A^{16} - 2A^{15} + A^{14} - A^{13} + 2A^{11} + A^{10} + 2A^{9} + A^{7} - A^{5} + A^{4} + A^{3} - A^{2} + A + 1$}
\end{minipage}

\noindent{\color{gray!40}\rule{\textwidth}{0.4pt}}
\vspace{0.9\baselineskip}
\noindent \begin{minipage}[t]{0.25\textwidth}
\vspace{0pt}
\centering
\includegraphics[page=27,width=\linewidth]{knotoids.pdf}
\end{minipage}
\hfill
\begin{minipage}[t]{0.73\textwidth}
\vspace{0pt}
\raggedright
\textbf{Name:} {\large{$\mathbf{K5_{15}}$}} (chiral, rotatable) \\ \textbf{EM:} {\small\texttt{[0],[0,1,2,3],[1,4,5,2],[6,7,4,3],[7,8,9,5],[6,9,8,10],[10]}} \\ \textbf{EM:} {\small\texttt{(B0, A0C0C3D3, B1D2E3B2, F0E0C1B3, D1F2F1C2, D0E2E1G0, F3)}} \\ \textbf{Kauffman bracket:} {\scriptsize $2A^{15} + A^{13} - 2A^{11} - 3A^{9} + 2A^{5} - A$} \\ \textbf{Arrow:} {\scriptsize $-2A^{18}L_1 - A^{16} + 2A^{14}L_1 + 3A^{12} - 2A^{8} + A^{4}$} \\ \textbf{Mock:} {\scriptsize $w^{2} - 2w - 3 + 2/w + 3/w^{2}$} \\ \textbf{Affine:} {\scriptsize $-2t + 4 - 2/t$} \\ \textbf{Yamada:} {\scriptsize $-A^{21} - A^{20} + A^{19} + A^{18} - A^{17} + 2A^{16} + 3A^{15} - A^{14} - A^{13} + 2A^{12} - 2A^{11} - 2A^{10} + A^{9} - A^{8} - 2A^{6} + 2A^{5} - A^{4} - 3A^{3} + 2A^{2} - A - 3$}
\end{minipage}

\noindent{\color{gray!40}\rule{\textwidth}{0.4pt}}
\vspace{0.9\baselineskip}
\noindent \begin{minipage}[t]{0.25\textwidth}
\vspace{0pt}
\centering
\includegraphics[page=28,width=\linewidth]{knotoids.pdf}
\end{minipage}
\hfill
\begin{minipage}[t]{0.73\textwidth}
\vspace{0pt}
\raggedright
\textbf{Name:} {\large{$\mathbf{K5_{16}}$}} (chiral, non-rotatable$^{*}$) \\ \textbf{EM:} {\small\texttt{[0],[0,1,2,3],[1,4,5,2],[3,6,7,4],[5,8,9,10],[6,10,8,7],[9]}} \\ \textbf{EM:} {\small\texttt{(B0, A0C0C3D0, B1D3E0B2, B3F0F3C1, C2F2G0F1, D1E3E1D2, E2)}} \\ \textbf{Kauffman bracket:} {\scriptsize $A^{14} - A^{12} - 2A^{10} + 3A^{6} + 2A^{4} - A^{2} - 1$} \\ \textbf{Arrow:} {\scriptsize $A^{8} - A^{6}L_1 - 2A^{4} + 3 + 2L_1/A^{2} - 1/A^{4} - L_1/A^{6}$} \\ \textbf{Mock:} {\scriptsize $-w^{3} - 2w^{2} + 4 + 2/w - 1/w^{2} - 1/w^{3}$} \\ \textbf{Affine:} {\scriptsize $-t + 2 - 1/t$} \\ \textbf{Yamada:} {\scriptsize $-A^{21} + A^{20} + 2A^{19} - 2A^{18} + A^{17} + 2A^{16} - 3A^{15} - 3A^{12} - A^{11} - 2A^{10} + A^{9} - A^{8} - A^{7} + 3A^{6} - A^{5} - 2A^{4} + 2A^{3} - A^{2} - A + 1$}
\end{minipage}

\noindent{\color{gray!40}\rule{\textwidth}{0.4pt}}
\vspace{0.9\baselineskip}
\noindent \begin{minipage}[t]{0.25\textwidth}
\vspace{0pt}
\centering
\includegraphics[page=29,width=\linewidth]{knotoids.pdf}
\end{minipage}
\hfill
\begin{minipage}[t]{0.73\textwidth}
\vspace{0pt}
\raggedright
\textbf{Name:} {\large{$\mathbf{K5_{17}}$}} (chiral, non-rotatable$^{*}$) \\ \textbf{EM:} {\small\texttt{[0],[0,1,2,3],[1,4,5,2],[3,6,7,4],[8,9,10,5],[6,10,8,7],[9]}} \\ \textbf{EM:} {\small\texttt{(B0, A0C0C3D0, B1D3E3B2, B3F0F3C1, F2G0F1C2, D1E2E0D2, E1)}} \\ \textbf{Kauffman bracket:} {\scriptsize $A^{16} - A^{12} - 2A^{10} + A^{8} + 3A^{6} + A^{4} - A^{2} - 1$} \\ \textbf{Arrow:} {\scriptsize $A^{4} - 1 - 2L_1/A^{2} + A^{-4} + 3L_1/A^{6} + A^{-8} - L_1/A^{10} - 1/A^{12}$} \\ \textbf{Mock:} {\scriptsize $-w^{4} - 2w^{3} + w^{2} + 3w + 1 - 1/w$} \\ \textbf{Affine:} {\scriptsize $-t + 2 - 1/t$} \\ \textbf{Yamada:} {\scriptsize $A^{20} - 2A^{18} + 2A^{17} - 4A^{15} - A^{13} - 3A^{12} - A^{10} + A^{9} - A^{8} + A^{7} + 3A^{6} - 2A^{5} + 2A^{3} - 2A^{2} - A + 1$}
\end{minipage}

\noindent{\color{gray!40}\rule{\textwidth}{0.4pt}}
\vspace{0.9\baselineskip}
\noindent \begin{minipage}[t]{0.25\textwidth}
\vspace{0pt}
\centering
\includegraphics[page=30,width=\linewidth]{knotoids.pdf}
\end{minipage}
\hfill
\begin{minipage}[t]{0.73\textwidth}
\vspace{0pt}
\raggedright
\textbf{Name:} {\large{$\mathbf{K5_{18}}$}} (chiral, non-rotatable$^{*}$) \\ \textbf{EM:} {\small\texttt{[0],[0,1,2,3],[1,4,5,2],[6,7,4,3],[5,8,9,10],[10,8,7,6],[9]}} \\ \textbf{EM:} {\small\texttt{(B0, A0C0C3D3, B1D2E0B2, F3F2C1B3, C2F1G0F0, E3E1D1D0, E2)}} \\ \textbf{Kauffman bracket:} {\scriptsize $-A^{13} - 2A^{11} + 2A^{7} + 2A^{5} - A^{3} - A$} \\ \textbf{Arrow:} {\scriptsize $L_1/A^{2} + 2/A^{4} - 2/A^{8} - 2L_1/A^{10} + A^{-12} + L_1/A^{14}$} \\ \textbf{Mock:} {\scriptsize $w^{3} + 2w^{2} - 2 - 2/w + w^{-2} + w^{-3}$} \\ \textbf{Affine:} {\scriptsize $t - 2 + 1/t$} \\ \textbf{Yamada:} {\scriptsize $A^{21} - A^{19} + A^{18} + A^{17} - A^{16} + 2A^{15} + 2A^{14} + A^{12} + A^{11} - A^{10} - A^{8} + A^{7} - A^{6} - A^{5} + 3A^{4} - A^{3} + A - 1$}
\end{minipage}

\noindent{\color{gray!40}\rule{\textwidth}{0.4pt}}
\vspace{0.9\baselineskip}
\noindent \begin{minipage}[t]{0.25\textwidth}
\vspace{0pt}
\centering
\includegraphics[page=31,width=\linewidth]{knotoids.pdf}
\end{minipage}
\hfill
\begin{minipage}[t]{0.73\textwidth}
\vspace{0pt}
\raggedright
\textbf{Name:} {\large{$\mathbf{K5_{19}}$}} (chiral, non-rotatable$^{*}$) \\ \textbf{EM:} {\small\texttt{[0],[0,1,2,3],[1,4,5,2],[6,7,4,3],[8,9,10,5],[10,8,7,6],[9]}} \\ \textbf{EM:} {\small\texttt{(B0, A0C0C3D3, B1D2E3B2, F3F2C1B3, F1G0F0C2, E2E0D1D0, E1)}} \\ \textbf{Kauffman bracket:} {\scriptsize $A^{16} + 2A^{14} + A^{12} - 3A^{10} - 2A^{8} + A^{6} + 2A^{4} - 1$} \\ \textbf{Arrow:} {\scriptsize $A^{-8} + 2L_1/A^{10} + A^{-12} - 3L_1/A^{14} - 2/A^{16} + L_1/A^{18} + 2/A^{20} - 1/A^{24}$} \\ \textbf{Mock:} {\scriptsize $w^{4} + 2w^{3} + w^{2} - 3w - 3 + 1/w + 2/w^{2}$} \\ \textbf{Affine:} {\scriptsize $t - 2 + 1/t$} \\ \textbf{Yamada:} {\scriptsize $-A^{22} + A^{20} - 2A^{19} - 2A^{18} + 2A^{17} - 3A^{16} - 4A^{15} + A^{14} - A^{13} - 2A^{12} + 2A^{11} + A^{10} + 2A^{9} - A^{8} + 2A^{7} + 2A^{6} - 4A^{5} + 2A^{4} + A^{3} - 3A^{2} + 1$}
\end{minipage}

\noindent{\color{gray!40}\rule{\textwidth}{0.4pt}}
\vspace{0.9\baselineskip}
\noindent \begin{minipage}[t]{0.25\textwidth}
\vspace{0pt}
\centering
\includegraphics[page=32,width=\linewidth]{knotoids.pdf}
\end{minipage}
\hfill
\begin{minipage}[t]{0.73\textwidth}
\vspace{0pt}
\raggedright
\textbf{Name:} {\large{$\mathbf{K5_{20}}$}} (chiral, non-rotatable$^{*}$) \\ \textbf{EM:} {\small\texttt{[0],[0,1,2,3],[1,4,5,2],[3,6,7,4],[5,8,9,6],[7,9,10,8],[10]}} \\ \textbf{EM:} {\small\texttt{(B0, A0C0C3D0, B1D3E0B2, B3E3F0C1, C2F3F1D1, D2E2G0E1, F2)}} \\ \textbf{Kauffman bracket:} {\scriptsize $A^{14} + A^{12} - A^{10} - 2A^{8} + A^{4} + A^{2}$} \\ \textbf{Arrow:} {\scriptsize $A^{8} + A^{6}L_1 - A^{4} - 2A^{2}L_1 - L_2 + 1 + L_1/A^{2} + L_2/A^{4}$} \\ \textbf{Mock:} {\scriptsize $w^{3} - 2w + 1 + 1/w$} \\ \textbf{Affine:} {\scriptsize $-t^{2} + 2 - 1/t^{2}$} \\ \textbf{Yamada:} {\scriptsize $A^{21} + A^{18} - 2A^{17} - A^{16} - 2A^{14} - A^{12} - A^{10} - A^{9} + A^{8} + A^{5} - A^{3} - 1$}
\end{minipage}

\noindent{\color{gray!40}\rule{\textwidth}{0.4pt}}
\vspace{0.9\baselineskip}
\noindent \begin{minipage}[t]{0.25\textwidth}
\vspace{0pt}
\centering
\includegraphics[page=33,width=\linewidth]{knotoids.pdf}
\end{minipage}
\hfill
\begin{minipage}[t]{0.73\textwidth}
\vspace{0pt}
\raggedright
\textbf{Name:} {\large{$\mathbf{K5_{21}}$}} (chiral, non-rotatable$^{*}$) \\ \textbf{EM:} {\small\texttt{[0],[0,1,2,3],[1,4,5,2],[3,6,7,4],[8,9,6,5],[7,9,10,8],[10]}} \\ \textbf{EM:} {\small\texttt{(B0, A0C0C3D0, B1D3E3B2, B3E2F0C1, F3F1D1C2, D2E1G0E0, F2)}} \\ \textbf{Kauffman bracket:} {\scriptsize $A^{16} + A^{14} - A^{12} - 2A^{10} + 2A^{6} + A^{4} - A^{2}$} \\ \textbf{Arrow:} {\scriptsize $A^{4} + A^{2}L_1 - 1 - 2L_1/A^{2} - L_2/A^{4} + A^{-4} + 2L_1/A^{6} + L_2/A^{8} - L_1/A^{10}$} \\ \textbf{Mock:} {\scriptsize $-w^{4} - w^{3} + 2w^{2} + 3w + 1 - 2/w - 1/w^{2}$} \\ \textbf{Affine:} {\scriptsize $-t^{2} + t + 1/t - 1/t^{2}$} \\ \textbf{Yamada:} {\scriptsize $A^{21} - A^{19} + A^{18} - 2A^{16} - A^{14} - 3A^{13} - A^{12} - A^{10} + A^{9} + 2A^{7} - A^{6} - A^{5} + 2A^{4} - 2A^{3} + A - 1$}
\end{minipage}

\noindent{\color{gray!40}\rule{\textwidth}{0.4pt}}
\vspace{0.9\baselineskip}
\noindent \begin{minipage}[t]{0.25\textwidth}
\vspace{0pt}
\centering
\includegraphics[page=34,width=\linewidth]{knotoids.pdf}
\end{minipage}
\hfill
\begin{minipage}[t]{0.73\textwidth}
\vspace{0pt}
\raggedright
\textbf{Name:} {\large{$\mathbf{K5_{22}}$}} (chiral, non-rotatable$^{*}$) \\ \textbf{EM:} {\small\texttt{[0],[0,1,2,3],[1,4,5,2],[3,6,7,4],[8,9,6,5],[9,10,8,7],[10]}} \\ \textbf{EM:} {\small\texttt{(B0, A0C0C3D0, B1D3E3B2, B3E2F3C1, F2F0D1C2, E1G0E0D2, F1)}} \\ \textbf{Kauffman bracket:} {\scriptsize $A^{16} + A^{14} - A^{12} - 3A^{10} - A^{8} + 3A^{6} + 2A^{4} - 1$} \\ \textbf{Arrow:} {\scriptsize $A^{10}L_1 + A^{8} - A^{6}L_1 - A^{4}L_2 - 2A^{4} - A^{2}L_1 + L_2 + 2 + 2L_1/A^{2} - L_1/A^{6}$} \\ \textbf{Mock:} {\scriptsize $-w^{2} + 4 + 1/w - 2/w^{2} - 1/w^{3}$} \\ \textbf{Affine:} {\scriptsize $-t^{2} + t + 1/t - 1/t^{2}$} \\ \textbf{Yamada:} {\scriptsize $-A^{22} + 2A^{20} - A^{19} - 2A^{18} + 3A^{17} - 3A^{15} + 3A^{14} + A^{13} - A^{12} + 2A^{11} + A^{10} + 2A^{9} - 2A^{8} + A^{7} + 2A^{6} - 4A^{5} + A^{4} + 3A^{3} - 2A^{2} + 1$}
\end{minipage}

\noindent{\color{gray!40}\rule{\textwidth}{0.4pt}}
\vspace{0.9\baselineskip}
\noindent \begin{minipage}[t]{0.25\textwidth}
\vspace{0pt}
\centering
\includegraphics[page=35,width=\linewidth]{knotoids.pdf}
\end{minipage}
\hfill
\begin{minipage}[t]{0.73\textwidth}
\vspace{0pt}
\raggedright
\textbf{Name:} {\large{$\mathbf{K5_{23}}$}} (chiral, non-rotatable$^{*}$) \\ \textbf{EM:} {\small\texttt{[0],[0,1,2,3],[1,4,5,2],[6,7,4,3],[5,8,9,6],[7,9,10,8],[10]}} \\ \textbf{EM:} {\small\texttt{(B0, A0C0C3D3, B1D2E0B2, E3F0C1B3, C2F3F1D0, D1E2G0E1, F2)}} \\ \textbf{Kauffman bracket:} {\scriptsize $-A^{13} - A^{11} + A^{9} + 3A^{7} + A^{5} - 2A^{3} - 2A$} \\ \textbf{Arrow:} {\scriptsize $A^{16} + A^{14}L_1 - A^{12} - 3A^{10}L_1 - A^{8}L_2 + 2A^{6}L_1 + A^{4}L_2 + A^{4}$} \\ \textbf{Mock:} {\scriptsize $w^{3} - 3w + 2/w + w^{-2}$} \\ \textbf{Affine:} {\scriptsize $-t^{2} - t + 4 - 1/t - 1/t^{2}$} \\ \textbf{Yamada:} {\scriptsize $A^{23} - A^{21} + A^{20} + A^{19} - 2A^{18} - A^{17} + 2A^{16} - 2A^{15} - A^{14} + A^{13} - A^{12} - A^{10} + 2A^{9} - A^{8} - 2A^{7} + 2A^{6} - A^{5} - 2A^{4} + A^{3} - A - 1$}
\end{minipage}

\noindent{\color{gray!40}\rule{\textwidth}{0.4pt}}
\vspace{0.9\baselineskip}
\noindent \begin{minipage}[t]{0.25\textwidth}
\vspace{0pt}
\centering
\includegraphics[page=36,width=\linewidth]{knotoids.pdf}
\end{minipage}
\hfill
\begin{minipage}[t]{0.73\textwidth}
\vspace{0pt}
\raggedright
\textbf{Name:} {\large{$\mathbf{K5_{24}}$}} (chiral, rotatable) \\ \textbf{EM:} {\small\texttt{[0],[0,1,2,3],[1,3,4,5],[6,7,8,2],[9,6,5,4],[7,9,8,10],[10]}} \\ \textbf{EM:} {\small\texttt{(B0, A0C0D3C1, B1B3E3E2, E1F0F2B2, F1D0C3C2, D1E0D2G0, F3)}} \\ \textbf{Kauffman bracket:} {\scriptsize $-A^{18} - A^{16} + 2A^{14} + 3A^{12} - 3A^{8} - A^{6} + A^{4} + A^{2}$} \\ \textbf{Arrow:} {\scriptsize $-A^{12} - A^{10}L_1 + 2A^{8} + 3A^{6}L_1 - 3A^{2}L_1 - 1 + L_1/A^{2} + A^{-4}$} \\ \textbf{Mock:} {\scriptsize $-w^{2} - 2w + 2 + 4/w + w^{-2} - 2/w^{3} - 1/w^{4}$} \\ \textbf{Affine:} {\scriptsize $0$} \\ \textbf{Yamada:} {\scriptsize $-A^{23} + A^{22} + A^{21} - 3A^{20} + 2A^{19} + 2A^{18} - 4A^{17} + 2A^{16} + 2A^{15} - 2A^{14} + A^{13} + 2A^{11} - A^{10} + 4A^{8} - 2A^{7} - A^{6} + 4A^{5} - A^{4} - A^{3} + 2A^{2} - 1$}
\end{minipage}

\noindent{\color{gray!40}\rule{\textwidth}{0.4pt}}
\vspace{0.9\baselineskip}
\noindent \begin{minipage}[t]{0.25\textwidth}
\vspace{0pt}
\centering
\includegraphics[page=37,width=\linewidth]{knotoids.pdf}
\end{minipage}
\hfill
\begin{minipage}[t]{0.73\textwidth}
\vspace{0pt}
\raggedright
\textbf{Name:} {\large{$\mathbf{K5_{25}}$}} (chiral, rotatable) \\ \textbf{EM:} {\small\texttt{[0],[0,1,2,3],[1,3,4,5],[5,6,7,2],[8,9,6,4],[9,10,8,7],[10]}} \\ \textbf{EM:} {\small\texttt{(B0, A0C0D3C1, B1B3E3D0, C3E2F3B2, F2F0D1C2, E1G0E0D2, F1)}} \\ \textbf{Kauffman bracket:} {\scriptsize $A^{16} + 2A^{14} - 4A^{10} - 2A^{8} + 2A^{6} + 3A^{4} - 1$} \\ \textbf{Arrow:} {\scriptsize $A^{4} + 2A^{2}L_1 - 4L_1/A^{2} - 2/A^{4} + 2L_1/A^{6} + 3/A^{8} - 1/A^{12}$} \\ \textbf{Mock:} {\scriptsize $2w + 2 - 4/w - 2/w^{2} + 2/w^{3} + w^{-4}$} \\ \textbf{Affine:} {\scriptsize $0$} \\ \textbf{Yamada:} {\scriptsize $A^{23} - A^{22} - 2A^{21} + 3A^{20} + A^{19} - 4A^{18} + 4A^{17} + 3A^{16} - 4A^{15} + 2A^{14} + 2A^{13} - 2A^{12} + 3A^{9} - 2A^{8} + 6A^{6} - 4A^{5} - 2A^{4} + 4A^{3} - 2A^{2} - A + 1$}
\end{minipage}

\noindent{\color{gray!40}\rule{\textwidth}{0.4pt}}
\vspace{0.9\baselineskip}

\section*{Number of crossings: 6}

\noindent{\color{gray!40}\rule{\textwidth}{0.4pt}}
\vspace{0.9\baselineskip}
\noindent \begin{minipage}[t]{0.25\textwidth}
\vspace{0pt}
\centering
\includegraphics[page=38,width=\linewidth]{knotoids.pdf}
\end{minipage}
\hfill
\begin{minipage}[t]{0.73\textwidth}
\vspace{0pt}
\raggedright
\textbf{Name:} {\large{$\mathbf{K6_{1}}$}} (chiral, rotatable) \\ \textbf{PD:} {\scriptsize\texttt{[0],[0,1,2,3],[1,4,5,2],[3,5,6,7],[4,8,9,10],[11,12,7,6],[8],[12,11,10,9]}} \\ \textbf{EM:} {\small\texttt{(B0, A0C0C3D0, B1E0D1B2, B3C2F3F2, C1G0H3H2, H1H0D3D2, E1, F1F0E3E2)}} \\ \textbf{Kauffman bracket:} {\scriptsize $A^{26} - A^{22} + 2A^{18} - 2A^{14} + A^{10} - A^{6} + A^{2}$} \\ \textbf{Arrow:} {\scriptsize $A^{8} - A^{4} + 2 - 2/A^{4} + A^{-8} - 1/A^{12} + A^{-16}$} \\ \textbf{Mock:} {\scriptsize $-2w^{2} + 5 - 2/w^{2}$} \\ \textbf{Affine:} {\scriptsize $0$} \\ \textbf{Yamada:} {\scriptsize $-A^{22} - A^{20} - A^{19} + A^{18} + A^{15} - A^{14} - A^{12} - A^{11} - A^{10} - 2A^{9} + A^{5} + A^{4} - 1$}
\end{minipage}

\noindent{\color{gray!40}\rule{\textwidth}{0.4pt}}
\vspace{0.9\baselineskip}
\noindent \begin{minipage}[t]{0.25\textwidth}
\vspace{0pt}
\centering
\includegraphics[page=39,width=\linewidth]{knotoids.pdf}
\end{minipage}
\hfill
\begin{minipage}[t]{0.73\textwidth}
\vspace{0pt}
\raggedright
\textbf{Name:} {\large{$\mathbf{K6_{2}}$}} (chiral, rotatable) \\ \textbf{PD:} {\scriptsize\texttt{[0],[0,1,2,3],[1,4,5,2],[3,5,6,7],[4,8,9,6],[7,10,11,12],[8,11,10,9],[12]}} \\ \textbf{EM:} {\small\texttt{(B0, A0C0C3D0, B1E0D1B2, B3C2E3F0, C1G0G3D2, D3G2G1H0, E1F2F1E2, F3)}} \\ \textbf{Kauffman bracket:} {\scriptsize $-A^{25} + A^{21} - 2A^{17} + 2A^{13} - 2A^{9} + 2A^{5} - A$} \\ \textbf{Arrow:} {\scriptsize $A^{4} - 1 + 2/A^{4} - 2/A^{8} + 2/A^{12} - 2/A^{16} + A^{-20}$} \\ \textbf{Mock:} {\scriptsize $-w^{4} + 3w^{2} - 3 + 3/w^{2} - 1/w^{4}$} \\ \textbf{Affine:} {\scriptsize $0$} \\ \textbf{Yamada:} {\scriptsize $-A^{22} - A^{19} - 2A^{16} - A^{13} - A^{10} - A^{8} + A^{7} - A^{5} + 2A^{4} - A^{3} + A - 1$}
\end{minipage}

\noindent{\color{gray!40}\rule{\textwidth}{0.4pt}}
\vspace{0.9\baselineskip}
\noindent \begin{minipage}[t]{0.25\textwidth}
\vspace{0pt}
\centering
\includegraphics[page=40,width=\linewidth]{knotoids.pdf}
\end{minipage}
\hfill
\begin{minipage}[t]{0.73\textwidth}
\vspace{0pt}
\raggedright
\textbf{Name:} {\large{$\mathbf{K6_{3}}$}} (achiral, rotatable) \\ \textbf{PD:} {\scriptsize\texttt{[0],[0,1,2,3],[1,4,5,2],[3,5,6,7],[4,8,9,6],[7,9,10,11],[8,12,11,10],[12]}} \\ \textbf{EM:} {\small\texttt{(B0, A0C0C3D0, B1E0D1B2, B3C2E3F0, C1G0F1D2, D3E2G3G2, E1H0F3F2, G1)}} \\ \textbf{Kauffman bracket:} {\scriptsize $-A^{24} + 2A^{20} - 2A^{16} + 3A^{12} - 2A^{8} + 2A^{4} - 1$} \\ \textbf{Arrow:} {\scriptsize $-A^{12} + 2A^{8} - 2A^{4} + 3 - 2/A^{4} + 2/A^{8} - 1/A^{12}$} \\ \textbf{Mock:} {\scriptsize $w^{4} - 3w^{2} + 5 - 3/w^{2} + w^{-4}$} \\ \textbf{Affine:} {\scriptsize $0$} \\ \textbf{Yamada:} {\scriptsize $A^{22} - A^{21} - A^{20} + A^{19} - 2A^{18} + 2A^{16} - A^{15} + A^{14} + 2A^{13} + 2A^{11} + 2A^{9} + A^{8} - A^{7} + 2A^{6} - 2A^{4} + A^{3} - A^{2} - A + 1$}
\end{minipage}

\noindent{\color{gray!40}\rule{\textwidth}{0.4pt}}
\vspace{0.9\baselineskip}
\noindent \begin{minipage}[t]{0.25\textwidth}
\vspace{0pt}
\centering
\includegraphics[page=41,width=\linewidth]{knotoids.pdf}
\end{minipage}
\hfill
\begin{minipage}[t]{0.73\textwidth}
\vspace{0pt}
\raggedright
\textbf{Name:} {\large{$\mathbf{K6_{4}}$}} (chiral, non-rotatable$^{*}$) \\ \textbf{PD:} {\scriptsize\texttt{[0],[0,1,2,3],[1,4,5,2],[6,7,8,3],[4,9,10,5],[6],[7,11,12,8],[9,12,11,10]}} \\ \textbf{EM:} {\small\texttt{(B0, A0C0C3D3, B1E0E3B2, F0G0G3B3, C1H0H3C2, D0, D1H2H1D2, E1G2G1E2)}} \\ \textbf{Kauffman bracket:} {\scriptsize $A^{24} + A^{22} - A^{18} + A^{14} - A^{10} + A^{6} - A^{2}$} \\ \textbf{Arrow:} {\scriptsize $A^{-12} + L_1/A^{14} - L_1/A^{18} + L_1/A^{22} - L_1/A^{26} + L_1/A^{30} - L_1/A^{34}$} \\ \textbf{Mock:} {\scriptsize $w^{6} + w^{5} - w^{3} + w - 1/w + w^{-3} - 1/w^{5}$} \\ \textbf{Affine:} {\scriptsize $3t - 6 + 3/t$} \\ \textbf{Yamada:} {\scriptsize $A^{23} + A^{22} + A^{21} + A^{20} + 2A^{19} + A^{18} + A^{17} + A^{16} + A^{15} - A^{11} - A^{10} - A^{8} - A^{6} - A^{4} - A^{3} + A^{2} + 1$}
\end{minipage}

\noindent{\color{gray!40}\rule{\textwidth}{0.4pt}}
\vspace{0.9\baselineskip}
\noindent \begin{minipage}[t]{0.25\textwidth}
\vspace{0pt}
\centering
\includegraphics[page=42,width=\linewidth]{knotoids.pdf}
\end{minipage}
\hfill
\begin{minipage}[t]{0.73\textwidth}
\vspace{0pt}
\raggedright
\textbf{Name:} {\large{$\mathbf{K6_{5}}$}} (chiral, non-rotatable$^{*}$) \\ \textbf{PD:} {\scriptsize\texttt{[0],[0,1,2,3],[1,4,5,2],[3,5,6,7],[8,9,10,4],[11,12,7,6],[12,11,9,8],[10]}} \\ \textbf{EM:} {\small\texttt{(B0, A0C0C3D0, B1E3D1B2, B3C2F3F2, G3G2H0C1, G1G0D3D2, F1F0E1E0, E2)}} \\ \textbf{Kauffman bracket:} {\scriptsize $-A^{23} + A^{19} - 2A^{15} + 2A^{11} - 2A^{7} - A^{5} + A^{3} + A$} \\ \textbf{Arrow:} {\scriptsize $A^{-4} - 1/A^{8} + 2/A^{12} - 2/A^{16} + 2/A^{20} + L_1/A^{22} - 1/A^{24} - L_1/A^{26}$} \\ \textbf{Mock:} {\scriptsize $3w^{2} + w - 4 - 1/w + 2/w^{2}$} \\ \textbf{Affine:} {\scriptsize $t - 2 + 1/t$} \\ \textbf{Yamada:} {\scriptsize $A^{23} + 2A^{21} + 2A^{20} - A^{19} + 2A^{18} + 2A^{17} - A^{16} + A^{15} + A^{13} - A^{12} + A^{10} - 2A^{9} - A^{8} + A^{7} - A^{6} - A^{5} + A^{4} - A^{2} + 1$}
\end{minipage}

\noindent{\color{gray!40}\rule{\textwidth}{0.4pt}}
\vspace{0.9\baselineskip}
\noindent \begin{minipage}[t]{0.25\textwidth}
\vspace{0pt}
\centering
\includegraphics[page=43,width=\linewidth]{knotoids.pdf}
\end{minipage}
\hfill
\begin{minipage}[t]{0.73\textwidth}
\vspace{0pt}
\raggedright
\textbf{Name:} {\large{$\mathbf{K6_{6}}$}} (chiral, non-rotatable$^{*}$) \\ \textbf{PD:} {\scriptsize\texttt{[0],[0,1,2,3],[1,4,5,2],[3,6,7,8],[4,9,6,5],[10,11,8,7],[11,10,12,9],[12]}} \\ \textbf{EM:} {\small\texttt{(B0, A0C0C3D0, B1E0E3B2, B3E2F3F2, C1G3D1C2, G1G0D3D2, F1F0H0E1, G2)}} \\ \textbf{Kauffman bracket:} {\scriptsize $-A^{24} + A^{20} - 2A^{16} + 3A^{12} + A^{10} - 2A^{8} - A^{6} + A^{4} + A^{2}$} \\ \textbf{Arrow:} {\scriptsize $-A^{30}L_1 + A^{26}L_1 - 2A^{22}L_1 + 3A^{18}L_1 + A^{16} - 2A^{14}L_1 - A^{12} + A^{10}L_1 + A^{8}$} \\ \textbf{Mock:} {\scriptsize $-2w^{3} + 3w - 3/w - 1/w^{2} + 2/w^{3} + 2/w^{4}$} \\ \textbf{Affine:} {\scriptsize $-3t + 6 - 3/t$} \\ \textbf{Yamada:} {\scriptsize $-A^{23} + A^{22} - 2A^{20} + A^{19} + 2A^{18} - 2A^{17} + 2A^{16} + 3A^{15} - A^{14} + A^{13} - A^{12} + A^{11} - 2A^{10} - A^{9} + 2A^{8} - 3A^{7} - 2A^{6} + A^{5} - 2A^{4} - 2A^{3} - 1$}
\end{minipage}

\noindent{\color{gray!40}\rule{\textwidth}{0.4pt}}
\vspace{0.9\baselineskip}
\noindent \begin{minipage}[t]{0.25\textwidth}
\vspace{0pt}
\centering
\includegraphics[page=44,width=\linewidth]{knotoids.pdf}
\end{minipage}
\hfill
\begin{minipage}[t]{0.73\textwidth}
\vspace{0pt}
\raggedright
\textbf{Name:} {\large{$\mathbf{K6_{7}}$}} (chiral, non-rotatable$^{*}$) \\ \textbf{PD:} {\scriptsize\texttt{[0],[0,1,2,3],[1,4,5,2],[3,6,7,8],[4,9,10,5],[6,11,12,8],[7],[9,12,11,10]}} \\ \textbf{EM:} {\small\texttt{(B0, A0C0C3D0, B1E0E3B2, B3F0G0F3, C1H0H3C2, D1H2H1D3, D2, E1F2F1E2)}} \\ \textbf{Kauffman bracket:} {\scriptsize $A^{23} - A^{19} - A^{17} + A^{15} + A^{13} - A^{11} - A^{9} + A^{7} + A^{5} - A^{3} - A$} \\ \textbf{Arrow:} {\scriptsize $-A^{26}L_1 + A^{22}L_1 + A^{20} - A^{18}L_1 - A^{16} + A^{14}L_1 + A^{12} - A^{10}L_1 - A^{8} + A^{6}L_1 + A^{4}$} \\ \textbf{Mock:} {\scriptsize $-3w - 2 + 3/w + 3/w^{2}$} \\ \textbf{Affine:} {\scriptsize $-3t + 6 - 3/t$} \\ \textbf{Yamada:} {\scriptsize $A^{24} - A^{22} - A^{19} + A^{17} + A^{14} + A^{11} - 2A^{8} - A^{6} - 2A^{5} - A^{3} - A^{2} - 1$}
\end{minipage}

\noindent{\color{gray!40}\rule{\textwidth}{0.4pt}}
\vspace{0.9\baselineskip}
\noindent \begin{minipage}[t]{0.25\textwidth}
\vspace{0pt}
\centering
\includegraphics[page=45,width=\linewidth]{knotoids.pdf}
\end{minipage}
\hfill
\begin{minipage}[t]{0.73\textwidth}
\vspace{0pt}
\raggedright
\textbf{Name:} {\large{$\mathbf{K6_{8}}$}} (chiral, non-rotatable$^{*}$) \\ \textbf{PD:} {\scriptsize\texttt{[0],[0,1,2,3],[1,4,5,2],[6,7,4,3],[8,9,10,5],[11,12,7,6],[12,11,9,8],[10]}} \\ \textbf{EM:} {\small\texttt{(B0, A0C0C3D3, B1D2E3B2, F3F2C1B3, G3G2H0C2, G1G0D1D0, F1F0E1E0, E2)}} \\ \textbf{Kauffman bracket:} {\scriptsize $-A^{24} - A^{22} + A^{20} + 2A^{18} - A^{16} - 2A^{14} + A^{12} + 2A^{10} - A^{6} + A^{2}$} \\ \textbf{Arrow:} {\scriptsize $-A^{30}L_1 - A^{28} + A^{26}L_1 + 2A^{24} - A^{22}L_1 - 2A^{20} + A^{18}L_1 + 2A^{16} - A^{12} + A^{8}$} \\ \textbf{Mock:} {\scriptsize $w^{4} - w^{3} - 2w^{2} + w + 2 - 1/w - 2/w^{2} + w^{-3} + 2/w^{4}$} \\ \textbf{Affine:} {\scriptsize $-2t + 4 - 2/t$} \\ \textbf{Yamada:} {\scriptsize $-A^{23} - A^{20} + A^{19} + 2A^{18} - 2A^{17} + 2A^{16} + 2A^{15} - A^{14} + 2A^{13} + A^{11} - 2A^{10} - A^{9} + A^{8} - 2A^{7} - A^{6} + A^{5} - 2A^{4} - 2A^{3} - A - 2$}
\end{minipage}

\noindent{\color{gray!40}\rule{\textwidth}{0.4pt}}
\vspace{0.9\baselineskip}
\noindent \begin{minipage}[t]{0.25\textwidth}
\vspace{0pt}
\centering
\includegraphics[page=46,width=\linewidth]{knotoids.pdf}
\end{minipage}
\hfill
\begin{minipage}[t]{0.73\textwidth}
\vspace{0pt}
\raggedright
\textbf{Name:} {\large{$\mathbf{K6_{9}}$}} (chiral, non-rotatable$^{*}$) \\ \textbf{PD:} {\scriptsize\texttt{[0],[0,1,2,3],[1,4,5,2],[3,5,6,7],[4,8,9,10],[10,11,7,6],[8,11,12,9],[12]}} \\ \textbf{EM:} {\small\texttt{(B0, A0C0C3D0, B1E0D1B2, B3C2F3F2, C1G0G3F0, E3G1D3D2, E1F1H0E2, G2)}} \\ \textbf{Kauffman bracket:} {\scriptsize $A^{22} + A^{20} - A^{18} - 2A^{16} + 2A^{14} + 3A^{12} - A^{10} - 3A^{8} + 2A^{4} - 1$} \\ \textbf{Arrow:} {\scriptsize $A^{-8} + L_1/A^{10} - 1/A^{12} - 2L_1/A^{14} + 2/A^{16} + 3L_1/A^{18} - 1/A^{20} - 3L_1/A^{22} + 2L_1/A^{26} - L_1/A^{30}$} \\ \textbf{Mock:} {\scriptsize $2w^{4} + 2w^{3} - 2w^{2} - 4w + 1 + 4/w - 2/w^{3}$} \\ \textbf{Affine:} {\scriptsize $2t - 4 + 2/t$} \\ \textbf{Yamada:} {\scriptsize $A^{23} + 3A^{20} + A^{19} - 2A^{18} + 5A^{17} + 2A^{16} - 3A^{15} + 3A^{14} + A^{13} - 3A^{12} - A^{10} + A^{9} - 3A^{8} + 4A^{6} - 4A^{5} + 3A^{3} - 3A^{2} - A + 2$}
\end{minipage}

\noindent{\color{gray!40}\rule{\textwidth}{0.4pt}}
\vspace{0.9\baselineskip}
\noindent \begin{minipage}[t]{0.25\textwidth}
\vspace{0pt}
\centering
\includegraphics[page=47,width=\linewidth]{knotoids.pdf}
\end{minipage}
\hfill
\begin{minipage}[t]{0.73\textwidth}
\vspace{0pt}
\raggedright
\textbf{Name:} {\large{$\mathbf{K6_{10}}$}} (chiral, non-rotatable$^{*}$) \\ \textbf{PD:} {\scriptsize\texttt{[0],[0,1,2,3],[1,4,5,2],[6,7,4,3],[7,8,9,5],[10,11,12,6],[8,12,11,9],[10]}} \\ \textbf{EM:} {\small\texttt{(B0, A0C0C3D3, B1D2E3B2, F3E0C1B3, D1G0G3C2, H0G2G1D0, E1F2F1E2, F0)}} \\ \textbf{Kauffman bracket:} {\scriptsize $A^{23} + A^{21} - A^{19} - 2A^{17} + A^{15} + 2A^{13} - A^{11} - 3A^{9} + 2A^{5} - A$} \\ \textbf{Arrow:} {\scriptsize $-A^{26}L_1 - A^{24} + A^{22}L_1 + 2A^{20} - A^{18}L_1 - 2A^{16} + A^{14}L_1 + 3A^{12} - 2A^{8} + A^{4}$} \\ \textbf{Mock:} {\scriptsize $2w^{2} - 2w - 5 + 2/w + 4/w^{2}$} \\ \textbf{Affine:} {\scriptsize $-2t + 4 - 2/t$} \\ \textbf{Yamada:} {\scriptsize $-A^{23} + A^{21} - A^{20} - A^{19} + 3A^{18} - A^{16} + 3A^{15} + A^{14} - 2A^{13} + A^{12} + A^{11} - 2A^{10} - A^{8} + A^{7} - 4A^{6} - A^{5} + 2A^{4} - 4A^{3} + A - 2$}
\end{minipage}

\noindent{\color{gray!40}\rule{\textwidth}{0.4pt}}
\vspace{0.9\baselineskip}
\noindent \begin{minipage}[t]{0.25\textwidth}
\vspace{0pt}
\centering
\includegraphics[page=48,width=\linewidth]{knotoids.pdf}
\end{minipage}
\hfill
\begin{minipage}[t]{0.73\textwidth}
\vspace{0pt}
\raggedright
\textbf{Name:} {\large{$\mathbf{K6_{11}}$}} (chiral, non-rotatable$^{*}$) \\ \textbf{PD:} {\scriptsize\texttt{[0],[0,1,2,3],[1,4,5,2],[3,6,7,8],[4,9,6,5],[10,11,8,7],[9,11,12,10],[12]}} \\ \textbf{EM:} {\small\texttt{(B0, A0C0C3D0, B1E0E3B2, B3E2F3F2, C1G0D1C2, G3G1D3D2, E1F1H0F0, G2)}} \\ \textbf{Kauffman bracket:} {\scriptsize $A^{22} + A^{20} - A^{18} - A^{16} + 2A^{14} + 2A^{12} - 2A^{10} - 3A^{8} + A^{6} + 2A^{4} - 1$} \\ \textbf{Arrow:} {\scriptsize $A^{-8} + L_1/A^{10} - 1/A^{12} - L_1/A^{14} + 2/A^{16} + 2L_1/A^{18} - 2/A^{20} - 3L_1/A^{22} + A^{-24} + 2L_1/A^{26} - L_1/A^{30}$} \\ \textbf{Mock:} {\scriptsize $2w^{4} + 2w^{3} - 2w^{2} - 3w + 2 + 3/w - 1/w^{2} - 2/w^{3}$} \\ \textbf{Affine:} {\scriptsize $3t - 6 + 3/t$} \\ \textbf{Yamada:} {\scriptsize $-A^{24} - A^{22} - 3A^{21} - 5A^{18} + A^{17} + 2A^{16} - 3A^{15} + A^{14} + 2A^{13} - A^{12} + A^{11} + 3A^{9} - A^{8} - 2A^{7} + 4A^{6} - 3A^{5} - 2A^{4} + 3A^{3} - A^{2} - A + 1$}
\end{minipage}

\noindent{\color{gray!40}\rule{\textwidth}{0.4pt}}
\vspace{0.9\baselineskip}
\noindent \begin{minipage}[t]{0.25\textwidth}
\vspace{0pt}
\centering
\includegraphics[page=49,width=\linewidth]{knotoids.pdf}
\end{minipage}
\hfill
\begin{minipage}[t]{0.73\textwidth}
\vspace{0pt}
\raggedright
\textbf{Name:} {\large{$\mathbf{K6_{12}}$}} (chiral, non-rotatable$^{*}$) \\ \textbf{PD:} {\scriptsize\texttt{[0],[0,1,2,3],[1,4,5,2],[3,5,6,7],[4,8,9,6],[7,10,11,12],[8,12,10,9],[11]}} \\ \textbf{EM:} {\small\texttt{(B0, A0C0C3D0, B1E0D1B2, B3C2E3F0, C1G0G3D2, D3G2H0G1, E1F3F1E2, F2)}} \\ \textbf{Kauffman bracket:} {\scriptsize $A^{22} - 2A^{18} - A^{16} + 3A^{14} + 2A^{12} - 2A^{10} - 2A^{8} + 2A^{6} + 2A^{4} - A^{2} - 1$} \\ \textbf{Arrow:} {\scriptsize $A^{4} - 2 - L_1/A^{2} + 3/A^{4} + 2L_1/A^{6} - 2/A^{8} - 2L_1/A^{10} + 2/A^{12} + 2L_1/A^{14} - 1/A^{16} - L_1/A^{18}$} \\ \textbf{Mock:} {\scriptsize $-w^{4} - w^{3} + 4w^{2} + 3w - 4 - 3/w + 2/w^{2} + w^{-3}$} \\ \textbf{Affine:} {\scriptsize $0$} \\ \textbf{Yamada:} {\scriptsize $-A^{24} + A^{23} + 2A^{22} - 3A^{21} + 4A^{19} - 3A^{18} + 5A^{16} - A^{15} + 3A^{13} + A^{11} - 2A^{10} + 3A^{9} - 5A^{7} + 4A^{6} - 4A^{4} + 2A^{3} - A + 1$}
\end{minipage}

\noindent{\color{gray!40}\rule{\textwidth}{0.4pt}}
\vspace{0.9\baselineskip}
\noindent \begin{minipage}[t]{0.25\textwidth}
\vspace{0pt}
\centering
\includegraphics[page=50,width=\linewidth]{knotoids.pdf}
\end{minipage}
\hfill
\begin{minipage}[t]{0.73\textwidth}
\vspace{0pt}
\raggedright
\textbf{Name:} {\large{$\mathbf{K6_{13}}$}} (chiral, non-rotatable$^{*}$) \\ \textbf{PD:} {\scriptsize\texttt{[0],[0,1,2,3],[1,4,5,2],[3,5,6,7],[4,8,9,10],[10,11,7,6],[11,9,12,8],[12]}} \\ \textbf{EM:} {\small\texttt{(B0, A0C0C3D0, B1E0D1B2, B3C2F3F2, C1G3G1F0, E3G0D3D2, F1E2H0E1, G2)}} \\ \textbf{Kauffman bracket:} {\scriptsize $-A^{24} + 2A^{20} + A^{18} - 3A^{16} - A^{14} + 3A^{12} + 2A^{10} - 2A^{8} - 2A^{6} + A^{4} + A^{2}$} \\ \textbf{Arrow:} {\scriptsize $-A^{18}L_1 + 2A^{14}L_1 + A^{12} - 3A^{10}L_1 - A^{8} + 3A^{6}L_1 + 2A^{4} - 2A^{2}L_1 - 2 + L_1/A^{2} + A^{-4}$} \\ \textbf{Mock:} {\scriptsize $w^{3} - 5w - 2 + 5/w + 4/w^{2} - 1/w^{3} - 1/w^{4}$} \\ \textbf{Affine:} {\scriptsize $-2t + 4 - 2/t$} \\ \textbf{Yamada:} {\scriptsize $A^{24} - A^{23} + 3A^{21} - 4A^{20} - A^{19} + 5A^{18} - 4A^{17} - A^{16} + 3A^{15} - 2A^{14} - A^{12} + 2A^{11} - A^{10} - 4A^{9} + 3A^{8} - A^{7} - 5A^{6} + 3A^{5} + A^{4} - 3A^{3} + A^{2} + A - 1$}
\end{minipage}

\noindent{\color{gray!40}\rule{\textwidth}{0.4pt}}
\vspace{0.9\baselineskip}
\noindent \begin{minipage}[t]{0.25\textwidth}
\vspace{0pt}
\centering
\includegraphics[page=51,width=\linewidth]{knotoids.pdf}
\end{minipage}
\hfill
\begin{minipage}[t]{0.73\textwidth}
\vspace{0pt}
\raggedright
\textbf{Name:} {\large{$\mathbf{K6_{14}}$}} (chiral, non-rotatable$^{*}$) \\ \textbf{PD:} {\scriptsize\texttt{[0],[0,1,2,3],[1,4,5,2],[3,5,6,7],[4,8,9,6],[7,9,10,11],[8,11,12,10],[12]}} \\ \textbf{EM:} {\small\texttt{(B0, A0C0C3D0, B1E0D1B2, B3C2E3F0, C1G0F1D2, D3E2G3G1, E1F3H0F2, G2)}} \\ \textbf{Kauffman bracket:} {\scriptsize $A^{22} + A^{20} - 2A^{18} - 2A^{16} + 2A^{14} + 3A^{12} - 2A^{10} - 3A^{8} + A^{6} + 3A^{4} - 1$} \\ \textbf{Arrow:} {\scriptsize $A^{10}L_1 + A^{8} - 2A^{6}L_1 - 2A^{4} + 2A^{2}L_1 + 3 - 2L_1/A^{2} - 3/A^{4} + L_1/A^{6} + 3/A^{8} - 1/A^{12}$} \\ \textbf{Mock:} {\scriptsize $-w^{3} - 2w^{2} + 3w + 6 - 3/w - 4/w^{2} + w^{-3} + w^{-4}$} \\ \textbf{Affine:} {\scriptsize $0$} \\ \textbf{Yamada:} {\scriptsize $-A^{24} + A^{23} + A^{22} - 3A^{21} + 2A^{20} + 3A^{19} - 6A^{18} + A^{17} + 4A^{16} - 5A^{15} + 2A^{13} - 2A^{12} - A^{11} - 2A^{10} + 3A^{9} - 2A^{8} - 4A^{7} + 6A^{6} - 2A^{5} - 4A^{4} + 5A^{3} - A^{2} - 2A + 1$}
\end{minipage}

\noindent{\color{gray!40}\rule{\textwidth}{0.4pt}}
\vspace{0.9\baselineskip}
\noindent \begin{minipage}[t]{0.25\textwidth}
\vspace{0pt}
\centering
\includegraphics[page=52,width=\linewidth]{knotoids.pdf}
\end{minipage}
\hfill
\begin{minipage}[t]{0.73\textwidth}
\vspace{0pt}
\raggedright
\textbf{Name:} {\large{$\mathbf{K6_{15}}$}} (chiral, non-rotatable$^{*}$) \\ \textbf{PD:} {\scriptsize\texttt{[0],[0,1,2,3],[1,4,5,2],[3,5,6,7],[4,8,9,6],[10,11,8,7],[11,10,12,9],[12]}} \\ \textbf{EM:} {\small\texttt{(B0, A0C0C3D0, B1E0D1B2, B3C2E3F3, C1F2G3D2, G1G0E1D3, F1F0H0E2, G2)}} \\ \textbf{Kauffman bracket:} {\scriptsize $A^{22} + A^{20} - 2A^{18} - 2A^{16} + 2A^{14} + 4A^{12} - A^{10} - 3A^{8} + 2A^{4} - 1$} \\ \textbf{Arrow:} {\scriptsize $L_1/A^{2} + A^{-4} - 2L_1/A^{6} - 2/A^{8} + 2L_1/A^{10} + 4/A^{12} - L_1/A^{14} - 3/A^{16} + 2/A^{20} - 1/A^{24}$} \\ \textbf{Mock:} {\scriptsize $w^{3} + 3w^{2} - 2w - 5 + 2/w + 4/w^{2} - 1/w^{3} - 1/w^{4}$} \\ \textbf{Affine:} {\scriptsize $t - 2 + 1/t$} \\ \textbf{Yamada:} {\scriptsize $A^{23} - A^{22} - A^{21} + 4A^{20} - 4A^{18} + 5A^{17} + 2A^{16} - 4A^{15} + 4A^{14} + 3A^{13} - 2A^{12} + A^{11} + 2A^{9} - 4A^{8} + 5A^{6} - 6A^{5} - A^{4} + 4A^{3} - 3A^{2} - A + 2$}
\end{minipage}

\noindent{\color{gray!40}\rule{\textwidth}{0.4pt}}
\vspace{0.9\baselineskip}
\noindent \begin{minipage}[t]{0.25\textwidth}
\vspace{0pt}
\centering
\includegraphics[page=53,width=\linewidth]{knotoids.pdf}
\end{minipage}
\hfill
\begin{minipage}[t]{0.73\textwidth}
\vspace{0pt}
\raggedright
\textbf{Name:} {\large{$\mathbf{K6_{16}}$}} (chiral, non-rotatable$^{*}$) \\ \textbf{PD:} {\scriptsize\texttt{[0],[0,1,2,3],[1,4,5,2],[3,6,7,4],[8,9,10,5],[6,10,11,7],[12,11,9,8],[12]}} \\ \textbf{EM:} {\small\texttt{(B0, A0C0C3D0, B1D3E3B2, B3F0F3C1, G3G2F1C2, D1E2G1D2, H0F2E1E0, G0)}} \\ \textbf{Kauffman bracket:} {\scriptsize $A^{22} - 2A^{18} - A^{16} + 3A^{14} + 3A^{12} - 2A^{10} - 3A^{8} + A^{6} + 2A^{4} - 1$} \\ \textbf{Arrow:} {\scriptsize $A^{4} - 2 - L_1/A^{2} + 3/A^{4} + 3L_1/A^{6} - 2/A^{8} - 3L_1/A^{10} + A^{-12} + 2L_1/A^{14} - L_1/A^{18}$} \\ \textbf{Mock:} {\scriptsize $-w^{4} - w^{3} + 4w^{2} + 4w - 3 - 4/w + w^{-2} + w^{-3}$} \\ \textbf{Affine:} {\scriptsize $t - 2 + 1/t$} \\ \textbf{Yamada:} {\scriptsize $A^{23} - A^{22} - 3A^{21} + 3A^{20} + 2A^{19} - 6A^{18} + 2A^{17} + 3A^{16} - 6A^{15} + 2A^{13} - 3A^{12} - A^{11} - A^{10} + 3A^{9} - 3A^{8} - 2A^{7} + 7A^{6} - 3A^{5} - 3A^{4} + 5A^{3} - A^{2} - 2A + 1$}
\end{minipage}

\noindent{\color{gray!40}\rule{\textwidth}{0.4pt}}
\vspace{0.9\baselineskip}
\noindent \begin{minipage}[t]{0.25\textwidth}
\vspace{0pt}
\centering
\includegraphics[page=54,width=\linewidth]{knotoids.pdf}
\end{minipage}
\hfill
\begin{minipage}[t]{0.73\textwidth}
\vspace{0pt}
\raggedright
\textbf{Name:} {\large{$\mathbf{K6_{17}}$}} (chiral, non-rotatable$^{*}$) \\ \textbf{PD:} {\scriptsize\texttt{[0],[0,1,2,3],[1,4,5,2],[3,6,7,8],[4,9,10,5],[6,10,11,7],[8,11,12,9],[12]}} \\ \textbf{EM:} {\small\texttt{(B0, A0C0C3D0, B1E0E3B2, B3F0F3G0, C1G3F1C2, D1E2G1D2, D3F2H0E1, G2)}} \\ \textbf{Kauffman bracket:} {\scriptsize $A^{22} - A^{18} + 3A^{14} + A^{12} - 3A^{10} - 2A^{8} + 2A^{6} + 2A^{4} - A^{2} - 1$} \\ \textbf{Arrow:} {\scriptsize $A^{-8} - 1/A^{12} + 3/A^{16} + L_1/A^{18} - 3/A^{20} - 2L_1/A^{22} + 2/A^{24} + 2L_1/A^{26} - 1/A^{28} - L_1/A^{30}$} \\ \textbf{Mock:} {\scriptsize $2w^{4} + w^{3} - 3w^{2} - 2w + 3 + 2/w - 2/w^{2} - 1/w^{3} + w^{-4}$} \\ \textbf{Affine:} {\scriptsize $t - 2 + 1/t$} \\ \textbf{Yamada:} {\scriptsize $2A^{23} + A^{22} - A^{21} + 3A^{20} + 4A^{19} - 3A^{18} + 2A^{17} + 4A^{16} - 3A^{15} + 2A^{13} - 2A^{12} - A^{11} - 2A^{10} + 2A^{9} - 2A^{8} - 3A^{7} + 5A^{6} - 2A^{5} - 3A^{4} + 3A^{3} - A + 1$}
\end{minipage}

\noindent{\color{gray!40}\rule{\textwidth}{0.4pt}}
\vspace{0.9\baselineskip}
\noindent \begin{minipage}[t]{0.25\textwidth}
\vspace{0pt}
\centering
\includegraphics[page=55,width=\linewidth]{knotoids.pdf}
\end{minipage}
\hfill
\begin{minipage}[t]{0.73\textwidth}
\vspace{0pt}
\raggedright
\textbf{Name:} {\large{$\mathbf{K6_{18}}$}} (chiral, non-rotatable$^{*}$) \\ \textbf{PD:} {\scriptsize\texttt{[0],[0,1,2,3],[1,4,5,2],[3,6,7,4],[8,9,6,5],[9,10,11,7],[12,11,10,8],[12]}} \\ \textbf{EM:} {\small\texttt{(B0, A0C0C3D0, B1D3E3B2, B3E2F3C1, G3F0D1C2, E1G2G1D2, H0F2F1E0, G0)}} \\ \textbf{Kauffman bracket:} {\scriptsize $A^{22} - 2A^{18} - A^{16} + 3A^{14} + 2A^{12} - 3A^{10} - 3A^{8} + 2A^{6} + 3A^{4} - 1$} \\ \textbf{Arrow:} {\scriptsize $A^{16} - 2A^{12} - A^{10}L_1 + 3A^{8} + 2A^{6}L_1 - 3A^{4} - 3A^{2}L_1 + 2 + 3L_1/A^{2} - L_1/A^{6}$} \\ \textbf{Mock:} {\scriptsize $w^{4} + w^{3} - 4w^{2} - 4w + 5 + 4/w - 1/w^{2} - 1/w^{3}$} \\ \textbf{Affine:} {\scriptsize $-t + 2 - 1/t$} \\ \textbf{Yamada:} {\scriptsize $A^{23} - A^{22} - 2A^{21} + 4A^{20} - 6A^{18} + 5A^{17} + 2A^{16} - 7A^{15} + 3A^{14} + 2A^{13} - 3A^{12} + 2A^{9} - 5A^{8} - A^{7} + 6A^{6} - 6A^{5} - 2A^{4} + 6A^{3} - 4A^{2} - 2A + 2$}
\end{minipage}

\noindent{\color{gray!40}\rule{\textwidth}{0.4pt}}
\vspace{0.9\baselineskip}
\noindent \begin{minipage}[t]{0.25\textwidth}
\vspace{0pt}
\centering
\includegraphics[page=56,width=\linewidth]{knotoids.pdf}
\end{minipage}
\hfill
\begin{minipage}[t]{0.73\textwidth}
\vspace{0pt}
\raggedright
\textbf{Name:} {\large{$\mathbf{K6_{19}}$}} (chiral, non-rotatable$^{*}$) \\ \textbf{PD:} {\scriptsize\texttt{[0],[0,1,2,3],[1,4,5,2],[3,5,6,7],[4,8,9,6],[10,11,8,7],[9,12,11,10],[12]}} \\ \textbf{EM:} {\small\texttt{(B0, A0C0C3D0, B1E0D1B2, B3C2E3F3, C1F2G0D2, G3G2E1D3, E2H0F1F0, G1)}} \\ \textbf{Kauffman bracket:} {\scriptsize $A^{22} + A^{20} - A^{18} - 2A^{16} + A^{14} + 3A^{12} - 3A^{8} - A^{6} + A^{4} + A^{2}$} \\ \textbf{Arrow:} {\scriptsize $A^{4} + A^{2}L_1 - 1 - 2L_1/A^{2} + A^{-4} + 3L_1/A^{6} - 3L_1/A^{10} - 1/A^{12} + L_1/A^{14} + A^{-16}$} \\ \textbf{Mock:} {\scriptsize $-w^{4} - w^{3} + 2w^{2} + 4w + 1 - 4/w - 1/w^{2} + w^{-3}$} \\ \textbf{Affine:} {\scriptsize $t - 2 + 1/t$} \\ \textbf{Yamada:} {\scriptsize $-A^{24} + A^{22} - 2A^{21} - 2A^{20} + 4A^{19} - A^{18} - 4A^{17} + 3A^{16} - 3A^{14} + A^{13} - 3A^{10} + A^{9} + 2A^{8} - 5A^{7} + 2A^{6} + 3A^{5} - 3A^{4} + 2A^{2} - 1$}
\end{minipage}

\noindent{\color{gray!40}\rule{\textwidth}{0.4pt}}
\vspace{0.9\baselineskip}
\noindent \begin{minipage}[t]{0.25\textwidth}
\vspace{0pt}
\centering
\includegraphics[page=57,width=\linewidth]{knotoids.pdf}
\end{minipage}
\hfill
\begin{minipage}[t]{0.73\textwidth}
\vspace{0pt}
\raggedright
\textbf{Name:} {\large{$\mathbf{K6_{20}}$}} (chiral, non-rotatable$^{*}$) \\ \textbf{PD:} {\scriptsize\texttt{[0],[0,1,2,3],[1,4,5,2],[3,6,7,8],[4,8,9,5],[6,10,11,7],[12,11,10,9],[12]}} \\ \textbf{EM:} {\small\texttt{(B0, A0C0C3D0, B1E0E3B2, B3F0F3E1, C1D3G3C2, D1G2G1D2, H0F2F1E2, G0)}} \\ \textbf{Kauffman bracket:} {\scriptsize $A^{23} - A^{19} - A^{17} + 2A^{15} + 2A^{13} - 2A^{11} - 3A^{9} + 2A^{5} - A$} \\ \textbf{Arrow:} {\scriptsize $-A^{26}L_1 + A^{22}L_1 + A^{20} - 2A^{18}L_1 - 2A^{16} + 2A^{14}L_1 + 3A^{12} - 2A^{8} + A^{4}$} \\ \textbf{Mock:} {\scriptsize $w^{2} - 3w - 4 + 3/w + 4/w^{2}$} \\ \textbf{Affine:} {\scriptsize $-3t + 6 - 3/t$} \\ \textbf{Yamada:} {\scriptsize $A^{22} - 2A^{20} + 2A^{18} - 3A^{17} - A^{16} + 4A^{15} - A^{13} + 4A^{12} + A^{11} - 2A^{10} + A^{9} - A^{8} - A^{7} - 5A^{6} + A^{5} + A^{4} - 5A^{3} + 2A^{2} + A - 3$}
\end{minipage}

\noindent{\color{gray!40}\rule{\textwidth}{0.4pt}}
\vspace{0.9\baselineskip}
\noindent \begin{minipage}[t]{0.25\textwidth}
\vspace{0pt}
\centering
\includegraphics[page=58,width=\linewidth]{knotoids.pdf}
\end{minipage}
\hfill
\begin{minipage}[t]{0.73\textwidth}
\vspace{0pt}
\raggedright
\textbf{Name:} {\large{$\mathbf{K6_{21}}$}} (chiral, non-rotatable$^{*}$) \\ \textbf{PD:} {\scriptsize\texttt{[0],[0,1,2,3],[1,4,5,2],[3,5,6,7],[4,7,8,9],[9,10,11,6],[12,11,10,8],[12]}} \\ \textbf{EM:} {\small\texttt{(B0, A0C0C3D0, B1E0D1B2, B3C2F3E1, C1D3G3F0, E3G2G1D2, H0F2F1E2, G0)}} \\ \textbf{Kauffman bracket:} {\scriptsize $A^{22} - 2A^{18} - 2A^{16} + 3A^{14} + 3A^{12} - 2A^{10} - 3A^{8} + A^{6} + 3A^{4} - 1$} \\ \textbf{Arrow:} {\scriptsize $A^{16} - 2A^{12} - 2A^{10}L_1 + 3A^{8} + 3A^{6}L_1 - 2A^{4} - 3A^{2}L_1 + 1 + 3L_1/A^{2} - L_1/A^{6}$} \\ \textbf{Mock:} {\scriptsize $w^{4} + w^{3} - 4w^{2} - 5w + 4 + 5/w - 1/w^{3}$} \\ \textbf{Affine:} {\scriptsize $-2t + 4 - 2/t$} \\ \textbf{Yamada:} {\scriptsize $2A^{23} - 4A^{21} + 3A^{20} + 2A^{19} - 7A^{18} + 2A^{17} + 4A^{16} - 6A^{15} + A^{14} + 4A^{13} - A^{12} - A^{10} + 3A^{9} - 5A^{8} - 5A^{7} + 6A^{6} - 4A^{5} - 4A^{4} + 6A^{3} - A^{2} - 2A + 1$}
\end{minipage}

\noindent{\color{gray!40}\rule{\textwidth}{0.4pt}}
\vspace{0.9\baselineskip}
\noindent \begin{minipage}[t]{0.25\textwidth}
\vspace{0pt}
\centering
\includegraphics[page=59,width=\linewidth]{knotoids.pdf}
\end{minipage}
\hfill
\begin{minipage}[t]{0.73\textwidth}
\vspace{0pt}
\raggedright
\textbf{Name:} {\large{$\mathbf{K6_{22}}$}} (achiral, non-rotatable$^{*}$) \\ \textbf{PD:} {\scriptsize\texttt{[0],[0,1,2,3],[1,4,5,2],[3,5,6,7],[7,8,9,4],[6,9,10,11],[8,12,11,10],[12]}} \\ \textbf{EM:} {\small\texttt{(B0, A0C0C3D0, B1E3D1B2, B3C2F0E0, D3G0F1C1, D2E2G3G2, E1H0F3F2, G1)}} \\ \textbf{Kauffman bracket:} {\scriptsize $A^{22} + A^{20} - A^{18} - 2A^{16} + 3A^{12} - 2A^{8} - A^{6} + A^{4} + A^{2}$} \\ \textbf{Arrow:} {\scriptsize $A^{10}L_1 + A^{8} - A^{6}L_1 - 2A^{4} + 3 - 2/A^{4} - L_1/A^{6} + A^{-8} + L_1/A^{10}$} \\ \textbf{Mock:} {\scriptsize $-2w^{2} + 5 - 2/w^{2}$} \\ \textbf{Affine:} {\scriptsize $0$} \\ \textbf{Yamada:} {\scriptsize $-A^{24} + A^{22} - A^{21} - 2A^{20} + 3A^{19} + A^{18} - 3A^{17} + 3A^{16} + A^{15} - 3A^{14} - A^{13} - 2A^{12} - A^{11} - 3A^{10} + A^{9} + 3A^{8} - 3A^{7} + A^{6} + 3A^{5} - 2A^{4} - A^{3} + A^{2} - 1$}
\end{minipage}

\noindent{\color{gray!40}\rule{\textwidth}{0.4pt}}
\vspace{0.9\baselineskip}
\noindent \begin{minipage}[t]{0.25\textwidth}
\vspace{0pt}
\centering
\includegraphics[page=60,width=\linewidth]{knotoids.pdf}
\end{minipage}
\hfill
\begin{minipage}[t]{0.73\textwidth}
\vspace{0pt}
\raggedright
\textbf{Name:} {\large{$\mathbf{K6_{23}}$}} (chiral, non-rotatable$^{*}$) \\ \textbf{PD:} {\scriptsize\texttt{[0],[0,1,2,3],[1,4,5,2],[3,6,7,8],[9,10,11,4],[5,11,12,6],[10,9,8,7],[12]}} \\ \textbf{EM:} {\small\texttt{(B0, A0C0C3D0, B1E3F0B2, B3F3G3G2, G1G0F1C1, C2E2H0D1, E1E0D3D2, F2)}} \\ \textbf{Kauffman bracket:} {\scriptsize $A^{19} - A^{15} + A^{11} - A^{9} - 2A^{7} - A^{5} + A^{3} + A$} \\ \textbf{Arrow:} {\scriptsize $-A^{10}L_1 + A^{6}L_1 - A^{2}L_1 + 1 + 2L_1/A^{2} + L_2/A^{4} - L_1/A^{6} - L_2/A^{8}$} \\ \textbf{Mock:} {\scriptsize $-w^{5} - w^{4} + w^{3} + w^{2} - w + 1 + 2/w - 1/w^{3}$} \\ \textbf{Affine:} {\scriptsize $t^{2} - t - 1/t + t^{-2}$} \\ \textbf{Yamada:} {\scriptsize $A^{20} - A^{19} + A^{17} - A^{16} - A^{14} - 2A^{12} - A^{11} - A^{9} - A^{8} - A^{6} - A^{5} + A^{3} + 1$}
\end{minipage}

\noindent{\color{gray!40}\rule{\textwidth}{0.4pt}}
\vspace{0.9\baselineskip}
\noindent \begin{minipage}[t]{0.25\textwidth}
\vspace{0pt}
\centering
\includegraphics[page=61,width=\linewidth]{knotoids.pdf}
\end{minipage}
\hfill
\begin{minipage}[t]{0.73\textwidth}
\vspace{0pt}
\raggedright
\textbf{Name:} {\large{$\mathbf{K6_{24}}$}} (chiral, non-rotatable$^{*}$) \\ \textbf{PD:} {\scriptsize\texttt{[0],[0,1,2,3],[1,4,5,2],[3,6,7,8],[9,10,11,4],[11,12,6,5],[10,9,8,7],[12]}} \\ \textbf{EM:} {\small\texttt{(B0, A0C0C3D0, B1E3F3B2, B3F2G3G2, G1G0F0C1, E2H0D1C2, E1E0D3D2, F1)}} \\ \textbf{Kauffman bracket:} {\scriptsize $-A^{18} + A^{14} - 2A^{10} + 3A^{6} + 2A^{4} - A^{2} - 1$} \\ \textbf{Arrow:} {\scriptsize $-A^{18}L_1 + A^{14}L_1 - 2A^{10}L_1 + 3A^{6}L_1 + A^{4}L_2 + A^{4} - A^{2}L_1 - L_2$} \\ \textbf{Mock:} {\scriptsize $-w^{5} - w^{4} + w^{3} + w^{2} - 2w + 3/w + w^{-2} - 1/w^{3}$} \\ \textbf{Affine:} {\scriptsize $t^{2} - 2t + 2 - 2/t + t^{-2}$} \\ \textbf{Yamada:} {\scriptsize $-A^{23} + A^{22} + A^{21} - 2A^{20} + 3A^{18} - A^{17} + 2A^{15} - A^{14} - 2A^{12} + A^{11} - 2A^{10} - 3A^{9} + A^{8} - 2A^{7} - 3A^{6} + A^{5} - A^{3} + A + 1$}
\end{minipage}

\noindent{\color{gray!40}\rule{\textwidth}{0.4pt}}
\vspace{0.9\baselineskip}
\noindent \begin{minipage}[t]{0.25\textwidth}
\vspace{0pt}
\centering
\includegraphics[page=62,width=\linewidth]{knotoids.pdf}
\end{minipage}
\hfill
\begin{minipage}[t]{0.73\textwidth}
\vspace{0pt}
\raggedright
\textbf{Name:} {\large{$\mathbf{K6_{25}}$}} (chiral, non-rotatable$^{*}$) \\ \textbf{PD:} {\scriptsize\texttt{[0],[0,1,2,3],[1,4,5,2],[6,7,8,3],[4,9,10,11],[5,11,12,6],[7,10,9,8],[12]}} \\ \textbf{EM:} {\small\texttt{(B0, A0C0C3D3, B1E0F0B2, F3G0G3B3, C1G2G1F1, C2E3H0D0, D1E2E1D2, F2)}} \\ \textbf{Kauffman bracket:} {\scriptsize $A^{12} + A^{10} + A^{8} - A^{4} - A^{2}$} \\ \textbf{Arrow:} {\scriptsize $A^{-12} + L_1/A^{14} + L_2/A^{16} - L_2/A^{20} - L_1/A^{22}$} \\ \textbf{Mock:} {\scriptsize $w^{6} + w^{5} - w - 1/w^{2} + w^{-4}$} \\ \textbf{Affine:} {\scriptsize $t^{2} + 2t - 6 + 2/t + t^{-2}$} \\ \textbf{Yamada:} {\scriptsize $-A^{20} - A^{19} - A^{18} - A^{17} - A^{16} - A^{15} - A^{14} - A^{12} + A^{7} + A^{5} - A + 1$}
\end{minipage}

\noindent{\color{gray!40}\rule{\textwidth}{0.4pt}}
\vspace{0.9\baselineskip}
\noindent \begin{minipage}[t]{0.25\textwidth}
\vspace{0pt}
\centering
\includegraphics[page=63,width=\linewidth]{knotoids.pdf}
\end{minipage}
\hfill
\begin{minipage}[t]{0.73\textwidth}
\vspace{0pt}
\raggedright
\textbf{Name:} {\large{$\mathbf{K6_{26}}$}} (chiral, non-rotatable$^{*}$) \\ \textbf{PD:} {\scriptsize\texttt{[0],[0,1,2,3],[1,4,5,2],[6,7,8,3],[4,9,10,11],[11,12,6,5],[7,10,9,8],[12]}} \\ \textbf{EM:} {\small\texttt{(B0, A0C0C3D3, B1E0F3B2, F2G0G3B3, C1G2G1F0, E3H0D0C2, D1E2E1D2, F1)}} \\ \textbf{Kauffman bracket:} {\scriptsize $-2A^{19} - 2A^{17} + A^{15} + 3A^{13} - 2A^{9} + 2A^{5} - A$} \\ \textbf{Arrow:} {\scriptsize $L_2/A^{8} + A^{-8} + 2L_1/A^{10} - L_2/A^{12} - 3L_1/A^{14} + 2L_1/A^{18} - 2L_1/A^{22} + L_1/A^{26}$} \\ \textbf{Mock:} {\scriptsize $w^{4} + 2w^{3} - 3w + 2/w - 1/w^{2} - 2/w^{3} + w^{-4} + w^{-5}$} \\ \textbf{Affine:} {\scriptsize $t^{2} + t - 4 + 1/t + t^{-2}$} \\ \textbf{Yamada:} {\scriptsize $-A^{24} - A^{23} - 2A^{20} - A^{19} + 2A^{18} - 3A^{16} + 2A^{15} - 3A^{13} + A^{12} + A^{11} - 2A^{10} + A^{7} - 2A^{6} + 4A^{4} - 3A^{3} + 2A - 1$}
\end{minipage}

\noindent{\color{gray!40}\rule{\textwidth}{0.4pt}}
\vspace{0.9\baselineskip}
\noindent \begin{minipage}[t]{0.25\textwidth}
\vspace{0pt}
\centering
\includegraphics[page=64,width=\linewidth]{knotoids.pdf}
\end{minipage}
\hfill
\begin{minipage}[t]{0.73\textwidth}
\vspace{0pt}
\raggedright
\textbf{Name:} {\large{$\mathbf{K6_{27}}$}} (chiral, non-rotatable$^{*}$) \\ \textbf{PD:} {\scriptsize\texttt{[0],[0,1,2,3],[1,3,4,5],[2,6,7,4],[8,9,6,5],[10,11,12,7],[11,10,9,8],[12]}} \\ \textbf{EM:} {\small\texttt{(B0, A0C0D0C1, B1B3D3E3, B2E2F3C2, G3G2D1C3, G1G0H0D2, F1F0E1E0, F2)}} \\ \textbf{Kauffman bracket:} {\scriptsize $A^{16} + A^{14} - 2A^{10} + 2A^{6} + A^{4} - A^{2} - 1$} \\ \textbf{Arrow:} {\scriptsize $A^{-8} + L_1/A^{10} + L_2/A^{12} - 1/A^{12} - 2L_1/A^{14} - L_2/A^{16} + A^{-16} + 2L_1/A^{18} + L_2/A^{20} - L_1/A^{22} - L_2/A^{24}$} \\ \textbf{Mock:} {\scriptsize $2w^{4} + 2w^{3} - w^{2} - 3w - 2 + 1/w + 2/w^{2}$} \\ \textbf{Affine:} {\scriptsize $2t^{2} + t - 6 + 1/t + 2/t^{2}$} \\ \textbf{Yamada:} {\scriptsize $-A^{24} - A^{22} - 2A^{21} - 2A^{18} + A^{16} - 2A^{15} + A^{13} - A^{12} + A^{11} - A^{10} + 2A^{9} - A^{8} - A^{7} + 2A^{6} - 2A^{5} - A^{4} + A^{3} + 1$}
\end{minipage}

\noindent{\color{gray!40}\rule{\textwidth}{0.4pt}}
\vspace{0.9\baselineskip}
\noindent \begin{minipage}[t]{0.25\textwidth}
\vspace{0pt}
\centering
\includegraphics[page=65,width=\linewidth]{knotoids.pdf}
\end{minipage}
\hfill
\begin{minipage}[t]{0.73\textwidth}
\vspace{0pt}
\raggedright
\textbf{Name:} {\large{$\mathbf{K6_{28}}$}} (chiral, non-rotatable$^{*}$) \\ \textbf{PD:} {\scriptsize\texttt{[0],[0,1,2,3],[1,3,4,5],[6,7,4,2],[5,8,9,6],[7,10,11,12],[8,11,10,9],[12]}} \\ \textbf{EM:} {\small\texttt{(B0, A0C0D3C1, B1B3D2E0, E3F0C2B2, C3G0G3D0, D1G2G1H0, E1F2F1E2, F3)}} \\ \textbf{Kauffman bracket:} {\scriptsize $-A^{18} + 2A^{14} + A^{12} - 2A^{10} - 2A^{8} + 2A^{6} + 2A^{4} - 1$} \\ \textbf{Arrow:} {\scriptsize $-A^{12}L_2 + A^{8}L_2 + A^{8} + A^{6}L_1 - A^{4}L_2 - A^{4} - 2A^{2}L_1 + L_2 + 1 + 2L_1/A^{2} - L_1/A^{6}$} \\ \textbf{Mock:} {\scriptsize $-w^{2} - w + 4 + 3/w - 2/w^{2} - 2/w^{3}$} \\ \textbf{Affine:} {\scriptsize $-2t^{2} + t + 2 + 1/t - 2/t^{2}$} \\ \textbf{Yamada:} {\scriptsize $A^{22} + A^{21} - A^{20} + 2A^{18} - 2A^{17} - 2A^{16} + 2A^{15} - A^{14} - 3A^{13} + A^{12} - 2A^{10} + A^{9} + A^{7} - 3A^{6} + A^{5} + 2A^{4} - 4A^{3} + A^{2} + A - 2$}
\end{minipage}

\noindent{\color{gray!40}\rule{\textwidth}{0.4pt}}
\vspace{0.9\baselineskip}
\noindent \begin{minipage}[t]{0.25\textwidth}
\vspace{0pt}
\centering
\includegraphics[page=66,width=\linewidth]{knotoids.pdf}
\end{minipage}
\hfill
\begin{minipage}[t]{0.73\textwidth}
\vspace{0pt}
\raggedright
\textbf{Name:} {\large{$\mathbf{K6_{29}}$}} (chiral, non-rotatable$^{*}$) \\ \textbf{PD:} {\scriptsize\texttt{[0],[0,1,2,3],[1,4,5,2],[3,6,7,4],[5,7,8,9],[10,11,12,6],[11,10,9,8],[12]}} \\ \textbf{EM:} {\small\texttt{(B0, A0C0C3D0, B1D3E0B2, B3F3E1C1, C2D2G3G2, G1G0H0D1, F1F0E3E2, F2)}} \\ \textbf{Kauffman bracket:} {\scriptsize $-A^{17} - A^{15} + A^{13} + 2A^{11} - A^{9} - 2A^{7} - A^{5} + A^{3} + A$} \\ \textbf{Arrow:} {\scriptsize $A^{-4} + L_1/A^{6} - 1/A^{8} - 2L_1/A^{10} + A^{-12} + 2L_1/A^{14} + L_2/A^{16} - L_1/A^{18} - L_2/A^{20}$} \\ \textbf{Mock:} {\scriptsize $-w^{3} + w^{2} + 3w - 2/w$} \\ \textbf{Affine:} {\scriptsize $t^{2} + t - 4 + 1/t + t^{-2}$} \\ \textbf{Yamada:} {\scriptsize $-2A^{22} - 3A^{19} + A^{18} + A^{17} - 2A^{16} - A^{14} - 2A^{12} + 2A^{10} - A^{9} + 2A^{7} - A^{6} - A^{5} + A^{4} - A^{2} + 1$}
\end{minipage}

\noindent{\color{gray!40}\rule{\textwidth}{0.4pt}}
\vspace{0.9\baselineskip}
\noindent \begin{minipage}[t]{0.25\textwidth}
\vspace{0pt}
\centering
\includegraphics[page=67,width=\linewidth]{knotoids.pdf}
\end{minipage}
\hfill
\begin{minipage}[t]{0.73\textwidth}
\vspace{0pt}
\raggedright
\textbf{Name:} {\large{$\mathbf{K6_{30}}$}} (chiral, non-rotatable$^{*}$) \\ \textbf{PD:} {\scriptsize\texttt{[0],[0,1,2,3],[1,4,5,2],[3,6,7,8],[8,7,9,4],[5,9,10,11],[11,10,12,6],[12]}} \\ \textbf{EM:} {\small\texttt{(B0, A0C0C3D0, B1E3F0B2, B3G3E1E0, D3D2F1C1, C2E2G1G0, F3F2H0D1, G2)}} \\ \textbf{Kauffman bracket:} {\scriptsize $-A^{15} + A^{13} + A^{11} - A^{9} - 2A^{7} - A^{5} + A^{3} + A$} \\ \textbf{Arrow:} {\scriptsize $1 - L_1/A^{2} - 1/A^{4} + L_1/A^{6} + 2/A^{8} + L_1/A^{10} - 1/A^{12} - L_1/A^{14}$} \\ \textbf{Mock:} {\scriptsize $-w^{3} - w^{2} + w + 3 + 1/w - 1/w^{2} - 1/w^{3}$} \\ \textbf{Affine:} {\scriptsize $0$} \\ \textbf{Yamada:} {\scriptsize $A^{21} + A^{20} - A^{19} - A^{17} - 4A^{16} - 2A^{15} - 2A^{14} - A^{12} + A^{11} + 2A^{10} + A^{8} + A^{7} - A^{6} - A^{5} - A^{2} + 1$}
\end{minipage}

\noindent{\color{gray!40}\rule{\textwidth}{0.4pt}}
\vspace{0.9\baselineskip}
\noindent \begin{minipage}[t]{0.25\textwidth}
\vspace{0pt}
\centering
\includegraphics[page=68,width=\linewidth]{knotoids.pdf}
\end{minipage}
\hfill
\begin{minipage}[t]{0.73\textwidth}
\vspace{0pt}
\raggedright
\textbf{Name:} {\large{$\mathbf{K6_{31}}$}} (chiral, rotatable) \\ \textbf{PD:} {\scriptsize\texttt{[0],[0,1,2,3],[1,4,5,2],[3,6,7,8],[8,7,9,4],[9,10,11,5],[6,11,10,12],[12]}} \\ \textbf{EM:} {\small\texttt{(B0, A0C0C3D0, B1E3F3B2, B3G0E1E0, D3D2F0C1, E2G2G1C2, D1F2F1H0, G3)}} \\ \textbf{Kauffman bracket:} {\scriptsize $A^{16} - A^{12} + 2A^{8} + 2A^{6} - A^{4} - 2A^{2}$} \\ \textbf{Arrow:} {\scriptsize $A^{-8} - 1/A^{12} + 2/A^{16} + 2L_1/A^{18} - 1/A^{20} - 2L_1/A^{22}$} \\ \textbf{Mock:} {\scriptsize $2w^{4} + 2w^{3} - w^{2} - 2w - 1/w^{2} + w^{-4}$} \\ \textbf{Affine:} {\scriptsize $2t - 4 + 2/t$} \\ \textbf{Yamada:} {\scriptsize $-2A^{22} - 2A^{21} - A^{20} - A^{19} - 3A^{18} - A^{17} + 2A^{16} + 2A^{13} - A^{12} - A^{10} + A^{9} - A^{8} - A^{7} + 3A^{6} + A^{3} + A^{2} - A - 1$}
\end{minipage}

\noindent{\color{gray!40}\rule{\textwidth}{0.4pt}}
\vspace{0.9\baselineskip}
\noindent \begin{minipage}[t]{0.25\textwidth}
\vspace{0pt}
\centering
\includegraphics[page=69,width=\linewidth]{knotoids.pdf}
\end{minipage}
\hfill
\begin{minipage}[t]{0.73\textwidth}
\vspace{0pt}
\raggedright
\textbf{Name:} {\large{$\mathbf{K6_{32}}$}} (chiral, non-rotatable$^{*}$) \\ \textbf{PD:} {\scriptsize\texttt{[0],[0,1,2,3],[1,4,5,2],[6,7,4,3],[5,7,8,9],[10,11,12,6],[11,10,9,8],[12]}} \\ \textbf{EM:} {\small\texttt{(B0, A0C0C3D3, B1D2E0B2, F3E1C1B3, C2D1G3G2, G1G0H0D0, F1F0E3E2, F2)}} \\ \textbf{Kauffman bracket:} {\scriptsize $A^{16} + A^{14} - A^{12} - 2A^{10} + A^{8} + 3A^{6} + A^{4} - 2A^{2} - 1$} \\ \textbf{Arrow:} {\scriptsize $A^{-8} + L_1/A^{10} - 1/A^{12} - 2L_1/A^{14} + A^{-16} + 3L_1/A^{18} + L_2/A^{20} - 2L_1/A^{22} - L_2/A^{24}$} \\ \textbf{Mock:} {\scriptsize $2w^{4} + 2w^{3} - w^{2} - 3w - 1 + 2/w + w^{-2} - 1/w^{3}$} \\ \textbf{Affine:} {\scriptsize $t^{2} + 2t - 6 + 2/t + t^{-2}$} \\ \textbf{Yamada:} {\scriptsize $-2A^{22} - A^{21} + A^{20} - A^{19} - 3A^{18} + A^{17} - 4A^{15} + A^{14} + A^{13} - 2A^{12} + A^{11} + 2A^{9} - 2A^{8} + A^{7} + 3A^{6} - 3A^{5} + 2A^{3} - A^{2} - A + 1$}
\end{minipage}

\noindent{\color{gray!40}\rule{\textwidth}{0.4pt}}
\vspace{0.9\baselineskip}
\noindent \begin{minipage}[t]{0.25\textwidth}
\vspace{0pt}
\centering
\includegraphics[page=70,width=\linewidth]{knotoids.pdf}
\end{minipage}
\hfill
\begin{minipage}[t]{0.73\textwidth}
\vspace{0pt}
\raggedright
\textbf{Name:} {\large{$\mathbf{K6_{33}}$}} (chiral, non-rotatable$^{*}$) \\ \textbf{PD:} {\scriptsize\texttt{[0],[0,1,2,3],[1,4,5,2],[6,7,4,3],[7,8,9,5],[6,10,11,12],[8,11,10,9],[12]}} \\ \textbf{EM:} {\small\texttt{(B0, A0C0C3D3, B1D2E3B2, F0E0C1B3, D1G0G3C2, D0G2G1H0, E1F2F1E2, F3)}} \\ \textbf{Kauffman bracket:} {\scriptsize $-2A^{19} - A^{17} + 2A^{15} + 3A^{13} - A^{11} - 3A^{9} + 2A^{5} - A$} \\ \textbf{Arrow:} {\scriptsize $A^{4}L_2 + A^{4} + A^{2}L_1 - L_2 - 1 - 3L_1/A^{2} + A^{-4} + 3L_1/A^{6} - 2L_1/A^{10} + L_1/A^{14}$} \\ \textbf{Mock:} {\scriptsize $-w^{4} - w^{3} + 2w^{2} + 3w - 1 - 4/w + w^{-2} + 2/w^{3}$} \\ \textbf{Affine:} {\scriptsize $t^{2} - t - 1/t + t^{-2}$} \\ \textbf{Yamada:} {\scriptsize $-A^{24} - A^{23} + A^{22} + A^{21} - 2A^{20} + A^{19} + 3A^{18} - 2A^{17} - 2A^{16} + 3A^{15} - A^{14} - 4A^{13} + 2A^{12} - 3A^{10} - 4A^{6} + 2A^{5} + 3A^{4} - 4A^{3} + 2A^{2} + 2A - 2$}
\end{minipage}

\noindent{\color{gray!40}\rule{\textwidth}{0.4pt}}
\vspace{0.9\baselineskip}
\noindent \begin{minipage}[t]{0.25\textwidth}
\vspace{0pt}
\centering
\includegraphics[page=71,width=\linewidth]{knotoids.pdf}
\end{minipage}
\hfill
\begin{minipage}[t]{0.73\textwidth}
\vspace{0pt}
\raggedright
\textbf{Name:} {\large{$\mathbf{K6_{34}}$}} (chiral, rotatable) \\ \textbf{PD:} {\scriptsize\texttt{[0],[0,1,2,3],[1,4,5,2],[6,7,8,3],[4,8,7,9],[9,10,11,5],[6,11,10,12],[12]}} \\ \textbf{EM:} {\small\texttt{(B0, A0C0C3D3, B1E0F3B2, G0E2E1B3, C1D2D1F0, E3G2G1C2, D0F2F1H0, G3)}} \\ \textbf{Kauffman bracket:} {\scriptsize $-2A^{19} - 2A^{17} + 2A^{15} + 3A^{13} - 3A^{9} + 2A^{5} - A$} \\ \textbf{Arrow:} {\scriptsize $2L_1/A^{2} + 2/A^{4} - 2L_1/A^{6} - 3/A^{8} + 3/A^{12} - 2/A^{16} + A^{-20}$} \\ \textbf{Mock:} {\scriptsize $2w^{3} + 3w^{2} - 2w - 4 + 3/w^{2} - 1/w^{4}$} \\ \textbf{Affine:} {\scriptsize $2t - 4 + 2/t$} \\ \textbf{Yamada:} {\scriptsize $-A^{24} - A^{23} + A^{22} + A^{21} - 2A^{20} + 3A^{18} - A^{17} - 4A^{16} + 2A^{15} - A^{14} - 5A^{13} + A^{12} + A^{11} - 2A^{10} + A^{9} + 2A^{8} + 2A^{7} - 3A^{6} + 3A^{4} - 5A^{3} + 3A - 1$}
\end{minipage}

\noindent{\color{gray!40}\rule{\textwidth}{0.4pt}}
\vspace{0.9\baselineskip}
\noindent \begin{minipage}[t]{0.25\textwidth}
\vspace{0pt}
\centering
\includegraphics[page=72,width=\linewidth]{knotoids.pdf}
\end{minipage}
\hfill
\begin{minipage}[t]{0.73\textwidth}
\vspace{0pt}
\raggedright
\textbf{Name:} {\large{$\mathbf{K6_{35}}$}} (chiral, non-rotatable$^{*}$) \\ \textbf{PD:} {\scriptsize\texttt{[0],[0,1,2,3],[1,4,5,2],[3,6,7,8],[8,7,9,4],[5,10,11,6],[12,11,10,9],[12]}} \\ \textbf{EM:} {\small\texttt{(B0, A0C0C3D0, B1E3F0B2, B3F3E1E0, D3D2G3C1, C2G2G1D1, H0F2F1E2, G0)}} \\ \textbf{Kauffman bracket:} {\scriptsize $A^{20} - 2A^{16} + A^{14} + 3A^{12} + A^{10} - 3A^{8} - 2A^{6} + A^{4} + A^{2}$} \\ \textbf{Arrow:} {\scriptsize $A^{2}L_1 - 2L_1/A^{2} + A^{-4} + 3L_1/A^{6} + L_2/A^{8} - 3L_1/A^{10} - 2L_2/A^{12} + L_1/A^{14} + L_2/A^{16}$} \\ \textbf{Mock:} {\scriptsize $w^{5} + w^{4} - 2w^{3} - w^{2} + 3w + 1 - 3/w + w^{-3}$} \\ \textbf{Affine:} {\scriptsize $t - 2 + 1/t$} \\ \textbf{Yamada:} {\scriptsize $-A^{21} - 3A^{20} + 2A^{19} + A^{18} - 4A^{17} + 3A^{16} + 2A^{15} - 3A^{14} + A^{13} - A^{12} - A^{11} - 4A^{10} - A^{9} + 2A^{8} - 3A^{7} + A^{6} + 4A^{5} - A^{4} - A^{3} + 3A^{2} - 2$}
\end{minipage}

\noindent{\color{gray!40}\rule{\textwidth}{0.4pt}}
\vspace{0.9\baselineskip}
\noindent \begin{minipage}[t]{0.25\textwidth}
\vspace{0pt}
\centering
\includegraphics[page=73,width=\linewidth]{knotoids.pdf}
\end{minipage}
\hfill
\begin{minipage}[t]{0.73\textwidth}
\vspace{0pt}
\raggedright
\textbf{Name:} {\large{$\mathbf{K6_{36}}$}} (chiral, non-rotatable$^{*}$) \\ \textbf{PD:} {\scriptsize\texttt{[0],[0,1,2,3],[1,4,5,2],[3,6,7,8],[8,7,9,4],[10,11,6,5],[9,12,11,10],[12]}} \\ \textbf{EM:} {\small\texttt{(B0, A0C0C3D0, B1E3F3B2, B3F2E1E0, D3D2G0C1, G3G2D1C2, E2H0F1F0, G1)}} \\ \textbf{Kauffman bracket:} {\scriptsize $-A^{18} + 2A^{14} - 3A^{10} - A^{8} + 3A^{6} + 3A^{4} - A^{2} - 1$} \\ \textbf{Arrow:} {\scriptsize $-A^{18}L_1 + 2A^{14}L_1 - 3A^{10}L_1 - A^{8}L_2 + 3A^{6}L_1 + 2A^{4}L_2 + A^{4} - A^{2}L_1 - L_2$} \\ \textbf{Mock:} {\scriptsize $-w^{5} - w^{4} + 2w^{3} + 2w^{2} - 3w - 1 + 3/w + w^{-2} - 1/w^{3}$} \\ \textbf{Affine:} {\scriptsize $-t + 2 - 1/t$} \\ \textbf{Yamada:} {\scriptsize $-2A^{23} + A^{22} + 4A^{21} - 2A^{20} - A^{19} + 5A^{18} - 2A^{17} - 3A^{16} + A^{15} - 2A^{14} - A^{13} - 2A^{12} + 3A^{11} - 4A^{9} + 3A^{8} - 5A^{6} + A^{5} + A^{4} - 2A^{3} - A^{2} + A + 1$}
\end{minipage}

\noindent{\color{gray!40}\rule{\textwidth}{0.4pt}}
\vspace{0.9\baselineskip}
\noindent \begin{minipage}[t]{0.25\textwidth}
\vspace{0pt}
\centering
\includegraphics[page=74,width=\linewidth]{knotoids.pdf}
\end{minipage}
\hfill
\begin{minipage}[t]{0.73\textwidth}
\vspace{0pt}
\raggedright
\textbf{Name:} {\large{$\mathbf{K6_{37}}$}} (chiral, non-rotatable$^{*}$) \\ \textbf{PD:} {\scriptsize\texttt{[0],[0,1,2,3],[1,4,5,2],[6,7,8,3],[4,8,7,9],[5,10,11,6],[12,11,10,9],[12]}} \\ \textbf{EM:} {\small\texttt{(B0, A0C0C3D3, B1E0F0B2, F3E2E1B3, C1D2D1G3, C2G2G1D0, H0F2F1E3, G0)}} \\ \textbf{Kauffman bracket:} {\scriptsize $A^{18} + A^{16} + A^{14} - 2A^{10} - A^{8} + A^{6} + A^{4} - 1$} \\ \textbf{Arrow:} {\scriptsize $A^{-12} + L_1/A^{14} + L_2/A^{16} - 2L_2/A^{20} - L_1/A^{22} + L_2/A^{24} + L_1/A^{26} - L_1/A^{30}$} \\ \textbf{Mock:} {\scriptsize $w^{6} + w^{5} - w + 1 + 1/w - 2/w^{2} - 1/w^{3} + w^{-4}$} \\ \textbf{Affine:} {\scriptsize $3t - 6 + 3/t$} \\ \textbf{Yamada:} {\scriptsize $-A^{24} - A^{23} - A^{20} - A^{19} - A^{18} - 3A^{16} - 2A^{15} - 2A^{13} - A^{12} + 2A^{11} + A^{10} + 3A^{9} + 2A^{8} + 2A^{7} - 2A^{5} + A^{4} - A^{3} - 2A^{2} + 1$}
\end{minipage}

\noindent{\color{gray!40}\rule{\textwidth}{0.4pt}}
\vspace{0.9\baselineskip}
\noindent \begin{minipage}[t]{0.25\textwidth}
\vspace{0pt}
\centering
\includegraphics[page=75,width=\linewidth]{knotoids.pdf}
\end{minipage}
\hfill
\begin{minipage}[t]{0.73\textwidth}
\vspace{0pt}
\raggedright
\textbf{Name:} {\large{$\mathbf{K6_{38}}$}} (chiral, non-rotatable$^{*}$) \\ \textbf{PD:} {\scriptsize\texttt{[0],[0,1,2,3],[1,4,5,2],[3,6,7,4],[5,8,9,6],[10,11,12,7],[8,12,11,9],[10]}} \\ \textbf{EM:} {\small\texttt{(B0, A0C0C3D0, B1D3E0B2, B3E3F3C1, C2G0G3D1, H0G2G1D2, E1F2F1E2, F0)}} \\ \textbf{Kauffman bracket:} {\scriptsize $-A^{18} + A^{16} + 3A^{14} + A^{12} - 3A^{10} - 2A^{8} + 2A^{6} + 2A^{4} - A^{2} - 1$} \\ \textbf{Arrow:} {\scriptsize $-L_1/A^{6} + A^{-8} + 3L_1/A^{10} + L_2/A^{12} - 3L_1/A^{14} - 2L_2/A^{16} + 2L_1/A^{18} + 2L_2/A^{20} - L_1/A^{22} - L_2/A^{24}$} \\ \textbf{Mock:} {\scriptsize $-w^{5} + 3w^{3} + 2w^{2} - 3w - 2 + 2/w + w^{-2} - 1/w^{3}$} \\ \textbf{Affine:} {\scriptsize $t^{2} + t - 4 + 1/t + t^{-2}$} \\ \textbf{Yamada:} {\scriptsize $A^{23} - A^{21} + 3A^{20} + 2A^{19} - 4A^{18} + 3A^{17} + 4A^{16} - 3A^{15} + 2A^{14} + 3A^{13} - 2A^{12} - A^{10} + 2A^{9} - 3A^{8} - 2A^{7} + 5A^{6} - 3A^{5} - 2A^{4} + 4A^{3} - A^{2} - 2A + 1$}
\end{minipage}

\noindent{\color{gray!40}\rule{\textwidth}{0.4pt}}
\vspace{0.9\baselineskip}
\noindent \begin{minipage}[t]{0.25\textwidth}
\vspace{0pt}
\centering
\includegraphics[page=76,width=\linewidth]{knotoids.pdf}
\end{minipage}
\hfill
\begin{minipage}[t]{0.73\textwidth}
\vspace{0pt}
\raggedright
\textbf{Name:} {\large{$\mathbf{K6_{39}}$}} (chiral, non-rotatable$^{*}$) \\ \textbf{PD:} {\scriptsize\texttt{[0],[0,1,2,3],[1,4,5,2],[3,6,7,4],[5,8,9,10],[6,11,12,7],[8,12,11,10],[9]}} \\ \textbf{EM:} {\small\texttt{(B0, A0C0C3D0, B1D3E0B2, B3F0F3C1, C2G0H0G3, D1G2G1D2, E1F2F1E3, E2)}} \\ \textbf{Kauffman bracket:} {\scriptsize $-A^{18} + A^{16} + 2A^{14} - 3A^{10} - A^{8} + 3A^{6} + 2A^{4} - A^{2} - 1$} \\ \textbf{Arrow:} {\scriptsize $-A^{6}L_1 + A^{4} + 2A^{2}L_1 + L_2 - 1 - 3L_1/A^{2} - 2L_2/A^{4} + A^{-4} + 3L_1/A^{6} + 2L_2/A^{8} - L_1/A^{10} - L_2/A^{12}$} \\ \textbf{Mock:} {\scriptsize $-w^{4} - 3w^{3} + 5w + 3 - 2/w - 1/w^{2}$} \\ \textbf{Affine:} {\scriptsize $t^{2} - t - 1/t + t^{-2}$} \\ \textbf{Yamada:} {\scriptsize $-A^{24} + A^{23} + 2A^{22} - 3A^{21} + 3A^{19} - 4A^{18} - A^{17} + 3A^{16} - 3A^{15} - A^{14} + 2A^{13} - A^{12} - 3A^{10} + 2A^{9} - A^{8} - 5A^{7} + 4A^{6} - 3A^{4} + 3A^{3} - A + 1$}
\end{minipage}

\noindent{\color{gray!40}\rule{\textwidth}{0.4pt}}
\vspace{0.9\baselineskip}
\noindent \begin{minipage}[t]{0.25\textwidth}
\vspace{0pt}
\centering
\includegraphics[page=77,width=\linewidth]{knotoids.pdf}
\end{minipage}
\hfill
\begin{minipage}[t]{0.73\textwidth}
\vspace{0pt}
\raggedright
\textbf{Name:} {\large{$\mathbf{K6_{40}}$}} (chiral, non-rotatable$^{*}$) \\ \textbf{PD:} {\scriptsize\texttt{[0],[0,1,2,3],[1,4,5,2],[3,6,7,4],[8,9,6,5],[7,10,11,12],[12,11,9,8],[10]}} \\ \textbf{EM:} {\small\texttt{(B0, A0C0C3D0, B1D3E3B2, B3E2F0C1, G3G2D1C2, D2H0G1G0, F3F2E1E0, F1)}} \\ \textbf{Kauffman bracket:} {\scriptsize $-A^{18} + 2A^{14} + A^{12} - 3A^{10} - 2A^{8} + 2A^{6} + 3A^{4} - 1$} \\ \textbf{Arrow:} {\scriptsize $-A^{18}L_1 + 2A^{14}L_1 + A^{12}L_2 - 3A^{10}L_1 - 2A^{8}L_2 + 2A^{6}L_1 + 2A^{4}L_2 + A^{4} - L_2$} \\ \textbf{Mock:} {\scriptsize $-w^{5} - w^{4} + 2w^{3} + 2w^{2} - 3w - 2 + 2/w + 2/w^{2}$} \\ \textbf{Affine:} {\scriptsize $t^{2} - 2t + 2 - 2/t + t^{-2}$} \\ \textbf{Yamada:} {\scriptsize $A^{24} - A^{23} - 2A^{22} + 3A^{21} - A^{20} - 4A^{19} + 5A^{18} + A^{17} - 2A^{16} + 2A^{15} - 3A^{12} + A^{11} + A^{10} - 5A^{9} + 2A^{7} - 5A^{6} - A^{5} + 2A^{4} - A^{3} - A^{2} + A + 1$}
\end{minipage}

\noindent{\color{gray!40}\rule{\textwidth}{0.4pt}}
\vspace{0.9\baselineskip}
\noindent \begin{minipage}[t]{0.25\textwidth}
\vspace{0pt}
\centering
\includegraphics[page=78,width=\linewidth]{knotoids.pdf}
\end{minipage}
\hfill
\begin{minipage}[t]{0.73\textwidth}
\vspace{0pt}
\raggedright
\textbf{Name:} {\large{$\mathbf{K6_{41}}$}} (chiral, non-rotatable$^{*}$) \\ \textbf{PD:} {\scriptsize\texttt{[0],[0,1,2,3],[1,4,5,2],[3,6,7,4],[8,9,10,5],[6,11,12,7],[8,12,11,10],[9]}} \\ \textbf{EM:} {\small\texttt{(B0, A0C0C3D0, B1D3E3B2, B3F0F3C1, G0H0G3C2, D1G2G1D2, E0F2F1E2, E1)}} \\ \textbf{Kauffman bracket:} {\scriptsize $-A^{20} + A^{16} + 2A^{14} - A^{12} - 3A^{10} + 3A^{6} + 2A^{4} - A^{2} - 1$} \\ \textbf{Arrow:} {\scriptsize $-A^{14}L_1 + A^{10}L_1 + A^{8}L_2 + A^{8} - A^{6}L_1 - 2A^{4}L_2 - A^{4} + 2L_2 + 1 + 2L_1/A^{2} - L_2/A^{4} - L_1/A^{6}$} \\ \textbf{Mock:} {\scriptsize $-2w^{3} - 3w^{2} + w + 5 + 2/w - 1/w^{2} - 1/w^{3}$} \\ \textbf{Affine:} {\scriptsize $t^{2} - 2t + 2 - 2/t + t^{-2}$} \\ \textbf{Yamada:} {\scriptsize $A^{23} - 3A^{21} + 2A^{20} + 2A^{19} - 5A^{18} + A^{17} + 4A^{16} - 3A^{15} + A^{14} + 3A^{13} - 2A^{12} - A^{11} - 3A^{10} + A^{9} - 4A^{8} - 4A^{7} + 5A^{6} - 2A^{5} - 3A^{4} + 4A^{3} - A + 1$}
\end{minipage}

\noindent{\color{gray!40}\rule{\textwidth}{0.4pt}}
\vspace{0.9\baselineskip}
\noindent \begin{minipage}[t]{0.25\textwidth}
\vspace{0pt}
\centering
\includegraphics[page=79,width=\linewidth]{knotoids.pdf}
\end{minipage}
\hfill
\begin{minipage}[t]{0.73\textwidth}
\vspace{0pt}
\raggedright
\textbf{Name:} {\large{$\mathbf{K6_{42}}$}} (chiral, non-rotatable$^{*}$) \\ \textbf{PD:} {\scriptsize\texttt{[0],[0,1,2,3],[1,4,5,2],[6,7,4,3],[5,8,9,6],[10,11,12,7],[8,12,11,9],[10]}} \\ \textbf{EM:} {\small\texttt{(B0, A0C0C3D3, B1D2E0B2, E3F3C1B3, C2G0G3D0, H0G2G1D1, E1F2F1E2, F0)}} \\ \textbf{Kauffman bracket:} {\scriptsize $A^{18} + A^{16} + A^{14} - A^{12} - 2A^{10} + 2A^{6} + A^{4} - A^{2} - 1$} \\ \textbf{Arrow:} {\scriptsize $A^{-12} + L_1/A^{14} + L_2/A^{16} - L_1/A^{18} - 2L_2/A^{20} + 2L_2/A^{24} + L_1/A^{26} - L_2/A^{28} - L_1/A^{30}$} \\ \textbf{Mock:} {\scriptsize $w^{6} + w^{5} - w^{3} - w^{2} + 2 + 1/w - 2/w^{2} - 1/w^{3} + w^{-4}$} \\ \textbf{Affine:} {\scriptsize $t^{2} + 2t - 6 + 2/t + t^{-2}$} \\ \textbf{Yamada:} {\scriptsize $A^{25} + A^{24} + 2A^{21} + A^{20} + A^{18} + 3A^{17} - A^{15} + 2A^{14} - A^{13} - 2A^{12} + A^{11} - A^{10} - A^{8} + 2A^{7} - 3A^{5} + A^{4} - A^{2} + 1$}
\end{minipage}

\noindent{\color{gray!40}\rule{\textwidth}{0.4pt}}
\vspace{0.9\baselineskip}
\noindent \begin{minipage}[t]{0.25\textwidth}
\vspace{0pt}
\centering
\includegraphics[page=80,width=\linewidth]{knotoids.pdf}
\end{minipage}
\hfill
\begin{minipage}[t]{0.73\textwidth}
\vspace{0pt}
\raggedright
\textbf{Name:} {\large{$\mathbf{K6_{43}}$}} (chiral, non-rotatable$^{*}$) \\ \textbf{PD:} {\scriptsize\texttt{[0],[0,1,2,3],[1,4,5,2],[6,7,4,3],[5,8,9,10],[11,12,7,6],[12,11,10,8],[9]}} \\ \textbf{EM:} {\small\texttt{(B0, A0C0C3D3, B1D2E0B2, F3F2C1B3, C2G3H0G2, G1G0D1D0, F1F0E3E1, E2)}} \\ \textbf{Kauffman bracket:} {\scriptsize $A^{16} + 2A^{14} - 3A^{10} - A^{8} + 2A^{6} + 2A^{4} - A^{2} - 1$} \\ \textbf{Arrow:} {\scriptsize $A^{-8} + 2L_1/A^{10} + L_2/A^{12} - 1/A^{12} - 3L_1/A^{14} - 2L_2/A^{16} + A^{-16} + 2L_1/A^{18} + 2L_2/A^{20} - L_1/A^{22} - L_2/A^{24}$} \\ \textbf{Mock:} {\scriptsize $2w^{4} + 3w^{3} - 4w - 3 + 1/w + 2/w^{2}$} \\ \textbf{Affine:} {\scriptsize $t^{2} + 2t - 6 + 2/t + t^{-2}$} \\ \textbf{Yamada:} {\scriptsize $-A^{24} - 3A^{21} - A^{20} + 2A^{19} - 4A^{18} - A^{17} + 3A^{16} - 2A^{15} + 2A^{13} - A^{12} - 2A^{10} + 2A^{9} - 3A^{7} + 4A^{6} - 3A^{4} + 2A^{3} - A + 1$}
\end{minipage}

\noindent{\color{gray!40}\rule{\textwidth}{0.4pt}}
\vspace{0.9\baselineskip}
\noindent \begin{minipage}[t]{0.25\textwidth}
\vspace{0pt}
\centering
\includegraphics[page=81,width=\linewidth]{knotoids.pdf}
\end{minipage}
\hfill
\begin{minipage}[t]{0.73\textwidth}
\vspace{0pt}
\raggedright
\textbf{Name:} {\large{$\mathbf{K6_{44}}$}} (chiral, non-rotatable$^{*}$) \\ \textbf{PD:} {\scriptsize\texttt{[0],[0,1,2,3],[1,4,5,2],[6,7,4,3],[8,9,10,5],[11,12,7,6],[12,11,10,8],[9]}} \\ \textbf{EM:} {\small\texttt{(B0, A0C0C3D3, B1D2E3B2, F3F2C1B3, G3H0G2C2, G1G0D1D0, F1F0E2E0, E1)}} \\ \textbf{Kauffman bracket:} {\scriptsize $A^{22} + 2A^{20} - 3A^{16} - 2A^{14} + 3A^{12} + 2A^{10} - A^{8} - 2A^{6} + A^{2}$} \\ \textbf{Arrow:} {\scriptsize $L_1/A^{2} + L_2/A^{4} + A^{-4} - 2L_2/A^{8} - 1/A^{8} - 2L_1/A^{10} + 2L_2/A^{12} + A^{-12} + 2L_1/A^{14} - L_2/A^{16} - 2L_1/A^{18} + L_1/A^{22}$} \\ \textbf{Mock:} {\scriptsize $w^{3} + 3w^{2} + w - 4 - 4/w + 2/w^{2} + 2/w^{3}$} \\ \textbf{Affine:} {\scriptsize $t^{2} + t - 4 + 1/t + t^{-2}$} \\ \textbf{Yamada:} {\scriptsize $A^{25} + 3A^{22} - 3A^{20} + 5A^{19} - 6A^{17} + 3A^{16} + A^{15} - 3A^{14} + 2A^{13} + 2A^{12} + 3A^{11} - 2A^{10} + A^{9} + 4A^{8} - 6A^{7} + 4A^{5} - 4A^{4} - A^{3} + 3A^{2} - 1$}
\end{minipage}

\noindent{\color{gray!40}\rule{\textwidth}{0.4pt}}
\vspace{0.9\baselineskip}
\noindent \begin{minipage}[t]{0.25\textwidth}
\vspace{0pt}
\centering
\includegraphics[page=82,width=\linewidth]{knotoids.pdf}
\end{minipage}
\hfill
\begin{minipage}[t]{0.73\textwidth}
\vspace{0pt}
\raggedright
\textbf{Name:} {\large{$\mathbf{K6_{45}}$}} (chiral, non-rotatable$^{*}$) \\ \textbf{PD:} {\scriptsize\texttt{[0],[0,1,2,3],[1,4,5,2],[3,6,7,4],[5,7,8,9],[6,10,11,8],[12,11,10,9],[12]}} \\ \textbf{EM:} {\small\texttt{(B0, A0C0C3D0, B1D3E0B2, B3F0E1C1, C2D2F3G3, D1G2G1E2, H0F2F1E3, G0)}} \\ \textbf{Kauffman bracket:} {\scriptsize $-A^{18} + A^{16} + 3A^{14} + 2A^{12} - 3A^{10} - 3A^{8} + A^{6} + 2A^{4} - 1$} \\ \textbf{Arrow:} {\scriptsize $-L_1/A^{6} + A^{-8} + 3L_1/A^{10} + 2L_2/A^{12} - 3L_1/A^{14} - 3L_2/A^{16} + L_1/A^{18} + 2L_2/A^{20} - L_2/A^{24}$} \\ \textbf{Mock:} {\scriptsize $-w^{5} + 3w^{3} + 2w^{2} - 3w - 3 + 1/w + 2/w^{2}$} \\ \textbf{Affine:} {\scriptsize $2t^{2} - 4 + 2/t^{2}$} \\ \textbf{Yamada:} {\scriptsize $-2A^{22} - A^{21} + 2A^{20} - 2A^{19} - 4A^{18} + 4A^{17} - 5A^{15} + 3A^{14} + A^{13} - 3A^{12} + A^{11} + A^{9} - 4A^{8} + A^{7} + 4A^{6} - 5A^{5} + A^{4} + 4A^{3} - 3A^{2} - A + 2$}
\end{minipage}

\noindent{\color{gray!40}\rule{\textwidth}{0.4pt}}
\vspace{0.9\baselineskip}
\noindent \begin{minipage}[t]{0.25\textwidth}
\vspace{0pt}
\centering
\includegraphics[page=83,width=\linewidth]{knotoids.pdf}
\end{minipage}
\hfill
\begin{minipage}[t]{0.73\textwidth}
\vspace{0pt}
\raggedright
\textbf{Name:} {\large{$\mathbf{K6_{46}}$}} (chiral, non-rotatable$^{*}$) \\ \textbf{PD:} {\scriptsize\texttt{[0],[0,1,2,3],[1,4,5,2],[3,6,7,4],[7,8,9,5],[6,10,11,8],[12,11,10,9],[12]}} \\ \textbf{EM:} {\small\texttt{(B0, A0C0C3D0, B1D3E3B2, B3F0E0C1, D2F3G3C2, D1G2G1E1, H0F2F1E2, G0)}} \\ \textbf{Kauffman bracket:} {\scriptsize $-A^{20} + A^{16} + 2A^{14} - A^{10} + A^{6} - A^{2}$} \\ \textbf{Arrow:} {\scriptsize $-A^{2}L_1 + L_1/A^{2} + L_2/A^{4} + A^{-4} - L_2/A^{8} + L_2/A^{12} - L_2/A^{16}$} \\ \textbf{Mock:} {\scriptsize $-w^{5} - w^{4} + w^{3} + 2w^{2} - 1 + w^{-2}$} \\ \textbf{Affine:} {\scriptsize $2t^{2} - t - 2 - 1/t + 2/t^{2}$} \\ \textbf{Yamada:} {\scriptsize $-A^{20} - A^{19} - A^{16} + A^{14} - A^{13} - A^{12} - 2A^{10} - A^{8} + A^{7} - A^{6} - A^{5} + A^{4} - A^{3} + A + 1$}
\end{minipage}

\noindent{\color{gray!40}\rule{\textwidth}{0.4pt}}
\vspace{0.9\baselineskip}
\noindent \begin{minipage}[t]{0.25\textwidth}
\vspace{0pt}
\centering
\includegraphics[page=84,width=\linewidth]{knotoids.pdf}
\end{minipage}
\hfill
\begin{minipage}[t]{0.73\textwidth}
\vspace{0pt}
\raggedright
\textbf{Name:} {\large{$\mathbf{K6_{47}}$}} (chiral, non-rotatable$^{*}$) \\ \textbf{PD:} {\scriptsize\texttt{[0],[0,1,2,3],[1,4,5,2],[6,7,4,3],[5,7,8,9],[6,10,11,8],[12,11,10,9],[12]}} \\ \textbf{EM:} {\small\texttt{(B0, A0C0C3D3, B1D2E0B2, F0E1C1B3, C2D1F3G3, D0G2G1E2, H0F2F1E3, G0)}} \\ \textbf{Kauffman bracket:} {\scriptsize $A^{18} + A^{16} + A^{14} - A^{12} - A^{10} + A^{6} - A^{2}$} \\ \textbf{Arrow:} {\scriptsize $A^{-12} + L_1/A^{14} + L_2/A^{16} - L_1/A^{18} - L_2/A^{20} + L_2/A^{24} - L_2/A^{28}$} \\ \textbf{Mock:} {\scriptsize $w^{6} + w^{5} - w^{3} - w^{2} + 1 - 1/w^{2} + w^{-4}$} \\ \textbf{Affine:} {\scriptsize $2t^{2} + t - 6 + 1/t + 2/t^{2}$} \\ \textbf{Yamada:} {\scriptsize $-A^{22} - A^{21} - A^{20} - A^{19} - 2A^{18} - A^{17} + A^{13} - A^{11} - A^{8} + A^{7} + A^{5} - A^{4} + A^{2} - A + 1$}
\end{minipage}

\noindent{\color{gray!40}\rule{\textwidth}{0.4pt}}
\vspace{0.9\baselineskip}
\noindent \begin{minipage}[t]{0.25\textwidth}
\vspace{0pt}
\centering
\includegraphics[page=85,width=\linewidth]{knotoids.pdf}
\end{minipage}
\hfill
\begin{minipage}[t]{0.73\textwidth}
\vspace{0pt}
\raggedright
\textbf{Name:} {\large{$\mathbf{K6_{48}}$}} (achiral, non-rotatable$^{*}$) \\ \textbf{PD:} {\scriptsize\texttt{[0],[0,1,2,3],[1,4,5,2],[6,7,4,3],[7,8,9,5],[10,11,8,6],[9,12,11,10],[12]}} \\ \textbf{EM:} {\small\texttt{(B0, A0C0C3D3, B1D2E3B2, F3E0C1B3, D1F2G0C2, G3G2E1D0, E2H0F1F0, G1)}} \\ \textbf{Kauffman bracket:} {\scriptsize $A^{22} + A^{20} - 2A^{18} - 3A^{16} + A^{14} + 5A^{12} + A^{10} - 3A^{8} - 2A^{6} + A^{4} + A^{2}$} \\ \textbf{Arrow:} {\scriptsize $A^{10}L_1 + A^{8}L_2 - 2A^{6}L_1 - 3A^{4}L_2 + A^{2}L_1 + 4L_2 + 1 + L_1/A^{2} - 3L_2/A^{4} - 2L_1/A^{6} + L_2/A^{8} + L_1/A^{10}$} \\ \textbf{Mock:} {\scriptsize $w^{5} + w^{4} - 2w^{3} - 3w^{2} + w + 5 + 1/w - 3/w^{2} - 2/w^{3} + w^{-4} + w^{-5}$} \\ \textbf{Affine:} {\scriptsize $0$} \\ \textbf{Yamada:} {\scriptsize $-2A^{24} + 4A^{22} - 2A^{21} - 4A^{20} + 7A^{19} + A^{18} - 8A^{17} + 5A^{16} + 2A^{15} - 6A^{14} - 6A^{10} + 2A^{9} + 5A^{8} - 8A^{7} + A^{6} + 7A^{5} - 4A^{4} - 2A^{3} + 4A^{2} - 2$}
\end{minipage}

\noindent{\color{gray!40}\rule{\textwidth}{0.4pt}}
\vspace{0.9\baselineskip}
\noindent \begin{minipage}[t]{0.25\textwidth}
\vspace{0pt}
\centering
\includegraphics[page=86,width=\linewidth]{knotoids.pdf}
\end{minipage}
\hfill
\begin{minipage}[t]{0.73\textwidth}
\vspace{0pt}
\raggedright
\textbf{Name:} {\large{$\mathbf{K6_{49}}$}} (chiral, non-rotatable$^{*}$) \\ \textbf{PD:} {\scriptsize\texttt{[0],[0,1,2,3],[1,4,5,2],[3,6,7,4],[5,8,9,6],[10,11,8,7],[11,10,12,9],[12]}} \\ \textbf{EM:} {\small\texttt{(B0, A0C0C3D0, B1D3E0B2, B3E3F3C1, C2F2G3D1, G1G0E1D2, F1F0H0E2, G2)}} \\ \textbf{Kauffman bracket:} {\scriptsize $A^{20} - A^{18} - 3A^{16} + 4A^{12} + 2A^{10} - 2A^{8} - 2A^{6} + A^{4} + A^{2}$} \\ \textbf{Arrow:} {\scriptsize $A^{20} - A^{18}L_1 - 3A^{16} + 4A^{12} + 2A^{10}L_1 - 2A^{8} - 2A^{6}L_1 + A^{4} + A^{2}L_1$} \\ \textbf{Mock:} {\scriptsize $-w^{4} - w^{3} + 3w^{2} + w - 4 - 1/w + 3/w^{2} + w^{-3}$} \\ \textbf{Affine:} {\scriptsize $-2t + 4 - 2/t$} \\ \textbf{Yamada:} {\scriptsize $-A^{23} + A^{22} + A^{21} - 4A^{20} + 3A^{19} + 4A^{18} - 6A^{17} + 3A^{16} + 3A^{15} - 3A^{14} + A^{13} + A^{11} - 4A^{10} - 2A^{9} + 3A^{8} - 4A^{7} - 3A^{6} + 5A^{5} - 2A^{4} - 3A^{3} + 3A^{2} - 2$}
\end{minipage}

\noindent{\color{gray!40}\rule{\textwidth}{0.4pt}}
\vspace{0.9\baselineskip}
\noindent \begin{minipage}[t]{0.25\textwidth}
\vspace{0pt}
\centering
\includegraphics[page=87,width=\linewidth]{knotoids.pdf}
\end{minipage}
\hfill
\begin{minipage}[t]{0.73\textwidth}
\vspace{0pt}
\raggedright
\textbf{Name:} {\large{$\mathbf{K6_{50}}$}} (chiral, non-rotatable$^{*}$) \\ \textbf{PD:} {\scriptsize\texttt{[0],[0,1,2,3],[1,4,5,2],[3,6,7,4],[5,8,9,10],[6,11,8,7],[11,10,12,9],[12]}} \\ \textbf{EM:} {\small\texttt{(B0, A0C0C3D0, B1D3E0B2, B3F0F3C1, C2F2G3G1, D1G0E1D2, F1E3H0E2, G2)}} \\ \textbf{Kauffman bracket:} {\scriptsize $-A^{18} + A^{16} + 3A^{14} + A^{12} - 4A^{10} - 3A^{8} + 2A^{6} + 3A^{4} - 1$} \\ \textbf{Arrow:} {\scriptsize $-A^{6}L_1 + A^{4}L_2 + 3A^{2}L_1 - L_2 + 2 - 4L_1/A^{2} - 3/A^{4} + 2L_1/A^{6} + 3/A^{8} - 1/A^{12}$} \\ \textbf{Mock:} {\scriptsize $-w^{3} - w^{2} + 3w + 3 - 4/w - 2/w^{2} + 2/w^{3} + w^{-4}$} \\ \textbf{Affine:} {\scriptsize $t^{2} - t - 1/t + t^{-2}$} \\ \textbf{Yamada:} {\scriptsize $2A^{23} - 4A^{21} + 3A^{20} + 3A^{19} - 7A^{18} + 2A^{17} + 5A^{16} - 5A^{15} + A^{14} + 4A^{13} - A^{12} + 5A^{9} - 2A^{8} - 3A^{7} + 8A^{6} - 3A^{5} - 5A^{4} + 5A^{3} - A^{2} - 2A + 1$}
\end{minipage}

\noindent{\color{gray!40}\rule{\textwidth}{0.4pt}}
\vspace{0.9\baselineskip}
\noindent \begin{minipage}[t]{0.25\textwidth}
\vspace{0pt}
\centering
\includegraphics[page=88,width=\linewidth]{knotoids.pdf}
\end{minipage}
\hfill
\begin{minipage}[t]{0.73\textwidth}
\vspace{0pt}
\raggedright
\textbf{Name:} {\large{$\mathbf{K6_{51}}$}} (chiral, non-rotatable$^{*}$) \\ \textbf{PD:} {\scriptsize\texttt{[0],[0,1,2,3],[1,4,5,2],[3,6,7,4],[8,9,6,5],[7,10,11,8],[9,11,10,12],[12]}} \\ \textbf{EM:} {\small\texttt{(B0, A0C0C3D0, B1D3E3B2, B3E2F0C1, F3G0D1C2, D2G2G1E0, E1F2F1H0, G3)}} \\ \textbf{Kauffman bracket:} {\scriptsize $-A^{19} + 2A^{15} + A^{13} - 3A^{11} - 3A^{9} + A^{7} + 3A^{5} - A$} \\ \textbf{Arrow:} {\scriptsize $A^{4} - 2 - L_1/A^{2} + 3/A^{4} + 3L_1/A^{6} - 1/A^{8} - 3L_1/A^{10} + L_1/A^{14}$} \\ \textbf{Mock:} {\scriptsize $-w^{4} - w^{3} + 3w^{2} + 3w - 2 - 3/w + w^{-2} + w^{-3}$} \\ \textbf{Affine:} {\scriptsize $0$} \\ \textbf{Yamada:} {\scriptsize $-A^{24} + 2A^{22} - A^{21} - 2A^{20} + 5A^{19} + 2A^{18} - 3A^{17} + 4A^{16} + 2A^{15} - 3A^{14} + A^{13} + A^{11} - 3A^{10} + A^{9} + 3A^{8} - 5A^{7} + A^{6} + 4A^{5} - 2A^{4} - A^{3} + 2A^{2} - 1$}
\end{minipage}

\noindent{\color{gray!40}\rule{\textwidth}{0.4pt}}
\vspace{0.9\baselineskip}
\noindent \begin{minipage}[t]{0.25\textwidth}
\vspace{0pt}
\centering
\includegraphics[page=89,width=\linewidth]{knotoids.pdf}
\end{minipage}
\hfill
\begin{minipage}[t]{0.73\textwidth}
\vspace{0pt}
\raggedright
\textbf{Name:} {\large{$\mathbf{K6_{52}}$}} (chiral, non-rotatable$^{*}$) \\ \textbf{PD:} {\scriptsize\texttt{[0],[0,1,2,3],[1,4,5,2],[3,5,6,7],[4,7,8,9],[9,10,11,6],[8,11,10,12],[12]}} \\ \textbf{EM:} {\small\texttt{(B0, A0C0C3D0, B1E0D1B2, B3C2F3E1, C1D3G0F0, E3G2G1D2, E2F2F1H0, G3)}} \\ \textbf{Kauffman bracket:} {\scriptsize $-A^{18} + 3A^{14} + 2A^{12} - 3A^{10} - 3A^{8} + A^{6} + 3A^{4} - 1$} \\ \textbf{Arrow:} {\scriptsize $-A^{6}L_1 + 3A^{2}L_1 + 2 - 3L_1/A^{2} - 3/A^{4} + L_1/A^{6} + 3/A^{8} - 1/A^{12}$} \\ \textbf{Mock:} {\scriptsize $-w^{3} - w^{2} + 3w + 4 - 3/w - 3/w^{2} + w^{-3} + w^{-4}$} \\ \textbf{Affine:} {\scriptsize $0$} \\ \textbf{Yamada:} {\scriptsize $A^{24} + A^{23} - 2A^{22} - A^{21} + 3A^{20} - 2A^{19} - 5A^{18} + 4A^{17} + A^{16} - 5A^{15} + 3A^{14} + 2A^{13} - 3A^{12} + A^{9} - 5A^{8} + 4A^{6} - 6A^{5} + 5A^{3} - 2A^{2} - A + 1$}
\end{minipage}

\noindent{\color{gray!40}\rule{\textwidth}{0.4pt}}
\vspace{0.9\baselineskip}
\noindent \begin{minipage}[t]{0.25\textwidth}
\vspace{0pt}
\centering
\includegraphics[page=90,width=\linewidth]{knotoids.pdf}
\end{minipage}
\hfill
\begin{minipage}[t]{0.73\textwidth}
\vspace{0pt}
\raggedright
\textbf{Name:} {\large{$\mathbf{K6_{53}}$}} (chiral, non-rotatable$^{*}$) \\ \textbf{PD:} {\scriptsize\texttt{[0],[0,1,2,3],[1,4,5,2],[3,5,6,7],[4,7,8,9],[10,11,12,6],[12,10,9,8],[11]}} \\ \textbf{EM:} {\small\texttt{(B0, A0C0C3D0, B1E0D1B2, B3C2F3E1, C1D3G3G2, G1H0G0D2, F2F0E3E2, F1)}} \\ \textbf{Kauffman bracket:} {\scriptsize $A^{22} + A^{20} - A^{18} - 3A^{16} + 4A^{12} + 2A^{10} - 2A^{8} - 2A^{6} + A^{2}$} \\ \textbf{Arrow:} {\scriptsize $A^{4} + A^{2}L_1 - 1 - 3L_1/A^{2} + 4L_1/A^{6} + 2/A^{8} - 2L_1/A^{10} - 2/A^{12} + A^{-16}$} \\ \textbf{Mock:} {\scriptsize $-w^{4} - 2w^{3} + w^{2} + 5w + 3 - 3/w - 2/w^{2}$} \\ \textbf{Affine:} {\scriptsize $t - 2 + 1/t$} \\ \textbf{Yamada:} {\scriptsize $A^{23} - A^{22} - 3A^{21} + 3A^{20} + A^{19} - 6A^{18} + 4A^{17} + 3A^{16} - 6A^{15} + A^{14} + 2A^{13} - 3A^{12} - A^{11} - A^{10} + 2A^{9} - 4A^{8} - 2A^{7} + 7A^{6} - 4A^{5} - 2A^{4} + 6A^{3} - 2A^{2} - 2A + 1$}
\end{minipage}

\noindent{\color{gray!40}\rule{\textwidth}{0.4pt}}
\vspace{0.9\baselineskip}
\noindent \begin{minipage}[t]{0.25\textwidth}
\vspace{0pt}
\centering
\includegraphics[page=91,width=\linewidth]{knotoids.pdf}
\end{minipage}
\hfill
\begin{minipage}[t]{0.73\textwidth}
\vspace{0pt}
\raggedright
\textbf{Name:} {\large{$\mathbf{K6_{54}}$}} (chiral, non-rotatable$^{*}$) \\ \textbf{PD:} {\scriptsize\texttt{[0],[0,1,2,3],[1,4,5,2],[3,5,6,7],[7,8,9,4],[6,9,10,11],[11,10,12,8],[12]}} \\ \textbf{EM:} {\small\texttt{(B0, A0C0C3D0, B1E3D1B2, B3C2F0E0, D3G3F1C1, D2E2G1G0, F3F2H0E1, G2)}} \\ \textbf{Kauffman bracket:} {\scriptsize $-A^{19} - A^{17} + 2A^{13} + A^{11} - 2A^{9} - 2A^{7} + A^{3} + A$} \\ \textbf{Arrow:} {\scriptsize $A^{-8} + L_1/A^{10} - 2L_1/A^{14} - 1/A^{16} + 2L_1/A^{18} + 2/A^{20} - 1/A^{24} - L_1/A^{26}$} \\ \textbf{Mock:} {\scriptsize $w^{4} + w^{3} + w^{2} - w - 2 + 1/w + w^{-2} - 1/w^{3}$} \\ \textbf{Affine:} {\scriptsize $2t - 4 + 2/t$} \\ \textbf{Yamada:} {\scriptsize $A^{25} + A^{24} + A^{23} + 2A^{21} + 2A^{20} - 2A^{19} + A^{18} + 2A^{17} - 2A^{16} + A^{15} + A^{14} + 2A^{13} - A^{12} + A^{10} - 3A^{9} - A^{8} + A^{7} - A^{6} - 2A^{5} + A^{4} + A^{3} - A^{2} + 1$}
\end{minipage}

\noindent{\color{gray!40}\rule{\textwidth}{0.4pt}}
\vspace{0.9\baselineskip}
\noindent \begin{minipage}[t]{0.25\textwidth}
\vspace{0pt}
\centering
\includegraphics[page=92,width=\linewidth]{knotoids.pdf}
\end{minipage}
\hfill
\begin{minipage}[t]{0.73\textwidth}
\vspace{0pt}
\raggedright
\textbf{Name:} {\large{$\mathbf{K6_{55}}$}} (chiral, non-rotatable$^{*}$) \\ \textbf{PD:} {\scriptsize\texttt{[0],[0,1,2,3],[1,4,5,2],[3,5,6,7],[7,8,9,4],[10,11,12,6],[8,12,10,9],[11]}} \\ \textbf{EM:} {\small\texttt{(B0, A0C0C3D0, B1E3D1B2, B3C2F3E0, D3G0G3C1, G2H0G1D2, E1F2F0E2, F1)}} \\ \textbf{Kauffman bracket:} {\scriptsize $A^{22} - 2A^{18} - A^{16} + 3A^{14} + 3A^{12} - 2A^{10} - 4A^{8} + 2A^{4} + A^{2}$} \\ \textbf{Arrow:} {\scriptsize $A^{16} - 2A^{12} - A^{10}L_1 + 3A^{8} + 3A^{6}L_1 - 2A^{4} - 4A^{2}L_1 + 2L_1/A^{2} + A^{-4}$} \\ \textbf{Mock:} {\scriptsize $w^{4} + 2w^{3} - 3w^{2} - 5w + 3 + 3/w$} \\ \textbf{Affine:} {\scriptsize $-t + 2 - 1/t$} \\ \textbf{Yamada:} {\scriptsize $A^{25} - 3A^{23} + A^{22} + 4A^{21} - 6A^{20} + 7A^{18} - 6A^{17} - A^{16} + 4A^{15} - 2A^{14} - A^{13} - A^{12} + 3A^{11} - 3A^{10} - 5A^{9} + 6A^{8} - A^{7} - 6A^{6} + 5A^{5} + A^{4} - 4A^{3} + A^{2} + A - 1$}
\end{minipage}

\noindent{\color{gray!40}\rule{\textwidth}{0.4pt}}
\vspace{0.9\baselineskip}
\noindent \begin{minipage}[t]{0.25\textwidth}
\vspace{0pt}
\centering
\includegraphics[page=93,width=\linewidth]{knotoids.pdf}
\end{minipage}
\hfill
\begin{minipage}[t]{0.73\textwidth}
\vspace{0pt}
\raggedright
\textbf{Name:} {\large{$\mathbf{K6_{56}}$}} (chiral, non-rotatable$^{*}$) \\ \textbf{PD:} {\scriptsize\texttt{[0],[0,1,2,3],[1,4,5,2],[3,6,7,8],[8,9,10,4],[5,11,12,6],[11,10,9,7],[12]}} \\ \textbf{EM:} {\small\texttt{(B0, A0C0C3D0, B1E3F0B2, B3F3G3E0, D3G2G1C1, C2G0H0D1, F1E2E1D2, F2)}} \\ \textbf{Kauffman bracket:} {\scriptsize $A^{20} - A^{16} + A^{14} + 3A^{12} - 3A^{8} - 2A^{6} + A^{4} + A^{2}$} \\ \textbf{Arrow:} {\scriptsize $A^{-4} - 1/A^{8} + L_1/A^{10} + 3/A^{12} - 3/A^{16} - 2L_1/A^{18} + A^{-20} + L_1/A^{22}$} \\ \textbf{Mock:} {\scriptsize $w^{3} + 3w^{2} - 4 - 2/w + 2/w^{2} + w^{-3}$} \\ \textbf{Affine:} {\scriptsize $t - 2 + 1/t$} \\ \textbf{Yamada:} {\scriptsize $A^{22} + 3A^{21} - A^{20} + A^{19} + 4A^{18} - 2A^{17} + 3A^{15} - 2A^{14} - A^{12} + A^{11} - A^{10} - 2A^{9} + 3A^{8} - A^{7} - 2A^{6} + 2A^{5} + A^{4} - 2A^{3} + A^{2} + A - 1$}
\end{minipage}

\noindent{\color{gray!40}\rule{\textwidth}{0.4pt}}
\vspace{0.9\baselineskip}
\noindent \begin{minipage}[t]{0.25\textwidth}
\vspace{0pt}
\centering
\includegraphics[page=94,width=\linewidth]{knotoids.pdf}
\end{minipage}
\hfill
\begin{minipage}[t]{0.73\textwidth}
\vspace{0pt}
\raggedright
\textbf{Name:} {\large{$\mathbf{K6_{57}}$}} (chiral, non-rotatable$^{*}$) \\ \textbf{PD:} {\scriptsize\texttt{[0],[0,1,2,3],[1,4,5,2],[3,6,7,8],[8,7,9,4],[5,9,10,11],[6,11,12,10],[12]}} \\ \textbf{EM:} {\small\texttt{(B0, A0C0C3D0, B1E3F0B2, B3G0E1E0, D3D2F1C1, C2E2G3G1, D1F3H0F2, G2)}} \\ \textbf{Kauffman bracket:} {\scriptsize $A^{20} - 2A^{16} + 3A^{12} + 2A^{10} - 2A^{8} - 2A^{6} + A^{2}$} \\ \textbf{Arrow:} {\scriptsize $A^{2}L_1 - 2L_1/A^{2} - L_2/A^{4} + A^{-4} + 3L_1/A^{6} + 2L_2/A^{8} - 2L_1/A^{10} - 2L_2/A^{12} + L_2/A^{16}$} \\ \textbf{Mock:} {\scriptsize $w^{5} + w^{4} - 2w^{3} - w^{2} + 3w + 2 - 2/w - 1/w^{2}$} \\ \textbf{Affine:} {\scriptsize $-t^{2} + 2t - 2 + 2/t - 1/t^{2}$} \\ \textbf{Yamada:} {\scriptsize $A^{22} + A^{21} - 3A^{20} + 2A^{18} - 5A^{17} - A^{16} + 2A^{15} - 3A^{14} - A^{12} + A^{11} - 2A^{10} - A^{9} + 4A^{8} - A^{7} - A^{6} + 3A^{5} - 3A^{3} + A^{2} + A - 1$}
\end{minipage}

\noindent{\color{gray!40}\rule{\textwidth}{0.4pt}}
\vspace{0.9\baselineskip}
\noindent \begin{minipage}[t]{0.25\textwidth}
\vspace{0pt}
\centering
\includegraphics[page=95,width=\linewidth]{knotoids.pdf}
\end{minipage}
\hfill
\begin{minipage}[t]{0.73\textwidth}
\vspace{0pt}
\raggedright
\textbf{Name:} {\large{$\mathbf{K6_{58}}$}} (chiral, non-rotatable$^{*}$) \\ \textbf{PD:} {\scriptsize\texttt{[0],[0,1,2,3],[1,4,5,2],[3,6,7,8],[8,9,10,4],[11,12,6,5],[11,10,9,7],[12]}} \\ \textbf{EM:} {\small\texttt{(B0, A0C0C3D0, B1E3F3B2, B3F2G3E0, D3G2G1C1, G0H0D1C2, F0E2E1D2, F1)}} \\ \textbf{Kauffman bracket:} {\scriptsize $-A^{19} + A^{15} - 3A^{11} - 2A^{9} + 2A^{7} + 3A^{5} - A$} \\ \textbf{Arrow:} {\scriptsize $A^{-8} - 1/A^{12} + 3/A^{16} + 2L_1/A^{18} - 2/A^{20} - 3L_1/A^{22} + L_1/A^{26}$} \\ \textbf{Mock:} {\scriptsize $2w^{4} + 2w^{3} - 2w^{2} - 3w + 1 + 1/w - 1/w^{2} + w^{-4}$} \\ \textbf{Affine:} {\scriptsize $t - 2 + 1/t$} \\ \textbf{Yamada:} {\scriptsize $2A^{23} + 2A^{22} + A^{20} + 5A^{19} - 2A^{17} + 3A^{16} - 3A^{14} + A^{13} - A^{12} - A^{11} - 3A^{10} + A^{9} + A^{8} - 4A^{7} + 3A^{6} + 2A^{5} - 2A^{4} + 2A^{2} - 1$}
\end{minipage}

\noindent{\color{gray!40}\rule{\textwidth}{0.4pt}}
\vspace{0.9\baselineskip}
\noindent \begin{minipage}[t]{0.25\textwidth}
\vspace{0pt}
\centering
\includegraphics[page=96,width=\linewidth]{knotoids.pdf}
\end{minipage}
\hfill
\begin{minipage}[t]{0.73\textwidth}
\vspace{0pt}
\raggedright
\textbf{Name:} {\large{$\mathbf{K6_{59}}$}} (chiral, non-rotatable$^{*}$) \\ \textbf{PD:} {\scriptsize\texttt{[0],[0,1,2,3],[1,4,5,2],[3,6,7,8],[8,7,9,4],[9,10,11,5],[6,11,12,10],[12]}} \\ \textbf{EM:} {\small\texttt{(B0, A0C0C3D0, B1E3F3B2, B3G0E1E0, D3D2F0C1, E2G3G1C2, D1F2H0F1, G2)}} \\ \textbf{Kauffman bracket:} {\scriptsize $-A^{19} + 2A^{15} + A^{13} - 2A^{11} - 3A^{9} + A^{7} + 2A^{5} - A$} \\ \textbf{Arrow:} {\scriptsize $A^{10}L_1 - 2A^{6}L_1 - A^{4}L_2 + 2A^{2}L_1 + 2L_2 + 1 - L_1/A^{2} - 2L_2/A^{4} + L_2/A^{8}$} \\ \textbf{Mock:} {\scriptsize $w^{5} + w^{4} - 2w^{3} - 2w^{2} + 2w + 3 - 1/w - 1/w^{2}$} \\ \textbf{Affine:} {\scriptsize $-t^{2} + t + 1/t - 1/t^{2}$} \\ \textbf{Yamada:} {\scriptsize $A^{23} - 2A^{21} + A^{20} + 2A^{19} - 3A^{18} - A^{17} + 3A^{16} - 2A^{15} - 2A^{14} + A^{13} - 2A^{12} - A^{11} - 2A^{10} + 2A^{9} - A^{8} - 3A^{7} + 4A^{6} - 2A^{4} + 2A^{3} + A^{2} - A - 1$}
\end{minipage}

\noindent{\color{gray!40}\rule{\textwidth}{0.4pt}}
\vspace{0.9\baselineskip}
\noindent \begin{minipage}[t]{0.25\textwidth}
\vspace{0pt}
\centering
\includegraphics[page=97,width=\linewidth]{knotoids.pdf}
\end{minipage}
\hfill
\begin{minipage}[t]{0.73\textwidth}
\vspace{0pt}
\raggedright
\textbf{Name:} {\large{$\mathbf{K6_{60}}$}} (chiral, non-rotatable$^{*}$) \\ \textbf{PD:} {\scriptsize\texttt{[0],[0,1,2,3],[1,4,5,2],[6,7,4,3],[5,8,9,10],[11,8,7,6],[11,10,12,9],[12]}} \\ \textbf{EM:} {\small\texttt{(B0, A0C0C3D3, B1D2E0B2, F3F2C1B3, C2F1G3G1, G0E1D1D0, F0E3H0E2, G2)}} \\ \textbf{Kauffman bracket:} {\scriptsize $A^{17} + 2A^{15} - 3A^{11} - 3A^{9} + A^{7} + 2A^{5} - A$} \\ \textbf{Arrow:} {\scriptsize $-A^{20}L_2 - 2A^{18}L_1 + A^{16}L_2 - A^{16} + 3A^{14}L_1 + 3A^{12} - A^{10}L_1 - 2A^{8} + A^{4}$} \\ \textbf{Mock:} {\scriptsize $w^{2} - 2w - 2 + 3/w + 2/w^{2} - 1/w^{3}$} \\ \textbf{Affine:} {\scriptsize $-t^{2} - t + 4 - 1/t - 1/t^{2}$} \\ \textbf{Yamada:} {\scriptsize $-A^{23} + 2A^{21} - 2A^{19} + 3A^{18} + A^{17} - 3A^{16} + A^{15} + 2A^{14} - 2A^{13} + 2A^{11} - 2A^{10} - A^{9} - 2A^{8} + 2A^{7} - 3A^{6} - 2A^{5} + 4A^{4} - 3A^{3} - 2A^{2} + A - 1$}
\end{minipage}

\noindent{\color{gray!40}\rule{\textwidth}{0.4pt}}
\vspace{0.9\baselineskip}
\noindent \begin{minipage}[t]{0.25\textwidth}
\vspace{0pt}
\centering
\includegraphics[page=98,width=\linewidth]{knotoids.pdf}
\end{minipage}
\hfill
\begin{minipage}[t]{0.73\textwidth}
\vspace{0pt}
\raggedright
\textbf{Name:} {\large{$\mathbf{K6_{61}}$}} (chiral, non-rotatable$^{*}$) \\ \textbf{PD:} {\scriptsize\texttt{[0],[0,1,2,3],[1,4,5,2],[6,7,8,3],[4,8,9,10],[11,12,6,5],[7,11,10,9],[12]}} \\ \textbf{EM:} {\small\texttt{(B0, A0C0C3D3, B1E0F3B2, F2G0E1B3, C1D2G3G2, G1H0D0C2, D1F0E3E2, F1)}} \\ \textbf{Kauffman bracket:} {\scriptsize $-2A^{18} - A^{16} + 3A^{14} + 4A^{12} - A^{10} - 3A^{8} + 2A^{4} - 1$} \\ \textbf{Arrow:} {\scriptsize $-2A^{6}L_1 - A^{4} + 3A^{2}L_1 + 4 - L_1/A^{2} - 3/A^{4} + 2/A^{8} - 1/A^{12}$} \\ \textbf{Mock:} {\scriptsize $-2w^{3} - 2w^{2} + 3w + 5 - 1/w - 3/w^{2} + w^{-4}$} \\ \textbf{Affine:} {\scriptsize $-t + 2 - 1/t$} \\ \textbf{Yamada:} {\scriptsize $A^{24} - 2A^{22} + 2A^{20} - 3A^{19} - 3A^{18} + 5A^{17} - 4A^{15} + 5A^{14} + 2A^{13} - 2A^{12} + 3A^{11} + 2A^{10} + 2A^{9} - 3A^{8} + 3A^{7} + 3A^{6} - 6A^{5} + A^{4} + 3A^{3} - 4A^{2} - A + 2$}
\end{minipage}

\noindent{\color{gray!40}\rule{\textwidth}{0.4pt}}
\vspace{0.9\baselineskip}
\noindent \begin{minipage}[t]{0.25\textwidth}
\vspace{0pt}
\centering
\includegraphics[page=99,width=\linewidth]{knotoids.pdf}
\end{minipage}
\hfill
\begin{minipage}[t]{0.73\textwidth}
\vspace{0pt}
\raggedright
\textbf{Name:} {\large{$\mathbf{K6_{62}}$}} (chiral, non-rotatable$^{*}$) \\ \textbf{PD:} {\scriptsize\texttt{[0],[0,1,2,3],[1,3,4,5],[5,6,7,2],[4,7,8,9],[10,11,12,6],[11,10,9,8],[12]}} \\ \textbf{EM:} {\small\texttt{(B0, A0C0D3C1, B1B3E0D0, C3F3E1B2, C2D2G3G2, G1G0H0D1, F1F0E3E2, F2)}} \\ \textbf{Kauffman bracket:} {\scriptsize $-2A^{17} - A^{15} + 2A^{13} + 3A^{11} - A^{9} - 3A^{7} - A^{5} + A^{3} + A$} \\ \textbf{Arrow:} {\scriptsize $L_2/A^{4} + A^{-4} + L_1/A^{6} - L_2/A^{8} - 1/A^{8} - 3L_1/A^{10} + A^{-12} + 3L_1/A^{14} + L_2/A^{16} - L_1/A^{18} - L_2/A^{20}$} \\ \textbf{Mock:} {\scriptsize $-w^{3} + w^{2} + 3w - 1 - 3/w + w^{-2} + w^{-3}$} \\ \textbf{Affine:} {\scriptsize $2t^{2} - 4 + 2/t^{2}$} \\ \textbf{Yamada:} {\scriptsize $-A^{26} + A^{24} - A^{23} - 3A^{22} + 2A^{21} + A^{20} - 5A^{19} + 2A^{18} + 2A^{17} - 3A^{16} + A^{13} - 3A^{12} + A^{11} + 3A^{10} - 4A^{9} - A^{8} + 3A^{7} - 2A^{6} - 2A^{5} + 2A^{4} + A^{3} - A^{2} + 1$}
\end{minipage}

\noindent{\color{gray!40}\rule{\textwidth}{0.4pt}}
\vspace{0.9\baselineskip}
\noindent \begin{minipage}[t]{0.25\textwidth}
\vspace{0pt}
\centering
\includegraphics[page=100,width=\linewidth]{knotoids.pdf}
\end{minipage}
\hfill
\begin{minipage}[t]{0.73\textwidth}
\vspace{0pt}
\raggedright
\textbf{Name:} {\large{$\mathbf{K6_{63}}$}} (chiral, non-rotatable$^{*}$) \\ \textbf{PD:} {\scriptsize\texttt{[0],[0,1,2,3],[1,4,5,2],[3,6,7,4],[5,7,8,9],[6,10,11,8],[9,11,12,10],[12]}} \\ \textbf{EM:} {\small\texttt{(B0, A0C0C3D0, B1D3E0B2, B3F0E1C1, C2D2F3G0, D1G3G1E2, E3F2H0F1, G2)}} \\ \textbf{Kauffman bracket:} {\scriptsize $-A^{18} - 2A^{16} + 2A^{14} + 4A^{12} + A^{10} - 3A^{8} - 2A^{6} + A^{4} + A^{2}$} \\ \textbf{Arrow:} {\scriptsize $-A^{12}L_2 - 2A^{10}L_1 + 2A^{8}L_2 + 4A^{6}L_1 - A^{4}L_2 + 2A^{4} - 3A^{2}L_1 - 2 + L_1/A^{2} + A^{-4}$} \\ \textbf{Mock:} {\scriptsize $-w^{2} - 3w + 1 + 5/w + 2/w^{2} - 2/w^{3} - 1/w^{4}$} \\ \textbf{Affine:} {\scriptsize $-t + 2 - 1/t$} \\ \textbf{Yamada:} {\scriptsize $-2A^{23} + A^{22} + 3A^{21} - 4A^{20} + A^{19} + 5A^{18} - 5A^{17} + 3A^{15} - 3A^{14} + 4A^{11} - A^{9} + 6A^{8} - A^{7} - 4A^{6} + 4A^{5} - 3A^{3} + 2A^{2} + A - 1$}
\end{minipage}

\noindent{\color{gray!40}\rule{\textwidth}{0.4pt}}
\vspace{0.9\baselineskip}
\noindent \begin{minipage}[t]{0.25\textwidth}
\vspace{0pt}
\centering
\includegraphics[page=101,width=\linewidth]{knotoids.pdf}
\end{minipage}
\hfill
\begin{minipage}[t]{0.73\textwidth}
\vspace{0pt}
\raggedright
\textbf{Name:} {\large{$\mathbf{K6_{64}}$}} (chiral, non-rotatable$^{*}$) \\ \textbf{PD:} {\scriptsize\texttt{[0],[0,1,2,3],[1,4,5,2],[3,6,7,4],[7,8,9,5],[6,10,11,8],[11,12,10,9],[12]}} \\ \textbf{EM:} {\small\texttt{(B0, A0C0C3D0, B1D3E3B2, B3F0E0C1, D2F3G3C2, D1G2G0E1, F2H0F1E2, G1)}} \\ \textbf{Kauffman bracket:} {\scriptsize $A^{17} - A^{13} - 2A^{11} + 2A^{7} + A^{5} - A^{3} - A$} \\ \textbf{Arrow:} {\scriptsize $-A^{14}L_1 + A^{10}L_1 + A^{8}L_2 + A^{8} - A^{4}L_2 - A^{4} - A^{2}L_1 + 1 + L_1/A^{2}$} \\ \textbf{Mock:} {\scriptsize $-w^{3} - 2w^{2} + 3 + 1/w$} \\ \textbf{Affine:} {\scriptsize $t^{2} - 2t + 2 - 2/t + t^{-2}$} \\ \textbf{Yamada:} {\scriptsize $-A^{19} + 2A^{17} - A^{16} - A^{15} + 2A^{14} - A^{13} - A^{12} + A^{11} - A^{10} - 2A^{8} - A^{6} - 3A^{5} + A^{4} - 2A^{2} + A + 1$}
\end{minipage}

\noindent{\color{gray!40}\rule{\textwidth}{0.4pt}}
\vspace{0.9\baselineskip}
\noindent \begin{minipage}[t]{0.25\textwidth}
\vspace{0pt}
\centering
\includegraphics[page=102,width=\linewidth]{knotoids.pdf}
\end{minipage}
\hfill
\begin{minipage}[t]{0.73\textwidth}
\vspace{0pt}
\raggedright
\textbf{Name:} {\large{$\mathbf{K6_{65}}$}} (chiral, non-rotatable$^{*}$) \\ \textbf{PD:} {\scriptsize\texttt{[0],[0,1,2,3],[1,4,5,2],[3,6,7,8],[8,7,9,4],[5,10,11,6],[9,11,12,10],[12]}} \\ \textbf{EM:} {\small\texttt{(B0, A0C0C3D0, B1E3F0B2, B3F3E1E0, D3D2G0C1, C2G3G1D1, E2F2H0F1, G2)}} \\ \textbf{Kauffman bracket:} {\scriptsize $A^{13} - 2A^{9} - 2A^{7} + A^{3} + A$} \\ \textbf{Arrow:} {\scriptsize $-L_1/A^{2} + 2L_1/A^{6} + 2/A^{8} - 1/A^{12} - L_1/A^{14}$} \\ \textbf{Mock:} {\scriptsize $-w^{3} + 2w + 2 - 1/w^{2} - 1/w^{3}$} \\ \textbf{Affine:} {\scriptsize $t - 2 + 1/t$} \\ \textbf{Yamada:} {\scriptsize $-A^{19} + A^{18} + 2A^{17} + 2A^{15} + A^{14} + 2A^{13} + A^{11} + A^{10} - A^{9} - A^{6} - A^{5} - A^{2} + 1$}
\end{minipage}

\noindent{\color{gray!40}\rule{\textwidth}{0.4pt}}
\vspace{0.9\baselineskip}
\noindent \begin{minipage}[t]{0.25\textwidth}
\vspace{0pt}
\centering
\includegraphics[page=103,width=\linewidth]{knotoids.pdf}
\end{minipage}
\hfill
\begin{minipage}[t]{0.73\textwidth}
\vspace{0pt}
\raggedright
\textbf{Name:} {\large{$\mathbf{K6_{66}}$}} (chiral, non-rotatable$^{*}$) \\ \textbf{PD:} {\scriptsize\texttt{[0],[0,1,2,3],[1,4,5,2],[6,7,8,3],[4,8,7,9],[5,10,11,6],[9,11,12,10],[12]}} \\ \textbf{EM:} {\small\texttt{(B0, A0C0C3D3, B1E0F0B2, F3E2E1B3, C1D2D1G0, C2G3G1D0, E3F2H0F1, G2)}} \\ \textbf{Kauffman bracket:} {\scriptsize $A^{17} + A^{15} - A^{13} - 3A^{11} - A^{9} + 2A^{7} + 2A^{5} - A^{3} - A$} \\ \textbf{Arrow:} {\scriptsize $-A^{14}L_1 - A^{12} + A^{10}L_1 + 3A^{8} + A^{6}L_1 - 2A^{4} - 2A^{2}L_1 + 1 + L_1/A^{2}$} \\ \textbf{Mock:} {\scriptsize $-w^{3} - 2w^{2} + 4 + 2/w - 1/w^{2} - 1/w^{3}$} \\ \textbf{Affine:} {\scriptsize $-t + 2 - 1/t$} \\ \textbf{Yamada:} {\scriptsize $-A^{25} + A^{23} - A^{22} - A^{21} + 2A^{20} + 2A^{19} - 2A^{18} + 3A^{16} - 2A^{15} - A^{14} + 2A^{13} - A^{12} + 3A^{9} - A^{8} - A^{7} + 4A^{6} - A^{4} + 2A^{3} + A^{2} - A - 1$}
\end{minipage}

\noindent{\color{gray!40}\rule{\textwidth}{0.4pt}}
\vspace{0.9\baselineskip}
\noindent \begin{minipage}[t]{0.25\textwidth}
\vspace{0pt}
\centering
\includegraphics[page=104,width=\linewidth]{knotoids.pdf}
\end{minipage}
\hfill
\begin{minipage}[t]{0.73\textwidth}
\vspace{0pt}
\raggedright
\textbf{Name:} {\large{$\mathbf{K6_{67}}$}} (chiral, non-rotatable$^{*}$) \\ \textbf{PD:} {\scriptsize\texttt{[0],[0,1,2,3],[1,3,4,5],[2,6,7,8],[9,6,5,4],[10,11,8,7],[9,11,12,10],[12]}} \\ \textbf{EM:} {\small\texttt{(B0, A0C0D0C1, B1B3E3E2, B2E1F3F2, G0D1C3C2, G3G1D3D2, E0F1H0F0, G2)}} \\ \textbf{Kauffman bracket:} {\scriptsize $A^{16} + 2A^{14} - A^{12} - 4A^{10} - A^{8} + 3A^{6} + 3A^{4} - A^{2} - 1$} \\ \textbf{Arrow:} {\scriptsize $A^{4}L_2 + 2A^{2}L_1 - 2L_2 + 1 - 4L_1/A^{2} + L_2/A^{4} - 2/A^{4} + 3L_1/A^{6} + 3/A^{8} - L_1/A^{10} - 1/A^{12}$} \\ \textbf{Mock:} {\scriptsize $w^{2} + 3w + 1 - 5/w - 2/w^{2} + 2/w^{3} + w^{-4}$} \\ \textbf{Affine:} {\scriptsize $t - 2 + 1/t$} \\ \textbf{Yamada:} {\scriptsize $-A^{24} + A^{23} + A^{22} - 4A^{21} + A^{20} + 4A^{19} - 5A^{18} + A^{17} + 5A^{16} - 4A^{15} - A^{14} + A^{13} - 3A^{12} - 2A^{11} - 3A^{10} + 3A^{9} - A^{8} - 4A^{7} + 7A^{6} - 4A^{4} + 4A^{3} - A^{2} - 2A + 1$}
\end{minipage}

\noindent{\color{gray!40}\rule{\textwidth}{0.4pt}}
\vspace{0.9\baselineskip}
\noindent \begin{minipage}[t]{0.25\textwidth}
\vspace{0pt}
\centering
\includegraphics[page=105,width=\linewidth]{knotoids.pdf}
\end{minipage}
\hfill
\begin{minipage}[t]{0.73\textwidth}
\vspace{0pt}
\raggedright
\textbf{Name:} {\large{$\mathbf{K6_{68}}$}} (chiral, non-rotatable$^{*}$) \\ \textbf{PD:} {\scriptsize\texttt{[0],[0,1,2,3],[1,3,4,5],[6,7,8,2],[8,9,5,4],[10,11,7,6],[9,11,12,10],[12]}} \\ \textbf{EM:} {\small\texttt{(B0, A0C0D3C1, B1B3E3E2, F3F2E0B2, D2G0C3C2, G3G1D1D0, E1F1H0F0, G2)}} \\ \textbf{Kauffman bracket:} {\scriptsize $-A^{19} - 2A^{17} - A^{15} + 2A^{13} + 3A^{11} - 2A^{7} - A^{5} + A$} \\ \textbf{Arrow:} {\scriptsize $A^{-8} + 2L_1/A^{10} + L_2/A^{12} - 2L_1/A^{14} - 2L_2/A^{16} - 1/A^{16} + L_2/A^{20} + A^{-20} + L_1/A^{22} - L_1/A^{26}$} \\ \textbf{Mock:} {\scriptsize $w^{4} + 2w^{3} + 2w^{2} - w - 3 - 1/w + w^{-2}$} \\ \textbf{Affine:} {\scriptsize $3t - 6 + 3/t$} \\ \textbf{Yamada:} {\scriptsize $-A^{22} - A^{21} + A^{20} - 4A^{18} + A^{17} - 5A^{15} + A^{13} - A^{12} + 2A^{11} + A^{10} + 2A^{9} - 2A^{8} - A^{7} + 3A^{6} - 3A^{5} - A^{4} + 3A^{3} - A^{2} - A + 1$}
\end{minipage}

\noindent{\color{gray!40}\rule{\textwidth}{0.4pt}}
\vspace{0.9\baselineskip}
\noindent \begin{minipage}[t]{0.25\textwidth}
\vspace{0pt}
\centering
\includegraphics[page=106,width=\linewidth]{knotoids.pdf}
\end{minipage}
\hfill
\begin{minipage}[t]{0.73\textwidth}
\vspace{0pt}
\raggedright
\textbf{Name:} {\large{$\mathbf{K6_{69}}$}} (chiral, non-rotatable$^{*}$) \\ \textbf{PD:} {\scriptsize\texttt{[0],[0,1,2,3],[1,4,5,2],[3,6,7,4],[5,8,9,10],[6,10,11,7],[8,11,12,9],[12]}} \\ \textbf{EM:} {\small\texttt{(B0, A0C0C3D0, B1D3E0B2, B3F0F3C1, C2G0G3F1, D1E3G1D2, E1F2H0E2, G2)}} \\ \textbf{Kauffman bracket:} {\scriptsize $-A^{18} + A^{16} + 3A^{14} - 4A^{10} - 3A^{8} + 3A^{6} + 3A^{4} - 1$} \\ \textbf{Arrow:} {\scriptsize $-A^{12} + A^{10}L_1 + 3A^{8} - 4A^{4} - 3A^{2}L_1 + 3 + 3L_1/A^{2} - L_1/A^{6}$} \\ \textbf{Mock:} {\scriptsize $w^{4} + w^{3} - 3w^{2} - 2w + 5 + 2/w - 2/w^{2} - 1/w^{3}$} \\ \textbf{Affine:} {\scriptsize $t - 2 + 1/t$} \\ \textbf{Yamada:} {\scriptsize $A^{23} - A^{22} - 3A^{21} + 4A^{20} + 2A^{19} - 7A^{18} + 4A^{17} + 5A^{16} - 6A^{15} + 2A^{14} + 3A^{13} - 3A^{12} - A^{11} + 4A^{9} - 3A^{8} - A^{7} + 9A^{6} - 3A^{5} - 3A^{4} + 7A^{3} - 2A^{2} - 3A + 1$}
\end{minipage}

\noindent{\color{gray!40}\rule{\textwidth}{0.4pt}}
\vspace{0.9\baselineskip}
\noindent \begin{minipage}[t]{0.25\textwidth}
\vspace{0pt}
\centering
\includegraphics[page=107,width=\linewidth]{knotoids.pdf}
\end{minipage}
\hfill
\begin{minipage}[t]{0.73\textwidth}
\vspace{0pt}
\raggedright
\textbf{Name:} {\large{$\mathbf{K6_{70}}$}} (chiral, non-rotatable$^{*}$) \\ \textbf{PD:} {\scriptsize\texttt{[0],[0,1,2,3],[1,4,5,2],[3,6,7,4],[8,9,10,5],[6,10,11,7],[11,12,9,8],[12]}} \\ \textbf{EM:} {\small\texttt{(B0, A0C0C3D0, B1D3E3B2, B3F0F3C1, G3G2F1C2, D1E2G0D2, F2H0E1E0, G1)}} \\ \textbf{Kauffman bracket:} {\scriptsize $A^{22} - 2A^{18} - A^{16} + 3A^{14} + 4A^{12} - A^{10} - 4A^{8} - A^{6} + A^{4} + A^{2}$} \\ \textbf{Arrow:} {\scriptsize $A^{4} - 2 - L_1/A^{2} + 3/A^{4} + 4L_1/A^{6} - 1/A^{8} - 4L_1/A^{10} - 1/A^{12} + L_1/A^{14} + A^{-16}$} \\ \textbf{Mock:} {\scriptsize $-w^{4} - w^{3} + 3w^{2} + 4w - 1 - 4/w + w^{-3}$} \\ \textbf{Affine:} {\scriptsize $t - 2 + 1/t$} \\ \textbf{Yamada:} {\scriptsize $-A^{24} + 4A^{22} - 5A^{20} + 6A^{19} + 3A^{18} - 7A^{17} + 3A^{16} + 3A^{15} - 4A^{14} + A^{13} + 2A^{12} + 3A^{11} - 3A^{10} + A^{9} + 6A^{8} - 7A^{7} - 2A^{6} + 6A^{5} - 3A^{4} - 3A^{3} + 3A^{2} + A - 1$}
\end{minipage}

\noindent{\color{gray!40}\rule{\textwidth}{0.4pt}}
\vspace{0.9\baselineskip}
\noindent \begin{minipage}[t]{0.25\textwidth}
\vspace{0pt}
\centering
\includegraphics[page=108,width=\linewidth]{knotoids.pdf}
\end{minipage}
\hfill
\begin{minipage}[t]{0.73\textwidth}
\vspace{0pt}
\raggedright
\textbf{Name:} {\large{$\mathbf{K6_{71}}$}} (chiral, non-rotatable$^{*}$) \\ \textbf{PD:} {\scriptsize\texttt{[0],[0,1,2,3],[1,4,5,2],[6,7,4,3],[5,8,9,10],[10,11,7,6],[8,11,12,9],[12]}} \\ \textbf{EM:} {\small\texttt{(B0, A0C0C3D3, B1D2E0B2, F3F2C1B3, C2G0G3F0, E3G1D1D0, E1F1H0E2, G2)}} \\ \textbf{Kauffman bracket:} {\scriptsize $-A^{20} - 2A^{18} + 3A^{14} + 3A^{12} - A^{10} - 3A^{8} + A^{4} + A^{2}$} \\ \textbf{Arrow:} {\scriptsize $-A^{26}L_1 - 2A^{24} + 3A^{20} + 3A^{18}L_1 - A^{16} - 3A^{14}L_1 + A^{10}L_1 + A^{8}$} \\ \textbf{Mock:} {\scriptsize $-w^{3} + 2w^{2} + 2w - 3 - 2/w + w^{-2} + w^{-3} + w^{-4}$} \\ \textbf{Affine:} {\scriptsize $-t + 2 - 1/t$} \\ \textbf{Yamada:} {\scriptsize $-A^{25} + 2A^{23} - A^{22} - 3A^{21} + 3A^{20} + A^{19} - 4A^{18} + A^{17} + 5A^{16} - A^{15} + A^{14} + 4A^{13} - A^{12} - A^{11} - 2A^{10} + 2A^{9} - 5A^{8} - 4A^{7} + 4A^{6} - 3A^{5} - 3A^{4} + 2A^{3} - A - 1$}
\end{minipage}

\noindent{\color{gray!40}\rule{\textwidth}{0.4pt}}
\vspace{0.9\baselineskip}
\noindent \begin{minipage}[t]{0.25\textwidth}
\vspace{0pt}
\centering
\includegraphics[page=109,width=\linewidth]{knotoids.pdf}
\end{minipage}
\hfill
\begin{minipage}[t]{0.73\textwidth}
\vspace{0pt}
\raggedright
\textbf{Name:} {\large{$\mathbf{K6_{72}}$}} (chiral, non-rotatable$^{*}$) \\ \textbf{PD:} {\scriptsize\texttt{[0],[0,1,2,3],[1,4,5,2],[3,6,7,4],[5,8,9,6],[7,9,10,11],[11,10,12,8],[12]}} \\ \textbf{EM:} {\small\texttt{(B0, A0C0C3D0, B1D3E0B2, B3E3F0C1, C2G3F1D1, D2E2G1G0, F3F2H0E1, G2)}} \\ \textbf{Kauffman bracket:} {\scriptsize $-A^{17} + A^{13} + A^{11} - 2A^{9} - 2A^{7} + A^{3} + A$} \\ \textbf{Arrow:} {\scriptsize $A^{-4} - 1/A^{8} - L_1/A^{10} + 2/A^{12} + 2L_1/A^{14} + L_2/A^{16} - 1/A^{16} - L_1/A^{18} - L_2/A^{20}$} \\ \textbf{Mock:} {\scriptsize $-w^{3} + w^{2} + 2w - 1 - 1/w + w^{-2}$} \\ \textbf{Affine:} {\scriptsize $t^{2} - 2 + t^{-2}$} \\ \textbf{Yamada:} {\scriptsize $A^{23} + A^{21} + 2A^{20} - A^{19} + A^{18} + 2A^{17} - A^{16} + A^{15} + A^{13} - A^{12} + 2A^{10} + A^{7} - A^{6} - 2A^{5} - A^{2} + 1$}
\end{minipage}

\noindent{\color{gray!40}\rule{\textwidth}{0.4pt}}
\vspace{0.9\baselineskip}
\noindent \begin{minipage}[t]{0.25\textwidth}
\vspace{0pt}
\centering
\includegraphics[page=110,width=\linewidth]{knotoids.pdf}
\end{minipage}
\hfill
\begin{minipage}[t]{0.73\textwidth}
\vspace{0pt}
\raggedright
\textbf{Name:} {\large{$\mathbf{K6_{73}}$}} (chiral, non-rotatable$^{*}$) \\ \textbf{PD:} {\scriptsize\texttt{[0],[0,1,2,3],[1,4,5,2],[3,6,7,4],[5,7,8,9],[6,9,10,11],[11,10,12,8],[12]}} \\ \textbf{EM:} {\small\texttt{(B0, A0C0C3D0, B1D3E0B2, B3F0E1C1, C2D2G3F1, D1E3G1G0, F3F2H0E2, G2)}} \\ \textbf{Kauffman bracket:} {\scriptsize $A^{20} - A^{18} - 3A^{16} + 4A^{12} + 3A^{10} - 2A^{8} - 2A^{6} + A^{2}$} \\ \textbf{Arrow:} {\scriptsize $A^{14}L_1 - A^{12}L_2 - 3A^{10}L_1 + A^{8}L_2 - A^{8} + 4A^{6}L_1 + 3A^{4} - 2A^{2}L_1 - 2 + A^{-4}$} \\ \textbf{Mock:} {\scriptsize $w^{3} + w^{2} - 3w - 1 + 4/w + 2/w^{2} - 2/w^{3} - 1/w^{4}$} \\ \textbf{Affine:} {\scriptsize $-t^{2} + t + 1/t - 1/t^{2}$} \\ \textbf{Yamada:} {\scriptsize $-A^{23} + 2A^{22} + 2A^{21} - 5A^{20} + 2A^{19} + 5A^{18} - 6A^{17} + A^{16} + 3A^{15} - 4A^{14} - A^{13} - A^{12} + 2A^{11} - 3A^{10} - 2A^{9} + 5A^{8} - 3A^{7} - 4A^{6} + 5A^{5} - A^{4} - 4A^{3} + 2A^{2} + A - 1$}
\end{minipage}

\noindent{\color{gray!40}\rule{\textwidth}{0.4pt}}
\vspace{0.9\baselineskip}
\noindent \begin{minipage}[t]{0.25\textwidth}
\vspace{0pt}
\centering
\includegraphics[page=111,width=\linewidth]{knotoids.pdf}
\end{minipage}
\hfill
\begin{minipage}[t]{0.73\textwidth}
\vspace{0pt}
\raggedright
\textbf{Name:} {\large{$\mathbf{K6_{74}}$}} (chiral, non-rotatable$^{*}$) \\ \textbf{PD:} {\scriptsize\texttt{[0],[0,1,2,3],[1,4,5,2],[3,6,7,4],[8,9,6,5],[9,10,11,7],[8,11,10,12],[12]}} \\ \textbf{EM:} {\small\texttt{(B0, A0C0C3D0, B1D3E3B2, B3E2F3C1, G0F0D1C2, E1G2G1D2, E0F2F1H0, G3)}} \\ \textbf{Kauffman bracket:} {\scriptsize $-A^{20} - A^{18} + A^{16} + 4A^{14} + A^{12} - 4A^{10} - 3A^{8} + 2A^{6} + 3A^{4} - 1$} \\ \textbf{Arrow:} {\scriptsize $-A^{14}L_1 - A^{12} + A^{10}L_1 + A^{8}L_2 + 3A^{8} + A^{6}L_1 - A^{4}L_2 - 3A^{4} - 3A^{2}L_1 + 2 + 3L_1/A^{2} - L_1/A^{6}$} \\ \textbf{Mock:} {\scriptsize $w^{4} - 4w^{2} - 2w + 5 + 3/w - 1/w^{2} - 1/w^{3}$} \\ \textbf{Affine:} {\scriptsize $t^{2} - t - 1/t + t^{-2}$} \\ \textbf{Yamada:} {\scriptsize $A^{24} - 3A^{22} + A^{21} + 4A^{20} - 5A^{19} - 3A^{18} + 8A^{17} - 2A^{16} - 6A^{15} + 7A^{14} + A^{13} - 4A^{12} + 2A^{11} + A^{10} - A^{9} - 7A^{8} + 3A^{7} + 2A^{6} - 10A^{5} + 3A^{4} + 5A^{3} - 6A^{2} + 3$}
\end{minipage}

\noindent{\color{gray!40}\rule{\textwidth}{0.4pt}}
\vspace{0.9\baselineskip}
\noindent \begin{minipage}[t]{0.25\textwidth}
\vspace{0pt}
\centering
\includegraphics[page=112,width=\linewidth]{knotoids.pdf}
\end{minipage}
\hfill
\begin{minipage}[t]{0.73\textwidth}
\vspace{0pt}
\raggedright
\textbf{Name:} {\large{$\mathbf{K6_{75}}$}} (chiral, non-rotatable$^{*}$) \\ \textbf{PD:} {\scriptsize\texttt{[0],[0,1,2,3],[1,4,5,2],[3,6,7,4],[7,8,9,5],[6,9,10,11],[11,10,12,8],[12]}} \\ \textbf{EM:} {\small\texttt{(B0, A0C0C3D0, B1D3E3B2, B3F0E0C1, D2G3F1C2, D1E2G1G0, F3F2H0E1, G2)}} \\ \textbf{Kauffman bracket:} {\scriptsize $-A^{19} + 2A^{15} + 2A^{13} - A^{11} - 3A^{9} + A^{5} - A$} \\ \textbf{Arrow:} {\scriptsize $A^{22}L_1 - 2A^{18}L_1 - A^{16}L_2 - A^{16} + A^{14}L_1 + A^{12}L_2 + 2A^{12} - A^{8} + A^{4}$} \\ \textbf{Mock:} {\scriptsize $w^{3} + 2w^{2} - 2w - 3 + 1/w + 2/w^{2}$} \\ \textbf{Affine:} {\scriptsize $-t^{2} + 2 - 1/t^{2}$} \\ \textbf{Yamada:} {\scriptsize $-A^{21} + 3A^{19} - 2A^{17} + 3A^{16} - A^{15} - 2A^{14} + A^{13} - A^{12} - 2A^{10} + A^{9} - 3A^{7} + A^{6} + A^{5} - 3A^{4} + A^{2} - A - 1$}
\end{minipage}

\noindent{\color{gray!40}\rule{\textwidth}{0.4pt}}
\vspace{0.9\baselineskip}
\noindent \begin{minipage}[t]{0.25\textwidth}
\vspace{0pt}
\centering
\includegraphics[page=113,width=\linewidth]{knotoids.pdf}
\end{minipage}
\hfill
\begin{minipage}[t]{0.73\textwidth}
\vspace{0pt}
\raggedright
\textbf{Name:} {\large{$\mathbf{K6_{76}}$}} (chiral, non-rotatable$^{*}$) \\ \textbf{PD:} {\scriptsize\texttt{[0],[0,1,2,3],[1,4,5,2],[6,7,4,3],[5,8,9,6],[7,9,10,11],[11,10,12,8],[12]}} \\ \textbf{EM:} {\small\texttt{(B0, A0C0C3D3, B1D2E0B2, E3F0C1B3, C2G3F1D0, D1E2G1G0, F3F2H0E1, G2)}} \\ \textbf{Kauffman bracket:} {\scriptsize $A^{16} + A^{14} - A^{12} - 3A^{10} + A^{8} + 4A^{6} + 2A^{4} - 2A^{2} - 2$} \\ \textbf{Arrow:} {\scriptsize $A^{-8} + L_1/A^{10} - 1/A^{12} - 3L_1/A^{14} + A^{-16} + 4L_1/A^{18} + L_2/A^{20} + A^{-20} - 2L_1/A^{22} - L_2/A^{24} - 1/A^{24}$} \\ \textbf{Mock:} {\scriptsize $2w^{4} + 2w^{3} - w^{2} - 4w - 2 + 3/w + 2/w^{2} - 1/w^{3}$} \\ \textbf{Affine:} {\scriptsize $t^{2} + t - 4 + 1/t + t^{-2}$} \\ \textbf{Yamada:} {\scriptsize $2A^{25} + A^{24} - 2A^{23} + A^{22} + 5A^{21} - 2A^{20} - 2A^{19} + 6A^{18} - A^{17} - 3A^{16} + 3A^{15} - 2A^{12} + 4A^{11} - 5A^{9} + 4A^{8} + A^{7} - 5A^{6} + 2A^{4} - 2A^{3} - A^{2} + A + 1$}
\end{minipage}

\noindent{\color{gray!40}\rule{\textwidth}{0.4pt}}
\vspace{0.9\baselineskip}
\noindent \begin{minipage}[t]{0.25\textwidth}
\vspace{0pt}
\centering
\includegraphics[page=114,width=\linewidth]{knotoids.pdf}
\end{minipage}
\hfill
\begin{minipage}[t]{0.73\textwidth}
\vspace{0pt}
\raggedright
\textbf{Name:} {\large{$\mathbf{K6_{77}}$}} (chiral, non-rotatable$^{*}$) \\ \textbf{PD:} {\scriptsize\texttt{[0],[0,1,2,3],[1,4,5,2],[3,6,7,4],[5,8,9,6],[7,10,11,12],[8,12,10,9],[11]}} \\ \textbf{EM:} {\small\texttt{(B0, A0C0C3D0, B1D3E0B2, B3E3F0C1, C2G0G3D1, D2G2H0G1, E1F3F1E2, F2)}} \\ \textbf{Kauffman bracket:} {\scriptsize $-A^{18} - A^{16} + 2A^{14} + 3A^{12} - A^{10} - 3A^{8} - A^{6} + 2A^{4} + A^{2}$} \\ \textbf{Arrow:} {\scriptsize $-A^{18}L_1 - A^{16} + 2A^{14}L_1 + 3A^{12} - A^{10}L_1 - 3A^{8} - A^{6}L_1 + 2A^{4} + A^{2}L_1$} \\ \textbf{Mock:} {\scriptsize $w^{3} + 2w^{2} - 2w - 4 + 3/w^{2} + w^{-3}$} \\ \textbf{Affine:} {\scriptsize $-t + 2 - 1/t$} \\ \textbf{Yamada:} {\scriptsize $A^{24} - A^{23} - A^{22} + 3A^{21} - 3A^{20} - 2A^{19} + 5A^{18} - 3A^{17} + 4A^{15} - A^{14} - 2A^{12} + A^{11} - A^{10} - 4A^{9} + 2A^{8} - 5A^{6} + 2A^{5} + A^{4} - 3A^{3} + A^{2} + A - 1$}
\end{minipage}

\noindent{\color{gray!40}\rule{\textwidth}{0.4pt}}
\vspace{0.9\baselineskip}
\noindent \begin{minipage}[t]{0.25\textwidth}
\vspace{0pt}
\centering
\includegraphics[page=115,width=\linewidth]{knotoids.pdf}
\end{minipage}
\hfill
\begin{minipage}[t]{0.73\textwidth}
\vspace{0pt}
\raggedright
\textbf{Name:} {\large{$\mathbf{K6_{78}}$}} (chiral, non-rotatable$^{*}$) \\ \textbf{PD:} {\scriptsize\texttt{[0],[0,1,2,3],[1,4,5,2],[3,6,7,4],[8,9,6,5],[7,10,11,12],[12,10,9,8],[11]}} \\ \textbf{EM:} {\small\texttt{(B0, A0C0C3D0, B1D3E3B2, B3E2F0C1, G3G2D1C2, D2G1H0G0, F3F1E1E0, F2)}} \\ \textbf{Kauffman bracket:} {\scriptsize $-A^{19} - A^{17} + 2A^{15} + 3A^{13} - A^{11} - 4A^{9} - A^{7} + 2A^{5} + A^{3} - A$} \\ \textbf{Arrow:} {\scriptsize $A^{10}L_1 + A^{8} - 2A^{6}L_1 - 3A^{4} + A^{2}L_1 + 4 + L_1/A^{2} - 2/A^{4} - L_1/A^{6} + A^{-8}$} \\ \textbf{Mock:} {\scriptsize $-w^{3} - 2w^{2} + 2w + 6 - 3/w^{2} - 1/w^{3}$} \\ \textbf{Affine:} {\scriptsize $t - 2 + 1/t$} \\ \textbf{Yamada:} {\scriptsize $-A^{24} + A^{23} + 2A^{22} - 3A^{21} + A^{20} + 4A^{19} - 5A^{18} + 4A^{16} - 4A^{15} + 3A^{13} + A^{11} - A^{10} + 4A^{9} - 4A^{7} + 6A^{6} - 4A^{4} + 4A^{3} - A^{2} - 2A + 1$}
\end{minipage}

\noindent{\color{gray!40}\rule{\textwidth}{0.4pt}}
\vspace{0.9\baselineskip}
\noindent \begin{minipage}[t]{0.25\textwidth}
\vspace{0pt}
\centering
\includegraphics[page=116,width=\linewidth]{knotoids.pdf}
\end{minipage}
\hfill
\begin{minipage}[t]{0.73\textwidth}
\vspace{0pt}
\raggedright
\textbf{Name:} {\large{$\mathbf{K6_{79}}$}} (chiral, non-rotatable$^{*}$) \\ \textbf{PD:} {\scriptsize\texttt{[0],[0,1,2,3],[1,4,5,2],[3,6,7,4],[8,9,6,5],[10,11,12,7],[12,10,9,8],[11]}} \\ \textbf{EM:} {\small\texttt{(B0, A0C0C3D0, B1D3E3B2, B3E2F3C1, G3G2D1C2, G1H0G0D2, F2F0E1E0, F1)}} \\ \textbf{Kauffman bracket:} {\scriptsize $A^{22} + A^{20} - 2A^{18} - 3A^{16} + A^{14} + 5A^{12} + 2A^{10} - 3A^{8} - 2A^{6} + A^{2}$} \\ \textbf{Arrow:} {\scriptsize $A^{4} + A^{2}L_1 - 2 - 3L_1/A^{2} + A^{-4} + 5L_1/A^{6} + 2/A^{8} - 3L_1/A^{10} - 2/A^{12} + A^{-16}$} \\ \textbf{Mock:} {\scriptsize $-w^{4} - 2w^{3} + 2w^{2} + 6w + 2 - 4/w - 2/w^{2}$} \\ \textbf{Affine:} {\scriptsize $2t - 4 + 2/t$} \\ \textbf{Yamada:} {\scriptsize $A^{25} - A^{24} - 2A^{23} + 4A^{22} - A^{21} - 7A^{20} + 6A^{19} + A^{18} - 9A^{17} + 4A^{16} + 3A^{15} - 4A^{14} + A^{13} + A^{12} + 2A^{11} - 6A^{10} - A^{9} + 6A^{8} - 8A^{7} - A^{6} + 8A^{5} - 4A^{4} - 2A^{3} + 4A^{2} - 1$}
\end{minipage}

\noindent{\color{gray!40}\rule{\textwidth}{0.4pt}}
\vspace{0.9\baselineskip}
\noindent \begin{minipage}[t]{0.25\textwidth}
\vspace{0pt}
\centering
\includegraphics[page=117,width=\linewidth]{knotoids.pdf}
\end{minipage}
\hfill
\begin{minipage}[t]{0.73\textwidth}
\vspace{0pt}
\raggedright
\textbf{Name:} {\large{$\mathbf{K6_{80}}$}} (chiral, non-rotatable$^{*}$) \\ \textbf{PD:} {\scriptsize\texttt{[0],[0,1,2,3],[1,4,5,2],[6,7,4,3],[5,8,9,6],[7,10,11,12],[8,12,10,9],[11]}} \\ \textbf{EM:} {\small\texttt{(B0, A0C0C3D3, B1D2E0B2, E3F0C1B3, C2G0G3D0, D1G2H0G1, E1F3F1E2, F2)}} \\ \textbf{Kauffman bracket:} {\scriptsize $-A^{20} - A^{18} + A^{16} + 3A^{14} + A^{12} - 3A^{10} - 3A^{8} + A^{6} + 2A^{4} + A^{2}$} \\ \textbf{Arrow:} {\scriptsize $-A^{26}L_1 - A^{24} + A^{22}L_1 + 3A^{20} + A^{18}L_1 - 3A^{16} - 3A^{14}L_1 + A^{12} + 2A^{10}L_1 + A^{8}$} \\ \textbf{Mock:} {\scriptsize $2w^{2} - 4 - 2/w + 2/w^{2} + 2/w^{3} + w^{-4}$} \\ \textbf{Affine:} {\scriptsize $-2t + 4 - 2/t$} \\ \textbf{Yamada:} {\scriptsize $A^{26} - 2A^{24} + 2A^{22} - 2A^{21} - 3A^{20} + 4A^{19} + 2A^{18} - 4A^{17} + 4A^{16} + 2A^{15} - 3A^{14} + A^{13} + A^{12} + A^{11} - 4A^{10} + A^{9} + 2A^{8} - 6A^{7} - A^{6} + 3A^{5} - 3A^{4} - 2A^{3} + A^{2} - 1$}
\end{minipage}

\noindent{\color{gray!40}\rule{\textwidth}{0.4pt}}
\vspace{0.9\baselineskip}
\noindent \begin{minipage}[t]{0.25\textwidth}
\vspace{0pt}
\centering
\includegraphics[page=118,width=\linewidth]{knotoids.pdf}
\end{minipage}
\hfill
\begin{minipage}[t]{0.73\textwidth}
\vspace{0pt}
\raggedright
\textbf{Name:} {\large{$\mathbf{K6_{81}}$}} (chiral, rotatable) \\ \textbf{PD:} {\scriptsize\texttt{[0],[0,1,2,3],[1,4,5,2],[6,7,4,3],[8,9,6,5],[7,10,11,12],[12,10,9,8],[11]}} \\ \textbf{EM:} {\small\texttt{(B0, A0C0C3D3, B1D2E3B2, E2F0C1B3, G3G2D0C2, D1G1H0G0, F3F1E1E0, F2)}} \\ \textbf{Kauffman bracket:} {\scriptsize $-A^{18} - 2A^{16} + A^{14} + 4A^{12} + 2A^{10} - 2A^{8} - 2A^{6} + A^{2}$} \\ \textbf{Arrow:} {\scriptsize $-A^{12} - 2A^{10}L_1 + A^{8} + 4A^{6}L_1 + 2A^{4} - 2A^{2}L_1 - 2 + A^{-4}$} \\ \textbf{Mock:} {\scriptsize $-2w + 4/w + 2/w^{2} - 2/w^{3} - 1/w^{4}$} \\ \textbf{Affine:} {\scriptsize $0$} \\ \textbf{Yamada:} {\scriptsize $-A^{22} - A^{21} + 3A^{20} - 3A^{18} + 4A^{17} + 2A^{16} - 4A^{15} + 2A^{14} + A^{13} - 3A^{12} - A^{10} + A^{9} - 4A^{8} - A^{7} + 3A^{6} - 5A^{5} - A^{4} + 4A^{3} - 2A^{2} - A + 1$}
\end{minipage}

\noindent{\color{gray!40}\rule{\textwidth}{0.4pt}}
\vspace{0.9\baselineskip}
\noindent \begin{minipage}[t]{0.25\textwidth}
\vspace{0pt}
\centering
\includegraphics[page=119,width=\linewidth]{knotoids.pdf}
\end{minipage}
\hfill
\begin{minipage}[t]{0.73\textwidth}
\vspace{0pt}
\raggedright
\textbf{Name:} {\large{$\mathbf{K6_{82}}$}} (chiral, rotatable) \\ \textbf{PD:} {\scriptsize\texttt{[0],[0,1,2,3],[1,4,5,6],[6,7,8,2],[9,10,4,3],[10,11,7,5],[11,12,9,8],[12]}} \\ \textbf{EM:} {\small\texttt{(B0, A0C0D3E3, B1E2F3D0, C3F2G3B2, G2F0C1B3, E1G0D1C2, F1H0E0D2, G1)}} \\ \textbf{Kauffman bracket:} {\scriptsize $-2A^{18} - A^{16} + 4A^{14} + 4A^{12} - 2A^{10} - 4A^{8} + 3A^{4} - 1$} \\ \textbf{Arrow:} {\scriptsize $-2A^{6}L_1 - A^{4} + 4A^{2}L_1 + 4 - 2L_1/A^{2} - 4/A^{4} + 3/A^{8} - 1/A^{12}$} \\ \textbf{Mock:} {\scriptsize $-2w^{3} - 2w^{2} + 4w + 6 - 2/w - 4/w^{2} + w^{-4}$} \\ \textbf{Affine:} {\scriptsize $0$} \\ \textbf{Yamada:} {\scriptsize $A^{24} - A^{23} - 3A^{22} + 2A^{21} + 4A^{20} - 4A^{19} - 2A^{18} + 8A^{17} - 3A^{16} - 7A^{15} + 6A^{14} - 4A^{12} + 4A^{11} + 4A^{10} + 2A^{9} - 3A^{8} + 6A^{7} + 3A^{6} - 9A^{5} + 3A^{4} + 4A^{3} - 7A^{2} - A + 3$}
\end{minipage}

\noindent{\color{gray!40}\rule{\textwidth}{0.4pt}}
\vspace{0.9\baselineskip}

\section*{Number of crossings: 7}

\noindent{\color{gray!40}\rule{\textwidth}{0.4pt}}
\vspace{0.9\baselineskip}
\noindent \begin{minipage}[t]{0.25\textwidth}
\vspace{0pt}
\centering
\includegraphics[page=120,width=\linewidth]{knotoids.pdf}
\end{minipage}
\hfill
\begin{minipage}[t]{0.73\textwidth}
\vspace{0pt}
\raggedright
\textbf{Name:} {\large{$\mathbf{K7_{1}}$}} (chiral, rotatable) \\ \textbf{PD:} {\scriptsize\texttt{[0],[0,1,2,3],[1,4,5,2],[3,6,7,8],[4,9,10,5],[6,11,12,7],[8],[9,13,14,10],[11,14,13,12]}} \\ \textbf{EM:} {\scriptsize\texttt{(B0, A0C0C3D0, B1E0E3B2, B3F0F3G0, C1H0H3C2, D1I0I3D2, D3, E1I2I1E2, F1H2H1F2)}} \\ \textbf{Kauffman bracket:} {\scriptsize $A^{30} + A^{22} - A^{18} + A^{14} - A^{10} + A^{6} - A^{2}$} \\ \textbf{Arrow:} {\scriptsize $A^{-12} + A^{-20} - 1/A^{24} + A^{-28} - 1/A^{32} + A^{-36} - 1/A^{40}$} \\ \textbf{Mock:} {\scriptsize $w^{6} - w^{4} + w^{2} - 1 + w^{-2} - 1/w^{4} + w^{-6}$} \\ \textbf{Affine:} {\scriptsize $0$} \\ \textbf{Yamada:} {\scriptsize $A^{25} + A^{24} + 2A^{23} + 2A^{22} + 2A^{21} + A^{20} + A^{19} - A^{10} - A^{9} - 2A^{8} - A^{7} - A^{6} + A^{2} + 1$}
\end{minipage}

\noindent{\color{gray!40}\rule{\textwidth}{0.4pt}}
\vspace{0.9\baselineskip}
\noindent \begin{minipage}[t]{0.25\textwidth}
\vspace{0pt}
\centering
\includegraphics[page=121,width=\linewidth]{knotoids.pdf}
\end{minipage}
\hfill
\begin{minipage}[t]{0.73\textwidth}
\vspace{0pt}
\raggedright
\textbf{Name:} {\large{$\mathbf{K7_{2}}$}} (chiral, rotatable) \\ \textbf{PD:} {\scriptsize\texttt{[0],[0,1,2,3],[1,4,5,2],[3,5,6,7],[4,8,9,10],[11,12,7,6],[8],[13,14,10,9],[14,13,12,11]}} \\ \textbf{EM:} {\scriptsize\texttt{(B0, A0C0C3D0, B1E0D1B2, B3C2F3F2, C1G0H3H2, I3I2D3D2, E1, I1I0E3E2, H1H0F1F0)}} \\ \textbf{Kauffman bracket:} {\scriptsize $-A^{29} + A^{25} - 2A^{21} + 2A^{17} - 2A^{13} + A^{9} - A^{5} + A$} \\ \textbf{Arrow:} {\scriptsize $A^{-4} - 1/A^{8} + 2/A^{12} - 2/A^{16} + 2/A^{20} - 1/A^{24} + A^{-28} - 1/A^{32}$} \\ \textbf{Mock:} {\scriptsize $3w^{2} - 5 + 3/w^{2}$} \\ \textbf{Affine:} {\scriptsize $0$} \\ \textbf{Yamada:} {\scriptsize $A^{25} + A^{23} + A^{22} + 2A^{20} + 2A^{19} + A^{18} + 2A^{17} - A^{14} - A^{13} - A^{12} - A^{11} + A^{9} - A^{5} - A^{4} + 1$}
\end{minipage}

\noindent{\color{gray!40}\rule{\textwidth}{0.4pt}}
\vspace{0.9\baselineskip}
\noindent \begin{minipage}[t]{0.25\textwidth}
\vspace{0pt}
\centering
\includegraphics[page=122,width=\linewidth]{knotoids.pdf}
\end{minipage}
\hfill
\begin{minipage}[t]{0.73\textwidth}
\vspace{0pt}
\raggedright
\textbf{Name:} {\large{$\mathbf{K7_{3}}$}} (chiral, rotatable) \\ \textbf{PD:} {\scriptsize\texttt{[0],[0,1,2,3],[1,4,5,2],[3,6,7,8],[4,9,6,5],[10,11,8,7],[9,12,13,14],[14,13,11,10],[12]}} \\ \textbf{EM:} {\scriptsize\texttt{(B0, A0C0C3D0, B1E0E3B2, B3E2F3F2, C1G0D1C2, H3H2D3D2, E1I0H1H0, G3G2F1F0, G1)}} \\ \textbf{Kauffman bracket:} {\scriptsize $-A^{30} + A^{26} - 2A^{22} + 3A^{18} - 2A^{14} + 2A^{10} - A^{6} + A^{2}$} \\ \textbf{Arrow:} {\scriptsize $-A^{36} + A^{32} - 2A^{28} + 3A^{24} - 2A^{20} + 2A^{16} - A^{12} + A^{8}$} \\ \textbf{Mock:} {\scriptsize $2w^{4} - 3w^{2} + 3 - 3/w^{2} + 2/w^{4}$} \\ \textbf{Affine:} {\scriptsize $0$} \\ \textbf{Yamada:} {\scriptsize $-A^{25} - A^{23} - A^{22} + A^{21} + A^{20} + 3A^{18} + A^{17} + 2A^{15} - A^{14} + A^{13} - A^{12} - 2A^{9} - 2A^{6} - A^{5} - A^{4} - 2A^{3} - A^{2} - 1$}
\end{minipage}

\noindent{\color{gray!40}\rule{\textwidth}{0.4pt}}
\vspace{0.9\baselineskip}
\noindent \begin{minipage}[t]{0.25\textwidth}
\vspace{0pt}
\centering
\includegraphics[page=123,width=\linewidth]{knotoids.pdf}
\end{minipage}
\hfill
\begin{minipage}[t]{0.73\textwidth}
\vspace{0pt}
\raggedright
\textbf{Name:} {\large{$\mathbf{K7_{4}}$}} (chiral, rotatable) \\ \textbf{PD:} {\scriptsize\texttt{[0],[0,1,2,3],[1,4,5,2],[3,5,6,7],[4,8,9,10],[10,11,7,6],[8,12,13,9],[11,13,12,14],[14]}} \\ \textbf{EM:} {\scriptsize\texttt{(B0, A0C0C3D0, B1E0D1B2, B3C2F3F2, C1G0G3F0, E3H0D3D2, E1H2H1E2, F1G2G1I0, H3)}} \\ \textbf{Kauffman bracket:} {\scriptsize $A^{28} - A^{24} + 3A^{20} - 3A^{16} + 3A^{12} - 3A^{8} + 2A^{4} - 1$} \\ \textbf{Arrow:} {\scriptsize $A^{-8} - 1/A^{12} + 3/A^{16} - 3/A^{20} + 3/A^{24} - 3/A^{28} + 2/A^{32} - 1/A^{36}$} \\ \textbf{Mock:} {\scriptsize $2w^{4} - 4w^{2} + 5 - 4/w^{2} + 2/w^{4}$} \\ \textbf{Affine:} {\scriptsize $0$} \\ \textbf{Yamada:} {\scriptsize $A^{25} + A^{23} + 3A^{22} + A^{21} + A^{20} + 4A^{19} - A^{18} + 2A^{16} - 2A^{15} - A^{14} - 2A^{12} - 2A^{10} + A^{9} + A^{8} - 3A^{7} + 2A^{6} - 2A^{4} + 2A^{3} - A + 1$}
\end{minipage}

\noindent{\color{gray!40}\rule{\textwidth}{0.4pt}}
\vspace{0.9\baselineskip}
\noindent \begin{minipage}[t]{0.25\textwidth}
\vspace{0pt}
\centering
\includegraphics[page=124,width=\linewidth]{knotoids.pdf}
\end{minipage}
\hfill
\begin{minipage}[t]{0.73\textwidth}
\vspace{0pt}
\raggedright
\textbf{Name:} {\large{$\mathbf{K7_{5}}$}} (chiral, rotatable) \\ \textbf{PD:} {\scriptsize\texttt{[0],[0,1,2,3],[1,4,5,2],[3,6,7,8],[4,9,6,5],[9,10,11,7],[8,12,13,14],[10,13,12,11],[14]}} \\ \textbf{EM:} {\scriptsize\texttt{(B0, A0C0C3D0, B1E0E3B2, B3E2F3G0, C1F0D1C2, E1H0H3D2, D3H2H1I0, F1G2G1F2, G3)}} \\ \textbf{Kauffman bracket:} {\scriptsize $A^{29} - A^{25} + 2A^{21} - 3A^{17} + 2A^{13} - 3A^{9} + 2A^{5} - A$} \\ \textbf{Arrow:} {\scriptsize $-A^{32} + A^{28} - 2A^{24} + 3A^{20} - 2A^{16} + 3A^{12} - 2A^{8} + A^{4}$} \\ \textbf{Mock:} {\scriptsize $4w^{2} - 7 + 4/w^{2}$} \\ \textbf{Affine:} {\scriptsize $0$} \\ \textbf{Yamada:} {\scriptsize $-A^{25} - A^{22} + A^{20} - A^{19} + 2A^{18} + 3A^{17} + A^{16} + 2A^{15} + A^{14} - A^{13} - 2A^{12} - A^{11} - 2A^{10} - 2A^{8} + A^{7} - A^{6} - 3A^{5} + A^{4} - 2A^{3} - A^{2} + A - 1$}
\end{minipage}

\noindent{\color{gray!40}\rule{\textwidth}{0.4pt}}
\vspace{0.9\baselineskip}
\noindent \begin{minipage}[t]{0.25\textwidth}
\vspace{0pt}
\centering
\includegraphics[page=125,width=\linewidth]{knotoids.pdf}
\end{minipage}
\hfill
\begin{minipage}[t]{0.73\textwidth}
\vspace{0pt}
\raggedright
\textbf{Name:} {\large{$\mathbf{K7_{6}}$}} (chiral, rotatable) \\ \textbf{PD:} {\scriptsize\texttt{[0],[0,1,2,3],[1,4,5,2],[3,5,6,7],[4,8,9,6],[7,10,11,12],[8,13,10,9],[13,14,12,11],[14]}} \\ \textbf{EM:} {\scriptsize\texttt{(B0, A0C0C3D0, B1E0D1B2, B3C2E3F0, C1G0G3D2, D3G2H3H2, E1H0F1E2, G1I0F3F2, H1)}} \\ \textbf{Kauffman bracket:} {\scriptsize $A^{28} - 2A^{24} + 3A^{20} - 3A^{16} + 4A^{12} - 3A^{8} + 2A^{4} - 1$} \\ \textbf{Arrow:} {\scriptsize $A^{4} - 2 + 3/A^{4} - 3/A^{8} + 4/A^{12} - 3/A^{16} + 2/A^{20} - 1/A^{24}$} \\ \textbf{Mock:} {\scriptsize $-w^{4} + 5w^{2} - 7 + 5/w^{2} - 1/w^{4}$} \\ \textbf{Affine:} {\scriptsize $0$} \\ \textbf{Yamada:} {\scriptsize $A^{25} - A^{24} - A^{23} + 2A^{22} - A^{21} + 4A^{19} - A^{18} + 3A^{16} - A^{15} + 2A^{13} + 2A^{11} - 2A^{10} + A^{9} + A^{8} - 4A^{7} + 2A^{6} - 3A^{4} + 2A^{3} - A + 1$}
\end{minipage}

\noindent{\color{gray!40}\rule{\textwidth}{0.4pt}}
\vspace{0.9\baselineskip}
\noindent \begin{minipage}[t]{0.25\textwidth}
\vspace{0pt}
\centering
\includegraphics[page=126,width=\linewidth]{knotoids.pdf}
\end{minipage}
\hfill
\begin{minipage}[t]{0.73\textwidth}
\vspace{0pt}
\raggedright
\textbf{Name:} {\large{$\mathbf{K7_{7}}$}} (chiral, rotatable) \\ \textbf{PD:} {\scriptsize\texttt{[0],[0,1,2,3],[1,4,5,2],[3,5,6,7],[4,8,9,6],[7,9,10,11],[8,12,13,10],[11,13,12,14],[14]}} \\ \textbf{EM:} {\scriptsize\texttt{(B0, A0C0C3D0, B1E0D1B2, B3C2E3F0, C1G0F1D2, D3E2G3H0, E1H2H1F2, F3G2G1I0, H3)}} \\ \textbf{Kauffman bracket:} {\scriptsize $A^{28} - 2A^{24} + 3A^{20} - 4A^{16} + 4A^{12} - 3A^{8} + 3A^{4} - 1$} \\ \textbf{Arrow:} {\scriptsize $A^{16} - 2A^{12} + 3A^{8} - 4A^{4} + 4 - 3/A^{4} + 3/A^{8} - 1/A^{12}$} \\ \textbf{Mock:} {\scriptsize $w^{4} - 5w^{2} + 9 - 5/w^{2} + w^{-4}$} \\ \textbf{Affine:} {\scriptsize $0$} \\ \textbf{Yamada:} {\scriptsize $A^{25} - A^{24} - A^{23} + 2A^{22} - 2A^{21} - A^{20} + 4A^{19} - 2A^{18} + 4A^{16} - 2A^{15} + A^{13} - A^{12} + A^{11} - 2A^{10} + 2A^{9} + 2A^{8} - 3A^{7} + 4A^{6} + 2A^{5} - 3A^{4} + 3A^{3} - A^{2} - 2A + 1$}
\end{minipage}

\noindent{\color{gray!40}\rule{\textwidth}{0.4pt}}
\vspace{0.9\baselineskip}
\noindent \begin{minipage}[t]{0.25\textwidth}
\vspace{0pt}
\centering
\includegraphics[page=127,width=\linewidth]{knotoids.pdf}
\end{minipage}
\hfill
\begin{minipage}[t]{0.73\textwidth}
\vspace{0pt}
\raggedright
\textbf{Name:} {\large{$\mathbf{K7_{8}}$}} (chiral, non-rotatable$^{*}$) \\ \textbf{PD:} {\scriptsize\texttt{[0],[0,1,2,3],[1,4,5,2],[3,5,6,7],[8,9,10,4],[11,12,7,6],[13,14,9,8],[10],[14,13,12,11]}} \\ \textbf{EM:} {\scriptsize\texttt{(B0, A0C0C3D0, B1E3D1B2, B3C2F3F2, G3G2H0C1, I3I2D3D2, I1I0E1E0, E2, G1G0F1F0)}} \\ \textbf{Kauffman bracket:} {\scriptsize $A^{26} - A^{22} + 2A^{18} - 2A^{14} + 2A^{10} - 2A^{6} - A^{4} + A^{2} + 1$} \\ \textbf{Arrow:} {\scriptsize $A^{8} - A^{4} + 2 - 2/A^{4} + 2/A^{8} - 2/A^{12} - L_1/A^{14} + A^{-16} + L_1/A^{18}$} \\ \textbf{Mock:} {\scriptsize $-3w^{2} - w + 6 + 1/w - 2/w^{2}$} \\ \textbf{Affine:} {\scriptsize $-t + 2 - 1/t$} \\ \textbf{Yamada:} {\scriptsize $-A^{26} - A^{24} + A^{22} - 2A^{21} + A^{19} - 2A^{18} + A^{15} - A^{14} - A^{13} + A^{12} - 2A^{11} - A^{10} + A^{9} - A^{7} + A^{6} + A^{5} - A^{4} + A^{2} - 1$}
\end{minipage}

\noindent{\color{gray!40}\rule{\textwidth}{0.4pt}}
\vspace{0.9\baselineskip}
\noindent \begin{minipage}[t]{0.25\textwidth}
\vspace{0pt}
\centering
\includegraphics[page=128,width=\linewidth]{knotoids.pdf}
\end{minipage}
\hfill
\begin{minipage}[t]{0.73\textwidth}
\vspace{0pt}
\raggedright
\textbf{Name:} {\large{$\mathbf{K7_{9}}$}} (chiral, non-rotatable$^{*}$) \\ \textbf{PD:} {\scriptsize\texttt{[0],[0,1,2,3],[1,4,5,2],[3,5,6,7],[4,8,9,6],[10,11,12,7],[8,13,14,9],[10],[11,14,13,12]}} \\ \textbf{EM:} {\scriptsize\texttt{(B0, A0C0C3D0, B1E0D1B2, B3C2E3F3, C1G0G3D2, H0I0I3D3, E1I2I1E2, F0, F1G2G1F2)}} \\ \textbf{Kauffman bracket:} {\scriptsize $-A^{27} - A^{25} + A^{23} + 2A^{21} - A^{19} - 3A^{17} + 3A^{13} - 2A^{9} + 2A^{5} - A$} \\ \textbf{Arrow:} {\scriptsize $1 + L_1/A^{2} - 1/A^{4} - 2L_1/A^{6} + A^{-8} + 3L_1/A^{10} - 3L_1/A^{14} + 2L_1/A^{18} - 2L_1/A^{22} + L_1/A^{26}$} \\ \textbf{Mock:} {\scriptsize $-w^{6} - w^{5} + 2w^{4} + 3w^{3} - 3w + 3/w - 3/w^{3} + w^{-5}$} \\ \textbf{Affine:} {\scriptsize $t - 2 + 1/t$} \\ \textbf{Yamada:} {\scriptsize $-A^{26} + A^{24} - A^{23} - 2A^{22} + A^{21} - 3A^{19} + 2A^{17} - 3A^{16} - A^{15} + 2A^{14} - A^{13} - A^{12} + 2A^{11} - A^{9} - 2A^{8} + 2A^{7} - 3A^{5} + 4A^{4} - A^{2} + A - 1$}
\end{minipage}

\noindent{\color{gray!40}\rule{\textwidth}{0.4pt}}
\vspace{0.9\baselineskip}
\noindent \begin{minipage}[t]{0.25\textwidth}
\vspace{0pt}
\centering
\includegraphics[page=129,width=\linewidth]{knotoids.pdf}
\end{minipage}
\hfill
\begin{minipage}[t]{0.73\textwidth}
\vspace{0pt}
\raggedright
\textbf{Name:} {\large{$\mathbf{K7_{10}}$}} (chiral, non-rotatable$^{*}$) \\ \textbf{PD:} {\scriptsize\texttt{[0],[0,1,2,3],[1,4,5,2],[6,7,8,3],[4,8,9,5],[10,11,7,6],[12,13,14,9],[13,12,11,10],[14]}} \\ \textbf{EM:} {\scriptsize\texttt{(B0, A0C0C3D3, B1E0E3B2, F3F2E1B3, C1D2G3C2, H3H2D1D0, H1H0I0E2, G1G0F1F0, G2)}} \\ \textbf{Kauffman bracket:} {\scriptsize $A^{28} + A^{26} - A^{24} - 2A^{22} + A^{20} + 3A^{18} - A^{16} - 3A^{14} + 2A^{10} - A^{6} + A^{2}$} \\ \textbf{Arrow:} {\scriptsize $A^{10}L_1 + A^{8} - A^{6}L_1 - 2A^{4} + A^{2}L_1 + 3 - L_1/A^{2} - 3/A^{4} + 2/A^{8} - 1/A^{12} + A^{-16}$} \\ \textbf{Mock:} {\scriptsize $-2w^{2} + 2w + 7 - 2/w - 4/w^{2}$} \\ \textbf{Affine:} {\scriptsize $2t - 4 + 2/t$} \\ \textbf{Yamada:} {\scriptsize $-A^{26} - A^{23} + 2A^{21} - 2A^{20} - A^{19} + 4A^{18} - 2A^{17} - A^{16} + 2A^{15} - 2A^{14} - A^{13} - 2A^{12} + A^{11} - A^{10} - 3A^{9} + 3A^{8} - 2A^{6} + 2A^{5} + A^{4} - 2A^{3} + A - 1$}
\end{minipage}

\noindent{\color{gray!40}\rule{\textwidth}{0.4pt}}
\vspace{0.9\baselineskip}
\noindent \begin{minipage}[t]{0.25\textwidth}
\vspace{0pt}
\centering
\includegraphics[page=130,width=\linewidth]{knotoids.pdf}
\end{minipage}
\hfill
\begin{minipage}[t]{0.73\textwidth}
\vspace{0pt}
\raggedright
\textbf{Name:} {\large{$\mathbf{K7_{11}}$}} (chiral, non-rotatable$^{*}$) \\ \textbf{PD:} {\scriptsize\texttt{[0],[0,1,2,3],[1,4,5,2],[3,6,7,8],[4,9,6,5],[9,10,11,7],[12,13,14,8],[10,14,13,11],[12]}} \\ \textbf{EM:} {\scriptsize\texttt{(B0, A0C0C3D0, B1E0E3B2, B3E2F3G3, C1F0D1C2, E1H0H3D2, I0H2H1D3, F1G2G1F2, G0)}} \\ \textbf{Kauffman bracket:} {\scriptsize $-A^{27} - A^{25} + A^{23} + 2A^{21} - A^{19} - 3A^{17} + A^{15} + 3A^{13} - 3A^{9} + 2A^{5} - A$} \\ \textbf{Arrow:} {\scriptsize $A^{6}L_1 + A^{4} - A^{2}L_1 - 2 + L_1/A^{2} + 3/A^{4} - L_1/A^{6} - 3/A^{8} + 3/A^{12} - 2/A^{16} + A^{-20}$} \\ \textbf{Mock:} {\scriptsize $-w^{4} + w^{3} + 4w^{2} - w - 4 + 1/w + 4/w^{2} - 1/w^{3} - 2/w^{4}$} \\ \textbf{Affine:} {\scriptsize $2t - 4 + 2/t$} \\ \textbf{Yamada:} {\scriptsize $-A^{26} + A^{24} - A^{23} - 2A^{22} + 2A^{21} - 3A^{19} + 2A^{18} + 2A^{17} - 3A^{16} + 2A^{14} - 3A^{13} - 2A^{12} + A^{11} - A^{10} - A^{9} - A^{8} + 4A^{7} - A^{6} - 3A^{5} + 4A^{4} - 2A^{3} - A^{2} + 2A - 1$}
\end{minipage}

\noindent{\color{gray!40}\rule{\textwidth}{0.4pt}}
\vspace{0.9\baselineskip}
\noindent \begin{minipage}[t]{0.25\textwidth}
\vspace{0pt}
\centering
\includegraphics[page=131,width=\linewidth]{knotoids.pdf}
\end{minipage}
\hfill
\begin{minipage}[t]{0.73\textwidth}
\vspace{0pt}
\raggedright
\textbf{Name:} {\large{$\mathbf{K7_{12}}$}} (chiral, non-rotatable$^{*}$) \\ \textbf{PD:} {\scriptsize\texttt{[0],[0,1,2,3],[1,3,4,5],[5,6,7,2],[8,9,10,4],[6,11,12,7],[8],[9,13,14,10],[11,14,13,12]}} \\ \textbf{EM:} {\scriptsize\texttt{(B0, A0C0D3C1, B1B3E3D0, C3F0F3B2, G0H0H3C2, D1I0I3D2, E0, E1I2I1E2, F1H2H1F2)}} \\ \textbf{Kauffman bracket:} {\scriptsize $-A^{27} - A^{25} + A^{23} + A^{21} - A^{19} - 2A^{17} + A^{15} + 2A^{13} - A^{11} - 2A^{9} + A^{7} + 2A^{5} - A$} \\ \textbf{Arrow:} {\scriptsize $L_1/A^{6} + A^{-8} - L_1/A^{10} - 1/A^{12} + L_1/A^{14} + 2/A^{16} - L_1/A^{18} - 2/A^{20} + L_1/A^{22} + 2/A^{24} - L_1/A^{26} - 2/A^{28} + A^{-32}$} \\ \textbf{Mock:} {\scriptsize $w^{5} + 2w^{4} - w^{3} - 2w^{2} + w + 2 - 1/w - 2/w^{2} + w^{-3} + 2/w^{4} - 1/w^{5} - 1/w^{6}$} \\ \textbf{Affine:} {\scriptsize $3t - 6 + 3/t$} \\ \textbf{Yamada:} {\scriptsize $A^{27} + 2A^{24} + A^{23} - A^{22} + 2A^{21} + 2A^{20} + A^{18} + 2A^{17} - A^{16} - A^{15} - A^{13} - 2A^{12} - A^{8} + 2A^{7} + A^{6} - 2A^{5} + 2A^{4} - A^{3} + A - 1$}
\end{minipage}

\noindent{\color{gray!40}\rule{\textwidth}{0.4pt}}
\vspace{0.9\baselineskip}
\noindent \begin{minipage}[t]{0.25\textwidth}
\vspace{0pt}
\centering
\includegraphics[page=132,width=\linewidth]{knotoids.pdf}
\end{minipage}
\hfill
\begin{minipage}[t]{0.73\textwidth}
\vspace{0pt}
\raggedright
\textbf{Name:} {\large{$\mathbf{K7_{13}}$}} (chiral, non-rotatable$^{*}$) \\ \textbf{PD:} {\scriptsize\texttt{[0],[0,1,2,3],[1,4,5,2],[3,5,6,7],[4,8,9,6],[7,9,10,11],[12,13,14,8],[13,12,11,10],[14]}} \\ \textbf{EM:} {\scriptsize\texttt{(B0, A0C0C3D0, B1E0D1B2, B3C2E3F0, C1G3F1D2, D3E2H3H2, H1H0I0E1, G1G0F3F2, G2)}} \\ \textbf{Kauffman bracket:} {\scriptsize $A^{27} - 2A^{23} + 3A^{19} - 4A^{15} - A^{13} + 4A^{11} + A^{9} - 3A^{7} - 2A^{5} + A^{3} + A$} \\ \textbf{Arrow:} {\scriptsize $-A^{18}L_1 + 2A^{14}L_1 - 3A^{10}L_1 + 4A^{6}L_1 + A^{4} - 4A^{2}L_1 - 1 + 3L_1/A^{2} + 2/A^{4} - L_1/A^{6} - 1/A^{8}$} \\ \textbf{Mock:} {\scriptsize $-w^{5} + 3w^{3} - 5w + 5/w + 2/w^{2} - 3/w^{3} - 2/w^{4} + w^{-5} + w^{-6}$} \\ \textbf{Affine:} {\scriptsize $-t + 2 - 1/t$} \\ \textbf{Yamada:} {\scriptsize $A^{26} - 2A^{25} + 3A^{23} - 5A^{22} + 6A^{20} - 5A^{19} - A^{18} + 5A^{17} - 2A^{16} - A^{15} + 4A^{13} - 2A^{11} + 7A^{10} - 5A^{8} + 5A^{7} + A^{6} - 4A^{5} + A^{4} + 2A^{3} - 2A^{2} - A + 1$}
\end{minipage}

\noindent{\color{gray!40}\rule{\textwidth}{0.4pt}}
\vspace{0.9\baselineskip}
\noindent \begin{minipage}[t]{0.25\textwidth}
\vspace{0pt}
\centering
\includegraphics[page=133,width=\linewidth]{knotoids.pdf}
\end{minipage}
\hfill
\begin{minipage}[t]{0.73\textwidth}
\vspace{0pt}
\raggedright
\textbf{Name:} {\large{$\mathbf{K7_{14}}$}} (chiral, non-rotatable$^{*}$) \\ \textbf{PD:} {\scriptsize\texttt{[0],[0,1,2,3],[1,4,5,2],[3,5,6,7],[4,8,9,10],[11,12,7,6],[8,13,14,9],[10,13,12,11],[14]}} \\ \textbf{EM:} {\scriptsize\texttt{(B0, A0C0C3D0, B1E0D1B2, B3C2F3F2, C1G0G3H0, H3H2D3D2, E1H1I0E2, E3G1F1F0, G2)}} \\ \textbf{Kauffman bracket:} {\scriptsize $A^{28} + A^{26} - A^{24} - 2A^{22} + 2A^{20} + 3A^{18} - 2A^{16} - 4A^{14} + A^{12} + 3A^{10} - 2A^{6} + A^{2}$} \\ \textbf{Arrow:} {\scriptsize $A^{4} + A^{2}L_1 - 1 - 2L_1/A^{2} + 2/A^{4} + 3L_1/A^{6} - 2/A^{8} - 4L_1/A^{10} + A^{-12} + 3L_1/A^{14} - 2L_1/A^{18} + L_1/A^{22}$} \\ \textbf{Mock:} {\scriptsize $-2w^{4} - 2w^{3} + 4w^{2} + 6w - 1 - 6/w + 2/w^{3}$} \\ \textbf{Affine:} {\scriptsize $0$} \\ \textbf{Yamada:} {\scriptsize $-A^{26} - 2A^{23} + A^{22} + 2A^{21} - 5A^{20} + 2A^{19} + 4A^{18} - 6A^{17} + A^{16} + 3A^{15} - 4A^{14} - A^{13} - A^{12} + 2A^{11} - 2A^{10} - 2A^{9} + 6A^{8} - 3A^{7} - 3A^{6} + 5A^{5} - A^{4} - 3A^{3} + 3A^{2} + A - 2$}
\end{minipage}

\noindent{\color{gray!40}\rule{\textwidth}{0.4pt}}
\vspace{0.9\baselineskip}
\noindent \begin{minipage}[t]{0.25\textwidth}
\vspace{0pt}
\centering
\includegraphics[page=134,width=\linewidth]{knotoids.pdf}
\end{minipage}
\hfill
\begin{minipage}[t]{0.73\textwidth}
\vspace{0pt}
\raggedright
\textbf{Name:} {\large{$\mathbf{K7_{15}}$}} (chiral, non-rotatable$^{*}$) \\ \textbf{PD:} {\scriptsize\texttt{[0],[0,1,2,3],[1,4,5,2],[6,7,4,3],[7,8,9,5],[10,11,12,6],[8,13,14,9],[10],[11,14,13,12]}} \\ \textbf{EM:} {\scriptsize\texttt{(B0, A0C0C3D3, B1D2E3B2, F3E0C1B3, D1G0G3C2, H0I0I3D0, E1I2I1E2, F0, F1G2G1F2)}} \\ \textbf{Kauffman bracket:} {\scriptsize $-A^{27} - A^{25} + A^{23} + 2A^{21} - A^{19} - 2A^{17} + A^{15} + 3A^{13} - A^{11} - 3A^{9} + 2A^{5} - A$} \\ \textbf{Arrow:} {\scriptsize $A^{12} + A^{10}L_1 - A^{8} - 2A^{6}L_1 + A^{4} + 2A^{2}L_1 - 1 - 3L_1/A^{2} + A^{-4} + 3L_1/A^{6} - 2L_1/A^{10} + L_1/A^{14}$} \\ \textbf{Mock:} {\scriptsize $-2w^{4} - 2w^{3} + 3w^{2} + 5w - 5/w + 2/w^{3}$} \\ \textbf{Affine:} {\scriptsize $-t + 2 - 1/t$} \\ \textbf{Yamada:} {\scriptsize $-A^{26} + A^{24} - A^{23} - A^{22} + 3A^{21} - 2A^{19} + 2A^{18} - 3A^{16} + A^{14} - 3A^{13} - A^{12} + A^{11} - A^{10} - A^{9} - A^{8} + 3A^{7} - 2A^{6} - A^{5} + 4A^{4} - 2A^{3} + A - 2$}
\end{minipage}

\noindent{\color{gray!40}\rule{\textwidth}{0.4pt}}
\vspace{0.9\baselineskip}
\noindent \begin{minipage}[t]{0.25\textwidth}
\vspace{0pt}
\centering
\includegraphics[page=135,width=\linewidth]{knotoids.pdf}
\end{minipage}
\hfill
\begin{minipage}[t]{0.73\textwidth}
\vspace{0pt}
\raggedright
\textbf{Name:} {\large{$\mathbf{K7_{16}}$}} (chiral, non-rotatable$^{*}$) \\ \textbf{PD:} {\scriptsize\texttt{[0],[0,1,2,3],[1,4,5,2],[3,5,6,7],[4,8,9,6],[7,10,11,12],[8,13,10,9],[14,13,12,11],[14]}} \\ \textbf{EM:} {\scriptsize\texttt{(B0, A0C0C3D0, B1E0D1B2, B3C2E3F0, C1G0G3D2, D3G2H3H2, E1H1F1E2, I0G1F3F2, H0)}} \\ \textbf{Kauffman bracket:} {\scriptsize $A^{28} - 2A^{24} + 4A^{20} + A^{18} - 4A^{16} - 2A^{14} + 4A^{12} + 2A^{10} - 3A^{8} - 2A^{6} + A^{4} + A^{2}$} \\ \textbf{Arrow:} {\scriptsize $A^{4} - 2 + 4/A^{4} + L_1/A^{6} - 4/A^{8} - 2L_1/A^{10} + 4/A^{12} + 2L_1/A^{14} - 3/A^{16} - 2L_1/A^{18} + A^{-20} + L_1/A^{22}$} \\ \textbf{Mock:} {\scriptsize $-2w^{4} - w^{3} + 6w^{2} + 3w - 6 - 3/w + 4/w^{2} + w^{-3} - 1/w^{4}$} \\ \textbf{Affine:} {\scriptsize $0$} \\ \textbf{Yamada:} {\scriptsize $-2A^{26} + A^{25} + 3A^{24} - 5A^{23} - A^{22} + 6A^{21} - 7A^{20} - 3A^{19} + 7A^{18} - 4A^{17} - 3A^{16} + 4A^{15} - A^{14} - A^{13} - 3A^{12} + 5A^{11} + A^{10} - 8A^{9} + 7A^{8} + A^{7} - 7A^{6} + 4A^{5} + 2A^{4} - 3A^{3} + A^{2} + A - 1$}
\end{minipage}

\noindent{\color{gray!40}\rule{\textwidth}{0.4pt}}
\vspace{0.9\baselineskip}
\noindent \begin{minipage}[t]{0.25\textwidth}
\vspace{0pt}
\centering
\includegraphics[page=136,width=\linewidth]{knotoids.pdf}
\end{minipage}
\hfill
\begin{minipage}[t]{0.73\textwidth}
\vspace{0pt}
\raggedright
\textbf{Name:} {\large{$\mathbf{K7_{17}}$}} (chiral, non-rotatable$^{*}$) \\ \textbf{PD:} {\scriptsize\texttt{[0],[0,1,2,3],[1,4,5,2],[6,7,4,3],[7,8,9,5],[10,11,8,6],[12,13,14,9],[13,12,11,10],[14]}} \\ \textbf{EM:} {\scriptsize\texttt{(B0, A0C0C3D3, B1D2E3B2, F3E0C1B3, D1F2G3C2, H3H2E1D0, H1H0I0E2, G1G0F1F0, G2)}} \\ \textbf{Kauffman bracket:} {\scriptsize $A^{28} + A^{26} - 2A^{24} - 3A^{22} + 2A^{20} + 4A^{18} - A^{16} - 4A^{14} + A^{12} + 3A^{10} - 2A^{6} + A^{2}$} \\ \textbf{Arrow:} {\scriptsize $A^{16} + A^{14}L_1 - 2A^{12} - 3A^{10}L_1 + 2A^{8} + 4A^{6}L_1 - A^{4} - 4A^{2}L_1 + 1 + 3L_1/A^{2} - 2L_1/A^{6} + L_1/A^{10}$} \\ \textbf{Mock:} {\scriptsize $-w^{6} - w^{5} + 2w^{4} + 3w^{3} - 2w^{2} - 5w + 2 + 5/w - 3/w^{3} + w^{-5}$} \\ \textbf{Affine:} {\scriptsize $-t + 2 - 1/t$} \\ \textbf{Yamada:} {\scriptsize $-A^{26} + A^{25} + A^{24} - 3A^{23} + A^{22} + 4A^{21} - 5A^{20} + 7A^{18} - 6A^{17} - A^{16} + 4A^{15} - 3A^{14} - A^{13} - 2A^{12} + 3A^{11} - 3A^{10} - 5A^{9} + 6A^{8} - A^{7} - 5A^{6} + 5A^{5} + A^{4} - 4A^{3} + 2A^{2} + A - 2$}
\end{minipage}

\noindent{\color{gray!40}\rule{\textwidth}{0.4pt}}
\vspace{0.9\baselineskip}
\noindent \begin{minipage}[t]{0.25\textwidth}
\vspace{0pt}
\centering
\includegraphics[page=137,width=\linewidth]{knotoids.pdf}
\end{minipage}
\hfill
\begin{minipage}[t]{0.73\textwidth}
\vspace{0pt}
\raggedright
\textbf{Name:} {\large{$\mathbf{K7_{18}}$}} (chiral, non-rotatable$^{*}$) \\ \textbf{PD:} {\scriptsize\texttt{[0],[0,1,2,3],[1,4,5,2],[3,5,6,7],[4,8,9,6],[7,9,10,11],[8,12,13,10],[14,13,12,11],[14]}} \\ \textbf{EM:} {\scriptsize\texttt{(B0, A0C0C3D0, B1E0D1B2, B3C2E3F0, C1G0F1D2, D3E2G3H3, E1H2H1F2, I0G2G1F3, H0)}} \\ \textbf{Kauffman bracket:} {\scriptsize $-A^{26} - A^{24} + 2A^{22} + 3A^{20} - 2A^{18} - 4A^{16} + 2A^{14} + 5A^{12} - A^{10} - 4A^{8} + 3A^{4} - 1$} \\ \textbf{Arrow:} {\scriptsize $-A^{14}L_1 - A^{12} + 2A^{10}L_1 + 3A^{8} - 2A^{6}L_1 - 4A^{4} + 2A^{2}L_1 + 5 - L_1/A^{2} - 4/A^{4} + 3/A^{8} - 1/A^{12}$} \\ \textbf{Mock:} {\scriptsize $w^{4} - w^{3} - 4w^{2} + 3w + 8 - 3/w - 6/w^{2} + w^{-3} + 2/w^{4}$} \\ \textbf{Affine:} {\scriptsize $0$} \\ \textbf{Yamada:} {\scriptsize $A^{26} - A^{25} - 2A^{24} + 3A^{23} - 6A^{21} + 5A^{20} + 4A^{19} - 9A^{18} + 4A^{17} + 6A^{16} - 7A^{15} + A^{14} + 5A^{13} - A^{12} - A^{11} + 7A^{9} - 4A^{8} - 4A^{7} + 11A^{6} - 4A^{5} - 5A^{4} + 7A^{3} - 3A^{2} - 3A + 2$}
\end{minipage}

\noindent{\color{gray!40}\rule{\textwidth}{0.4pt}}
\vspace{0.9\baselineskip}
\noindent \begin{minipage}[t]{0.25\textwidth}
\vspace{0pt}
\centering
\includegraphics[page=138,width=\linewidth]{knotoids.pdf}
\end{minipage}
\hfill
\begin{minipage}[t]{0.73\textwidth}
\vspace{0pt}
\raggedright
\textbf{Name:} {\large{$\mathbf{K7_{19}}$}} (chiral, non-rotatable$^{*}$) \\ \textbf{PD:} {\scriptsize\texttt{[0],[0,1,2,3],[1,4,5,2],[3,6,7,8],[4,9,6,5],[10,11,8,7],[9,12,13,10],[11,13,14,12],[14]}} \\ \textbf{EM:} {\scriptsize\texttt{(B0, A0C0C3D0, B1E0E3B2, B3E2F3F2, C1G0D1C2, G3H0D3D2, E1H3H1F0, F1G2I0G1, H2)}} \\ \textbf{Kauffman bracket:} {\scriptsize $A^{28} - 2A^{24} - A^{22} + 3A^{20} + A^{18} - 4A^{16} - 2A^{14} + 4A^{12} + 3A^{10} - 2A^{8} - 2A^{6} + A^{4} + A^{2}$} \\ \textbf{Arrow:} {\scriptsize $A^{28} - 2A^{24} - A^{22}L_1 + 3A^{20} + A^{18}L_1 - 4A^{16} - 2A^{14}L_1 + 4A^{12} + 3A^{10}L_1 - 2A^{8} - 2A^{6}L_1 + A^{4} + A^{2}L_1$} \\ \textbf{Mock:} {\scriptsize $-2w^{4} - 2w^{3} + 5w^{2} + 3w - 6 - 3/w + 4/w^{2} + 2/w^{3}$} \\ \textbf{Affine:} {\scriptsize $-3t + 6 - 3/t$} \\ \textbf{Yamada:} {\scriptsize $A^{27} - A^{26} + 3A^{24} - 4A^{23} - 2A^{22} + 6A^{21} - 4A^{20} - 4A^{19} + 8A^{18} - 3A^{17} - 3A^{16} + 4A^{15} + A^{13} - 3A^{12} + 4A^{11} + A^{10} - 8A^{9} + 4A^{8} + 2A^{7} - 8A^{6} + 2A^{5} + 2A^{4} - 4A^{3} + A - 1$}
\end{minipage}

\noindent{\color{gray!40}\rule{\textwidth}{0.4pt}}
\vspace{0.9\baselineskip}
\noindent \begin{minipage}[t]{0.25\textwidth}
\vspace{0pt}
\centering
\includegraphics[page=139,width=\linewidth]{knotoids.pdf}
\end{minipage}
\hfill
\begin{minipage}[t]{0.73\textwidth}
\vspace{0pt}
\raggedright
\textbf{Name:} {\large{$\mathbf{K7_{20}}$}} (chiral, non-rotatable$^{*}$) \\ \textbf{PD:} {\scriptsize\texttt{[0],[0,1,2,3],[1,4,5,2],[3,6,7,4],[8,9,10,5],[6,11,12,7],[8],[9,13,14,10],[11,14,13,12]}} \\ \textbf{EM:} {\scriptsize\texttt{(B0, A0C0C3D0, B1D3E3B2, B3F0F3C1, G0H0H3C2, D1I0I3D2, E0, E1I2I1E2, F1H2H1F2)}} \\ \textbf{Kauffman bracket:} {\scriptsize $-A^{27} + A^{23} + A^{21} - 2A^{19} - 2A^{17} + A^{15} + 2A^{13} - A^{11} - 2A^{9} + A^{7} + 2A^{5} - A$} \\ \textbf{Arrow:} {\scriptsize $1 - 1/A^{4} - L_1/A^{6} + 2/A^{8} + 2L_1/A^{10} - 1/A^{12} - 2L_1/A^{14} + A^{-16} + 2L_1/A^{18} - 1/A^{20} - 2L_1/A^{22} + L_1/A^{26}$} \\ \textbf{Mock:} {\scriptsize $-w^{6} - w^{5} + 2w^{4} + 2w^{3} - w^{2} - 2w + 1 + 2/w - 1/w^{2} - 2/w^{3} + w^{-4} + w^{-5}$} \\ \textbf{Affine:} {\scriptsize $-t + 2 - 1/t$} \\ \textbf{Yamada:} {\scriptsize $-A^{26} + 2A^{24} - A^{22} + 3A^{21} + 2A^{20} - A^{19} + A^{18} + 3A^{17} - A^{16} + 2A^{14} - A^{13} - 2A^{12} + A^{11} - A^{9} - 2A^{8} + 2A^{7} - 3A^{5} + 3A^{4} + A - 1$}
\end{minipage}

\noindent{\color{gray!40}\rule{\textwidth}{0.4pt}}
\vspace{0.9\baselineskip}
\noindent \begin{minipage}[t]{0.25\textwidth}
\vspace{0pt}
\centering
\includegraphics[page=140,width=\linewidth]{knotoids.pdf}
\end{minipage}
\hfill
\begin{minipage}[t]{0.73\textwidth}
\vspace{0pt}
\raggedright
\textbf{Name:} {\large{$\mathbf{K7_{21}}$}} (chiral, non-rotatable$^{*}$) \\ \textbf{PD:} {\scriptsize\texttt{[0],[0,1,2,3],[1,4,5,2],[3,6,7,8],[4,9,6,5],[10,11,8,7],[11,12,13,9],[14,13,12,10],[14]}} \\ \textbf{EM:} {\scriptsize\texttt{(B0, A0C0C3D0, B1E0E3B2, B3E2F3F2, C1G3D1C2, H3G0D3D2, F1H2H1E1, I0G2G1F0, H0)}} \\ \textbf{Kauffman bracket:} {\scriptsize $A^{28} - A^{24} + 3A^{20} - 4A^{16} - A^{14} + 4A^{12} + 2A^{10} - 3A^{8} - 2A^{6} + A^{4} + A^{2}$} \\ \textbf{Arrow:} {\scriptsize $A^{4} - 1 + 3/A^{4} - 4/A^{8} - L_1/A^{10} + 4/A^{12} + 2L_1/A^{14} - 3/A^{16} - 2L_1/A^{18} + A^{-20} + L_1/A^{22}$} \\ \textbf{Mock:} {\scriptsize $-2w^{4} - w^{3} + 5w^{2} + 2w - 5 - 2/w + 4/w^{2} + w^{-3} - 1/w^{4}$} \\ \textbf{Affine:} {\scriptsize $-t + 2 - 1/t$} \\ \textbf{Yamada:} {\scriptsize $-A^{26} + A^{25} - 4A^{23} + A^{22} + 3A^{21} - 7A^{20} + 6A^{18} - 5A^{17} - A^{16} + 4A^{15} - 2A^{14} - A^{13} - A^{12} + 4A^{11} - 2A^{10} - 5A^{9} + 7A^{8} - 2A^{7} - 5A^{6} + 5A^{5} + A^{4} - 3A^{3} + A^{2} + A - 1$}
\end{minipage}

\noindent{\color{gray!40}\rule{\textwidth}{0.4pt}}
\vspace{0.9\baselineskip}
\noindent \begin{minipage}[t]{0.25\textwidth}
\vspace{0pt}
\centering
\includegraphics[page=141,width=\linewidth]{knotoids.pdf}
\end{minipage}
\hfill
\begin{minipage}[t]{0.73\textwidth}
\vspace{0pt}
\raggedright
\textbf{Name:} {\large{$\mathbf{K7_{22}}$}} (chiral, non-rotatable$^{*}$) \\ \textbf{PD:} {\scriptsize\texttt{[0],[0,1,2,3],[1,4,5,2],[3,5,6,7],[4,8,9,6],[10,11,12,7],[8,12,13,9],[14,13,11,10],[14]}} \\ \textbf{EM:} {\scriptsize\texttt{(B0, A0C0C3D0, B1E0D1B2, B3C2E3F3, C1G0G3D2, H3H2G1D3, E1F2H1E2, I0G2F1F0, H0)}} \\ \textbf{Kauffman bracket:} {\scriptsize $-A^{26} - A^{24} + 2A^{22} + 3A^{20} - 2A^{18} - 4A^{16} + A^{14} + 5A^{12} - 3A^{8} + 2A^{4} - 1$} \\ \textbf{Arrow:} {\scriptsize $-A^{14}L_1 - A^{12} + 2A^{10}L_1 + 3A^{8} - 2A^{6}L_1 - 4A^{4} + A^{2}L_1 + 5 - 3/A^{4} + 2/A^{8} - 1/A^{12}$} \\ \textbf{Mock:} {\scriptsize $w^{4} - w^{3} - 4w^{2} + 2w + 7 - 2/w - 5/w^{2} + w^{-3} + 2/w^{4}$} \\ \textbf{Affine:} {\scriptsize $-t + 2 - 1/t$} \\ \textbf{Yamada:} {\scriptsize $A^{26} - A^{25} - 2A^{24} + 3A^{23} + A^{22} - 6A^{21} + 2A^{20} + 5A^{19} - 7A^{18} + 7A^{16} - 4A^{15} - A^{14} + 5A^{13} + A^{12} - A^{10} + 7A^{9} - 2A^{8} - 6A^{7} + 8A^{6} - A^{5} - 6A^{4} + 4A^{3} - 2A + 1$}
\end{minipage}

\noindent{\color{gray!40}\rule{\textwidth}{0.4pt}}
\vspace{0.9\baselineskip}
\noindent \begin{minipage}[t]{0.25\textwidth}
\vspace{0pt}
\centering
\includegraphics[page=142,width=\linewidth]{knotoids.pdf}
\end{minipage}
\hfill
\begin{minipage}[t]{0.73\textwidth}
\vspace{0pt}
\raggedright
\textbf{Name:} {\large{$\mathbf{K7_{23}}$}} (chiral, non-rotatable$^{*}$) \\ \textbf{PD:} {\scriptsize\texttt{[0],[0,1,2,3],[1,4,5,2],[3,6,7,8],[4,9,10,5],[6,10,11,7],[8,12,13,14],[14,13,11,9],[12]}} \\ \textbf{EM:} {\scriptsize\texttt{(B0, A0C0C3D0, B1E0E3B2, B3F0F3G0, C1H3F1C2, D1E2H2D2, D3I0H1H0, G3G2F2E1, G1)}} \\ \textbf{Kauffman bracket:} {\scriptsize $-A^{20} + A^{16} + A^{14} - A^{12} - A^{10} + A^{8} + 2A^{6} - A^{2}$} \\ \textbf{Arrow:} {\scriptsize $-A^{14}L_1 + A^{10}L_1 + A^{8} - A^{6}L_1 - A^{4} + A^{2}L_1 + 2 - 1/A^{4}$} \\ \textbf{Mock:} {\scriptsize $-w^{3} - w^{2} + w + 2 - 1/w - 1/w^{2} + w^{-3} + w^{-4}$} \\ \textbf{Affine:} {\scriptsize $-2t + 4 - 2/t$} \\ \textbf{Yamada:} {\scriptsize $-A^{21} + A^{19} - A^{17} + A^{15} - A^{14} + 2A^{12} + 2A^{9} + A^{7} + 2A^{5} - A^{3} + 2A^{2} - 1$}
\end{minipage}

\noindent{\color{gray!40}\rule{\textwidth}{0.4pt}}
\vspace{0.9\baselineskip}
\noindent \begin{minipage}[t]{0.25\textwidth}
\vspace{0pt}
\centering
\includegraphics[page=143,width=\linewidth]{knotoids.pdf}
\end{minipage}
\hfill
\begin{minipage}[t]{0.73\textwidth}
\vspace{0pt}
\raggedright
\textbf{Name:} {\large{$\mathbf{K7_{24}}$}} (chiral, non-rotatable$^{*}$) \\ \textbf{PD:} {\scriptsize\texttt{[0],[0,1,2,3],[1,4,5,2],[3,6,7,8],[4,9,10,5],[6,10,11,7],[12,13,14,8],[9,14,13,11],[12]}} \\ \textbf{EM:} {\scriptsize\texttt{(B0, A0C0C3D0, B1E0E3B2, B3F0F3G3, C1H0F1C2, D1E2H3D2, I0H2H1D3, E1G2G1F2, G0)}} \\ \textbf{Kauffman bracket:} {\scriptsize $-A^{27} - A^{25} + A^{23} + 2A^{21} - 2A^{19} - 3A^{17} + A^{15} + 4A^{13} - 3A^{9} + 2A^{5} - A$} \\ \textbf{Arrow:} {\scriptsize $1 + L_1/A^{2} - 1/A^{4} - 2L_1/A^{6} + 2/A^{8} + 3L_1/A^{10} - 1/A^{12} - 4L_1/A^{14} + 3L_1/A^{18} - 2L_1/A^{22} + L_1/A^{26}$} \\ \textbf{Mock:} {\scriptsize $-w^{6} - w^{5} + 2w^{4} + 3w^{3} - w^{2} - 4w + 1 + 4/w - 3/w^{3} + w^{-5}$} \\ \textbf{Affine:} {\scriptsize $0$} \\ \textbf{Yamada:} {\scriptsize $-A^{26} + A^{24} - 2A^{23} - 2A^{22} + 3A^{21} - A^{20} - 4A^{19} + 2A^{18} + 2A^{17} - 5A^{16} + 3A^{14} - 3A^{13} - A^{12} + 3A^{11} - 2A^{9} - A^{8} + 4A^{7} - 2A^{6} - 3A^{5} + 6A^{4} - 2A^{3} - 2A^{2} + 2A - 1$}
\end{minipage}

\noindent{\color{gray!40}\rule{\textwidth}{0.4pt}}
\vspace{0.9\baselineskip}
\noindent \begin{minipage}[t]{0.25\textwidth}
\vspace{0pt}
\centering
\includegraphics[page=144,width=\linewidth]{knotoids.pdf}
\end{minipage}
\hfill
\begin{minipage}[t]{0.73\textwidth}
\vspace{0pt}
\raggedright
\textbf{Name:} {\large{$\mathbf{K7_{25}}$}} (chiral, non-rotatable$^{*}$) \\ \textbf{PD:} {\scriptsize\texttt{[0],[0,1,2,3],[1,4,5,2],[3,5,6,7],[4,8,9,10],[10,11,7,6],[8,11,12,13],[13,12,14,9],[14]}} \\ \textbf{EM:} {\scriptsize\texttt{(B0, A0C0C3D0, B1E0D1B2, B3C2F3F2, C1G0H3F0, E3G1D3D2, E1F1H1H0, G3G2I0E2, H2)}} \\ \textbf{Kauffman bracket:} {\scriptsize $A^{28} + A^{26} - 2A^{24} - 2A^{22} + 3A^{20} + 4A^{18} - 3A^{16} - 5A^{14} + A^{12} + 4A^{10} - 2A^{6} + A^{2}$} \\ \textbf{Arrow:} {\scriptsize $A^{10}L_1 + A^{8} - 2A^{6}L_1 - 2A^{4} + 3A^{2}L_1 + 4 - 3L_1/A^{2} - 5/A^{4} + L_1/A^{6} + 4/A^{8} - 2/A^{12} + A^{-16}$} \\ \textbf{Mock:} {\scriptsize $-w^{3} - 3w^{2} + 4w + 9 - 4/w - 6/w^{2} + w^{-3} + w^{-4}$} \\ \textbf{Affine:} {\scriptsize $t - 2 + 1/t$} \\ \textbf{Yamada:} {\scriptsize $-A^{26} + A^{25} - 4A^{23} + 4A^{22} + 3A^{21} - 10A^{20} + 5A^{19} + 7A^{18} - 10A^{17} + 2A^{16} + 5A^{15} - 5A^{14} - 2A^{13} + 4A^{11} - 5A^{10} - 4A^{9} + 10A^{8} - 6A^{7} - 6A^{6} + 10A^{5} - 2A^{4} - 5A^{3} + 4A^{2} + A - 2$}
\end{minipage}

\noindent{\color{gray!40}\rule{\textwidth}{0.4pt}}
\vspace{0.9\baselineskip}
\noindent \begin{minipage}[t]{0.25\textwidth}
\vspace{0pt}
\centering
\includegraphics[page=145,width=\linewidth]{knotoids.pdf}
\end{minipage}
\hfill
\begin{minipage}[t]{0.73\textwidth}
\vspace{0pt}
\raggedright
\textbf{Name:} {\large{$\mathbf{K7_{26}}$}} (achiral, non-rotatable$^{*}$, Possible duplicate [7\_28]) \\ \textbf{PD:} {\scriptsize\texttt{[0],[0,1,2,3],[1,4,5,2],[3,5,6,7],[4,8,9,10],[10,9,11,6],[7,11,12,13],[13,12,14,8],[14]}} \\ \textbf{EM:} {\scriptsize\texttt{(B0, A0C0C3D0, B1E0D1B2, B3C2F3G0, C1H3F1F0, E3E2G1D2, D3F2H1H0, G3G2I0E1, H2)}} \\ \textbf{Kauffman bracket:} {\scriptsize $A^{27} - 2A^{23} + 4A^{19} - 5A^{15} - 2A^{13} + 4A^{11} + 2A^{9} - 3A^{7} - 2A^{5} + A^{3} + A$} \\ \textbf{Arrow:} {\scriptsize $-A^{18}L_1 + 2A^{14}L_1 - 4A^{10}L_1 + 5A^{6}L_1 + 2A^{4} - 4A^{2}L_1 - 2 + 3L_1/A^{2} + 2/A^{4} - L_1/A^{6} - 1/A^{8}$} \\ \textbf{Mock:} {\scriptsize $-w^{5} + 3w^{3} - 6w - 1 + 6/w + 3/w^{2} - 3/w^{3} - 2/w^{4} + w^{-5} + w^{-6}$} \\ \textbf{Affine:} {\scriptsize $-2t + 4 - 2/t$} \\ \textbf{Yamada:} {\scriptsize $A^{26} - 3A^{25} + 6A^{23} - 6A^{22} - 3A^{21} + 9A^{20} - 6A^{19} - 4A^{18} + 7A^{17} - A^{16} - A^{15} + 6A^{13} - 5A^{11} + 8A^{10} + 2A^{9} - 8A^{8} + 5A^{7} + 4A^{6} - 6A^{5} + 3A^{3} - 2A^{2} - A + 1$}
\end{minipage}

\noindent{\color{gray!40}\rule{\textwidth}{0.4pt}}
\vspace{0.9\baselineskip}
\noindent \begin{minipage}[t]{0.25\textwidth}
\vspace{0pt}
\centering
\includegraphics[page=146,width=\linewidth]{knotoids.pdf}
\end{minipage}
\hfill
\begin{minipage}[t]{0.73\textwidth}
\vspace{0pt}
\raggedright
\textbf{Name:} {\large{$\mathbf{K7_{27}}$}} (chiral, non-rotatable$^{*}$) \\ \textbf{PD:} {\scriptsize\texttt{[0],[0,1,2,3],[1,4,5,2],[3,5,6,7],[4,8,9,10],[10,9,11,6],[11,12,13,7],[8,13,12,14],[14]}} \\ \textbf{EM:} {\scriptsize\texttt{(B0, A0C0C3D0, B1E0D1B2, B3C2F3G3, C1H0F1F0, E3E2G0D2, F2H2H1D3, E1G2G1I0, H3)}} \\ \textbf{Kauffman bracket:} {\scriptsize $2A^{20} + A^{18} - 2A^{16} - 2A^{14} + 2A^{12} + 2A^{10} - A^{8} - 2A^{6} + A^{2}$} \\ \textbf{Arrow:} {\scriptsize $2/A^{4} + L_1/A^{6} - 2/A^{8} - 2L_1/A^{10} + 2/A^{12} + 2L_1/A^{14} - 1/A^{16} - 2L_1/A^{18} + L_1/A^{22}$} \\ \textbf{Mock:} {\scriptsize $-w^{4} - w^{3} + 3w^{2} + 3w - 2 - 3/w + w^{-2} + w^{-3}$} \\ \textbf{Affine:} {\scriptsize $0$} \\ \textbf{Yamada:} {\scriptsize $A^{22} + A^{21} - 2A^{20} + 3A^{19} + 3A^{18} - 3A^{17} + 2A^{16} + 3A^{15} - 2A^{14} + A^{13} + A^{11} - 2A^{10} - A^{9} + 3A^{8} - 3A^{7} - A^{6} + 3A^{5} - A^{4} - 2A^{3} + 2A^{2} + A - 1$}
\end{minipage}

\noindent{\color{gray!40}\rule{\textwidth}{0.4pt}}
\vspace{0.9\baselineskip}
\noindent \begin{minipage}[t]{0.25\textwidth}
\vspace{0pt}
\centering
\includegraphics[page=147,width=\linewidth]{knotoids.pdf}
\end{minipage}
\hfill
\begin{minipage}[t]{0.73\textwidth}
\vspace{0pt}
\raggedright
\textbf{Name:} {\large{$\mathbf{K7_{28}}$}} (achiral, non-rotatable$^{*}$, Possible duplicate [7\_26]) \\ \textbf{PD:} {\scriptsize\texttt{[0],[0,1,2,3],[1,4,5,2],[3,6,7,8],[4,9,10,11],[11,12,6,5],[12,13,8,7],[9,13,14,10],[14]}} \\ \textbf{EM:} {\scriptsize\texttt{(B0, A0C0C3D0, B1E0F3B2, B3F2G3G2, C1H0H3F0, E3G0D1C2, F1H1D3D2, E1G1I0E2, H2)}} \\ \textbf{Kauffman bracket:} {\scriptsize $-A^{26} - A^{24} + 2A^{22} + 3A^{20} - 2A^{18} - 4A^{16} + 2A^{14} + 5A^{12} - 4A^{8} + 2A^{4} - 1$} \\ \textbf{Arrow:} {\scriptsize $-A^{8} - A^{6}L_1 + 2A^{4} + 3A^{2}L_1 - 2 - 4L_1/A^{2} + 2/A^{4} + 5L_1/A^{6} - 4L_1/A^{10} + 2L_1/A^{14} - L_1/A^{18}$} \\ \textbf{Mock:} {\scriptsize $w^{6} + w^{5} - 2w^{4} - 3w^{3} + 3w^{2} + 6w - 1 - 6/w + 3/w^{3} - 1/w^{5}$} \\ \textbf{Affine:} {\scriptsize $2t - 4 + 2/t$} \\ \textbf{Yamada:} {\scriptsize $A^{26} - A^{25} - 2A^{24} + 3A^{23} - 6A^{21} + 4A^{20} + 5A^{19} - 8A^{18} + 2A^{17} + 8A^{16} - 5A^{15} + 6A^{13} - A^{11} - A^{10} + 7A^{9} - 4A^{8} - 6A^{7} + 9A^{6} - 3A^{5} - 6A^{4} + 6A^{3} - 3A + 1$}
\end{minipage}

\noindent{\color{gray!40}\rule{\textwidth}{0.4pt}}
\vspace{0.9\baselineskip}
\noindent \begin{minipage}[t]{0.25\textwidth}
\vspace{0pt}
\centering
\includegraphics[page=148,width=\linewidth]{knotoids.pdf}
\end{minipage}
\hfill
\begin{minipage}[t]{0.73\textwidth}
\vspace{0pt}
\raggedright
\textbf{Name:} {\large{$\mathbf{K7_{29}}$}} (chiral, non-rotatable$^{*}$) \\ \textbf{PD:} {\scriptsize\texttt{[0],[0,1,2,3],[1,4,5,2],[3,5,6,7],[4,8,9,6],[7,9,10,11],[11,12,13,8],[14,13,12,10],[14]}} \\ \textbf{EM:} {\scriptsize\texttt{(B0, A0C0C3D0, B1E0D1B2, B3C2E3F0, C1G3F1D2, D3E2H3G0, F3H2H1E1, I0G2G1F2, H0)}} \\ \textbf{Kauffman bracket:} {\scriptsize $A^{28} - 2A^{24} + 4A^{20} + A^{18} - 5A^{16} - 2A^{14} + 5A^{12} + 3A^{10} - 3A^{8} - 3A^{6} + A^{4} + A^{2}$} \\ \textbf{Arrow:} {\scriptsize $A^{4} - 2 + 4/A^{4} + L_1/A^{6} - 5/A^{8} - 2L_1/A^{10} + 5/A^{12} + 3L_1/A^{14} - 3/A^{16} - 3L_1/A^{18} + A^{-20} + L_1/A^{22}$} \\ \textbf{Mock:} {\scriptsize $-2w^{4} - w^{3} + 7w^{2} + 4w - 7 - 4/w + 4/w^{2} + w^{-3} - 1/w^{4}$} \\ \textbf{Affine:} {\scriptsize $t - 2 + 1/t$} \\ \textbf{Yamada:} {\scriptsize $-2A^{26} + A^{25} + 3A^{24} - 6A^{23} + 9A^{21} - 10A^{20} - 4A^{19} + 11A^{18} - 7A^{17} - 5A^{16} + 6A^{15} - A^{14} - 2A^{13} - 3A^{12} + 7A^{11} - 10A^{9} + 10A^{8} + 2A^{7} - 11A^{6} + 6A^{5} + 4A^{4} - 6A^{3} + A^{2} + 2A - 1$}
\end{minipage}

\noindent{\color{gray!40}\rule{\textwidth}{0.4pt}}
\vspace{0.9\baselineskip}
\noindent \begin{minipage}[t]{0.25\textwidth}
\vspace{0pt}
\centering
\includegraphics[page=149,width=\linewidth]{knotoids.pdf}
\end{minipage}
\hfill
\begin{minipage}[t]{0.73\textwidth}
\vspace{0pt}
\raggedright
\textbf{Name:} {\large{$\mathbf{K7_{30}}$}} (chiral, non-rotatable$^{*}$) \\ \textbf{PD:} {\scriptsize\texttt{[0],[0,1,2,3],[1,4,5,2],[3,6,7,4],[8,9,6,5],[10,11,12,7],[13,14,9,8],[14,13,11,10],[12]}} \\ \textbf{EM:} {\scriptsize\texttt{(B0, A0C0C3D0, B1D3E3B2, B3E2F3C1, G3G2D1C2, H3H2I0D2, H1H0E1E0, G1G0F1F0, F2)}} \\ \textbf{Kauffman bracket:} {\scriptsize $A^{27} - 2A^{23} - A^{21} + 3A^{19} + 2A^{17} - 3A^{15} - 2A^{13} + 2A^{11} + 2A^{9} - 2A^{7} - 2A^{5} + A$} \\ \textbf{Arrow:} {\scriptsize $-A^{24} + 2A^{20} + A^{18}L_1 - 3A^{16} - 2A^{14}L_1 + 3A^{12} + 2A^{10}L_1 - 2A^{8} - 2A^{6}L_1 + 2A^{4} + 2A^{2}L_1 - L_1/A^{2}$} \\ \textbf{Mock:} {\scriptsize $w^{6} + w^{5} - 2w^{4} - 2w^{3} + 3w^{2} + 2w - 3 - 2/w + 3/w^{2} + 2/w^{3} - 1/w^{4} - 1/w^{5}$} \\ \textbf{Affine:} {\scriptsize $t - 2 + 1/t$} \\ \textbf{Yamada:} {\scriptsize $A^{26} - A^{25} - A^{24} + 2A^{23} - 3A^{22} + 5A^{20} - 5A^{19} + 5A^{17} - 2A^{16} + A^{14} + 2A^{13} - 2A^{12} - 3A^{11} + 3A^{10} - 2A^{9} - 5A^{8} + 3A^{7} - A^{6} - 5A^{5} + A^{4} + 2A^{3} - 3A^{2} + 2$}
\end{minipage}

\noindent{\color{gray!40}\rule{\textwidth}{0.4pt}}
\vspace{0.9\baselineskip}
\noindent \begin{minipage}[t]{0.25\textwidth}
\vspace{0pt}
\centering
\includegraphics[page=150,width=\linewidth]{knotoids.pdf}
\end{minipage}
\hfill
\begin{minipage}[t]{0.73\textwidth}
\vspace{0pt}
\raggedright
\textbf{Name:} {\large{$\mathbf{K7_{31}}$}} (chiral, non-rotatable$^{*}$) \\ \textbf{PD:} {\scriptsize\texttt{[0],[0,1,2,3],[1,4,5,2],[3,6,7,4],[8,9,10,5],[6,10,11,7],[12,13,9,8],[13,12,14,11],[14]}} \\ \textbf{EM:} {\scriptsize\texttt{(B0, A0C0C3D0, B1D3E3B2, B3F0F3C1, G3G2F1C2, D1E2H3D2, H1H0E1E0, G1G0I0F2, H2)}} \\ \textbf{Kauffman bracket:} {\scriptsize $A^{28} - 2A^{24} - A^{22} + 4A^{20} + 3A^{18} - 3A^{16} - 4A^{14} + 2A^{12} + 3A^{10} - A^{8} - 2A^{6} + A^{2}$} \\ \textbf{Arrow:} {\scriptsize $A^{4} - 2 - L_1/A^{2} + 4/A^{4} + 3L_1/A^{6} - 3/A^{8} - 4L_1/A^{10} + 2/A^{12} + 3L_1/A^{14} - 1/A^{16} - 2L_1/A^{18} + L_1/A^{22}$} \\ \textbf{Mock:} {\scriptsize $-2w^{4} - 2w^{3} + 5w^{2} + 5w - 4 - 5/w + 2/w^{2} + 2/w^{3}$} \\ \textbf{Affine:} {\scriptsize $-t + 2 - 1/t$} \\ \textbf{Yamada:} {\scriptsize $-A^{26} + A^{25} + 2A^{24} - 4A^{23} + A^{22} + 8A^{21} - 6A^{20} - A^{19} + 10A^{18} - 5A^{17} - 3A^{16} + 6A^{15} - A^{14} - A^{13} - 2A^{12} + 5A^{11} - A^{10} - 8A^{9} + 8A^{8} + A^{7} - 8A^{6} + 5A^{5} + 3A^{4} - 5A^{3} + A^{2} + 2A - 1$}
\end{minipage}

\noindent{\color{gray!40}\rule{\textwidth}{0.4pt}}
\vspace{0.9\baselineskip}
\noindent \begin{minipage}[t]{0.25\textwidth}
\vspace{0pt}
\centering
\includegraphics[page=151,width=\linewidth]{knotoids.pdf}
\end{minipage}
\hfill
\begin{minipage}[t]{0.73\textwidth}
\vspace{0pt}
\raggedright
\textbf{Name:} {\large{$\mathbf{K7_{32}}$}} (chiral, non-rotatable$^{*}$) \\ \textbf{PD:} {\scriptsize\texttt{[0],[0,1,2,3],[1,4,5,2],[3,6,7,8],[4,9,6,5],[10,11,8,7],[11,12,13,9],[10,13,12,14],[14]}} \\ \textbf{EM:} {\scriptsize\texttt{(B0, A0C0C3D0, B1E0E3B2, B3E2F3F2, C1G3D1C2, H0G0D3D2, F1H2H1E1, F0G2G1I0, H3)}} \\ \textbf{Kauffman bracket:} {\scriptsize $-A^{26} - A^{24} + A^{22} + 2A^{20} - A^{18} - 3A^{16} + A^{14} + 4A^{12} + A^{10} - 3A^{8} - A^{6} + A^{4} + A^{2}$} \\ \textbf{Arrow:} {\scriptsize $-A^{32} - A^{30}L_1 + A^{28} + 2A^{26}L_1 - A^{24} - 3A^{22}L_1 + A^{20} + 4A^{18}L_1 + A^{16} - 3A^{14}L_1 - A^{12} + A^{10}L_1 + A^{8}$} \\ \textbf{Mock:} {\scriptsize $-2w^{3} + 2w^{2} + 5w - 2 - 5/w - 1/w^{2} + 2/w^{3} + 2/w^{4}$} \\ \textbf{Affine:} {\scriptsize $-t + 2 - 1/t$} \\ \textbf{Yamada:} {\scriptsize $-A^{27} + A^{25} - 2A^{24} - 2A^{23} + 5A^{22} - 4A^{20} + 5A^{19} + 3A^{18} - 5A^{17} + 3A^{16} + 4A^{15} - 3A^{14} + 2A^{11} - 6A^{10} - A^{9} + 5A^{8} - 7A^{7} - 2A^{6} + 4A^{5} - 3A^{4} - 2A^{3} + A^{2} - 1$}
\end{minipage}

\noindent{\color{gray!40}\rule{\textwidth}{0.4pt}}
\vspace{0.9\baselineskip}
\noindent \begin{minipage}[t]{0.25\textwidth}
\vspace{0pt}
\centering
\includegraphics[page=152,width=\linewidth]{knotoids.pdf}
\end{minipage}
\hfill
\begin{minipage}[t]{0.73\textwidth}
\vspace{0pt}
\raggedright
\textbf{Name:} {\large{$\mathbf{K7_{33}}$}} (chiral, non-rotatable$^{*}$) \\ \textbf{PD:} {\scriptsize\texttt{[0],[0,1,2,3],[1,4,5,2],[6,7,4,3],[5,8,9,6],[7,10,11,12],[8,13,14,9],[10,14,13,11],[12]}} \\ \textbf{EM:} {\scriptsize\texttt{(B0, A0C0C3D3, B1D2E0B2, E3F0C1B3, C2G0G3D0, D1H0H3I0, E1H2H1E2, F1G2G1F2, F3)}} \\ \textbf{Kauffman bracket:} {\scriptsize $A^{24} + A^{22} - 2A^{18} + 2A^{14} + A^{12} - 2A^{10} - A^{8} + 2A^{6} + A^{4} - A^{2} - 1$} \\ \textbf{Arrow:} {\scriptsize $A^{-12} + L_1/A^{14} - 2L_1/A^{18} + 2L_1/A^{22} + A^{-24} - 2L_1/A^{26} - 1/A^{28} + 2L_1/A^{30} + A^{-32} - L_1/A^{34} - 1/A^{36}$} \\ \textbf{Mock:} {\scriptsize $w^{6} + w^{5} - 2w^{3} - w^{2} + 2w + 1 - 2/w - 1/w^{2} + 2/w^{3} + w^{-4} - 1/w^{5}$} \\ \textbf{Affine:} {\scriptsize $t - 2 + 1/t$} \\ \textbf{Yamada:} {\scriptsize $A^{27} + A^{26} + 3A^{23} + A^{22} - A^{21} + 2A^{20} + 3A^{19} - A^{18} + A^{17} + 3A^{16} - 2A^{15} - A^{14} + A^{13} - 2A^{12} - A^{11} - 2A^{10} + 2A^{9} - 2A^{8} - 3A^{7} + 3A^{6} - A^{5} - A^{4} + A^{3} + 1$}
\end{minipage}

\noindent{\color{gray!40}\rule{\textwidth}{0.4pt}}
\vspace{0.9\baselineskip}
\noindent \begin{minipage}[t]{0.25\textwidth}
\vspace{0pt}
\centering
\includegraphics[page=153,width=\linewidth]{knotoids.pdf}
\end{minipage}
\hfill
\begin{minipage}[t]{0.73\textwidth}
\vspace{0pt}
\raggedright
\textbf{Name:} {\large{$\mathbf{K7_{34}}$}} (achiral, non-rotatable$^{*}$, Possible duplicate [7\_59]) \\ \textbf{PD:} {\scriptsize\texttt{[0],[0,1,2,3],[1,4,5,2],[3,6,7,8],[4,9,10,5],[6,10,11,12],[12,13,8,7],[9,13,14,11],[14]}} \\ \textbf{EM:} {\scriptsize\texttt{(B0, A0C0C3D0, B1E0E3B2, B3F0G3G2, C1H0F1C2, D1E2H3G0, F3H1D3D2, E1G1I0F2, H2)}} \\ \textbf{Kauffman bracket:} {\scriptsize $-A^{26} - A^{24} + 2A^{22} + 2A^{20} - 2A^{18} - 3A^{16} + 3A^{14} + 5A^{12} - A^{10} - 4A^{8} + 2A^{4} - 1$} \\ \textbf{Arrow:} {\scriptsize $-A^{8} - A^{6}L_1 + 2A^{4} + 2A^{2}L_1 - 2 - 3L_1/A^{2} + 3/A^{4} + 5L_1/A^{6} - 1/A^{8} - 4L_1/A^{10} + 2L_1/A^{14} - L_1/A^{18}$} \\ \textbf{Mock:} {\scriptsize $w^{6} + w^{5} - 2w^{4} - 3w^{3} + 3w^{2} + 5w - 2 - 5/w + w^{-2} + 3/w^{3} - 1/w^{5}$} \\ \textbf{Affine:} {\scriptsize $t - 2 + 1/t$} \\ \textbf{Yamada:} {\scriptsize $A^{26} - A^{25} - A^{24} + 3A^{23} - 2A^{22} - 5A^{21} + 5A^{20} + A^{19} - 8A^{18} + 5A^{17} + 6A^{16} - 5A^{15} + 3A^{14} + 6A^{13} - A^{12} + A^{10} + 6A^{9} - 6A^{8} - 3A^{7} + 8A^{6} - 7A^{5} - 4A^{4} + 6A^{3} - 2A^{2} - 2A + 2$}
\end{minipage}

\noindent{\color{gray!40}\rule{\textwidth}{0.4pt}}
\vspace{0.9\baselineskip}
\noindent \begin{minipage}[t]{0.25\textwidth}
\vspace{0pt}
\centering
\includegraphics[page=154,width=\linewidth]{knotoids.pdf}
\end{minipage}
\hfill
\begin{minipage}[t]{0.73\textwidth}
\vspace{0pt}
\raggedright
\textbf{Name:} {\large{$\mathbf{K7_{35}}$}} (achiral, non-rotatable$^{*}$, Possible duplicate [7\_60]) \\ \textbf{PD:} {\scriptsize\texttt{[0],[0,1,2,3],[1,4,5,2],[3,6,7,8],[4,9,10,5],[6,10,11,12],[12,13,8,7],[13,14,11,9],[14]}} \\ \textbf{EM:} {\scriptsize\texttt{(B0, A0C0C3D0, B1E0E3B2, B3F0G3G2, C1H3F1C2, D1E2H2G0, F3H0D3D2, G1I0F2E1, H1)}} \\ \textbf{Kauffman bracket:} {\scriptsize $-A^{20} - A^{18} + A^{16} + 2A^{14} - 2A^{10} + 3A^{6} + A^{4} - A^{2} - 1$} \\ \textbf{Arrow:} {\scriptsize $-A^{14}L_1 - A^{12} + A^{10}L_1 + 2A^{8} - 2A^{4} + 3 + L_1/A^{2} - 1/A^{4} - L_1/A^{6}$} \\ \textbf{Mock:} {\scriptsize $w^{4} - 2w^{2} + 3 - 2/w^{2} + w^{-4}$} \\ \textbf{Affine:} {\scriptsize $0$} \\ \textbf{Yamada:} {\scriptsize $A^{24} - A^{22} + A^{20} - 2A^{19} - 3A^{18} + 2A^{17} - 3A^{15} + 3A^{14} + 2A^{13} - A^{12} + 3A^{11} + 2A^{10} + 3A^{9} - A^{8} + 2A^{7} + 2A^{6} - 4A^{5} + 2A^{3} - 2A^{2} - A + 1$}
\end{minipage}

\noindent{\color{gray!40}\rule{\textwidth}{0.4pt}}
\vspace{0.9\baselineskip}
\noindent \begin{minipage}[t]{0.25\textwidth}
\vspace{0pt}
\centering
\includegraphics[page=155,width=\linewidth]{knotoids.pdf}
\end{minipage}
\hfill
\begin{minipage}[t]{0.73\textwidth}
\vspace{0pt}
\raggedright
\textbf{Name:} {\large{$\mathbf{K7_{36}}$}} (chiral, non-rotatable$^{*}$) \\ \textbf{PD:} {\scriptsize\texttt{[0],[0,1,2,3],[1,4,5,2],[3,6,7,4],[8,9,6,5],[9,10,11,7],[12,13,14,8],[10,14,13,11],[12]}} \\ \textbf{EM:} {\scriptsize\texttt{(B0, A0C0C3D0, B1D3E3B2, B3E2F3C1, G3F0D1C2, E1H0H3D2, I0H2H1E0, F1G2G1F2, G0)}} \\ \textbf{Kauffman bracket:} {\scriptsize $-A^{26} + 2A^{22} + A^{20} - 3A^{18} - 2A^{16} + 4A^{14} + 3A^{12} - 3A^{10} - 4A^{8} + 2A^{6} + 3A^{4} - 1$} \\ \textbf{Arrow:} {\scriptsize $-A^{20} + 2A^{16} + A^{14}L_1 - 3A^{12} - 2A^{10}L_1 + 4A^{8} + 3A^{6}L_1 - 3A^{4} - 4A^{2}L_1 + 2 + 3L_1/A^{2} - L_1/A^{6}$} \\ \textbf{Mock:} {\scriptsize $2w^{4} + 2w^{3} - 5w^{2} - 5w + 6 + 5/w - 2/w^{2} - 2/w^{3}$} \\ \textbf{Affine:} {\scriptsize $t - 2 + 1/t$} \\ \textbf{Yamada:} {\scriptsize $A^{26} - A^{25} - 2A^{24} + 4A^{23} - 7A^{21} + 5A^{20} + 4A^{19} - 9A^{18} + 4A^{17} + 7A^{16} - 7A^{15} + 4A^{13} - 2A^{12} - 3A^{11} - 2A^{10} + 6A^{9} - 6A^{8} - 6A^{7} + 10A^{6} - 5A^{5} - 6A^{4} + 7A^{3} - 2A^{2} - 2A + 2$}
\end{minipage}

\noindent{\color{gray!40}\rule{\textwidth}{0.4pt}}
\vspace{0.9\baselineskip}
\noindent \begin{minipage}[t]{0.25\textwidth}
\vspace{0pt}
\centering
\includegraphics[page=156,width=\linewidth]{knotoids.pdf}
\end{minipage}
\hfill
\begin{minipage}[t]{0.73\textwidth}
\vspace{0pt}
\raggedright
\textbf{Name:} {\large{$\mathbf{K7_{37}}$}} (chiral, non-rotatable$^{*}$) \\ \textbf{PD:} {\scriptsize\texttt{[0],[0,1,2,3],[1,4,5,2],[3,6,7,4],[8,9,6,5],[10,11,12,7],[13,10,9,8],[11,13,14,12],[14]}} \\ \textbf{EM:} {\scriptsize\texttt{(B0, A0C0C3D0, B1D3E3B2, B3E2F3C1, G3G2D1C2, G1H0H3D2, H1F0E1E0, F1G0I0F2, H2)}} \\ \textbf{Kauffman bracket:} {\scriptsize $A^{28} - 2A^{24} - A^{22} + 4A^{20} + 2A^{18} - 4A^{16} - 3A^{14} + 4A^{12} + 4A^{10} - 2A^{8} - 3A^{6} + A^{2}$} \\ \textbf{Arrow:} {\scriptsize $A^{4} - 2 - L_1/A^{2} + 4/A^{4} + 2L_1/A^{6} - 4/A^{8} - 3L_1/A^{10} + 4/A^{12} + 4L_1/A^{14} - 2/A^{16} - 3L_1/A^{18} + L_1/A^{22}$} \\ \textbf{Mock:} {\scriptsize $-2w^{4} - 2w^{3} + 6w^{2} + 5w - 6 - 5/w + 3/w^{2} + 2/w^{3}$} \\ \textbf{Affine:} {\scriptsize $-t + 2 - 1/t$} \\ \textbf{Yamada:} {\scriptsize $-A^{26} + A^{25} + A^{24} - 4A^{23} + 4A^{22} + 5A^{21} - 10A^{20} + 5A^{19} + 10A^{18} - 10A^{17} + 2A^{16} + 8A^{15} - 4A^{14} - A^{13} + 5A^{11} - 5A^{10} - 5A^{9} + 11A^{8} - 5A^{7} - 8A^{6} + 10A^{5} - A^{4} - 7A^{3} + 5A^{2} + 2A - 2$}
\end{minipage}

\noindent{\color{gray!40}\rule{\textwidth}{0.4pt}}
\vspace{0.9\baselineskip}
\noindent \begin{minipage}[t]{0.25\textwidth}
\vspace{0pt}
\centering
\includegraphics[page=157,width=\linewidth]{knotoids.pdf}
\end{minipage}
\hfill
\begin{minipage}[t]{0.73\textwidth}
\vspace{0pt}
\raggedright
\textbf{Name:} {\large{$\mathbf{K7_{38}}$}} (chiral, non-rotatable$^{*}$) \\ \textbf{PD:} {\scriptsize\texttt{[0],[0,1,2,3],[1,4,5,2],[3,5,6,7],[4,8,9,10],[10,11,7,6],[8,11,12,13],[9,14,13,12],[14]}} \\ \textbf{EM:} {\scriptsize\texttt{(B0, A0C0C3D0, B1E0D1B2, B3C2F3F2, C1G0H0F0, E3G1D3D2, E1F1H3H2, E2I0G3G2, H1)}} \\ \textbf{Kauffman bracket:} {\scriptsize $-A^{25} - A^{23} + A^{21} + 2A^{19} - 2A^{17} - 3A^{15} + A^{13} + 4A^{11} - 3A^{7} - A^{5} + A^{3} + A$} \\ \textbf{Arrow:} {\scriptsize $A^{-8} + L_1/A^{10} - 1/A^{12} - 2L_1/A^{14} + 2/A^{16} + 3L_1/A^{18} - 1/A^{20} - 4L_1/A^{22} + 3L_1/A^{26} + A^{-28} - L_1/A^{30} - 1/A^{32}$} \\ \textbf{Mock:} {\scriptsize $2w^{4} + 2w^{3} - 2w^{2} - 5w + 5/w + w^{-2} - 2/w^{3}$} \\ \textbf{Affine:} {\scriptsize $t - 2 + 1/t$} \\ \textbf{Yamada:} {\scriptsize $A^{27} + 3A^{24} - 3A^{22} + 7A^{21} + 2A^{20} - 6A^{19} + 5A^{18} + 3A^{17} - 4A^{16} + A^{15} + A^{14} + A^{13} - 5A^{12} + 4A^{10} - 7A^{9} + 6A^{7} - 4A^{6} - 2A^{5} + 4A^{4} - 2A^{2} + 1$}
\end{minipage}

\noindent{\color{gray!40}\rule{\textwidth}{0.4pt}}
\vspace{0.9\baselineskip}
\noindent \begin{minipage}[t]{0.25\textwidth}
\vspace{0pt}
\centering
\includegraphics[page=158,width=\linewidth]{knotoids.pdf}
\end{minipage}
\hfill
\begin{minipage}[t]{0.73\textwidth}
\vspace{0pt}
\raggedright
\textbf{Name:} {\large{$\mathbf{K7_{39}}$}} (chiral, non-rotatable$^{*}$) \\ \textbf{PD:} {\scriptsize\texttt{[0],[0,1,2,3],[1,4,5,2],[3,6,7,8],[4,9,6,5],[7,10,11,12],[8,11,13,14],[14,13,10,9],[12]}} \\ \textbf{EM:} {\scriptsize\texttt{(B0, A0C0C3D0, B1E0E3B2, B3E2F0G0, C1H3D1C2, D2H2G1I0, D3F2H1H0, G3G2F1E1, F3)}} \\ \textbf{Kauffman bracket:} {\scriptsize $-A^{26} - A^{24} + A^{22} + 2A^{20} - 3A^{16} + 3A^{12} + A^{10} - 2A^{8} - A^{6} + A^{4} + A^{2}$} \\ \textbf{Arrow:} {\scriptsize $-A^{32} - A^{30}L_1 + A^{28} + 2A^{26}L_1 - 3A^{22}L_1 + 3A^{18}L_1 + A^{16} - 2A^{14}L_1 - A^{12} + A^{10}L_1 + A^{8}$} \\ \textbf{Mock:} {\scriptsize $-2w^{3} + w^{2} + 4w - 1 - 4/w - 1/w^{2} + 2/w^{3} + 2/w^{4}$} \\ \textbf{Affine:} {\scriptsize $-2t + 4 - 2/t$} \\ \textbf{Yamada:} {\scriptsize $-A^{27} + A^{25} - A^{24} - 2A^{23} + 3A^{22} + 2A^{21} - 3A^{20} + 2A^{19} + 4A^{18} - 4A^{17} + 3A^{15} - 2A^{14} + 3A^{11} - 3A^{10} - 2A^{9} + 4A^{8} - 4A^{7} - 4A^{6} + 2A^{5} - A^{4} - 2A^{3} - 1$}
\end{minipage}

\noindent{\color{gray!40}\rule{\textwidth}{0.4pt}}
\vspace{0.9\baselineskip}
\noindent \begin{minipage}[t]{0.25\textwidth}
\vspace{0pt}
\centering
\includegraphics[page=159,width=\linewidth]{knotoids.pdf}
\end{minipage}
\hfill
\begin{minipage}[t]{0.73\textwidth}
\vspace{0pt}
\raggedright
\textbf{Name:} {\large{$\mathbf{K7_{40}}$}} (chiral, non-rotatable$^{*}$) \\ \textbf{PD:} {\scriptsize\texttt{[0],[0,1,2,3],[1,4,5,2],[3,6,7,8],[4,9,10,5],[6,10,11,7],[12,13,14,8],[14,12,11,9],[13]}} \\ \textbf{EM:} {\scriptsize\texttt{(B0, A0C0C3D0, B1E0E3B2, B3F0F3G3, C1H3F1C2, D1E2H2D2, H1I0H0D3, G2G0F2E1, G1)}} \\ \textbf{Kauffman bracket:} {\scriptsize $A^{18} + A^{16} - A^{12} + A^{8} + A^{6} - A^{4} - A^{2}$} \\ \textbf{Arrow:} {\scriptsize $A^{-12} + L_1/A^{14} - L_1/A^{18} + L_1/A^{22} + A^{-24} - L_1/A^{26} - 1/A^{28}$} \\ \textbf{Mock:} {\scriptsize $w^{6} + w^{5} - w^{3} - w^{2} + 1 + w^{-3} - 1/w^{5}$} \\ \textbf{Affine:} {\scriptsize $2t - 4 + 2/t$} \\ \textbf{Yamada:} {\scriptsize $-A^{22} - A^{21} - A^{20} - A^{19} - 2A^{18} - 2A^{17} - A^{16} - 2A^{14} + A^{12} + A^{10} + 2A^{9} + A^{7} + A^{5} - A^{3} + 2A^{2} - A - 1$}
\end{minipage}

\noindent{\color{gray!40}\rule{\textwidth}{0.4pt}}
\vspace{0.9\baselineskip}
\noindent \begin{minipage}[t]{0.25\textwidth}
\vspace{0pt}
\centering
\includegraphics[page=160,width=\linewidth]{knotoids.pdf}
\end{minipage}
\hfill
\begin{minipage}[t]{0.73\textwidth}
\vspace{0pt}
\raggedright
\textbf{Name:} {\large{$\mathbf{K7_{41}}$}} (chiral, non-rotatable$^{*}$) \\ \textbf{PD:} {\scriptsize\texttt{[0],[0,1,2,3],[1,4,5,2],[3,5,6,7],[4,8,9,10],[10,9,11,6],[11,12,13,7],[13,14,12,8],[14]}} \\ \textbf{EM:} {\scriptsize\texttt{(B0, A0C0C3D0, B1E0D1B2, B3C2F3G3, C1H3F1F0, E3E2G0D2, F2H2H0D3, G2I0G1E1, H1)}} \\ \textbf{Kauffman bracket:} {\scriptsize $A^{19} - A^{17} - 2A^{15} + 2A^{11} - 2A^{7} - A^{5} + A^{3} + A$} \\ \textbf{Arrow:} {\scriptsize $-L_1/A^{2} + A^{-4} + 2L_1/A^{6} - 2L_1/A^{10} + 2L_1/A^{14} + A^{-16} - L_1/A^{18} - 1/A^{20}$} \\ \textbf{Mock:} {\scriptsize $-w^{3} + w^{2} + 3w - 3/w - 1/w^{2} + w^{-3} + w^{-4}$} \\ \textbf{Affine:} {\scriptsize $0$} \\ \textbf{Yamada:} {\scriptsize $-A^{23} - A^{22} + 2A^{21} + A^{20} - A^{19} + 3A^{18} + 2A^{17} - A^{16} + 2A^{15} + A^{14} + 2A^{13} - A^{12} + A^{11} + A^{10} - 3A^{9} + A^{7} - 2A^{6} - A^{5} + A^{4} - A^{2} + 1$}
\end{minipage}

\noindent{\color{gray!40}\rule{\textwidth}{0.4pt}}
\vspace{0.9\baselineskip}
\noindent \begin{minipage}[t]{0.25\textwidth}
\vspace{0pt}
\centering
\includegraphics[page=161,width=\linewidth]{knotoids.pdf}
\end{minipage}
\hfill
\begin{minipage}[t]{0.73\textwidth}
\vspace{0pt}
\raggedright
\textbf{Name:} {\large{$\mathbf{K7_{42}}$}} (chiral, non-rotatable$^{*}$) \\ \textbf{PD:} {\scriptsize\texttt{[0],[0,1,2,3],[1,4,5,2],[3,6,7,8],[4,8,9,5],[6,10,11,7],[12,13,14,9],[10,14,13,11],[12]}} \\ \textbf{EM:} {\scriptsize\texttt{(B0, A0C0C3D0, B1E0E3B2, B3F0F3E1, C1D3G3C2, D1H0H3D2, I0H2H1E2, F1G2G1F2, G0)}} \\ \textbf{Kauffman bracket:} {\scriptsize $-A^{27} + A^{23} + A^{21} - 2A^{19} - 2A^{17} + 2A^{15} + 3A^{13} - A^{11} - 3A^{9} + 2A^{5} - A$} \\ \textbf{Arrow:} {\scriptsize $A^{12} - A^{8} - A^{6}L_1 + 2A^{4} + 2A^{2}L_1 - 2 - 3L_1/A^{2} + A^{-4} + 3L_1/A^{6} - 2L_1/A^{10} + L_1/A^{14}$} \\ \textbf{Mock:} {\scriptsize $-2w^{4} - 2w^{3} + 3w^{2} + 4w - 1 - 4/w + w^{-2} + 2/w^{3}$} \\ \textbf{Affine:} {\scriptsize $-2t + 4 - 2/t$} \\ \textbf{Yamada:} {\scriptsize $A^{25} - 2A^{23} + 2A^{21} - 3A^{20} - 2A^{19} + 4A^{18} - 3A^{16} + 3A^{15} + A^{14} - 5A^{13} + A^{11} - 3A^{10} - A^{9} + A^{8} + 3A^{7} - 4A^{6} + 5A^{4} - 4A^{3} + 2A - 2$}
\end{minipage}

\noindent{\color{gray!40}\rule{\textwidth}{0.4pt}}
\vspace{0.9\baselineskip}
\noindent \begin{minipage}[t]{0.25\textwidth}
\vspace{0pt}
\centering
\includegraphics[page=162,width=\linewidth]{knotoids.pdf}
\end{minipage}
\hfill
\begin{minipage}[t]{0.73\textwidth}
\vspace{0pt}
\raggedright
\textbf{Name:} {\large{$\mathbf{K7_{43}}$}} (chiral, non-rotatable$^{*}$) \\ \textbf{PD:} {\scriptsize\texttt{[0],[0,1,2,3],[1,4,5,2],[3,5,6,7],[4,8,9,6],[10,11,8,7],[12,13,14,9],[13,12,11,10],[14]}} \\ \textbf{EM:} {\scriptsize\texttt{(B0, A0C0C3D0, B1E0D1B2, B3C2E3F3, C1F2G3D2, H3H2E1D3, H1H0I0E2, G1G0F1F0, G2)}} \\ \textbf{Kauffman bracket:} {\scriptsize $A^{28} + A^{26} - 2A^{24} - 2A^{22} + 3A^{20} + 4A^{18} - 2A^{16} - 5A^{14} + A^{12} + 3A^{10} - 2A^{6} + A^{2}$} \\ \textbf{Arrow:} {\scriptsize $A^{4} + A^{2}L_1 - 2 - 2L_1/A^{2} + 3/A^{4} + 4L_1/A^{6} - 2/A^{8} - 5L_1/A^{10} + A^{-12} + 3L_1/A^{14} - 2L_1/A^{18} + L_1/A^{22}$} \\ \textbf{Mock:} {\scriptsize $-2w^{4} - 2w^{3} + 5w^{2} + 7w - 2 - 7/w + 2/w^{3}$} \\ \textbf{Affine:} {\scriptsize $t - 2 + 1/t$} \\ \textbf{Yamada:} {\scriptsize $-A^{26} + A^{25} - 4A^{23} + 3A^{22} + 3A^{21} - 9A^{20} + 2A^{19} + 6A^{18} - 8A^{17} + A^{16} + 6A^{15} - 3A^{14} - A^{13} - A^{12} + 3A^{11} - 5A^{10} - 5A^{9} + 9A^{8} - 4A^{7} - 5A^{6} + 8A^{5} - 4A^{3} + 3A^{2} + A - 2$}
\end{minipage}

\noindent{\color{gray!40}\rule{\textwidth}{0.4pt}}
\vspace{0.9\baselineskip}
\noindent \begin{minipage}[t]{0.25\textwidth}
\vspace{0pt}
\centering
\includegraphics[page=163,width=\linewidth]{knotoids.pdf}
\end{minipage}
\hfill
\begin{minipage}[t]{0.73\textwidth}
\vspace{0pt}
\raggedright
\textbf{Name:} {\large{$\mathbf{K7_{44}}$}} (chiral, non-rotatable$^{*}$) \\ \textbf{PD:} {\scriptsize\texttt{[0],[0,1,2,3],[1,4,5,2],[3,5,6,7],[4,7,8,9],[10,11,12,6],[13,14,9,8],[14,13,11,10],[12]}} \\ \textbf{EM:} {\scriptsize\texttt{(B0, A0C0C3D0, B1E0D1B2, B3C2F3E1, C1D3G3G2, H3H2I0D2, H1H0E3E2, G1G0F1F0, F2)}} \\ \textbf{Kauffman bracket:} {\scriptsize $A^{27} - 2A^{23} - A^{21} + 3A^{19} + 2A^{17} - 3A^{15} - 3A^{13} + 2A^{11} + 2A^{9} - A^{7} - 2A^{5} + A$} \\ \textbf{Arrow:} {\scriptsize $-A^{18}L_1 + 2A^{14}L_1 + A^{12} - 3A^{10}L_1 - 2A^{8} + 3A^{6}L_1 + 3A^{4} - 2A^{2}L_1 - 2 + L_1/A^{2} + 2/A^{4} - 1/A^{8}$} \\ \textbf{Mock:} {\scriptsize $-w^{5} - w^{4} + 2w^{3} + 2w^{2} - 3w - 2 + 3/w + 3/w^{2} - 2/w^{3} - 2/w^{4} + w^{-5} + w^{-6}$} \\ \textbf{Affine:} {\scriptsize $-2t + 4 - 2/t$} \\ \textbf{Yamada:} {\scriptsize $A^{26} - A^{25} + 3A^{23} - 4A^{22} - A^{21} + 6A^{20} - 5A^{19} - 2A^{18} + 5A^{17} - 2A^{16} - A^{15} + 3A^{13} - 2A^{12} - 5A^{11} + 3A^{10} - 2A^{9} - 6A^{8} + 4A^{7} + 2A^{6} - 4A^{5} + A^{4} + 3A^{3} - 2A^{2} - A + 1$}
\end{minipage}

\noindent{\color{gray!40}\rule{\textwidth}{0.4pt}}
\vspace{0.9\baselineskip}
\noindent \begin{minipage}[t]{0.25\textwidth}
\vspace{0pt}
\centering
\includegraphics[page=164,width=\linewidth]{knotoids.pdf}
\end{minipage}
\hfill
\begin{minipage}[t]{0.73\textwidth}
\vspace{0pt}
\raggedright
\textbf{Name:} {\large{$\mathbf{K7_{45}}$}} (chiral, non-rotatable$^{*}$) \\ \textbf{PD:} {\scriptsize\texttt{[0],[0,1,2,3],[1,4,5,2],[3,6,7,8],[4,9,10,5],[6,11,12,7],[8,13,14,9],[14,12,11,10],[13]}} \\ \textbf{EM:} {\scriptsize\texttt{(B0, A0C0C3D0, B1E0E3B2, B3F0F3G0, C1G3H3C2, D1H2H1D2, D3I0H0E1, G2F2F1E2, G1)}} \\ \textbf{Kauffman bracket:} {\scriptsize $-A^{21} - A^{19} + A^{15} - 2A^{11} - A^{9} + 2A^{7} + 2A^{5} - A$} \\ \textbf{Arrow:} {\scriptsize $A^{-12} + L_1/A^{14} - L_1/A^{18} + 2L_1/A^{22} + A^{-24} - 2L_1/A^{26} - 2/A^{28} + A^{-32}$} \\ \textbf{Mock:} {\scriptsize $w^{6} + w^{5} - w^{3} - w^{2} + w + 2 - 1/w - 1/w^{2} + w^{-3} - 1/w^{5}$} \\ \textbf{Affine:} {\scriptsize $3t - 6 + 3/t$} \\ \textbf{Yamada:} {\scriptsize $A^{23} + A^{22} + A^{21} + A^{20} + 3A^{19} + 3A^{18} + 3A^{15} - 2A^{14} - 2A^{13} + A^{12} - 2A^{11} - 3A^{10} - A^{8} - A^{6} + 2A^{5} + 2A^{4} - 3A^{3} + 3A^{2} + A - 2$}
\end{minipage}

\noindent{\color{gray!40}\rule{\textwidth}{0.4pt}}
\vspace{0.9\baselineskip}
\noindent \begin{minipage}[t]{0.25\textwidth}
\vspace{0pt}
\centering
\includegraphics[page=165,width=\linewidth]{knotoids.pdf}
\end{minipage}
\hfill
\begin{minipage}[t]{0.73\textwidth}
\vspace{0pt}
\raggedright
\textbf{Name:} {\large{$\mathbf{K7_{46}}$}} (chiral, non-rotatable$^{*}$) \\ \textbf{PD:} {\scriptsize\texttt{[0],[0,1,2,3],[1,4,5,2],[3,6,7,8],[4,9,10,5],[6,11,12,7],[13,14,9,8],[14,12,11,10],[13]}} \\ \textbf{EM:} {\scriptsize\texttt{(B0, A0C0C3D0, B1E0E3B2, B3F0F3G3, C1G2H3C2, D1H2H1D2, I0H0E1D3, G1F2F1E2, G0)}} \\ \textbf{Kauffman bracket:} {\scriptsize $-A^{27} - A^{25} + A^{23} + A^{21} - 2A^{19} - 3A^{17} + 2A^{15} + 4A^{13} - 3A^{9} + 2A^{5} - A$} \\ \textbf{Arrow:} {\scriptsize $L_1/A^{6} + A^{-8} - L_1/A^{10} - 1/A^{12} + 2L_1/A^{14} + 3/A^{16} - 2L_1/A^{18} - 4/A^{20} + 3/A^{24} - 2/A^{28} + A^{-32}$} \\ \textbf{Mock:} {\scriptsize $w^{5} + 2w^{4} - w^{3} - 3w^{2} + w + 4 - 1/w - 3/w^{2} + w^{-3} + 2/w^{4} - 1/w^{5} - 1/w^{6}$} \\ \textbf{Affine:} {\scriptsize $3t - 6 + 3/t$} \\ \textbf{Yamada:} {\scriptsize $-A^{26} - 3A^{23} - 2A^{22} + 2A^{21} - 3A^{20} - 4A^{19} + 3A^{18} + 2A^{17} - 4A^{16} + 3A^{15} + 4A^{14} - 3A^{13} + 3A^{11} - A^{10} - 2A^{9} + 3A^{7} - 4A^{6} - 3A^{5} + 6A^{4} - 3A^{3} - A^{2} + 3A - 1$}
\end{minipage}

\noindent{\color{gray!40}\rule{\textwidth}{0.4pt}}
\vspace{0.9\baselineskip}
\noindent \begin{minipage}[t]{0.25\textwidth}
\vspace{0pt}
\centering
\includegraphics[page=166,width=\linewidth]{knotoids.pdf}
\end{minipage}
\hfill
\begin{minipage}[t]{0.73\textwidth}
\vspace{0pt}
\raggedright
\textbf{Name:} {\large{$\mathbf{K7_{47}}$}} (chiral, non-rotatable$^{*}$) \\ \textbf{PD:} {\scriptsize\texttt{[0],[0,1,2,3],[1,4,5,2],[3,6,7,8],[4,9,6,5],[7,10,11,12],[8,13,14,9],[10,14,13,11],[12]}} \\ \textbf{EM:} {\scriptsize\texttt{(B0, A0C0C3D0, B1E0E3B2, B3E2F0G0, C1G3D1C2, D2H0H3I0, D3H2H1E1, F1G2G1F2, F3)}} \\ \textbf{Kauffman bracket:} {\scriptsize $A^{24} + A^{22} - 2A^{18} + 3A^{14} + A^{12} - 3A^{10} - 2A^{8} + 2A^{6} + 2A^{4} - A^{2} - 1$} \\ \textbf{Arrow:} {\scriptsize $A^{-12} + L_1/A^{14} - 2L_1/A^{18} + 3L_1/A^{22} + A^{-24} - 3L_1/A^{26} - 2/A^{28} + 2L_1/A^{30} + 2/A^{32} - L_1/A^{34} - 1/A^{36}$} \\ \textbf{Mock:} {\scriptsize $w^{6} + w^{5} - 2w^{3} - w^{2} + 3w + 2 - 3/w - 2/w^{2} + 2/w^{3} + w^{-4} - 1/w^{5}$} \\ \textbf{Affine:} {\scriptsize $2t - 4 + 2/t$} \\ \textbf{Yamada:} {\scriptsize $A^{27} + A^{26} + 4A^{23} + 2A^{22} - 2A^{21} + 3A^{20} + 4A^{19} - 5A^{18} + 4A^{16} - 4A^{15} - A^{14} + 3A^{13} - A^{12} - A^{11} - 2A^{10} + 3A^{9} - 3A^{8} - 5A^{7} + 6A^{6} - A^{5} - 3A^{4} + 3A^{3} - A + 1$}
\end{minipage}

\noindent{\color{gray!40}\rule{\textwidth}{0.4pt}}
\vspace{0.9\baselineskip}
\noindent \begin{minipage}[t]{0.25\textwidth}
\vspace{0pt}
\centering
\includegraphics[page=167,width=\linewidth]{knotoids.pdf}
\end{minipage}
\hfill
\begin{minipage}[t]{0.73\textwidth}
\vspace{0pt}
\raggedright
\textbf{Name:} {\large{$\mathbf{K7_{48}}$}} (chiral, non-rotatable$^{*}$) \\ \textbf{PD:} {\scriptsize\texttt{[0],[0,1,2,3],[1,4,5,2],[3,6,7,8],[4,9,6,5],[10,11,12,7],[13,14,9,8],[14,13,11,10],[12]}} \\ \textbf{EM:} {\scriptsize\texttt{(B0, A0C0C3D0, B1E0E3B2, B3E2F3G3, C1G2D1C2, H3H2I0D2, H1H0E1D3, G1G0F1F0, F2)}} \\ \textbf{Kauffman bracket:} {\scriptsize $A^{28} - 2A^{24} - A^{22} + 3A^{20} + 3A^{18} - 3A^{16} - 3A^{14} + 2A^{12} + 3A^{10} - A^{8} - 2A^{6} + A^{2}$} \\ \textbf{Arrow:} {\scriptsize $A^{10}L_1 - 2A^{6}L_1 - A^{4} + 3A^{2}L_1 + 3 - 3L_1/A^{2} - 3/A^{4} + 2L_1/A^{6} + 3/A^{8} - L_1/A^{10} - 2/A^{12} + A^{-16}$} \\ \textbf{Mock:} {\scriptsize $w^{5} + w^{4} - 2w^{3} - 2w^{2} + 3w + 4 - 3/w - 3/w^{2} + 2/w^{3} + 2/w^{4} - 1/w^{5} - 1/w^{6}$} \\ \textbf{Affine:} {\scriptsize $2t - 4 + 2/t$} \\ \textbf{Yamada:} {\scriptsize $-A^{26} + 2A^{25} + 2A^{24} - 4A^{23} + A^{22} + 5A^{21} - 6A^{20} - 2A^{19} + 7A^{18} - 4A^{17} - A^{16} + 6A^{15} - 2A^{12} + 5A^{11} - A^{10} - 6A^{9} + 7A^{8} - 6A^{6} + 4A^{5} + 2A^{4} - 3A^{3} + A^{2} + A - 1$}
\end{minipage}

\noindent{\color{gray!40}\rule{\textwidth}{0.4pt}}
\vspace{0.9\baselineskip}
\noindent \begin{minipage}[t]{0.25\textwidth}
\vspace{0pt}
\centering
\includegraphics[page=168,width=\linewidth]{knotoids.pdf}
\end{minipage}
\hfill
\begin{minipage}[t]{0.73\textwidth}
\vspace{0pt}
\raggedright
\textbf{Name:} {\large{$\mathbf{K7_{49}}$}} (chiral, non-rotatable$^{*}$) \\ \textbf{PD:} {\scriptsize\texttt{[0],[0,1,2,3],[1,4,5,2],[3,6,7,8],[4,9,10,5],[6,10,11,12],[7],[13,14,9,8],[14,13,12,11]}} \\ \textbf{EM:} {\scriptsize\texttt{(B0, A0C0C3D0, B1E0E3B2, B3F0G0H3, C1H2F1C2, D1E2I3I2, D2, I1I0E1D3, H1H0F3F2)}} \\ \textbf{Kauffman bracket:} {\scriptsize $A^{28} - 2A^{24} - A^{22} + 3A^{20} + 2A^{18} - 3A^{16} - 3A^{14} + 2A^{12} + 2A^{10} - A^{8} - A^{6} + A^{4} + A^{2}$} \\ \textbf{Arrow:} {\scriptsize $A^{28} - 2A^{24} - A^{22}L_1 + 3A^{20} + 2A^{18}L_1 - 3A^{16} - 3A^{14}L_1 + 2A^{12} + 2A^{10}L_1 - A^{8} - A^{6}L_1 + A^{4} + A^{2}L_1$} \\ \textbf{Mock:} {\scriptsize $-2w^{4} - 2w^{3} + 4w^{2} + 3w - 4 - 3/w + 3/w^{2} + 2/w^{3}$} \\ \textbf{Affine:} {\scriptsize $-3t + 6 - 3/t$} \\ \textbf{Yamada:} {\scriptsize $A^{27} - A^{26} - A^{25} + 3A^{24} - A^{23} - 4A^{22} + 4A^{21} + A^{20} - 5A^{19} + 4A^{18} - 3A^{16} + A^{15} + 2A^{13} - 3A^{12} + A^{11} + 4A^{10} - 5A^{9} + 4A^{7} - 4A^{6} - A^{5} + A^{4} - 2A^{3} - A^{2} - 1$}
\end{minipage}

\noindent{\color{gray!40}\rule{\textwidth}{0.4pt}}
\vspace{0.9\baselineskip}
\noindent \begin{minipage}[t]{0.25\textwidth}
\vspace{0pt}
\centering
\includegraphics[page=169,width=\linewidth]{knotoids.pdf}
\end{minipage}
\hfill
\begin{minipage}[t]{0.73\textwidth}
\vspace{0pt}
\raggedright
\textbf{Name:} {\large{$\mathbf{K7_{50}}$}} (achiral, non-rotatable$^{*}$, Possible duplicate [7\_52]) \\ \textbf{PD:} {\scriptsize\texttt{[0],[0,1,2,3],[1,4,5,2],[3,6,7,8],[4,8,9,10],[11,12,13,5],[6,13,12,7],[14,11,10,9],[14]}} \\ \textbf{EM:} {\scriptsize\texttt{(B0, A0C0C3D0, B1E0F3B2, B3G0G3E1, C1D3H3H2, H1G2G1C2, D1F2F1D2, I0F0E3E2, H0)}} \\ \textbf{Kauffman bracket:} {\scriptsize $-A^{26} + 2A^{22} + 2A^{20} - 3A^{18} - 3A^{16} + 3A^{14} + 5A^{12} - A^{10} - 4A^{8} + 2A^{4} - 1$} \\ \textbf{Arrow:} {\scriptsize $-A^{2}L_1 + 2L_1/A^{2} + 2/A^{4} - 3L_1/A^{6} - 3/A^{8} + 3L_1/A^{10} + 5/A^{12} - L_1/A^{14} - 4/A^{16} + 2/A^{20} - 1/A^{24}$} \\ \textbf{Mock:} {\scriptsize $-w^{5} - w^{4} + 2w^{3} + 4w^{2} - 2w - 5 + 2/w + 4/w^{2} - 2/w^{3} - 2/w^{4} + w^{-5} + w^{-6}$} \\ \textbf{Affine:} {\scriptsize $-t + 2 - 1/t$} \\ \textbf{Yamada:} {\scriptsize $-2A^{25} + 4A^{23} - 3A^{22} - 4A^{21} + 8A^{20} + 2A^{19} - 8A^{18} + 7A^{17} + 6A^{16} - 7A^{15} + 3A^{14} + 6A^{13} - 3A^{12} - 2A^{11} + A^{10} + 5A^{9} - 8A^{8} - 2A^{7} + 10A^{6} - 8A^{5} - 4A^{4} + 7A^{3} - 2A^{2} - 2A + 2$}
\end{minipage}

\noindent{\color{gray!40}\rule{\textwidth}{0.4pt}}
\vspace{0.9\baselineskip}
\noindent \begin{minipage}[t]{0.25\textwidth}
\vspace{0pt}
\centering
\includegraphics[page=170,width=\linewidth]{knotoids.pdf}
\end{minipage}
\hfill
\begin{minipage}[t]{0.73\textwidth}
\vspace{0pt}
\raggedright
\textbf{Name:} {\large{$\mathbf{K7_{51}}$}} (achiral, non-rotatable$^{*}$, Possible duplicate [7\_75]) \\ \textbf{PD:} {\scriptsize\texttt{[0],[0,1,2,3],[1,4,5,2],[3,6,7,8],[9,10,11,4],[5,12,13,6],[10,9,8,7],[11,13,12,14],[14]}} \\ \textbf{EM:} {\scriptsize\texttt{(B0, A0C0C3D0, B1E3F0B2, B3F3G3G2, G1G0H0C1, C2H2H1D1, E1E0D3D2, E2F2F1I0, H3)}} \\ \textbf{Kauffman bracket:} {\scriptsize $A^{21} + A^{19} - 2A^{17} - 3A^{15} + A^{13} + 3A^{11} - 3A^{7} - A^{5} + A^{3} + A$} \\ \textbf{Arrow:} {\scriptsize $-A^{6}L_1 - A^{4} + 2A^{2}L_1 + 3 - L_1/A^{2} - 3/A^{4} + 3/A^{8} + L_1/A^{10} - 1/A^{12} - L_1/A^{14}$} \\ \textbf{Mock:} {\scriptsize $w^{4} - 2w^{2} + w + 4 - 1/w - 3/w^{2} + w^{-4}$} \\ \textbf{Affine:} {\scriptsize $t - 2 + 1/t$} \\ \textbf{Yamada:} {\scriptsize $A^{24} - 2A^{23} - 2A^{22} + 5A^{21} - A^{20} - 5A^{19} + 5A^{18} + A^{17} - 3A^{16} + 3A^{15} + 2A^{14} + 2A^{13} - 2A^{12} + 3A^{11} + 3A^{10} - 5A^{9} + 2A^{8} + 4A^{7} - 5A^{6} - A^{5} + 3A^{4} - A^{3} - 2A^{2} + 1$}
\end{minipage}

\noindent{\color{gray!40}\rule{\textwidth}{0.4pt}}
\vspace{0.9\baselineskip}
\noindent \begin{minipage}[t]{0.25\textwidth}
\vspace{0pt}
\centering
\includegraphics[page=171,width=\linewidth]{knotoids.pdf}
\end{minipage}
\hfill
\begin{minipage}[t]{0.73\textwidth}
\vspace{0pt}
\raggedright
\textbf{Name:} {\large{$\mathbf{K7_{52}}$}} (achiral, non-rotatable$^{*}$, Possible duplicate [7\_50]) \\ \textbf{PD:} {\scriptsize\texttt{[0],[0,1,2,3],[1,4,5,2],[3,6,7,8],[8,7,9,4],[10,11,6,5],[12,13,14,9],[13,12,11,10],[14]}} \\ \textbf{EM:} {\scriptsize\texttt{(B0, A0C0C3D0, B1E3F3B2, B3F2E1E0, D3D2G3C1, H3H2D1C2, H1H0I0E2, G1G0F1F0, G2)}} \\ \textbf{Kauffman bracket:} {\scriptsize $A^{27} - 2A^{23} + 4A^{19} + A^{17} - 5A^{15} - 3A^{13} + 3A^{11} + 3A^{9} - 2A^{7} - 2A^{5} + A$} \\ \textbf{Arrow:} {\scriptsize $-A^{24} + 2A^{20} - 4A^{16} - A^{14}L_1 + 5A^{12} + 3A^{10}L_1 - 3A^{8} - 3A^{6}L_1 + 2A^{4} + 2A^{2}L_1 - L_1/A^{2}$} \\ \textbf{Mock:} {\scriptsize $w^{6} + w^{5} - 2w^{4} - 2w^{3} + 4w^{2} + 2w - 5 - 2/w + 4/w^{2} + 2/w^{3} - 1/w^{4} - 1/w^{5}$} \\ \textbf{Affine:} {\scriptsize $t - 2 + 1/t$} \\ \textbf{Yamada:} {\scriptsize $-2A^{25} + 2A^{24} + 2A^{23} - 7A^{22} + 4A^{21} + 8A^{20} - 10A^{19} + 2A^{18} + 8A^{17} - 5A^{16} - A^{15} + 2A^{14} + 3A^{13} - 6A^{12} - 3A^{11} + 7A^{10} - 6A^{9} - 7A^{8} + 8A^{7} - 2A^{6} - 8A^{5} + 4A^{4} + 3A^{3} - 4A^{2} + 2$}
\end{minipage}

\noindent{\color{gray!40}\rule{\textwidth}{0.4pt}}
\vspace{0.9\baselineskip}
\noindent \begin{minipage}[t]{0.25\textwidth}
\vspace{0pt}
\centering
\includegraphics[page=172,width=\linewidth]{knotoids.pdf}
\end{minipage}
\hfill
\begin{minipage}[t]{0.73\textwidth}
\vspace{0pt}
\raggedright
\textbf{Name:} {\large{$\mathbf{K7_{53}}$}} (achiral, non-rotatable$^{*}$, Possible duplicate [7\_54]) \\ \textbf{PD:} {\scriptsize\texttt{[0],[0,1,2,3],[1,4,5,2],[3,6,7,8],[8,9,10,4],[11,12,13,5],[6,13,12,7],[9,14,11,10],[14]}} \\ \textbf{EM:} {\scriptsize\texttt{(B0, A0C0C3D0, B1E3F3B2, B3G0G3E0, D3H0H3C1, H2G2G1C2, D1F2F1D2, E1I0F0E2, H1)}} \\ \textbf{Kauffman bracket:} {\scriptsize $A^{24} + A^{22} - 2A^{18} + 3A^{14} + 2A^{12} - 3A^{10} - 3A^{8} + A^{6} + 2A^{4} - 1$} \\ \textbf{Arrow:} {\scriptsize $A^{-12} + L_1/A^{14} - 2L_1/A^{18} + 3L_1/A^{22} + 2/A^{24} - 3L_1/A^{26} - 3/A^{28} + L_1/A^{30} + 2/A^{32} - 1/A^{36}$} \\ \textbf{Mock:} {\scriptsize $w^{6} + w^{5} - 2w^{3} - 2w^{2} + 2w + 3 - 2/w - 2/w^{2} + 2/w^{3} + w^{-4} - 1/w^{5}$} \\ \textbf{Affine:} {\scriptsize $t - 2 + 1/t$} \\ \textbf{Yamada:} {\scriptsize $-A^{26} - A^{25} - 4A^{22} - 3A^{21} + 3A^{20} - A^{19} - 7A^{18} + 4A^{17} + 2A^{16} - 6A^{15} + 3A^{14} + 4A^{13} - 2A^{12} + A^{11} + 2A^{10} + 3A^{9} - 4A^{8} + 7A^{6} - 6A^{5} - 2A^{4} + 5A^{3} - 2A^{2} - 2A + 1$}
\end{minipage}

\noindent{\color{gray!40}\rule{\textwidth}{0.4pt}}
\vspace{0.9\baselineskip}
\noindent \begin{minipage}[t]{0.25\textwidth}
\vspace{0pt}
\centering
\includegraphics[page=173,width=\linewidth]{knotoids.pdf}
\end{minipage}
\hfill
\begin{minipage}[t]{0.73\textwidth}
\vspace{0pt}
\raggedright
\textbf{Name:} {\large{$\mathbf{K7_{54}}$}} (achiral, non-rotatable$^{*}$, Possible duplicate [7\_53]) \\ \textbf{PD:} {\scriptsize\texttt{[0],[0,1,2,3],[1,4,5,2],[6,7,8,3],[4,8,7,9],[5,10,11,6],[9,12,13,14],[10,13,12,11],[14]}} \\ \textbf{EM:} {\scriptsize\texttt{(B0, A0C0C3D3, B1E0F0B2, F3E2E1B3, C1D2D1G0, C2H0H3D0, E3H2H1I0, F1G2G1F2, G3)}} \\ \textbf{Kauffman bracket:} {\scriptsize $A^{24} + A^{22} - 2A^{18} + 3A^{14} + 2A^{12} - 3A^{10} - 3A^{8} + A^{6} + 2A^{4} - 1$} \\ \textbf{Arrow:} {\scriptsize $A^{-12} + L_1/A^{14} - 2L_1/A^{18} + 3L_1/A^{22} + 2/A^{24} - 3L_1/A^{26} - 3/A^{28} + L_1/A^{30} + 2/A^{32} - 1/A^{36}$} \\ \textbf{Mock:} {\scriptsize $w^{6} + w^{5} - 2w^{3} - 2w^{2} + 2w + 3 - 2/w - 2/w^{2} + 2/w^{3} + w^{-4} - 1/w^{5}$} \\ \textbf{Affine:} {\scriptsize $t - 2 + 1/t$} \\ \textbf{Yamada:} {\scriptsize $-A^{26} - A^{25} - 4A^{22} - 3A^{21} + 3A^{20} - A^{19} - 7A^{18} + 4A^{17} + 2A^{16} - 6A^{15} + 3A^{14} + 4A^{13} - 2A^{12} + A^{11} + 2A^{10} + 3A^{9} - 4A^{8} + 7A^{6} - 6A^{5} - 2A^{4} + 5A^{3} - 2A^{2} - 2A + 1$}
\end{minipage}

\noindent{\color{gray!40}\rule{\textwidth}{0.4pt}}
\vspace{0.9\baselineskip}
\noindent \begin{minipage}[t]{0.25\textwidth}
\vspace{0pt}
\centering
\includegraphics[page=174,width=\linewidth]{knotoids.pdf}
\end{minipage}
\hfill
\begin{minipage}[t]{0.73\textwidth}
\vspace{0pt}
\raggedright
\textbf{Name:} {\large{$\mathbf{K7_{55}}$}} (chiral, non-rotatable$^{*}$) \\ \textbf{PD:} {\scriptsize\texttt{[0],[0,1,2,3],[1,4,5,2],[3,6,7,4],[8,9,10,5],[6,11,12,7],[8],[9,12,13,14],[14,13,11,10]}} \\ \textbf{EM:} {\scriptsize\texttt{(B0, A0C0C3D0, B1D3E3B2, B3F0F3C1, G0H0I3C2, D1I2H1D2, E0, E1F2I1I0, H3H2F1E2)}} \\ \textbf{Kauffman bracket:} {\scriptsize $-A^{26} + 2A^{22} + A^{20} - 3A^{18} - 3A^{16} + 3A^{14} + 4A^{12} - A^{10} - 3A^{8} + A^{6} + 2A^{4} - 1$} \\ \textbf{Arrow:} {\scriptsize $-A^{20} + 2A^{16} + A^{14}L_1 - 3A^{12} - 3A^{10}L_1 + 3A^{8} + 4A^{6}L_1 - A^{4} - 3A^{2}L_1 + 1 + 2L_1/A^{2} - L_1/A^{6}$} \\ \textbf{Mock:} {\scriptsize $2w^{4} + 2w^{3} - 4w^{2} - 5w + 4 + 5/w - 1/w^{2} - 2/w^{3}$} \\ \textbf{Affine:} {\scriptsize $t - 2 + 1/t$} \\ \textbf{Yamada:} {\scriptsize $A^{26} - A^{25} - 3A^{24} + 3A^{23} + 3A^{22} - 6A^{21} + A^{20} + 7A^{19} - 5A^{18} - 2A^{17} + 7A^{16} - 3A^{15} - 3A^{14} + 3A^{13} - 2A^{11} - 5A^{10} + 5A^{9} - 2A^{8} - 9A^{7} + 6A^{6} + A^{5} - 6A^{4} + 3A^{3} + A^{2} - A + 1$}
\end{minipage}

\noindent{\color{gray!40}\rule{\textwidth}{0.4pt}}
\vspace{0.9\baselineskip}
\noindent \begin{minipage}[t]{0.25\textwidth}
\vspace{0pt}
\centering
\includegraphics[page=175,width=\linewidth]{knotoids.pdf}
\end{minipage}
\hfill
\begin{minipage}[t]{0.73\textwidth}
\vspace{0pt}
\raggedright
\textbf{Name:} {\large{$\mathbf{K7_{56}}$}} (chiral, non-rotatable$^{*}$) \\ \textbf{PD:} {\scriptsize\texttt{[0],[0,1,2,3],[1,4,5,2],[3,5,6,7],[4,8,9,6],[7,10,11,8],[12,13,14,9],[10,14,13,11],[12]}} \\ \textbf{EM:} {\scriptsize\texttt{(B0, A0C0C3D0, B1E0D1B2, B3C2E3F0, C1F3G3D2, D3H0H3E1, I0H2H1E2, F1G2G1F2, G0)}} \\ \textbf{Kauffman bracket:} {\scriptsize $-A^{26} + 2A^{22} + A^{20} - 3A^{18} - 2A^{16} + 4A^{14} + 4A^{12} - 3A^{10} - 4A^{8} + A^{6} + 3A^{4} - 1$} \\ \textbf{Arrow:} {\scriptsize $-A^{14}L_1 + 2A^{10}L_1 + A^{8} - 3A^{6}L_1 - 2A^{4} + 4A^{2}L_1 + 4 - 3L_1/A^{2} - 4/A^{4} + L_1/A^{6} + 3/A^{8} - 1/A^{12}$} \\ \textbf{Mock:} {\scriptsize $-2w^{3} - 2w^{2} + 5w + 6 - 5/w - 5/w^{2} + 2/w^{3} + 2/w^{4}$} \\ \textbf{Affine:} {\scriptsize $-t + 2 - 1/t$} \\ \textbf{Yamada:} {\scriptsize $A^{26} - A^{25} - 2A^{24} + 4A^{23} + A^{22} - 7A^{21} + 5A^{20} + 6A^{19} - 11A^{18} + A^{17} + 7A^{16} - 8A^{15} - 2A^{14} + 5A^{13} - A^{12} - 3A^{11} - 2A^{10} + 7A^{9} - 5A^{8} - 8A^{7} + 11A^{6} - 3A^{5} - 7A^{4} + 8A^{3} - 3A + 1$}
\end{minipage}

\noindent{\color{gray!40}\rule{\textwidth}{0.4pt}}
\vspace{0.9\baselineskip}
\noindent \begin{minipage}[t]{0.25\textwidth}
\vspace{0pt}
\centering
\includegraphics[page=176,width=\linewidth]{knotoids.pdf}
\end{minipage}
\hfill
\begin{minipage}[t]{0.73\textwidth}
\vspace{0pt}
\raggedright
\textbf{Name:} {\large{$\mathbf{K7_{57}}$}} (chiral, non-rotatable$^{*}$) \\ \textbf{PD:} {\scriptsize\texttt{[0],[0,1,2,3],[1,4,5,2],[3,5,6,7],[4,8,9,6],[10,11,8,7],[9,12,13,14],[14,13,11,10],[12]}} \\ \textbf{EM:} {\scriptsize\texttt{(B0, A0C0C3D0, B1E0D1B2, B3C2E3F3, C1F2G0D2, H3H2E1D3, E2I0H1H0, G3G2F1F0, G1)}} \\ \textbf{Kauffman bracket:} {\scriptsize $-A^{25} - A^{23} + 2A^{21} + 2A^{19} - 2A^{17} - 4A^{15} + A^{13} + 4A^{11} - 3A^{7} - A^{5} + A^{3} + A$} \\ \textbf{Arrow:} {\scriptsize $L_1/A^{2} + A^{-4} - 2L_1/A^{6} - 2/A^{8} + 2L_1/A^{10} + 4/A^{12} - L_1/A^{14} - 4/A^{16} + 3/A^{20} + L_1/A^{22} - 1/A^{24} - L_1/A^{26}$} \\ \textbf{Mock:} {\scriptsize $w^{3} + 4w^{2} - w - 6 + 1/w + 4/w^{2} - 1/w^{3} - 1/w^{4}$} \\ \textbf{Affine:} {\scriptsize $2t - 4 + 2/t$} \\ \textbf{Yamada:} {\scriptsize $A^{27} - A^{26} - A^{25} + 4A^{24} - A^{23} - 5A^{22} + 8A^{21} + 3A^{20} - 7A^{19} + 6A^{18} + 4A^{17} - 5A^{16} + A^{14} + 2A^{13} - 5A^{12} + 2A^{11} + 7A^{10} - 8A^{9} - A^{8} + 7A^{7} - 5A^{6} - 3A^{5} + 4A^{4} - 2A^{2} + 1$}
\end{minipage}

\noindent{\color{gray!40}\rule{\textwidth}{0.4pt}}
\vspace{0.9\baselineskip}
\noindent \begin{minipage}[t]{0.25\textwidth}
\vspace{0pt}
\centering
\includegraphics[page=177,width=\linewidth]{knotoids.pdf}
\end{minipage}
\hfill
\begin{minipage}[t]{0.73\textwidth}
\vspace{0pt}
\raggedright
\textbf{Name:} {\large{$\mathbf{K7_{58}}$}} (chiral, non-rotatable$^{*}$) \\ \textbf{PD:} {\scriptsize\texttt{[0],[0,1,2,3],[1,4,5,2],[3,5,6,7],[4,7,8,9],[10,11,12,6],[13,10,9,8],[11,13,14,12],[14]}} \\ \textbf{EM:} {\scriptsize\texttt{(B0, A0C0C3D0, B1E0D1B2, B3C2F3E1, C1D3G3G2, G1H0H3D2, H1F0E3E2, F1G0I0F2, H2)}} \\ \textbf{Kauffman bracket:} {\scriptsize $A^{28} - 2A^{24} - A^{22} + 4A^{20} + 3A^{18} - 4A^{16} - 4A^{14} + 3A^{12} + 4A^{10} - A^{8} - 3A^{6} + A^{2}$} \\ \textbf{Arrow:} {\scriptsize $A^{4} - 2 - L_1/A^{2} + 4/A^{4} + 3L_1/A^{6} - 4/A^{8} - 4L_1/A^{10} + 3/A^{12} + 4L_1/A^{14} - 1/A^{16} - 3L_1/A^{18} + L_1/A^{22}$} \\ \textbf{Mock:} {\scriptsize $-2w^{4} - 2w^{3} + 6w^{2} + 6w - 5 - 6/w + 2/w^{2} + 2/w^{3}$} \\ \textbf{Affine:} {\scriptsize $0$} \\ \textbf{Yamada:} {\scriptsize $-A^{26} + 2A^{25} + 2A^{24} - 6A^{23} + 3A^{22} + 8A^{21} - 11A^{20} + A^{19} + 12A^{18} - 8A^{17} - A^{16} + 9A^{15} - A^{14} - 2A^{13} - 2A^{12} + 6A^{11} - 4A^{10} - 9A^{9} + 12A^{8} - A^{7} - 10A^{6} + 9A^{5} + 2A^{4} - 7A^{3} + 2A^{2} + 2A - 1$}
\end{minipage}

\noindent{\color{gray!40}\rule{\textwidth}{0.4pt}}
\vspace{0.9\baselineskip}
\noindent \begin{minipage}[t]{0.25\textwidth}
\vspace{0pt}
\centering
\includegraphics[page=178,width=\linewidth]{knotoids.pdf}
\end{minipage}
\hfill
\begin{minipage}[t]{0.73\textwidth}
\vspace{0pt}
\raggedright
\textbf{Name:} {\large{$\mathbf{K7_{59}}$}} (achiral, non-rotatable$^{*}$, Possible duplicate [7\_34]) \\ \textbf{PD:} {\scriptsize\texttt{[0],[0,1,2,3],[1,4,5,2],[3,5,6,7],[4,8,9,10],[11,12,13,6],[7,13,14,8],[12,11,10,9],[14]}} \\ \textbf{EM:} {\scriptsize\texttt{(B0, A0C0C3D0, B1E0D1B2, B3C2F3G0, C1G3H3H2, H1H0G1D2, D3F2I0E1, F1F0E3E2, G2)}} \\ \textbf{Kauffman bracket:} {\scriptsize $A^{27} - 2A^{23} + 4A^{19} + A^{17} - 5A^{15} - 3A^{13} + 3A^{11} + 2A^{9} - 2A^{7} - 2A^{5} + A^{3} + A$} \\ \textbf{Arrow:} {\scriptsize $-A^{18}L_1 + 2A^{14}L_1 - 4A^{10}L_1 - A^{8} + 5A^{6}L_1 + 3A^{4} - 3A^{2}L_1 - 2 + 2L_1/A^{2} + 2/A^{4} - L_1/A^{6} - 1/A^{8}$} \\ \textbf{Mock:} {\scriptsize $-w^{5} + 3w^{3} + w^{2} - 5w - 2 + 5/w + 3/w^{2} - 3/w^{3} - 2/w^{4} + w^{-5} + w^{-6}$} \\ \textbf{Affine:} {\scriptsize $-t + 2 - 1/t$} \\ \textbf{Yamada:} {\scriptsize $2A^{26} - 2A^{25} - 2A^{24} + 6A^{23} - 4A^{22} - 7A^{21} + 8A^{20} - 3A^{19} - 6A^{18} + 6A^{17} + A^{16} - A^{14} + 6A^{13} + 3A^{12} - 5A^{11} + 6A^{10} + 5A^{9} - 8A^{8} + A^{7} + 5A^{6} - 5A^{5} - 2A^{4} + 3A^{3} - A^{2} - A + 1$}
\end{minipage}

\noindent{\color{gray!40}\rule{\textwidth}{0.4pt}}
\vspace{0.9\baselineskip}
\noindent \begin{minipage}[t]{0.25\textwidth}
\vspace{0pt}
\centering
\includegraphics[page=179,width=\linewidth]{knotoids.pdf}
\end{minipage}
\hfill
\begin{minipage}[t]{0.73\textwidth}
\vspace{0pt}
\raggedright
\textbf{Name:} {\large{$\mathbf{K7_{60}}$}} (achiral, non-rotatable$^{*}$, Possible duplicate [7\_35]) \\ \textbf{PD:} {\scriptsize\texttt{[0],[0,1,2,3],[1,4,5,2],[3,5,6,7],[4,8,9,10],[11,12,13,6],[13,14,8,7],[12,11,10,9],[14]}} \\ \textbf{EM:} {\scriptsize\texttt{(B0, A0C0C3D0, B1E0D1B2, B3C2F3G3, C1G2H3H2, H1H0G0D2, F2I0E1D3, F1F0E3E2, G1)}} \\ \textbf{Kauffman bracket:} {\scriptsize $A^{21} + A^{19} - A^{17} - 3A^{15} + 2A^{11} - 2A^{7} - A^{5} + A^{3} + A$} \\ \textbf{Arrow:} {\scriptsize $-A^{6}L_1 - A^{4} + A^{2}L_1 + 3 - 2/A^{4} + 2/A^{8} + L_1/A^{10} - 1/A^{12} - L_1/A^{14}$} \\ \textbf{Mock:} {\scriptsize $w^{4} - 2w^{2} + 3 - 2/w^{2} + w^{-4}$} \\ \textbf{Affine:} {\scriptsize $0$} \\ \textbf{Yamada:} {\scriptsize $A^{24} - A^{23} - 2A^{22} + 2A^{21} - 4A^{19} + 2A^{18} + 2A^{17} - A^{16} + 3A^{15} + 2A^{14} + 3A^{13} - A^{12} + 2A^{11} + 3A^{10} - 3A^{9} + 2A^{7} - 3A^{6} - 2A^{5} + A^{4} - A^{2} + 1$}
\end{minipage}

\noindent{\color{gray!40}\rule{\textwidth}{0.4pt}}
\vspace{0.9\baselineskip}
\noindent \begin{minipage}[t]{0.25\textwidth}
\vspace{0pt}
\centering
\includegraphics[page=180,width=\linewidth]{knotoids.pdf}
\end{minipage}
\hfill
\begin{minipage}[t]{0.73\textwidth}
\vspace{0pt}
\raggedright
\textbf{Name:} {\large{$\mathbf{K7_{61}}$}} (chiral, non-rotatable$^{*}$) \\ \textbf{PD:} {\scriptsize\texttt{[0],[0,1,2,3],[1,4,5,2],[3,6,7,8],[4,9,10,5],[6,11,12,13],[7],[13,14,9,8],[14,12,11,10]}} \\ \textbf{EM:} {\scriptsize\texttt{(B0, A0C0C3D0, B1E0E3B2, B3F0G0H3, C1H2I3C2, D1I2I1H0, D2, F3I0E1D3, H1F2F1E2)}} \\ \textbf{Kauffman bracket:} {\scriptsize $-A^{26} + 2A^{22} + A^{20} - 3A^{18} - 2A^{16} + 3A^{14} + 3A^{12} - 2A^{10} - 2A^{8} + 2A^{6} + 2A^{4} - A^{2} - 1$} \\ \textbf{Arrow:} {\scriptsize $-A^{14}L_1 + 2A^{10}L_1 + A^{8} - 3A^{6}L_1 - 2A^{4} + 3A^{2}L_1 + 3 - 2L_1/A^{2} - 2/A^{4} + 2L_1/A^{6} + 2/A^{8} - L_1/A^{10} - 1/A^{12}$} \\ \textbf{Mock:} {\scriptsize $-2w^{3} - w^{2} + 5w + 4 - 5/w - 4/w^{2} + 2/w^{3} + 2/w^{4}$} \\ \textbf{Affine:} {\scriptsize $-t + 2 - 1/t$} \\ \textbf{Yamada:} {\scriptsize $-A^{27} + A^{26} + 2A^{25} - 3A^{24} - A^{23} + 4A^{22} - 3A^{21} - 3A^{20} + 5A^{19} - A^{18} - 3A^{17} + 4A^{16} + A^{15} - A^{14} + 2A^{13} + 2A^{12} + 3A^{11} - 4A^{10} + 3A^{9} + 4A^{8} - 6A^{7} + 2A^{6} + 2A^{5} - 4A^{4} + A^{3} - A + 1$}
\end{minipage}

\noindent{\color{gray!40}\rule{\textwidth}{0.4pt}}
\vspace{0.9\baselineskip}
\noindent \begin{minipage}[t]{0.25\textwidth}
\vspace{0pt}
\centering
\includegraphics[page=181,width=\linewidth]{knotoids.pdf}
\end{minipage}
\hfill
\begin{minipage}[t]{0.73\textwidth}
\vspace{0pt}
\raggedright
\textbf{Name:} {\large{$\mathbf{K7_{62}}$}} (chiral, non-rotatable$^{*}$) \\ \textbf{PD:} {\scriptsize\texttt{[0],[0,1,2,3],[1,4,5,2],[3,5,6,7],[4,7,8,9],[9,10,11,6],[12,13,14,8],[10,14,13,11],[12]}} \\ \textbf{EM:} {\scriptsize\texttt{(B0, A0C0C3D0, B1E0D1B2, B3C2F3E1, C1D3G3F0, E3H0H3D2, I0H2H1E2, F1G2G1F2, G0)}} \\ \textbf{Kauffman bracket:} {\scriptsize $-A^{26} + 2A^{22} + 2A^{20} - 3A^{18} - 3A^{16} + 3A^{14} + 4A^{12} - 2A^{10} - 4A^{8} + A^{6} + 3A^{4} - 1$} \\ \textbf{Arrow:} {\scriptsize $-A^{14}L_1 + 2A^{10}L_1 + 2A^{8} - 3A^{6}L_1 - 3A^{4} + 3A^{2}L_1 + 4 - 2L_1/A^{2} - 4/A^{4} + L_1/A^{6} + 3/A^{8} - 1/A^{12}$} \\ \textbf{Mock:} {\scriptsize $-2w^{3} - 3w^{2} + 4w + 7 - 4/w - 5/w^{2} + 2/w^{3} + 2/w^{4}$} \\ \textbf{Affine:} {\scriptsize $-2t + 4 - 2/t$} \\ \textbf{Yamada:} {\scriptsize $2A^{26} - 4A^{24} + 3A^{23} + 3A^{22} - 8A^{21} + A^{20} + 7A^{19} - 8A^{18} - A^{17} + 8A^{16} - 4A^{15} - 3A^{14} + 3A^{13} - 3A^{11} - 5A^{10} + 7A^{9} - A^{8} - 9A^{7} + 9A^{6} - 7A^{4} + 5A^{3} - 2A + 1$}
\end{minipage}

\noindent{\color{gray!40}\rule{\textwidth}{0.4pt}}
\vspace{0.9\baselineskip}
\noindent \begin{minipage}[t]{0.25\textwidth}
\vspace{0pt}
\centering
\includegraphics[page=182,width=\linewidth]{knotoids.pdf}
\end{minipage}
\hfill
\begin{minipage}[t]{0.73\textwidth}
\vspace{0pt}
\raggedright
\textbf{Name:} {\large{$\mathbf{K7_{63}}$}} (achiral, non-rotatable$^{*}$, Possible duplicate [7\_76]) \\ \textbf{PD:} {\scriptsize\texttt{[0],[0,1,2,3],[1,4,5,2],[3,6,7,8],[8,9,10,4],[5,11,12,6],[13,10,9,7],[11,14,13,12],[14]}} \\ \textbf{EM:} {\scriptsize\texttt{(B0, A0C0C3D0, B1E3F0B2, B3F3G3E0, D3G2G1C1, C2H0H3D1, H2E2E1D2, F1I0G0F2, H1)}} \\ \textbf{Kauffman bracket:} {\scriptsize $A^{22} + 2A^{20} - A^{18} - 3A^{16} + 4A^{12} + A^{10} - 3A^{8} - 2A^{6} + A^{4} + A^{2}$} \\ \textbf{Arrow:} {\scriptsize $L_1/A^{2} + 2/A^{4} - L_1/A^{6} - 3/A^{8} + 4/A^{12} + L_1/A^{14} - 3/A^{16} - 2L_1/A^{18} + A^{-20} + L_1/A^{22}$} \\ \textbf{Mock:} {\scriptsize $-w^{4} + 4w^{2} + w - 4 - 1/w + 3/w^{2} - 1/w^{4}$} \\ \textbf{Affine:} {\scriptsize $t - 2 + 1/t$} \\ \textbf{Yamada:} {\scriptsize $-2A^{24} - A^{23} + 3A^{22} - A^{21} - 6A^{20} + 4A^{19} + 2A^{18} - 7A^{17} + 2A^{16} + 3A^{15} - 4A^{14} + A^{12} + 2A^{11} - 4A^{10} + 6A^{8} - 6A^{7} - A^{6} + 6A^{5} - 2A^{4} - 3A^{3} + 2A^{2} + A - 1$}
\end{minipage}

\noindent{\color{gray!40}\rule{\textwidth}{0.4pt}}
\vspace{0.9\baselineskip}
\noindent \begin{minipage}[t]{0.25\textwidth}
\vspace{0pt}
\centering
\includegraphics[page=183,width=\linewidth]{knotoids.pdf}
\end{minipage}
\hfill
\begin{minipage}[t]{0.73\textwidth}
\vspace{0pt}
\raggedright
\textbf{Name:} {\large{$\mathbf{K7_{64}}$}} (achiral, non-rotatable$^{*}$, Possible duplicate [7\_80]) \\ \textbf{PD:} {\scriptsize\texttt{[0],[0,1,2,3],[1,4,5,2],[3,6,7,8],[8,7,9,4],[10,11,6,5],[11,12,13,9],[14,13,12,10],[14]}} \\ \textbf{EM:} {\scriptsize\texttt{(B0, A0C0C3D0, B1E3F3B2, B3F2E1E0, D3D2G3C1, H3G0D1C2, F1H2H1E2, I0G2G1F0, H0)}} \\ \textbf{Kauffman bracket:} {\scriptsize $A^{28} - 2A^{24} + 5A^{20} + 2A^{18} - 5A^{16} - 4A^{14} + 4A^{12} + 4A^{10} - 2A^{8} - 3A^{6} + A^{2}$} \\ \textbf{Arrow:} {\scriptsize $A^{-8} - 2/A^{12} + 5/A^{16} + 2L_1/A^{18} - 5/A^{20} - 4L_1/A^{22} + 4/A^{24} + 4L_1/A^{26} - 2/A^{28} - 3L_1/A^{30} + L_1/A^{34}$} \\ \textbf{Mock:} {\scriptsize $3w^{4} + 2w^{3} - 6w^{2} - 5w + 6 + 5/w - 3/w^{2} - 2/w^{3} + w^{-4}$} \\ \textbf{Affine:} {\scriptsize $t - 2 + 1/t$} \\ \textbf{Yamada:} {\scriptsize $3A^{25} - 5A^{23} + 8A^{22} + 7A^{21} - 13A^{20} + 7A^{19} + 12A^{18} - 13A^{17} + A^{16} + 9A^{15} - 6A^{14} - 4A^{13} + 6A^{11} - 8A^{10} - 6A^{9} + 15A^{8} - 7A^{7} - 10A^{6} + 13A^{5} - 8A^{3} + 5A^{2} + 2A - 2$}
\end{minipage}

\noindent{\color{gray!40}\rule{\textwidth}{0.4pt}}
\vspace{0.9\baselineskip}
\noindent \begin{minipage}[t]{0.25\textwidth}
\vspace{0pt}
\centering
\includegraphics[page=184,width=\linewidth]{knotoids.pdf}
\end{minipage}
\hfill
\begin{minipage}[t]{0.73\textwidth}
\vspace{0pt}
\raggedright
\textbf{Name:} {\large{$\mathbf{K7_{65}}$}} (achiral, non-rotatable$^{*}$, Possible duplicate [7\_81]) \\ \textbf{PD:} {\scriptsize\texttt{[0],[0,1,2,3],[1,4,5,2],[6,7,8,3],[4,8,7,9],[5,10,11,6],[9,11,12,13],[10,14,13,12],[14]}} \\ \textbf{EM:} {\scriptsize\texttt{(B0, A0C0C3D3, B1E0F0B2, F3E2E1B3, C1D2D1G0, C2H0G1D0, E3F2H3H2, F1I0G3G2, H1)}} \\ \textbf{Kauffman bracket:} {\scriptsize $-A^{26} - A^{24} + A^{22} + 3A^{20} - 4A^{16} - A^{14} + 4A^{12} + 3A^{10} - 2A^{8} - 2A^{6} + A^{2}$} \\ \textbf{Arrow:} {\scriptsize $-A^{20} - A^{18}L_1 + A^{16} + 3A^{14}L_1 - 4A^{10}L_1 - A^{8} + 4A^{6}L_1 + 3A^{4} - 2A^{2}L_1 - 2 + A^{-4}$} \\ \textbf{Mock:} {\scriptsize $w^{4} + 2w^{3} - 5w - 2 + 5/w + 3/w^{2} - 2/w^{3} - 1/w^{4}$} \\ \textbf{Affine:} {\scriptsize $t - 2 + 1/t$} \\ \textbf{Yamada:} {\scriptsize $A^{26} - 2A^{24} + 4A^{22} - 3A^{21} - 6A^{20} + 8A^{19} + A^{18} - 10A^{17} + 6A^{16} + 3A^{15} - 6A^{14} + 2A^{13} + 3A^{12} + 2A^{11} - 5A^{10} + 4A^{9} + 7A^{8} - 8A^{7} + A^{6} + 9A^{5} - 5A^{4} - 3A^{3} + 5A^{2} - 2$}
\end{minipage}

\noindent{\color{gray!40}\rule{\textwidth}{0.4pt}}
\vspace{0.9\baselineskip}
\noindent \begin{minipage}[t]{0.25\textwidth}
\vspace{0pt}
\centering
\includegraphics[page=185,width=\linewidth]{knotoids.pdf}
\end{minipage}
\hfill
\begin{minipage}[t]{0.73\textwidth}
\vspace{0pt}
\raggedright
\textbf{Name:} {\large{$\mathbf{K7_{66}}$}} (chiral, non-rotatable$^{*}$) \\ \textbf{PD:} {\scriptsize\texttt{[0],[0,1,2,3],[1,4,5,2],[3,6,7,8],[4,9,10,11],[11,12,6,5],[7],[13,14,9,8],[14,13,12,10]}} \\ \textbf{EM:} {\scriptsize\texttt{(B0, A0C0C3D0, B1E0F3B2, B3F2G0H3, C1H2I3F0, E3I2D1C2, D2, I1I0E1D3, H1H0F1E2)}} \\ \textbf{Kauffman bracket:} {\scriptsize $A^{28} - 3A^{24} - A^{22} + 4A^{20} + 3A^{18} - 4A^{16} - 3A^{14} + 3A^{12} + 3A^{10} - 2A^{8} - 2A^{6} + A^{4} + A^{2}$} \\ \textbf{Arrow:} {\scriptsize $A^{22}L_1 - 3A^{18}L_1 - A^{16} + 4A^{14}L_1 + 3A^{12} - 4A^{10}L_1 - 3A^{8} + 3A^{6}L_1 + 3A^{4} - 2A^{2}L_1 - 2 + L_1/A^{2} + A^{-4}$} \\ \textbf{Mock:} {\scriptsize $2w^{3} + w^{2} - 7w - 4 + 7/w + 6/w^{2} - 2/w^{3} - 2/w^{4}$} \\ \textbf{Affine:} {\scriptsize $-t + 2 - 1/t$} \\ \textbf{Yamada:} {\scriptsize $A^{27} - 2A^{26} - 2A^{25} + 6A^{24} - 2A^{23} - 6A^{22} + 9A^{21} - A^{20} - 8A^{19} + 8A^{18} - 5A^{16} + 2A^{15} + A^{13} - 6A^{12} + 2A^{11} + 6A^{10} - 10A^{9} + A^{8} + 7A^{7} - 8A^{6} + 4A^{4} - 3A^{3} + A - 1$}
\end{minipage}

\noindent{\color{gray!40}\rule{\textwidth}{0.4pt}}
\vspace{0.9\baselineskip}
\noindent \begin{minipage}[t]{0.25\textwidth}
\vspace{0pt}
\centering
\includegraphics[page=186,width=\linewidth]{knotoids.pdf}
\end{minipage}
\hfill
\begin{minipage}[t]{0.73\textwidth}
\vspace{0pt}
\raggedright
\textbf{Name:} {\large{$\mathbf{K7_{67}}$}} (chiral, non-rotatable$^{*}$) \\ \textbf{PD:} {\scriptsize\texttt{[0],[0,1,2,3],[1,4,5,6],[6,7,8,2],[3,9,10,11],[4,12,13,5],[7,14,9,8],[10],[14,13,12,11]}} \\ \textbf{EM:} {\scriptsize\texttt{(B0, A0C0D3E0, B1F0F3D0, C3G0G3B2, B3G2H0I3, C1I2I1C2, D1I0E1D2, E2, G1F2F1E3)}} \\ \textbf{Kauffman bracket:} {\scriptsize $-A^{26} + 2A^{22} + A^{20} - 3A^{18} - A^{16} + 4A^{14} + 3A^{12} - 3A^{10} - 3A^{8} + 2A^{6} + 2A^{4} - A^{2} - 1$} \\ \textbf{Arrow:} {\scriptsize $-A^{2}L_1 + 2L_1/A^{2} + A^{-4} - 3L_1/A^{6} - 1/A^{8} + 4L_1/A^{10} + 3/A^{12} - 3L_1/A^{14} - 3/A^{16} + 2L_1/A^{18} + 2/A^{20} - L_1/A^{22} - 1/A^{24}$} \\ \textbf{Mock:} {\scriptsize $-w^{5} + 3w^{3} + 2w^{2} - 4w - 3 + 4/w + 3/w^{2} - 3/w^{3} - 2/w^{4} + w^{-5} + w^{-6}$} \\ \textbf{Affine:} {\scriptsize $0$} \\ \textbf{Yamada:} {\scriptsize $-A^{27} + A^{26} + A^{25} - 4A^{24} + A^{23} + 5A^{22} - 4A^{21} - A^{20} + 8A^{19} - 2A^{18} - 2A^{17} + 7A^{16} + A^{15} - 2A^{14} + 2A^{13} + A^{12} - 6A^{10} + 4A^{9} + 2A^{8} - 9A^{7} + 4A^{6} + 2A^{5} - 5A^{4} + 2A^{3} + A^{2} - A + 1$}
\end{minipage}

\noindent{\color{gray!40}\rule{\textwidth}{0.4pt}}
\vspace{0.9\baselineskip}
\noindent \begin{minipage}[t]{0.25\textwidth}
\vspace{0pt}
\centering
\includegraphics[page=187,width=\linewidth]{knotoids.pdf}
\end{minipage}
\hfill
\begin{minipage}[t]{0.73\textwidth}
\vspace{0pt}
\raggedright
\textbf{Name:} {\large{$\mathbf{K7_{68}}$}} (chiral, non-rotatable$^{*}$) \\ \textbf{PD:} {\scriptsize\texttt{[0],[0,1,2,3],[1,4,5,6],[6,7,8,2],[3,9,10,4],[5],[7,11,12,8],[9,12,13,14],[14,13,11,10]}} \\ \textbf{EM:} {\scriptsize\texttt{(B0, A0C0D3E0, B1E3F0D0, C3G0G3B2, B3H0I3C1, C2, D1I2H1D2, E1G2I1I0, H3H2G1E2)}} \\ \textbf{Kauffman bracket:} {\scriptsize $-A^{26} - A^{24} + 2A^{22} + 2A^{20} - 2A^{18} - 3A^{16} + 2A^{14} + 4A^{12} - A^{10} - 3A^{8} + A^{6} + 2A^{4} - 1$} \\ \textbf{Arrow:} {\scriptsize $-A^{20} - A^{18}L_1 + 2A^{16} + 2A^{14}L_1 - 2A^{12} - 3A^{10}L_1 + 2A^{8} + 4A^{6}L_1 - A^{4} - 3A^{2}L_1 + 1 + 2L_1/A^{2} - L_1/A^{6}$} \\ \textbf{Mock:} {\scriptsize $2w^{4} + 2w^{3} - 4w^{2} - 6w + 3 + 6/w - 2/w^{3}$} \\ \textbf{Affine:} {\scriptsize $0$} \\ \textbf{Yamada:} {\scriptsize $-A^{27} + A^{26} + A^{25} - 3A^{24} + A^{23} + 4A^{22} - 4A^{21} - A^{20} + 6A^{19} - 2A^{18} - 2A^{17} + 5A^{16} - A^{15} - 3A^{14} - A^{12} - A^{11} - 6A^{10} + 3A^{9} + A^{8} - 7A^{7} + 4A^{6} + 2A^{5} - 4A^{4} + 2A^{3} - A + 1$}
\end{minipage}

\noindent{\color{gray!40}\rule{\textwidth}{0.4pt}}
\vspace{0.9\baselineskip}
\noindent \begin{minipage}[t]{0.25\textwidth}
\vspace{0pt}
\centering
\includegraphics[page=188,width=\linewidth]{knotoids.pdf}
\end{minipage}
\hfill
\begin{minipage}[t]{0.73\textwidth}
\vspace{0pt}
\raggedright
\textbf{Name:} {\large{$\mathbf{K7_{69}}$}} (chiral, non-rotatable$^{*}$) \\ \textbf{PD:} {\scriptsize\texttt{[0],[0,1,2,3],[1,4,5,2],[3,5,6,7],[4,8,9,10],[10,11,12,6],[7,13,14,8],[13,12,11,9],[14]}} \\ \textbf{EM:} {\scriptsize\texttt{(B0, A0C0C3D0, B1E0D1B2, B3C2F3G0, C1G3H3F0, E3H2H1D2, D3H0I0E1, G1F2F1E2, G2)}} \\ \textbf{Kauffman bracket:} {\scriptsize $A^{28} - 2A^{24} + 5A^{20} + 2A^{18} - 5A^{16} - 3A^{14} + 4A^{12} + 3A^{10} - 3A^{8} - 3A^{6} + A^{4} + A^{2}$} \\ \textbf{Arrow:} {\scriptsize $A^{-8} - 2/A^{12} + 5/A^{16} + 2L_1/A^{18} - 5/A^{20} - 3L_1/A^{22} + 4/A^{24} + 3L_1/A^{26} - 3/A^{28} - 3L_1/A^{30} + A^{-32} + L_1/A^{34}$} \\ \textbf{Mock:} {\scriptsize $3w^{4} + 2w^{3} - 6w^{2} - 4w + 7 + 4/w - 4/w^{2} - 2/w^{3} + w^{-4}$} \\ \textbf{Affine:} {\scriptsize $2t - 4 + 2/t$} \\ \textbf{Yamada:} {\scriptsize $-3A^{26} + 4A^{24} - 8A^{23} - 5A^{22} + 10A^{21} - 9A^{20} - 7A^{19} + 12A^{18} - 3A^{17} - 5A^{16} + 7A^{15} + 2A^{14} - 4A^{12} + 7A^{11} + 3A^{10} - 13A^{9} + 8A^{8} + 5A^{7} - 12A^{6} + 3A^{5} + 5A^{4} - 5A^{3} + A^{2} + 2A - 1$}
\end{minipage}

\noindent{\color{gray!40}\rule{\textwidth}{0.4pt}}
\vspace{0.9\baselineskip}
\noindent \begin{minipage}[t]{0.25\textwidth}
\vspace{0pt}
\centering
\includegraphics[page=189,width=\linewidth]{knotoids.pdf}
\end{minipage}
\hfill
\begin{minipage}[t]{0.73\textwidth}
\vspace{0pt}
\raggedright
\textbf{Name:} {\large{$\mathbf{K7_{70}}$}} (chiral, non-rotatable$^{*}$) \\ \textbf{PD:} {\scriptsize\texttt{[0],[0,1,2,3],[1,4,5,2],[3,5,6,7],[4,8,9,10],[10,11,12,6],[13,14,8,7],[13,12,11,9],[14]}} \\ \textbf{EM:} {\scriptsize\texttt{(B0, A0C0C3D0, B1E0D1B2, B3C2F3G3, C1G2H3F0, E3H2H1D2, H0I0E1D3, G0F2F1E2, G1)}} \\ \textbf{Kauffman bracket:} {\scriptsize $A^{22} + 2A^{20} - 3A^{16} + 3A^{12} + A^{10} - 3A^{8} - 2A^{6} + A^{4} + A^{2}$} \\ \textbf{Arrow:} {\scriptsize $A^{-8} + 2L_1/A^{10} - 3L_1/A^{14} + 3L_1/A^{18} + A^{-20} - 3L_1/A^{22} - 2/A^{24} + L_1/A^{26} + A^{-28}$} \\ \textbf{Mock:} {\scriptsize $w^{4} + 2w^{3} - 4w - 1 + 4/w + 2/w^{2} - 2/w^{3} - 1/w^{4}$} \\ \textbf{Affine:} {\scriptsize $2t - 4 + 2/t$} \\ \textbf{Yamada:} {\scriptsize $-2A^{24} - 2A^{23} + A^{22} - A^{21} - 6A^{20} + 2A^{19} + 2A^{18} - 5A^{17} + 2A^{16} + 4A^{15} - 2A^{14} + A^{13} + A^{12} + 2A^{11} - 3A^{10} - A^{9} + 5A^{8} - 5A^{7} - 2A^{6} + 4A^{5} - A^{4} - 2A^{3} + 2A^{2} + A - 1$}
\end{minipage}

\noindent{\color{gray!40}\rule{\textwidth}{0.4pt}}
\vspace{0.9\baselineskip}
\noindent \begin{minipage}[t]{0.25\textwidth}
\vspace{0pt}
\centering
\includegraphics[page=190,width=\linewidth]{knotoids.pdf}
\end{minipage}
\hfill
\begin{minipage}[t]{0.73\textwidth}
\vspace{0pt}
\raggedright
\textbf{Name:} {\large{$\mathbf{K7_{71}}$}} (chiral, non-rotatable$^{*}$) \\ \textbf{PD:} {\scriptsize\texttt{[0],[0,1,2,3],[1,4,5,2],[3,6,7,8],[4,9,10,11],[11,12,6,5],[12,13,14,7],[13,10,9,8],[14]}} \\ \textbf{EM:} {\scriptsize\texttt{(B0, A0C0C3D0, B1E0F3B2, B3F2G3H3, C1H2H1F0, E3G0D1C2, F1H0I0D2, G1E2E1D3, G2)}} \\ \textbf{Kauffman bracket:} {\scriptsize $-A^{26} + 3A^{22} + 2A^{20} - 4A^{18} - 4A^{16} + 4A^{14} + 5A^{12} - 2A^{10} - 5A^{8} + A^{6} + 3A^{4} - 1$} \\ \textbf{Arrow:} {\scriptsize $-A^{8} + 3A^{4} + 2A^{2}L_1 - 4 - 4L_1/A^{2} + 4/A^{4} + 5L_1/A^{6} - 2/A^{8} - 5L_1/A^{10} + A^{-12} + 3L_1/A^{14} - L_1/A^{18}$} \\ \textbf{Mock:} {\scriptsize $w^{6} + w^{5} - 3w^{4} - 3w^{3} + 5w^{2} + 6w - 3 - 6/w + 2/w^{2} + 3/w^{3} - 1/w^{4} - 1/w^{5}$} \\ \textbf{Affine:} {\scriptsize $2t - 4 + 2/t$} \\ \textbf{Yamada:} {\scriptsize $2A^{26} - A^{25} - 6A^{24} + 6A^{23} + 5A^{22} - 13A^{21} + 2A^{20} + 12A^{19} - 12A^{18} - 4A^{17} + 12A^{16} - 5A^{15} - 6A^{14} + 5A^{13} + 2A^{12} - 5A^{11} - 7A^{10} + 12A^{9} - A^{8} - 14A^{7} + 14A^{6} + 2A^{5} - 12A^{4} + 5A^{3} + 2A^{2} - 2A + 1$}
\end{minipage}

\noindent{\color{gray!40}\rule{\textwidth}{0.4pt}}
\vspace{0.9\baselineskip}
\noindent \begin{minipage}[t]{0.25\textwidth}
\vspace{0pt}
\centering
\includegraphics[page=191,width=\linewidth]{knotoids.pdf}
\end{minipage}
\hfill
\begin{minipage}[t]{0.73\textwidth}
\vspace{0pt}
\raggedright
\textbf{Name:} {\large{$\mathbf{K7_{72}}$}} (achiral, non-rotatable$^{*}$, Possible duplicate [7\_85]) \\ \textbf{PD:} {\scriptsize\texttt{[0],[0,1,2,3],[1,4,5,6],[7,8,3,2],[4,9,10,5],[6,11,12,7],[8,12,13,9],[14,13,11,10],[14]}} \\ \textbf{EM:} {\scriptsize\texttt{(B0, A0C0D3D2, B1E0E3F0, F3G0B3B2, C1G3H3C2, C3H2G1D0, D1F2H1E1, I0G2F1E2, H0)}} \\ \textbf{Kauffman bracket:} {\scriptsize $A^{28} - 3A^{24} - A^{22} + 5A^{20} + 4A^{18} - 5A^{16} - 5A^{14} + 3A^{12} + 5A^{10} - A^{8} - 3A^{6} + A^{2}$} \\ \textbf{Arrow:} {\scriptsize $A^{10}L_1 - 3A^{6}L_1 - A^{4} + 5A^{2}L_1 + 4 - 5L_1/A^{2} - 5/A^{4} + 3L_1/A^{6} + 5/A^{8} - L_1/A^{10} - 3/A^{12} + A^{-16}$} \\ \textbf{Mock:} {\scriptsize $w^{5} + w^{4} - 3w^{3} - 3w^{2} + 5w + 6 - 5/w - 5/w^{2} + 3/w^{3} + 3/w^{4} - 1/w^{5} - 1/w^{6}$} \\ \textbf{Affine:} {\scriptsize $t - 2 + 1/t$} \\ \textbf{Yamada:} {\scriptsize $-A^{26} + 3A^{25} + 2A^{24} - 10A^{23} + 4A^{22} + 13A^{21} - 16A^{20} - 2A^{19} + 18A^{18} - 11A^{17} - 5A^{16} + 12A^{15} - 4A^{13} - 2A^{12} + 11A^{11} - 4A^{10} - 13A^{9} + 17A^{8} + A^{7} - 17A^{6} + 11A^{5} + 6A^{4} - 10A^{3} + A^{2} + 3A - 1$}
\end{minipage}

\noindent{\color{gray!40}\rule{\textwidth}{0.4pt}}
\vspace{0.9\baselineskip}
\noindent \begin{minipage}[t]{0.25\textwidth}
\vspace{0pt}
\centering
\includegraphics[page=192,width=\linewidth]{knotoids.pdf}
\end{minipage}
\hfill
\begin{minipage}[t]{0.73\textwidth}
\vspace{0pt}
\raggedright
\textbf{Name:} {\large{$\mathbf{K7_{73}}$}} (chiral, non-rotatable$^{*}$) \\ \textbf{PD:} {\scriptsize\texttt{[0],[0,1,2,3],[1,4,5,6],[7,8,9,2],[3,9,10,4],[5],[6,11,12,7],[8,13,14,10],[11,14,13,12]}} \\ \textbf{EM:} {\scriptsize\texttt{(B0, A0C0D3E0, B1E3F0G0, G3H0E1B2, B3D2H3C1, C2, C3I0I3D0, D1I2I1E2, G1H2H1G2)}} \\ \textbf{Kauffman bracket:} {\scriptsize $-A^{26} - A^{24} + 2A^{22} + 2A^{20} - 2A^{18} - 3A^{16} + 3A^{14} + 4A^{12} - 2A^{10} - 4A^{8} + A^{6} + 3A^{4} - 1$} \\ \textbf{Arrow:} {\scriptsize $-A^{20} - A^{18}L_1 + 2A^{16} + 2A^{14}L_1 - 2A^{12} - 3A^{10}L_1 + 3A^{8} + 4A^{6}L_1 - 2A^{4} - 4A^{2}L_1 + 1 + 3L_1/A^{2} - L_1/A^{6}$} \\ \textbf{Mock:} {\scriptsize $2w^{4} + 2w^{3} - 5w^{2} - 7w + 4 + 7/w - 2/w^{3}$} \\ \textbf{Affine:} {\scriptsize $-t + 2 - 1/t$} \\ \textbf{Yamada:} {\scriptsize $-A^{27} + A^{26} + A^{25} - 3A^{24} + 2A^{23} + 4A^{22} - 6A^{21} + A^{20} + 7A^{19} - 6A^{18} - A^{17} + 6A^{16} - 4A^{15} - 3A^{14} + 2A^{13} - A^{11} - 5A^{10} + 5A^{9} - A^{8} - 9A^{7} + 7A^{6} - 6A^{4} + 5A^{3} - 2A + 1$}
\end{minipage}

\noindent{\color{gray!40}\rule{\textwidth}{0.4pt}}
\vspace{0.9\baselineskip}
\noindent \begin{minipage}[t]{0.25\textwidth}
\vspace{0pt}
\centering
\includegraphics[page=193,width=\linewidth]{knotoids.pdf}
\end{minipage}
\hfill
\begin{minipage}[t]{0.73\textwidth}
\vspace{0pt}
\raggedright
\textbf{Name:} {\large{$\mathbf{K7_{74}}$}} (achiral, non-rotatable$^{*}$, Possible duplicate [7\_79]) \\ \textbf{PD:} {\scriptsize\texttt{[0],[0,1,2,3],[1,4,5,2],[3,6,7,8],[4,9,10,11],[11,12,6,5],[12,10,13,7],[8,13,14,9],[14]}} \\ \textbf{EM:} {\scriptsize\texttt{(B0, A0C0C3D0, B1E0F3B2, B3F2G3H0, C1H3G1F0, E3G0D1C2, F1E2H1D2, D3G2I0E1, H2)}} \\ \textbf{Kauffman bracket:} {\scriptsize $A^{28} - 3A^{24} + 5A^{20} + 2A^{18} - 6A^{16} - 3A^{14} + 5A^{12} + 4A^{10} - 3A^{8} - 3A^{6} + A^{4} + A^{2}$} \\ \textbf{Arrow:} {\scriptsize $A^{10}L_1 - 3A^{6}L_1 + 5A^{2}L_1 + 2 - 6L_1/A^{2} - 3/A^{4} + 5L_1/A^{6} + 4/A^{8} - 3L_1/A^{10} - 3/A^{12} + L_1/A^{14} + A^{-16}$} \\ \textbf{Mock:} {\scriptsize $w^{5} - 4w^{3} - w^{2} + 7w + 4 - 7/w - 4/w^{2} + 4/w^{3} + 3/w^{4} - 1/w^{5} - 1/w^{6}$} \\ \textbf{Affine:} {\scriptsize $0$} \\ \textbf{Yamada:} {\scriptsize $-3A^{26} + 2A^{25} + 7A^{24} - 9A^{23} - 3A^{22} + 15A^{21} - 10A^{20} - 9A^{19} + 15A^{18} - 5A^{17} - 9A^{16} + 6A^{15} + A^{14} - 3A^{13} - 7A^{12} + 9A^{11} + 3A^{10} - 16A^{9} + 11A^{8} + 8A^{7} - 15A^{6} + 4A^{5} + 7A^{4} - 6A^{3} + 2A - 1$}
\end{minipage}

\noindent{\color{gray!40}\rule{\textwidth}{0.4pt}}
\vspace{0.9\baselineskip}
\noindent \begin{minipage}[t]{0.25\textwidth}
\vspace{0pt}
\centering
\includegraphics[page=194,width=\linewidth]{knotoids.pdf}
\end{minipage}
\hfill
\begin{minipage}[t]{0.73\textwidth}
\vspace{0pt}
\raggedright
\textbf{Name:} {\large{$\mathbf{K7_{75}}$}} (achiral, non-rotatable$^{*}$, Possible duplicate [7\_51]) \\ \textbf{PD:} {\scriptsize\texttt{[0],[0,1,2,3],[1,4,5,2],[3,6,7,8],[4,9,10,11],[11,12,6,5],[12,10,13,7],[8,14,13,9],[14]}} \\ \textbf{EM:} {\scriptsize\texttt{(B0, A0C0C3D0, B1E0F3B2, B3F2G3H0, C1H3G1F0, E3G0D1C2, F1E2H2D2, D3I0G2E1, H1)}} \\ \textbf{Kauffman bracket:} {\scriptsize $-A^{20} - A^{18} + A^{16} + 3A^{14} - 3A^{10} - A^{8} + 3A^{6} + 2A^{4} - A^{2} - 1$} \\ \textbf{Arrow:} {\scriptsize $-A^{14}L_1 - A^{12} + A^{10}L_1 + 3A^{8} - 3A^{4} - A^{2}L_1 + 3 + 2L_1/A^{2} - 1/A^{4} - L_1/A^{6}$} \\ \textbf{Mock:} {\scriptsize $w^{4} - 3w^{2} - w + 4 + 1/w - 2/w^{2} + w^{-4}$} \\ \textbf{Affine:} {\scriptsize $-t + 2 - 1/t$} \\ \textbf{Yamada:} {\scriptsize $A^{24} - 2A^{22} - A^{21} + 3A^{20} - A^{19} - 5A^{18} + 4A^{17} + 2A^{16} - 5A^{15} + 3A^{14} + 3A^{13} - 2A^{12} + 2A^{11} + 2A^{10} + 3A^{9} - 3A^{8} + A^{7} + 5A^{6} - 5A^{5} - A^{4} + 5A^{3} - 2A^{2} - 2A + 1$}
\end{minipage}

\noindent{\color{gray!40}\rule{\textwidth}{0.4pt}}
\vspace{0.9\baselineskip}
\noindent \begin{minipage}[t]{0.25\textwidth}
\vspace{0pt}
\centering
\includegraphics[page=195,width=\linewidth]{knotoids.pdf}
\end{minipage}
\hfill
\begin{minipage}[t]{0.73\textwidth}
\vspace{0pt}
\raggedright
\textbf{Name:} {\large{$\mathbf{K7_{76}}$}} (achiral, non-rotatable$^{*}$, Possible duplicate [7\_63]) \\ \textbf{PD:} {\scriptsize\texttt{[0],[0,1,2,3],[1,4,5,2],[3,6,7,8],[4,9,10,11],[11,12,6,5],[12,10,13,7],[13,14,9,8],[14]}} \\ \textbf{EM:} {\scriptsize\texttt{(B0, A0C0C3D0, B1E0F3B2, B3F2G3H3, C1H2G1F0, E3G0D1C2, F1E2H0D2, G2I0E1D3, H1)}} \\ \textbf{Kauffman bracket:} {\scriptsize $A^{22} + 2A^{20} - A^{18} - 3A^{16} + 4A^{12} + A^{10} - 3A^{8} - 2A^{6} + A^{4} + A^{2}$} \\ \textbf{Arrow:} {\scriptsize $L_1/A^{2} + 2/A^{4} - L_1/A^{6} - 3/A^{8} + 4/A^{12} + L_1/A^{14} - 3/A^{16} - 2L_1/A^{18} + A^{-20} + L_1/A^{22}$} \\ \textbf{Mock:} {\scriptsize $-w^{4} + 4w^{2} + w - 4 - 1/w + 3/w^{2} - 1/w^{4}$} \\ \textbf{Affine:} {\scriptsize $t - 2 + 1/t$} \\ \textbf{Yamada:} {\scriptsize $-2A^{24} - A^{23} + 3A^{22} - A^{21} - 6A^{20} + 4A^{19} + 2A^{18} - 7A^{17} + 2A^{16} + 3A^{15} - 4A^{14} + A^{12} + 2A^{11} - 4A^{10} + 6A^{8} - 6A^{7} - A^{6} + 6A^{5} - 2A^{4} - 3A^{3} + 2A^{2} + A - 1$}
\end{minipage}

\noindent{\color{gray!40}\rule{\textwidth}{0.4pt}}
\vspace{0.9\baselineskip}
\noindent \begin{minipage}[t]{0.25\textwidth}
\vspace{0pt}
\centering
\includegraphics[page=196,width=\linewidth]{knotoids.pdf}
\end{minipage}
\hfill
\begin{minipage}[t]{0.73\textwidth}
\vspace{0pt}
\raggedright
\textbf{Name:} {\large{$\mathbf{K7_{77}}$}} (achiral, non-rotatable$^{*}$, Possible duplicate [7\_78]) \\ \textbf{PD:} {\scriptsize\texttt{[0],[0,1,2,3],[1,4,5,2],[3,6,7,8],[4,9,10,11],[11,12,6,5],[12,10,13,7],[14,13,9,8],[14]}} \\ \textbf{EM:} {\scriptsize\texttt{(B0, A0C0C3D0, B1E0F3B2, B3F2G3H3, C1H2G1F0, E3G0D1C2, F1E2H1D2, I0G2E1D3, H0)}} \\ \textbf{Kauffman bracket:} {\scriptsize $-A^{26} - A^{24} + 3A^{22} + 3A^{20} - 3A^{18} - 5A^{16} + 3A^{14} + 6A^{12} - A^{10} - 5A^{8} + 3A^{4} - 1$} \\ \textbf{Arrow:} {\scriptsize $-A^{8} - A^{6}L_1 + 3A^{4} + 3A^{2}L_1 - 3 - 5L_1/A^{2} + 3/A^{4} + 6L_1/A^{6} - 1/A^{8} - 5L_1/A^{10} + 3L_1/A^{14} - L_1/A^{18}$} \\ \textbf{Mock:} {\scriptsize $w^{6} + w^{5} - 3w^{4} - 4w^{3} + 4w^{2} + 7w - 2 - 7/w + w^{-2} + 4/w^{3} - 1/w^{5}$} \\ \textbf{Affine:} {\scriptsize $0$} \\ \textbf{Yamada:} {\scriptsize $A^{26} - 2A^{25} - 2A^{24} + 6A^{23} - 2A^{22} - 9A^{21} + 10A^{20} + 5A^{19} - 14A^{18} + 7A^{17} + 10A^{16} - 10A^{15} + A^{14} + 8A^{13} - 2A^{12} - 4A^{11} + A^{10} + 10A^{9} - 9A^{8} - 6A^{7} + 16A^{6} - 8A^{5} - 8A^{4} + 11A^{3} - 2A^{2} - 4A + 2$}
\end{minipage}

\noindent{\color{gray!40}\rule{\textwidth}{0.4pt}}
\vspace{0.9\baselineskip}
\noindent \begin{minipage}[t]{0.25\textwidth}
\vspace{0pt}
\centering
\includegraphics[page=197,width=\linewidth]{knotoids.pdf}
\end{minipage}
\hfill
\begin{minipage}[t]{0.73\textwidth}
\vspace{0pt}
\raggedright
\textbf{Name:} {\large{$\mathbf{K7_{78}}$}} (achiral, non-rotatable$^{*}$, Possible duplicate [7\_77]) \\ \textbf{PD:} {\scriptsize\texttt{[0],[0,1,2,3],[1,4,5,6],[7,8,3,2],[9,10,11,4],[10,12,13,5],[6,13,14,7],[8,14,12,9],[11]}} \\ \textbf{EM:} {\scriptsize\texttt{(B0, A0C0D3D2, B1E3F3G0, G3H0B3B2, H3F0I0C1, E1H2G1C2, C3F2H1D0, D1G2F1E0, E2)}} \\ \textbf{Kauffman bracket:} {\scriptsize $A^{27} - 3A^{23} + 5A^{19} + A^{17} - 6A^{15} - 3A^{13} + 5A^{11} + 3A^{9} - 3A^{7} - 3A^{5} + A^{3} + A$} \\ \textbf{Arrow:} {\scriptsize $-A^{18}L_1 + 3A^{14}L_1 - 5A^{10}L_1 - A^{8} + 6A^{6}L_1 + 3A^{4} - 5A^{2}L_1 - 3 + 3L_1/A^{2} + 3/A^{4} - L_1/A^{6} - 1/A^{8}$} \\ \textbf{Mock:} {\scriptsize $-w^{5} + 4w^{3} + w^{2} - 7w - 2 + 7/w + 4/w^{2} - 4/w^{3} - 3/w^{4} + w^{-5} + w^{-6}$} \\ \textbf{Affine:} {\scriptsize $0$} \\ \textbf{Yamada:} {\scriptsize $2A^{26} - 4A^{25} - 2A^{24} + 11A^{23} - 8A^{22} - 8A^{21} + 16A^{20} - 6A^{19} - 9A^{18} + 10A^{17} + A^{16} - 4A^{15} - 2A^{14} + 8A^{13} + A^{12} - 10A^{11} + 10A^{10} + 7A^{9} - 14A^{8} + 5A^{7} + 10A^{6} - 9A^{5} - 2A^{4} + 6A^{3} - 2A^{2} - 2A + 1$}
\end{minipage}

\noindent{\color{gray!40}\rule{\textwidth}{0.4pt}}
\vspace{0.9\baselineskip}
\noindent \begin{minipage}[t]{0.25\textwidth}
\vspace{0pt}
\centering
\includegraphics[page=198,width=\linewidth]{knotoids.pdf}
\end{minipage}
\hfill
\begin{minipage}[t]{0.73\textwidth}
\vspace{0pt}
\raggedright
\textbf{Name:} {\large{$\mathbf{K7_{79}}$}} (achiral, non-rotatable$^{*}$, Possible duplicate [7\_74]) \\ \textbf{PD:} {\scriptsize\texttt{[0],[0,1,2,3],[1,4,5,6],[6,7,8,2],[3,8,9,10],[4,11,12,5],[7,12,13,9],[10,13,14,11],[14]}} \\ \textbf{EM:} {\scriptsize\texttt{(B0, A0C0D3E0, B1F0F3D0, C3G0E1B2, B3D2G3H0, C1H3G1C2, D1F2H1E2, E3G2I0F1, H2)}} \\ \textbf{Kauffman bracket:} {\scriptsize $A^{28} - 3A^{24} + 5A^{20} + 2A^{18} - 6A^{16} - 3A^{14} + 5A^{12} + 4A^{10} - 3A^{8} - 3A^{6} + A^{4} + A^{2}$} \\ \textbf{Arrow:} {\scriptsize $A^{10}L_1 - 3A^{6}L_1 + 5A^{2}L_1 + 2 - 6L_1/A^{2} - 3/A^{4} + 5L_1/A^{6} + 4/A^{8} - 3L_1/A^{10} - 3/A^{12} + L_1/A^{14} + A^{-16}$} \\ \textbf{Mock:} {\scriptsize $w^{5} - 4w^{3} - w^{2} + 7w + 4 - 7/w - 4/w^{2} + 4/w^{3} + 3/w^{4} - 1/w^{5} - 1/w^{6}$} \\ \textbf{Affine:} {\scriptsize $0$} \\ \textbf{Yamada:} {\scriptsize $-3A^{26} + 2A^{25} + 7A^{24} - 9A^{23} - 3A^{22} + 15A^{21} - 10A^{20} - 9A^{19} + 15A^{18} - 5A^{17} - 9A^{16} + 6A^{15} + A^{14} - 3A^{13} - 7A^{12} + 9A^{11} + 3A^{10} - 16A^{9} + 11A^{8} + 8A^{7} - 15A^{6} + 4A^{5} + 7A^{4} - 6A^{3} + 2A - 1$}
\end{minipage}

\noindent{\color{gray!40}\rule{\textwidth}{0.4pt}}
\vspace{0.9\baselineskip}
\noindent \begin{minipage}[t]{0.25\textwidth}
\vspace{0pt}
\centering
\includegraphics[page=199,width=\linewidth]{knotoids.pdf}
\end{minipage}
\hfill
\begin{minipage}[t]{0.73\textwidth}
\vspace{0pt}
\raggedright
\textbf{Name:} {\large{$\mathbf{K7_{80}}$}} (achiral, non-rotatable$^{*}$, Possible duplicate [7\_64]) \\ \textbf{PD:} {\scriptsize\texttt{[0],[0,1,2,3],[1,4,5,2],[3,6,7,8],[4,8,9,10],[11,12,6,5],[12,11,13,7],[14,13,10,9],[14]}} \\ \textbf{EM:} {\scriptsize\texttt{(B0, A0C0C3D0, B1E0F3B2, B3F2G3E1, C1D3H3H2, G1G0D1C2, F1F0H1D2, I0G2E3E2, H0)}} \\ \textbf{Kauffman bracket:} {\scriptsize $A^{28} - 3A^{24} - 2A^{22} + 4A^{20} + 4A^{18} - 4A^{16} - 5A^{14} + 2A^{12} + 5A^{10} - 2A^{6} + A^{2}$} \\ \textbf{Arrow:} {\scriptsize $A^{34}L_1 - 3A^{30}L_1 - 2A^{28} + 4A^{26}L_1 + 4A^{24} - 4A^{22}L_1 - 5A^{20} + 2A^{18}L_1 + 5A^{16} - 2A^{12} + A^{8}$} \\ \textbf{Mock:} {\scriptsize $w^{4} - 2w^{3} - 3w^{2} + 5w + 6 - 5/w - 6/w^{2} + 2/w^{3} + 3/w^{4}$} \\ \textbf{Affine:} {\scriptsize $-t + 2 - 1/t$} \\ \textbf{Yamada:} {\scriptsize $2A^{25} - 2A^{24} - 5A^{23} + 8A^{22} - 13A^{20} + 10A^{19} + 7A^{18} - 15A^{17} + 6A^{16} + 8A^{15} - 6A^{14} + 4A^{12} + 6A^{11} - 9A^{10} - A^{9} + 13A^{8} - 12A^{7} - 7A^{6} + 13A^{5} - 7A^{4} - 8A^{3} + 5A^{2} - 3$}
\end{minipage}

\noindent{\color{gray!40}\rule{\textwidth}{0.4pt}}
\vspace{0.9\baselineskip}
\noindent \begin{minipage}[t]{0.25\textwidth}
\vspace{0pt}
\centering
\includegraphics[page=200,width=\linewidth]{knotoids.pdf}
\end{minipage}
\hfill
\begin{minipage}[t]{0.73\textwidth}
\vspace{0pt}
\raggedright
\textbf{Name:} {\large{$\mathbf{K7_{81}}$}} (achiral, non-rotatable$^{*}$, Possible duplicate [7\_65]) \\ \textbf{PD:} {\scriptsize\texttt{[0],[0,1,2,3],[1,4,5,2],[3,6,7,8],[8,9,10,4],[11,12,6,5],[12,11,13,7],[9,14,13,10],[14]}} \\ \textbf{EM:} {\scriptsize\texttt{(B0, A0C0C3D0, B1E3F3B2, B3F2G3E0, D3H0H3C1, G1G0D1C2, F1F0H2D2, E1I0G2E2, H1)}} \\ \textbf{Kauffman bracket:} {\scriptsize $-A^{26} - A^{24} + A^{22} + 3A^{20} - 4A^{16} - A^{14} + 4A^{12} + 3A^{10} - 2A^{8} - 2A^{6} + A^{2}$} \\ \textbf{Arrow:} {\scriptsize $-A^{20} - A^{18}L_1 + A^{16} + 3A^{14}L_1 - 4A^{10}L_1 - A^{8} + 4A^{6}L_1 + 3A^{4} - 2A^{2}L_1 - 2 + A^{-4}$} \\ \textbf{Mock:} {\scriptsize $w^{4} + 2w^{3} - 5w - 2 + 5/w + 3/w^{2} - 2/w^{3} - 1/w^{4}$} \\ \textbf{Affine:} {\scriptsize $t - 2 + 1/t$} \\ \textbf{Yamada:} {\scriptsize $A^{26} - 2A^{24} + 4A^{22} - 3A^{21} - 6A^{20} + 8A^{19} + A^{18} - 10A^{17} + 6A^{16} + 3A^{15} - 6A^{14} + 2A^{13} + 3A^{12} + 2A^{11} - 5A^{10} + 4A^{9} + 7A^{8} - 8A^{7} + A^{6} + 9A^{5} - 5A^{4} - 3A^{3} + 5A^{2} - 2$}
\end{minipage}

\noindent{\color{gray!40}\rule{\textwidth}{0.4pt}}
\vspace{0.9\baselineskip}
\noindent \begin{minipage}[t]{0.25\textwidth}
\vspace{0pt}
\centering
\includegraphics[page=201,width=\linewidth]{knotoids.pdf}
\end{minipage}
\hfill
\begin{minipage}[t]{0.73\textwidth}
\vspace{0pt}
\raggedright
\textbf{Name:} {\large{$\mathbf{K7_{82}}$}} (chiral, non-rotatable$^{*}$) \\ \textbf{PD:} {\scriptsize\texttt{[0],[0,1,2,3],[1,4,5,2],[3,6,7,8],[4,8,9,10],[10,11,12,5],[6,12,13,7],[14,13,11,9],[14]}} \\ \textbf{EM:} {\scriptsize\texttt{(B0, A0C0C3D0, B1E0F3B2, B3G0G3E1, C1D3H3F0, E3H2G1C2, D1F2H1D2, I0G2F1E2, H0)}} \\ \textbf{Kauffman bracket:} {\scriptsize $-A^{26} + 3A^{22} + 2A^{20} - 4A^{18} - 4A^{16} + 4A^{14} + 6A^{12} - 2A^{10} - 5A^{8} + 3A^{4} - 1$} \\ \textbf{Arrow:} {\scriptsize $-A^{2}L_1 + 3L_1/A^{2} + 2/A^{4} - 4L_1/A^{6} - 4/A^{8} + 4L_1/A^{10} + 6/A^{12} - 2L_1/A^{14} - 5/A^{16} + 3/A^{20} - 1/A^{24}$} \\ \textbf{Mock:} {\scriptsize $-w^{5} - w^{4} + 3w^{3} + 5w^{2} - 3w - 6 + 3/w + 5/w^{2} - 3/w^{3} - 3/w^{4} + w^{-5} + w^{-6}$} \\ \textbf{Affine:} {\scriptsize $t - 2 + 1/t$} \\ \textbf{Yamada:} {\scriptsize $-2A^{25} + A^{24} + 6A^{23} - 6A^{22} - 6A^{21} + 14A^{20} - 14A^{18} + 13A^{17} + 7A^{16} - 13A^{15} + 5A^{14} + 8A^{13} - 6A^{12} - 4A^{11} + 4A^{10} + 7A^{9} - 13A^{8} + 17A^{6} - 14A^{5} - 5A^{4} + 13A^{3} - 5A^{2} - 4A + 3$}
\end{minipage}

\noindent{\color{gray!40}\rule{\textwidth}{0.4pt}}
\vspace{0.9\baselineskip}
\noindent \begin{minipage}[t]{0.25\textwidth}
\vspace{0pt}
\centering
\includegraphics[page=202,width=\linewidth]{knotoids.pdf}
\end{minipage}
\hfill
\begin{minipage}[t]{0.73\textwidth}
\vspace{0pt}
\raggedright
\textbf{Name:} {\large{$\mathbf{K7_{83}}$}} (chiral, non-rotatable$^{*}$) \\ \textbf{PD:} {\scriptsize\texttt{[0],[0,1,2,3],[1,4,5,2],[3,6,7,8],[8,9,10,4],[5,10,11,12],[6,12,13,7],[9,14,13,11],[14]}} \\ \textbf{EM:} {\scriptsize\texttt{(B0, A0C0C3D0, B1E3F0B2, B3G0G3E0, D3H0F1C1, C2E2H3G1, D1F3H2D2, E1I0G2F2, H1)}} \\ \textbf{Kauffman bracket:} {\scriptsize $A^{22} + 2A^{20} - 3A^{16} - A^{14} + 3A^{12} + 2A^{10} - 2A^{8} - 2A^{6} + A^{2}$} \\ \textbf{Arrow:} {\scriptsize $A^{-8} + 2L_1/A^{10} - 3L_1/A^{14} - 1/A^{16} + 3L_1/A^{18} + 2/A^{20} - 2L_1/A^{22} - 2/A^{24} + A^{-28}$} \\ \textbf{Mock:} {\scriptsize $w^{4} + 2w^{3} + w^{2} - 3w - 2 + 3/w + 2/w^{2} - 2/w^{3} - 1/w^{4}$} \\ \textbf{Affine:} {\scriptsize $3t - 6 + 3/t$} \\ \textbf{Yamada:} {\scriptsize $2A^{23} + 2A^{22} - A^{21} - A^{20} + 6A^{19} + A^{18} - 4A^{17} + 5A^{16} + A^{15} - 5A^{14} + A^{13} - 3A^{10} + 2A^{9} + 2A^{8} - 6A^{7} + 2A^{6} + 4A^{5} - 3A^{4} - A^{3} + 4A^{2} - 2$}
\end{minipage}

\noindent{\color{gray!40}\rule{\textwidth}{0.4pt}}
\vspace{0.9\baselineskip}
\noindent \begin{minipage}[t]{0.25\textwidth}
\vspace{0pt}
\centering
\includegraphics[page=203,width=\linewidth]{knotoids.pdf}
\end{minipage}
\hfill
\begin{minipage}[t]{0.73\textwidth}
\vspace{0pt}
\raggedright
\textbf{Name:} {\large{$\mathbf{K7_{84}}$}} (chiral, non-rotatable$^{*}$) \\ \textbf{PD:} {\scriptsize\texttt{[0],[0,1,2,3],[1,4,5,2],[3,6,7,8],[4,8,9,10],[11,12,6,5],[7],[12,13,14,9],[10,14,13,11]}} \\ \textbf{EM:} {\scriptsize\texttt{(B0, A0C0C3D0, B1E0F3B2, B3F2G0E1, C1D3H3I0, I3H0D1C2, D2, F1I2I1E2, E3H2H1F0)}} \\ \textbf{Kauffman bracket:} {\scriptsize $A^{28} - 2A^{24} - A^{22} + 4A^{20} + 2A^{18} - 4A^{16} - 3A^{14} + 3A^{12} + 3A^{10} - 2A^{8} - 2A^{6} + A^{4} + A^{2}$} \\ \textbf{Arrow:} {\scriptsize $A^{16} - 2A^{12} - A^{10}L_1 + 4A^{8} + 2A^{6}L_1 - 4A^{4} - 3A^{2}L_1 + 3 + 3L_1/A^{2} - 2/A^{4} - 2L_1/A^{6} + A^{-8} + L_1/A^{10}$} \\ \textbf{Mock:} {\scriptsize $w^{4} + w^{3} - 6w^{2} - 5w + 8 + 5/w - 2/w^{2} - 1/w^{3}$} \\ \textbf{Affine:} {\scriptsize $-2t + 4 - 2/t$} \\ \textbf{Yamada:} {\scriptsize $A^{27} - A^{26} - A^{25} + 4A^{24} - 3A^{23} - 5A^{22} + 8A^{21} - 3A^{20} - 7A^{19} + 8A^{18} - A^{17} - 4A^{16} + 3A^{15} - 5A^{12} + 2A^{11} + 4A^{10} - 9A^{9} + 3A^{8} + 6A^{7} - 8A^{6} + A^{5} + 4A^{4} - 3A^{3} + A - 1$}
\end{minipage}

\noindent{\color{gray!40}\rule{\textwidth}{0.4pt}}
\vspace{0.9\baselineskip}
\noindent \begin{minipage}[t]{0.25\textwidth}
\vspace{0pt}
\centering
\includegraphics[page=204,width=\linewidth]{knotoids.pdf}
\end{minipage}
\hfill
\begin{minipage}[t]{0.73\textwidth}
\vspace{0pt}
\raggedright
\textbf{Name:} {\large{$\mathbf{K7_{85}}$}} (achiral, non-rotatable$^{*}$, Possible duplicate [7\_72]) \\ \textbf{PD:} {\scriptsize\texttt{[0],[0,1,2,3],[1,4,5,2],[3,6,7,8],[4,8,9,10],[11,12,6,5],[12,13,14,7],[13,11,10,9],[14]}} \\ \textbf{EM:} {\scriptsize\texttt{(B0, A0C0C3D0, B1E0F3B2, B3F2G3E1, C1D3H3H2, H1G0D1C2, F1H0I0D2, G1F0E3E2, G2)}} \\ \textbf{Kauffman bracket:} {\scriptsize $A^{28} - 3A^{24} - A^{22} + 5A^{20} + 4A^{18} - 5A^{16} - 5A^{14} + 3A^{12} + 5A^{10} - A^{8} - 3A^{6} + A^{2}$} \\ \textbf{Arrow:} {\scriptsize $A^{10}L_1 - 3A^{6}L_1 - A^{4} + 5A^{2}L_1 + 4 - 5L_1/A^{2} - 5/A^{4} + 3L_1/A^{6} + 5/A^{8} - L_1/A^{10} - 3/A^{12} + A^{-16}$} \\ \textbf{Mock:} {\scriptsize $w^{5} + w^{4} - 3w^{3} - 3w^{2} + 5w + 6 - 5/w - 5/w^{2} + 3/w^{3} + 3/w^{4} - 1/w^{5} - 1/w^{6}$} \\ \textbf{Affine:} {\scriptsize $t - 2 + 1/t$} \\ \textbf{Yamada:} {\scriptsize $-A^{26} + 3A^{25} + 2A^{24} - 10A^{23} + 4A^{22} + 13A^{21} - 16A^{20} - 2A^{19} + 18A^{18} - 11A^{17} - 5A^{16} + 12A^{15} - 4A^{13} - 2A^{12} + 11A^{11} - 4A^{10} - 13A^{9} + 17A^{8} + A^{7} - 17A^{6} + 11A^{5} + 6A^{4} - 10A^{3} + A^{2} + 3A - 1$}
\end{minipage}

\noindent{\color{gray!40}\rule{\textwidth}{0.4pt}}
\vspace{0.9\baselineskip}
\noindent \begin{minipage}[t]{0.25\textwidth}
\vspace{0pt}
\centering
\includegraphics[page=205,width=\linewidth]{knotoids.pdf}
\end{minipage}
\hfill
\begin{minipage}[t]{0.73\textwidth}
\vspace{0pt}
\raggedright
\textbf{Name:} {\large{$\mathbf{K7_{86}}$}} (chiral, non-rotatable$^{*}$) \\ \textbf{PD:} {\scriptsize\texttt{[0],[0,1,2,3],[1,4,5,6],[6,7,8,2],[3,9,10,11],[4,11,12,5],[7,12,13,14],[14,13,9,8],[10]}} \\ \textbf{EM:} {\scriptsize\texttt{(B0, A0C0D3E0, B1F0F3D0, C3G0H3B2, B3H2I0F1, C1E3G1C2, D1F2H1H0, G3G2E1D2, E2)}} \\ \textbf{Kauffman bracket:} {\scriptsize $A^{28} - 3A^{24} - A^{22} + 4A^{20} + 2A^{18} - 5A^{16} - 3A^{14} + 4A^{12} + 4A^{10} - 2A^{8} - 2A^{6} + A^{4} + A^{2}$} \\ \textbf{Arrow:} {\scriptsize $A^{34}L_1 - 3A^{30}L_1 - A^{28} + 4A^{26}L_1 + 2A^{24} - 5A^{22}L_1 - 3A^{20} + 4A^{18}L_1 + 4A^{16} - 2A^{14}L_1 - 2A^{12} + A^{10}L_1 + A^{8}$} \\ \textbf{Mock:} {\scriptsize $-3w^{3} - w^{2} + 7w + 4 - 7/w - 5/w^{2} + 3/w^{3} + 3/w^{4}$} \\ \textbf{Affine:} {\scriptsize $-2t + 4 - 2/t$} \\ \textbf{Yamada:} {\scriptsize $A^{27} - 2A^{26} - A^{25} + 6A^{24} - 5A^{23} - 5A^{22} + 11A^{21} - 4A^{20} - 9A^{19} + 11A^{18} - A^{17} - 5A^{16} + 5A^{15} + 2A^{14} + A^{13} - 6A^{12} + 5A^{11} + 5A^{10} - 12A^{9} + 4A^{8} + 7A^{7} - 12A^{6} + 4A^{4} - 5A^{3} - A^{2} + A - 1$}
\end{minipage}

\noindent{\color{gray!40}\rule{\textwidth}{0.4pt}}
\vspace{0.9\baselineskip}
\noindent \begin{minipage}[t]{0.25\textwidth}
\vspace{0pt}
\centering
\includegraphics[page=206,width=\linewidth]{knotoids.pdf}
\end{minipage}
\hfill
\begin{minipage}[t]{0.73\textwidth}
\vspace{0pt}
\raggedright
\textbf{Name:} {\large{$\mathbf{K7_{87}}$}} (chiral, non-rotatable$^{*}$) \\ \textbf{PD:} {\scriptsize\texttt{[0],[0,1,2,3],[1,4,5,6],[7,8,9,2],[3,9,10,4],[5],[6,10,11,12],[13,14,8,7],[14,13,12,11]}} \\ \textbf{EM:} {\scriptsize\texttt{(B0, A0C0D3E0, B1E3F0G0, H3H2E1B2, B3D2G1C1, C2, C3E2I3I2, I1I0D1D0, H1H0G3G2)}} \\ \textbf{Kauffman bracket:} {\scriptsize $A^{27} + A^{25} - 2A^{23} - 2A^{21} + 2A^{19} + 3A^{17} - 2A^{15} - 3A^{13} + A^{11} + 2A^{9} - A^{7} - 2A^{5} + A$} \\ \textbf{Arrow:} {\scriptsize $-A^{24} - A^{22}L_1 + 2A^{20} + 2A^{18}L_1 - 2A^{16} - 3A^{14}L_1 + 2A^{12} + 3A^{10}L_1 - A^{8} - 2A^{6}L_1 + A^{4} + 2A^{2}L_1 - L_1/A^{2}$} \\ \textbf{Mock:} {\scriptsize $w^{6} + w^{5} - 2w^{4} - 3w^{3} + 2w^{2} + 3w - 2 - 3/w + 2/w^{2} + 3/w^{3} - 1/w^{5}$} \\ \textbf{Affine:} {\scriptsize $-t + 2 - 1/t$} \\ \textbf{Yamada:} {\scriptsize $-A^{27} + A^{26} - 2A^{24} + 3A^{23} - 3A^{21} + 5A^{20} - A^{19} - 3A^{18} + 4A^{17} - A^{14} + A^{13} - 5A^{11} + A^{9} - 5A^{8} + 2A^{6} - 3A^{5} - A^{4} + 3A^{3} - A^{2} - A + 1$}
\end{minipage}

\noindent{\color{gray!40}\rule{\textwidth}{0.4pt}}
\vspace{0.9\baselineskip}
\noindent \begin{minipage}[t]{0.25\textwidth}
\vspace{0pt}
\centering
\includegraphics[page=207,width=\linewidth]{knotoids.pdf}
\end{minipage}
\hfill
\begin{minipage}[t]{0.73\textwidth}
\vspace{0pt}
\raggedright
\textbf{Name:} {\large{$\mathbf{K7_{88}}$}} (achiral, non-rotatable$^{*}$, Possible duplicate [7\_90]) \\ \textbf{PD:} {\scriptsize\texttt{[0],[0,1,2,3],[1,4,5,2],[3,6,7,8],[4,8,9,10],[10,11,12,5],[6,12,11,13],[13,9,14,7],[14]}} \\ \textbf{EM:} {\scriptsize\texttt{(B0, A0C0C3D0, B1E0F3B2, B3G0H3E1, C1D3H1F0, E3G2G1C2, D1F2F1H0, G3E2I0D2, H2)}} \\ \textbf{Kauffman bracket:} {\scriptsize $-A^{26} + 3A^{22} + 2A^{20} - 4A^{18} - 3A^{16} + 4A^{14} + 5A^{12} - 3A^{10} - 5A^{8} + A^{6} + 3A^{4} - 1$} \\ \textbf{Arrow:} {\scriptsize $-A^{2}L_1 + 3L_1/A^{2} + 2/A^{4} - 4L_1/A^{6} - 3/A^{8} + 4L_1/A^{10} + 5/A^{12} - 3L_1/A^{14} - 5/A^{16} + L_1/A^{18} + 3/A^{20} - 1/A^{24}$} \\ \textbf{Mock:} {\scriptsize $-w^{5} - w^{4} + 3w^{3} + 4w^{2} - 4w - 5 + 4/w + 5/w^{2} - 3/w^{3} - 3/w^{4} + w^{-5} + w^{-6}$} \\ \textbf{Affine:} {\scriptsize $0$} \\ \textbf{Yamada:} {\scriptsize $2A^{26} - A^{25} - 5A^{24} + 6A^{23} + 3A^{22} - 13A^{21} + 3A^{20} + 10A^{19} - 14A^{18} - 2A^{17} + 11A^{16} - 7A^{15} - 5A^{14} + 6A^{13} + 2A^{12} - 4A^{11} - 4A^{10} + 12A^{9} - 3A^{8} - 13A^{7} + 15A^{6} - 11A^{4} + 7A^{3} + A^{2} - 3A + 1$}
\end{minipage}

\noindent{\color{gray!40}\rule{\textwidth}{0.4pt}}
\vspace{0.9\baselineskip}
\noindent \begin{minipage}[t]{0.25\textwidth}
\vspace{0pt}
\centering
\includegraphics[page=208,width=\linewidth]{knotoids.pdf}
\end{minipage}
\hfill
\begin{minipage}[t]{0.73\textwidth}
\vspace{0pt}
\raggedright
\textbf{Name:} {\large{$\mathbf{K7_{89}}$}} (chiral, non-rotatable$^{*}$) \\ \textbf{PD:} {\scriptsize\texttt{[0],[0,1,2,3],[1,4,5,6],[7,8,9,2],[3,9,10,4],[5],[6,11,12,13],[13,14,8,7],[14,12,11,10]}} \\ \textbf{EM:} {\scriptsize\texttt{(B0, A0C0D3E0, B1E3F0G0, H3H2E1B2, B3D2I3C1, C2, C3I2I1H0, G3I0D1D0, H1G2G1E2)}} \\ \textbf{Kauffman bracket:} {\scriptsize $A^{28} + A^{26} - 2A^{24} - 2A^{22} + 3A^{20} + 4A^{18} - 3A^{16} - 4A^{14} + 2A^{12} + 4A^{10} - A^{8} - 3A^{6} + A^{2}$} \\ \textbf{Arrow:} {\scriptsize $L_1/A^{2} + A^{-4} - 2L_1/A^{6} - 2/A^{8} + 3L_1/A^{10} + 4/A^{12} - 3L_1/A^{14} - 4/A^{16} + 2L_1/A^{18} + 4/A^{20} - L_1/A^{22} - 3/A^{24} + A^{-28}$} \\ \textbf{Mock:} {\scriptsize $2w^{3} + 4w^{2} - 4w - 7 + 4/w + 6/w^{2} - 2/w^{3} - 2/w^{4}$} \\ \textbf{Affine:} {\scriptsize $2t - 4 + 2/t$} \\ \textbf{Yamada:} {\scriptsize $A^{27} - A^{26} + 4A^{24} - 3A^{23} - A^{22} + 9A^{21} - 6A^{20} - 3A^{19} + 11A^{18} - 4A^{17} - 3A^{16} + 6A^{15} - A^{14} - A^{13} - 4A^{12} + 4A^{11} + A^{10} - 9A^{9} + 7A^{8} + 3A^{7} - 9A^{6} + 5A^{5} + 4A^{4} - 6A^{3} + A^{2} + 2A - 1$}
\end{minipage}

\noindent{\color{gray!40}\rule{\textwidth}{0.4pt}}
\vspace{0.9\baselineskip}
\noindent \begin{minipage}[t]{0.25\textwidth}
\vspace{0pt}
\centering
\includegraphics[page=209,width=\linewidth]{knotoids.pdf}
\end{minipage}
\hfill
\begin{minipage}[t]{0.73\textwidth}
\vspace{0pt}
\raggedright
\textbf{Name:} {\large{$\mathbf{K7_{90}}$}} (achiral, non-rotatable$^{*}$, Possible duplicate [7\_88]) \\ \textbf{PD:} {\scriptsize\texttt{[0],[0,1,2,3],[1,4,5,6],[6,7,8,2],[3,9,10,4],[11,12,13,5],[7,13,14,8],[9,14,12,10],[11]}} \\ \textbf{EM:} {\scriptsize\texttt{(B0, A0C0D3E0, B1E3F3D0, C3G0G3B2, B3H0H3C1, I0H2G1C2, D1F2H1D2, E1G2F1E2, F0)}} \\ \textbf{Kauffman bracket:} {\scriptsize $-A^{26} + 3A^{22} + 2A^{20} - 4A^{18} - 3A^{16} + 4A^{14} + 5A^{12} - 3A^{10} - 5A^{8} + A^{6} + 3A^{4} - 1$} \\ \textbf{Arrow:} {\scriptsize $-A^{2}L_1 + 3L_1/A^{2} + 2/A^{4} - 4L_1/A^{6} - 3/A^{8} + 4L_1/A^{10} + 5/A^{12} - 3L_1/A^{14} - 5/A^{16} + L_1/A^{18} + 3/A^{20} - 1/A^{24}$} \\ \textbf{Mock:} {\scriptsize $-w^{5} - w^{4} + 3w^{3} + 4w^{2} - 4w - 5 + 4/w + 5/w^{2} - 3/w^{3} - 3/w^{4} + w^{-5} + w^{-6}$} \\ \textbf{Affine:} {\scriptsize $0$} \\ \textbf{Yamada:} {\scriptsize $2A^{26} - A^{25} - 5A^{24} + 6A^{23} + 3A^{22} - 13A^{21} + 3A^{20} + 10A^{19} - 14A^{18} - 2A^{17} + 11A^{16} - 7A^{15} - 5A^{14} + 6A^{13} + 2A^{12} - 4A^{11} - 4A^{10} + 12A^{9} - 3A^{8} - 13A^{7} + 15A^{6} - 11A^{4} + 7A^{3} + A^{2} - 3A + 1$}
\end{minipage}

\noindent{\color{gray!40}\rule{\textwidth}{0.4pt}}
\vspace{0.9\baselineskip}
\noindent \begin{minipage}[t]{0.25\textwidth}
\vspace{0pt}
\centering
\includegraphics[page=210,width=\linewidth]{knotoids.pdf}
\end{minipage}
\hfill
\begin{minipage}[t]{0.73\textwidth}
\vspace{0pt}
\raggedright
\textbf{Name:} {\large{$\mathbf{K7_{91}}$}} (achiral, non-rotatable$^{*}$, Possible duplicate [7\_92]) \\ \textbf{PD:} {\scriptsize\texttt{[0],[0,1,2,3],[1,4,5,2],[3,6,7,8],[4,8,9,10],[11,12,6,5],[13,9,14,7],[10,13,12,11],[14]}} \\ \textbf{EM:} {\scriptsize\texttt{(B0, A0C0C3D0, B1E0F3B2, B3F2G3E1, C1D3G1H0, H3H2D1C2, H1E2I0D2, E3G0F1F0, G2)}} \\ \textbf{Kauffman bracket:} {\scriptsize $A^{28} - 3A^{24} - 2A^{22} + 5A^{20} + 4A^{18} - 4A^{16} - 5A^{14} + 3A^{12} + 5A^{10} - A^{8} - 3A^{6} + A^{2}$} \\ \textbf{Arrow:} {\scriptsize $A^{16} - 3A^{12} - 2A^{10}L_1 + 5A^{8} + 4A^{6}L_1 - 4A^{4} - 5A^{2}L_1 + 3 + 5L_1/A^{2} - 1/A^{4} - 3L_1/A^{6} + L_1/A^{10}$} \\ \textbf{Mock:} {\scriptsize $-w^{6} - w^{5} + 3w^{4} + 3w^{3} - 5w^{2} - 6w + 5 + 6/w - 2/w^{2} - 3/w^{3} + w^{-4} + w^{-5}$} \\ \textbf{Affine:} {\scriptsize $-2t + 4 - 2/t$} \\ \textbf{Yamada:} {\scriptsize $-2A^{26} + 2A^{25} + 4A^{24} - 9A^{23} + 2A^{22} + 13A^{21} - 14A^{20} - 4A^{19} + 17A^{18} - 9A^{17} - 6A^{16} + 10A^{15} - 3A^{13} - 3A^{12} + 11A^{11} - A^{10} - 13A^{9} + 15A^{8} + 3A^{7} - 16A^{6} + 9A^{5} + 6A^{4} - 9A^{3} + A^{2} + 3A - 1$}
\end{minipage}

\noindent{\color{gray!40}\rule{\textwidth}{0.4pt}}
\vspace{0.9\baselineskip}
\noindent \begin{minipage}[t]{0.25\textwidth}
\vspace{0pt}
\centering
\includegraphics[page=211,width=\linewidth]{knotoids.pdf}
\end{minipage}
\hfill
\begin{minipage}[t]{0.73\textwidth}
\vspace{0pt}
\raggedright
\textbf{Name:} {\large{$\mathbf{K7_{92}}$}} (achiral, non-rotatable$^{*}$, Possible duplicate [7\_91]) \\ \textbf{PD:} {\scriptsize\texttt{[0],[0,1,2,3],[1,4,5,6],[7,8,3,2],[8,9,10,4],[11,10,12,5],[6,13,14,7],[9,14,13,12],[11]}} \\ \textbf{EM:} {\scriptsize\texttt{(B0, A0C0D3D2, B1E3F3G0, G3E0B3B2, D1H0F1C1, I0E2H3C2, C3H2H1D0, E1G2G1F2, F0)}} \\ \textbf{Kauffman bracket:} {\scriptsize $A^{28} - 3A^{24} - A^{22} + 5A^{20} + 3A^{18} - 5A^{16} - 4A^{14} + 4A^{12} + 5A^{10} - 2A^{8} - 3A^{6} + A^{2}$} \\ \textbf{Arrow:} {\scriptsize $A^{10}L_1 - 3A^{6}L_1 - A^{4} + 5A^{2}L_1 + 3 - 5L_1/A^{2} - 4/A^{4} + 4L_1/A^{6} + 5/A^{8} - 2L_1/A^{10} - 3/A^{12} + A^{-16}$} \\ \textbf{Mock:} {\scriptsize $w^{5} + w^{4} - 3w^{3} - 2w^{2} + 6w + 5 - 6/w - 5/w^{2} + 3/w^{3} + 3/w^{4} - 1/w^{5} - 1/w^{6}$} \\ \textbf{Affine:} {\scriptsize $2t - 4 + 2/t$} \\ \textbf{Yamada:} {\scriptsize $-A^{26} + 3A^{25} + A^{24} - 9A^{23} + 6A^{22} + 9A^{21} - 16A^{20} + 3A^{19} + 15A^{18} - 13A^{17} - A^{16} + 11A^{15} - 3A^{14} - 3A^{13} + 10A^{11} - 6A^{10} - 9A^{9} + 17A^{8} - 4A^{7} - 14A^{6} + 13A^{5} + 2A^{4} - 9A^{3} + 4A^{2} + 2A - 2$}
\end{minipage}

\noindent{\color{gray!40}\rule{\textwidth}{0.4pt}}
\vspace{0.9\baselineskip}
\noindent \begin{minipage}[t]{0.25\textwidth}
\vspace{0pt}
\centering
\includegraphics[page=212,width=\linewidth]{knotoids.pdf}
\end{minipage}
\hfill
\begin{minipage}[t]{0.73\textwidth}
\vspace{0pt}
\raggedright
\textbf{Name:} {\large{$\mathbf{K7_{93}}$}} (chiral, non-rotatable$^{*}$) \\ \textbf{PD:} {\scriptsize\texttt{[0],[0,1,2,3],[1,4,5,6],[6,7,8,2],[3,9,10,4],[11,12,7,5],[12,13,9,8],[13,11,14,10],[14]}} \\ \textbf{EM:} {\scriptsize\texttt{(B0, A0C0D3E0, B1E3F3D0, C3F2G3B2, B3G2H3C1, H1G0D1C2, F1H0E1D2, G1F0I0E2, H2)}} \\ \textbf{Kauffman bracket:} {\scriptsize $A^{28} - 4A^{24} - A^{22} + 6A^{20} + 4A^{18} - 6A^{16} - 5A^{14} + 5A^{12} + 6A^{10} - 2A^{8} - 4A^{6} + A^{2}$} \\ \textbf{Arrow:} {\scriptsize $A^{10}L_1 - 4A^{6}L_1 - A^{4} + 6A^{2}L_1 + 4 - 6L_1/A^{2} - 5/A^{4} + 5L_1/A^{6} + 6/A^{8} - 2L_1/A^{10} - 4/A^{12} + A^{-16}$} \\ \textbf{Mock:} {\scriptsize $w^{5} + w^{4} - 4w^{3} - 3w^{2} + 7w + 6 - 7/w - 6/w^{2} + 4/w^{3} + 4/w^{4} - 1/w^{5} - 1/w^{6}$} \\ \textbf{Affine:} {\scriptsize $0$} \\ \textbf{Yamada:} {\scriptsize $-A^{26} + 4A^{25} + A^{24} - 14A^{23} + 9A^{22} + 14A^{21} - 25A^{20} + 4A^{19} + 23A^{18} - 19A^{17} - 3A^{16} + 17A^{15} - 3A^{14} - 6A^{13} + A^{12} + 15A^{11} - 9A^{10} - 14A^{9} + 25A^{8} - 5A^{7} - 23A^{6} + 19A^{5} + 4A^{4} - 15A^{3} + 5A^{2} + 4A - 2$}
\end{minipage}

\noindent{\color{gray!40}\rule{\textwidth}{0.4pt}}
\vspace{0.9\baselineskip}
\noindent \begin{minipage}[t]{0.25\textwidth}
\vspace{0pt}
\centering
\includegraphics[page=213,width=\linewidth]{knotoids.pdf}
\end{minipage}
\hfill
\begin{minipage}[t]{0.73\textwidth}
\vspace{0pt}
\raggedright
\textbf{Name:} {\large{$\mathbf{K7_{94}}$}} (chiral, non-rotatable$^{*}$) \\ \textbf{PD:} {\scriptsize\texttt{[0],[0,1,2,3],[1,4,5,2],[3,5,6,7],[4,8,9,10],[10,11,12,6],[12,13,14,7],[8,14,13,9],[11]}} \\ \textbf{EM:} {\scriptsize\texttt{(B0, A0C0C3D0, B1E0D1B2, B3C2F3G3, C1H0H3F0, E3I0G0D2, F2H2H1D3, E1G2G1E2, F1)}} \\ \textbf{Kauffman bracket:} {\scriptsize $A^{22} + A^{20} - 2A^{18} - 3A^{16} + 2A^{14} + 4A^{12} - 3A^{8} + 2A^{4} - 1$} \\ \textbf{Arrow:} {\scriptsize $A^{4} + A^{2}L_1 - L_2 - 1 - 3L_1/A^{2} + L_2/A^{4} + A^{-4} + 4L_1/A^{6} - 3L_1/A^{10} + 2L_1/A^{14} - L_1/A^{18}$} \\ \textbf{Mock:} {\scriptsize $-w^{4} - 2w^{3} + 2w^{2} + 5w - 4/w + w^{-2} + 2/w^{3} - 1/w^{4} - 1/w^{5}$} \\ \textbf{Affine:} {\scriptsize $-t^{2} + t + 1/t - 1/t^{2}$} \\ \textbf{Yamada:} {\scriptsize $-A^{26} - A^{25} + A^{24} + A^{23} - A^{22} - A^{21} + 4A^{20} + 2A^{19} - 5A^{18} + 4A^{17} + 4A^{16} - 5A^{15} + 2A^{14} + 3A^{13} - 3A^{12} + A^{10} + 3A^{9} - 3A^{8} - A^{7} + 6A^{6} - 4A^{5} - 3A^{4} + 5A^{3} - A^{2} - 2A + 1$}
\end{minipage}

\noindent{\color{gray!40}\rule{\textwidth}{0.4pt}}
\vspace{0.9\baselineskip}
\noindent \begin{minipage}[t]{0.25\textwidth}
\vspace{0pt}
\centering
\includegraphics[page=214,width=\linewidth]{knotoids.pdf}
\end{minipage}
\hfill
\begin{minipage}[t]{0.73\textwidth}
\vspace{0pt}
\raggedright
\textbf{Name:} {\large{$\mathbf{K7_{95}}$}} (chiral, non-rotatable$^{*}$) \\ \textbf{PD:} {\scriptsize\texttt{[0],[0,1,2,3],[1,4,5,2],[3,5,6,7],[8,9,10,4],[10,11,12,6],[7,12,13,14],[14,13,9,8],[11]}} \\ \textbf{EM:} {\scriptsize\texttt{(B0, A0C0C3D0, B1E3D1B2, B3C2F3G0, H3H2F0C1, E2I0G1D2, D3F2H1H0, G3G2E1E0, F1)}} \\ \textbf{Kauffman bracket:} {\scriptsize $-A^{24} + 2A^{20} - 3A^{16} + 4A^{12} + A^{10} - 4A^{8} - 2A^{6} + 2A^{4} + 2A^{2}$} \\ \textbf{Arrow:} {\scriptsize $-A^{30}L_1 + 2A^{26}L_1 - 3A^{22}L_1 + 4A^{18}L_1 + A^{16} - 4A^{14}L_1 - A^{12}L_2 - A^{12} + 2A^{10}L_1 + A^{8}L_2 + A^{8}$} \\ \textbf{Mock:} {\scriptsize $w^{5} + w^{4} - 3w^{3} - w^{2} + 4w - 5/w - 1/w^{2} + 3/w^{3} + 2/w^{4}$} \\ \textbf{Affine:} {\scriptsize $-t^{2} - 2t + 6 - 2/t - 1/t^{2}$} \\ \textbf{Yamada:} {\scriptsize $-A^{27} + 2A^{26} - 4A^{24} + 5A^{23} + 2A^{22} - 7A^{21} + 4A^{20} + 3A^{19} - 4A^{18} + A^{16} + 2A^{15} - 4A^{14} + 5A^{12} - 6A^{11} - 2A^{10} + 6A^{9} - 4A^{8} - 3A^{7} + 3A^{6} - 3A^{4} - A^{3} + A^{2} - 1$}
\end{minipage}

\noindent{\color{gray!40}\rule{\textwidth}{0.4pt}}
\vspace{0.9\baselineskip}
\noindent \begin{minipage}[t]{0.25\textwidth}
\vspace{0pt}
\centering
\includegraphics[page=215,width=\linewidth]{knotoids.pdf}
\end{minipage}
\hfill
\begin{minipage}[t]{0.73\textwidth}
\vspace{0pt}
\raggedright
\textbf{Name:} {\large{$\mathbf{K7_{96}}$}} (chiral, non-rotatable$^{*}$) \\ \textbf{PD:} {\scriptsize\texttt{[0],[0,1,2,3],[1,4,5,2],[3,6,7,4],[5,7,8,9],[10,11,12,6],[13,14,9,8],[14,13,11,10],[12]}} \\ \textbf{EM:} {\scriptsize\texttt{(B0, A0C0C3D0, B1D3E0B2, B3F3E1C1, C2D2G3G2, H3H2I0D1, H1H0E3E2, G1G0F1F0, F2)}} \\ \textbf{Kauffman bracket:} {\scriptsize $A^{20} + A^{18} - A^{16} - A^{14} + A^{12} + 2A^{10} - A^{8} - 2A^{6} - A^{4} + A^{2} + 1$} \\ \textbf{Arrow:} {\scriptsize $A^{2}L_1 + 1 - L_1/A^{2} - 1/A^{4} + L_1/A^{6} + 2/A^{8} - L_1/A^{10} - 2/A^{12} - L_1/A^{14} + A^{-16} + L_1/A^{18}$} \\ \textbf{Mock:} {\scriptsize $-w^{2} + w + 4 - 1/w - 2/w^{2}$} \\ \textbf{Affine:} {\scriptsize $t - 2 + 1/t$} \\ \textbf{Yamada:} {\scriptsize $A^{25} - A^{24} + A^{23} + 2A^{22} - 2A^{21} + 2A^{20} + 3A^{19} - A^{18} + A^{17} + A^{15} - 2A^{14} - A^{13} + 2A^{12} - 2A^{11} - A^{10} + 2A^{9} - A^{7} + 2A^{6} + A^{5} - A^{4} + A^{2} - 1$}
\end{minipage}

\noindent{\color{gray!40}\rule{\textwidth}{0.4pt}}
\vspace{0.9\baselineskip}
\noindent \begin{minipage}[t]{0.25\textwidth}
\vspace{0pt}
\centering
\includegraphics[page=216,width=\linewidth]{knotoids.pdf}
\end{minipage}
\hfill
\begin{minipage}[t]{0.73\textwidth}
\vspace{0pt}
\raggedright
\textbf{Name:} {\large{$\mathbf{K7_{97}}$}} (chiral, non-rotatable$^{*}$) \\ \textbf{PD:} {\scriptsize\texttt{[0],[0,1,2,3],[1,4,5,2],[3,6,7,8],[9,10,11,4],[5,11,12,13],[13,12,14,6],[10,9,8,7],[14]}} \\ \textbf{EM:} {\scriptsize\texttt{(B0, A0C0C3D0, B1E3F0B2, B3G3H3H2, H1H0F1C1, C2E2G1G0, F3F2I0D1, E1E0D3D2, G2)}} \\ \textbf{Kauffman bracket:} {\scriptsize $A^{18} + A^{12} + A^{10} - A^{8} - 2A^{6} - A^{4} + A^{2} + 1$} \\ \textbf{Arrow:} {\scriptsize $1 + L_1/A^{6} + A^{-8} - L_1/A^{10} - 2/A^{12} - L_1/A^{14} + A^{-16} + L_1/A^{18}$} \\ \textbf{Mock:} {\scriptsize $-w^{2} + 3 - 1/w^{2}$} \\ \textbf{Affine:} {\scriptsize $0$} \\ \textbf{Yamada:} {\scriptsize $A^{23} + 2A^{20} + 2A^{19} + A^{17} + A^{15} - A^{14} + A^{12} - 2A^{11} - A^{8} + A^{6} + A^{5} + A^{2} - 1$}
\end{minipage}

\noindent{\color{gray!40}\rule{\textwidth}{0.4pt}}
\vspace{0.9\baselineskip}
\noindent \begin{minipage}[t]{0.25\textwidth}
\vspace{0pt}
\centering
\includegraphics[page=217,width=\linewidth]{knotoids.pdf}
\end{minipage}
\hfill
\begin{minipage}[t]{0.73\textwidth}
\vspace{0pt}
\raggedright
\textbf{Name:} {\large{$\mathbf{K7_{98}}$}} (chiral, rotatable) \\ \textbf{PD:} {\scriptsize\texttt{[0],[0,1,2,3],[1,4,5,2],[3,6,7,8],[9,10,11,4],[11,12,13,5],[6,13,12,14],[10,9,8,7],[14]}} \\ \textbf{EM:} {\scriptsize\texttt{(B0, A0C0C3D0, B1E3F3B2, B3G0H3H2, H1H0F0C1, E2G2G1C2, D1F2F1I0, E1E0D3D2, G3)}} \\ \textbf{Kauffman bracket:} {\scriptsize $-A^{19} + A^{15} - 2A^{11} + 3A^{7} + 2A^{5} - 2A^{3} - 2A$} \\ \textbf{Arrow:} {\scriptsize $A^{16} - A^{12} + 2A^{8} - 3A^{4} - 2A^{2}L_1 + 2 + 2L_1/A^{2}$} \\ \textbf{Mock:} {\scriptsize $-3w^{2} - 2w + 5 + 2/w - 1/w^{2}$} \\ \textbf{Affine:} {\scriptsize $-2t + 4 - 2/t$} \\ \textbf{Yamada:} {\scriptsize $A^{25} - A^{24} - A^{23} + 2A^{22} + A^{21} - 3A^{20} + 3A^{18} - 2A^{17} + 2A^{15} + 3A^{11} - A^{10} - A^{9} + 4A^{8} - A^{7} - 2A^{6} + 2A^{5} - A^{3} - A^{2} + A + 1$}
\end{minipage}

\noindent{\color{gray!40}\rule{\textwidth}{0.4pt}}
\vspace{0.9\baselineskip}
\noindent \begin{minipage}[t]{0.25\textwidth}
\vspace{0pt}
\centering
\includegraphics[page=218,width=\linewidth]{knotoids.pdf}
\end{minipage}
\hfill
\begin{minipage}[t]{0.73\textwidth}
\vspace{0pt}
\raggedright
\textbf{Name:} {\large{$\mathbf{K7_{99}}$}} (chiral, non-rotatable$^{*}$) \\ \textbf{PD:} {\scriptsize\texttt{[0],[0,1,2,3],[1,4,5,2],[6,7,4,3],[5,7,8,9],[10,11,12,6],[13,14,9,8],[14,13,11,10],[12]}} \\ \textbf{EM:} {\scriptsize\texttt{(B0, A0C0C3D3, B1D2E0B2, F3E1C1B3, C2D1G3G2, H3H2I0D0, H1H0E3E2, G1G0F1F0, F2)}} \\ \textbf{Kauffman bracket:} {\scriptsize $A^{22} + A^{20} - A^{18} - 2A^{16} + A^{14} + 3A^{12} - A^{10} - 3A^{8} - A^{6} + 2A^{4} + A^{2}$} \\ \textbf{Arrow:} {\scriptsize $A^{10}L_1 + A^{8} - A^{6}L_1 - 2A^{4} + A^{2}L_1 + 3 - L_1/A^{2} - 3/A^{4} - L_1/A^{6} + 2/A^{8} + L_1/A^{10}$} \\ \textbf{Mock:} {\scriptsize $-2w^{2} + w + 6 - 1/w - 3/w^{2}$} \\ \textbf{Affine:} {\scriptsize $t - 2 + 1/t$} \\ \textbf{Yamada:} {\scriptsize $A^{25} - A^{24} - A^{23} + 2A^{22} + A^{21} - 3A^{20} + 2A^{19} + 4A^{18} - 4A^{17} + A^{16} + 3A^{15} - 2A^{14} + 2A^{11} - 3A^{10} - A^{9} + 5A^{8} - 2A^{7} - A^{6} + 4A^{5} - 2A^{3} + A^{2} + A - 1$}
\end{minipage}

\noindent{\color{gray!40}\rule{\textwidth}{0.4pt}}
\vspace{0.9\baselineskip}
\noindent \begin{minipage}[t]{0.25\textwidth}
\vspace{0pt}
\centering
\includegraphics[page=219,width=\linewidth]{knotoids.pdf}
\end{minipage}
\hfill
\begin{minipage}[t]{0.73\textwidth}
\vspace{0pt}
\raggedright
\textbf{Name:} {\large{$\mathbf{K7_{100}}$}} (chiral, non-rotatable$^{*}$) \\ \textbf{PD:} {\scriptsize\texttt{[0],[0,1,2,3],[1,4,5,2],[6,7,4,3],[7,8,9,5],[6,10,11,12],[8,13,14,9],[10,14,13,11],[12]}} \\ \textbf{EM:} {\scriptsize\texttt{(B0, A0C0C3D3, B1D2E3B2, F0E0C1B3, D1G0G3C2, D0H0H3I0, E1H2H1E2, F1G2G1F2, F3)}} \\ \textbf{Kauffman bracket:} {\scriptsize $2A^{23} + A^{21} - 2A^{19} - 3A^{17} + A^{15} + 3A^{13} - A^{11} - 3A^{9} + 2A^{5} - A$} \\ \textbf{Arrow:} {\scriptsize $-2A^{26}L_1 - A^{24} + 2A^{22}L_1 + 3A^{20} - A^{18}L_1 - 3A^{16} + A^{14}L_1 + 3A^{12} - 2A^{8} + A^{4}$} \\ \textbf{Mock:} {\scriptsize $2w^{2} - 3w - 6 + 3/w + 5/w^{2}$} \\ \textbf{Affine:} {\scriptsize $-3t + 6 - 3/t$} \\ \textbf{Yamada:} {\scriptsize $-A^{27} - A^{26} + A^{25} + A^{24} - 2A^{23} + A^{22} + 4A^{21} - A^{20} - 2A^{19} + 4A^{18} - 4A^{16} + 2A^{15} + 2A^{14} - 4A^{13} + 3A^{11} - 2A^{10} - A^{9} + 3A^{7} - 5A^{6} - 2A^{5} + 4A^{4} - 4A^{3} - A^{2} + A - 2$}
\end{minipage}

\noindent{\color{gray!40}\rule{\textwidth}{0.4pt}}
\vspace{0.9\baselineskip}
\noindent \begin{minipage}[t]{0.25\textwidth}
\vspace{0pt}
\centering
\includegraphics[page=220,width=\linewidth]{knotoids.pdf}
\end{minipage}
\hfill
\begin{minipage}[t]{0.73\textwidth}
\vspace{0pt}
\raggedright
\textbf{Name:} {\large{$\mathbf{K7_{101}}$}} (chiral, rotatable) \\ \textbf{PD:} {\scriptsize\texttt{[0],[0,1,2,3],[1,4,5,2],[6,7,8,3],[4,9,10,11],[11,12,13,5],[6,13,12,14],[7,10,9,8],[14]}} \\ \textbf{EM:} {\scriptsize\texttt{(B0, A0C0C3D3, B1E0F3B2, G0H0H3B3, C1H2H1F0, E3G2G1C2, D0F2F1I0, D1E2E1D2, G3)}} \\ \textbf{Kauffman bracket:} {\scriptsize $2A^{23} + 2A^{21} - 2A^{19} - 4A^{17} + 3A^{13} - 3A^{9} + 2A^{5} - A$} \\ \textbf{Arrow:} {\scriptsize $-2A^{26}L_1 - 2A^{24} + 2A^{22}L_1 + 4A^{20} - 3A^{16} + 3A^{12} - 2A^{8} + A^{4}$} \\ \textbf{Mock:} {\scriptsize $3w^{2} - 2w - 7 + 2/w + 5/w^{2}$} \\ \textbf{Affine:} {\scriptsize $-2t + 4 - 2/t$} \\ \textbf{Yamada:} {\scriptsize $-A^{27} - A^{26} + A^{25} + A^{24} - 2A^{23} + 4A^{21} + A^{20} - 4A^{19} + 3A^{18} + 3A^{17} - 4A^{16} + A^{15} + 4A^{14} - 3A^{13} - 2A^{12} + 3A^{11} - 2A^{10} - 2A^{9} - 2A^{8} + 4A^{7} - 3A^{6} - 5A^{5} + 5A^{4} - 2A^{3} - 3A^{2} + A - 1$}
\end{minipage}

\noindent{\color{gray!40}\rule{\textwidth}{0.4pt}}
\vspace{0.9\baselineskip}
\noindent \begin{minipage}[t]{0.25\textwidth}
\vspace{0pt}
\centering
\includegraphics[page=221,width=\linewidth]{knotoids.pdf}
\end{minipage}
\hfill
\begin{minipage}[t]{0.73\textwidth}
\vspace{0pt}
\raggedright
\textbf{Name:} {\large{$\mathbf{K7_{102}}$}} (chiral, non-rotatable$^{*}$) \\ \textbf{PD:} {\scriptsize\texttt{[0],[0,1,2,3],[1,4,5,2],[3,6,7,4],[5,8,9,10],[6,11,8,7],[12,13,10,9],[13,12,14,11],[14]}} \\ \textbf{EM:} {\scriptsize\texttt{(B0, A0C0C3D0, B1D3E0B2, B3F0F3C1, C2F2G3G2, D1H3E1D2, H1H0E3E2, G1G0I0F1, H2)}} \\ \textbf{Kauffman bracket:} {\scriptsize $A^{21} + A^{19} - 2A^{17} - 3A^{15} + 2A^{13} + 4A^{11} - 4A^{7} - 2A^{5} + A^{3} + A$} \\ \textbf{Arrow:} {\scriptsize $-A^{18}L_1 - A^{16} + 2A^{14}L_1 + 3A^{12} - 2A^{10}L_1 - 4A^{8} + 4A^{4} + 2A^{2}L_1 - 1 - L_1/A^{2}$} \\ \textbf{Mock:} {\scriptsize $-w^{5} - w^{4} + 2w^{3} + 3w^{2} - 2w - 4 + 4/w^{2} + 2/w^{3} - 1/w^{4} - 1/w^{5}$} \\ \textbf{Affine:} {\scriptsize $-t + 2 - 1/t$} \\ \textbf{Yamada:} {\scriptsize $-2A^{25} + 2A^{24} + 3A^{23} - 6A^{22} + 2A^{21} + 7A^{20} - 7A^{19} + 6A^{17} - 3A^{16} - 2A^{15} + 3A^{13} - 4A^{12} - 4A^{11} + 6A^{10} - 4A^{9} - 7A^{8} + 6A^{7} - 5A^{5} + 2A^{4} + 3A^{3} - 2A^{2} - A + 1$}
\end{minipage}

\noindent{\color{gray!40}\rule{\textwidth}{0.4pt}}
\vspace{0.9\baselineskip}
\noindent \begin{minipage}[t]{0.25\textwidth}
\vspace{0pt}
\centering
\includegraphics[page=222,width=\linewidth]{knotoids.pdf}
\end{minipage}
\hfill
\begin{minipage}[t]{0.73\textwidth}
\vspace{0pt}
\raggedright
\textbf{Name:} {\large{$\mathbf{K7_{103}}$}} (chiral, non-rotatable$^{*}$) \\ \textbf{PD:} {\scriptsize\texttt{[0],[0,1,2,3],[1,4,5,2],[3,6,7,4],[8,9,10,5],[6,11,8,7],[9,12,13,10],[11,13,12,14],[14]}} \\ \textbf{EM:} {\scriptsize\texttt{(B0, A0C0C3D0, B1D3E3B2, B3F0F3C1, F2G0G3C2, D1H0E0D2, E1H2H1E2, F1G2G1I0, H3)}} \\ \textbf{Kauffman bracket:} {\scriptsize $A^{24} - A^{20} - 2A^{18} + A^{16} + 4A^{14} + 2A^{12} - 3A^{10} - 3A^{8} + A^{6} + 2A^{4} - 1$} \\ \textbf{Arrow:} {\scriptsize $1 - 1/A^{4} - 2L_1/A^{6} + A^{-8} + 4L_1/A^{10} + 2/A^{12} - 3L_1/A^{14} - 3/A^{16} + L_1/A^{18} + 2/A^{20} - 1/A^{24}$} \\ \textbf{Mock:} {\scriptsize $-w^{6} - 2w^{5} + w^{4} + 4w^{3} + 2w^{2} - 3w - 3 + 1/w + 2/w^{2}$} \\ \textbf{Affine:} {\scriptsize $-t + 2 - 1/t$} \\ \textbf{Yamada:} {\scriptsize $A^{24} + A^{23} - 3A^{22} - A^{21} + 4A^{20} - 4A^{19} - 7A^{18} + 6A^{17} - A^{16} - 8A^{15} + 4A^{14} + 2A^{13} - 4A^{12} + A^{11} + 2A^{10} + 2A^{9} - 5A^{8} + 3A^{7} + 7A^{6} - 8A^{5} + A^{4} + 6A^{3} - 5A^{2} - 2A + 2$}
\end{minipage}

\noindent{\color{gray!40}\rule{\textwidth}{0.4pt}}
\vspace{0.9\baselineskip}
\noindent \begin{minipage}[t]{0.25\textwidth}
\vspace{0pt}
\centering
\includegraphics[page=223,width=\linewidth]{knotoids.pdf}
\end{minipage}
\hfill
\begin{minipage}[t]{0.73\textwidth}
\vspace{0pt}
\raggedright
\textbf{Name:} {\large{$\mathbf{K7_{104}}$}} (chiral, non-rotatable$^{*}$) \\ \textbf{PD:} {\scriptsize\texttt{[0],[0,1,2,3],[1,4,5,2],[3,6,7,8],[4,8,9,5],[6,10,11,7],[9,11,12,13],[13,12,14,10],[14]}} \\ \textbf{EM:} {\scriptsize\texttt{(B0, A0C0C3D0, B1E0E3B2, B3F0F3E1, C1D3G0C2, D1H3G1D2, E2F2H1H0, G3G2I0F1, H2)}} \\ \textbf{Kauffman bracket:} {\scriptsize $A^{22} + 2A^{20} - 2A^{18} - 4A^{16} + 5A^{12} + 2A^{10} - 3A^{8} - 2A^{6} + A^{4} + A^{2}$} \\ \textbf{Arrow:} {\scriptsize $A^{10}L_1 + 2A^{8} - 2A^{6}L_1 - 4A^{4} + 5 + 2L_1/A^{2} - 3/A^{4} - 2L_1/A^{6} + A^{-8} + L_1/A^{10}$} \\ \textbf{Mock:} {\scriptsize $w^{5} + w^{4} - 2w^{3} - 4w^{2} + 6 + 2/w - 3/w^{2} - 2/w^{3} + w^{-4} + w^{-5}$} \\ \textbf{Affine:} {\scriptsize $-t + 2 - 1/t$} \\ \textbf{Yamada:} {\scriptsize $A^{25} - 2A^{24} - 2A^{23} + 5A^{22} - 7A^{20} + 6A^{19} + 4A^{18} - 10A^{17} + 3A^{16} + 4A^{15} - 6A^{14} - A^{13} + A^{12} + 3A^{11} - 6A^{10} + 8A^{8} - 8A^{7} - 3A^{6} + 8A^{5} - 3A^{4} - 4A^{3} + 4A^{2} + A - 2$}
\end{minipage}

\noindent{\color{gray!40}\rule{\textwidth}{0.4pt}}
\vspace{0.9\baselineskip}
\noindent \begin{minipage}[t]{0.25\textwidth}
\vspace{0pt}
\centering
\includegraphics[page=224,width=\linewidth]{knotoids.pdf}
\end{minipage}
\hfill
\begin{minipage}[t]{0.73\textwidth}
\vspace{0pt}
\raggedright
\textbf{Name:} {\large{$\mathbf{K7_{105}}$}} (chiral, non-rotatable$^{*}$) \\ \textbf{PD:} {\scriptsize\texttt{[0],[0,1,2,3],[1,4,5,2],[6,7,4,3],[8,9,10,5],[11,8,7,6],[9,12,13,10],[11,13,12,14],[14]}} \\ \textbf{EM:} {\scriptsize\texttt{(B0, A0C0C3D3, B1D2E3B2, F3F2C1B3, F1G0G3C2, H0E0D1D0, E1H2H1E2, F0G2G1I0, H3)}} \\ \textbf{Kauffman bracket:} {\scriptsize $A^{24} + 2A^{22} + A^{20} - 4A^{18} - 3A^{16} + 3A^{14} + 5A^{12} - A^{10} - 4A^{8} + 2A^{4} - 1$} \\ \textbf{Arrow:} {\scriptsize $A^{-12} + 2L_1/A^{14} + A^{-16} - 4L_1/A^{18} - 3/A^{20} + 3L_1/A^{22} + 5/A^{24} - L_1/A^{26} - 4/A^{28} + 2/A^{32} - 1/A^{36}$} \\ \textbf{Mock:} {\scriptsize $w^{6} + 2w^{5} + w^{4} - 4w^{3} - 4w^{2} + 3w + 5 - 1/w - 4/w^{2} + 2/w^{4}$} \\ \textbf{Affine:} {\scriptsize $t - 2 + 1/t$} \\ \textbf{Yamada:} {\scriptsize $A^{27} - 2A^{25} + 2A^{24} + 5A^{23} - 4A^{22} - A^{21} + 10A^{20} - A^{19} - 8A^{18} + 9A^{17} + 4A^{16} - 9A^{15} + 4A^{14} + 5A^{13} - 6A^{12} - 3A^{11} + 2A^{10} + 3A^{9} - 10A^{8} + A^{7} + 10A^{6} - 10A^{5} - A^{4} + 8A^{3} - 3A^{2} - 2A + 2$}
\end{minipage}

\noindent{\color{gray!40}\rule{\textwidth}{0.4pt}}
\vspace{0.9\baselineskip}
\noindent \begin{minipage}[t]{0.25\textwidth}
\vspace{0pt}
\centering
\includegraphics[page=225,width=\linewidth]{knotoids.pdf}
\end{minipage}
\hfill
\begin{minipage}[t]{0.73\textwidth}
\vspace{0pt}
\raggedright
\textbf{Name:} {\large{$\mathbf{K7_{106}}$}} (chiral, non-rotatable$^{*}$) \\ \textbf{PD:} {\scriptsize\texttt{[0],[0,1,2,3],[1,4,5,2],[3,5,6,7],[4,8,9,10],[11,12,13,6],[13,9,8,7],[10,14,12,11],[14]}} \\ \textbf{EM:} {\scriptsize\texttt{(B0, A0C0C3D0, B1E0D1B2, B3C2F3G3, C1G2G1H0, H3H2G0D2, F2E2E1D3, E3I0F1F0, H1)}} \\ \textbf{Kauffman bracket:} {\scriptsize $2A^{20} + A^{18} - 3A^{16} - 3A^{14} + 2A^{12} + 3A^{10} - 2A^{6} + A^{2}$} \\ \textbf{Arrow:} {\scriptsize $A^{8}L_2 + A^{8} + A^{6}L_1 - 2A^{4}L_2 - A^{4} - 3A^{2}L_1 + L_2 + 1 + 3L_1/A^{2} - 2L_1/A^{6} + L_1/A^{10}$} \\ \textbf{Mock:} {\scriptsize $w^{3} - w^{2} - 3w + 3 + 3/w - 2/w^{2} - 2/w^{3} + w^{-4} + w^{-5}$} \\ \textbf{Affine:} {\scriptsize $-t + 2 - 1/t$} \\ \textbf{Yamada:} {\scriptsize $-A^{25} - A^{24} + 2A^{23} + A^{22} - 2A^{21} + 5A^{19} - A^{18} - 5A^{17} + 4A^{16} - A^{15} - 5A^{14} + A^{13} - A^{11} - 3A^{10} + 3A^{9} + A^{8} - 5A^{7} + 3A^{6} + 3A^{5} - 4A^{4} - A^{3} + 3A^{2} - 2$}
\end{minipage}

\noindent{\color{gray!40}\rule{\textwidth}{0.4pt}}
\vspace{0.9\baselineskip}
\noindent \begin{minipage}[t]{0.25\textwidth}
\vspace{0pt}
\centering
\includegraphics[page=226,width=\linewidth]{knotoids.pdf}
\end{minipage}
\hfill
\begin{minipage}[t]{0.73\textwidth}
\vspace{0pt}
\raggedright
\textbf{Name:} {\large{$\mathbf{K7_{107}}$}} (chiral, non-rotatable$^{*}$) \\ \textbf{PD:} {\scriptsize\texttt{[0],[0,1,2,3],[1,4,5,2],[3,5,6,7],[8,9,10,4],[11,12,13,6],[7,13,9,8],[10,14,12,11],[14]}} \\ \textbf{EM:} {\scriptsize\texttt{(B0, A0C0C3D0, B1E3D1B2, B3C2F3G0, G3G2H0C1, H3H2G1D2, D3F2E1E0, E2I0F1F0, H1)}} \\ \textbf{Kauffman bracket:} {\scriptsize $-A^{24} + 2A^{20} - 4A^{16} + 5A^{12} + 2A^{10} - 4A^{8} - 3A^{6} + 2A^{4} + 2A^{2}$} \\ \textbf{Arrow:} {\scriptsize $-A^{30}L_1 + 2A^{26}L_1 - 4A^{22}L_1 + 5A^{18}L_1 + A^{16}L_2 + A^{16} - 4A^{14}L_1 - 2A^{12}L_2 - A^{12} + 2A^{10}L_1 + A^{8}L_2 + A^{8}$} \\ \textbf{Mock:} {\scriptsize $w^{5} + w^{4} - 4w^{3} - 2w^{2} + 5w + 1 - 5/w - 1/w^{2} + 3/w^{3} + 2/w^{4}$} \\ \textbf{Affine:} {\scriptsize $-3t + 6 - 3/t$} \\ \textbf{Yamada:} {\scriptsize $-A^{27} + 3A^{26} - 7A^{24} + 5A^{23} + 4A^{22} - 10A^{21} + 4A^{20} + 7A^{19} - 4A^{18} + 3A^{16} + 3A^{15} - 6A^{14} - A^{13} + 7A^{12} - 8A^{11} - 5A^{10} + 8A^{9} - 4A^{8} - 6A^{7} + 4A^{6} + 2A^{5} - 4A^{4} - A^{3} + 2A^{2} - 1$}
\end{minipage}

\noindent{\color{gray!40}\rule{\textwidth}{0.4pt}}
\vspace{0.9\baselineskip}
\noindent \begin{minipage}[t]{0.25\textwidth}
\vspace{0pt}
\centering
\includegraphics[page=227,width=\linewidth]{knotoids.pdf}
\end{minipage}
\hfill
\begin{minipage}[t]{0.73\textwidth}
\vspace{0pt}
\raggedright
\textbf{Name:} {\large{$\mathbf{K7_{108}}$}} (chiral, non-rotatable$^{*}$) \\ \textbf{PD:} {\scriptsize\texttt{[0],[0,1,2,3],[1,4,5,2],[3,6,7,8],[9,10,11,4],[5,12,13,6],[14,9,8,7],[10,14,12,11],[13]}} \\ \textbf{EM:} {\scriptsize\texttt{(B0, A0C0C3D0, B1E3F0B2, B3F3G3G2, G1H0H3C1, C2H2I0D1, H1E0D3D2, E1G0F1E2, F2)}} \\ \textbf{Kauffman bracket:} {\scriptsize $-A^{23} + A^{19} - 2A^{15} + A^{13} + 3A^{11} - 3A^{7} - 2A^{5} + A^{3} + A$} \\ \textbf{Arrow:} {\scriptsize $A^{8} - A^{4} + 2 - L_1/A^{2} - 3/A^{4} + 3/A^{8} + 2L_1/A^{10} - 1/A^{12} - L_1/A^{14}$} \\ \textbf{Mock:} {\scriptsize $-w^{3} - 3w^{2} + 6 + 2/w - 2/w^{2} - 1/w^{3}$} \\ \textbf{Affine:} {\scriptsize $-t + 2 - 1/t$} \\ \textbf{Yamada:} {\scriptsize $A^{23} - 2A^{22} - A^{21} + 3A^{20} - 2A^{19} - 2A^{18} + 2A^{17} - 2A^{16} - 2A^{15} - A^{14} + A^{13} - A^{12} - 2A^{11} + 3A^{10} - 2A^{8} + 2A^{7} + A^{6} - 2A^{5} - A^{4} + 2A^{3} - A^{2} - A + 1$}
\end{minipage}

\noindent{\color{gray!40}\rule{\textwidth}{0.4pt}}
\vspace{0.9\baselineskip}
\noindent \begin{minipage}[t]{0.25\textwidth}
\vspace{0pt}
\centering
\includegraphics[page=228,width=\linewidth]{knotoids.pdf}
\end{minipage}
\hfill
\begin{minipage}[t]{0.73\textwidth}
\vspace{0pt}
\raggedright
\textbf{Name:} {\large{$\mathbf{K7_{109}}$}} (chiral, non-rotatable$^{*}$) \\ \textbf{PD:} {\scriptsize\texttt{[0],[0,1,2,3],[1,4,5,2],[3,6,7,8],[8,7,9,4],[10,11,6,5],[9,12,13,10],[11,13,12,14],[14]}} \\ \textbf{EM:} {\scriptsize\texttt{(B0, A0C0C3D0, B1E3F3B2, B3F2E1E0, D3D2G0C1, G3H0D1C2, E2H2H1F0, F1G2G1I0, H3)}} \\ \textbf{Kauffman bracket:} {\scriptsize $A^{22} - 2A^{18} + 4A^{14} + A^{12} - 4A^{10} - 3A^{8} + 3A^{6} + 3A^{4} - A^{2} - 1$} \\ \textbf{Arrow:} {\scriptsize $A^{4} - 2 + 4/A^{4} + L_1/A^{6} - 4/A^{8} - 3L_1/A^{10} + 3/A^{12} + 3L_1/A^{14} - 1/A^{16} - L_1/A^{18}$} \\ \textbf{Mock:} {\scriptsize $-2w^{4} - 2w^{3} + 5w^{2} + 4w - 4 - 2/w + 3/w^{2} - 1/w^{4}$} \\ \textbf{Affine:} {\scriptsize $0$} \\ \textbf{Yamada:} {\scriptsize $-2A^{26} + A^{25} + 4A^{24} - 4A^{23} - 3A^{22} + 7A^{21} - 3A^{20} - 8A^{19} + 4A^{18} - 6A^{16} + A^{15} + 2A^{14} + A^{13} - 4A^{12} + 5A^{11} + 5A^{10} - 9A^{9} + 3A^{8} + 5A^{7} - 6A^{6} - A^{5} + 3A^{4} - 2A^{2} + 1$}
\end{minipage}

\noindent{\color{gray!40}\rule{\textwidth}{0.4pt}}
\vspace{0.9\baselineskip}
\noindent \begin{minipage}[t]{0.25\textwidth}
\vspace{0pt}
\centering
\includegraphics[page=229,width=\linewidth]{knotoids.pdf}
\end{minipage}
\hfill
\begin{minipage}[t]{0.73\textwidth}
\vspace{0pt}
\raggedright
\textbf{Name:} {\large{$\mathbf{K7_{110}}$}} (chiral, non-rotatable$^{*}$) \\ \textbf{PD:} {\scriptsize\texttt{[0],[0,1,2,3],[1,4,5,2],[3,6,7,8],[9,10,11,4],[12,13,6,5],[14,9,8,7],[10,14,12,11],[13]}} \\ \textbf{EM:} {\scriptsize\texttt{(B0, A0C0C3D0, B1E3F3B2, B3F2G3G2, G1H0H3C1, H2I0D1C2, H1E0D3D2, E1G0F0E2, F1)}} \\ \textbf{Kauffman bracket:} {\scriptsize $A^{22} - A^{18} + 3A^{14} - 4A^{10} - 2A^{8} + 3A^{6} + 3A^{4} - A^{2} - 1$} \\ \textbf{Arrow:} {\scriptsize $A^{4} - 1 + 3/A^{4} - 4/A^{8} - 2L_1/A^{10} + 3/A^{12} + 3L_1/A^{14} - 1/A^{16} - L_1/A^{18}$} \\ \textbf{Mock:} {\scriptsize $-2w^{4} - 2w^{3} + 4w^{2} + 3w - 3 - 1/w + 3/w^{2} - 1/w^{4}$} \\ \textbf{Affine:} {\scriptsize $-t + 2 - 1/t$} \\ \textbf{Yamada:} {\scriptsize $-A^{26} + A^{25} + A^{24} - 4A^{23} - A^{22} + 5A^{21} - 4A^{20} - 4A^{19} + 5A^{18} - A^{17} - 4A^{16} + A^{15} - A^{13} - 3A^{12} + 4A^{11} + 2A^{10} - 6A^{9} + 5A^{8} + 3A^{7} - 5A^{6} + A^{5} + 2A^{4} - A^{3} - 2A^{2} + 1$}
\end{minipage}

\noindent{\color{gray!40}\rule{\textwidth}{0.4pt}}
\vspace{0.9\baselineskip}
\noindent \begin{minipage}[t]{0.25\textwidth}
\vspace{0pt}
\centering
\includegraphics[page=230,width=\linewidth]{knotoids.pdf}
\end{minipage}
\hfill
\begin{minipage}[t]{0.73\textwidth}
\vspace{0pt}
\raggedright
\textbf{Name:} {\large{$\mathbf{K7_{111}}$}} (chiral, non-rotatable$^{*}$) \\ \textbf{PD:} {\scriptsize\texttt{[0],[0,1,2,3],[1,4,5,2],[6,7,8,3],[4,8,7,9],[5,10,11,6],[12,13,10,9],[13,12,14,11],[14]}} \\ \textbf{EM:} {\scriptsize\texttt{(B0, A0C0C3D3, B1E0F0B2, F3E2E1B3, C1D2D1G3, C2G2H3D0, H1H0F1E3, G1G0I0F2, H2)}} \\ \textbf{Kauffman bracket:} {\scriptsize $A^{24} + A^{22} - A^{18} - 2A^{16} + 2A^{12} + A^{10} - A^{8} - A^{6} + A^{2}$} \\ \textbf{Arrow:} {\scriptsize $A^{6}L_1 + A^{4} - 1 - 2L_1/A^{2} + 2L_1/A^{6} + A^{-8} - L_1/A^{10} - 1/A^{12} + A^{-16}$} \\ \textbf{Mock:} {\scriptsize $w^{3} + w^{2} + w + 1 - 3/w - 1/w^{2} + w^{-3}$} \\ \textbf{Affine:} {\scriptsize $2t - 4 + 2/t$} \\ \textbf{Yamada:} {\scriptsize $A^{27} - A^{25} + A^{24} + A^{23} - A^{22} - A^{21} + 2A^{19} - 2A^{18} - A^{17} + 4A^{16} + 3A^{13} + A^{9} - A^{8} - 3A^{7} + 2A^{6} + A^{5} - 2A^{4} + A^{3} + 2A^{2} - 1$}
\end{minipage}

\noindent{\color{gray!40}\rule{\textwidth}{0.4pt}}
\vspace{0.9\baselineskip}
\noindent \begin{minipage}[t]{0.25\textwidth}
\vspace{0pt}
\centering
\includegraphics[page=231,width=\linewidth]{knotoids.pdf}
\end{minipage}
\hfill
\begin{minipage}[t]{0.73\textwidth}
\vspace{0pt}
\raggedright
\textbf{Name:} {\large{$\mathbf{K7_{112}}$}} (chiral, non-rotatable$^{*}$) \\ \textbf{PD:} {\scriptsize\texttt{[0],[0,1,2,3],[1,4,5,2],[6,7,8,3],[4,8,7,9],[10,11,6,5],[9,12,13,10],[11,13,12,14],[14]}} \\ \textbf{EM:} {\scriptsize\texttt{(B0, A0C0C3D3, B1E0F3B2, F2E2E1B3, C1D2D1G0, G3H0D0C2, E3H2H1F0, F1G2G1I0, H3)}} \\ \textbf{Kauffman bracket:} {\scriptsize $A^{22} + A^{20} - 3A^{18} - 3A^{16} + 3A^{14} + 5A^{12} - A^{10} - 4A^{8} + 3A^{4} - 1$} \\ \textbf{Arrow:} {\scriptsize $A^{10}L_1 + A^{8} - 3A^{6}L_1 - 3A^{4} + 3A^{2}L_1 + 5 - L_1/A^{2} - 4/A^{4} + 3/A^{8} - 1/A^{12}$} \\ \textbf{Mock:} {\scriptsize $-2w^{3} - 3w^{2} + 4w + 8 - 2/w - 5/w^{2} + w^{-4}$} \\ \textbf{Affine:} {\scriptsize $0$} \\ \textbf{Yamada:} {\scriptsize $-A^{26} - A^{25} + 2A^{24} + 2A^{23} - 3A^{22} - A^{21} + 7A^{20} - A^{19} - 8A^{18} + 7A^{17} + 2A^{16} - 10A^{15} + 4A^{14} + 4A^{13} - 4A^{12} + A^{11} + 4A^{10} + 4A^{9} - 6A^{8} + 2A^{7} + 9A^{6} - 8A^{5} - 2A^{4} + 8A^{3} - 4A^{2} - 3A + 2$}
\end{minipage}

\noindent{\color{gray!40}\rule{\textwidth}{0.4pt}}
\vspace{0.9\baselineskip}
\noindent \begin{minipage}[t]{0.25\textwidth}
\vspace{0pt}
\centering
\includegraphics[page=232,width=\linewidth]{knotoids.pdf}
\end{minipage}
\hfill
\begin{minipage}[t]{0.73\textwidth}
\vspace{0pt}
\raggedright
\textbf{Name:} {\large{$\mathbf{K7_{113}}$}} (chiral, non-rotatable$^{*}$) \\ \textbf{PD:} {\scriptsize\texttt{[0],[0,1,2,3],[1,4,5,2],[6,7,8,3],[4,9,10,11],[12,13,6,5],[7,14,9,8],[14,12,11,10],[13]}} \\ \textbf{EM:} {\scriptsize\texttt{(B0, A0C0C3D3, B1E0F3B2, F2G0G3B3, C1G2H3H2, H1I0D0C2, D1H0E1D2, G1F0E3E2, F1)}} \\ \textbf{Kauffman bracket:} {\scriptsize $2A^{22} + 2A^{20} - 3A^{18} - 4A^{16} + A^{14} + 5A^{12} - 3A^{8} + 2A^{4} - 1$} \\ \textbf{Arrow:} {\scriptsize $2L_1/A^{2} + 2/A^{4} - 3L_1/A^{6} - 4/A^{8} + L_1/A^{10} + 5/A^{12} - 3/A^{16} + 2/A^{20} - 1/A^{24}$} \\ \textbf{Mock:} {\scriptsize $2w^{3} + 4w^{2} - 3w - 7 + 1/w + 5/w^{2} - 1/w^{4}$} \\ \textbf{Affine:} {\scriptsize $t - 2 + 1/t$} \\ \textbf{Yamada:} {\scriptsize $A^{27} - 2A^{25} + 3A^{23} - A^{22} - 4A^{21} + 5A^{20} + 5A^{19} - 7A^{18} + 3A^{17} + 7A^{16} - 6A^{15} + 5A^{13} - 2A^{12} - 2A^{11} + 5A^{9} - 4A^{8} - 4A^{7} + 9A^{6} - 3A^{5} - 6A^{4} + 5A^{3} - 2A + 1$}
\end{minipage}

\noindent{\color{gray!40}\rule{\textwidth}{0.4pt}}
\vspace{0.9\baselineskip}
\noindent \begin{minipage}[t]{0.25\textwidth}
\vspace{0pt}
\centering
\includegraphics[page=233,width=\linewidth]{knotoids.pdf}
\end{minipage}
\hfill
\begin{minipage}[t]{0.73\textwidth}
\vspace{0pt}
\raggedright
\textbf{Name:} {\large{$\mathbf{K7_{114}}$}} (chiral, non-rotatable$^{*}$) \\ \textbf{PD:} {\scriptsize\texttt{[0],[0,1,2,3],[1,4,5,2],[3,6,7,8],[4,9,10,11],[11,12,6,5],[12,13,14,7],[14,10,9,8],[13]}} \\ \textbf{EM:} {\scriptsize\texttt{(B0, A0C0C3D0, B1E0F3B2, B3F2G3H3, C1H2H1F0, E3G0D1C2, F1I0H0D2, G2E2E1D3, G1)}} \\ \textbf{Kauffman bracket:} {\scriptsize $A^{22} + A^{20} - 3A^{18} - 3A^{16} + 3A^{14} + 5A^{12} - 4A^{8} + 2A^{4} - 1$} \\ \textbf{Arrow:} {\scriptsize $A^{4} + A^{2}L_1 - L_2 - 2 - 3L_1/A^{2} + L_2/A^{4} + 2/A^{4} + 5L_1/A^{6} - 4L_1/A^{10} + 2L_1/A^{14} - L_1/A^{18}$} \\ \textbf{Mock:} {\scriptsize $-w^{4} - 2w^{3} + 3w^{2} + 6w - 1 - 5/w + w^{-2} + 2/w^{3} - 1/w^{4} - 1/w^{5}$} \\ \textbf{Affine:} {\scriptsize $-t^{2} + 2t - 2 + 2/t - 1/t^{2}$} \\ \textbf{Yamada:} {\scriptsize $-A^{26} - A^{25} + 2A^{24} + 2A^{23} - 3A^{22} - 2A^{21} + 6A^{20} - 8A^{18} + 6A^{17} + 5A^{16} - 7A^{15} + 4A^{14} + 6A^{13} - 3A^{12} + 2A^{10} + 3A^{9} - 7A^{8} - A^{7} + 8A^{6} - 7A^{5} - 3A^{4} + 8A^{3} - A^{2} - 3A + 1$}
\end{minipage}

\noindent{\color{gray!40}\rule{\textwidth}{0.4pt}}
\vspace{0.9\baselineskip}
\noindent \begin{minipage}[t]{0.25\textwidth}
\vspace{0pt}
\centering
\includegraphics[page=234,width=\linewidth]{knotoids.pdf}
\end{minipage}
\hfill
\begin{minipage}[t]{0.73\textwidth}
\vspace{0pt}
\raggedright
\textbf{Name:} {\large{$\mathbf{K7_{115}}$}} (chiral, non-rotatable$^{*}$) \\ \textbf{PD:} {\scriptsize\texttt{[0],[0,1,2,3],[1,4,5,2],[3,6,7,8],[9,10,11,4],[11,12,6,5],[12,13,14,7],[8,14,10,9],[13]}} \\ \textbf{EM:} {\scriptsize\texttt{(B0, A0C0C3D0, B1E3F3B2, B3F2G3H0, H3H2F0C1, E2G0D1C2, F1I0H1D2, D3G2E1E0, G1)}} \\ \textbf{Kauffman bracket:} {\scriptsize $-A^{19} + 2A^{15} - 3A^{11} - A^{9} + 4A^{7} + 3A^{5} - 2A^{3} - 3A$} \\ \textbf{Arrow:} {\scriptsize $A^{22}L_1 - 2A^{18}L_1 + 3A^{14}L_1 + A^{12} - 4A^{10}L_1 - A^{8}L_2 - 2A^{8} + 2A^{6}L_1 + A^{4}L_2 + 2A^{4}$} \\ \textbf{Mock:} {\scriptsize $2w^{3} + w^{2} - 5w - 2 + 4/w + 3/w^{2} - 1/w^{3} - 1/w^{4}$} \\ \textbf{Affine:} {\scriptsize $-t^{2} + 2 - 1/t^{2}$} \\ \textbf{Yamada:} {\scriptsize $2A^{25} - A^{24} - 4A^{23} + 2A^{22} + 3A^{21} - 5A^{20} - A^{19} + 7A^{18} - 3A^{17} - 2A^{16} + 4A^{15} - A^{13} - A^{12} + 5A^{11} - 2A^{10} - 4A^{9} + 7A^{8} - 5A^{6} + 3A^{5} + 2A^{4} - 2A^{3} - A^{2} + 2A + 1$}
\end{minipage}

\noindent{\color{gray!40}\rule{\textwidth}{0.4pt}}
\vspace{0.9\baselineskip}
\noindent \begin{minipage}[t]{0.25\textwidth}
\vspace{0pt}
\centering
\includegraphics[page=235,width=\linewidth]{knotoids.pdf}
\end{minipage}
\hfill
\begin{minipage}[t]{0.73\textwidth}
\vspace{0pt}
\raggedright
\textbf{Name:} {\large{$\mathbf{K7_{116}}$}} (chiral, non-rotatable$^{*}$) \\ \textbf{PD:} {\scriptsize\texttt{[0],[0,1,2,3],[1,4,5,2],[6,7,8,3],[4,9,10,11],[5,11,12,6],[7,12,13,14],[14,10,9,8],[13]}} \\ \textbf{EM:} {\scriptsize\texttt{(B0, A0C0C3D3, B1E0F0B2, F3G0H3B3, C1H2H1F1, C2E3G1D0, D1F2I0H0, G3E2E1D2, G2)}} \\ \textbf{Kauffman bracket:} {\scriptsize $A^{12} + A^{10} - A^{8} - 2A^{6} - A^{4} + A^{2} + 2$} \\ \textbf{Arrow:} {\scriptsize $A^{12} + A^{10}L_1 - A^{8} - 2A^{6}L_1 - A^{4}L_2 + A^{2}L_1 + L_2 + 1$} \\ \textbf{Mock:} {\scriptsize $-w^{4} - w^{3} + w^{2} + 2w + 2 - 1/w - 1/w^{2}$} \\ \textbf{Affine:} {\scriptsize $-t^{2} + 2 - 1/t^{2}$} \\ \textbf{Yamada:} {\scriptsize $A^{23} - A^{21} - A^{17} + 2A^{15} - A^{14} - A^{13} + A^{12} - 2A^{11} - 2A^{10} - 2A^{8} + 2A^{5} + A^{4} - 2A^{3} + A^{2} - 2$}
\end{minipage}

\noindent{\color{gray!40}\rule{\textwidth}{0.4pt}}
\vspace{0.9\baselineskip}
\noindent \begin{minipage}[t]{0.25\textwidth}
\vspace{0pt}
\centering
\includegraphics[page=236,width=\linewidth]{knotoids.pdf}
\end{minipage}
\hfill
\begin{minipage}[t]{0.73\textwidth}
\vspace{0pt}
\raggedright
\textbf{Name:} {\large{$\mathbf{K7_{117}}$}} (chiral, non-rotatable$^{*}$) \\ \textbf{PD:} {\scriptsize\texttt{[0],[0,1,2,3],[1,4,5,2],[3,6,7,4],[5,8,9,10],[6,11,12,7],[8,12,11,13],[14,13,10,9],[14]}} \\ \textbf{EM:} {\scriptsize\texttt{(B0, A0C0C3D0, B1D3E0B2, B3F0F3C1, C2G0H3H2, D1G2G1D2, E1F2F1H1, I0G3E3E2, H0)}} \\ \textbf{Kauffman bracket:} {\scriptsize $A^{22} + 2A^{20} - 2A^{18} - 3A^{16} + A^{14} + 5A^{12} + A^{10} - 4A^{8} - 2A^{6} + A^{4} + A^{2}$} \\ \textbf{Arrow:} {\scriptsize $L_1/A^{2} + 2/A^{4} - 2L_1/A^{6} - 3/A^{8} + L_1/A^{10} + 5/A^{12} + L_1/A^{14} - 4/A^{16} - 2L_1/A^{18} + A^{-20} + L_1/A^{22}$} \\ \textbf{Mock:} {\scriptsize $-w^{4} + 4w^{2} - 5 + 4/w^{2} - 1/w^{4}$} \\ \textbf{Affine:} {\scriptsize $0$} \\ \textbf{Yamada:} {\scriptsize $-A^{26} + 2A^{25} + 2A^{24} - 5A^{23} + 2A^{22} + 7A^{21} - 7A^{20} + A^{19} + 10A^{18} - 4A^{17} - A^{16} + 6A^{15} - 2A^{14} - 3A^{13} - 4A^{12} + 4A^{11} - 2A^{10} - 7A^{9} + 10A^{8} - 7A^{6} + 6A^{5} + 2A^{4} - 4A^{3} + A^{2} + A - 1$}
\end{minipage}

\noindent{\color{gray!40}\rule{\textwidth}{0.4pt}}
\vspace{0.9\baselineskip}
\noindent \begin{minipage}[t]{0.25\textwidth}
\vspace{0pt}
\centering
\includegraphics[page=237,width=\linewidth]{knotoids.pdf}
\end{minipage}
\hfill
\begin{minipage}[t]{0.73\textwidth}
\vspace{0pt}
\raggedright
\textbf{Name:} {\large{$\mathbf{K7_{118}}$}} (chiral, rotatable) \\ \textbf{PD:} {\scriptsize\texttt{[0],[0,1,2,3],[1,4,5,2],[3,6,7,4],[8,9,10,5],[6,11,12,7],[8,12,11,13],[9,14,13,10],[14]}} \\ \textbf{EM:} {\scriptsize\texttt{(B0, A0C0C3D0, B1D3E3B2, B3F0F3C1, G0H0H3C2, D1G2G1D2, E0F2F1H2, E1I0G3E2, H1)}} \\ \textbf{Kauffman bracket:} {\scriptsize $A^{24} - A^{20} - 2A^{18} + A^{16} + 4A^{14} + A^{12} - 4A^{10} - 3A^{8} + 2A^{6} + 3A^{4} - 1$} \\ \textbf{Arrow:} {\scriptsize $A^{12} - A^{8} - 2A^{6}L_1 + A^{4} + 4A^{2}L_1 + 1 - 4L_1/A^{2} - 3/A^{4} + 2L_1/A^{6} + 3/A^{8} - 1/A^{12}$} \\ \textbf{Mock:} {\scriptsize $-w^{4} - 2w^{3} + 4w + 3 - 4/w - 2/w^{2} + 2/w^{3} + w^{-4}$} \\ \textbf{Affine:} {\scriptsize $-2t + 4 - 2/t$} \\ \textbf{Yamada:} {\scriptsize $-A^{25} - A^{24} + 4A^{23} - 6A^{21} + 6A^{20} + 5A^{19} - 10A^{18} + 2A^{17} + 7A^{16} - 7A^{15} + 6A^{13} - A^{11} + 7A^{9} - 5A^{8} - 6A^{7} + 11A^{6} - 4A^{5} - 7A^{4} + 7A^{3} - 2A + 1$}
\end{minipage}

\noindent{\color{gray!40}\rule{\textwidth}{0.4pt}}
\vspace{0.9\baselineskip}
\noindent \begin{minipage}[t]{0.25\textwidth}
\vspace{0pt}
\centering
\includegraphics[page=238,width=\linewidth]{knotoids.pdf}
\end{minipage}
\hfill
\begin{minipage}[t]{0.73\textwidth}
\vspace{0pt}
\raggedright
\textbf{Name:} {\large{$\mathbf{K7_{119}}$}} (chiral, rotatable) \\ \textbf{PD:} {\scriptsize\texttt{[0],[0,1,2,3],[1,4,5,2],[6,7,4,3],[8,9,10,5],[11,12,7,6],[12,11,13,8],[9,14,13,10],[14]}} \\ \textbf{EM:} {\scriptsize\texttt{(B0, A0C0C3D3, B1D2E3B2, F3F2C1B3, G3H0H3C2, G1G0D1D0, F1F0H2E0, E1I0G2E2, H1)}} \\ \textbf{Kauffman bracket:} {\scriptsize $-A^{26} - 2A^{24} + 4A^{20} + 3A^{18} - 4A^{16} - 4A^{14} + 2A^{12} + 4A^{10} - 2A^{6} + A^{2}$} \\ \textbf{Arrow:} {\scriptsize $-A^{32} - 2A^{30}L_1 + 4A^{26}L_1 + 3A^{24} - 4A^{22}L_1 - 4A^{20} + 2A^{18}L_1 + 4A^{16} - 2A^{12} + A^{8}$} \\ \textbf{Mock:} {\scriptsize $w^{4} - 2w^{3} - 2w^{2} + 4w + 3 - 4/w - 4/w^{2} + 2/w^{3} + 3/w^{4}$} \\ \textbf{Affine:} {\scriptsize $-2t + 4 - 2/t$} \\ \textbf{Yamada:} {\scriptsize $-A^{27} + A^{25} - 3A^{24} - A^{23} + 7A^{22} - 4A^{21} - 5A^{20} + 12A^{19} + A^{18} - 10A^{17} + 8A^{16} + 3A^{15} - 8A^{14} + A^{13} + 3A^{12} + 2A^{11} - 8A^{10} + 4A^{9} + 7A^{8} - 13A^{7} + A^{6} + 8A^{5} - 9A^{4} - 3A^{3} + 5A^{2} - A - 3$}
\end{minipage}

\noindent{\color{gray!40}\rule{\textwidth}{0.4pt}}
\vspace{0.9\baselineskip}
\noindent \begin{minipage}[t]{0.25\textwidth}
\vspace{0pt}
\centering
\includegraphics[page=239,width=\linewidth]{knotoids.pdf}
\end{minipage}
\hfill
\begin{minipage}[t]{0.73\textwidth}
\vspace{0pt}
\raggedright
\textbf{Name:} {\large{$\mathbf{K7_{120}}$}} (chiral, non-rotatable$^{*}$) \\ \textbf{PD:} {\scriptsize\texttt{[0],[0,1,2,3],[1,4,5,2],[3,6,7,4],[5,8,9,6],[10,11,12,7],[8,12,13,9],[14,13,11,10],[14]}} \\ \textbf{EM:} {\scriptsize\texttt{(B0, A0C0C3D0, B1D3E0B2, B3E3F3C1, C2G0G3D1, H3H2G1D2, E1F2H1E2, I0G2F1F0, H0)}} \\ \textbf{Kauffman bracket:} {\scriptsize $-A^{24} + A^{22} + 4A^{20} - 5A^{16} - 3A^{14} + 4A^{12} + 3A^{10} - 2A^{8} - 2A^{6} + A^{4} + A^{2}$} \\ \textbf{Arrow:} {\scriptsize $-A^{12} + A^{10}L_1 + 4A^{8} - 5A^{4} - 3A^{2}L_1 + 4 + 3L_1/A^{2} - 2/A^{4} - 2L_1/A^{6} + A^{-8} + L_1/A^{10}$} \\ \textbf{Mock:} {\scriptsize $w^{4} + w^{3} - 5w^{2} - 3w + 8 + 3/w - 3/w^{2} - 1/w^{3}$} \\ \textbf{Affine:} {\scriptsize $0$} \\ \textbf{Yamada:} {\scriptsize $-A^{26} + A^{25} + 2A^{24} - 5A^{23} + A^{22} + 9A^{21} - 9A^{20} - 3A^{19} + 10A^{18} - 8A^{17} - 4A^{16} + 6A^{15} - A^{14} - 2A^{13} - 2A^{12} + 6A^{11} - A^{10} - 9A^{9} + 9A^{8} + A^{7} - 11A^{6} + 6A^{5} + 3A^{4} - 6A^{3} + A^{2} + 2A - 1$}
\end{minipage}

\noindent{\color{gray!40}\rule{\textwidth}{0.4pt}}
\vspace{0.9\baselineskip}
\noindent \begin{minipage}[t]{0.25\textwidth}
\vspace{0pt}
\centering
\includegraphics[page=240,width=\linewidth]{knotoids.pdf}
\end{minipage}
\hfill
\begin{minipage}[t]{0.73\textwidth}
\vspace{0pt}
\raggedright
\textbf{Name:} {\large{$\mathbf{K7_{121}}$}} (chiral, non-rotatable$^{*}$) \\ \textbf{PD:} {\scriptsize\texttt{[0],[0,1,2,3],[1,4,5,2],[3,6,7,4],[5,7,8,9],[6,10,11,12],[12,13,9,8],[10,13,14,11],[14]}} \\ \textbf{EM:} {\scriptsize\texttt{(B0, A0C0C3D0, B1D3E0B2, B3F0E1C1, C2D2G3G2, D1H0H3G0, F3H1E3E2, F1G1I0F2, H2)}} \\ \textbf{Kauffman bracket:} {\scriptsize $A^{22} + 3A^{20} - 4A^{16} - 3A^{14} + 3A^{12} + 3A^{10} - A^{8} - 2A^{6} + A^{2}$} \\ \textbf{Arrow:} {\scriptsize $L_1/A^{2} + 3/A^{4} - 4/A^{8} - 3L_1/A^{10} + 3/A^{12} + 3L_1/A^{14} - 1/A^{16} - 2L_1/A^{18} + L_1/A^{22}$} \\ \textbf{Mock:} {\scriptsize $-w^{4} + 4w^{2} + 2w - 4 - 4/w + 2/w^{2} + 2/w^{3}$} \\ \textbf{Affine:} {\scriptsize $0$} \\ \textbf{Yamada:} {\scriptsize $A^{25} - A^{24} - A^{23} + 4A^{22} - 5A^{20} + 7A^{19} + 4A^{18} - 7A^{17} + 4A^{16} + 4A^{15} - 5A^{14} + A^{12} + 2A^{11} - 5A^{10} + 7A^{8} - 7A^{7} - A^{6} + 7A^{5} - 3A^{4} - 3A^{3} + 4A^{2} + A - 2$}
\end{minipage}

\noindent{\color{gray!40}\rule{\textwidth}{0.4pt}}
\vspace{0.9\baselineskip}
\noindent \begin{minipage}[t]{0.25\textwidth}
\vspace{0pt}
\centering
\includegraphics[page=241,width=\linewidth]{knotoids.pdf}
\end{minipage}
\hfill
\begin{minipage}[t]{0.73\textwidth}
\vspace{0pt}
\raggedright
\textbf{Name:} {\large{$\mathbf{K7_{122}}$}} (chiral, non-rotatable$^{*}$) \\ \textbf{PD:} {\scriptsize\texttt{[0],[0,1,2,3],[1,4,5,2],[3,6,7,4],[5,8,9,10],[6,11,8,7],[11,12,13,9],[14,13,12,10],[14]}} \\ \textbf{EM:} {\scriptsize\texttt{(B0, A0C0C3D0, B1D3E0B2, B3F0F3C1, C2F2G3H3, D1G0E1D2, F1H2H1E2, I0G2G1E3, H0)}} \\ \textbf{Kauffman bracket:} {\scriptsize $A^{22} - A^{20} - 3A^{18} - A^{16} + 5A^{14} + 4A^{12} - 3A^{10} - 4A^{8} + A^{6} + 3A^{4} - 1$} \\ \textbf{Arrow:} {\scriptsize $A^{16} - A^{14}L_1 - 3A^{12} - A^{10}L_1 + 5A^{8} + 4A^{6}L_1 - 3A^{4} - 4A^{2}L_1 + 1 + 3L_1/A^{2} - L_1/A^{6}$} \\ \textbf{Mock:} {\scriptsize $w^{4} - 5w^{2} - 4w + 6 + 6/w - 1/w^{2} - 2/w^{3}$} \\ \textbf{Affine:} {\scriptsize $-2t + 4 - 2/t$} \\ \textbf{Yamada:} {\scriptsize $-2A^{25} + 5A^{23} - 2A^{22} - 5A^{21} + 9A^{20} + 2A^{19} - 11A^{18} + 6A^{17} + 5A^{16} - 11A^{15} + A^{14} + 5A^{13} - 4A^{12} - 3A^{11} + A^{10} + 6A^{9} - 9A^{8} - 3A^{7} + 12A^{6} - 8A^{5} - 5A^{4} + 9A^{3} - 3A^{2} - 3A + 2$}
\end{minipage}

\noindent{\color{gray!40}\rule{\textwidth}{0.4pt}}
\vspace{0.9\baselineskip}
\noindent \begin{minipage}[t]{0.25\textwidth}
\vspace{0pt}
\centering
\includegraphics[page=242,width=\linewidth]{knotoids.pdf}
\end{minipage}
\hfill
\begin{minipage}[t]{0.73\textwidth}
\vspace{0pt}
\raggedright
\textbf{Name:} {\large{$\mathbf{K7_{123}}$}} (chiral, non-rotatable$^{*}$) \\ \textbf{PD:} {\scriptsize\texttt{[0],[0,1,2,3],[1,4,5,2],[3,6,7,4],[5,7,8,9],[9,10,11,6],[12,13,14,8],[10,14,13,11],[12]}} \\ \textbf{EM:} {\scriptsize\texttt{(B0, A0C0C3D0, B1D3E0B2, B3F3E1C1, C2D2G3F0, E3H0H3D1, I0H2H1E2, F1G2G1F2, G0)}} \\ \textbf{Kauffman bracket:} {\scriptsize $A^{22} + 2A^{20} - A^{18} - 3A^{16} + 4A^{12} + A^{10} - 3A^{8} - 2A^{6} + A^{4} + A^{2}$} \\ \textbf{Arrow:} {\scriptsize $L_1/A^{2} + 2/A^{4} - L_1/A^{6} - 3/A^{8} + 4/A^{12} + L_1/A^{14} - 3/A^{16} - 2L_1/A^{18} + A^{-20} + L_1/A^{22}$} \\ \textbf{Mock:} {\scriptsize $-w^{4} + 4w^{2} + w - 4 - 1/w + 3/w^{2} - 1/w^{4}$} \\ \textbf{Affine:} {\scriptsize $t - 2 + 1/t$} \\ \textbf{Yamada:} {\scriptsize $-A^{26} + A^{25} + 2A^{24} - 2A^{23} + A^{22} + 7A^{21} - 3A^{20} - 2A^{19} + 7A^{18} - 3A^{17} - 3A^{16} + 4A^{15} - A^{14} - A^{13} - 3A^{12} + 3A^{11} - 6A^{9} + 6A^{8} + 2A^{7} - 5A^{6} + 3A^{5} + 3A^{4} - 3A^{3} + A - 1$}
\end{minipage}

\noindent{\color{gray!40}\rule{\textwidth}{0.4pt}}
\vspace{0.9\baselineskip}
\noindent \begin{minipage}[t]{0.25\textwidth}
\vspace{0pt}
\centering
\includegraphics[page=243,width=\linewidth]{knotoids.pdf}
\end{minipage}
\hfill
\begin{minipage}[t]{0.73\textwidth}
\vspace{0pt}
\raggedright
\textbf{Name:} {\large{$\mathbf{K7_{124}}$}} (chiral, non-rotatable$^{*}$) \\ \textbf{PD:} {\scriptsize\texttt{[0],[0,1,2,3],[1,4,5,2],[3,6,7,4],[8,9,6,5],[7,10,11,12],[12,13,9,8],[10,14,13,11],[14]}} \\ \textbf{EM:} {\scriptsize\texttt{(B0, A0C0C3D0, B1D3E3B2, B3E2F0C1, G3G2D1C2, D2H0H3G0, F3H2E1E0, F1I0G1F2, H1)}} \\ \textbf{Kauffman bracket:} {\scriptsize $A^{22} - 2A^{18} - A^{16} + 4A^{14} + 2A^{12} - 4A^{10} - 4A^{8} + 2A^{6} + 4A^{4} - 1$} \\ \textbf{Arrow:} {\scriptsize $A^{16} - 2A^{12} - A^{10}L_1 + 4A^{8} + 2A^{6}L_1 - 4A^{4} - 4A^{2}L_1 + 2 + 4L_1/A^{2} - L_1/A^{6}$} \\ \textbf{Mock:} {\scriptsize $w^{4} + w^{3} - 5w^{2} - 5w + 6 + 5/w - 1/w^{2} - 1/w^{3}$} \\ \textbf{Affine:} {\scriptsize $-2t + 4 - 2/t$} \\ \textbf{Yamada:} {\scriptsize $A^{27} - A^{26} - 2A^{25} + 4A^{24} - A^{23} - 7A^{22} + 8A^{21} + 3A^{20} - 10A^{19} + 5A^{18} + 4A^{17} - 6A^{16} - A^{15} + A^{14} + A^{13} - 7A^{12} + 8A^{10} - 8A^{9} - A^{8} + 10A^{7} - 5A^{6} - 4A^{5} + 4A^{4} - 3A^{2} + 1$}
\end{minipage}

\noindent{\color{gray!40}\rule{\textwidth}{0.4pt}}
\vspace{0.9\baselineskip}
\noindent \begin{minipage}[t]{0.25\textwidth}
\vspace{0pt}
\centering
\includegraphics[page=244,width=\linewidth]{knotoids.pdf}
\end{minipage}
\hfill
\begin{minipage}[t]{0.73\textwidth}
\vspace{0pt}
\raggedright
\textbf{Name:} {\large{$\mathbf{K7_{125}}$}} (chiral, non-rotatable$^{*}$) \\ \textbf{PD:} {\scriptsize\texttt{[0],[0,1,2,3],[1,4,5,2],[3,6,7,4],[7,8,9,5],[6,9,10,11],[8,12,13,14],[14,13,11,10],[12]}} \\ \textbf{EM:} {\scriptsize\texttt{(B0, A0C0C3D0, B1D3E3B2, B3F0E0C1, D2G0F1C2, D1E2H3H2, E1I0H1H0, G3G2F3F2, G1)}} \\ \textbf{Kauffman bracket:} {\scriptsize $A^{22} - 2A^{18} - 2A^{16} + 2A^{14} + 3A^{12} - A^{10} - 3A^{8} + 2A^{4} + A^{2}$} \\ \textbf{Arrow:} {\scriptsize $A^{28} - 2A^{24} - 2A^{22}L_1 + 2A^{20} + 3A^{18}L_1 - A^{16} - 3A^{14}L_1 + 2A^{10}L_1 + A^{8}$} \\ \textbf{Mock:} {\scriptsize $-w^{4} - 2w^{3} + 2w^{2} + 3w - 2 - 3/w + w^{-2} + 2/w^{3} + w^{-4}$} \\ \textbf{Affine:} {\scriptsize $-3t + 6 - 3/t$} \\ \textbf{Yamada:} {\scriptsize $A^{24} - 3A^{22} + 3A^{21} + 2A^{20} - 7A^{19} + 3A^{18} + 2A^{17} - 4A^{16} + A^{15} + A^{14} + A^{13} - 3A^{12} + A^{11} + 4A^{10} - 3A^{9} + 6A^{7} - 2A^{6} - 2A^{5} + 4A^{4} + A^{3} - 2A^{2} + A + 1$}
\end{minipage}

\noindent{\color{gray!40}\rule{\textwidth}{0.4pt}}
\vspace{0.9\baselineskip}
\noindent \begin{minipage}[t]{0.25\textwidth}
\vspace{0pt}
\centering
\includegraphics[page=245,width=\linewidth]{knotoids.pdf}
\end{minipage}
\hfill
\begin{minipage}[t]{0.73\textwidth}
\vspace{0pt}
\raggedright
\textbf{Name:} {\large{$\mathbf{K7_{126}}$}} (chiral, non-rotatable$^{*}$) \\ \textbf{PD:} {\scriptsize\texttt{[0],[0,1,2,3],[1,4,5,2],[3,5,6,7],[4,8,9,10],[11,8,7,6],[9,12,13,14],[10,13,12,11],[14]}} \\ \textbf{EM:} {\scriptsize\texttt{(B0, A0C0C3D0, B1E0D1B2, B3C2F3F2, C1F1G0H0, H3E1D3D2, E2H2H1I0, E3G2G1F0, G3)}} \\ \textbf{Kauffman bracket:} {\scriptsize $-A^{24} + 2A^{20} + A^{18} - 4A^{16} - 2A^{14} + 4A^{12} + 4A^{10} - A^{8} - 3A^{6} + A^{2}$} \\ \textbf{Arrow:} {\scriptsize $-A^{18}L_1 + 2A^{14}L_1 + A^{12} - 4A^{10}L_1 - 2A^{8} + 4A^{6}L_1 + 4A^{4} - A^{2}L_1 - 3 + A^{-4}$} \\ \textbf{Mock:} {\scriptsize $w^{3} + w^{2} - 5w - 4 + 5/w + 5/w^{2} - 1/w^{3} - 1/w^{4}$} \\ \textbf{Affine:} {\scriptsize $-2t + 4 - 2/t$} \\ \textbf{Yamada:} {\scriptsize $A^{26} + A^{25} - 2A^{24} - 2A^{23} + 4A^{22} - A^{21} - 8A^{20} + 5A^{19} + 4A^{18} - 9A^{17} + 4A^{16} + 6A^{15} - 3A^{14} + A^{13} + 2A^{12} + 3A^{11} - 6A^{10} + 8A^{8} - 6A^{7} - 2A^{6} + 9A^{5} - 2A^{4} - 4A^{3} + 3A^{2} + A - 1$}
\end{minipage}

\noindent{\color{gray!40}\rule{\textwidth}{0.4pt}}
\vspace{0.9\baselineskip}
\noindent \begin{minipage}[t]{0.25\textwidth}
\vspace{0pt}
\centering
\includegraphics[page=246,width=\linewidth]{knotoids.pdf}
\end{minipage}
\hfill
\begin{minipage}[t]{0.73\textwidth}
\vspace{0pt}
\raggedright
\textbf{Name:} {\large{$\mathbf{K7_{127}}$}} (chiral, non-rotatable$^{*}$) \\ \textbf{PD:} {\scriptsize\texttt{[0],[0,1,2,3],[1,4,5,2],[3,5,6,7],[8,9,10,4],[11,8,7,6],[12,13,14,9],[13,12,11,10],[14]}} \\ \textbf{EM:} {\scriptsize\texttt{(B0, A0C0C3D0, B1E3D1B2, B3C2F3F2, F1G3H3C1, H2E0D3D2, H1H0I0E1, G1G0F0E2, G2)}} \\ \textbf{Kauffman bracket:} {\scriptsize $A^{22} + A^{20} - A^{18} - 2A^{16} + A^{14} + 3A^{12} + A^{10} - 3A^{8} - 2A^{6} + A^{2} + 1$} \\ \textbf{Arrow:} {\scriptsize $A^{4} + A^{2}L_1 - 1 - 2L_1/A^{2} + A^{-4} + 3L_1/A^{6} + A^{-8} - 3L_1/A^{10} - 2/A^{12} + A^{-16} + L_1/A^{18}$} \\ \textbf{Mock:} {\scriptsize $-w^{4} - w^{3} + w^{2} + 3w + 2 - 3/w - 1/w^{2} + w^{-3}$} \\ \textbf{Affine:} {\scriptsize $0$} \\ \textbf{Yamada:} {\scriptsize $-A^{28} + A^{26} - A^{25} - 2A^{24} + 3A^{23} - 6A^{21} + A^{20} + 2A^{19} - 4A^{18} + 2A^{16} + 2A^{15} - 3A^{14} + A^{13} + 4A^{12} - 5A^{11} + 3A^{9} - 3A^{8} - 3A^{7} + 2A^{6} + 2A^{5} - 2A^{4} + 2A^{2} - 1$}
\end{minipage}

\noindent{\color{gray!40}\rule{\textwidth}{0.4pt}}
\vspace{0.9\baselineskip}
\noindent \begin{minipage}[t]{0.25\textwidth}
\vspace{0pt}
\centering
\includegraphics[page=247,width=\linewidth]{knotoids.pdf}
\end{minipage}
\hfill
\begin{minipage}[t]{0.73\textwidth}
\vspace{0pt}
\raggedright
\textbf{Name:} {\large{$\mathbf{K7_{128}}$}} (chiral, non-rotatable$^{*}$) \\ \textbf{PD:} {\scriptsize\texttt{[0],[0,1,2,3],[1,4,5,2],[6,7,4,3],[5,7,8,9],[6,10,11,12],[12,13,9,8],[10,13,14,11],[14]}} \\ \textbf{EM:} {\scriptsize\texttt{(B0, A0C0C3D3, B1D2E0B2, F0E1C1B3, C2D1G3G2, D0H0H3G0, F3H1E3E2, F1G1I0F2, H2)}} \\ \textbf{Kauffman bracket:} {\scriptsize $-A^{21} - 2A^{19} + 3A^{15} + 3A^{13} - 2A^{11} - 3A^{9} + 2A^{5} - A$} \\ \textbf{Arrow:} {\scriptsize $A^{6}L_1 + 2A^{4} - 3 - 3L_1/A^{2} + 2/A^{4} + 3L_1/A^{6} - 2L_1/A^{10} + L_1/A^{14}$} \\ \textbf{Mock:} {\scriptsize $-w^{4} + 3w^{2} + 2w - 2 - 4/w + w^{-2} + 2/w^{3}$} \\ \textbf{Affine:} {\scriptsize $0$} \\ \textbf{Yamada:} {\scriptsize $A^{25} - 2A^{23} + 3A^{21} - 2A^{20} - 2A^{19} + 5A^{18} + A^{17} - 4A^{16} + 3A^{15} + 3A^{14} - 4A^{13} + A^{12} + 3A^{11} - 2A^{10} - A^{9} + A^{8} + 4A^{7} - 4A^{6} - A^{5} + 6A^{4} - 4A^{3} - A^{2} + 3A - 1$}
\end{minipage}

\noindent{\color{gray!40}\rule{\textwidth}{0.4pt}}
\vspace{0.9\baselineskip}
\noindent \begin{minipage}[t]{0.25\textwidth}
\vspace{0pt}
\centering
\includegraphics[page=248,width=\linewidth]{knotoids.pdf}
\end{minipage}
\hfill
\begin{minipage}[t]{0.73\textwidth}
\vspace{0pt}
\raggedright
\textbf{Name:} {\large{$\mathbf{K7_{129}}$}} (chiral, non-rotatable$^{*}$) \\ \textbf{PD:} {\scriptsize\texttt{[0],[0,1,2,3],[1,4,5,2],[6,7,4,3],[7,8,9,5],[6,9,10,11],[8,12,13,14],[14,13,11,10],[12]}} \\ \textbf{EM:} {\scriptsize\texttt{(B0, A0C0C3D3, B1D2E3B2, F0E0C1B3, D1G0F1C2, D0E2H3H2, E1I0H1H0, G3G2F3F2, G1)}} \\ \textbf{Kauffman bracket:} {\scriptsize $-2A^{24} - A^{22} + 3A^{20} + 3A^{18} - 3A^{16} - 4A^{14} + 2A^{12} + 4A^{10} - 2A^{6} + A^{2}$} \\ \textbf{Arrow:} {\scriptsize $-2A^{30}L_1 - A^{28} + 3A^{26}L_1 + 3A^{24} - 3A^{22}L_1 - 4A^{20} + 2A^{18}L_1 + 4A^{16} - 2A^{12} + A^{8}$} \\ \textbf{Mock:} {\scriptsize $w^{4} - 2w^{3} - 3w^{2} + 3w + 4 - 3/w - 4/w^{2} + 2/w^{3} + 3/w^{4}$} \\ \textbf{Affine:} {\scriptsize $-3t + 6 - 3/t$} \\ \textbf{Yamada:} {\scriptsize $-A^{27} + A^{25} - 2A^{24} - A^{23} + 5A^{22} - A^{21} - 5A^{20} + 7A^{19} + 3A^{18} - 9A^{17} + 5A^{16} + 4A^{15} - 5A^{14} + A^{13} + 2A^{12} + 2A^{11} - 7A^{10} + A^{9} + 6A^{8} - 8A^{7} - A^{6} + 7A^{5} - 5A^{4} - 4A^{3} + 3A^{2} - A - 3$}
\end{minipage}

\noindent{\color{gray!40}\rule{\textwidth}{0.4pt}}
\vspace{0.9\baselineskip}
\noindent \begin{minipage}[t]{0.25\textwidth}
\vspace{0pt}
\centering
\includegraphics[page=249,width=\linewidth]{knotoids.pdf}
\end{minipage}
\hfill
\begin{minipage}[t]{0.73\textwidth}
\vspace{0pt}
\raggedright
\textbf{Name:} {\large{$\mathbf{K7_{130}}$}} (chiral, non-rotatable$^{*}$) \\ \textbf{PD:} {\scriptsize\texttt{[0],[0,1,2,3],[1,4,5,2],[6,7,4,3],[7,8,9,5],[10,11,12,6],[8,12,13,9],[13,14,11,10],[14]}} \\ \textbf{EM:} {\scriptsize\texttt{(B0, A0C0C3D3, B1D2E3B2, F3E0C1B3, D1G0G3C2, H3H2G1D0, E1F2H0E2, G2I0F1F0, H1)}} \\ \textbf{Kauffman bracket:} {\scriptsize $-A^{26} - A^{24} + 2A^{22} + 3A^{20} - 2A^{18} - 6A^{16} + 6A^{12} + 3A^{10} - 3A^{8} - 2A^{6} + A^{4} + A^{2}$} \\ \textbf{Arrow:} {\scriptsize $-A^{32} - A^{30}L_1 + 2A^{28} + 3A^{26}L_1 - 2A^{24} - 6A^{22}L_1 + 6A^{18}L_1 + 3A^{16} - 3A^{14}L_1 - 2A^{12} + A^{10}L_1 + A^{8}$} \\ \textbf{Mock:} {\scriptsize $-2w^{3} + 2w^{2} + 6w - 2 - 8/w - 2/w^{2} + 4/w^{3} + 3/w^{4}$} \\ \textbf{Affine:} {\scriptsize $-2t + 4 - 2/t$} \\ \textbf{Yamada:} {\scriptsize $-2A^{27} + 5A^{25} - 3A^{24} - 7A^{23} + 10A^{22} + 3A^{21} - 14A^{20} + 7A^{19} + 11A^{18} - 13A^{17} + 2A^{16} + 10A^{15} - 5A^{14} - 4A^{13} + 2A^{12} + 8A^{11} - 11A^{10} - 4A^{9} + 15A^{8} - 10A^{7} - 10A^{6} + 11A^{5} - 2A^{4} - 7A^{3} + 3A^{2} + A - 2$}
\end{minipage}

\noindent{\color{gray!40}\rule{\textwidth}{0.4pt}}
\vspace{0.9\baselineskip}
\noindent \begin{minipage}[t]{0.25\textwidth}
\vspace{0pt}
\centering
\includegraphics[page=250,width=\linewidth]{knotoids.pdf}
\end{minipage}
\hfill
\begin{minipage}[t]{0.73\textwidth}
\vspace{0pt}
\raggedright
\textbf{Name:} {\large{$\mathbf{K7_{131}}$}} (chiral, non-rotatable$^{*}$) \\ \textbf{PD:} {\scriptsize\texttt{[0],[0,1,2,3],[1,4,5,2],[3,6,7,4],[5,7,8,9],[9,10,11,6],[12,13,10,8],[13,12,14,11],[14]}} \\ \textbf{EM:} {\scriptsize\texttt{(B0, A0C0C3D0, B1D3E0B2, B3F3E1C1, C2D2G3F0, E3G2H3D1, H1H0F1E2, G1G0I0F2, H2)}} \\ \textbf{Kauffman bracket:} {\scriptsize $A^{21} + 2A^{19} - 2A^{17} - 4A^{15} + 4A^{11} + A^{9} - 3A^{7} - 2A^{5} + A^{3} + A$} \\ \textbf{Arrow:} {\scriptsize $-A^{12} - 2A^{10}L_1 + 2A^{8} + 4A^{6}L_1 + A^{4}L_2 - A^{4} - 4A^{2}L_1 - 2L_2 + 1 + 3L_1/A^{2} + 2L_2/A^{4} - L_1/A^{6} - L_2/A^{8}$} \\ \textbf{Mock:} {\scriptsize $-w^{5} + 3w^{3} - 5w + 4/w + w^{-2} - 1/w^{3}$} \\ \textbf{Affine:} {\scriptsize $t^{2} - t - 1/t + t^{-2}$} \\ \textbf{Yamada:} {\scriptsize $A^{26} - 2A^{25} - A^{24} + 5A^{23} - 4A^{22} - 3A^{21} + 9A^{20} - 4A^{19} - 4A^{18} + 6A^{17} - A^{16} - 2A^{15} - 2A^{14} + 3A^{13} - 2A^{12} - 7A^{11} + 4A^{10} + A^{9} - 8A^{8} + 3A^{7} + 4A^{6} - 5A^{5} + 4A^{3} - A^{2} - A + 1$}
\end{minipage}

\noindent{\color{gray!40}\rule{\textwidth}{0.4pt}}
\vspace{0.9\baselineskip}
\noindent \begin{minipage}[t]{0.25\textwidth}
\vspace{0pt}
\centering
\includegraphics[page=251,width=\linewidth]{knotoids.pdf}
\end{minipage}
\hfill
\begin{minipage}[t]{0.73\textwidth}
\vspace{0pt}
\raggedright
\textbf{Name:} {\large{$\mathbf{K7_{132}}$}} (chiral, non-rotatable$^{*}$) \\ \textbf{PD:} {\scriptsize\texttt{[0],[0,1,2,3],[1,4,5,2],[3,6,7,4],[7,8,9,5],[6,9,10,11],[8,12,13,10],[11,13,12,14],[14]}} \\ \textbf{EM:} {\scriptsize\texttt{(B0, A0C0C3D0, B1D3E3B2, B3F0E0C1, D2G0F1C2, D1E2G3H0, E1H2H1F2, F3G2G1I0, H3)}} \\ \textbf{Kauffman bracket:} {\scriptsize $-A^{20} + 2A^{16} + 2A^{14} - 2A^{12} - 3A^{10} + A^{8} + 4A^{6} + A^{4} - 2A^{2} - 1$} \\ \textbf{Arrow:} {\scriptsize $-A^{2}L_1 + 2L_1/A^{2} + L_2/A^{4} + A^{-4} - 2L_1/A^{6} - 2L_2/A^{8} - 1/A^{8} + L_1/A^{10} + 2L_2/A^{12} + 2/A^{12} + L_1/A^{14} - L_2/A^{16} - 1/A^{16} - L_1/A^{18}$} \\ \textbf{Mock:} {\scriptsize $-w^{5} - w^{4} + 2w^{3} + 4w^{2} - w - 4 + 2/w^{2}$} \\ \textbf{Affine:} {\scriptsize $t^{2} - 2 + t^{-2}$} \\ \textbf{Yamada:} {\scriptsize $-A^{24} - A^{23} + 3A^{22} - 6A^{20} + 4A^{19} + 2A^{18} - 8A^{17} + 2A^{16} + 2A^{15} - 5A^{14} + 2A^{11} - 4A^{10} + 6A^{8} - 5A^{7} - A^{6} + 7A^{5} - 2A^{4} - 3A^{3} + 3A^{2} - 1$}
\end{minipage}

\noindent{\color{gray!40}\rule{\textwidth}{0.4pt}}
\vspace{0.9\baselineskip}
\noindent \begin{minipage}[t]{0.25\textwidth}
\vspace{0pt}
\centering
\includegraphics[page=252,width=\linewidth]{knotoids.pdf}
\end{minipage}
\hfill
\begin{minipage}[t]{0.73\textwidth}
\vspace{0pt}
\raggedright
\textbf{Name:} {\large{$\mathbf{K7_{133}}$}} (chiral, non-rotatable$^{*}$) \\ \textbf{PD:} {\scriptsize\texttt{[0],[0,1,2,3],[1,4,5,2],[3,5,6,7],[4,8,9,6],[7,10,11,8],[9,11,12,13],[13,12,14,10],[14]}} \\ \textbf{EM:} {\scriptsize\texttt{(B0, A0C0C3D0, B1E0D1B2, B3C2E3F0, C1F3G0D2, D3H3G1E1, E2F2H1H0, G3G2I0F1, H2)}} \\ \textbf{Kauffman bracket:} {\scriptsize $A^{21} + A^{19} - 2A^{17} - 3A^{15} + A^{13} + 4A^{11} - 3A^{7} - 2A^{5} + A^{3} + A$} \\ \textbf{Arrow:} {\scriptsize $-A^{12}L_2 - A^{10}L_1 + 2A^{8}L_2 + 3A^{6}L_1 - 2A^{4}L_2 + A^{4} - 4A^{2}L_1 + L_2 - 1 + 3L_1/A^{2} + 2/A^{4} - L_1/A^{6} - 1/A^{8}$} \\ \textbf{Mock:} {\scriptsize $w^{4} + w^{3} - 2w^{2} - 3w + 2 + 4/w + w^{-2} - 3/w^{3} - 2/w^{4} + w^{-5} + w^{-6}$} \\ \textbf{Affine:} {\scriptsize $-t^{2} + 2 - 1/t^{2}$} \\ \textbf{Yamada:} {\scriptsize $A^{26} - A^{25} + A^{24} + 3A^{23} - 5A^{22} + 6A^{20} - 6A^{19} - 2A^{18} + 5A^{17} - 3A^{16} - 2A^{15} - A^{14} + 2A^{13} - 3A^{12} - 4A^{11} + 6A^{10} - A^{9} - 5A^{8} + 5A^{7} + A^{6} - 5A^{5} + A^{4} + 3A^{3} - 2A^{2} - A + 1$}
\end{minipage}

\noindent{\color{gray!40}\rule{\textwidth}{0.4pt}}
\vspace{0.9\baselineskip}
\noindent \begin{minipage}[t]{0.25\textwidth}
\vspace{0pt}
\centering
\includegraphics[page=253,width=\linewidth]{knotoids.pdf}
\end{minipage}
\hfill
\begin{minipage}[t]{0.73\textwidth}
\vspace{0pt}
\raggedright
\textbf{Name:} {\large{$\mathbf{K7_{134}}$}} (chiral, non-rotatable$^{*}$) \\ \textbf{PD:} {\scriptsize\texttt{[0],[0,1,2,3],[1,4,5,2],[3,5,6,7],[4,8,9,6],[10,11,8,7],[11,12,13,9],[10,13,12,14],[14]}} \\ \textbf{EM:} {\scriptsize\texttt{(B0, A0C0C3D0, B1E0D1B2, B3C2E3F3, C1F2G3D2, H0G0E1D3, F1H2H1E2, F0G2G1I0, H3)}} \\ \textbf{Kauffman bracket:} {\scriptsize $2A^{22} + A^{20} - 3A^{18} - 3A^{16} + 4A^{14} + 5A^{12} - 2A^{10} - 5A^{8} + 3A^{4} - 1$} \\ \textbf{Arrow:} {\scriptsize $L_2/A^{8} + A^{-8} + L_1/A^{10} - 2L_2/A^{12} - 1/A^{12} - 3L_1/A^{14} + 2L_2/A^{16} + 2/A^{16} + 5L_1/A^{18} - L_2/A^{20} - 1/A^{20} - 5L_1/A^{22} + 3L_1/A^{26} - L_1/A^{30}$} \\ \textbf{Mock:} {\scriptsize $2w^{4} + 2w^{3} - 3w^{2} - 5w + 3 + 6/w - 2/w^{2} - 4/w^{3} + w^{-4} + w^{-5}$} \\ \textbf{Affine:} {\scriptsize $t^{2} + t - 4 + 1/t + t^{-2}$} \\ \textbf{Yamada:} {\scriptsize $A^{27} - 2A^{25} + 2A^{24} + 4A^{23} - 5A^{22} - A^{21} + 10A^{20} - 3A^{19} - 8A^{18} + 11A^{17} + 2A^{16} - 10A^{15} + 6A^{14} + 5A^{13} - 6A^{12} - A^{11} + 4A^{10} + 3A^{9} - 10A^{8} + 3A^{7} + 10A^{6} - 12A^{5} + 8A^{3} - 6A^{2} - 2A + 3$}
\end{minipage}

\noindent{\color{gray!40}\rule{\textwidth}{0.4pt}}
\vspace{0.9\baselineskip}
\noindent \begin{minipage}[t]{0.25\textwidth}
\vspace{0pt}
\centering
\includegraphics[page=254,width=\linewidth]{knotoids.pdf}
\end{minipage}
\hfill
\begin{minipage}[t]{0.73\textwidth}
\vspace{0pt}
\raggedright
\textbf{Name:} {\large{$\mathbf{K7_{135}}$}} (chiral, non-rotatable$^{*}$) \\ \textbf{PD:} {\scriptsize\texttt{[0],[0,1,2,3],[1,4,5,2],[3,6,7,8],[9,10,11,4],[5,12,13,6],[10,9,8,7],[11,13,14,12],[14]}} \\ \textbf{EM:} {\scriptsize\texttt{(B0, A0C0C3D0, B1E3F0B2, B3F3G3G2, G1G0H0C1, C2H3H1D1, E1E0D3D2, E2F2I0F1, H2)}} \\ \textbf{Kauffman bracket:} {\scriptsize $A^{14} + 2A^{12} - 2A^{8} - 2A^{6} + A^{2} + 1$} \\ \textbf{Arrow:} {\scriptsize $A^{2}L_1 + 2 - 2/A^{4} - 2L_1/A^{6} - L_2/A^{8} + A^{-8} + L_1/A^{10} + L_2/A^{12}$} \\ \textbf{Mock:} {\scriptsize $w^{3} - w + 2 - 1/w^{2}$} \\ \textbf{Affine:} {\scriptsize $-t^{2} + t + 1/t - 1/t^{2}$} \\ \textbf{Yamada:} {\scriptsize $A^{23} - A^{22} - 2A^{21} + A^{20} - 2A^{18} - A^{16} - 2A^{14} + A^{12} - 2A^{11} + A^{10} + A^{9} - A^{8} + A^{6} - A^{4} + A^{2} - 1$}
\end{minipage}

\noindent{\color{gray!40}\rule{\textwidth}{0.4pt}}
\vspace{0.9\baselineskip}
\noindent \begin{minipage}[t]{0.25\textwidth}
\vspace{0pt}
\centering
\includegraphics[page=255,width=\linewidth]{knotoids.pdf}
\end{minipage}
\hfill
\begin{minipage}[t]{0.73\textwidth}
\vspace{0pt}
\raggedright
\textbf{Name:} {\large{$\mathbf{K7_{136}}$}} (chiral, non-rotatable$^{*}$) \\ \textbf{PD:} {\scriptsize\texttt{[0],[0,1,2,3],[1,4,5,2],[3,6,7,4],[5,8,9,6],[10,11,8,7],[11,12,13,9],[14,13,12,10],[14]}} \\ \textbf{EM:} {\scriptsize\texttt{(B0, A0C0C3D0, B1D3E0B2, B3E3F3C1, C2F2G3D1, H3G0E1D2, F1H2H1E2, I0G2G1F0, H0)}} \\ \textbf{Kauffman bracket:} {\scriptsize $-A^{24} + A^{22} + 3A^{20} - 5A^{16} - A^{14} + 5A^{12} + 3A^{10} - 3A^{8} - 3A^{6} + A^{4} + A^{2}$} \\ \textbf{Arrow:} {\scriptsize $-A^{6}L_1 + A^{4} + 3A^{2}L_1 + L_2 - 1 - 5L_1/A^{2} - 2L_2/A^{4} + A^{-4} + 5L_1/A^{6} + 3L_2/A^{8} - 3L_1/A^{10} - 3L_2/A^{12} + L_1/A^{14} + L_2/A^{16}$} \\ \textbf{Mock:} {\scriptsize $w^{5} - 5w^{3} - 2w^{2} + 7w + 4 - 4/w - 1/w^{2} + w^{-3}$} \\ \textbf{Affine:} {\scriptsize $t^{2} - t - 1/t + t^{-2}$} \\ \textbf{Yamada:} {\scriptsize $-A^{26} + A^{25} + A^{24} - 5A^{23} + 3A^{22} + 6A^{21} - 11A^{20} + 2A^{19} + 10A^{18} - 11A^{17} - A^{16} + 7A^{15} - 4A^{14} - 2A^{13} - A^{12} + 5A^{11} - 5A^{10} - 7A^{9} + 11A^{8} - 3A^{7} - 9A^{6} + 10A^{5} + A^{4} - 7A^{3} + 4A^{2} + 2A - 2$}
\end{minipage}

\noindent{\color{gray!40}\rule{\textwidth}{0.4pt}}
\vspace{0.9\baselineskip}
\noindent \begin{minipage}[t]{0.25\textwidth}
\vspace{0pt}
\centering
\includegraphics[page=256,width=\linewidth]{knotoids.pdf}
\end{minipage}
\hfill
\begin{minipage}[t]{0.73\textwidth}
\vspace{0pt}
\raggedright
\textbf{Name:} {\large{$\mathbf{K7_{137}}$}} (chiral, non-rotatable$^{*}$) \\ \textbf{PD:} {\scriptsize\texttt{[0],[0,1,2,3],[1,4,5,2],[3,6,7,4],[8,9,6,5],[7,10,11,8],[9,11,12,13],[10,14,13,12],[14]}} \\ \textbf{EM:} {\scriptsize\texttt{(B0, A0C0C3D0, B1D3E3B2, B3E2F0C1, F3G0D1C2, D2H0G1E0, E1F2H3H2, F1I0G3G2, H1)}} \\ \textbf{Kauffman bracket:} {\scriptsize $A^{22} - 3A^{18} - A^{16} + 4A^{14} + 4A^{12} - 3A^{10} - 4A^{8} + A^{6} + 3A^{4} - 1$} \\ \textbf{Arrow:} {\scriptsize $A^{10}L_1 - 3A^{6}L_1 - A^{4}L_2 + 4A^{2}L_1 + 3L_2 + 1 - 3L_1/A^{2} - 3L_2/A^{4} - 1/A^{4} + L_1/A^{6} + 2L_2/A^{8} + A^{-8} - L_2/A^{12}$} \\ \textbf{Mock:} {\scriptsize $w^{5} + w^{4} - 4w^{3} - 4w^{2} + 5w + 6 - 2/w - 2/w^{2}$} \\ \textbf{Affine:} {\scriptsize $t^{2} - 2 + t^{-2}$} \\ \textbf{Yamada:} {\scriptsize $A^{27} - A^{26} - 3A^{25} + 4A^{24} + 2A^{23} - 9A^{22} + 4A^{21} + 7A^{20} - 10A^{19} + A^{18} + 6A^{17} - 5A^{16} - 2A^{15} + A^{14} + 4A^{13} - 6A^{12} - 3A^{11} + 9A^{10} - 6A^{9} - 6A^{8} + 9A^{7} - A^{6} - 6A^{5} + 3A^{4} + 2A^{3} - 2A^{2} + 1$}
\end{minipage}

\noindent{\color{gray!40}\rule{\textwidth}{0.4pt}}
\vspace{0.9\baselineskip}
\noindent \begin{minipage}[t]{0.25\textwidth}
\vspace{0pt}
\centering
\includegraphics[page=257,width=\linewidth]{knotoids.pdf}
\end{minipage}
\hfill
\begin{minipage}[t]{0.73\textwidth}
\vspace{0pt}
\raggedright
\textbf{Name:} {\large{$\mathbf{K7_{138}}$}} (chiral, non-rotatable$^{*}$) \\ \textbf{PD:} {\scriptsize\texttt{[0],[0,1,2,3],[1,4,5,2],[3,6,7,4],[7,8,9,5],[6,10,11,8],[12,13,14,9],[10,14,12,11],[13]}} \\ \textbf{EM:} {\scriptsize\texttt{(B0, A0C0C3D0, B1D3E3B2, B3F0E0C1, D2F3G3C2, D1H0H3E1, H2I0H1E2, F1G2G0F2, G1)}} \\ \textbf{Kauffman bracket:} {\scriptsize $-A^{21} + A^{17} + 2A^{15} - 2A^{11} - A^{9} + A^{7} + A^{5} - A^{3} - A$} \\ \textbf{Arrow:} {\scriptsize $A^{24} - A^{20} - 2A^{18}L_1 + 2A^{14}L_1 + A^{12} - A^{10}L_1 - A^{8} + A^{6}L_1 + A^{4}$} \\ \textbf{Mock:} {\scriptsize $-w^{2} - 3w + 3/w + 2/w^{2}$} \\ \textbf{Affine:} {\scriptsize $-3t + 6 - 3/t$} \\ \textbf{Yamada:} {\scriptsize $A^{21} - 2A^{19} + A^{17} - 2A^{16} - A^{15} + 3A^{14} - A^{13} + 2A^{11} + A^{9} - A^{8} + 2A^{7} - A^{6} - 2A^{5} + 3A^{4} + A^{3} - A^{2} + 2A + 1$}
\end{minipage}

\noindent{\color{gray!40}\rule{\textwidth}{0.4pt}}
\vspace{0.9\baselineskip}
\noindent \begin{minipage}[t]{0.25\textwidth}
\vspace{0pt}
\centering
\includegraphics[page=258,width=\linewidth]{knotoids.pdf}
\end{minipage}
\hfill
\begin{minipage}[t]{0.73\textwidth}
\vspace{0pt}
\raggedright
\textbf{Name:} {\large{$\mathbf{K7_{139}}$}} (chiral, non-rotatable$^{*}$) \\ \textbf{PD:} {\scriptsize\texttt{[0],[0,1,2,3],[1,4,5,2],[3,5,6,7],[4,8,9,6],[7,10,11,8],[12,13,10,9],[11,14,13,12],[14]}} \\ \textbf{EM:} {\scriptsize\texttt{(B0, A0C0C3D0, B1E0D1B2, B3C2E3F0, C1F3G3D2, D3G2H0E1, H3H2F1E2, F2I0G1G0, H1)}} \\ \textbf{Kauffman bracket:} {\scriptsize $A^{22} - 3A^{18} - A^{16} + 4A^{14} + 3A^{12} - 4A^{10} - 4A^{8} + 2A^{6} + 4A^{4} - 1$} \\ \textbf{Arrow:} {\scriptsize $A^{22}L_1 - 3A^{18}L_1 - A^{16}L_2 + 4A^{14}L_1 + 2A^{12}L_2 + A^{12} - 4A^{10}L_1 - 3A^{8}L_2 - A^{8} + 2A^{6}L_1 + 3A^{4}L_2 + A^{4} - L_2$} \\ \textbf{Mock:} {\scriptsize $-w^{5} - w^{4} + 4w^{3} + 4w^{2} - 6w - 5 + 3/w + 3/w^{2}$} \\ \textbf{Affine:} {\scriptsize $-t^{2} - t + 4 - 1/t - 1/t^{2}$} \\ \textbf{Yamada:} {\scriptsize $A^{27} - A^{26} - 2A^{25} + 6A^{24} - 9A^{22} + 7A^{21} + 4A^{20} - 12A^{19} + 4A^{18} + 6A^{17} - 6A^{16} - A^{15} + 2A^{14} + 3A^{13} - 7A^{12} + 10A^{10} - 9A^{9} - 4A^{8} + 10A^{7} - 5A^{6} - 6A^{5} + 4A^{4} + A^{3} - 3A^{2} + 1$}
\end{minipage}

\noindent{\color{gray!40}\rule{\textwidth}{0.4pt}}
\vspace{0.9\baselineskip}
\noindent \begin{minipage}[t]{0.25\textwidth}
\vspace{0pt}
\centering
\includegraphics[page=259,width=\linewidth]{knotoids.pdf}
\end{minipage}
\hfill
\begin{minipage}[t]{0.73\textwidth}
\vspace{0pt}
\raggedright
\textbf{Name:} {\large{$\mathbf{K7_{140}}$}} (chiral, non-rotatable$^{*}$) \\ \textbf{PD:} {\scriptsize\texttt{[0],[0,1,2,3],[1,4,5,2],[3,5,6,7],[4,8,9,6],[10,11,8,7],[9,12,13,10],[14,13,12,11],[14]}} \\ \textbf{EM:} {\scriptsize\texttt{(B0, A0C0C3D0, B1E0D1B2, B3C2E3F3, C1F2G0D2, G3H3E1D3, E2H2H1F0, I0G2G1F1, H0)}} \\ \textbf{Kauffman bracket:} {\scriptsize $A^{24} + A^{22} - 2A^{18} - A^{16} + 2A^{14} + 3A^{12} - A^{10} - 3A^{8} - A^{6} + A^{4} + A^{2}$} \\ \textbf{Arrow:} {\scriptsize $1 + L_1/A^{2} + L_2/A^{4} - 1/A^{4} - 2L_1/A^{6} - 2L_2/A^{8} + A^{-8} + 2L_1/A^{10} + 3L_2/A^{12} - L_1/A^{14} - 3L_2/A^{16} - L_1/A^{18} + L_2/A^{20} + L_1/A^{22}$} \\ \textbf{Mock:} {\scriptsize $-w^{6} - w^{5} + 2w^{4} + 3w^{3} + w^{2} - 2w - 4 - 1/w + 4/w^{2} + w^{-3} - 1/w^{4}$} \\ \textbf{Affine:} {\scriptsize $t^{2} - 2 + t^{-2}$} \\ \textbf{Yamada:} {\scriptsize $-A^{28} + A^{26} - 2A^{25} - 2A^{24} + 3A^{23} - 4A^{21} + A^{20} + 4A^{19} - 3A^{18} - 2A^{17} + 4A^{16} - 2A^{15} - 4A^{14} + 2A^{13} - 2A^{11} - 2A^{10} + 4A^{9} - A^{8} - 5A^{7} + 5A^{6} + A^{5} - 2A^{4} + A^{3} + A^{2} - 1$}
\end{minipage}

\noindent{\color{gray!40}\rule{\textwidth}{0.4pt}}
\vspace{0.9\baselineskip}
\noindent \begin{minipage}[t]{0.25\textwidth}
\vspace{0pt}
\centering
\includegraphics[page=260,width=\linewidth]{knotoids.pdf}
\end{minipage}
\hfill
\begin{minipage}[t]{0.73\textwidth}
\vspace{0pt}
\raggedright
\textbf{Name:} {\large{$\mathbf{K7_{141}}$}} (chiral, non-rotatable$^{*}$) \\ \textbf{PD:} {\scriptsize\texttt{[0],[0,1,2,3],[1,4,5,2],[3,5,6,7],[4,8,9,10],[10,11,12,6],[13,14,8,7],[14,13,12,9],[11]}} \\ \textbf{EM:} {\scriptsize\texttt{(B0, A0C0C3D0, B1E0D1B2, B3C2F3G3, C1G2H3F0, E3I0H2D2, H1H0E1D3, G1G0F2E2, F1)}} \\ \textbf{Kauffman bracket:} {\scriptsize $-A^{24} + 3A^{20} + 2A^{18} - 4A^{16} - 4A^{14} + 2A^{12} + 4A^{10} - 2A^{6} + A^{2}$} \\ \textbf{Arrow:} {\scriptsize $-A^{18}L_1 + 3A^{14}L_1 + 2A^{12} - 4A^{10}L_1 - 4A^{8} + 2A^{6}L_1 + 4A^{4} - 2 + A^{-4}$} \\ \textbf{Mock:} {\scriptsize $2w^{3} + 2w^{2} - 5w - 5 + 3/w + 5/w^{2} - 1/w^{4}$} \\ \textbf{Affine:} {\scriptsize $-t + 2 - 1/t$} \\ \textbf{Yamada:} {\scriptsize $A^{26} - 3A^{24} + A^{23} + 4A^{22} - 4A^{21} - 4A^{20} + 8A^{19} - 8A^{17} + 7A^{16} + 2A^{15} - 6A^{14} + A^{13} + 2A^{12} - 6A^{10} + 4A^{9} + 4A^{8} - 10A^{7} + A^{6} + 6A^{5} - 6A^{4} - 2A^{3} + 4A^{2} - 2$}
\end{minipage}

\noindent{\color{gray!40}\rule{\textwidth}{0.4pt}}
\vspace{0.9\baselineskip}
\noindent \begin{minipage}[t]{0.25\textwidth}
\vspace{0pt}
\centering
\includegraphics[page=261,width=\linewidth]{knotoids.pdf}
\end{minipage}
\hfill
\begin{minipage}[t]{0.73\textwidth}
\vspace{0pt}
\raggedright
\textbf{Name:} {\large{$\mathbf{K7_{142}}$}} (chiral, non-rotatable$^{*}$) \\ \textbf{PD:} {\scriptsize\texttt{[0],[0,1,2,3],[1,4,5,2],[3,6,7,8],[4,9,10,11],[11,12,6,5],[13,9,8,7],[13,12,14,10],[14]}} \\ \textbf{EM:} {\scriptsize\texttt{(B0, A0C0C3D0, B1E0F3B2, B3F2G3G2, C1G1H3F0, E3H1D1C2, H0E1D3D2, G0F1I0E2, H2)}} \\ \textbf{Kauffman bracket:} {\scriptsize $A^{23} + A^{21} - 2A^{19} - 3A^{17} + 2A^{15} + 4A^{13} - 4A^{9} - A^{7} + 2A^{5} - A$} \\ \textbf{Arrow:} {\scriptsize $-A^{26}L_1 - A^{24} + 2A^{22}L_1 + 3A^{20} - 2A^{18}L_1 - 4A^{16} + 4A^{12} + A^{10}L_1 - 2A^{8} + A^{4}$} \\ \textbf{Mock:} {\scriptsize $w^{3} + 3w^{2} - 3w - 7 + 1/w + 5/w^{2} + w^{-3}$} \\ \textbf{Affine:} {\scriptsize $-2t + 4 - 2/t$} \\ \textbf{Yamada:} {\scriptsize $2A^{24} - 4A^{22} + 2A^{21} + 3A^{20} - 6A^{19} - A^{18} + 5A^{17} - 3A^{16} - 3A^{15} + 6A^{14} - 3A^{12} + 2A^{11} + A^{10} - A^{9} - 5A^{8} + 6A^{7} + 2A^{6} - 6A^{5} + 6A^{4} + 3A^{3} - 3A^{2} + 2A + 1$}
\end{minipage}

\noindent{\color{gray!40}\rule{\textwidth}{0.4pt}}
\vspace{0.9\baselineskip}
\noindent \begin{minipage}[t]{0.25\textwidth}
\vspace{0pt}
\centering
\includegraphics[page=262,width=\linewidth]{knotoids.pdf}
\end{minipage}
\hfill
\begin{minipage}[t]{0.73\textwidth}
\vspace{0pt}
\raggedright
\textbf{Name:} {\large{$\mathbf{K7_{143}}$}} (chiral, non-rotatable$^{*}$) \\ \textbf{PD:} {\scriptsize\texttt{[0],[0,1,2,3],[1,4,5,2],[3,5,6,7],[8,9,10,4],[10,11,12,6],[7,13,14,8],[9,14,13,12],[11]}} \\ \textbf{EM:} {\scriptsize\texttt{(B0, A0C0C3D0, B1E3D1B2, B3C2F3G0, G3H0F0C1, E2I0H3D2, D3H2H1E0, E1G2G1F2, F1)}} \\ \textbf{Kauffman bracket:} {\scriptsize $A^{22} - 2A^{18} + 4A^{14} + A^{12} - 5A^{10} - 3A^{8} + 4A^{6} + 4A^{4} - A^{2} - 2$} \\ \textbf{Arrow:} {\scriptsize $A^{4} - 2 + 4/A^{4} + L_1/A^{6} - 5/A^{8} - 3L_1/A^{10} + 4/A^{12} + 4L_1/A^{14} - 1/A^{16} - 2L_1/A^{18}$} \\ \textbf{Mock:} {\scriptsize $-2w^{4} - 2w^{3} + 6w^{2} + 5w - 5 - 3/w + 3/w^{2} - 1/w^{4}$} \\ \textbf{Affine:} {\scriptsize $t - 2 + 1/t$} \\ \textbf{Yamada:} {\scriptsize $2A^{27} - A^{26} - 4A^{25} + 5A^{24} + 3A^{23} - 9A^{22} + 5A^{21} + 9A^{20} - 9A^{19} + 8A^{17} - 3A^{16} - 2A^{15} + 2A^{14} + 6A^{13} - 6A^{12} - 4A^{11} + 11A^{10} - 6A^{9} - 7A^{8} + 9A^{7} - 6A^{5} + 2A^{4} + 3A^{3} - 2A^{2} - A + 1$}
\end{minipage}

\noindent{\color{gray!40}\rule{\textwidth}{0.4pt}}
\vspace{0.9\baselineskip}
\noindent \begin{minipage}[t]{0.25\textwidth}
\vspace{0pt}
\centering
\includegraphics[page=263,width=\linewidth]{knotoids.pdf}
\end{minipage}
\hfill
\begin{minipage}[t]{0.73\textwidth}
\vspace{0pt}
\raggedright
\textbf{Name:} {\large{$\mathbf{K7_{144}}$}} (chiral, non-rotatable$^{*}$) \\ \textbf{PD:} {\scriptsize\texttt{[0],[0,1,2,3],[1,4,5,2],[3,6,7,8],[9,10,11,4],[11,12,6,5],[13,9,8,7],[10,13,12,14],[14]}} \\ \textbf{EM:} {\scriptsize\texttt{(B0, A0C0C3D0, B1E3F3B2, B3F2G3G2, G1H0F0C1, E2H2D1C2, H1E0D3D2, E1G0F1I0, H3)}} \\ \textbf{Kauffman bracket:} {\scriptsize $A^{22} - 2A^{18} - A^{16} + 3A^{14} + 2A^{12} - 3A^{10} - 4A^{8} + A^{6} + 3A^{4} + A^{2}$} \\ \textbf{Arrow:} {\scriptsize $A^{28} - 2A^{24} - A^{22}L_1 + 3A^{20} + 2A^{18}L_1 - 3A^{16} - 4A^{14}L_1 + A^{12} + 3A^{10}L_1 + A^{8}$} \\ \textbf{Mock:} {\scriptsize $-w^{4} - w^{3} + 3w^{2} + 2w - 4 - 4/w + 2/w^{2} + 3/w^{3} + w^{-4}$} \\ \textbf{Affine:} {\scriptsize $-3t + 6 - 3/t$} \\ \textbf{Yamada:} {\scriptsize $-A^{26} + A^{25} + 2A^{24} - 3A^{23} - A^{22} + 6A^{21} - 3A^{20} - 6A^{19} + 7A^{18} - A^{17} - 4A^{16} + 3A^{15} + A^{14} - A^{13} - 4A^{12} + 4A^{11} + 2A^{10} - 6A^{9} + 5A^{8} + 5A^{7} - 6A^{6} + A^{5} + 4A^{4} - A^{3} - A^{2} + 2A + 1$}
\end{minipage}

\noindent{\color{gray!40}\rule{\textwidth}{0.4pt}}
\vspace{0.9\baselineskip}
\noindent \begin{minipage}[t]{0.25\textwidth}
\vspace{0pt}
\centering
\includegraphics[page=264,width=\linewidth]{knotoids.pdf}
\end{minipage}
\hfill
\begin{minipage}[t]{0.73\textwidth}
\vspace{0pt}
\raggedright
\textbf{Name:} {\large{$\mathbf{K7_{145}}$}} (chiral, non-rotatable$^{*}$) \\ \textbf{PD:} {\scriptsize\texttt{[0],[0,1,2,3],[1,4,5,2],[6,7,8,3],[4,9,10,11],[5,11,12,6],[7,13,9,8],[13,12,14,10],[14]}} \\ \textbf{EM:} {\scriptsize\texttt{(B0, A0C0C3D3, B1E0F0B2, F3G0G3B3, C1G2H3F1, C2E3H1D0, D1H0E1D2, G1F2I0E2, H2)}} \\ \textbf{Kauffman bracket:} {\scriptsize $A^{18} + A^{16} - A^{14} - 2A^{12} - A^{10} + 2A^{8} + 2A^{6} - A^{2}$} \\ \textbf{Arrow:} {\scriptsize $A^{6}L_1 + A^{4} - A^{2}L_1 - 2 - L_1/A^{2} + 2/A^{4} + 2L_1/A^{6} - L_1/A^{10}$} \\ \textbf{Mock:} {\scriptsize $-w^{4} - w^{3} + 2w^{2} + 2w + w^{-2} - 1/w^{3} - 1/w^{4}$} \\ \textbf{Affine:} {\scriptsize $t - 2 + 1/t$} \\ \textbf{Yamada:} {\scriptsize $A^{25} - A^{23} - A^{20} - 2A^{19} + 2A^{17} - 2A^{16} - A^{15} + 2A^{14} - 2A^{13} - A^{12} + 2A^{11} - A^{10} - A^{8} + A^{7} - A^{6} - 3A^{5} + 2A^{4} - 2A^{2} + A + 1$}
\end{minipage}

\noindent{\color{gray!40}\rule{\textwidth}{0.4pt}}
\vspace{0.9\baselineskip}
\noindent \begin{minipage}[t]{0.25\textwidth}
\vspace{0pt}
\centering
\includegraphics[page=265,width=\linewidth]{knotoids.pdf}
\end{minipage}
\hfill
\begin{minipage}[t]{0.73\textwidth}
\vspace{0pt}
\raggedright
\textbf{Name:} {\large{$\mathbf{K7_{146}}$}} (chiral, non-rotatable$^{*}$) \\ \textbf{PD:} {\scriptsize\texttt{[0],[0,1,2,3],[1,4,5,2],[3,6,7,4],[5,7,8,9],[6,10,11,8],[11,12,13,9],[13,14,12,10],[14]}} \\ \textbf{EM:} {\scriptsize\texttt{(B0, A0C0C3D0, B1D3E0B2, B3F0E1C1, C2D2F3G3, D1H3G0E2, F2H2H0E3, G2I0G1F1, H1)}} \\ \textbf{Kauffman bracket:} {\scriptsize $-2A^{17} - A^{15} + 2A^{13} + 3A^{11} - A^{9} - 3A^{7} - A^{5} + A^{3} + A$} \\ \textbf{Arrow:} {\scriptsize $2/A^{4} + L_1/A^{6} - 2/A^{8} - 3L_1/A^{10} + A^{-12} + 3L_1/A^{14} + A^{-16} - L_1/A^{18} - 1/A^{20}$} \\ \textbf{Mock:} {\scriptsize $-w^{4} - 2w^{3} + 2w^{2} + 4w - 2/w$} \\ \textbf{Affine:} {\scriptsize $0$} \\ \textbf{Yamada:} {\scriptsize $A^{24} - 3A^{22} + 2A^{21} - 6A^{19} + A^{18} + A^{17} - 4A^{16} + A^{13} - 2A^{12} + A^{11} + 4A^{10} - 3A^{9} + 4A^{7} - 2A^{6} - 2A^{5} + 2A^{4} - 2A^{2} + 1$}
\end{minipage}

\noindent{\color{gray!40}\rule{\textwidth}{0.4pt}}
\vspace{0.9\baselineskip}
\noindent \begin{minipage}[t]{0.25\textwidth}
\vspace{0pt}
\centering
\includegraphics[page=266,width=\linewidth]{knotoids.pdf}
\end{minipage}
\hfill
\begin{minipage}[t]{0.73\textwidth}
\vspace{0pt}
\raggedright
\textbf{Name:} {\large{$\mathbf{K7_{147}}$}} (chiral, rotatable) \\ \textbf{PD:} {\scriptsize\texttt{[0],[0,1,2,3],[1,4,5,2],[3,6,7,4],[7,8,9,5],[6,10,11,8],[11,12,13,9],[10,13,14,12],[14]}} \\ \textbf{EM:} {\scriptsize\texttt{(B0, A0C0C3D0, B1D3E3B2, B3F0E0C1, D2F3G3C2, D1H0G0E1, F2H3H1E2, F1G2I0G1, H2)}} \\ \textbf{Kauffman bracket:} {\scriptsize $A^{18} - A^{14} - 2A^{12} + 2A^{8} + 2A^{6} - A^{2}$} \\ \textbf{Arrow:} {\scriptsize $A^{12} - A^{8} - 2A^{6}L_1 + 2A^{2}L_1 + 2 - 1/A^{4}$} \\ \textbf{Mock:} {\scriptsize $-w^{4} - 2w^{3} + 2w + 2$} \\ \textbf{Affine:} {\scriptsize $-2t + 4 - 2/t$} \\ \textbf{Yamada:} {\scriptsize $-A^{19} - A^{18} + 2A^{17} - A^{15} + 3A^{14} + A^{13} + 2A^{11} + A^{9} - A^{8} + A^{7} + A^{6} - 2A^{5} + A^{4} + A^{3} - 2A^{2} + 1$}
\end{minipage}

\noindent{\color{gray!40}\rule{\textwidth}{0.4pt}}
\vspace{0.9\baselineskip}
\noindent \begin{minipage}[t]{0.25\textwidth}
\vspace{0pt}
\centering
\includegraphics[page=267,width=\linewidth]{knotoids.pdf}
\end{minipage}
\hfill
\begin{minipage}[t]{0.73\textwidth}
\vspace{0pt}
\raggedright
\textbf{Name:} {\large{$\mathbf{K7_{148}}$}} (chiral, non-rotatable$^{*}$) \\ \textbf{PD:} {\scriptsize\texttt{[0],[0,1,2,3],[1,4,5,2],[3,5,6,7],[4,8,9,6],[10,11,8,7],[9,11,12,13],[10,13,14,12],[14]}} \\ \textbf{EM:} {\scriptsize\texttt{(B0, A0C0C3D0, B1E0D1B2, B3C2E3F3, C1F2G0D2, H0G1E1D3, E2F1H3H1, F0G3I0G2, H2)}} \\ \textbf{Kauffman bracket:} {\scriptsize $-A^{21} - A^{19} + 2A^{15} + A^{13} - 2A^{11} - 2A^{9} + A^{7} + 2A^{5} - A$} \\ \textbf{Arrow:} {\scriptsize $1 + L_1/A^{2} + L_2/A^{4} - 1/A^{4} - 2L_1/A^{6} - 2L_2/A^{8} + A^{-8} + 2L_1/A^{10} + 2L_2/A^{12} - L_1/A^{14} - 2L_2/A^{16} + L_2/A^{20}$} \\ \textbf{Mock:} {\scriptsize $-w^{6} - w^{5} + 2w^{4} + 3w^{3} + w^{2} - 2w - 3 + 3/w^{2} - 1/w^{4}$} \\ \textbf{Affine:} {\scriptsize $t - 2 + 1/t$} \\ \textbf{Yamada:} {\scriptsize $-A^{26} + A^{24} - A^{23} - A^{22} + 2A^{21} - 3A^{19} + A^{17} - 3A^{16} + 2A^{14} - 2A^{13} - 2A^{12} + A^{11} - A^{10} - A^{9} - A^{8} + 3A^{7} - 2A^{5} + 4A^{4} - A^{3} - A^{2} + A - 1$}
\end{minipage}

\noindent{\color{gray!40}\rule{\textwidth}{0.4pt}}
\vspace{0.9\baselineskip}
\noindent \begin{minipage}[t]{0.25\textwidth}
\vspace{0pt}
\centering
\includegraphics[page=268,width=\linewidth]{knotoids.pdf}
\end{minipage}
\hfill
\begin{minipage}[t]{0.73\textwidth}
\vspace{0pt}
\raggedright
\textbf{Name:} {\large{$\mathbf{K7_{149}}$}} (chiral, non-rotatable$^{*}$) \\ \textbf{PD:} {\scriptsize\texttt{[0],[0,1,2,3],[1,4,5,2],[3,5,6,7],[7,8,9,4],[10,11,8,6],[9,12,13,14],[14,13,11,10],[12]}} \\ \textbf{EM:} {\scriptsize\texttt{(B0, A0C0C3D0, B1E3D1B2, B3C2F3E0, D3F2G0C1, H3H2E1D2, E2I0H1H0, G3G2F1F0, G1)}} \\ \textbf{Kauffman bracket:} {\scriptsize $-A^{24} + 2A^{20} - 4A^{16} - A^{14} + 5A^{12} + 3A^{10} - 3A^{8} - 3A^{6} + A^{4} + 2A^{2}$} \\ \textbf{Arrow:} {\scriptsize $-A^{30}L_1 + 2A^{26}L_1 - 4A^{22}L_1 - A^{20}L_2 + 5A^{18}L_1 + 2A^{16}L_2 + A^{16} - 3A^{14}L_1 - 2A^{12}L_2 - A^{12} + A^{10}L_1 + A^{8}L_2 + A^{8}$} \\ \textbf{Mock:} {\scriptsize $w^{5} + w^{4} - 4w^{3} - 2w^{2} + 5w + 2 - 4/w - 2/w^{2} + 2/w^{3} + 2/w^{4}$} \\ \textbf{Affine:} {\scriptsize $-t^{2} - 2t + 6 - 2/t - 1/t^{2}$} \\ \textbf{Yamada:} {\scriptsize $-2A^{27} + 2A^{26} + 3A^{25} - 6A^{24} + 2A^{23} + 8A^{22} - 9A^{21} + 8A^{19} - 4A^{18} - 2A^{17} + A^{16} + 4A^{15} - 5A^{14} - 3A^{13} + 8A^{12} - 4A^{11} - 7A^{10} + 7A^{9} - 8A^{7} + 2A^{6} + 3A^{5} - 4A^{4} - A^{3} + 2A^{2} - 1$}
\end{minipage}

\noindent{\color{gray!40}\rule{\textwidth}{0.4pt}}
\vspace{0.9\baselineskip}
\noindent \begin{minipage}[t]{0.25\textwidth}
\vspace{0pt}
\centering
\includegraphics[page=269,width=\linewidth]{knotoids.pdf}
\end{minipage}
\hfill
\begin{minipage}[t]{0.73\textwidth}
\vspace{0pt}
\raggedright
\textbf{Name:} {\large{$\mathbf{K7_{150}}$}} (chiral, non-rotatable$^{*}$) \\ \textbf{PD:} {\scriptsize\texttt{[0],[0,1,2,3],[1,4,5,2],[3,6,7,4],[5,8,9,10],[6,10,11,7],[8,12,13,9],[14,13,12,11],[14]}} \\ \textbf{EM:} {\scriptsize\texttt{(B0, A0C0C3D0, B1D3E0B2, B3F0F3C1, C2G0G3F1, D1E3H3D2, E1H2H1E2, I0G2G1F2, H0)}} \\ \textbf{Kauffman bracket:} {\scriptsize $A^{22} - A^{20} - 3A^{18} + 5A^{14} + 3A^{12} - 4A^{10} - 4A^{8} + 2A^{6} + 3A^{4} - 1$} \\ \textbf{Arrow:} {\scriptsize $A^{16} - A^{14}L_1 - 3A^{12} + 5A^{8} + 3A^{6}L_1 - 4A^{4} - 4A^{2}L_1 + 2 + 3L_1/A^{2} - L_1/A^{6}$} \\ \textbf{Mock:} {\scriptsize $w^{4} - 5w^{2} - 3w + 7 + 5/w - 2/w^{2} - 2/w^{3}$} \\ \textbf{Affine:} {\scriptsize $-t + 2 - 1/t$} \\ \textbf{Yamada:} {\scriptsize $-A^{25} + A^{24} + 3A^{23} - 4A^{22} - 3A^{21} + 9A^{20} - 2A^{19} - 10A^{18} + 10A^{17} + 3A^{16} - 11A^{15} + 5A^{14} + 4A^{13} - 6A^{12} - 2A^{11} + 2A^{10} + 3A^{9} - 10A^{8} + A^{7} + 11A^{6} - 11A^{5} - 2A^{4} + 10A^{3} - 6A^{2} - 3A + 3$}
\end{minipage}

\noindent{\color{gray!40}\rule{\textwidth}{0.4pt}}
\vspace{0.9\baselineskip}
\noindent \begin{minipage}[t]{0.25\textwidth}
\vspace{0pt}
\centering
\includegraphics[page=270,width=\linewidth]{knotoids.pdf}
\end{minipage}
\hfill
\begin{minipage}[t]{0.73\textwidth}
\vspace{0pt}
\raggedright
\textbf{Name:} {\large{$\mathbf{K7_{151}}$}} (chiral, non-rotatable$^{*}$) \\ \textbf{PD:} {\scriptsize\texttt{[0],[0,1,2,3],[1,4,5,2],[3,6,7,4],[8,9,10,5],[6,10,11,7],[12,13,9,8],[11,14,13,12],[14]}} \\ \textbf{EM:} {\scriptsize\texttt{(B0, A0C0C3D0, B1D3E3B2, B3F0F3C1, G3G2F1C2, D1E2H0D2, H3H2E1E0, F2I0G1G0, H1)}} \\ \textbf{Kauffman bracket:} {\scriptsize $-A^{25} + 2A^{21} + A^{19} - 4A^{17} - 3A^{15} + 3A^{13} + 5A^{11} - A^{9} - 4A^{7} - A^{5} + A^{3} + A$} \\ \textbf{Arrow:} {\scriptsize $A^{4} - 2 - L_1/A^{2} + 4/A^{4} + 3L_1/A^{6} - 3/A^{8} - 5L_1/A^{10} + A^{-12} + 4L_1/A^{14} + A^{-16} - L_1/A^{18} - 1/A^{20}$} \\ \textbf{Mock:} {\scriptsize $-2w^{4} - 3w^{3} + 4w^{2} + 6w - 2 - 4/w + w^{-2} + w^{-3}$} \\ \textbf{Affine:} {\scriptsize $-t + 2 - 1/t$} \\ \textbf{Yamada:} {\scriptsize $A^{27} - A^{26} - 3A^{25} + 5A^{24} + 2A^{23} - 11A^{22} + 5A^{21} + 7A^{20} - 13A^{19} + 7A^{17} - 6A^{16} - 3A^{15} + 2A^{14} + 5A^{13} - 6A^{12} - 2A^{11} + 12A^{10} - 7A^{9} - 7A^{8} + 11A^{7} - A^{6} - 7A^{5} + 4A^{4} + 3A^{3} - 3A^{2} - A + 1$}
\end{minipage}

\noindent{\color{gray!40}\rule{\textwidth}{0.4pt}}
\vspace{0.9\baselineskip}
\noindent \begin{minipage}[t]{0.25\textwidth}
\vspace{0pt}
\centering
\includegraphics[page=271,width=\linewidth]{knotoids.pdf}
\end{minipage}
\hfill
\begin{minipage}[t]{0.73\textwidth}
\vspace{0pt}
\raggedright
\textbf{Name:} {\large{$\mathbf{K7_{152}}$}} (chiral, non-rotatable$^{*}$) \\ \textbf{PD:} {\scriptsize\texttt{[0],[0,1,2,3],[1,4,5,2],[3,6,7,8],[4,8,9,5],[6,10,11,7],[12,13,10,9],[11,14,13,12],[14]}} \\ \textbf{EM:} {\scriptsize\texttt{(B0, A0C0C3D0, B1E0E3B2, B3F0F3E1, C1D3G3C2, D1G2H0D2, H3H2F1E2, F2I0G1G0, H1)}} \\ \textbf{Kauffman bracket:} {\scriptsize $-A^{26} + 2A^{22} + A^{20} - 3A^{18} - 4A^{16} + 2A^{14} + 5A^{12} + A^{10} - 3A^{8} - A^{6} + A^{4} + A^{2}$} \\ \textbf{Arrow:} {\scriptsize $-A^{32} + 2A^{28} + A^{26}L_1 - 3A^{24} - 4A^{22}L_1 + 2A^{20} + 5A^{18}L_1 + A^{16} - 3A^{14}L_1 - A^{12} + A^{10}L_1 + A^{8}$} \\ \textbf{Mock:} {\scriptsize $-w^{3} + 3w^{2} + 4w - 4 - 6/w + 3/w^{3} + 2/w^{4}$} \\ \textbf{Affine:} {\scriptsize $-t + 2 - 1/t$} \\ \textbf{Yamada:} {\scriptsize $A^{26} - 4A^{24} + 6A^{22} - 5A^{21} - 5A^{20} + 10A^{19} + A^{18} - 9A^{17} + 8A^{16} + 4A^{15} - 6A^{14} + 2A^{13} + 3A^{12} + A^{11} - 8A^{10} + 4A^{9} + 5A^{8} - 13A^{7} + 7A^{5} - 7A^{4} - 3A^{3} + 4A^{2} - 2$}
\end{minipage}

\noindent{\color{gray!40}\rule{\textwidth}{0.4pt}}
\vspace{0.9\baselineskip}
\noindent \begin{minipage}[t]{0.25\textwidth}
\vspace{0pt}
\centering
\includegraphics[page=272,width=\linewidth]{knotoids.pdf}
\end{minipage}
\hfill
\begin{minipage}[t]{0.73\textwidth}
\vspace{0pt}
\raggedright
\textbf{Name:} {\large{$\mathbf{K7_{153}}$}} (chiral, non-rotatable$^{*}$) \\ \textbf{PD:} {\scriptsize\texttt{[0],[0,1,2,3],[1,4,5,2],[6,7,4,3],[5,8,9,10],[10,11,7,6],[8,12,13,9],[14,13,12,11],[14]}} \\ \textbf{EM:} {\scriptsize\texttt{(B0, A0C0C3D3, B1D2E0B2, F3F2C1B3, C2G0G3F0, E3H3D1D0, E1H2H1E2, I0G2G1F1, H0)}} \\ \textbf{Kauffman bracket:} {\scriptsize $A^{24} + 2A^{22} - 3A^{18} - 3A^{16} + 2A^{14} + 4A^{12} - 3A^{8} - A^{6} + A^{4} + A^{2}$} \\ \textbf{Arrow:} {\scriptsize $A^{6}L_1 + 2A^{4} - 3 - 3L_1/A^{2} + 2/A^{4} + 4L_1/A^{6} - 3L_1/A^{10} - 1/A^{12} + L_1/A^{14} + A^{-16}$} \\ \textbf{Mock:} {\scriptsize $-w^{4} + 3w^{2} + 3w - 1 - 5/w + 2/w^{3}$} \\ \textbf{Affine:} {\scriptsize $t - 2 + 1/t$} \\ \textbf{Yamada:} {\scriptsize $-A^{28} + 2A^{26} - A^{25} - 4A^{24} + 3A^{23} + 3A^{22} - 6A^{21} - A^{20} + 8A^{19} - 3A^{18} - 6A^{17} + 6A^{16} - A^{15} - 6A^{14} + 2A^{13} + 2A^{12} - 2A^{11} - 5A^{10} + 6A^{9} + 2A^{8} - 9A^{7} + 5A^{6} + 4A^{5} - 5A^{4} + 2A^{2} - 1$}
\end{minipage}

\noindent{\color{gray!40}\rule{\textwidth}{0.4pt}}
\vspace{0.9\baselineskip}
\noindent \begin{minipage}[t]{0.25\textwidth}
\vspace{0pt}
\centering
\includegraphics[page=273,width=\linewidth]{knotoids.pdf}
\end{minipage}
\hfill
\begin{minipage}[t]{0.73\textwidth}
\vspace{0pt}
\raggedright
\textbf{Name:} {\large{$\mathbf{K7_{154}}$}} (chiral, non-rotatable$^{*}$) \\ \textbf{PD:} {\scriptsize\texttt{[0],[0,1,2,3],[1,4,5,2],[3,6,7,4],[5,8,9,6],[7,9,10,11],[12,13,14,8],[13,12,11,10],[14]}} \\ \textbf{EM:} {\scriptsize\texttt{(B0, A0C0C3D0, B1D3E0B2, B3E3F0C1, C2G3F1D1, D2E2H3H2, H1H0I0E1, G1G0F3F2, G2)}} \\ \textbf{Kauffman bracket:} {\scriptsize $A^{20} - A^{16} + 2A^{12} + A^{10} - 2A^{8} - 2A^{6} + A^{2} + 1$} \\ \textbf{Arrow:} {\scriptsize $A^{8} - A^{4} + 2 + L_1/A^{2} - 2/A^{4} - 2L_1/A^{6} - L_2/A^{8} + A^{-8} + L_1/A^{10} + L_2/A^{12}$} \\ \textbf{Mock:} {\scriptsize $w^{3} - w^{2} - 2w + 3 + 1/w - 1/w^{2}$} \\ \textbf{Affine:} {\scriptsize $-t^{2} + 2 - 1/t^{2}$} \\ \textbf{Yamada:} {\scriptsize $A^{25} - A^{24} + A^{22} - 3A^{21} + A^{19} - 3A^{18} - A^{17} - A^{16} + A^{15} - A^{14} + 3A^{12} - A^{11} - A^{10} + A^{9} - A^{8} - 2A^{7} + A^{6} + A^{5} - A^{4} + A^{2} - 1$}
\end{minipage}

\noindent{\color{gray!40}\rule{\textwidth}{0.4pt}}
\vspace{0.9\baselineskip}
\noindent \begin{minipage}[t]{0.25\textwidth}
\vspace{0pt}
\centering
\includegraphics[page=274,width=\linewidth]{knotoids.pdf}
\end{minipage}
\hfill
\begin{minipage}[t]{0.73\textwidth}
\vspace{0pt}
\raggedright
\textbf{Name:} {\large{$\mathbf{K7_{155}}$}} (chiral, non-rotatable$^{*}$) \\ \textbf{PD:} {\scriptsize\texttt{[0],[0,1,2,3],[1,4,5,2],[3,6,7,4],[5,8,9,10],[6,11,12,13],[13,12,8,7],[14,11,10,9],[14]}} \\ \textbf{EM:} {\scriptsize\texttt{(B0, A0C0C3D0, B1D3E0B2, B3F0G3C1, C2G2H3H2, D1H1G1G0, F3F2E1D2, I0F1E3E2, H0)}} \\ \textbf{Kauffman bracket:} {\scriptsize $A^{21} + 2A^{19} - 2A^{17} - 4A^{15} + A^{13} + 4A^{11} + A^{9} - 4A^{7} - 2A^{5} + A^{3} + A$} \\ \textbf{Arrow:} {\scriptsize $-A^{18}L_1 - 2A^{16} + 2A^{14}L_1 + 4A^{12} - A^{10}L_1 - 4A^{8} - A^{6}L_1 + 4A^{4} + 2A^{2}L_1 - 1 - L_1/A^{2}$} \\ \textbf{Mock:} {\scriptsize $-w^{5} - w^{4} + 2w^{3} + 4w^{2} - w - 5 - 1/w + 4/w^{2} + 2/w^{3} - 1/w^{4} - 1/w^{5}$} \\ \textbf{Affine:} {\scriptsize $0$} \\ \textbf{Yamada:} {\scriptsize $A^{26} - 2A^{25} + 5A^{23} - 6A^{22} - 3A^{21} + 9A^{20} - 6A^{19} - 2A^{18} + 9A^{17} - A^{15} - A^{14} + 3A^{13} - 4A^{12} - 8A^{11} + 4A^{10} - A^{9} - 9A^{8} + 5A^{7} + 4A^{6} - 6A^{5} + A^{4} + 4A^{3} - 2A^{2} - A + 1$}
\end{minipage}

\noindent{\color{gray!40}\rule{\textwidth}{0.4pt}}
\vspace{0.9\baselineskip}
\noindent \begin{minipage}[t]{0.25\textwidth}
\vspace{0pt}
\centering
\includegraphics[page=275,width=\linewidth]{knotoids.pdf}
\end{minipage}
\hfill
\begin{minipage}[t]{0.73\textwidth}
\vspace{0pt}
\raggedright
\textbf{Name:} {\large{$\mathbf{K7_{156}}$}} (chiral, non-rotatable$^{*}$) \\ \textbf{PD:} {\scriptsize\texttt{[0],[0,1,2,3],[1,4,5,2],[3,6,7,4],[5,8,9,10],[11,12,13,6],[7,13,14,8],[12,11,10,9],[14]}} \\ \textbf{EM:} {\scriptsize\texttt{(B0, A0C0C3D0, B1D3E0B2, B3F3G0C1, C2G3H3H2, H1H0G1D1, D2F2I0E1, F1F0E3E2, G2)}} \\ \textbf{Kauffman bracket:} {\scriptsize $A^{12} - A^{8} - A^{6} + A^{2} + 1$} \\ \textbf{Arrow:} {\scriptsize $A^{12} - A^{8} - A^{6}L_1 - A^{4}L_2 + A^{4} + A^{2}L_1 + L_2$} \\ \textbf{Mock:} {\scriptsize $w^{5} + w^{4} - w^{3} - w^{2} - w - 1 + 2/w^{2} + w^{-3}$} \\ \textbf{Affine:} {\scriptsize $-t^{2} - t + 4 - 1/t - 1/t^{2}$} \\ \textbf{Yamada:} {\scriptsize $A^{23} - A^{21} + A^{20} - A^{16} - 2A^{14} - A^{13} - A^{11} + A^{9} - A^{8} - A^{3} - 1$}
\end{minipage}

\noindent{\color{gray!40}\rule{\textwidth}{0.4pt}}
\vspace{0.9\baselineskip}
\noindent \begin{minipage}[t]{0.25\textwidth}
\vspace{0pt}
\centering
\includegraphics[page=276,width=\linewidth]{knotoids.pdf}
\end{minipage}
\hfill
\begin{minipage}[t]{0.73\textwidth}
\vspace{0pt}
\raggedright
\textbf{Name:} {\large{$\mathbf{K7_{157}}$}} (chiral, non-rotatable$^{*}$) \\ \textbf{PD:} {\scriptsize\texttt{[0],[0,1,2,3],[1,4,5,2],[3,6,7,4],[5,8,9,10],[11,12,13,6],[13,14,8,7],[12,11,10,9],[14]}} \\ \textbf{EM:} {\scriptsize\texttt{(B0, A0C0C3D0, B1D3E0B2, B3F3G3C1, C2G2H3H2, H1H0G0D1, F2I0E1D2, F1F0E3E2, G1)}} \\ \textbf{Kauffman bracket:} {\scriptsize $-A^{18} - A^{16} + 2A^{14} + 3A^{12} - A^{10} - 4A^{8} - A^{6} + 2A^{4} + 2A^{2}$} \\ \textbf{Arrow:} {\scriptsize $-A^{24} - A^{22}L_1 + 2A^{20} + 3A^{18}L_1 - A^{16} - 4A^{14}L_1 - A^{12}L_2 + 2A^{10}L_1 + A^{8}L_2 + A^{8}$} \\ \textbf{Mock:} {\scriptsize $w^{5} + w^{4} - 2w^{3} + 2w - 2 - 3/w + w^{-2} + 2/w^{3} + w^{-4}$} \\ \textbf{Affine:} {\scriptsize $-t^{2} - t + 4 - 1/t - 1/t^{2}$} \\ \textbf{Yamada:} {\scriptsize $-2A^{25} + A^{24} + 4A^{23} - 4A^{22} - A^{21} + 7A^{20} - 3A^{19} - 2A^{18} + 3A^{17} - A^{16} - A^{15} - 2A^{14} + 3A^{13} - A^{12} - 5A^{11} + 4A^{10} - 6A^{8} + 3A^{7} + 2A^{6} - 3A^{5} - A^{4} + A^{3} - A - 1$}
\end{minipage}

\noindent{\color{gray!40}\rule{\textwidth}{0.4pt}}
\vspace{0.9\baselineskip}
\noindent \begin{minipage}[t]{0.25\textwidth}
\vspace{0pt}
\centering
\includegraphics[page=277,width=\linewidth]{knotoids.pdf}
\end{minipage}
\hfill
\begin{minipage}[t]{0.73\textwidth}
\vspace{0pt}
\raggedright
\textbf{Name:} {\large{$\mathbf{K7_{158}}$}} (chiral, non-rotatable$^{*}$) \\ \textbf{PD:} {\scriptsize\texttt{[0],[0,1,2,3],[1,4,5,2],[3,6,7,4],[8,9,6,5],[7,9,10,11],[12,13,14,8],[13,12,11,10],[14]}} \\ \textbf{EM:} {\scriptsize\texttt{(B0, A0C0C3D0, B1D3E3B2, B3E2F0C1, G3F1D1C2, D2E1H3H2, H1H0I0E0, G1G0F3F2, G2)}} \\ \textbf{Kauffman bracket:} {\scriptsize $A^{22} - 2A^{18} - A^{16} + 4A^{14} + 4A^{12} - 3A^{10} - 5A^{8} + 3A^{4} + A^{2} - 1$} \\ \textbf{Arrow:} {\scriptsize $A^{4} - 2 - L_1/A^{2} + 4/A^{4} + 4L_1/A^{6} - 3/A^{8} - 5L_1/A^{10} - L_2/A^{12} + A^{-12} + 3L_1/A^{14} + L_2/A^{16} - L_1/A^{18}$} \\ \textbf{Mock:} {\scriptsize $-2w^{4} - 2w^{3} + 5w^{2} + 6w - 2 - 5/w + w^{-3}$} \\ \textbf{Affine:} {\scriptsize $-t^{2} + t + 1/t - 1/t^{2}$} \\ \textbf{Yamada:} {\scriptsize $2A^{25} - A^{24} - 5A^{23} + 5A^{22} + 4A^{21} - 11A^{20} + 2A^{19} + 7A^{18} - 11A^{17} - 2A^{16} + 6A^{15} - 4A^{14} - 2A^{13} + A^{12} + 7A^{11} - 4A^{10} - 4A^{9} + 12A^{8} - 5A^{7} - 8A^{6} + 8A^{5} - 6A^{3} + 2A^{2} + 2A - 1$}
\end{minipage}

\noindent{\color{gray!40}\rule{\textwidth}{0.4pt}}
\vspace{0.9\baselineskip}
\noindent \begin{minipage}[t]{0.25\textwidth}
\vspace{0pt}
\centering
\includegraphics[page=278,width=\linewidth]{knotoids.pdf}
\end{minipage}
\hfill
\begin{minipage}[t]{0.73\textwidth}
\vspace{0pt}
\raggedright
\textbf{Name:} {\large{$\mathbf{K7_{159}}$}} (chiral, non-rotatable$^{*}$) \\ \textbf{PD:} {\scriptsize\texttt{[0],[0,1,2,3],[1,4,5,2],[3,6,7,4],[8,9,6,5],[9,10,11,7],[8,12,13,14],[10,13,12,11],[14]}} \\ \textbf{EM:} {\scriptsize\texttt{(B0, A0C0C3D0, B1D3E3B2, B3E2F3C1, G0F0D1C2, E1H0H3D2, E0H2H1I0, F1G2G1F2, G3)}} \\ \textbf{Kauffman bracket:} {\scriptsize $A^{24} + A^{22} - A^{20} - 4A^{18} - A^{16} + 5A^{14} + 3A^{12} - 3A^{10} - 4A^{8} + 2A^{6} + 3A^{4} - 1$} \\ \textbf{Arrow:} {\scriptsize $A^{18}L_1 + A^{16} - A^{14}L_1 - A^{12}L_2 - 3A^{12} - A^{10}L_1 + A^{8}L_2 + 4A^{8} + 3A^{6}L_1 - 3A^{4} - 4A^{2}L_1 + 2 + 3L_1/A^{2} - L_1/A^{6}$} \\ \textbf{Mock:} {\scriptsize $w^{4} + w^{3} - 4w^{2} - 3w + 7 + 4/w - 3/w^{2} - 2/w^{3}$} \\ \textbf{Affine:} {\scriptsize $-t^{2} + t + 1/t - 1/t^{2}$} \\ \textbf{Yamada:} {\scriptsize $A^{27} - 3A^{25} + A^{24} + 5A^{23} - 5A^{22} - 5A^{21} + 10A^{20} - 11A^{18} + 9A^{17} + 6A^{16} - 12A^{15} + 2A^{14} + 5A^{13} - 6A^{12} - 4A^{11} + 3A^{10} + 6A^{9} - 10A^{8} + 13A^{6} - 10A^{5} - 5A^{4} + 9A^{3} - 4A^{2} - 3A + 2$}
\end{minipage}

\noindent{\color{gray!40}\rule{\textwidth}{0.4pt}}
\vspace{0.9\baselineskip}
\noindent \begin{minipage}[t]{0.25\textwidth}
\vspace{0pt}
\centering
\includegraphics[page=279,width=\linewidth]{knotoids.pdf}
\end{minipage}
\hfill
\begin{minipage}[t]{0.73\textwidth}
\vspace{0pt}
\raggedright
\textbf{Name:} {\large{$\mathbf{K7_{160}}$}} (chiral, non-rotatable$^{*}$) \\ \textbf{PD:} {\scriptsize\texttt{[0],[0,1,2,3],[1,4,5,2],[3,6,7,4],[8,9,10,5],[6,11,8,7],[9,12,13,14],[14,13,11,10],[12]}} \\ \textbf{EM:} {\scriptsize\texttt{(B0, A0C0C3D0, B1D3E3B2, B3F0F3C1, F2G0H3C2, D1H2E0D2, E1I0H1H0, G3G2F1E2, G1)}} \\ \textbf{Kauffman bracket:} {\scriptsize $-A^{26} + 2A^{22} + 2A^{20} - 2A^{18} - 5A^{16} + 5A^{12} + 3A^{10} - 2A^{8} - 2A^{6} + A^{2}$} \\ \textbf{Arrow:} {\scriptsize $-A^{20} + 2A^{16} + 2A^{14}L_1 - 2A^{12} - 5A^{10}L_1 + 5A^{6}L_1 + 3A^{4} - 2A^{2}L_1 - 2 + A^{-4}$} \\ \textbf{Mock:} {\scriptsize $w^{4} + 2w^{3} - w^{2} - 5w + 5/w + 2/w^{2} - 2/w^{3} - 1/w^{4}$} \\ \textbf{Affine:} {\scriptsize $t - 2 + 1/t$} \\ \textbf{Yamada:} {\scriptsize $A^{25} - A^{24} - 4A^{23} + 6A^{22} + 3A^{21} - 12A^{20} + 6A^{19} + 9A^{18} - 12A^{17} + 2A^{16} + 7A^{15} - 5A^{14} - 2A^{13} + A^{12} + 5A^{11} - 7A^{10} - 4A^{9} + 11A^{8} - 8A^{7} - 9A^{6} + 11A^{5} - 2A^{4} - 7A^{3} + 4A^{2} + 2A - 1$}
\end{minipage}

\noindent{\color{gray!40}\rule{\textwidth}{0.4pt}}
\vspace{0.9\baselineskip}
\noindent \begin{minipage}[t]{0.25\textwidth}
\vspace{0pt}
\centering
\includegraphics[page=280,width=\linewidth]{knotoids.pdf}
\end{minipage}
\hfill
\begin{minipage}[t]{0.73\textwidth}
\vspace{0pt}
\raggedright
\textbf{Name:} {\large{$\mathbf{K7_{161}}$}} (chiral, non-rotatable$^{*}$) \\ \textbf{PD:} {\scriptsize\texttt{[0],[0,1,2,3],[1,4,5,2],[3,6,7,4],[8,9,10,5],[6,11,12,13],[7,13,14,8],[9,12,11,10],[14]}} \\ \textbf{EM:} {\scriptsize\texttt{(B0, A0C0C3D0, B1D3E3B2, B3F0G0C1, G3H0H3C2, D1H2H1G1, D2F3I0E0, E1F2F1E2, G2)}} \\ \textbf{Kauffman bracket:} {\scriptsize $A^{24} + A^{22} - A^{20} - 2A^{18} + 2A^{14} + A^{12} - A^{10} + A^{6} - A^{2}$} \\ \textbf{Arrow:} {\scriptsize $1 + L_1/A^{2} - 1/A^{4} - 2L_1/A^{6} - L_2/A^{8} + A^{-8} + 2L_1/A^{10} + L_2/A^{12} - L_1/A^{14} + L_1/A^{18} - L_1/A^{22}$} \\ \textbf{Mock:} {\scriptsize $-w^{6} - w^{5} + 2w^{4} + 2w^{3} - w + 2/w + w^{-2} - 2/w^{3} - 1/w^{4}$} \\ \textbf{Affine:} {\scriptsize $-t^{2} + 2t - 2 + 2/t - 1/t^{2}$} \\ \textbf{Yamada:} {\scriptsize $-A^{22} + A^{20} - A^{19} + 2A^{17} + A^{16} + A^{15} + A^{14} + 2A^{13} + 2A^{10} - A^{8} + A^{7} - A^{6} - A^{5} - A^{4} + A^{3} - A + 1$}
\end{minipage}

\noindent{\color{gray!40}\rule{\textwidth}{0.4pt}}
\vspace{0.9\baselineskip}
\noindent \begin{minipage}[t]{0.25\textwidth}
\vspace{0pt}
\centering
\includegraphics[page=281,width=\linewidth]{knotoids.pdf}
\end{minipage}
\hfill
\begin{minipage}[t]{0.73\textwidth}
\vspace{0pt}
\raggedright
\textbf{Name:} {\large{$\mathbf{K7_{162}}$}} (chiral, rotatable) \\ \textbf{PD:} {\scriptsize\texttt{[0],[0,1,2,3],[1,4,5,2],[3,6,7,4],[8,9,10,5],[6,11,12,13],[13,12,8,7],[9,14,11,10],[14]}} \\ \textbf{EM:} {\scriptsize\texttt{(B0, A0C0C3D0, B1D3E3B2, B3F0G3C1, G2H0H3C2, D1H2G1G0, F3F2E0D2, E1I0F1E2, H1)}} \\ \textbf{Kauffman bracket:} {\scriptsize $-A^{26} + 2A^{22} + 2A^{20} - A^{18} - 4A^{16} + 4A^{12} + 2A^{10} - 2A^{8} - 2A^{6} + A^{2}$} \\ \textbf{Arrow:} {\scriptsize $-A^{8} + 2A^{4} + 2A^{2}L_1 - 1 - 4L_1/A^{2} + 4L_1/A^{6} + 2/A^{8} - 2L_1/A^{10} - 2/A^{12} + A^{-16}$} \\ \textbf{Mock:} {\scriptsize $w^{6} + 2w^{5} - w^{4} - 4w^{3} + 4w + 2 - 2/w - 1/w^{2}$} \\ \textbf{Affine:} {\scriptsize $2t - 4 + 2/t$} \\ \textbf{Yamada:} {\scriptsize $-A^{24} - 2A^{23} + 4A^{22} + 3A^{21} - 8A^{20} + 4A^{19} + 6A^{18} - 11A^{17} - A^{16} + 4A^{15} - 5A^{14} - A^{13} + A^{12} + 4A^{11} - 4A^{10} - 2A^{9} + 9A^{8} - 5A^{7} - 6A^{6} + 8A^{5} - A^{4} - 6A^{3} + 3A^{2} + 2A - 1$}
\end{minipage}

\noindent{\color{gray!40}\rule{\textwidth}{0.4pt}}
\vspace{0.9\baselineskip}
\noindent \begin{minipage}[t]{0.25\textwidth}
\vspace{0pt}
\centering
\includegraphics[page=282,width=\linewidth]{knotoids.pdf}
\end{minipage}
\hfill
\begin{minipage}[t]{0.73\textwidth}
\vspace{0pt}
\raggedright
\textbf{Name:} {\large{$\mathbf{K7_{163}}$}} (chiral, non-rotatable$^{*}$) \\ \textbf{PD:} {\scriptsize\texttt{[0],[0,1,2,3],[1,4,5,2],[3,6,7,4],[8,9,10,5],[6,11,12,13],[13,14,8,7],[9,12,11,10],[14]}} \\ \textbf{EM:} {\scriptsize\texttt{(B0, A0C0C3D0, B1D3E3B2, B3F0G3C1, G2H0H3C2, D1H2H1G0, F3I0E0D2, E1F2F1E2, G1)}} \\ \textbf{Kauffman bracket:} {\scriptsize $A^{24} + A^{22} - A^{20} - 3A^{18} - A^{16} + 4A^{14} + 3A^{12} - 2A^{10} - 3A^{8} + A^{6} + 2A^{4} - 1$} \\ \textbf{Arrow:} {\scriptsize $A^{6}L_1 + A^{4} - A^{2}L_1 - L_2 - 2 - L_1/A^{2} + L_2/A^{4} + 3/A^{4} + 3L_1/A^{6} - 2/A^{8} - 3L_1/A^{10} + A^{-12} + 2L_1/A^{14} - L_1/A^{18}$} \\ \textbf{Mock:} {\scriptsize $-w^{4} - w^{3} + 3w^{2} + 3w - 2 - 2/w + 3/w^{2} + w^{-3} - 2/w^{4} - 1/w^{5}$} \\ \textbf{Affine:} {\scriptsize $-t^{2} + 2t - 2 + 2/t - 1/t^{2}$} \\ \textbf{Yamada:} {\scriptsize $A^{27} - 2A^{25} + A^{24} + 3A^{23} - 3A^{22} - 4A^{21} + 5A^{20} + A^{19} - 9A^{18} + 3A^{17} + 4A^{16} - 8A^{15} + A^{14} + 5A^{13} - 3A^{12} - 2A^{11} + A^{10} + 4A^{9} - 6A^{8} - 2A^{7} + 9A^{6} - 5A^{5} - 4A^{4} + 6A^{3} - A^{2} - 2A + 1$}
\end{minipage}

\noindent{\color{gray!40}\rule{\textwidth}{0.4pt}}
\vspace{0.9\baselineskip}
\noindent \begin{minipage}[t]{0.25\textwidth}
\vspace{0pt}
\centering
\includegraphics[page=283,width=\linewidth]{knotoids.pdf}
\end{minipage}
\hfill
\begin{minipage}[t]{0.73\textwidth}
\vspace{0pt}
\raggedright
\textbf{Name:} {\large{$\mathbf{K7_{164}}$}} (chiral, non-rotatable$^{*}$) \\ \textbf{PD:} {\scriptsize\texttt{[0],[0,1,2,3],[1,4,5,2],[3,5,6,7],[4,7,8,9],[6,9,10,11],[12,13,14,8],[13,12,11,10],[14]}} \\ \textbf{EM:} {\scriptsize\texttt{(B0, A0C0C3D0, B1E0D1B2, B3C2F0E1, C1D3G3F1, D2E3H3H2, H1H0I0E2, G1G0F3F2, G2)}} \\ \textbf{Kauffman bracket:} {\scriptsize $A^{22} + A^{20} - 2A^{18} - 2A^{16} + 3A^{14} + 4A^{12} - 2A^{10} - 5A^{8} + 3A^{4} + A^{2} - 1$} \\ \textbf{Arrow:} {\scriptsize $A^{4} + A^{2}L_1 - 2 - 2L_1/A^{2} + 3/A^{4} + 4L_1/A^{6} - 2/A^{8} - 5L_1/A^{10} - L_2/A^{12} + A^{-12} + 3L_1/A^{14} + L_2/A^{16} - L_1/A^{18}$} \\ \textbf{Mock:} {\scriptsize $-2w^{4} - 2w^{3} + 5w^{2} + 7w - 1 - 6/w - 1/w^{2} + w^{-3}$} \\ \textbf{Affine:} {\scriptsize $-t^{2} + 2t - 2 + 2/t - 1/t^{2}$} \\ \textbf{Yamada:} {\scriptsize $-A^{26} + 2A^{25} + A^{24} - 6A^{23} + 2A^{22} + 6A^{21} - 9A^{20} - A^{19} + 8A^{18} - 7A^{17} - 3A^{16} + 5A^{15} - 2A^{14} - 2A^{13} - 2A^{12} + 6A^{11} - 2A^{10} - 7A^{9} + 10A^{8} - A^{7} - 8A^{6} + 6A^{5} + 2A^{4} - 4A^{3} + A^{2} + A - 1$}
\end{minipage}

\noindent{\color{gray!40}\rule{\textwidth}{0.4pt}}
\vspace{0.9\baselineskip}
\noindent \begin{minipage}[t]{0.25\textwidth}
\vspace{0pt}
\centering
\includegraphics[page=284,width=\linewidth]{knotoids.pdf}
\end{minipage}
\hfill
\begin{minipage}[t]{0.73\textwidth}
\vspace{0pt}
\raggedright
\textbf{Name:} {\large{$\mathbf{K7_{165}}$}} (chiral, non-rotatable$^{*}$) \\ \textbf{PD:} {\scriptsize\texttt{[0],[0,1,2,3],[1,4,5,2],[3,5,6,7],[4,7,8,9],[9,10,11,6],[8,12,13,14],[10,13,12,11],[14]}} \\ \textbf{EM:} {\scriptsize\texttt{(B0, A0C0C3D0, B1E0D1B2, B3C2F3E1, C1D3G0F0, E3H0H3D2, E2H2H1I0, F1G2G1F2, G3)}} \\ \textbf{Kauffman bracket:} {\scriptsize $A^{22} - 3A^{18} - 2A^{16} + 4A^{14} + 4A^{12} - 2A^{10} - 4A^{8} + A^{6} + 3A^{4} - 1$} \\ \textbf{Arrow:} {\scriptsize $A^{16} - A^{12}L_2 - 2A^{12} - 2A^{10}L_1 + A^{8}L_2 + 3A^{8} + 4A^{6}L_1 - 2A^{4} - 4A^{2}L_1 + 1 + 3L_1/A^{2} - L_1/A^{6}$} \\ \textbf{Mock:} {\scriptsize $w^{4} + w^{3} - 4w^{2} - 5w + 5 + 6/w - 1/w^{2} - 2/w^{3}$} \\ \textbf{Affine:} {\scriptsize $-t^{2} - t + 4 - 1/t - 1/t^{2}$} \\ \textbf{Yamada:} {\scriptsize $A^{27} + A^{26} - 2A^{25} - A^{24} + 4A^{23} - A^{22} - 7A^{21} + 5A^{20} + 4A^{19} - 9A^{18} + 3A^{17} + 7A^{16} - 7A^{15} - A^{14} + 5A^{13} - 2A^{12} - 3A^{11} - A^{10} + 6A^{9} - 6A^{8} - 6A^{7} + 10A^{6} - 4A^{5} - 6A^{4} + 6A^{3} - A^{2} - 2A + 1$}
\end{minipage}

\noindent{\color{gray!40}\rule{\textwidth}{0.4pt}}
\vspace{0.9\baselineskip}
\noindent \begin{minipage}[t]{0.25\textwidth}
\vspace{0pt}
\centering
\includegraphics[page=285,width=\linewidth]{knotoids.pdf}
\end{minipage}
\hfill
\begin{minipage}[t]{0.73\textwidth}
\vspace{0pt}
\raggedright
\textbf{Name:} {\large{$\mathbf{K7_{166}}$}} (chiral, non-rotatable$^{*}$) \\ \textbf{PD:} {\scriptsize\texttt{[0],[0,1,2,3],[1,4,5,2],[3,6,7,8],[4,8,9,5],[6,10,11,7],[12,13,10,9],[11,13,14,12],[14]}} \\ \textbf{EM:} {\scriptsize\texttt{(B0, A0C0C3D0, B1E0E3B2, B3F0F3E1, C1D3G3C2, D1G2H0D2, H3H1F1E2, F2G1I0G0, H2)}} \\ \textbf{Kauffman bracket:} {\scriptsize $-A^{21} - A^{19} + A^{17} + 2A^{15} - 3A^{11} - A^{9} + 2A^{7} + 2A^{5} - A^{3} - A$} \\ \textbf{Arrow:} {\scriptsize $A^{6}L_1 + A^{4} - A^{2}L_1 - 2 + 3/A^{4} + L_1/A^{6} - 2/A^{8} - 2L_1/A^{10} + A^{-12} + L_1/A^{14}$} \\ \textbf{Mock:} {\scriptsize $-w^{4} + 3w^{2} + w - 2 - 1/w + 2/w^{2} - 1/w^{4}$} \\ \textbf{Affine:} {\scriptsize $t - 2 + 1/t$} \\ \textbf{Yamada:} {\scriptsize $-A^{26} + 2A^{24} - 2A^{22} + 3A^{21} + 2A^{20} - 3A^{19} + A^{18} + 3A^{17} - 2A^{16} - A^{15} + 3A^{14} - A^{13} - 2A^{12} + 2A^{11} + A^{10} - 2A^{8} + 4A^{7} - 4A^{5} + 4A^{4} - A^{2} + A - 1$}
\end{minipage}

\noindent{\color{gray!40}\rule{\textwidth}{0.4pt}}
\vspace{0.9\baselineskip}
\noindent \begin{minipage}[t]{0.25\textwidth}
\vspace{0pt}
\centering
\includegraphics[page=286,width=\linewidth]{knotoids.pdf}
\end{minipage}
\hfill
\begin{minipage}[t]{0.73\textwidth}
\vspace{0pt}
\raggedright
\textbf{Name:} {\large{$\mathbf{K7_{167}}$}} (chiral, non-rotatable$^{*}$) \\ \textbf{PD:} {\scriptsize\texttt{[0],[0,1,2,3],[1,4,5,2],[3,6,7,8],[4,8,9,5],[6,10,11,12],[7,13,12,14],[14,11,10,9],[13]}} \\ \textbf{EM:} {\scriptsize\texttt{(B0, A0C0C3D0, B1E0E3B2, B3F0G0E1, C1D3H3C2, D1H2H1G2, D2I0F3H0, G3F2F1E2, G1)}} \\ \textbf{Kauffman bracket:} {\scriptsize $A^{23} - 2A^{19} - A^{17} + 3A^{15} + 3A^{13} - 2A^{11} - 4A^{9} + 2A^{5} - A$} \\ \textbf{Arrow:} {\scriptsize $-A^{26}L_1 + 2A^{22}L_1 + A^{20} - 3A^{18}L_1 - A^{16}L_2 - 2A^{16} + 2A^{14}L_1 + A^{12}L_2 + 3A^{12} - 2A^{8} + A^{4}$} \\ \textbf{Mock:} {\scriptsize $w^{3} + 2w^{2} - 4w - 5 + 3/w + 4/w^{2}$} \\ \textbf{Affine:} {\scriptsize $-t^{2} - 2t + 6 - 2/t - 1/t^{2}$} \\ \textbf{Yamada:} {\scriptsize $A^{26} - 2A^{24} + A^{23} + 3A^{22} - 4A^{21} - 2A^{20} + 6A^{19} - A^{18} - 5A^{17} + 5A^{16} + A^{15} - 4A^{14} + 2A^{13} + A^{12} - 5A^{10} + 3A^{9} + 2A^{8} - 8A^{7} + 2A^{6} + 3A^{5} - 5A^{4} - A^{3} + 2A^{2} - 1$}
\end{minipage}

\noindent{\color{gray!40}\rule{\textwidth}{0.4pt}}
\vspace{0.9\baselineskip}
\noindent \begin{minipage}[t]{0.25\textwidth}
\vspace{0pt}
\centering
\includegraphics[page=287,width=\linewidth]{knotoids.pdf}
\end{minipage}
\hfill
\begin{minipage}[t]{0.73\textwidth}
\vspace{0pt}
\raggedright
\textbf{Name:} {\large{$\mathbf{K7_{168}}$}} (chiral, non-rotatable$^{*}$) \\ \textbf{PD:} {\scriptsize\texttt{[0],[0,1,2,3],[1,4,5,2],[3,5,6,7],[7,8,9,4],[6,9,10,11],[12,13,14,8],[13,12,11,10],[14]}} \\ \textbf{EM:} {\scriptsize\texttt{(B0, A0C0C3D0, B1E3D1B2, B3C2F0E0, D3G3F1C1, D2E2H3H2, H1H0I0E1, G1G0F3F2, G2)}} \\ \textbf{Kauffman bracket:} {\scriptsize $A^{22} + A^{20} - A^{18} - 2A^{16} + 3A^{12} + A^{10} - 2A^{8} - 2A^{6} + A^{2} + 1$} \\ \textbf{Arrow:} {\scriptsize $A^{10}L_1 + A^{8} - A^{6}L_1 - 2A^{4} + 3 + L_1/A^{2} - 2/A^{4} - 2L_1/A^{6} - L_2/A^{8} + A^{-8} + L_1/A^{10} + L_2/A^{12}$} \\ \textbf{Mock:} {\scriptsize $w^{3} - w^{2} - w + 4 - 2/w^{2}$} \\ \textbf{Affine:} {\scriptsize $-t^{2} + t + 1/t - 1/t^{2}$} \\ \textbf{Yamada:} {\scriptsize $-A^{28} + A^{26} - 2A^{24} + 2A^{23} + 2A^{22} - 5A^{21} + A^{20} + 3A^{19} - 3A^{18} - A^{17} + A^{15} - 4A^{14} - A^{13} + 3A^{12} - 4A^{11} - A^{10} + 4A^{9} - 2A^{7} + 2A^{6} + 2A^{5} - 2A^{4} - A^{3} + A^{2} - 1$}
\end{minipage}

\noindent{\color{gray!40}\rule{\textwidth}{0.4pt}}
\vspace{0.9\baselineskip}
\noindent \begin{minipage}[t]{0.25\textwidth}
\vspace{0pt}
\centering
\includegraphics[page=288,width=\linewidth]{knotoids.pdf}
\end{minipage}
\hfill
\begin{minipage}[t]{0.73\textwidth}
\vspace{0pt}
\raggedright
\textbf{Name:} {\large{$\mathbf{K7_{169}}$}} (chiral, non-rotatable$^{*}$) \\ \textbf{PD:} {\scriptsize\texttt{[0],[0,1,2,3],[1,4,5,2],[3,5,6,7],[7,8,9,4],[10,11,12,6],[8,12,13,9],[13,14,11,10],[14]}} \\ \textbf{EM:} {\scriptsize\texttt{(B0, A0C0C3D0, B1E3D1B2, B3C2F3E0, D3G0G3C1, H3H2G1D2, E1F2H0E2, G2I0F1F0, H1)}} \\ \textbf{Kauffman bracket:} {\scriptsize $-A^{25} + 2A^{21} - 5A^{17} - 2A^{15} + 5A^{13} + 5A^{11} - 3A^{9} - 5A^{7} + 2A^{3} + A$} \\ \textbf{Arrow:} {\scriptsize $A^{-8} - 2/A^{12} + 5/A^{16} + 2L_1/A^{18} - 5/A^{20} - 5L_1/A^{22} + 3/A^{24} + 5L_1/A^{26} - 2L_1/A^{30} - 1/A^{32}$} \\ \textbf{Mock:} {\scriptsize $3w^{4} + 2w^{3} - 5w^{2} - 5w + 4 + 5/w - 2/w^{2} - 2/w^{3} + w^{-4}$} \\ \textbf{Affine:} {\scriptsize $t - 2 + 1/t$} \\ \textbf{Yamada:} {\scriptsize $3A^{27} + A^{26} - 6A^{25} + 4A^{24} + 11A^{23} - 12A^{22} - A^{21} + 18A^{20} - 9A^{19} - 5A^{18} + 12A^{17} - 2A^{16} - 6A^{15} - 2A^{14} + 8A^{13} - 6A^{12} - 11A^{11} + 14A^{10} - 15A^{8} + 9A^{7} + 7A^{6} - 10A^{5} + A^{4} + 5A^{3} - 2A^{2} - A + 1$}
\end{minipage}

\noindent{\color{gray!40}\rule{\textwidth}{0.4pt}}
\vspace{0.9\baselineskip}
\noindent \begin{minipage}[t]{0.25\textwidth}
\vspace{0pt}
\centering
\includegraphics[page=289,width=\linewidth]{knotoids.pdf}
\end{minipage}
\hfill
\begin{minipage}[t]{0.73\textwidth}
\vspace{0pt}
\raggedright
\textbf{Name:} {\large{$\mathbf{K7_{170}}$}} (chiral, non-rotatable$^{*}$) \\ \textbf{PD:} {\scriptsize\texttt{[0],[0,1,2,3],[1,4,5,2],[3,6,7,8],[8,9,10,4],[5,11,12,6],[13,10,9,7],[11,13,14,12],[14]}} \\ \textbf{EM:} {\scriptsize\texttt{(B0, A0C0C3D0, B1E3F0B2, B3F3G3E0, D3G2G1C1, C2H0H3D1, H1E2E1D2, F1G0I0F2, H2)}} \\ \textbf{Kauffman bracket:} {\scriptsize $-A^{24} + 2A^{20} - A^{18} - 4A^{16} + 6A^{12} + 3A^{10} - 3A^{8} - 3A^{6} + A^{4} + A^{2}$} \\ \textbf{Arrow:} {\scriptsize $-A^{24} + 2A^{20} - A^{18}L_1 - 4A^{16} + 6A^{12} + 3A^{10}L_1 - 3A^{8} - 3A^{6}L_1 + A^{4} + A^{2}L_1$} \\ \textbf{Mock:} {\scriptsize $-w^{4} - w^{3} + 5w^{2} + 2w - 7 - 2/w + 4/w^{2} + w^{-3}$} \\ \textbf{Affine:} {\scriptsize $-t + 2 - 1/t$} \\ \textbf{Yamada:} {\scriptsize $-A^{24} - 3A^{23} + 4A^{22} + 3A^{21} - 8A^{20} + 5A^{19} + 9A^{18} - 10A^{17} + 2A^{16} + 7A^{15} - 5A^{14} - A^{13} + 3A^{11} - 7A^{10} - 4A^{9} + 8A^{8} - 6A^{7} - 7A^{6} + 9A^{5} - A^{4} - 7A^{3} + 4A^{2} + 2A - 2$}
\end{minipage}

\noindent{\color{gray!40}\rule{\textwidth}{0.4pt}}
\vspace{0.9\baselineskip}
\noindent \begin{minipage}[t]{0.25\textwidth}
\vspace{0pt}
\centering
\includegraphics[page=290,width=\linewidth]{knotoids.pdf}
\end{minipage}
\hfill
\begin{minipage}[t]{0.73\textwidth}
\vspace{0pt}
\raggedright
\textbf{Name:} {\large{$\mathbf{K7_{171}}$}} (chiral, non-rotatable$^{*}$) \\ \textbf{PD:} {\scriptsize\texttt{[0],[0,1,2,3],[1,4,5,2],[3,6,7,8],[8,7,9,4],[5,9,10,11],[6,12,13,14],[14,12,11,10],[13]}} \\ \textbf{EM:} {\scriptsize\texttt{(B0, A0C0C3D0, B1E3F0B2, B3G0E1E0, D3D2F1C1, C2E2H3H2, D1H1I0H0, G3G1F3F2, G2)}} \\ \textbf{Kauffman bracket:} {\scriptsize $A^{19} - A^{17} - 3A^{15} + 3A^{11} + 2A^{9} - 2A^{7} - 2A^{5} + A$} \\ \textbf{Arrow:} {\scriptsize $-A^{4} + A^{2}L_1 + 3 - 3/A^{4} - 2L_1/A^{6} + 2/A^{8} + 2L_1/A^{10} - L_1/A^{14}$} \\ \textbf{Mock:} {\scriptsize $w^{4} + w^{3} - 2w^{2} - w + 4 + 1/w - 2/w^{2} - 1/w^{3}$} \\ \textbf{Affine:} {\scriptsize $2t - 4 + 2/t$} \\ \textbf{Yamada:} {\scriptsize $A^{24} - 3A^{22} + 2A^{21} + 3A^{20} - 4A^{19} + A^{18} + 4A^{17} - 2A^{16} + A^{14} + 2A^{13} - 2A^{12} + 5A^{10} - 2A^{9} - A^{8} + 4A^{7} - A^{6} - 3A^{5} + A^{4} + 2A^{3} - 2A^{2} - A + 1$}
\end{minipage}

\noindent{\color{gray!40}\rule{\textwidth}{0.4pt}}
\vspace{0.9\baselineskip}
\noindent \begin{minipage}[t]{0.25\textwidth}
\vspace{0pt}
\centering
\includegraphics[page=291,width=\linewidth]{knotoids.pdf}
\end{minipage}
\hfill
\begin{minipage}[t]{0.73\textwidth}
\vspace{0pt}
\raggedright
\textbf{Name:} {\large{$\mathbf{K7_{172}}$}} (chiral, non-rotatable$^{*}$) \\ \textbf{PD:} {\scriptsize\texttt{[0],[0,1,2,3],[1,4,5,2],[3,6,7,8],[8,9,10,4],[11,12,6,5],[13,10,9,7],[13,14,12,11],[14]}} \\ \textbf{EM:} {\scriptsize\texttt{(B0, A0C0C3D0, B1E3F3B2, B3F2G3E0, D3G2G1C1, H3H2D1C2, H0E2E1D2, G0I0F1F0, H1)}} \\ \textbf{Kauffman bracket:} {\scriptsize $A^{22} - 2A^{18} + 4A^{14} + A^{12} - 5A^{10} - 4A^{8} + 3A^{6} + 4A^{4} - 1$} \\ \textbf{Arrow:} {\scriptsize $A^{16} - 2A^{12} + 4A^{8} + A^{6}L_1 - 5A^{4} - 4A^{2}L_1 + 3 + 4L_1/A^{2} - L_1/A^{6}$} \\ \textbf{Mock:} {\scriptsize $w^{4} + w^{3} - 5w^{2} - 4w + 7 + 4/w - 2/w^{2} - 1/w^{3}$} \\ \textbf{Affine:} {\scriptsize $-t + 2 - 1/t$} \\ \textbf{Yamada:} {\scriptsize $-2A^{26} + A^{25} + 5A^{24} - 4A^{23} - 5A^{22} + 10A^{21} - A^{20} - 11A^{19} + 7A^{18} + 2A^{17} - 7A^{16} + A^{15} + 3A^{14} - 7A^{12} + 4A^{11} + 6A^{10} - 11A^{9} + 2A^{8} + 9A^{7} - 8A^{6} - 3A^{5} + 5A^{4} - 3A^{2} + 1$}
\end{minipage}

\noindent{\color{gray!40}\rule{\textwidth}{0.4pt}}
\vspace{0.9\baselineskip}
\noindent \begin{minipage}[t]{0.25\textwidth}
\vspace{0pt}
\centering
\includegraphics[page=292,width=\linewidth]{knotoids.pdf}
\end{minipage}
\hfill
\begin{minipage}[t]{0.73\textwidth}
\vspace{0pt}
\raggedright
\textbf{Name:} {\large{$\mathbf{K7_{173}}$}} (chiral, non-rotatable$^{*}$) \\ \textbf{PD:} {\scriptsize\texttt{[0],[0,1,2,3],[1,4,5,2],[3,6,7,8],[8,7,9,4],[9,10,11,5],[6,12,13,14],[10,14,12,11],[13]}} \\ \textbf{EM:} {\scriptsize\texttt{(B0, A0C0C3D0, B1E3F3B2, B3G0E1E0, D3D2F0C1, E2H0H3C2, D1H2I0H1, F1G3G1F2, G2)}} \\ \textbf{Kauffman bracket:} {\scriptsize $-A^{20} + 2A^{16} + A^{14} - 2A^{12} - 3A^{10} + A^{8} + 3A^{6} + A^{4} - A^{2}$} \\ \textbf{Arrow:} {\scriptsize $-A^{20} + 2A^{16} + A^{14}L_1 - 2A^{12} - 3A^{10}L_1 + A^{8} + 3A^{6}L_1 + A^{4} - A^{2}L_1$} \\ \textbf{Mock:} {\scriptsize $w^{4} + w^{3} - 2w^{2} - 3w + 2 + 3/w - 1/w^{3}$} \\ \textbf{Affine:} {\scriptsize $0$} \\ \textbf{Yamada:} {\scriptsize $-A^{25} + 2A^{23} - 2A^{21} + A^{20} + 3A^{19} - 3A^{18} - 3A^{17} + 4A^{16} - 2A^{15} - 2A^{14} + 3A^{13} - A^{10} + 4A^{9} - 2A^{7} + 5A^{6} + 2A^{5} - 3A^{4} + A^{3} + 2A^{2} - A - 1$}
\end{minipage}

\noindent{\color{gray!40}\rule{\textwidth}{0.4pt}}
\vspace{0.9\baselineskip}
\noindent \begin{minipage}[t]{0.25\textwidth}
\vspace{0pt}
\centering
\includegraphics[page=293,width=\linewidth]{knotoids.pdf}
\end{minipage}
\hfill
\begin{minipage}[t]{0.73\textwidth}
\vspace{0pt}
\raggedright
\textbf{Name:} {\large{$\mathbf{K7_{174}}$}} (chiral, non-rotatable$^{*}$) \\ \textbf{PD:} {\scriptsize\texttt{[0],[0,1,2,3],[1,4,5,2],[6,7,4,3],[5,8,9,6],[7,9,10,11],[12,13,14,8],[13,12,11,10],[14]}} \\ \textbf{EM:} {\scriptsize\texttt{(B0, A0C0C3D3, B1D2E0B2, E3F0C1B3, C2G3F1D0, D1E2H3H2, H1H0I0E1, G1G0F3F2, G2)}} \\ \textbf{Kauffman bracket:} {\scriptsize $A^{22} + A^{20} - A^{18} - 3A^{16} + A^{14} + 4A^{12} - 4A^{8} - 2A^{6} + 2A^{4} + 2A^{2}$} \\ \textbf{Arrow:} {\scriptsize $A^{16} + A^{14}L_1 - A^{12} - 3A^{10}L_1 + A^{8} + 4A^{6}L_1 - 4A^{2}L_1 - L_2 - 1 + 2L_1/A^{2} + L_2/A^{4} + A^{-4}$} \\ \textbf{Mock:} {\scriptsize $2w^{3} - w^{2} - 6w + 1 + 5/w + 2/w^{2} - 1/w^{3} - 1/w^{4}$} \\ \textbf{Affine:} {\scriptsize $-t^{2} - t + 4 - 1/t - 1/t^{2}$} \\ \textbf{Yamada:} {\scriptsize $-A^{28} + A^{27} + 2A^{26} - 2A^{25} - 2A^{24} + 6A^{23} - A^{22} - 8A^{21} + 6A^{20} - 7A^{18} + 2A^{17} + 2A^{16} - 4A^{14} + 4A^{13} + 3A^{12} - 8A^{11} + 2A^{10} + 5A^{9} - 6A^{8} - 2A^{7} + 5A^{6} - A^{5} - 3A^{4} + 2A^{3} + A^{2} - A - 1$}
\end{minipage}

\noindent{\color{gray!40}\rule{\textwidth}{0.4pt}}
\vspace{0.9\baselineskip}
\noindent \begin{minipage}[t]{0.25\textwidth}
\vspace{0pt}
\centering
\includegraphics[page=294,width=\linewidth]{knotoids.pdf}
\end{minipage}
\hfill
\begin{minipage}[t]{0.73\textwidth}
\vspace{0pt}
\raggedright
\textbf{Name:} {\large{$\mathbf{K7_{175}}$}} (chiral, non-rotatable$^{*}$) \\ \textbf{PD:} {\scriptsize\texttt{[0],[0,1,2,3],[1,4,5,2],[6,7,4,3],[5,8,9,10],[11,12,13,6],[7,13,14,8],[12,11,10,9],[14]}} \\ \textbf{EM:} {\scriptsize\texttt{(B0, A0C0C3D3, B1D2E0B2, F3G0C1B3, C2G3H3H2, H1H0G1D0, D1F2I0E1, F1F0E3E2, G2)}} \\ \textbf{Kauffman bracket:} {\scriptsize $A^{22} + A^{20} - 2A^{18} - 3A^{16} + 2A^{14} + 4A^{12} - A^{10} - 5A^{8} - A^{6} + 3A^{4} + 2A^{2}$} \\ \textbf{Arrow:} {\scriptsize $A^{28} + A^{26}L_1 - 2A^{24} - 3A^{22}L_1 + 2A^{20} + 4A^{18}L_1 - A^{16} - 5A^{14}L_1 - A^{12}L_2 + 3A^{10}L_1 + A^{8}L_2 + A^{8}$} \\ \textbf{Mock:} {\scriptsize $w^{5} - 3w^{3} + w^{2} + 4w - 2 - 5/w + w^{-2} + 3/w^{3} + w^{-4}$} \\ \textbf{Affine:} {\scriptsize $-t^{2} - 2t + 6 - 2/t - 1/t^{2}$} \\ \textbf{Yamada:} {\scriptsize $-A^{28} + 2A^{27} + 2A^{26} - 5A^{25} + 8A^{23} - 6A^{22} - 6A^{21} + 11A^{20} - 3A^{19} - 6A^{18} + 4A^{17} + A^{16} - 2A^{15} - 4A^{14} + 7A^{13} + A^{12} - 9A^{11} + 7A^{10} + 4A^{9} - 10A^{8} + 2A^{7} + 5A^{6} - 4A^{5} - 3A^{4} + A^{3} - A - 1$}
\end{minipage}

\noindent{\color{gray!40}\rule{\textwidth}{0.4pt}}
\vspace{0.9\baselineskip}
\noindent \begin{minipage}[t]{0.25\textwidth}
\vspace{0pt}
\centering
\includegraphics[page=295,width=\linewidth]{knotoids.pdf}
\end{minipage}
\hfill
\begin{minipage}[t]{0.73\textwidth}
\vspace{0pt}
\raggedright
\textbf{Name:} {\large{$\mathbf{K7_{176}}$}} (chiral, rotatable) \\ \textbf{PD:} {\scriptsize\texttt{[0],[0,1,2,3],[1,4,5,2],[6,7,4,3],[8,9,10,5],[11,12,13,6],[7,13,12,8],[9,14,11,10],[14]}} \\ \textbf{EM:} {\scriptsize\texttt{(B0, A0C0C3D3, B1D2E3B2, F3G0C1B3, G3H0H3C2, H2G2G1D0, D1F2F1E0, E1I0F0E2, H1)}} \\ \textbf{Kauffman bracket:} {\scriptsize $A^{24} + 2A^{22} + A^{20} - 4A^{18} - 3A^{16} + 4A^{14} + 5A^{12} - 2A^{10} - 5A^{8} + 3A^{4} - 1$} \\ \textbf{Arrow:} {\scriptsize $A^{-12} + 2L_1/A^{14} + A^{-16} - 4L_1/A^{18} - 3/A^{20} + 4L_1/A^{22} + 5/A^{24} - 2L_1/A^{26} - 5/A^{28} + 3/A^{32} - 1/A^{36}$} \\ \textbf{Mock:} {\scriptsize $w^{6} + 2w^{5} + w^{4} - 4w^{3} - 4w^{2} + 4w + 6 - 2/w - 5/w^{2} + 2/w^{4}$} \\ \textbf{Affine:} {\scriptsize $2t - 4 + 2/t$} \\ \textbf{Yamada:} {\scriptsize $A^{27} - 2A^{25} + 3A^{24} + 5A^{23} - 6A^{22} + A^{21} + 13A^{20} - 4A^{19} - 9A^{18} + 12A^{17} - 13A^{15} + 6A^{14} + 5A^{13} - 7A^{12} - A^{11} + 6A^{10} + 3A^{9} - 12A^{8} + 4A^{7} + 10A^{6} - 15A^{5} + A^{4} + 10A^{3} - 6A^{2} - 2A + 3$}
\end{minipage}

\noindent{\color{gray!40}\rule{\textwidth}{0.4pt}}
\vspace{0.9\baselineskip}
\noindent \begin{minipage}[t]{0.25\textwidth}
\vspace{0pt}
\centering
\includegraphics[page=296,width=\linewidth]{knotoids.pdf}
\end{minipage}
\hfill
\begin{minipage}[t]{0.73\textwidth}
\vspace{0pt}
\raggedright
\textbf{Name:} {\large{$\mathbf{K7_{177}}$}} (chiral, non-rotatable$^{*}$) \\ \textbf{PD:} {\scriptsize\texttt{[0],[0,1,2,3],[1,4,5,2],[6,7,8,3],[4,8,9,10],[5,11,12,6],[7,13,10,9],[11,13,14,12],[14]}} \\ \textbf{EM:} {\scriptsize\texttt{(B0, A0C0C3D3, B1E0F0B2, F3G0E1B3, C1D2G3G2, C2H0H3D0, D1H1E3E2, F1G1I0F2, H2)}} \\ \textbf{Kauffman bracket:} {\scriptsize $A^{24} + A^{22} - 2A^{18} - 3A^{16} + 3A^{12} + 2A^{10} - A^{8} - A^{6} + A^{2}$} \\ \textbf{Arrow:} {\scriptsize $A^{18}L_1 + A^{16} - 2A^{12} - 3A^{10}L_1 + 3A^{6}L_1 + 2A^{4} - A^{2}L_1 - 1 + A^{-4}$} \\ \textbf{Mock:} {\scriptsize $w^{3} - 2w + 1 + 2/w + w^{-2} - 1/w^{3} - 1/w^{4}$} \\ \textbf{Affine:} {\scriptsize $t - 2 + 1/t$} \\ \textbf{Yamada:} {\scriptsize $A^{27} - 2A^{25} + 2A^{23} - A^{22} - 2A^{21} + A^{20} + 4A^{19} - 2A^{18} - 3A^{17} + 5A^{16} - 2A^{15} - 2A^{14} + 3A^{13} - A^{11} - A^{10} + 3A^{9} - A^{8} - 3A^{7} + 4A^{6} + 3A^{5} - 3A^{4} + A^{3} + 3A^{2} - 1$}
\end{minipage}

\noindent{\color{gray!40}\rule{\textwidth}{0.4pt}}
\vspace{0.9\baselineskip}
\noindent \begin{minipage}[t]{0.25\textwidth}
\vspace{0pt}
\centering
\includegraphics[page=297,width=\linewidth]{knotoids.pdf}
\end{minipage}
\hfill
\begin{minipage}[t]{0.73\textwidth}
\vspace{0pt}
\raggedright
\textbf{Name:} {\large{$\mathbf{K7_{178}}$}} (chiral, non-rotatable$^{*}$) \\ \textbf{PD:} {\scriptsize\texttt{[0],[0,1,2,3],[1,4,5,2],[3,6,7,4],[5,7,8,9],[6,10,11,8],[12,13,10,9],[13,12,14,11],[14]}} \\ \textbf{EM:} {\scriptsize\texttt{(B0, A0C0C3D0, B1D3E0B2, B3F0E1C1, C2D2F3G3, D1G2H3E2, H1H0F1E3, G1G0I0F2, H2)}} \\ \textbf{Kauffman bracket:} {\scriptsize $-A^{24} + A^{22} + 4A^{20} + A^{18} - 5A^{16} - 4A^{14} + 3A^{12} + 4A^{10} - A^{8} - 2A^{6} + A^{2}$} \\ \textbf{Arrow:} {\scriptsize $-A^{6}L_1 + A^{4}L_2 + 4A^{2}L_1 - L_2 + 2 - 5L_1/A^{2} - 4/A^{4} + 3L_1/A^{6} + 4/A^{8} - L_1/A^{10} - 2/A^{12} + A^{-16}$} \\ \textbf{Mock:} {\scriptsize $-w^{3} - w^{2} + 5w + 5 - 6/w - 4/w^{2} + 2/w^{3} + w^{-4}$} \\ \textbf{Affine:} {\scriptsize $t^{2} + t - 4 + 1/t + t^{-2}$} \\ \textbf{Yamada:} {\scriptsize $2A^{25} - A^{24} - 4A^{23} + 6A^{22} + 3A^{21} - 11A^{20} + 6A^{19} + 8A^{18} - 11A^{17} + 3A^{16} + 7A^{15} - 5A^{14} - 2A^{13} + A^{12} + 5A^{11} - 6A^{10} - 2A^{9} + 13A^{8} - 6A^{7} - 6A^{6} + 11A^{5} - 2A^{4} - 6A^{3} + 4A^{2} + A - 2$}
\end{minipage}

\noindent{\color{gray!40}\rule{\textwidth}{0.4pt}}
\vspace{0.9\baselineskip}
\noindent \begin{minipage}[t]{0.25\textwidth}
\vspace{0pt}
\centering
\includegraphics[page=298,width=\linewidth]{knotoids.pdf}
\end{minipage}
\hfill
\begin{minipage}[t]{0.73\textwidth}
\vspace{0pt}
\raggedright
\textbf{Name:} {\large{$\mathbf{K7_{179}}$}} (chiral, non-rotatable$^{*}$) \\ \textbf{PD:} {\scriptsize\texttt{[0],[0,1,2,3],[1,4,5,2],[3,6,7,4],[7,8,9,5],[6,10,11,8],[12,13,10,9],[13,12,14,11],[14]}} \\ \textbf{EM:} {\scriptsize\texttt{(B0, A0C0C3D0, B1D3E3B2, B3F0E0C1, D2F3G3C2, D1G2H3E1, H1H0F1E2, G1G0I0F2, H2)}} \\ \textbf{Kauffman bracket:} {\scriptsize $A^{23} - 2A^{19} - 2A^{17} + A^{15} + 3A^{13} - 2A^{9} + A^{5} - A$} \\ \textbf{Arrow:} {\scriptsize $-A^{14}L_1 + 2A^{10}L_1 + A^{8}L_2 + A^{8} - A^{6}L_1 - A^{4}L_2 - 2A^{4} + 2 - 1/A^{4} + A^{-8}$} \\ \textbf{Mock:} {\scriptsize $-w^{3} - 2w^{2} + 2w + 5 - 1/w - 2/w^{2}$} \\ \textbf{Affine:} {\scriptsize $t^{2} - 2 + t^{-2}$} \\ \textbf{Yamada:} {\scriptsize $A^{23} - 2A^{21} + 3A^{19} - A^{18} - A^{17} + 4A^{16} - A^{15} - 2A^{14} + A^{13} - A^{12} - A^{10} + 3A^{9} + A^{8} - 2A^{7} + 4A^{6} + 2A^{5} - 2A^{4} + A^{3} + A^{2} - A - 1$}
\end{minipage}

\noindent{\color{gray!40}\rule{\textwidth}{0.4pt}}
\vspace{0.9\baselineskip}
\noindent \begin{minipage}[t]{0.25\textwidth}
\vspace{0pt}
\centering
\includegraphics[page=299,width=\linewidth]{knotoids.pdf}
\end{minipage}
\hfill
\begin{minipage}[t]{0.73\textwidth}
\vspace{0pt}
\raggedright
\textbf{Name:} {\large{$\mathbf{K7_{180}}$}} (chiral, rotatable) \\ \textbf{PD:} {\scriptsize\texttt{[0],[0,1,2,3],[1,4,5,2],[3,5,6,7],[4,7,8,9],[10,11,8,6],[9,12,13,10],[11,13,12,14],[14]}} \\ \textbf{EM:} {\scriptsize\texttt{(B0, A0C0C3D0, B1E0D1B2, B3C2F3E1, C1D3F2G0, G3H0E2D2, E3H2H1F0, F1G2G1I0, H3)}} \\ \textbf{Kauffman bracket:} {\scriptsize $A^{22} + 2A^{20} - 4A^{16} - A^{14} + 4A^{12} + 3A^{10} - 2A^{8} - 3A^{6} + A^{2}$} \\ \textbf{Arrow:} {\scriptsize $A^{-8} + 2L_1/A^{10} - 4L_1/A^{14} - 1/A^{16} + 4L_1/A^{18} + 3/A^{20} - 2L_1/A^{22} - 3/A^{24} + A^{-28}$} \\ \textbf{Mock:} {\scriptsize $w^{4} + 2w^{3} + w^{2} - 4w - 3 + 4/w + 3/w^{2} - 2/w^{3} - 1/w^{4}$} \\ \textbf{Affine:} {\scriptsize $2t - 4 + 2/t$} \\ \textbf{Yamada:} {\scriptsize $A^{24} + A^{23} - 3A^{22} + A^{21} + 6A^{20} - 4A^{19} - 4A^{18} + 9A^{17} - 5A^{15} + 7A^{14} + 3A^{13} - 4A^{12} + A^{11} + A^{10} - 7A^{8} + 3A^{7} + 5A^{6} - 10A^{5} + 3A^{4} + 6A^{3} - 5A^{2} - A + 2$}
\end{minipage}

\noindent{\color{gray!40}\rule{\textwidth}{0.4pt}}
\vspace{0.9\baselineskip}
\noindent \begin{minipage}[t]{0.25\textwidth}
\vspace{0pt}
\centering
\includegraphics[page=300,width=\linewidth]{knotoids.pdf}
\end{minipage}
\hfill
\begin{minipage}[t]{0.73\textwidth}
\vspace{0pt}
\raggedright
\textbf{Name:} {\large{$\mathbf{K7_{181}}$}} (chiral, non-rotatable$^{*}$) \\ \textbf{PD:} {\scriptsize\texttt{[0],[0,1,2,3],[1,4,5,2],[3,6,7,8],[4,8,9,10],[10,11,12,5],[6,12,11,13],[13,14,9,7],[14]}} \\ \textbf{EM:} {\scriptsize\texttt{(B0, A0C0C3D0, B1E0F3B2, B3G0H3E1, C1D3H2F0, E3G2G1C2, D1F2F1H0, G3I0E2D2, H1)}} \\ \textbf{Kauffman bracket:} {\scriptsize $A^{22} - 3A^{18} - 2A^{16} + 4A^{14} + 5A^{12} - A^{10} - 4A^{8} + 2A^{4} - 1$} \\ \textbf{Arrow:} {\scriptsize $A^{4} - L_2 - 2 - 2L_1/A^{2} + L_2/A^{4} + 3/A^{4} + 5L_1/A^{6} - 1/A^{8} - 4L_1/A^{10} + 2L_1/A^{14} - L_1/A^{18}$} \\ \textbf{Mock:} {\scriptsize $-w^{4} - 2w^{3} + 3w^{2} + 5w - 2 - 4/w + 2/w^{2} + 2/w^{3} - 1/w^{4} - 1/w^{5}$} \\ \textbf{Affine:} {\scriptsize $-t^{2} + t + 1/t - 1/t^{2}$} \\ \textbf{Yamada:} {\scriptsize $-A^{26} - A^{25} + 2A^{24} + A^{23} - 4A^{22} + 7A^{20} - 3A^{19} - 6A^{18} + 9A^{17} + 2A^{16} - 7A^{15} + 6A^{14} + 4A^{13} - 5A^{12} + A^{11} + 3A^{10} + 2A^{9} - 7A^{8} + 3A^{7} + 7A^{6} - 10A^{5} + 7A^{3} - 4A^{2} - 2A + 2$}
\end{minipage}

\noindent{\color{gray!40}\rule{\textwidth}{0.4pt}}
\vspace{0.9\baselineskip}
\noindent \begin{minipage}[t]{0.25\textwidth}
\vspace{0pt}
\centering
\includegraphics[page=301,width=\linewidth]{knotoids.pdf}
\end{minipage}
\hfill
\begin{minipage}[t]{0.73\textwidth}
\vspace{0pt}
\raggedright
\textbf{Name:} {\large{$\mathbf{K7_{182}}$}} (chiral, rotatable) \\ \textbf{PD:} {\scriptsize\texttt{[0],[0,1,2,3],[1,4,5,2],[3,5,6,7],[7,8,9,4],[10,11,8,6],[9,12,13,10],[11,13,12,14],[14]}} \\ \textbf{EM:} {\scriptsize\texttt{(B0, A0C0C3D0, B1E3D1B2, B3C2F3E0, D3F2G0C1, G3H0E1D2, E2H2H1F0, F1G2G1I0, H3)}} \\ \textbf{Kauffman bracket:} {\scriptsize $A^{22} - 2A^{18} + 5A^{14} + 2A^{12} - 5A^{10} - 4A^{8} + 3A^{6} + 4A^{4} - A^{2} - 2$} \\ \textbf{Arrow:} {\scriptsize $A^{-8} - 2/A^{12} + 5/A^{16} + 2L_1/A^{18} - 5/A^{20} - 4L_1/A^{22} + 3/A^{24} + 4L_1/A^{26} - 1/A^{28} - 2L_1/A^{30}$} \\ \textbf{Mock:} {\scriptsize $3w^{4} + 2w^{3} - 5w^{2} - 4w + 5 + 4/w - 3/w^{2} - 2/w^{3} + w^{-4}$} \\ \textbf{Affine:} {\scriptsize $2t - 4 + 2/t$} \\ \textbf{Yamada:} {\scriptsize $3A^{27} + A^{26} - 4A^{25} + 5A^{24} + 7A^{23} - 10A^{22} + 2A^{21} + 12A^{20} - 9A^{19} - 3A^{18} + 9A^{17} - 3A^{16} - 4A^{15} + 7A^{13} - 5A^{12} - 7A^{11} + 12A^{10} - 4A^{9} - 11A^{8} + 8A^{7} + 3A^{6} - 7A^{5} + 2A^{4} + 4A^{3} - 2A^{2} - A + 1$}
\end{minipage}

\noindent{\color{gray!40}\rule{\textwidth}{0.4pt}}
\vspace{0.9\baselineskip}
\noindent \begin{minipage}[t]{0.25\textwidth}
\vspace{0pt}
\centering
\includegraphics[page=302,width=\linewidth]{knotoids.pdf}
\end{minipage}
\hfill
\begin{minipage}[t]{0.73\textwidth}
\vspace{0pt}
\raggedright
\textbf{Name:} {\large{$\mathbf{K7_{183}}$}} (chiral, non-rotatable$^{*}$) \\ \textbf{PD:} {\scriptsize\texttt{[0],[0,1,2,3],[1,4,5,2],[3,5,6,7],[7,8,9,4],[10,11,8,6],[11,12,13,9],[10,13,12,14],[14]}} \\ \textbf{EM:} {\scriptsize\texttt{(B0, A0C0C3D0, B1E3D1B2, B3C2F3E0, D3F2G3C1, H0G0E1D2, F1H2H1E2, F0G2G1I0, H3)}} \\ \textbf{Kauffman bracket:} {\scriptsize $-A^{26} - A^{24} + 2A^{22} + 3A^{20} - 2A^{18} - 6A^{16} + 7A^{12} + 3A^{10} - 3A^{8} - 3A^{6} + A^{4} + A^{2}$} \\ \textbf{Arrow:} {\scriptsize $-A^{26}L_1 - A^{24} + 2A^{22}L_1 + 3A^{20} - 2A^{18}L_1 - A^{16}L_2 - 5A^{16} + A^{12}L_2 + 6A^{12} + 3A^{10}L_1 - 3A^{8} - 3A^{6}L_1 + A^{4} + A^{2}L_1$} \\ \textbf{Mock:} {\scriptsize $-w^{4} + 6w^{2} - 9 - 1/w + 5/w^{2} + w^{-3}$} \\ \textbf{Affine:} {\scriptsize $-t^{2} - t + 4 - 1/t - 1/t^{2}$} \\ \textbf{Yamada:} {\scriptsize $-2A^{27} + 5A^{25} - 3A^{24} - 7A^{23} + 12A^{22} + 4A^{21} - 16A^{20} + 9A^{19} + 12A^{18} - 17A^{17} + 10A^{15} - 7A^{14} - 5A^{13} + 3A^{12} + 9A^{11} - 12A^{10} - 4A^{9} + 18A^{8} - 11A^{7} - 11A^{6} + 14A^{5} - 2A^{4} - 9A^{3} + 4A^{2} + 2A - 2$}
\end{minipage}

\noindent{\color{gray!40}\rule{\textwidth}{0.4pt}}
\vspace{0.9\baselineskip}
\noindent \begin{minipage}[t]{0.25\textwidth}
\vspace{0pt}
\centering
\includegraphics[page=303,width=\linewidth]{knotoids.pdf}
\end{minipage}
\hfill
\begin{minipage}[t]{0.73\textwidth}
\vspace{0pt}
\raggedright
\textbf{Name:} {\large{$\mathbf{K7_{184}}$}} (chiral, non-rotatable$^{*}$) \\ \textbf{PD:} {\scriptsize\texttt{[0],[0,1,2,3],[1,4,5,2],[3,6,7,8],[8,9,10,4],[5,10,11,12],[12,11,13,6],[13,14,9,7],[14]}} \\ \textbf{EM:} {\scriptsize\texttt{(B0, A0C0C3D0, B1E3F0B2, B3G3H3E0, D3H2F1C1, C2E2G1G0, F3F2H0D1, G2I0E1D2, H1)}} \\ \textbf{Kauffman bracket:} {\scriptsize $-A^{16} + A^{14} + 2A^{12} - A^{10} - 3A^{8} - A^{6} + 2A^{4} + 2A^{2}$} \\ \textbf{Arrow:} {\scriptsize $-A^{22}L_1 + A^{20} + 2A^{18}L_1 - A^{16} - 3A^{14}L_1 - A^{12}L_2 + 2A^{10}L_1 + A^{8}L_2 + A^{8}$} \\ \textbf{Mock:} {\scriptsize $w^{5} + w^{4} - 2w^{3} - w^{2} + w - 1 - 2/w + w^{-2} + 2/w^{3} + w^{-4}$} \\ \textbf{Affine:} {\scriptsize $-t^{2} - 2t + 6 - 2/t - 1/t^{2}$} \\ \textbf{Yamada:} {\scriptsize $A^{24} + 2A^{23} - 2A^{22} - 2A^{21} + 3A^{20} - 3A^{19} - 2A^{18} + 2A^{17} + A^{15} + 3A^{13} - 3A^{11} + A^{10} - A^{9} - 4A^{8} + A^{7} + A^{6} - 2A^{5} - A^{4} + A^{3} - A - 1$}
\end{minipage}

\noindent{\color{gray!40}\rule{\textwidth}{0.4pt}}
\vspace{0.9\baselineskip}
\noindent \begin{minipage}[t]{0.25\textwidth}
\vspace{0pt}
\centering
\includegraphics[page=304,width=\linewidth]{knotoids.pdf}
\end{minipage}
\hfill
\begin{minipage}[t]{0.73\textwidth}
\vspace{0pt}
\raggedright
\textbf{Name:} {\large{$\mathbf{K7_{185}}$}} (chiral, non-rotatable$^{*}$) \\ \textbf{PD:} {\scriptsize\texttt{[0],[0,1,2,3],[1,4,5,2],[3,6,7,8],[8,9,10,4],[10,11,12,5],[6,12,11,13],[13,14,9,7],[14]}} \\ \textbf{EM:} {\scriptsize\texttt{(B0, A0C0C3D0, B1E3F3B2, B3G0H3E0, D3H2F0C1, E2G2G1C2, D1F2F1H0, G3I0E1D2, H1)}} \\ \textbf{Kauffman bracket:} {\scriptsize $-A^{20} + 2A^{16} + A^{14} - 2A^{12} - 3A^{10} + A^{8} + 4A^{6} + A^{4} - A^{2} - 1$} \\ \textbf{Arrow:} {\scriptsize $-A^{14}L_1 + 2A^{10}L_1 + A^{8} - 2A^{6}L_1 - A^{4}L_2 - 2A^{4} + A^{2}L_1 + L_2 + 3 + L_1/A^{2} - 1/A^{4} - L_1/A^{6}$} \\ \textbf{Mock:} {\scriptsize $-w^{3} - w^{2} + 2w + 4 - 1/w - 3/w^{2} + w^{-4}$} \\ \textbf{Affine:} {\scriptsize $-t^{2} + 2 - 1/t^{2}$} \\ \textbf{Yamada:} {\scriptsize $-A^{25} + 2A^{23} - 4A^{21} + 4A^{19} - 4A^{18} - 2A^{17} + 7A^{16} - A^{14} + 5A^{13} + A^{12} - 2A^{10} + 3A^{9} - 2A^{8} - 6A^{7} + 5A^{6} + 2A^{5} - 4A^{4} + 2A^{3} + 3A^{2} - A - 1$}
\end{minipage}

\noindent{\color{gray!40}\rule{\textwidth}{0.4pt}}
\vspace{0.9\baselineskip}
\noindent \begin{minipage}[t]{0.25\textwidth}
\vspace{0pt}
\centering
\includegraphics[page=305,width=\linewidth]{knotoids.pdf}
\end{minipage}
\hfill
\begin{minipage}[t]{0.73\textwidth}
\vspace{0pt}
\raggedright
\textbf{Name:} {\large{$\mathbf{K7_{186}}$}} (chiral, non-rotatable$^{*}$) \\ \textbf{PD:} {\scriptsize\texttt{[0],[0,1,2,3],[1,4,5,2],[6,7,8,3],[4,8,9,10],[5,10,11,12],[12,11,13,6],[7,13,14,9],[14]}} \\ \textbf{EM:} {\scriptsize\texttt{(B0, A0C0C3D3, B1E0F0B2, G3H0E1B3, C1D2H3F1, C2E3G1G0, F3F2H1D0, D1G2I0E2, H2)}} \\ \textbf{Kauffman bracket:} {\scriptsize $-A^{20} - A^{18} + 2A^{16} + 3A^{14} - A^{12} - 4A^{10} - A^{8} + 3A^{6} + 2A^{4} - A^{2}$} \\ \textbf{Arrow:} {\scriptsize $-A^{20} - A^{18}L_1 + 2A^{16} + 3A^{14}L_1 - A^{12} - 4A^{10}L_1 - A^{8}L_2 + 3A^{6}L_1 + A^{4}L_2 + A^{4} - A^{2}L_1$} \\ \textbf{Mock:} {\scriptsize $w^{4} + 2w^{3} - w^{2} - 5w + 4/w + w^{-2} - 1/w^{3}$} \\ \textbf{Affine:} {\scriptsize $-t^{2} + 2 - 1/t^{2}$} \\ \textbf{Yamada:} {\scriptsize $-2A^{25} + 4A^{23} - 3A^{21} + 3A^{20} + 4A^{19} - 5A^{18} - 2A^{17} + 5A^{16} - 4A^{15} - 3A^{14} + 3A^{13} - A^{12} - 2A^{11} - A^{10} + 5A^{9} - 3A^{8} - 5A^{7} + 5A^{6} - A^{5} - 5A^{4} + 2A^{3} + 2A^{2} - A - 1$}
\end{minipage}

\noindent{\color{gray!40}\rule{\textwidth}{0.4pt}}
\vspace{0.9\baselineskip}
\noindent \begin{minipage}[t]{0.25\textwidth}
\vspace{0pt}
\centering
\includegraphics[page=306,width=\linewidth]{knotoids.pdf}
\end{minipage}
\hfill
\begin{minipage}[t]{0.73\textwidth}
\vspace{0pt}
\raggedright
\textbf{Name:} {\large{$\mathbf{K7_{187}}$}} (chiral, non-rotatable$^{*}$) \\ \textbf{PD:} {\scriptsize\texttt{[0],[0,1,2,3],[1,4,5,2],[3,6,7,8],[8,9,10,4],[5,11,12,6],[13,14,9,7],[14,13,11,10],[12]}} \\ \textbf{EM:} {\scriptsize\texttt{(B0, A0C0C3D0, B1E3F0B2, B3F3G3E0, D3G2H3C1, C2H2I0D1, H1H0E1D2, G1G0F1E2, F2)}} \\ \textbf{Kauffman bracket:} {\scriptsize $-A^{23} + 2A^{19} - A^{17} - 3A^{15} + 4A^{11} + A^{9} - 3A^{7} - 2A^{5} + A^{3} + A$} \\ \textbf{Arrow:} {\scriptsize $A^{2}L_1 - 2L_1/A^{2} + A^{-4} + 3L_1/A^{6} + L_2/A^{8} - 1/A^{8} - 4L_1/A^{10} - 2L_2/A^{12} + A^{-12} + 3L_1/A^{14} + 2L_2/A^{16} - L_1/A^{18} - L_2/A^{20}$} \\ \textbf{Mock:} {\scriptsize $w^{5} + w^{4} - 3w^{3} - w^{2} + 5w + 1 - 4/w + w^{-3}$} \\ \textbf{Affine:} {\scriptsize $t^{2} + t - 4 + 1/t + t^{-2}$} \\ \textbf{Yamada:} {\scriptsize $-A^{25} - 2A^{24} + 3A^{23} - 2A^{22} - 4A^{21} + 6A^{20} - 2A^{19} - 4A^{18} + 4A^{17} - A^{16} - 2A^{15} - 2A^{14} + 2A^{13} - 4A^{11} + 3A^{10} + 2A^{9} - 5A^{8} + A^{7} + 4A^{6} - 3A^{5} - A^{4} + 3A^{3} - A^{2} - A + 1$}
\end{minipage}

\noindent{\color{gray!40}\rule{\textwidth}{0.4pt}}
\vspace{0.9\baselineskip}
\noindent \begin{minipage}[t]{0.25\textwidth}
\vspace{0pt}
\centering
\includegraphics[page=307,width=\linewidth]{knotoids.pdf}
\end{minipage}
\hfill
\begin{minipage}[t]{0.73\textwidth}
\vspace{0pt}
\raggedright
\textbf{Name:} {\large{$\mathbf{K7_{188}}$}} (chiral, non-rotatable$^{*}$) \\ \textbf{PD:} {\scriptsize\texttt{[0],[0,1,2,3],[1,4,5,2],[3,6,7,8],[8,9,10,4],[5,10,11,12],[6,13,9,7],[13,12,14,11],[14]}} \\ \textbf{EM:} {\scriptsize\texttt{(B0, A0C0C3D0, B1E3F0B2, B3G0G3E0, D3G2F1C1, C2E2H3H1, D1H0E1D2, G1F3I0F2, H2)}} \\ \textbf{Kauffman bracket:} {\scriptsize $-A^{24} + 2A^{20} - 4A^{16} - A^{14} + 5A^{12} + 4A^{10} - 2A^{8} - 3A^{6} + A^{2}$} \\ \textbf{Arrow:} {\scriptsize $-A^{18}L_1 + 2A^{14}L_1 - A^{12}L_2 + A^{12} - 4A^{10}L_1 + A^{8}L_2 - 2A^{8} + 5A^{6}L_1 + 4A^{4} - 2A^{2}L_1 - 3 + A^{-4}$} \\ \textbf{Mock:} {\scriptsize $w^{3} + w^{2} - 5w - 3 + 6/w + 4/w^{2} - 2/w^{3} - 1/w^{4}$} \\ \textbf{Affine:} {\scriptsize $-t^{2} - t + 4 - 1/t - 1/t^{2}$} \\ \textbf{Yamada:} {\scriptsize $A^{25} + A^{24} - 4A^{23} + 2A^{22} + 6A^{21} - 8A^{20} - A^{19} + 10A^{18} - 7A^{17} - 3A^{16} + 6A^{15} - 3A^{14} - 2A^{13} - 2A^{12} + 5A^{11} - 3A^{10} - 7A^{9} + 8A^{8} - A^{7} - 9A^{6} + 6A^{5} + 3A^{4} - 6A^{3} + A^{2} + 2A - 1$}
\end{minipage}

\noindent{\color{gray!40}\rule{\textwidth}{0.4pt}}
\vspace{0.9\baselineskip}
\noindent \begin{minipage}[t]{0.25\textwidth}
\vspace{0pt}
\centering
\includegraphics[page=308,width=\linewidth]{knotoids.pdf}
\end{minipage}
\hfill
\begin{minipage}[t]{0.73\textwidth}
\vspace{0pt}
\raggedright
\textbf{Name:} {\large{$\mathbf{K7_{189}}$}} (chiral, non-rotatable$^{*}$) \\ \textbf{PD:} {\scriptsize\texttt{[0],[0,1,2,3],[1,4,5,2],[3,6,7,8],[8,9,10,4],[11,12,6,5],[13,14,9,7],[14,13,11,10],[12]}} \\ \textbf{EM:} {\scriptsize\texttt{(B0, A0C0C3D0, B1E3F3B2, B3F2G3E0, D3G2H3C1, H2I0D1C2, H1H0E1D2, G1G0F0E2, F1)}} \\ \textbf{Kauffman bracket:} {\scriptsize $A^{22} - 2A^{18} + 3A^{14} + 2A^{12} - 4A^{10} - 3A^{8} + 2A^{6} + 3A^{4} - 1$} \\ \textbf{Arrow:} {\scriptsize $A^{10}L_1 - 2A^{6}L_1 + 3A^{2}L_1 + L_2 + 1 - 4L_1/A^{2} - 2L_2/A^{4} - 1/A^{4} + 2L_1/A^{6} + 2L_2/A^{8} + A^{-8} - L_2/A^{12}$} \\ \textbf{Mock:} {\scriptsize $w^{5} + w^{4} - 3w^{3} - 3w^{2} + 4w + 4 - 3/w - 1/w^{2} + w^{-3}$} \\ \textbf{Affine:} {\scriptsize $t^{2} - 2 + t^{-2}$} \\ \textbf{Yamada:} {\scriptsize $-2A^{26} + 3A^{24} - 2A^{23} - 4A^{22} + 6A^{21} + 2A^{20} - 7A^{19} + 4A^{18} + 2A^{17} - 6A^{16} - A^{15} + A^{14} + A^{13} - 5A^{12} + 2A^{11} + 5A^{10} - 7A^{9} + 6A^{7} - 4A^{6} - 3A^{5} + 3A^{4} + A^{3} - 2A^{2} + 1$}
\end{minipage}

\noindent{\color{gray!40}\rule{\textwidth}{0.4pt}}
\vspace{0.9\baselineskip}
\noindent \begin{minipage}[t]{0.25\textwidth}
\vspace{0pt}
\centering
\includegraphics[page=309,width=\linewidth]{knotoids.pdf}
\end{minipage}
\hfill
\begin{minipage}[t]{0.73\textwidth}
\vspace{0pt}
\raggedright
\textbf{Name:} {\large{$\mathbf{K7_{190}}$}} (chiral, non-rotatable$^{*}$) \\ \textbf{PD:} {\scriptsize\texttt{[0],[0,1,2,3],[1,4,5,2],[6,7,8,3],[4,8,9,10],[5,11,12,6],[7,13,14,9],[10,14,13,11],[12]}} \\ \textbf{EM:} {\scriptsize\texttt{(B0, A0C0C3D3, B1E0F0B2, F3G0E1B3, C1D2G3H0, C2H3I0D0, D1H2H1E2, E3G2G1F1, F2)}} \\ \textbf{Kauffman bracket:} {\scriptsize $A^{18} + A^{16} - 2A^{12} - A^{10} + A^{8} + 2A^{6} - A^{2}$} \\ \textbf{Arrow:} {\scriptsize $1 + L_1/A^{2} + L_2/A^{4} - 1/A^{4} - 2L_1/A^{6} - 2L_2/A^{8} + A^{-8} + L_1/A^{10} + 2L_2/A^{12} - L_2/A^{16}$} \\ \textbf{Mock:} {\scriptsize $-w^{6} - w^{5} + 2w^{4} + 2w^{3} - w - 2 + 3/w^{2} - 1/w^{4}$} \\ \textbf{Affine:} {\scriptsize $t^{2} - 2 + t^{-2}$} \\ \textbf{Yamada:} {\scriptsize $A^{23} - A^{21} + A^{20} + A^{19} - A^{18} + A^{16} + A^{15} - A^{14} + A^{13} + 2A^{12} - A^{11} + A^{10} + 2A^{9} - A^{8} + A^{5} - A^{4} + 2A^{2} - A - 1$}
\end{minipage}

\noindent{\color{gray!40}\rule{\textwidth}{0.4pt}}
\vspace{0.9\baselineskip}
\noindent \begin{minipage}[t]{0.25\textwidth}
\vspace{0pt}
\centering
\includegraphics[page=310,width=\linewidth]{knotoids.pdf}
\end{minipage}
\hfill
\begin{minipage}[t]{0.73\textwidth}
\vspace{0pt}
\raggedright
\textbf{Name:} {\large{$\mathbf{K7_{191}}$}} (chiral, non-rotatable$^{*}$) \\ \textbf{PD:} {\scriptsize\texttt{[0],[0,1,2,3],[1,4,5,2],[6,7,8,3],[4,8,9,10],[11,12,6,5],[7,13,14,9],[10,14,13,11],[12]}} \\ \textbf{EM:} {\scriptsize\texttt{(B0, A0C0C3D3, B1E0F3B2, F2G0E1B3, C1D2G3H0, H3I0D0C2, D1H2H1E2, E3G2G1F0, F1)}} \\ \textbf{Kauffman bracket:} {\scriptsize $2A^{22} + A^{20} - 3A^{18} - 4A^{16} + 3A^{14} + 5A^{12} - A^{10} - 4A^{8} + 3A^{4} - 1$} \\ \textbf{Arrow:} {\scriptsize $A^{4}L_2 + A^{4} + A^{2}L_1 - 2L_2 - 1 - 4L_1/A^{2} + 2L_2/A^{4} + A^{-4} + 5L_1/A^{6} - L_2/A^{8} - 4L_1/A^{10} + 3L_1/A^{14} - L_1/A^{18}$} \\ \textbf{Mock:} {\scriptsize $-w^{4} - 2w^{3} + 2w^{2} + 5w - 2 - 6/w + 3/w^{2} + 4/w^{3} - 1/w^{4} - 1/w^{5}$} \\ \textbf{Affine:} {\scriptsize $t^{2} - t - 1/t + t^{-2}$} \\ \textbf{Yamada:} {\scriptsize $A^{27} - 2A^{25} + A^{24} + 3A^{23} - 4A^{22} - 3A^{21} + 8A^{20} - 9A^{18} + 8A^{17} + 5A^{16} - 9A^{15} + 4A^{14} + 5A^{13} - 5A^{12} - 2A^{11} + 3A^{10} + 5A^{9} - 7A^{8} + A^{7} + 11A^{6} - 8A^{5} - 3A^{4} + 8A^{3} - 4A^{2} - 3A + 2$}
\end{minipage}

\noindent{\color{gray!40}\rule{\textwidth}{0.4pt}}
\vspace{0.9\baselineskip}
\noindent \begin{minipage}[t]{0.25\textwidth}
\vspace{0pt}
\centering
\includegraphics[page=311,width=\linewidth]{knotoids.pdf}
\end{minipage}
\hfill
\begin{minipage}[t]{0.73\textwidth}
\vspace{0pt}
\raggedright
\textbf{Name:} {\large{$\mathbf{K7_{192}}$}} (chiral, rotatable) \\ \textbf{PD:} {\scriptsize\texttt{[0],[0,1,2,3],[1,4,5,6],[6,7,8,2],[9,10,11,3],[4,11,10,12],[12,13,7,5],[13,14,9,8],[14]}} \\ \textbf{EM:} {\scriptsize\texttt{(B0, A0C0D3E3, B1F0G3D0, C3G2H3B2, H2F2F1B3, C1E2E1G0, F3H0D1C2, G1I0E0D2, H1)}} \\ \textbf{Kauffman bracket:} {\scriptsize $2A^{22} + 2A^{20} - 4A^{18} - 5A^{16} + 2A^{14} + 6A^{12} - 4A^{8} + 3A^{4} - 1$} \\ \textbf{Arrow:} {\scriptsize $2A^{10}L_1 + 2A^{8} - 4A^{6}L_1 - 5A^{4} + 2A^{2}L_1 + 6 - 4/A^{4} + 3/A^{8} - 1/A^{12}$} \\ \textbf{Mock:} {\scriptsize $w^{4} - 4w^{2} + 2w + 8 - 4/w - 6/w^{2} + 2/w^{3} + 2/w^{4}$} \\ \textbf{Affine:} {\scriptsize $0$} \\ \textbf{Yamada:} {\scriptsize $A^{27} - A^{26} - 3A^{25} + 2A^{24} + 5A^{23} - 3A^{22} - 6A^{21} + 8A^{20} + 4A^{19} - 13A^{18} + 5A^{17} + 9A^{16} - 11A^{15} + A^{14} + 8A^{13} - 3A^{12} - 3A^{11} + 3A^{10} + 7A^{9} - 8A^{8} - 4A^{7} + 14A^{6} - 5A^{5} - 8A^{4} + 10A^{3} - 4A + 1$}
\end{minipage}

\noindent{\color{gray!40}\rule{\textwidth}{0.4pt}}
\vspace{0.9\baselineskip}
\noindent \begin{minipage}[t]{0.25\textwidth}
\vspace{0pt}
\centering
\includegraphics[page=312,width=\linewidth]{knotoids.pdf}
\end{minipage}
\hfill
\begin{minipage}[t]{0.73\textwidth}
\vspace{0pt}
\raggedright
\textbf{Name:} {\large{$\mathbf{K7_{193}}$}} (chiral, non-rotatable$^{*}$) \\ \textbf{PD:} {\scriptsize\texttt{[0],[0,1,2,3],[1,4,5,2],[3,6,7,4],[5,7,8,9],[6,10,11,12],[13,10,9,8],[14,13,12,11],[14]}} \\ \textbf{EM:} {\scriptsize\texttt{(B0, A0C0C3D0, B1D3E0B2, B3F0E1C1, C2D2G3G2, D1G1H3H2, H1F1E3E2, I0G0F3F2, H0)}} \\ \textbf{Kauffman bracket:} {\scriptsize $A^{21} + 2A^{19} - 2A^{17} - 5A^{15} - A^{13} + 4A^{11} + 3A^{9} - 2A^{7} - 2A^{5} + A$} \\ \textbf{Arrow:} {\scriptsize $-A^{6}L_1 - 2A^{4} + 2A^{2}L_1 + 5 + L_1/A^{2} - 4/A^{4} - 3L_1/A^{6} + 2/A^{8} + 2L_1/A^{10} - L_1/A^{14}$} \\ \textbf{Mock:} {\scriptsize $w^{4} - 3w^{2} + 6 + 2/w - 3/w^{2} - 2/w^{3}$} \\ \textbf{Affine:} {\scriptsize $2t - 4 + 2/t$} \\ \textbf{Yamada:} {\scriptsize $-2A^{25} + 2A^{24} + 3A^{23} - 7A^{22} + 3A^{21} + 8A^{20} - 9A^{19} + 7A^{17} - 4A^{16} - A^{15} + 2A^{14} + 5A^{13} - 3A^{12} - 2A^{11} + 9A^{10} - 3A^{9} - 6A^{8} + 8A^{7} - 7A^{5} + 3A^{4} + 3A^{3} - 3A^{2} - A + 1$}
\end{minipage}

\noindent{\color{gray!40}\rule{\textwidth}{0.4pt}}
\vspace{0.9\baselineskip}
\noindent \begin{minipage}[t]{0.25\textwidth}
\vspace{0pt}
\centering
\includegraphics[page=313,width=\linewidth]{knotoids.pdf}
\end{minipage}
\hfill
\begin{minipage}[t]{0.73\textwidth}
\vspace{0pt}
\raggedright
\textbf{Name:} {\large{$\mathbf{K7_{194}}$}} (chiral, non-rotatable$^{*}$) \\ \textbf{PD:} {\scriptsize\texttt{[0],[0,1,2,3],[1,4,5,2],[3,6,7,4],[5,8,9,10],[6,10,11,7],[12,13,14,8],[13,12,11,9],[14]}} \\ \textbf{EM:} {\scriptsize\texttt{(B0, A0C0C3D0, B1D3E0B2, B3F0F3C1, C2G3H3F1, D1E3H2D2, H1H0I0E1, G1G0F2E2, G2)}} \\ \textbf{Kauffman bracket:} {\scriptsize $A^{21} - 2A^{17} - 2A^{15} + 3A^{13} + 4A^{11} - A^{9} - 4A^{7} - 2A^{5} + A^{3} + A$} \\ \textbf{Arrow:} {\scriptsize $-A^{18}L_1 + 2A^{14}L_1 + 2A^{12} - 3A^{10}L_1 - 4A^{8} + A^{6}L_1 + 4A^{4} + 2A^{2}L_1 - 1 - L_1/A^{2}$} \\ \textbf{Mock:} {\scriptsize $-w^{5} - w^{4} + 2w^{3} + 2w^{2} - 3w - 3 + 1/w + 4/w^{2} + 2/w^{3} - 1/w^{4} - 1/w^{5}$} \\ \textbf{Affine:} {\scriptsize $-2t + 4 - 2/t$} \\ \textbf{Yamada:} {\scriptsize $-A^{25} + 2A^{24} + A^{23} - 5A^{22} + 4A^{21} + 4A^{20} - 8A^{19} + 3A^{18} + 5A^{17} - 5A^{16} - A^{15} + A^{14} + A^{13} - 5A^{12} - A^{11} + 6A^{10} - 5A^{9} - 3A^{8} + 7A^{7} - 3A^{6} - 5A^{5} + 3A^{4} + A^{3} - 4A^{2} + 2$}
\end{minipage}

\noindent{\color{gray!40}\rule{\textwidth}{0.4pt}}
\vspace{0.9\baselineskip}
\noindent \begin{minipage}[t]{0.25\textwidth}
\vspace{0pt}
\centering
\includegraphics[page=314,width=\linewidth]{knotoids.pdf}
\end{minipage}
\hfill
\begin{minipage}[t]{0.73\textwidth}
\vspace{0pt}
\raggedright
\textbf{Name:} {\large{$\mathbf{K7_{195}}$}} (chiral, non-rotatable$^{*}$) \\ \textbf{PD:} {\scriptsize\texttt{[0],[0,1,2,3],[1,4,5,2],[3,6,7,4],[5,8,9,10],[6,10,11,12],[12,11,13,7],[8,13,14,9],[14]}} \\ \textbf{EM:} {\scriptsize\texttt{(B0, A0C0C3D0, B1D3E0B2, B3F0G3C1, C2H0H3F1, D1E3G1G0, F3F2H1D2, E1G2I0E2, H2)}} \\ \textbf{Kauffman bracket:} {\scriptsize $-A^{24} + A^{22} + 4A^{20} - 5A^{16} - 2A^{14} + 5A^{12} + 4A^{10} - 3A^{8} - 3A^{6} + A^{2}$} \\ \textbf{Arrow:} {\scriptsize $-L_1/A^{6} + A^{-8} + 4L_1/A^{10} + L_2/A^{12} - 1/A^{12} - 5L_1/A^{14} - 3L_2/A^{16} + A^{-16} + 5L_1/A^{18} + 4L_2/A^{20} - 3L_1/A^{22} - 3L_2/A^{24} + L_2/A^{28}$} \\ \textbf{Mock:} {\scriptsize $-w^{5} + w^{4} + 6w^{3} + 3w^{2} - 7w - 5 + 3/w + 2/w^{2} - 1/w^{3}$} \\ \textbf{Affine:} {\scriptsize $3t - 6 + 3/t$} \\ \textbf{Yamada:} {\scriptsize $-A^{26} + A^{25} + A^{24} - 7A^{23} + 2A^{22} + 10A^{21} - 13A^{20} - 3A^{19} + 12A^{18} - 12A^{17} - 6A^{16} + 8A^{15} - 2A^{14} - 2A^{13} + 10A^{11} - A^{10} - 9A^{9} + 13A^{8} - 14A^{6} + 8A^{5} + 4A^{4} - 9A^{3} + 2A^{2} + 3A - 1$}
\end{minipage}

\noindent{\color{gray!40}\rule{\textwidth}{0.4pt}}
\vspace{0.9\baselineskip}
\noindent \begin{minipage}[t]{0.25\textwidth}
\vspace{0pt}
\centering
\includegraphics[page=315,width=\linewidth]{knotoids.pdf}
\end{minipage}
\hfill
\begin{minipage}[t]{0.73\textwidth}
\vspace{0pt}
\raggedright
\textbf{Name:} {\large{$\mathbf{K7_{196}}$}} (chiral, non-rotatable$^{*}$) \\ \textbf{PD:} {\scriptsize\texttt{[0],[0,1,2,3],[1,4,5,2],[3,6,7,4],[5,8,9,10],[10,9,11,6],[7,11,12,13],[13,12,14,8],[14]}} \\ \textbf{EM:} {\scriptsize\texttt{(B0, A0C0C3D0, B1D3E0B2, B3F3G0C1, C2H3F1F0, E3E2G1D1, D2F2H1H0, G3G2I0E1, H2)}} \\ \textbf{Kauffman bracket:} {\scriptsize $A^{14} + A^{12} - 2A^{8} - A^{6} + A^{2} + 1$} \\ \textbf{Arrow:} {\scriptsize $A^{8} + A^{6}L_1 - 2A^{2}L_1 - 1 + A^{-4} + L_1/A^{6}$} \\ \textbf{Mock:} {\scriptsize $w^{3} - 2w + w^{-2} + w^{-3}$} \\ \textbf{Affine:} {\scriptsize $-t + 2 - 1/t$} \\ \textbf{Yamada:} {\scriptsize $A^{23} - A^{22} - A^{21} + A^{20} - A^{19} - A^{18} + A^{15} - A^{14} - 2A^{11} - 2A^{8} + A^{2} - 1$}
\end{minipage}

\noindent{\color{gray!40}\rule{\textwidth}{0.4pt}}
\vspace{0.9\baselineskip}
\noindent \begin{minipage}[t]{0.25\textwidth}
\vspace{0pt}
\centering
\includegraphics[page=316,width=\linewidth]{knotoids.pdf}
\end{minipage}
\hfill
\begin{minipage}[t]{0.73\textwidth}
\vspace{0pt}
\raggedright
\textbf{Name:} {\large{$\mathbf{K7_{197}}$}} (chiral, non-rotatable$^{*}$) \\ \textbf{PD:} {\scriptsize\texttt{[0],[0,1,2,3],[1,4,5,2],[3,6,7,4],[5,8,9,10],[10,9,11,6],[11,12,13,7],[8,13,12,14],[14]}} \\ \textbf{EM:} {\scriptsize\texttt{(B0, A0C0C3D0, B1D3E0B2, B3F3G3C1, C2H0F1F0, E3E2G0D1, F2H2H1D2, E1G2G1I0, H3)}} \\ \textbf{Kauffman bracket:} {\scriptsize $A^{16} + 2A^{14} - A^{12} - 3A^{10} - A^{8} + 3A^{6} + 2A^{4} - A^{2} - 1$} \\ \textbf{Arrow:} {\scriptsize $L_1/A^{2} + 2/A^{4} - L_1/A^{6} - 3/A^{8} - L_1/A^{10} + 3/A^{12} + 2L_1/A^{14} - 1/A^{16} - L_1/A^{18}$} \\ \textbf{Mock:} {\scriptsize $-w^{4} - w^{3} + 3w^{2} + 2w - 2 + 2/w^{2} - 1/w^{3} - 1/w^{4}$} \\ \textbf{Affine:} {\scriptsize $t - 2 + 1/t$} \\ \textbf{Yamada:} {\scriptsize $A^{25} - A^{24} - 2A^{23} + 2A^{22} + A^{21} - 5A^{20} - A^{19} + 2A^{18} - 4A^{17} - 2A^{16} + 2A^{15} - A^{14} + 4A^{11} - A^{10} - 2A^{9} + 5A^{8} - 2A^{7} - 3A^{6} + 2A^{5} - 2A^{3} - A^{2} + A + 1$}
\end{minipage}

\noindent{\color{gray!40}\rule{\textwidth}{0.4pt}}
\vspace{0.9\baselineskip}
\noindent \begin{minipage}[t]{0.25\textwidth}
\vspace{0pt}
\centering
\includegraphics[page=317,width=\linewidth]{knotoids.pdf}
\end{minipage}
\hfill
\begin{minipage}[t]{0.73\textwidth}
\vspace{0pt}
\raggedright
\textbf{Name:} {\large{$\mathbf{K7_{198}}$}} (chiral, rotatable) \\ \textbf{PD:} {\scriptsize\texttt{[0],[0,1,2,3],[1,4,5,2],[3,6,7,4],[7,8,9,5],[6,10,11,12],[8,13,10,9],[14,13,12,11],[14]}} \\ \textbf{EM:} {\scriptsize\texttt{(B0, A0C0C3D0, B1D3E3B2, B3F0E0C1, D2G0G3C2, D1G2H3H2, E1H1F1E2, I0G1F3F2, H0)}} \\ \textbf{Kauffman bracket:} {\scriptsize $A^{17} - 2A^{13} - 2A^{11} + A^{9} + 4A^{7} + A^{5} - 2A^{3} - 2A$} \\ \textbf{Arrow:} {\scriptsize $-A^{20} + 2A^{16} + 2A^{14}L_1 - A^{12} - 4A^{10}L_1 - A^{8} + 2A^{6}L_1 + 2A^{4}$} \\ \textbf{Mock:} {\scriptsize $w^{4} + 2w^{3} - w^{2} - 4w + 2/w + w^{-2}$} \\ \textbf{Affine:} {\scriptsize $0$} \\ \textbf{Yamada:} {\scriptsize $A^{22} - A^{21} - 3A^{20} + 4A^{19} + A^{18} - 5A^{17} + 4A^{16} + A^{15} - 3A^{14} + A^{13} + A^{11} - 2A^{10} + 2A^{9} + 4A^{8} - 3A^{7} + A^{6} + 5A^{5} - 3A^{4} - A^{3} + 3A^{2} - 1$}
\end{minipage}

\noindent{\color{gray!40}\rule{\textwidth}{0.4pt}}
\vspace{0.9\baselineskip}
\noindent \begin{minipage}[t]{0.25\textwidth}
\vspace{0pt}
\centering
\includegraphics[page=318,width=\linewidth]{knotoids.pdf}
\end{minipage}
\hfill
\begin{minipage}[t]{0.73\textwidth}
\vspace{0pt}
\raggedright
\textbf{Name:} {\large{$\mathbf{K7_{199}}$}} (chiral, non-rotatable$^{*}$) \\ \textbf{PD:} {\scriptsize\texttt{[0],[0,1,2,3],[1,4,5,2],[3,6,7,4],[8,9,10,5],[6,10,11,7],[8,12,13,14],[9,13,12,11],[14]}} \\ \textbf{EM:} {\scriptsize\texttt{(B0, A0C0C3D0, B1D3E3B2, B3F0F3C1, G0H0F1C2, D1E2H3D2, E0H2H1I0, E1G2G1F2, G3)}} \\ \textbf{Kauffman bracket:} {\scriptsize $A^{24} + A^{22} - A^{20} - 4A^{18} + 6A^{14} + 4A^{12} - 4A^{10} - 5A^{8} + A^{6} + 3A^{4} - 1$} \\ \textbf{Arrow:} {\scriptsize $1 + L_1/A^{2} - 1/A^{4} - 4L_1/A^{6} + 6L_1/A^{10} + 4/A^{12} - 4L_1/A^{14} - 5/A^{16} + L_1/A^{18} + 3/A^{20} - 1/A^{24}$} \\ \textbf{Mock:} {\scriptsize $-w^{6} - 2w^{5} + w^{4} + 5w^{3} + 3w^{2} - 5w - 5 + 3/w + 4/w^{2} - 1/w^{3} - 1/w^{4}$} \\ \textbf{Affine:} {\scriptsize $0$} \\ \textbf{Yamada:} {\scriptsize $A^{27} - 4A^{25} + A^{24} + 7A^{23} - 7A^{22} - 7A^{21} + 13A^{20} - A^{19} - 16A^{18} + 10A^{17} + 6A^{16} - 15A^{15} + 2A^{14} + 8A^{13} - 6A^{12} - 5A^{11} + 5A^{10} + 8A^{9} - 13A^{8} + 18A^{6} - 12A^{5} - 5A^{4} + 12A^{3} - 4A^{2} - 4A + 2$}
\end{minipage}

\noindent{\color{gray!40}\rule{\textwidth}{0.4pt}}
\vspace{0.9\baselineskip}
\noindent \begin{minipage}[t]{0.25\textwidth}
\vspace{0pt}
\centering
\includegraphics[page=319,width=\linewidth]{knotoids.pdf}
\end{minipage}
\hfill
\begin{minipage}[t]{0.73\textwidth}
\vspace{0pt}
\raggedright
\textbf{Name:} {\large{$\mathbf{K7_{200}}$}} (chiral, non-rotatable$^{*}$) \\ \textbf{PD:} {\scriptsize\texttt{[0],[0,1,2,3],[1,4,5,2],[3,6,7,4],[8,9,10,5],[6,10,9,11],[7,11,12,13],[13,12,14,8],[14]}} \\ \textbf{EM:} {\scriptsize\texttt{(B0, A0C0C3D0, B1D3E3B2, B3F0G0C1, H3F2F1C2, D1E2E1G1, D2F3H1H0, G3G2I0E0, H2)}} \\ \textbf{Kauffman bracket:} {\scriptsize $A^{23} - 2A^{19} - 2A^{17} + 3A^{15} + 4A^{13} - A^{11} - 4A^{9} - A^{7} + 2A^{5} - A$} \\ \textbf{Arrow:} {\scriptsize $-A^{26}L_1 + 2A^{22}L_1 + 2A^{20} - 3A^{18}L_1 - 4A^{16} + A^{14}L_1 + 4A^{12} + A^{10}L_1 - 2A^{8} + A^{4}$} \\ \textbf{Mock:} {\scriptsize $w^{3} + 2w^{2} - 4w - 6 + 2/w + 5/w^{2} + w^{-3}$} \\ \textbf{Affine:} {\scriptsize $-3t + 6 - 3/t$} \\ \textbf{Yamada:} {\scriptsize $A^{24} - A^{23} - 2A^{22} + 4A^{21} + A^{20} - 6A^{19} + 2A^{18} + 4A^{17} - 6A^{16} - 2A^{15} + 6A^{14} - 3A^{13} - 2A^{12} + 4A^{11} + A^{10} - A^{9} - 2A^{8} + 7A^{7} - 2A^{6} - 6A^{5} + 7A^{4} - 4A^{2} + 4A + 2$}
\end{minipage}

\noindent{\color{gray!40}\rule{\textwidth}{0.4pt}}
\vspace{0.9\baselineskip}
\noindent \begin{minipage}[t]{0.25\textwidth}
\vspace{0pt}
\centering
\includegraphics[page=320,width=\linewidth]{knotoids.pdf}
\end{minipage}
\hfill
\begin{minipage}[t]{0.73\textwidth}
\vspace{0pt}
\raggedright
\textbf{Name:} {\large{$\mathbf{K7_{201}}$}} (chiral, non-rotatable$^{*}$) \\ \textbf{PD:} {\scriptsize\texttt{[0],[0,1,2,3],[1,4,5,2],[3,6,7,4],[8,9,10,5],[6,10,9,11],[11,12,13,7],[8,13,12,14],[14]}} \\ \textbf{EM:} {\scriptsize\texttt{(B0, A0C0C3D0, B1D3E3B2, B3F0G3C1, H0F2F1C2, D1E2E1G0, F3H2H1D2, E0G2G1I0, H3)}} \\ \textbf{Kauffman bracket:} {\scriptsize $A^{24} + A^{22} - A^{20} - 4A^{18} - A^{16} + 5A^{14} + 4A^{12} - 3A^{10} - 4A^{8} + A^{6} + 3A^{4} - 1$} \\ \textbf{Arrow:} {\scriptsize $A^{12} + A^{10}L_1 - A^{8} - 4A^{6}L_1 - A^{4} + 5A^{2}L_1 + 4 - 3L_1/A^{2} - 4/A^{4} + L_1/A^{6} + 3/A^{8} - 1/A^{12}$} \\ \textbf{Mock:} {\scriptsize $-w^{4} - 3w^{3} - w^{2} + 6w + 6 - 4/w - 4/w^{2} + w^{-3} + w^{-4}$} \\ \textbf{Affine:} {\scriptsize $-t + 2 - 1/t$} \\ \textbf{Yamada:} {\scriptsize $A^{27} - 3A^{25} + A^{24} + 5A^{23} - 5A^{22} - 6A^{21} + 10A^{20} + 2A^{19} - 13A^{18} + 7A^{17} + 8A^{16} - 12A^{15} + A^{14} + 7A^{13} - 5A^{12} - 5A^{11} + 2A^{10} + 6A^{9} - 11A^{8} - 3A^{7} + 14A^{6} - 8A^{5} - 6A^{4} + 11A^{3} - A^{2} - 4A + 1$}
\end{minipage}

\noindent{\color{gray!40}\rule{\textwidth}{0.4pt}}
\vspace{0.9\baselineskip}
\noindent \begin{minipage}[t]{0.25\textwidth}
\vspace{0pt}
\centering
\includegraphics[page=321,width=\linewidth]{knotoids.pdf}
\end{minipage}
\hfill
\begin{minipage}[t]{0.73\textwidth}
\vspace{0pt}
\raggedright
\textbf{Name:} {\large{$\mathbf{K7_{202}}$}} (chiral, non-rotatable$^{*}$) \\ \textbf{PD:} {\scriptsize\texttt{[0],[0,1,2,3],[1,4,5,2],[3,6,7,4],[8,9,10,5],[6,10,11,12],[12,11,13,7],[13,14,9,8],[14]}} \\ \textbf{EM:} {\scriptsize\texttt{(B0, A0C0C3D0, B1D3E3B2, B3F0G3C1, H3H2F1C2, D1E2G1G0, F3F2H0D2, G2I0E1E0, H1)}} \\ \textbf{Kauffman bracket:} {\scriptsize $-A^{25} + 3A^{21} + A^{19} - 4A^{17} - 4A^{15} + 3A^{13} + 5A^{11} - A^{9} - 4A^{7} - A^{5} + A^{3} + A$} \\ \textbf{Arrow:} {\scriptsize $A^{10}L_1 - 3A^{6}L_1 - A^{4}L_2 + 4A^{2}L_1 + 3L_2 + 1 - 3L_1/A^{2} - 4L_2/A^{4} - 1/A^{4} + L_1/A^{6} + 3L_2/A^{8} + A^{-8} + L_1/A^{10} - L_2/A^{12} - L_1/A^{14}$} \\ \textbf{Mock:} {\scriptsize $w^{5} + w^{4} - 4w^{3} - 4w^{2} + 5w + 7 - 1/w - 3/w^{2} - 1/w^{3}$} \\ \textbf{Affine:} {\scriptsize $t - 2 + 1/t$} \\ \textbf{Yamada:} {\scriptsize $A^{27} - A^{26} - 3A^{25} + 6A^{24} + 3A^{23} - 13A^{22} + 4A^{21} + 10A^{20} - 14A^{19} + 10A^{17} - 5A^{16} - 3A^{15} + 2A^{14} + 5A^{13} - 9A^{12} - 5A^{11} + 12A^{10} - 7A^{9} - 9A^{8} + 13A^{7} + A^{6} - 9A^{5} + 4A^{4} + 4A^{3} - 3A^{2} - A + 1$}
\end{minipage}

\noindent{\color{gray!40}\rule{\textwidth}{0.4pt}}
\vspace{0.9\baselineskip}
\noindent \begin{minipage}[t]{0.25\textwidth}
\vspace{0pt}
\centering
\includegraphics[page=322,width=\linewidth]{knotoids.pdf}
\end{minipage}
\hfill
\begin{minipage}[t]{0.73\textwidth}
\vspace{0pt}
\raggedright
\textbf{Name:} {\large{$\mathbf{K7_{203}}$}} (chiral, non-rotatable$^{*}$) \\ \textbf{PD:} {\scriptsize\texttt{[0],[0,1,2,3],[1,4,5,2],[6,7,4,3],[5,8,9,10],[10,11,7,6],[12,13,14,8],[13,12,11,9],[14]}} \\ \textbf{EM:} {\scriptsize\texttt{(B0, A0C0C3D3, B1D2E0B2, F3F2C1B3, C2G3H3F0, E3H2D1D0, H1H0I0E1, G1G0F1E2, G2)}} \\ \textbf{Kauffman bracket:} {\scriptsize $A^{22} + A^{20} - 2A^{18} - 3A^{16} + 2A^{14} + 6A^{12} - 5A^{8} - 3A^{6} + 2A^{4} + 2A^{2}$} \\ \textbf{Arrow:} {\scriptsize $A^{10}L_1 + A^{8} - 2A^{6}L_1 - 3A^{4} + 2A^{2}L_1 + 6 - 5/A^{4} - 3L_1/A^{6} + 2/A^{8} + 2L_1/A^{10}$} \\ \textbf{Mock:} {\scriptsize $w^{5} + w^{4} - 2w^{3} - 4w^{2} + w + 7 + 1/w - 4/w^{2} - 2/w^{3} + w^{-4} + w^{-5}$} \\ \textbf{Affine:} {\scriptsize $0$} \\ \textbf{Yamada:} {\scriptsize $-2A^{28} + A^{27} + 4A^{26} - 5A^{25} - 3A^{24} + 11A^{23} - 4A^{22} - 11A^{21} + 11A^{20} + A^{19} - 10A^{18} + 3A^{17} + 3A^{16} - 2A^{15} - 7A^{14} + 7A^{13} + 5A^{12} - 13A^{11} + 5A^{10} + 8A^{9} - 11A^{8} - 2A^{7} + 8A^{6} - 2A^{5} - 4A^{4} + 3A^{3} + 2A^{2} - A - 1$}
\end{minipage}

\noindent{\color{gray!40}\rule{\textwidth}{0.4pt}}
\vspace{0.9\baselineskip}
\noindent \begin{minipage}[t]{0.25\textwidth}
\vspace{0pt}
\centering
\includegraphics[page=323,width=\linewidth]{knotoids.pdf}
\end{minipage}
\hfill
\begin{minipage}[t]{0.73\textwidth}
\vspace{0pt}
\raggedright
\textbf{Name:} {\large{$\mathbf{K7_{204}}$}} (chiral, non-rotatable$^{*}$) \\ \textbf{PD:} {\scriptsize\texttt{[0],[0,1,2,3],[1,4,5,2],[6,7,4,3],[5,8,9,10],[10,11,12,6],[7,12,11,13],[8,13,14,9],[14]}} \\ \textbf{EM:} {\scriptsize\texttt{(B0, A0C0C3D3, B1D2E0B2, F3G0C1B3, C2H0H3F0, E3G2G1D0, D1F2F1H1, E1G3I0E2, H2)}} \\ \textbf{Kauffman bracket:} {\scriptsize $A^{24} + 2A^{22} - 3A^{18} - 2A^{16} + 2A^{14} + 4A^{12} - A^{10} - 3A^{8} - A^{6} + A^{4} + A^{2}$} \\ \textbf{Arrow:} {\scriptsize $1 + 2L_1/A^{2} + L_2/A^{4} - 1/A^{4} - 3L_1/A^{6} - 3L_2/A^{8} + A^{-8} + 2L_1/A^{10} + 4L_2/A^{12} - L_1/A^{14} - 3L_2/A^{16} - L_1/A^{18} + L_2/A^{20} + L_1/A^{22}$} \\ \textbf{Mock:} {\scriptsize $-w^{6} - w^{5} + 2w^{4} + 4w^{3} + 2w^{2} - 3w - 5 - 1/w + 4/w^{2} + w^{-3} - 1/w^{4}$} \\ \textbf{Affine:} {\scriptsize $t - 2 + 1/t$} \\ \textbf{Yamada:} {\scriptsize $-A^{28} + 2A^{26} - 2A^{25} - 4A^{24} + 4A^{23} + A^{22} - 7A^{21} + 7A^{19} - 3A^{18} - 3A^{17} + 7A^{16} - A^{15} - 6A^{14} + A^{13} - 4A^{11} - 4A^{10} + 7A^{9} + A^{8} - 7A^{7} + 7A^{6} + 3A^{5} - 4A^{4} + A^{2} - 1$}
\end{minipage}

\noindent{\color{gray!40}\rule{\textwidth}{0.4pt}}
\vspace{0.9\baselineskip}
\noindent \begin{minipage}[t]{0.25\textwidth}
\vspace{0pt}
\centering
\includegraphics[page=324,width=\linewidth]{knotoids.pdf}
\end{minipage}
\hfill
\begin{minipage}[t]{0.73\textwidth}
\vspace{0pt}
\raggedright
\textbf{Name:} {\large{$\mathbf{K7_{205}}$}} (chiral, non-rotatable$^{*}$) \\ \textbf{PD:} {\scriptsize\texttt{[0],[0,1,2,3],[1,4,5,2],[6,7,4,3],[5,8,9,10],[10,9,11,6],[7,11,12,13],[13,12,14,8],[14]}} \\ \textbf{EM:} {\scriptsize\texttt{(B0, A0C0C3D3, B1D2E0B2, F3G0C1B3, C2H3F1F0, E3E2G1D0, D1F2H1H0, G3G2I0E1, H2)}} \\ \textbf{Kauffman bracket:} {\scriptsize $A^{22} + A^{20} - 2A^{18} - 3A^{16} + 2A^{14} + 5A^{12} - A^{10} - 5A^{8} - 2A^{6} + 3A^{4} + 2A^{2}$} \\ \textbf{Arrow:} {\scriptsize $A^{22}L_1 + A^{20} - 2A^{18}L_1 - 3A^{16} + 2A^{14}L_1 + 5A^{12} - A^{10}L_1 - 5A^{8} - 2A^{6}L_1 + 3A^{4} + 2A^{2}L_1$} \\ \textbf{Mock:} {\scriptsize $-w^{4} + w^{3} + 4w^{2} - 2w - 6 + 5/w^{2} + w^{-3} - 1/w^{4}$} \\ \textbf{Affine:} {\scriptsize $-t + 2 - 1/t$} \\ \textbf{Yamada:} {\scriptsize $-A^{28} + 2A^{27} + 2A^{26} - 6A^{25} - A^{24} + 10A^{23} - 6A^{22} - 8A^{21} + 12A^{20} - 2A^{19} - 7A^{18} + 5A^{17} + 3A^{16} - 2A^{15} - 5A^{14} + 7A^{13} + A^{12} - 12A^{11} + 6A^{10} + 5A^{9} - 12A^{8} + A^{7} + 7A^{6} - 3A^{5} - 3A^{4} + 2A^{3} + A^{2} - A - 1$}
\end{minipage}

\noindent{\color{gray!40}\rule{\textwidth}{0.4pt}}
\vspace{0.9\baselineskip}
\noindent \begin{minipage}[t]{0.25\textwidth}
\vspace{0pt}
\centering
\includegraphics[page=325,width=\linewidth]{knotoids.pdf}
\end{minipage}
\hfill
\begin{minipage}[t]{0.73\textwidth}
\vspace{0pt}
\raggedright
\textbf{Name:} {\large{$\mathbf{K7_{206}}$}} (chiral, non-rotatable$^{*}$) \\ \textbf{PD:} {\scriptsize\texttt{[0],[0,1,2,3],[1,4,5,2],[3,6,7,4],[5,8,9,6],[7,10,11,8],[12,13,14,9],[10,14,13,11],[12]}} \\ \textbf{EM:} {\scriptsize\texttt{(B0, A0C0C3D0, B1D3E0B2, B3E3F0C1, C2F3G3D1, D2H0H3E1, I0H2H1E2, F1G2G1F2, G0)}} \\ \textbf{Kauffman bracket:} {\scriptsize $A^{22} + A^{20} - A^{18} - 3A^{16} + A^{14} + 4A^{12} + A^{10} - 3A^{8} - 2A^{6} + A^{4} + A^{2}$} \\ \textbf{Arrow:} {\scriptsize $A^{4} + A^{2}L_1 - 1 - 3L_1/A^{2} - L_2/A^{4} + 2/A^{4} + 4L_1/A^{6} + 2L_2/A^{8} - 1/A^{8} - 3L_1/A^{10} - 2L_2/A^{12} + L_1/A^{14} + L_2/A^{16}$} \\ \textbf{Mock:} {\scriptsize $w^{5} - 3w^{3} + 4w + 1 - 3/w + w^{-3}$} \\ \textbf{Affine:} {\scriptsize $-t^{2} + 2 - 1/t^{2}$} \\ \textbf{Yamada:} {\scriptsize $-A^{26} + A^{24} - 2A^{23} + A^{22} + 5A^{21} - 4A^{20} + 8A^{18} - 3A^{17} - A^{16} + 5A^{15} - 2A^{14} - A^{13} - 2A^{12} + 3A^{11} - A^{10} - 5A^{9} + 6A^{8} - 5A^{6} + 4A^{5} + 2A^{4} - 4A^{3} + A^{2} + 2A - 1$}
\end{minipage}

\noindent{\color{gray!40}\rule{\textwidth}{0.4pt}}
\vspace{0.9\baselineskip}
\noindent \begin{minipage}[t]{0.25\textwidth}
\vspace{0pt}
\centering
\includegraphics[page=326,width=\linewidth]{knotoids.pdf}
\end{minipage}
\hfill
\begin{minipage}[t]{0.73\textwidth}
\vspace{0pt}
\raggedright
\textbf{Name:} {\large{$\mathbf{K7_{207}}$}} (chiral, non-rotatable$^{*}$) \\ \textbf{PD:} {\scriptsize\texttt{[0],[0,1,2,3],[1,4,5,2],[3,6,7,4],[5,8,9,6],[10,11,8,7],[9,12,13,14],[14,13,11,10],[12]}} \\ \textbf{EM:} {\scriptsize\texttt{(B0, A0C0C3D0, B1D3E0B2, B3E3F3C1, C2F2G0D1, H3H2E1D2, E2I0H1H0, G3G2F1F0, G1)}} \\ \textbf{Kauffman bracket:} {\scriptsize $A^{20} - A^{18} - 3A^{16} + 4A^{12} + 2A^{10} - 3A^{8} - 2A^{6} + A^{4} + 2A^{2}$} \\ \textbf{Arrow:} {\scriptsize $A^{26}L_1 - A^{24} - 3A^{22}L_1 - A^{20}L_2 + A^{20} + 4A^{18}L_1 + 2A^{16}L_2 - 3A^{14}L_1 - 2A^{12}L_2 + A^{10}L_1 + A^{8}L_2 + A^{8}$} \\ \textbf{Mock:} {\scriptsize $w^{5} + w^{4} - 3w^{3} - w^{2} + 4w + 1 - 3/w - 1/w^{2} + w^{-3} + w^{-4}$} \\ \textbf{Affine:} {\scriptsize $-t^{2} + 2 - 1/t^{2}$} \\ \textbf{Yamada:} {\scriptsize $-A^{27} + A^{26} + A^{25} - 4A^{24} + 3A^{23} + 3A^{22} - 7A^{21} + 5A^{20} + 4A^{19} - 4A^{18} + A^{17} + A^{16} + A^{15} - 4A^{14} + 4A^{12} - 5A^{11} - 3A^{10} + 5A^{9} - 4A^{8} - 4A^{7} + 4A^{6} - 3A^{4} + A^{3} + A^{2} - A - 1$}
\end{minipage}

\noindent{\color{gray!40}\rule{\textwidth}{0.4pt}}
\vspace{0.9\baselineskip}
\noindent \begin{minipage}[t]{0.25\textwidth}
\vspace{0pt}
\centering
\includegraphics[page=327,width=\linewidth]{knotoids.pdf}
\end{minipage}
\hfill
\begin{minipage}[t]{0.73\textwidth}
\vspace{0pt}
\raggedright
\textbf{Name:} {\large{$\mathbf{K7_{208}}$}} (chiral, non-rotatable$^{*}$) \\ \textbf{PD:} {\scriptsize\texttt{[0],[0,1,2,3],[1,4,5,2],[3,6,7,4],[5,8,9,10],[6,10,11,12],[13,14,8,7],[9],[14,13,12,11]}} \\ \textbf{EM:} {\scriptsize\texttt{(B0, A0C0C3D0, B1D3E0B2, B3F0G3C1, C2G2H0F1, D1E3I3I2, I1I0E1D2, E2, G1G0F3F2)}} \\ \textbf{Kauffman bracket:} {\scriptsize $-A^{23} + A^{21} + 3A^{19} - 4A^{15} - A^{13} + 3A^{11} + A^{9} - 3A^{7} - 2A^{5} + A^{3} + A$} \\ \textbf{Arrow:} {\scriptsize $A^{20} - A^{18}L_1 - 3A^{16} + 4A^{12} + A^{10}L_1 - 3A^{8} - A^{6}L_1 + 3A^{4} + 2A^{2}L_1 - 1 - L_1/A^{2}$} \\ \textbf{Mock:} {\scriptsize $-w^{5} - 2w^{4} + 4w^{2} + w - 4 - 1/w + 4/w^{2} + 2/w^{3} - 1/w^{4} - 1/w^{5}$} \\ \textbf{Affine:} {\scriptsize $-2t + 4 - 2/t$} \\ \textbf{Yamada:} {\scriptsize $-A^{27} + A^{26} - 3A^{24} + 5A^{23} + A^{22} - 6A^{21} + 7A^{20} - 5A^{18} + 4A^{17} - A^{15} - 3A^{14} + A^{13} + A^{12} - 6A^{11} + 4A^{9} - 6A^{8} - A^{7} + 5A^{6} - 3A^{5} - 2A^{4} + 3A^{3} - A^{2} - A + 1$}
\end{minipage}

\noindent{\color{gray!40}\rule{\textwidth}{0.4pt}}
\vspace{0.9\baselineskip}
\noindent \begin{minipage}[t]{0.25\textwidth}
\vspace{0pt}
\centering
\includegraphics[page=328,width=\linewidth]{knotoids.pdf}
\end{minipage}
\hfill
\begin{minipage}[t]{0.73\textwidth}
\vspace{0pt}
\raggedright
\textbf{Name:} {\large{$\mathbf{K7_{209}}$}} (chiral, non-rotatable$^{*}$) \\ \textbf{PD:} {\scriptsize\texttt{[0],[0,1,2,3],[1,4,5,2],[3,6,7,4],[8,9,6,5],[7,10,11,8],[12,13,14,9],[10,14,13,11],[12]}} \\ \textbf{EM:} {\scriptsize\texttt{(B0, A0C0C3D0, B1D3E3B2, B3E2F0C1, F3G3D1C2, D2H0H3E0, I0H2H1E1, F1G2G1F2, G0)}} \\ \textbf{Kauffman bracket:} {\scriptsize $-A^{21} - A^{19} + A^{17} + 3A^{15} - 4A^{11} - 3A^{9} + 2A^{7} + 3A^{5} - A$} \\ \textbf{Arrow:} {\scriptsize $1 + L_1/A^{2} - 1/A^{4} - 3L_1/A^{6} - L_2/A^{8} + A^{-8} + 4L_1/A^{10} + 2L_2/A^{12} + A^{-12} - 2L_1/A^{14} - 2L_2/A^{16} - 1/A^{16} + L_2/A^{20}$} \\ \textbf{Mock:} {\scriptsize $-w^{6} - w^{5} + 2w^{4} + 3w^{3} + w^{2} - 3w - 3 + 2/w + 3/w^{2} - 1/w^{3} - 1/w^{4}$} \\ \textbf{Affine:} {\scriptsize $-t^{2} + t + 1/t - 1/t^{2}$} \\ \textbf{Yamada:} {\scriptsize $-A^{28} + 2A^{26} - A^{25} - 3A^{24} + 3A^{23} + 4A^{22} - 3A^{21} + 8A^{19} - A^{18} - 4A^{17} + 6A^{16} + A^{15} - 5A^{14} + 2A^{13} + A^{12} - 2A^{11} - 5A^{10} + 4A^{9} + 2A^{8} - 8A^{7} + 4A^{6} + 4A^{5} - 3A^{4} + 2A^{2} - 1$}
\end{minipage}

\noindent{\color{gray!40}\rule{\textwidth}{0.4pt}}
\vspace{0.9\baselineskip}
\noindent \begin{minipage}[t]{0.25\textwidth}
\vspace{0pt}
\centering
\includegraphics[page=329,width=\linewidth]{knotoids.pdf}
\end{minipage}
\hfill
\begin{minipage}[t]{0.73\textwidth}
\vspace{0pt}
\raggedright
\textbf{Name:} {\large{$\mathbf{K7_{210}}$}} (chiral, non-rotatable$^{*}$) \\ \textbf{PD:} {\scriptsize\texttt{[0],[0,1,2,3],[1,4,5,2],[3,6,7,4],[8,9,6,5],[10,11,8,7],[9,12,13,14],[14,13,11,10],[12]}} \\ \textbf{EM:} {\scriptsize\texttt{(B0, A0C0C3D0, B1D3E3B2, B3E2F3C1, F2G0D1C2, H3H2E0D2, E1I0H1H0, G3G2F1F0, G1)}} \\ \textbf{Kauffman bracket:} {\scriptsize $-A^{25} + 2A^{21} + 2A^{19} - 2A^{17} - 4A^{15} + A^{13} + 4A^{11} + A^{9} - 3A^{7} - 2A^{5} + A$} \\ \textbf{Arrow:} {\scriptsize $A^{22}L_1 - 2A^{18}L_1 - A^{16}L_2 - A^{16} + 2A^{14}L_1 + 2A^{12}L_2 + 2A^{12} - A^{10}L_1 - 2A^{8}L_2 - 2A^{8} - A^{6}L_1 + A^{4}L_2 + 2A^{4} + 2A^{2}L_1 - L_1/A^{2}$} \\ \textbf{Mock:} {\scriptsize $w^{3} + 2w^{2} - 2w - 4 + 1/w + 5/w^{2} + w^{-3} - 2/w^{4} - 1/w^{5}$} \\ \textbf{Affine:} {\scriptsize $-t^{2} + t + 1/t - 1/t^{2}$} \\ \textbf{Yamada:} {\scriptsize $-A^{26} + 4A^{24} - A^{23} - 7A^{22} + 6A^{21} + 4A^{20} - 11A^{19} + 4A^{18} + 5A^{17} - 6A^{16} + 3A^{14} + 3A^{13} - 5A^{12} + 2A^{11} + 8A^{10} - 7A^{9} - 2A^{8} + 10A^{7} - 3A^{6} - 5A^{5} + 5A^{4} + 2A^{3} - 3A^{2} + 1$}
\end{minipage}

\noindent{\color{gray!40}\rule{\textwidth}{0.4pt}}
\vspace{0.9\baselineskip}
\noindent \begin{minipage}[t]{0.25\textwidth}
\vspace{0pt}
\centering
\includegraphics[page=330,width=\linewidth]{knotoids.pdf}
\end{minipage}
\hfill
\begin{minipage}[t]{0.73\textwidth}
\vspace{0pt}
\raggedright
\textbf{Name:} {\large{$\mathbf{K7_{211}}$}} (chiral, non-rotatable$^{*}$) \\ \textbf{PD:} {\scriptsize\texttt{[0],[0,1,2,3],[1,4,5,2],[3,6,7,4],[8,9,10,5],[6,10,11,12],[13,14,8,7],[9],[14,13,12,11]}} \\ \textbf{EM:} {\scriptsize\texttt{(B0, A0C0C3D0, B1D3E3B2, B3F0G3C1, G2H0F1C2, D1E2I3I2, I1I0E0D2, E1, G1G0F3F2)}} \\ \textbf{Kauffman bracket:} {\scriptsize $-A^{25} + 2A^{21} + 2A^{19} - 2A^{17} - 3A^{15} + A^{13} + 3A^{11} - A^{9} - 3A^{7} - A^{5} + A^{3} + A$} \\ \textbf{Arrow:} {\scriptsize $A^{16} - 2A^{12} - 2A^{10}L_1 + 2A^{8} + 3A^{6}L_1 - A^{4} - 3A^{2}L_1 + 1 + 3L_1/A^{2} + A^{-4} - L_1/A^{6} - 1/A^{8}$} \\ \textbf{Mock:} {\scriptsize $-w^{6} - 2w^{5} + w^{4} + 3w^{3} - w^{2} - 3w + 2 + 3/w - 1/w^{3}$} \\ \textbf{Affine:} {\scriptsize $-2t + 4 - 2/t$} \\ \textbf{Yamada:} {\scriptsize $-A^{25} + A^{24} + 3A^{23} - 4A^{22} + 7A^{20} - 6A^{19} - 2A^{18} + 5A^{17} - 3A^{16} - 2A^{15} - A^{14} + 2A^{13} - 4A^{12} - 4A^{11} + 5A^{10} - A^{9} - 5A^{8} + 5A^{7} + 2A^{6} - 5A^{5} + A^{4} + 3A^{3} - 2A^{2} - A + 1$}
\end{minipage}

\noindent{\color{gray!40}\rule{\textwidth}{0.4pt}}
\vspace{0.9\baselineskip}
\noindent \begin{minipage}[t]{0.25\textwidth}
\vspace{0pt}
\centering
\includegraphics[page=331,width=\linewidth]{knotoids.pdf}
\end{minipage}
\hfill
\begin{minipage}[t]{0.73\textwidth}
\vspace{0pt}
\raggedright
\textbf{Name:} {\large{$\mathbf{K7_{212}}$}} (chiral, non-rotatable$^{*}$) \\ \textbf{PD:} {\scriptsize\texttt{[0],[0,1,2,3],[1,4,5,2],[3,6,7,8],[4,8,9,10],[11,12,13,5],[6,13,12,7],[9,11,10,14],[14]}} \\ \textbf{EM:} {\scriptsize\texttt{(B0, A0C0C3D0, B1E0F3B2, B3G0G3E1, C1D3H0H2, H1G2G1C2, D1F2F1D2, E2F0E3I0, H3)}} \\ \textbf{Kauffman bracket:} {\scriptsize $A^{22} - 2A^{18} - A^{16} + 4A^{14} + 3A^{12} - 3A^{10} - 4A^{8} + A^{6} + 3A^{4} - 1$} \\ \textbf{Arrow:} {\scriptsize $A^{4} - L_2 - 1 - L_1/A^{2} + 2L_2/A^{4} + 2/A^{4} + 3L_1/A^{6} - 2L_2/A^{8} - 1/A^{8} - 4L_1/A^{10} + L_2/A^{12} + 3L_1/A^{14} - L_1/A^{18}$} \\ \textbf{Mock:} {\scriptsize $-w^{4} - w^{3} + 3w^{2} + 3w - 3 - 4/w + 3/w^{2} + 3/w^{3} - 1/w^{4} - 1/w^{5}$} \\ \textbf{Affine:} {\scriptsize $-t^{2} + 2 - 1/t^{2}$} \\ \textbf{Yamada:} {\scriptsize $A^{27} + A^{26} - 2A^{25} + 4A^{23} - 3A^{22} - 5A^{21} + 5A^{20} - 9A^{18} + 3A^{17} + 3A^{16} - 7A^{15} + A^{14} + 5A^{13} - 2A^{12} - 2A^{11} + A^{10} + 5A^{9} - 6A^{8} - 2A^{7} + 9A^{6} - 6A^{5} - 3A^{4} + 5A^{3} - 2A^{2} - A + 1$}
\end{minipage}

\noindent{\color{gray!40}\rule{\textwidth}{0.4pt}}
\vspace{0.9\baselineskip}
\noindent \begin{minipage}[t]{0.25\textwidth}
\vspace{0pt}
\centering
\includegraphics[page=332,width=\linewidth]{knotoids.pdf}
\end{minipage}
\hfill
\begin{minipage}[t]{0.73\textwidth}
\vspace{0pt}
\raggedright
\textbf{Name:} {\large{$\mathbf{K7_{213}}$}} (chiral, non-rotatable$^{*}$) \\ \textbf{PD:} {\scriptsize\texttt{[0],[0,1,2,3],[1,4,5,2],[3,6,7,8],[8,7,9,4],[5,10,11,6],[9,12,13,10],[14,13,12,11],[14]}} \\ \textbf{EM:} {\scriptsize\texttt{(B0, A0C0C3D0, B1E3F0B2, B3F3E1E0, D3D2G0C1, C2G3H3D1, E2H2H1F1, I0G2G1F2, H0)}} \\ \textbf{Kauffman bracket:} {\scriptsize $-A^{17} - A^{15} + 2A^{13} + 3A^{11} - A^{9} - 3A^{7} - 2A^{5} + A^{3} + A$} \\ \textbf{Arrow:} {\scriptsize $A^{8} + A^{6}L_1 - 2A^{4} - 3A^{2}L_1 - L_2 + 2 + 3L_1/A^{2} + 2L_2/A^{4} - L_1/A^{6} - L_2/A^{8}$} \\ \textbf{Mock:} {\scriptsize $-w^{5} - w^{4} + 2w^{3} + w^{2} - 2w + 2 + 2/w - 1/w^{2} - 1/w^{3}$} \\ \textbf{Affine:} {\scriptsize $0$} \\ \textbf{Yamada:} {\scriptsize $-A^{22} + 2A^{21} + 2A^{20} - 4A^{19} + A^{18} + 2A^{17} - 4A^{16} - A^{15} - 3A^{12} + 3A^{10} - 2A^{9} + 3A^{7} - 2A^{6} - 3A^{5} + A^{4} + A^{3} - 3A^{2} + 2$}
\end{minipage}

\noindent{\color{gray!40}\rule{\textwidth}{0.4pt}}
\vspace{0.9\baselineskip}
\noindent \begin{minipage}[t]{0.25\textwidth}
\vspace{0pt}
\centering
\includegraphics[page=333,width=\linewidth]{knotoids.pdf}
\end{minipage}
\hfill
\begin{minipage}[t]{0.73\textwidth}
\vspace{0pt}
\raggedright
\textbf{Name:} {\large{$\mathbf{K7_{214}}$}} (chiral, non-rotatable$^{*}$) \\ \textbf{PD:} {\scriptsize\texttt{[0],[0,1,2,3],[1,4,5,2],[3,6,7,8],[8,7,9,4],[5,10,11,6],[12,13,10,9],[11,14,13,12],[14]}} \\ \textbf{EM:} {\scriptsize\texttt{(B0, A0C0C3D0, B1E3F0B2, B3F3E1E0, D3D2G3C1, C2G2H0D1, H3H2F1E2, F2I0G1G0, H1)}} \\ \textbf{Kauffman bracket:} {\scriptsize $A^{20} - 2A^{16} + A^{14} + 4A^{12} + A^{10} - 4A^{8} - 3A^{6} + A^{4} + 2A^{2}$} \\ \textbf{Arrow:} {\scriptsize $A^{2}L_1 - 2L_1/A^{2} + A^{-4} + 4L_1/A^{6} + L_2/A^{8} - 4L_1/A^{10} - 2L_2/A^{12} - 1/A^{12} + L_1/A^{14} + L_2/A^{16} + A^{-16}$} \\ \textbf{Mock:} {\scriptsize $w^{5} + w^{4} - 2w^{3} - w^{2} + 4w + 2 - 4/w - 1/w^{2} + w^{-3}$} \\ \textbf{Affine:} {\scriptsize $2t - 4 + 2/t$} \\ \textbf{Yamada:} {\scriptsize $-A^{25} - 3A^{24} + 3A^{23} - 7A^{21} + 4A^{20} + 3A^{19} - 5A^{18} + 2A^{17} + 2A^{16} - 4A^{14} + A^{13} + 3A^{12} - 6A^{11} + 4A^{9} - 4A^{8} - 3A^{7} + 5A^{6} + A^{5} - 3A^{4} + 2A^{3} + 2A^{2} - A - 1$}
\end{minipage}

\noindent{\color{gray!40}\rule{\textwidth}{0.4pt}}
\vspace{0.9\baselineskip}
\noindent \begin{minipage}[t]{0.25\textwidth}
\vspace{0pt}
\centering
\includegraphics[page=334,width=\linewidth]{knotoids.pdf}
\end{minipage}
\hfill
\begin{minipage}[t]{0.73\textwidth}
\vspace{0pt}
\raggedright
\textbf{Name:} {\large{$\mathbf{K7_{215}}$}} (chiral, non-rotatable$^{*}$) \\ \textbf{PD:} {\scriptsize\texttt{[0],[0,1,2,3],[1,4,5,2],[3,6,7,8],[8,7,9,4],[10,11,6,5],[9,12,13,10],[14,13,12,11],[14]}} \\ \textbf{EM:} {\scriptsize\texttt{(B0, A0C0C3D0, B1E3F3B2, B3F2E1E0, D3D2G0C1, G3H3D1C2, E2H2H1F0, I0G2G1F1, H0)}} \\ \textbf{Kauffman bracket:} {\scriptsize $-A^{20} - A^{18} + 2A^{16} + 3A^{14} - 2A^{12} - 5A^{10} + 5A^{6} + 3A^{4} - 2A^{2} - 1$} \\ \textbf{Arrow:} {\scriptsize $-A^{20} - A^{18}L_1 + 2A^{16} + 3A^{14}L_1 - 2A^{12} - 5A^{10}L_1 - A^{8}L_2 + A^{8} + 5A^{6}L_1 + 2A^{4}L_2 + A^{4} - 2A^{2}L_1 - L_2$} \\ \textbf{Mock:} {\scriptsize $-w^{5} + 3w^{3} - 5w + 1 + 5/w - 2/w^{3}$} \\ \textbf{Affine:} {\scriptsize $0$} \\ \textbf{Yamada:} {\scriptsize $-2A^{27} + 5A^{25} - 2A^{24} - 6A^{23} + 8A^{22} + 4A^{21} - 10A^{20} + 3A^{19} + 8A^{18} - 8A^{17} - 2A^{16} + 5A^{15} - 3A^{14} - 4A^{13} + 8A^{11} - 7A^{10} - 5A^{9} + 11A^{8} - 6A^{7} - 8A^{6} + 6A^{5} - 3A^{3} + A + 1$}
\end{minipage}

\noindent{\color{gray!40}\rule{\textwidth}{0.4pt}}
\vspace{0.9\baselineskip}
\noindent \begin{minipage}[t]{0.25\textwidth}
\vspace{0pt}
\centering
\includegraphics[page=335,width=\linewidth]{knotoids.pdf}
\end{minipage}
\hfill
\begin{minipage}[t]{0.73\textwidth}
\vspace{0pt}
\raggedright
\textbf{Name:} {\large{$\mathbf{K7_{216}}$}} (chiral, non-rotatable$^{*}$) \\ \textbf{PD:} {\scriptsize\texttt{[0],[0,1,2,3],[1,4,5,2],[3,6,7,8],[8,7,9,4],[10,11,6,5],[12,13,10,9],[11,14,13,12],[14]}} \\ \textbf{EM:} {\scriptsize\texttt{(B0, A0C0C3D0, B1E3F3B2, B3F2E1E0, D3D2G3C1, G2H0D1C2, H3H2F0E2, F1I0G1G0, H1)}} \\ \textbf{Kauffman bracket:} {\scriptsize $-A^{25} + 2A^{21} + A^{19} - 3A^{17} - 4A^{15} + 2A^{13} + 5A^{11} + A^{9} - 3A^{7} - 2A^{5} + A$} \\ \textbf{Arrow:} {\scriptsize $A^{10}L_1 - 2A^{6}L_1 - A^{4} + 3A^{2}L_1 + L_2 + 3 - 2L_1/A^{2} - 2L_2/A^{4} - 3/A^{4} - L_1/A^{6} + L_2/A^{8} + 2/A^{8} + 2L_1/A^{10} - L_1/A^{14}$} \\ \textbf{Mock:} {\scriptsize $w^{5} + 2w^{4} - w^{3} - 4w^{2} + w + 5 - 1/w - 2/w^{2}$} \\ \textbf{Affine:} {\scriptsize $2t - 4 + 2/t$} \\ \textbf{Yamada:} {\scriptsize $-2A^{26} + A^{25} + 6A^{24} - 3A^{23} - 6A^{22} + 10A^{21} + A^{20} - 12A^{19} + 6A^{18} + 4A^{17} - 7A^{16} + A^{15} + 4A^{14} + 2A^{13} - 6A^{12} + 5A^{11} + 9A^{10} - 10A^{9} + A^{8} + 10A^{7} - 7A^{6} - 5A^{5} + 5A^{4} + A^{3} - 3A^{2} + 1$}
\end{minipage}

\noindent{\color{gray!40}\rule{\textwidth}{0.4pt}}
\vspace{0.9\baselineskip}
\noindent \begin{minipage}[t]{0.25\textwidth}
\vspace{0pt}
\centering
\includegraphics[page=336,width=\linewidth]{knotoids.pdf}
\end{minipage}
\hfill
\begin{minipage}[t]{0.73\textwidth}
\vspace{0pt}
\raggedright
\textbf{Name:} {\large{$\mathbf{K7_{217}}$}} (chiral, non-rotatable$^{*}$) \\ \textbf{PD:} {\scriptsize\texttt{[0],[0,1,2,3],[1,4,5,2],[6,7,4,3],[5,8,9,6],[7,10,11,8],[12,13,14,9],[10,14,13,11],[12]}} \\ \textbf{EM:} {\scriptsize\texttt{(B0, A0C0C3D3, B1D2E0B2, E3F0C1B3, C2F3G3D0, D1H0H3E1, I0H2H1E2, F1G2G1F2, G0)}} \\ \textbf{Kauffman bracket:} {\scriptsize $-A^{21} - A^{19} + A^{17} + 3A^{15} - 3A^{11} - 2A^{9} + 2A^{7} + 2A^{5} - A^{3} - A$} \\ \textbf{Arrow:} {\scriptsize $A^{12} + A^{10}L_1 - A^{8} - 3A^{6}L_1 - A^{4}L_2 + A^{4} + 3A^{2}L_1 + 2L_2 - 2L_1/A^{2} - 2L_2/A^{4} + L_1/A^{6} + L_2/A^{8}$} \\ \textbf{Mock:} {\scriptsize $w^{5} - 3w^{3} - w^{2} + 3w + 2 - 2/w + w^{-3}$} \\ \textbf{Affine:} {\scriptsize $-t^{2} - t + 4 - 1/t - 1/t^{2}$} \\ \textbf{Yamada:} {\scriptsize $A^{27} - 2A^{25} + A^{24} + 3A^{23} - 3A^{22} - 2A^{21} + 4A^{20} - A^{19} - 5A^{18} + A^{17} + 2A^{16} - 4A^{15} + 3A^{13} - 2A^{12} - A^{11} + A^{10} + 3A^{9} - 5A^{8} - A^{7} + 5A^{6} - 4A^{5} - A^{4} + 2A^{3} - 1$}
\end{minipage}

\noindent{\color{gray!40}\rule{\textwidth}{0.4pt}}
\vspace{0.9\baselineskip}
\noindent \begin{minipage}[t]{0.25\textwidth}
\vspace{0pt}
\centering
\includegraphics[page=337,width=\linewidth]{knotoids.pdf}
\end{minipage}
\hfill
\begin{minipage}[t]{0.73\textwidth}
\vspace{0pt}
\raggedright
\textbf{Name:} {\large{$\mathbf{K7_{218}}$}} (chiral, non-rotatable$^{*}$) \\ \textbf{PD:} {\scriptsize\texttt{[0],[0,1,2,3],[1,4,5,2],[6,7,4,3],[5,8,9,10],[10,11,12,6],[7,13,14,8],[9],[11,14,13,12]}} \\ \textbf{EM:} {\scriptsize\texttt{(B0, A0C0C3D3, B1D2E0B2, F3G0C1B3, C2G3H0F0, E3I0I3D0, D1I2I1E1, E2, F1G2G1F2)}} \\ \textbf{Kauffman bracket:} {\scriptsize $-A^{21} - 2A^{19} + 2A^{15} + A^{13} - 2A^{11} - A^{9} + 2A^{7} + 2A^{5} - A^{3} - A$} \\ \textbf{Arrow:} {\scriptsize $L_1/A^{6} + 2/A^{8} - 2/A^{12} - L_1/A^{14} + 2/A^{16} + L_1/A^{18} - 2/A^{20} - 2L_1/A^{22} + A^{-24} + L_1/A^{26}$} \\ \textbf{Mock:} {\scriptsize $w^{5} + 2w^{4} - 2w^{2} - w + 2 + 1/w - 2/w^{2} - 2/w^{3} + w^{-4} + w^{-5}$} \\ \textbf{Affine:} {\scriptsize $2t - 4 + 2/t$} \\ \textbf{Yamada:} {\scriptsize $A^{27} - A^{25} + 2A^{24} + 2A^{23} - 2A^{22} + 2A^{21} + 3A^{20} - A^{19} + 3A^{17} - A^{16} - 2A^{15} + A^{14} - 2A^{12} + A^{10} - 3A^{8} + 2A^{7} + 2A^{6} - 3A^{5} + 3A^{4} - A^{2} + A - 1$}
\end{minipage}

\noindent{\color{gray!40}\rule{\textwidth}{0.4pt}}
\vspace{0.9\baselineskip}
\noindent \begin{minipage}[t]{0.25\textwidth}
\vspace{0pt}
\centering
\includegraphics[page=338,width=\linewidth]{knotoids.pdf}
\end{minipage}
\hfill
\begin{minipage}[t]{0.73\textwidth}
\vspace{0pt}
\raggedright
\textbf{Name:} {\large{$\mathbf{K7_{219}}$}} (chiral, non-rotatable$^{*}$) \\ \textbf{PD:} {\scriptsize\texttt{[0],[0,1,2,3],[1,4,5,2],[6,7,4,3],[8,9,6,5],[7,10,11,12],[13,14,9,8],[14,13,12,10],[11]}} \\ \textbf{EM:} {\scriptsize\texttt{(B0, A0C0C3D3, B1D2E3B2, E2F0C1B3, G3G2D0C2, D1H3I0H2, H1H0E1E0, G1G0F3F1, F2)}} \\ \textbf{Kauffman bracket:} {\scriptsize $A^{21} + A^{19} - 2A^{17} - 4A^{15} + 4A^{11} + 2A^{9} - 2A^{7} - 2A^{5} + A$} \\ \textbf{Arrow:} {\scriptsize $-A^{6}L_1 - A^{4}L_2 + 2A^{2}L_1 + 2L_2 + 2 - 2L_2/A^{4} - 2/A^{4} - 2L_1/A^{6} + L_2/A^{8} + A^{-8} + 2L_1/A^{10} - L_1/A^{14}$} \\ \textbf{Mock:} {\scriptsize $-w^{3} - 2w^{2} + w + 6 + 2/w - 3/w^{2} - 2/w^{3}$} \\ \textbf{Affine:} {\scriptsize $-t^{2} + t + 1/t - 1/t^{2}$} \\ \textbf{Yamada:} {\scriptsize $-A^{24} - 2A^{23} + 3A^{22} + A^{21} - 6A^{20} + 5A^{19} + 5A^{18} - 7A^{17} + 3A^{16} + 5A^{15} - 3A^{14} + A^{13} + A^{12} + 3A^{11} - 4A^{10} - A^{9} + 7A^{8} - 5A^{7} - 3A^{6} + 7A^{5} - 2A^{4} - 4A^{3} + 3A^{2} + A - 1$}
\end{minipage}

\noindent{\color{gray!40}\rule{\textwidth}{0.4pt}}
\vspace{0.9\baselineskip}
\noindent \begin{minipage}[t]{0.25\textwidth}
\vspace{0pt}
\centering
\includegraphics[page=339,width=\linewidth]{knotoids.pdf}
\end{minipage}
\hfill
\begin{minipage}[t]{0.73\textwidth}
\vspace{0pt}
\raggedright
\textbf{Name:} {\large{$\mathbf{K7_{220}}$}} (chiral, non-rotatable$^{*}$) \\ \textbf{PD:} {\scriptsize\texttt{[0],[0,1,2,3],[1,4,5,2],[6,7,4,3],[8,9,10,5],[10,11,12,6],[7,13,14,8],[9],[11,14,13,12]}} \\ \textbf{EM:} {\scriptsize\texttt{(B0, A0C0C3D3, B1D2E3B2, F3G0C1B3, G3H0F0C2, E2I0I3D0, D1I2I1E0, E1, F1G2G1F2)}} \\ \textbf{Kauffman bracket:} {\scriptsize $A^{24} + 2A^{22} + A^{20} - 3A^{18} - 2A^{16} + 3A^{14} + 3A^{12} - 3A^{10} - 4A^{8} + A^{6} + 3A^{4} - 1$} \\ \textbf{Arrow:} {\scriptsize $A^{-12} + 2L_1/A^{14} + A^{-16} - 3L_1/A^{18} - 2/A^{20} + 3L_1/A^{22} + 3/A^{24} - 3L_1/A^{26} - 4/A^{28} + L_1/A^{30} + 3/A^{32} - 1/A^{36}$} \\ \textbf{Mock:} {\scriptsize $w^{6} + 2w^{5} + w^{4} - 3w^{3} - 3w^{2} + 3w + 4 - 3/w - 4/w^{2} + w^{-3} + 2/w^{4}$} \\ \textbf{Affine:} {\scriptsize $2t - 4 + 2/t$} \\ \textbf{Yamada:} {\scriptsize $-A^{28} + A^{26} - 3A^{25} - 3A^{24} + 4A^{23} - 3A^{22} - 8A^{21} + 4A^{20} + 3A^{19} - 9A^{18} + 3A^{17} + 7A^{16} - 5A^{15} + 6A^{13} - 3A^{11} + 7A^{9} - 5A^{8} - 5A^{7} + 10A^{6} - 5A^{5} - 5A^{4} + 5A^{3} - A^{2} - A + 1$}
\end{minipage}

\noindent{\color{gray!40}\rule{\textwidth}{0.4pt}}
\vspace{0.9\baselineskip}
\noindent \begin{minipage}[t]{0.25\textwidth}
\vspace{0pt}
\centering
\includegraphics[page=340,width=\linewidth]{knotoids.pdf}
\end{minipage}
\hfill
\begin{minipage}[t]{0.73\textwidth}
\vspace{0pt}
\raggedright
\textbf{Name:} {\large{$\mathbf{K7_{221}}$}} (chiral, rotatable) \\ \textbf{PD:} {\scriptsize\texttt{[0],[0,1,2,3],[1,4,5,2],[6,7,8,3],[4,8,9,5],[6,9,10,11],[7,12,13,14],[14,12,11,10],[13]}} \\ \textbf{EM:} {\scriptsize\texttt{(B0, A0C0C3D3, B1E0E3B2, F0G0E1B3, C1D2F1C2, D0E2H3H2, D1H1I0H0, G3G1F3F2, G2)}} \\ \textbf{Kauffman bracket:} {\scriptsize $A^{22} + 2A^{20} - A^{18} - 4A^{16} - A^{14} + 4A^{12} + 3A^{10} - 2A^{8} - 2A^{6} + A^{2}$} \\ \textbf{Arrow:} {\scriptsize $A^{4} + 2A^{2}L_1 - 1 - 4L_1/A^{2} - 1/A^{4} + 4L_1/A^{6} + 3/A^{8} - 2L_1/A^{10} - 2/A^{12} + A^{-16}$} \\ \textbf{Mock:} {\scriptsize $-w^{4} - 2w^{3} + w^{2} + 6w + 4 - 4/w - 3/w^{2}$} \\ \textbf{Affine:} {\scriptsize $2t - 4 + 2/t$} \\ \textbf{Yamada:} {\scriptsize $-A^{25} - A^{24} + 3A^{23} - 5A^{21} + 4A^{20} + 4A^{19} - 8A^{18} + 2A^{17} + 5A^{16} - 7A^{15} - A^{14} + 3A^{13} - 2A^{12} - 2A^{11} - A^{10} + 5A^{9} - 4A^{8} - 5A^{7} + 9A^{6} - 3A^{5} - 5A^{4} + 6A^{3} - A^{2} - 2A + 1$}
\end{minipage}

\noindent{\color{gray!40}\rule{\textwidth}{0.4pt}}
\vspace{0.9\baselineskip}
\noindent \begin{minipage}[t]{0.25\textwidth}
\vspace{0pt}
\centering
\includegraphics[page=341,width=\linewidth]{knotoids.pdf}
\end{minipage}
\hfill
\begin{minipage}[t]{0.73\textwidth}
\vspace{0pt}
\raggedright
\textbf{Name:} {\large{$\mathbf{K7_{222}}$}} (chiral, non-rotatable$^{*}$) \\ \textbf{PD:} {\scriptsize\texttt{[0],[0,1,2,3],[1,4,5,2],[6,7,8,3],[4,8,7,9],[5,10,11,6],[9,12,13,10],[14,13,12,11],[14]}} \\ \textbf{EM:} {\scriptsize\texttt{(B0, A0C0C3D3, B1E0F0B2, F3E2E1B3, C1D2D1G0, C2G3H3D0, E3H2H1F1, I0G2G1F2, H0)}} \\ \textbf{Kauffman bracket:} {\scriptsize $A^{18} + A^{16} - A^{14} - 3A^{12} + 3A^{8} + 3A^{6} - A^{4} - 2A^{2}$} \\ \textbf{Arrow:} {\scriptsize $1 + L_1/A^{2} - 1/A^{4} - 3L_1/A^{6} - L_2/A^{8} + A^{-8} + 3L_1/A^{10} + 2L_2/A^{12} + A^{-12} - L_1/A^{14} - L_2/A^{16} - 1/A^{16}$} \\ \textbf{Mock:} {\scriptsize $-w^{6} - w^{5} + 2w^{4} + 2w^{3} - 2w - 2 + 2/w + 3/w^{2} - 1/w^{3} - 1/w^{4}$} \\ \textbf{Affine:} {\scriptsize $0$} \\ \textbf{Yamada:} {\scriptsize $A^{27} - 2A^{25} + A^{24} + 2A^{23} - 3A^{22} - 3A^{21} + A^{20} + A^{19} - 5A^{18} - A^{17} + 4A^{16} - 3A^{15} - A^{14} + 4A^{13} - A^{12} - A^{11} + 2A^{9} - 4A^{8} - 3A^{7} + 5A^{6} - 2A^{4} + 2A^{3} + 2A^{2} - A - 1$}
\end{minipage}

\noindent{\color{gray!40}\rule{\textwidth}{0.4pt}}
\vspace{0.9\baselineskip}
\noindent \begin{minipage}[t]{0.25\textwidth}
\vspace{0pt}
\centering
\includegraphics[page=342,width=\linewidth]{knotoids.pdf}
\end{minipage}
\hfill
\begin{minipage}[t]{0.73\textwidth}
\vspace{0pt}
\raggedright
\textbf{Name:} {\large{$\mathbf{K7_{223}}$}} (chiral, non-rotatable$^{*}$) \\ \textbf{PD:} {\scriptsize\texttt{[0],[0,1,2,3],[1,4,5,2],[6,7,8,3],[4,8,7,9],[5,10,11,6],[12,13,10,9],[11,14,13,12],[14]}} \\ \textbf{EM:} {\scriptsize\texttt{(B0, A0C0C3D3, B1E0F0B2, F3E2E1B3, C1D2D1G3, C2G2H0D0, H3H2F1E3, F2I0G1G0, H1)}} \\ \textbf{Kauffman bracket:} {\scriptsize $-A^{21} - A^{19} - A^{17} + 2A^{13} + 2A^{11} - A^{9} - 2A^{7} - A^{5} + A^{3} + A$} \\ \textbf{Arrow:} {\scriptsize $A^{-12} + L_1/A^{14} + L_2/A^{16} - 2L_2/A^{20} - 2L_1/A^{22} + L_2/A^{24} + 2L_1/A^{26} + A^{-28} - L_1/A^{30} - 1/A^{32}$} \\ \textbf{Mock:} {\scriptsize $w^{6} + w^{5} - 2w + 2/w - 1/w^{2} - 1/w^{3} + w^{-4}$} \\ \textbf{Affine:} {\scriptsize $2t - 4 + 2/t$} \\ \textbf{Yamada:} {\scriptsize $-A^{28} - A^{27} - A^{24} - A^{23} - A^{22} + A^{21} - 2A^{20} - 3A^{19} + A^{18} - A^{17} - 3A^{16} + A^{15} + A^{13} + 2A^{11} + 2A^{10} - A^{9} + 3A^{8} + 2A^{7} - 2A^{6} - A^{5} + A^{4} - A^{3} - 2A^{2} + 1$}
\end{minipage}

\noindent{\color{gray!40}\rule{\textwidth}{0.4pt}}
\vspace{0.9\baselineskip}
\noindent \begin{minipage}[t]{0.25\textwidth}
\vspace{0pt}
\centering
\includegraphics[page=343,width=\linewidth]{knotoids.pdf}
\end{minipage}
\hfill
\begin{minipage}[t]{0.73\textwidth}
\vspace{0pt}
\raggedright
\textbf{Name:} {\large{$\mathbf{K7_{224}}$}} (chiral, non-rotatable$^{*}$) \\ \textbf{PD:} {\scriptsize\texttt{[0],[0,1,2,3],[1,4,5,2],[6,7,8,3],[4,8,7,9],[10,11,6,5],[9,12,13,14],[14,12,11,10],[13]}} \\ \textbf{EM:} {\scriptsize\texttt{(B0, A0C0C3D3, B1E0F3B2, F2E2E1B3, C1D2D1G0, H3H2D0C2, E3H1I0H0, G3G1F1F0, G2)}} \\ \textbf{Kauffman bracket:} {\scriptsize $A^{22} + 3A^{20} - 5A^{16} - 3A^{14} + 4A^{12} + 4A^{10} - A^{8} - 3A^{6} + A^{2}$} \\ \textbf{Arrow:} {\scriptsize $L_1/A^{2} + L_2/A^{4} + 2/A^{4} - 2L_2/A^{8} - 3/A^{8} - 3L_1/A^{10} + L_2/A^{12} + 3/A^{12} + 4L_1/A^{14} - 1/A^{16} - 3L_1/A^{18} + L_1/A^{22}$} \\ \textbf{Mock:} {\scriptsize $-w^{4} + 5w^{2} + 3w - 5 - 5/w + 2/w^{2} + 2/w^{3}$} \\ \textbf{Affine:} {\scriptsize $t - 2 + 1/t$} \\ \textbf{Yamada:} {\scriptsize $A^{26} + 2A^{25} - 3A^{24} - 2A^{23} + 7A^{22} - 2A^{21} - 9A^{20} + 10A^{19} + 4A^{18} - 10A^{17} + 6A^{16} + 6A^{15} - 5A^{14} + 3A^{12} + 3A^{11} - 8A^{10} + A^{9} + 10A^{8} - 11A^{7} - 2A^{6} + 10A^{5} - 5A^{4} - 4A^{3} + 4A^{2} + A - 1$}
\end{minipage}

\noindent{\color{gray!40}\rule{\textwidth}{0.4pt}}
\vspace{0.9\baselineskip}
\noindent \begin{minipage}[t]{0.25\textwidth}
\vspace{0pt}
\centering
\includegraphics[page=344,width=\linewidth]{knotoids.pdf}
\end{minipage}
\hfill
\begin{minipage}[t]{0.73\textwidth}
\vspace{0pt}
\raggedright
\textbf{Name:} {\large{$\mathbf{K7_{225}}$}} (chiral, non-rotatable$^{*}$) \\ \textbf{PD:} {\scriptsize\texttt{[0],[0,1,2,3],[1,4,5,2],[6,7,8,3],[4,8,7,9],[10,11,6,5],[12,13,10,9],[11,14,13,12],[14]}} \\ \textbf{EM:} {\scriptsize\texttt{(B0, A0C0C3D3, B1E0F3B2, F2E2E1B3, C1D2D1G3, G2H0D0C2, H3H2F0E3, F1I0G1G0, H1)}} \\ \textbf{Kauffman bracket:} {\scriptsize $2A^{16} + 2A^{14} - 2A^{12} - 4A^{10} + 4A^{6} + 2A^{4} - 2A^{2} - 1$} \\ \textbf{Arrow:} {\scriptsize $L_2/A^{8} + A^{-8} + 2L_1/A^{10} - 2L_2/A^{12} - 4L_1/A^{14} + L_2/A^{16} - 1/A^{16} + 4L_1/A^{18} + 2/A^{20} - 2L_1/A^{22} - 1/A^{24}$} \\ \textbf{Mock:} {\scriptsize $w^{4} + 2w^{3} + w^{2} - 3w - 1 + 3/w - 1/w^{2} - 3/w^{3} + w^{-4} + w^{-5}$} \\ \textbf{Affine:} {\scriptsize $2t - 4 + 2/t$} \\ \textbf{Yamada:} {\scriptsize $-A^{26} - A^{25} + A^{24} + 2A^{23} - 2A^{22} - 3A^{21} + 4A^{20} + A^{19} - 8A^{18} + 2A^{17} + 2A^{16} - 8A^{15} + 3A^{13} - 3A^{12} + 2A^{10} + 5A^{9} - 3A^{8} - A^{7} + 8A^{6} - 5A^{5} - 4A^{4} + 5A^{3} - A^{2} - 2A + 1$}
\end{minipage}

\noindent{\color{gray!40}\rule{\textwidth}{0.4pt}}
\vspace{0.9\baselineskip}
\noindent \begin{minipage}[t]{0.25\textwidth}
\vspace{0pt}
\centering
\includegraphics[page=345,width=\linewidth]{knotoids.pdf}
\end{minipage}
\hfill
\begin{minipage}[t]{0.73\textwidth}
\vspace{0pt}
\raggedright
\textbf{Name:} {\large{$\mathbf{K7_{226}}$}} (chiral, non-rotatable$^{*}$) \\ \textbf{PD:} {\scriptsize\texttt{[0],[0,1,2,3],[1,4,5,2],[3,6,7,4],[5,8,9,6],[7,9,10,11],[8,12,13,10],[14,13,12,11],[14]}} \\ \textbf{EM:} {\scriptsize\texttt{(B0, A0C0C3D0, B1D3E0B2, B3E3F0C1, C2G0F1D1, D2E2G3H3, E1H2H1F2, I0G2G1F3, H0)}} \\ \textbf{Kauffman bracket:} {\scriptsize $A^{22} + 2A^{20} - A^{18} - 4A^{16} + 4A^{12} + 2A^{10} - 3A^{8} - 2A^{6} + A^{4} + A^{2}$} \\ \textbf{Arrow:} {\scriptsize $A^{4} + 2A^{2}L_1 - 1 - 4L_1/A^{2} - L_2/A^{4} + A^{-4} + 4L_1/A^{6} + 2L_2/A^{8} - 3L_1/A^{10} - 2L_2/A^{12} + L_1/A^{14} + L_2/A^{16}$} \\ \textbf{Mock:} {\scriptsize $w^{5} - 3w^{3} + 5w + 2 - 4/w - 1/w^{2} + w^{-3}$} \\ \textbf{Affine:} {\scriptsize $-t^{2} + t + 1/t - 1/t^{2}$} \\ \textbf{Yamada:} {\scriptsize $A^{25} - A^{24} - 2A^{23} + 3A^{22} - 7A^{20} + 4A^{19} + 4A^{18} - 8A^{17} + 2A^{16} + 4A^{15} - 5A^{14} - A^{13} + 2A^{11} - 5A^{10} - A^{9} + 8A^{8} - 5A^{7} - 2A^{6} + 7A^{5} - 2A^{4} - 4A^{3} + 3A^{2} + A - 2$}
\end{minipage}

\noindent{\color{gray!40}\rule{\textwidth}{0.4pt}}
\vspace{0.9\baselineskip}
\noindent \begin{minipage}[t]{0.25\textwidth}
\vspace{0pt}
\centering
\includegraphics[page=346,width=\linewidth]{knotoids.pdf}
\end{minipage}
\hfill
\begin{minipage}[t]{0.73\textwidth}
\vspace{0pt}
\raggedright
\textbf{Name:} {\large{$\mathbf{K7_{227}}$}} (chiral, non-rotatable$^{*}$) \\ \textbf{PD:} {\scriptsize\texttt{[0],[0,1,2,3],[1,4,5,2],[3,6,7,4],[5,7,8,9],[6,10,11,8],[9,12,13,10],[14,13,12,11],[14]}} \\ \textbf{EM:} {\scriptsize\texttt{(B0, A0C0C3D0, B1D3E0B2, B3F0E1C1, C2D2F3G0, D1G3H3E2, E3H2H1F1, I0G2G1F2, H0)}} \\ \textbf{Kauffman bracket:} {\scriptsize $A^{22} + 2A^{20} - 2A^{18} - 5A^{16} + 6A^{12} + 4A^{10} - 3A^{8} - 3A^{6} + A^{2}$} \\ \textbf{Arrow:} {\scriptsize $A^{4} + 2A^{2}L_1 - 2 - 5L_1/A^{2} - 2L_2/A^{4} + 2/A^{4} + 6L_1/A^{6} + 4L_2/A^{8} - 3L_1/A^{10} - 3L_2/A^{12} + L_2/A^{16}$} \\ \textbf{Mock:} {\scriptsize $w^{5} - 4w^{3} + 7w + 3 - 4/w - 2/w^{2}$} \\ \textbf{Affine:} {\scriptsize $-t^{2} + 2t - 2 + 2/t - 1/t^{2}$} \\ \textbf{Yamada:} {\scriptsize $2A^{25} - A^{24} - 5A^{23} + 6A^{22} + 3A^{21} - 13A^{20} + 6A^{19} + 9A^{18} - 14A^{17} + A^{16} + 7A^{15} - 7A^{14} - 3A^{13} + A^{12} + 6A^{11} - 7A^{10} - 3A^{9} + 14A^{8} - 7A^{7} - 8A^{6} + 12A^{5} - 2A^{4} - 8A^{3} + 5A^{2} + 2A - 2$}
\end{minipage}

\noindent{\color{gray!40}\rule{\textwidth}{0.4pt}}
\vspace{0.9\baselineskip}
\noindent \begin{minipage}[t]{0.25\textwidth}
\vspace{0pt}
\centering
\includegraphics[page=347,width=\linewidth]{knotoids.pdf}
\end{minipage}
\hfill
\begin{minipage}[t]{0.73\textwidth}
\vspace{0pt}
\raggedright
\textbf{Name:} {\large{$\mathbf{K7_{228}}$}} (chiral, non-rotatable$^{*}$) \\ \textbf{PD:} {\scriptsize\texttt{[0],[0,1,2,3],[1,4,5,2],[3,6,7,4],[5,7,8,9],[6,10,11,8],[12,13,10,9],[11,14,13,12],[14]}} \\ \textbf{EM:} {\scriptsize\texttt{(B0, A0C0C3D0, B1D3E0B2, B3F0E1C1, C2D2F3G3, D1G2H0E2, H3H2F1E3, F2I0G1G0, H1)}} \\ \textbf{Kauffman bracket:} {\scriptsize $A^{21} - A^{19} - 4A^{17} - 2A^{15} + 4A^{13} + 5A^{11} - A^{9} - 4A^{7} - A^{5} + A^{3} + A$} \\ \textbf{Arrow:} {\scriptsize $-L_1/A^{6} + A^{-8} + 4L_1/A^{10} + 2L_2/A^{12} - 4L_1/A^{14} - 4L_2/A^{16} - 1/A^{16} + L_1/A^{18} + 3L_2/A^{20} + A^{-20} + L_1/A^{22} - L_2/A^{24} - L_1/A^{26}$} \\ \textbf{Mock:} {\scriptsize $-w^{5} + 4w^{3} + 4w^{2} - 3w - 5 + 2/w^{2}$} \\ \textbf{Affine:} {\scriptsize $t^{2} + 2t - 6 + 2/t + t^{-2}$} \\ \textbf{Yamada:} {\scriptsize $-2A^{26} - A^{25} + 4A^{24} - 2A^{23} - 8A^{22} + 8A^{21} + 2A^{20} - 13A^{19} + 4A^{18} + 4A^{17} - 7A^{16} + A^{15} + 4A^{14} + 3A^{13} - 6A^{12} + 2A^{11} + 8A^{10} - 10A^{9} - 2A^{8} + 10A^{7} - 5A^{6} - 5A^{5} + 6A^{4} + A^{3} - 3A^{2} + 1$}
\end{minipage}

\noindent{\color{gray!40}\rule{\textwidth}{0.4pt}}
\vspace{0.9\baselineskip}
\noindent \begin{minipage}[t]{0.25\textwidth}
\vspace{0pt}
\centering
\includegraphics[page=348,width=\linewidth]{knotoids.pdf}
\end{minipage}
\hfill
\begin{minipage}[t]{0.73\textwidth}
\vspace{0pt}
\raggedright
\textbf{Name:} {\large{$\mathbf{K7_{229}}$}} (chiral, non-rotatable$^{*}$) \\ \textbf{PD:} {\scriptsize\texttt{[0],[0,1,2,3],[1,4,5,2],[3,6,7,4],[5,7,8,9],[6,9,10,11],[11,12,13,8],[14,13,12,10],[14]}} \\ \textbf{EM:} {\scriptsize\texttt{(B0, A0C0C3D0, B1D3E0B2, B3F0E1C1, C2D2G3F1, D1E3H3G0, F3H2H1E2, I0G2G1F2, H0)}} \\ \textbf{Kauffman bracket:} {\scriptsize $-A^{24} + A^{22} + 3A^{20} - 5A^{16} - 2A^{14} + 5A^{12} + 4A^{10} - 2A^{8} - 3A^{6} + A^{2}$} \\ \textbf{Arrow:} {\scriptsize $-A^{6}L_1 + A^{4} + 3A^{2}L_1 + L_2 - 1 - 5L_1/A^{2} - 3L_2/A^{4} + A^{-4} + 5L_1/A^{6} + 4L_2/A^{8} - 2L_1/A^{10} - 3L_2/A^{12} + L_2/A^{16}$} \\ \textbf{Mock:} {\scriptsize $w^{5} - 5w^{3} - 2w^{2} + 7w + 5 - 3/w - 2/w^{2}$} \\ \textbf{Affine:} {\scriptsize $0$} \\ \textbf{Yamada:} {\scriptsize $-A^{26} + 2A^{25} + 2A^{24} - 7A^{23} + 2A^{22} + 9A^{21} - 11A^{20} - 2A^{19} + 12A^{18} - 9A^{17} - 5A^{16} + 7A^{15} - 2A^{14} - 3A^{13} - 3A^{12} + 7A^{11} - 3A^{10} - 10A^{9} + 12A^{8} + A^{7} - 11A^{6} + 8A^{5} + 4A^{4} - 7A^{3} + A^{2} + 2A - 1$}
\end{minipage}

\noindent{\color{gray!40}\rule{\textwidth}{0.4pt}}
\vspace{0.9\baselineskip}
\noindent \begin{minipage}[t]{0.25\textwidth}
\vspace{0pt}
\centering
\includegraphics[page=349,width=\linewidth]{knotoids.pdf}
\end{minipage}
\hfill
\begin{minipage}[t]{0.73\textwidth}
\vspace{0pt}
\raggedright
\textbf{Name:} {\large{$\mathbf{K7_{230}}$}} (chiral, non-rotatable$^{*}$) \\ \textbf{PD:} {\scriptsize\texttt{[0],[0,1,2,3],[1,4,5,2],[3,6,7,4],[5,8,9,10],[6,11,12,13],[13,14,8,7],[9],[14,12,11,10]}} \\ \textbf{EM:} {\scriptsize\texttt{(B0, A0C0C3D0, B1D3E0B2, B3F0G3C1, C2G2H0I3, D1I2I1G0, F3I0E1D2, E2, G1F2F1E3)}} \\ \textbf{Kauffman bracket:} {\scriptsize $-A^{24} + A^{22} + 3A^{20} - 4A^{16} + 5A^{12} + 2A^{10} - 4A^{8} - 3A^{6} + A^{4} + A^{2}$} \\ \textbf{Arrow:} {\scriptsize $-L_1/A^{6} + A^{-8} + 3L_1/A^{10} + L_2/A^{12} - 1/A^{12} - 4L_1/A^{14} - 2L_2/A^{16} + 2/A^{16} + 5L_1/A^{18} + 3L_2/A^{20} - 1/A^{20} - 4L_1/A^{22} - 3L_2/A^{24} + L_1/A^{26} + L_2/A^{28}$} \\ \textbf{Mock:} {\scriptsize $-w^{5} + w^{4} + 5w^{3} + w^{2} - 7w - 3 + 4/w + 2/w^{2} - 1/w^{3}$} \\ \textbf{Affine:} {\scriptsize $t^{2} + t - 4 + 1/t + t^{-2}$} \\ \textbf{Yamada:} {\scriptsize $A^{27} - A^{26} + 6A^{24} - 3A^{23} - 4A^{22} + 11A^{21} - 4A^{20} - 6A^{19} + 11A^{18} - 2A^{17} - 5A^{16} + 5A^{15} - 5A^{12} + 4A^{11} + 4A^{10} - 11A^{9} + 5A^{8} + 6A^{7} - 10A^{6} + 2A^{5} + 5A^{4} - 5A^{3} + A^{2} + 2A - 1$}
\end{minipage}

\noindent{\color{gray!40}\rule{\textwidth}{0.4pt}}
\vspace{0.9\baselineskip}
\noindent \begin{minipage}[t]{0.25\textwidth}
\vspace{0pt}
\centering
\includegraphics[page=350,width=\linewidth]{knotoids.pdf}
\end{minipage}
\hfill
\begin{minipage}[t]{0.73\textwidth}
\vspace{0pt}
\raggedright
\textbf{Name:} {\large{$\mathbf{K7_{231}}$}} (chiral, non-rotatable$^{*}$) \\ \textbf{PD:} {\scriptsize\texttt{[0],[0,1,2,3],[1,4,5,2],[3,6,7,4],[8,9,6,5],[7,9,10,11],[8,12,13,10],[14,13,12,11],[14]}} \\ \textbf{EM:} {\scriptsize\texttt{(B0, A0C0C3D0, B1D3E3B2, B3E2F0C1, G0F1D1C2, D2E1G3H3, E0H2H1F2, I0G2G1F3, H0)}} \\ \textbf{Kauffman bracket:} {\scriptsize $-A^{21} - A^{19} + A^{17} + 2A^{15} - 3A^{11} - 2A^{9} + 2A^{7} + 2A^{5} - A$} \\ \textbf{Arrow:} {\scriptsize $1 + L_1/A^{2} - 1/A^{4} - 2L_1/A^{6} - L_2/A^{8} + A^{-8} + 3L_1/A^{10} + 2L_2/A^{12} - 2L_1/A^{14} - 2L_2/A^{16} + L_2/A^{20}$} \\ \textbf{Mock:} {\scriptsize $-w^{6} - w^{5} + 2w^{4} + 3w^{3} + w^{2} - 2w - 2 + 1/w + 2/w^{2} - 1/w^{3} - 1/w^{4}$} \\ \textbf{Affine:} {\scriptsize $-t^{2} + 2t - 2 + 2/t - 1/t^{2}$} \\ \textbf{Yamada:} {\scriptsize $A^{25} - 2A^{23} + 2A^{21} - 2A^{20} - 3A^{19} + 2A^{18} - 5A^{16} + A^{15} + A^{14} - 4A^{13} + 3A^{11} - A^{10} + A^{8} + 3A^{7} - 3A^{6} - 2A^{5} + 5A^{4} - 3A^{3} - A^{2} + 2A - 1$}
\end{minipage}

\noindent{\color{gray!40}\rule{\textwidth}{0.4pt}}
\vspace{0.9\baselineskip}
\noindent \begin{minipage}[t]{0.25\textwidth}
\vspace{0pt}
\centering
\includegraphics[page=351,width=\linewidth]{knotoids.pdf}
\end{minipage}
\hfill
\begin{minipage}[t]{0.73\textwidth}
\vspace{0pt}
\raggedright
\textbf{Name:} {\large{$\mathbf{K7_{232}}$}} (chiral, non-rotatable$^{*}$) \\ \textbf{PD:} {\scriptsize\texttt{[0],[0,1,2,3],[1,4,5,2],[3,6,7,4],[8,9,6,5],[9,10,11,7],[12,13,10,8],[11,14,13,12],[14]}} \\ \textbf{EM:} {\scriptsize\texttt{(B0, A0C0C3D0, B1D3E3B2, B3E2F3C1, G3F0D1C2, E1G2H0D2, H3H2F1E0, F2I0G1G0, H1)}} \\ \textbf{Kauffman bracket:} {\scriptsize $-A^{25} + 3A^{21} + A^{19} - 5A^{17} - 4A^{15} + 5A^{13} + 7A^{11} - A^{9} - 6A^{7} - 2A^{5} + A^{3} + A$} \\ \textbf{Arrow:} {\scriptsize $A^{22}L_1 - 3A^{18}L_1 - A^{16}L_2 + 5A^{14}L_1 + 3A^{12}L_2 + A^{12} - 5A^{10}L_1 - 5A^{8}L_2 - 2A^{8} + A^{6}L_1 + 4A^{4}L_2 + 2A^{4} + 2A^{2}L_1 - L_2 - L_1/A^{2}$} \\ \textbf{Mock:} {\scriptsize $-w^{5} - w^{4} + 4w^{3} + 4w^{2} - 6w - 6 + 3/w + 6/w^{2} + w^{-3} - 2/w^{4} - 1/w^{5}$} \\ \textbf{Affine:} {\scriptsize $-t^{2} + 2 - 1/t^{2}$} \\ \textbf{Yamada:} {\scriptsize $2A^{27} - 2A^{26} - 6A^{25} + 10A^{24} + 4A^{23} - 19A^{22} + 11A^{21} + 15A^{20} - 23A^{19} + 2A^{18} + 15A^{17} - 10A^{16} - 6A^{15} + 6A^{14} + 8A^{13} - 13A^{12} - 4A^{11} + 20A^{10} - 13A^{9} - 14A^{8} + 21A^{7} - 2A^{6} - 15A^{5} + 8A^{4} + 5A^{3} - 7A^{2} - A + 2$}
\end{minipage}

\noindent{\color{gray!40}\rule{\textwidth}{0.4pt}}
\vspace{0.9\baselineskip}
\noindent \begin{minipage}[t]{0.25\textwidth}
\vspace{0pt}
\centering
\includegraphics[page=352,width=\linewidth]{knotoids.pdf}
\end{minipage}
\hfill
\begin{minipage}[t]{0.73\textwidth}
\vspace{0pt}
\raggedright
\textbf{Name:} {\large{$\mathbf{K7_{233}}$}} (chiral, non-rotatable$^{*}$) \\ \textbf{PD:} {\scriptsize\texttt{[0],[0,1,2,3],[1,4,5,2],[3,6,7,4],[7,8,9,5],[6,10,11,8],[12,13,10,9],[11,14,13,12],[14]}} \\ \textbf{EM:} {\scriptsize\texttt{(B0, A0C0C3D0, B1D3E3B2, B3F0E0C1, D2F3G3C2, D1G2H0E1, H3H2F1E2, F2I0G1G0, H1)}} \\ \textbf{Kauffman bracket:} {\scriptsize $-A^{20} + 2A^{16} + 2A^{14} - A^{12} - 3A^{10} + 3A^{6} + A^{4} - A^{2} - 1$} \\ \textbf{Arrow:} {\scriptsize $-A^{2}L_1 + 2L_1/A^{2} + L_2/A^{4} + A^{-4} - L_1/A^{6} - 2L_2/A^{8} - 1/A^{8} + 2L_2/A^{12} + A^{-12} + L_1/A^{14} - L_2/A^{16} - L_1/A^{18}$} \\ \textbf{Mock:} {\scriptsize $-w^{5} - w^{4} + 2w^{3} + 4w^{2} - 3 - 1/w + w^{-2}$} \\ \textbf{Affine:} {\scriptsize $t^{2} + t - 4 + 1/t + t^{-2}$} \\ \textbf{Yamada:} {\scriptsize $A^{23} + A^{22} - 2A^{21} - A^{20} + 5A^{19} - A^{18} - 3A^{17} + 6A^{16} + A^{15} - 2A^{14} + 3A^{13} - 3A^{10} + 2A^{9} + A^{8} - 5A^{7} + 3A^{6} + 3A^{5} - 4A^{4} + 3A^{2} - A - 1$}
\end{minipage}

\noindent{\color{gray!40}\rule{\textwidth}{0.4pt}}
\vspace{0.9\baselineskip}
\noindent \begin{minipage}[t]{0.25\textwidth}
\vspace{0pt}
\centering
\includegraphics[page=353,width=\linewidth]{knotoids.pdf}
\end{minipage}
\hfill
\begin{minipage}[t]{0.73\textwidth}
\vspace{0pt}
\raggedright
\textbf{Name:} {\large{$\mathbf{K7_{234}}$}} (chiral, non-rotatable$^{*}$) \\ \textbf{PD:} {\scriptsize\texttt{[0],[0,1,2,3],[1,4,5,2],[3,6,7,4],[7,8,9,5],[6,9,10,11],[11,12,13,8],[14,13,12,10],[14]}} \\ \textbf{EM:} {\scriptsize\texttt{(B0, A0C0C3D0, B1D3E3B2, B3F0E0C1, D2G3F1C2, D1E2H3G0, F3H2H1E1, I0G2G1F2, H0)}} \\ \textbf{Kauffman bracket:} {\scriptsize $A^{23} - 2A^{19} - 2A^{17} + 2A^{15} + 3A^{13} - A^{11} - 3A^{9} + 2A^{5} - A$} \\ \textbf{Arrow:} {\scriptsize $-A^{14}L_1 + 2A^{10}L_1 + A^{8}L_2 + A^{8} - 2A^{6}L_1 - 2A^{4}L_2 - A^{4} + A^{2}L_1 + 2L_2 + 1 - 2L_2/A^{4} + L_2/A^{8}$} \\ \textbf{Mock:} {\scriptsize $w^{5} + w^{4} - 3w^{3} - 4w^{2} + 2w + 5 - 1/w^{2}$} \\ \textbf{Affine:} {\scriptsize $-t + 2 - 1/t$} \\ \textbf{Yamada:} {\scriptsize $-A^{24} + 3A^{22} - A^{21} - 3A^{20} + 4A^{19} - 5A^{17} + 3A^{16} + A^{15} - 4A^{14} - A^{12} - 4A^{10} + A^{9} + 4A^{8} - 5A^{7} + A^{6} + 4A^{5} - 3A^{4} - A^{3} + 2A^{2} - 1$}
\end{minipage}

\noindent{\color{gray!40}\rule{\textwidth}{0.4pt}}
\vspace{0.9\baselineskip}
\noindent \begin{minipage}[t]{0.25\textwidth}
\vspace{0pt}
\centering
\includegraphics[page=354,width=\linewidth]{knotoids.pdf}
\end{minipage}
\hfill
\begin{minipage}[t]{0.73\textwidth}
\vspace{0pt}
\raggedright
\textbf{Name:} {\large{$\mathbf{K7_{235}}$}} (chiral, non-rotatable$^{*}$) \\ \textbf{PD:} {\scriptsize\texttt{[0],[0,1,2,3],[1,4,5,2],[3,6,7,4],[8,9,10,5],[6,11,12,13],[13,14,8,7],[9],[14,12,11,10]}} \\ \textbf{EM:} {\scriptsize\texttt{(B0, A0C0C3D0, B1D3E3B2, B3F0G3C1, G2H0I3C2, D1I2I1G0, F3I0E0D2, E1, G1F2F1E2)}} \\ \textbf{Kauffman bracket:} {\scriptsize $-A^{26} + 2A^{22} + 2A^{20} - 2A^{18} - 3A^{16} + 2A^{14} + 5A^{12} - 4A^{8} - 2A^{6} + A^{4} + A^{2}$} \\ \textbf{Arrow:} {\scriptsize $-A^{2}L_1 + 2L_1/A^{2} + L_2/A^{4} + A^{-4} - 2L_1/A^{6} - 2L_2/A^{8} - 1/A^{8} + 2L_1/A^{10} + 3L_2/A^{12} + 2/A^{12} - 3L_2/A^{16} - 1/A^{16} - 2L_1/A^{18} + L_2/A^{20} + L_1/A^{22}$} \\ \textbf{Mock:} {\scriptsize $-w^{5} - w^{4} + 3w^{3} + 5w^{2} - 2w - 6 - 1/w + 3/w^{2} + w^{-3}$} \\ \textbf{Affine:} {\scriptsize $t^{2} - 2 + t^{-2}$} \\ \textbf{Yamada:} {\scriptsize $2A^{25} + A^{24} - 5A^{23} + 3A^{22} + 7A^{21} - 10A^{20} + A^{19} + 11A^{18} - 7A^{17} - A^{16} + 8A^{15} - 2A^{14} - 2A^{13} - A^{12} + 6A^{11} - 3A^{10} - 7A^{9} + 11A^{8} - 2A^{7} - 10A^{6} + 7A^{5} + 2A^{4} - 6A^{3} + 2A^{2} + 2A - 1$}
\end{minipage}

\noindent{\color{gray!40}\rule{\textwidth}{0.4pt}}
\vspace{0.9\baselineskip}
\noindent \begin{minipage}[t]{0.25\textwidth}
\vspace{0pt}
\centering
\includegraphics[page=355,width=\linewidth]{knotoids.pdf}
\end{minipage}
\hfill
\begin{minipage}[t]{0.73\textwidth}
\vspace{0pt}
\raggedright
\textbf{Name:} {\large{$\mathbf{K7_{236}}$}} (chiral, non-rotatable$^{*}$) \\ \textbf{PD:} {\scriptsize\texttt{[0],[0,1,2,3],[1,4,5,2],[3,5,6,7],[4,7,8,9],[9,10,11,6],[12,13,10,8],[11,14,13,12],[14]}} \\ \textbf{EM:} {\scriptsize\texttt{(B0, A0C0C3D0, B1E0D1B2, B3C2F3E1, C1D3G3F0, E3G2H0D2, H3H2F1E2, F2I0G1G0, H1)}} \\ \textbf{Kauffman bracket:} {\scriptsize $-A^{25} + 3A^{21} + 2A^{19} - 4A^{17} - 5A^{15} + 3A^{13} + 6A^{11} - 5A^{7} - 2A^{5} + A^{3} + A$} \\ \textbf{Arrow:} {\scriptsize $A^{22}L_1 - 3A^{18}L_1 - 2A^{16}L_2 + 4A^{14}L_1 + 4A^{12}L_2 + A^{12} - 3A^{10}L_1 - 5A^{8}L_2 - A^{8} + 4A^{4}L_2 + A^{4} + 2A^{2}L_1 - L_2 - L_1/A^{2}$} \\ \textbf{Mock:} {\scriptsize $-w^{5} - w^{4} + 4w^{3} + 5w^{2} - 5w - 7 + 1/w + 5/w^{2} + 2/w^{3} - 1/w^{4} - 1/w^{5}$} \\ \textbf{Affine:} {\scriptsize $-2t^{2} + 4 - 2/t^{2}$} \\ \textbf{Yamada:} {\scriptsize $3A^{27} - 8A^{25} + 7A^{24} + 8A^{23} - 17A^{22} + 2A^{21} + 16A^{20} - 15A^{19} - 4A^{18} + 13A^{17} - 4A^{16} - 6A^{15} + A^{14} + 10A^{13} - 8A^{12} - 9A^{11} + 18A^{10} - 5A^{9} - 16A^{8} + 13A^{7} + 3A^{6} - 12A^{5} + 2A^{4} + 5A^{3} - 3A^{2} - A + 1$}
\end{minipage}

\noindent{\color{gray!40}\rule{\textwidth}{0.4pt}}
\vspace{0.9\baselineskip}
\noindent \begin{minipage}[t]{0.25\textwidth}
\vspace{0pt}
\centering
\includegraphics[page=356,width=\linewidth]{knotoids.pdf}
\end{minipage}
\hfill
\begin{minipage}[t]{0.73\textwidth}
\vspace{0pt}
\raggedright
\textbf{Name:} {\large{$\mathbf{K7_{237}}$}} (chiral, non-rotatable$^{*}$) \\ \textbf{PD:} {\scriptsize\texttt{[0],[0,1,2,3],[1,4,5,2],[3,6,7,8],[4,8,9,10],[10,11,12,5],[6,12,13,7],[9,13,11,14],[14]}} \\ \textbf{EM:} {\scriptsize\texttt{(B0, A0C0C3D0, B1E0F3B2, B3G0G3E1, C1D3H0F0, E3H2G1C2, D1F2H1D2, E2G2F1I0, H3)}} \\ \textbf{Kauffman bracket:} {\scriptsize $A^{22} - 3A^{18} - 2A^{16} + 5A^{14} + 5A^{12} - 2A^{10} - 5A^{8} + 3A^{4} - 1$} \\ \textbf{Arrow:} {\scriptsize $A^{4} - L_2 - 2 - 2L_1/A^{2} + 2L_2/A^{4} + 3/A^{4} + 5L_1/A^{6} - L_2/A^{8} - 1/A^{8} - 5L_1/A^{10} + 3L_1/A^{14} - L_1/A^{18}$} \\ \textbf{Mock:} {\scriptsize $-w^{4} - 2w^{3} + 3w^{2} + 5w - 3 - 5/w + 3/w^{2} + 3/w^{3} - 1/w^{4} - 1/w^{5}$} \\ \textbf{Affine:} {\scriptsize $0$} \\ \textbf{Yamada:} {\scriptsize $-A^{26} - A^{25} + 3A^{24} + A^{23} - 6A^{22} + 2A^{21} + 9A^{20} - 6A^{19} - 6A^{18} + 12A^{17} - A^{16} - 10A^{15} + 8A^{14} + 3A^{13} - 7A^{12} + 2A^{11} + 6A^{10} + 2A^{9} - 8A^{8} + 7A^{7} + 7A^{6} - 14A^{5} + 2A^{4} + 8A^{3} - 7A^{2} - 2A + 3$}
\end{minipage}

\noindent{\color{gray!40}\rule{\textwidth}{0.4pt}}
\vspace{0.9\baselineskip}
\noindent \begin{minipage}[t]{0.25\textwidth}
\vspace{0pt}
\centering
\includegraphics[page=357,width=\linewidth]{knotoids.pdf}
\end{minipage}
\hfill
\begin{minipage}[t]{0.73\textwidth}
\vspace{0pt}
\raggedright
\textbf{Name:} {\large{$\mathbf{K7_{238}}$}} (chiral, non-rotatable$^{*}$) \\ \textbf{PD:} {\scriptsize\texttt{[0],[0,1,2,3],[1,4,5,2],[3,5,6,7],[7,8,9,4],[6,9,10,11],[8,12,13,10],[14,13,12,11],[14]}} \\ \textbf{EM:} {\scriptsize\texttt{(B0, A0C0C3D0, B1E3D1B2, B3C2F0E0, D3G0F1C1, D2E2G3H3, E1H2H1F2, I0G2G1F3, H0)}} \\ \textbf{Kauffman bracket:} {\scriptsize $2A^{22} + 2A^{20} - 2A^{18} - 5A^{16} + 5A^{12} + 2A^{10} - 3A^{8} - 2A^{6} + A^{4} + A^{2}$} \\ \textbf{Arrow:} {\scriptsize $A^{4}L_2 + A^{4} + 2A^{2}L_1 - L_2 - 1 - 5L_1/A^{2} - L_2/A^{4} + A^{-4} + 5L_1/A^{6} + 2L_2/A^{8} - 3L_1/A^{10} - 2L_2/A^{12} + L_1/A^{14} + L_2/A^{16}$} \\ \textbf{Mock:} {\scriptsize $w^{5} - 3w^{3} + 5w + 1 - 5/w + 2/w^{3}$} \\ \textbf{Affine:} {\scriptsize $0$} \\ \textbf{Yamada:} {\scriptsize $-A^{28} - A^{27} + A^{26} + 2A^{25} - 3A^{24} - 2A^{23} + 7A^{22} - 9A^{20} + 7A^{19} + 5A^{18} - 12A^{17} + 2A^{16} + 5A^{15} - 7A^{14} - 2A^{13} + 2A^{12} + 3A^{11} - 8A^{10} + 11A^{8} - 8A^{7} - 3A^{6} + 10A^{5} - 2A^{4} - 5A^{3} + 3A^{2} + A - 2$}
\end{minipage}

\noindent{\color{gray!40}\rule{\textwidth}{0.4pt}}
\vspace{0.9\baselineskip}
\noindent \begin{minipage}[t]{0.25\textwidth}
\vspace{0pt}
\centering
\includegraphics[page=358,width=\linewidth]{knotoids.pdf}
\end{minipage}
\hfill
\begin{minipage}[t]{0.73\textwidth}
\vspace{0pt}
\raggedright
\textbf{Name:} {\large{$\mathbf{K7_{239}}$}} (chiral, non-rotatable$^{*}$) \\ \textbf{PD:} {\scriptsize\texttt{[0],[0,1,2,3],[1,4,5,2],[3,6,7,8],[8,9,10,4],[10,11,12,5],[6,12,13,7],[9,13,11,14],[14]}} \\ \textbf{EM:} {\scriptsize\texttt{(B0, A0C0C3D0, B1E3F3B2, B3G0G3E0, D3H0F0C1, E2H2G1C2, D1F2H1D2, E1G2F1I0, H3)}} \\ \textbf{Kauffman bracket:} {\scriptsize $-A^{20} + 2A^{16} + A^{14} - 3A^{12} - 3A^{10} + 2A^{8} + 5A^{6} + A^{4} - 2A^{2} - 1$} \\ \textbf{Arrow:} {\scriptsize $-A^{14}L_1 + 2A^{10}L_1 + A^{8} - 3A^{6}L_1 - A^{4}L_2 - 2A^{4} + 2A^{2}L_1 + 2L_2 + 3 + L_1/A^{2} - L_2/A^{4} - 1/A^{4} - L_1/A^{6}$} \\ \textbf{Mock:} {\scriptsize $-2w^{3} - 2w^{2} + 3w + 5 - 1/w - 3/w^{2} + w^{-4}$} \\ \textbf{Affine:} {\scriptsize $-t + 2 - 1/t$} \\ \textbf{Yamada:} {\scriptsize $A^{26} - 2A^{24} + 4A^{22} - 2A^{21} - 5A^{20} + 7A^{19} + A^{18} - 8A^{17} + 5A^{16} + 2A^{15} - 6A^{14} - 6A^{10} + 2A^{9} + 5A^{8} - 8A^{7} + 7A^{5} - 4A^{4} - 2A^{3} + 4A^{2} - 1$}
\end{minipage}

\noindent{\color{gray!40}\rule{\textwidth}{0.4pt}}
\vspace{0.9\baselineskip}
\noindent \begin{minipage}[t]{0.25\textwidth}
\vspace{0pt}
\centering
\includegraphics[page=359,width=\linewidth]{knotoids.pdf}
\end{minipage}
\hfill
\begin{minipage}[t]{0.73\textwidth}
\vspace{0pt}
\raggedright
\textbf{Name:} {\large{$\mathbf{K7_{240}}$}} (chiral, non-rotatable$^{*}$) \\ \textbf{PD:} {\scriptsize\texttt{[0],[0,1,2,3],[1,4,5,2],[6,7,4,3],[5,8,9,6],[7,9,10,11],[8,12,13,10],[14,13,12,11],[14]}} \\ \textbf{EM:} {\scriptsize\texttt{(B0, A0C0C3D3, B1D2E0B2, E3F0C1B3, C2G0F1D0, D1E2G3H3, E1H2H1F2, I0G2G1F3, H0)}} \\ \textbf{Kauffman bracket:} {\scriptsize $A^{24} + 2A^{22} - 4A^{18} - 4A^{16} + 2A^{14} + 6A^{12} + A^{10} - 3A^{8} - 2A^{6} + A^{4} + A^{2}$} \\ \textbf{Arrow:} {\scriptsize $A^{12} + 2A^{10}L_1 + A^{8}L_2 - A^{8} - 4A^{6}L_1 - 4A^{4}L_2 + 2A^{2}L_1 + 5L_2 + 1 + L_1/A^{2} - 3L_2/A^{4} - 2L_1/A^{6} + L_2/A^{8} + L_1/A^{10}$} \\ \textbf{Mock:} {\scriptsize $w^{5} - 3w^{3} - 2w^{2} + 3w + 6 - 4/w^{2} - 2/w^{3} + w^{-4} + w^{-5}$} \\ \textbf{Affine:} {\scriptsize $-t^{2} + 2 - 1/t^{2}$} \\ \textbf{Yamada:} {\scriptsize $-A^{28} + 3A^{26} - 6A^{24} + 3A^{23} + 7A^{22} - 8A^{21} - 4A^{20} + 12A^{19} - 2A^{18} - 11A^{17} + 8A^{16} + A^{15} - 10A^{14} + A^{13} + 4A^{12} - A^{11} - 8A^{10} + 9A^{9} + 6A^{8} - 14A^{7} + 4A^{6} + 8A^{5} - 7A^{4} - 2A^{3} + 4A^{2} - 2$}
\end{minipage}

\noindent{\color{gray!40}\rule{\textwidth}{0.4pt}}
\vspace{0.9\baselineskip}
\noindent \begin{minipage}[t]{0.25\textwidth}
\vspace{0pt}
\centering
\includegraphics[page=360,width=\linewidth]{knotoids.pdf}
\end{minipage}
\hfill
\begin{minipage}[t]{0.73\textwidth}
\vspace{0pt}
\raggedright
\textbf{Name:} {\large{$\mathbf{K7_{241}}$}} (chiral, non-rotatable$^{*}$) \\ \textbf{PD:} {\scriptsize\texttt{[0],[0,1,2,3],[1,4,5,2],[6,7,4,3],[5,8,9,10],[11,12,13,6],[7,13,14,8],[9],[10,14,12,11]}} \\ \textbf{EM:} {\scriptsize\texttt{(B0, A0C0C3D3, B1D2E0B2, F3G0C1B3, C2G3H0I0, I3I2G1D0, D1F2I1E1, E2, E3G2F1F0)}} \\ \textbf{Kauffman bracket:} {\scriptsize $-A^{20} - 2A^{18} + A^{16} + 4A^{14} + 2A^{12} - 3A^{10} - 2A^{8} + 2A^{6} + 2A^{4} - A^{2} - 1$} \\ \textbf{Arrow:} {\scriptsize $-A^{8} - 2A^{6}L_1 - A^{4}L_2 + 2A^{4} + 4A^{2}L_1 + 3L_2 - 1 - 3L_1/A^{2} - 3L_2/A^{4} + A^{-4} + 2L_1/A^{6} + 2L_2/A^{8} - L_1/A^{10} - L_2/A^{12}$} \\ \textbf{Mock:} {\scriptsize $w^{6} + w^{5} - 2w^{4} - 4w^{3} + 5w + 4 - 2/w - 3/w^{2} + w^{-4}$} \\ \textbf{Affine:} {\scriptsize $t^{2} - 2 + t^{-2}$} \\ \textbf{Yamada:} {\scriptsize $-A^{27} + A^{26} + 2A^{25} - 3A^{24} + 5A^{22} - 4A^{21} - 3A^{20} + 6A^{19} - 2A^{18} - 4A^{17} + 4A^{16} - A^{15} - 4A^{14} + A^{11} - 6A^{10} + 3A^{9} + 3A^{8} - 8A^{7} + 3A^{6} + 3A^{5} - 4A^{4} + 2A^{3} + A^{2} - A + 1$}
\end{minipage}

\noindent{\color{gray!40}\rule{\textwidth}{0.4pt}}
\vspace{0.9\baselineskip}
\noindent \begin{minipage}[t]{0.25\textwidth}
\vspace{0pt}
\centering
\includegraphics[page=361,width=\linewidth]{knotoids.pdf}
\end{minipage}
\hfill
\begin{minipage}[t]{0.73\textwidth}
\vspace{0pt}
\raggedright
\textbf{Name:} {\large{$\mathbf{K7_{242}}$}} (chiral, non-rotatable$^{*}$) \\ \textbf{PD:} {\scriptsize\texttt{[0],[0,1,2,3],[1,4,5,2],[6,7,4,3],[8,9,10,5],[11,12,13,6],[7,13,14,8],[9],[10,14,12,11]}} \\ \textbf{EM:} {\scriptsize\texttt{(B0, A0C0C3D3, B1D2E3B2, F3G0C1B3, G3H0I0C2, I3I2G1D0, D1F2I1E0, E1, E2G2F1F0)}} \\ \textbf{Kauffman bracket:} {\scriptsize $-A^{26} - 2A^{24} + 5A^{20} + 3A^{18} - 4A^{16} - 4A^{14} + 3A^{12} + 4A^{10} - A^{8} - 3A^{6} + A^{2}$} \\ \textbf{Arrow:} {\scriptsize $-A^{14}L_1 - A^{12}L_2 - A^{12} + 3A^{8}L_2 + 2A^{8} + 3A^{6}L_1 - 3A^{4}L_2 - A^{4} - 4A^{2}L_1 + 2L_2 + 1 + 4L_1/A^{2} - L_2/A^{4} - 3L_1/A^{6} + L_1/A^{10}$} \\ \textbf{Mock:} {\scriptsize $w^{4} - 4w^{2} - 3w + 6 + 6/w - 3/w^{2} - 4/w^{3} + w^{-4} + w^{-5}$} \\ \textbf{Affine:} {\scriptsize $t^{2} - t - 1/t + t^{-2}$} \\ \textbf{Yamada:} {\scriptsize $A^{28} - A^{27} - 2A^{26} + 4A^{25} - 8A^{23} + 7A^{22} + 6A^{21} - 14A^{20} + 3A^{19} + 11A^{18} - 11A^{17} - A^{16} + 9A^{15} - 2A^{14} - 3A^{13} + 10A^{11} - 6A^{10} - 6A^{9} + 16A^{8} - 4A^{7} - 9A^{6} + 11A^{5} - 7A^{3} + 2A^{2} + A - 1$}
\end{minipage}

\noindent{\color{gray!40}\rule{\textwidth}{0.4pt}}
\vspace{0.9\baselineskip}
\noindent \begin{minipage}[t]{0.25\textwidth}
\vspace{0pt}
\centering
\includegraphics[page=362,width=\linewidth]{knotoids.pdf}
\end{minipage}
\hfill
\begin{minipage}[t]{0.73\textwidth}
\vspace{0pt}
\raggedright
\textbf{Name:} {\large{$\mathbf{K7_{243}}$}} (chiral, non-rotatable$^{*}$) \\ \textbf{PD:} {\scriptsize\texttt{[0],[0,1,2,3],[1,4,5,2],[6,7,8,3],[4,8,9,10],[5,10,11,12],[12,13,7,6],[13,11,14,9],[14]}} \\ \textbf{EM:} {\scriptsize\texttt{(B0, A0C0C3D3, B1E0F0B2, G3G2E1B3, C1D2H3F1, C2E3H1G0, F3H0D1D0, G1F2I0E2, H2)}} \\ \textbf{Kauffman bracket:} {\scriptsize $A^{17} + A^{15} - 2A^{13} - 3A^{11} + 4A^{7} + 2A^{5} - 2A^{3} - 2A$} \\ \textbf{Arrow:} {\scriptsize $-A^{20} - A^{18}L_1 + 2A^{16} + 3A^{14}L_1 + A^{12}L_2 - A^{12} - 4A^{10}L_1 - 2A^{8}L_2 + 2A^{6}L_1 + A^{4}L_2 + A^{4}$} \\ \textbf{Mock:} {\scriptsize $w^{4} + 2w^{3} - w^{2} - 5w - 1 + 3/w + 2/w^{2}$} \\ \textbf{Affine:} {\scriptsize $-t + 2 - 1/t$} \\ \textbf{Yamada:} {\scriptsize $A^{26} - A^{25} - 2A^{24} + 2A^{23} + A^{22} - 3A^{21} + 5A^{19} - 2A^{18} - 3A^{17} + 6A^{16} - A^{15} - 3A^{14} + 2A^{13} - 2A^{11} - 4A^{10} + 3A^{9} - A^{8} - 6A^{7} + 3A^{6} + 2A^{5} - 4A^{4} + A^{3} + 2A^{2} - A - 1$}
\end{minipage}

\noindent{\color{gray!40}\rule{\textwidth}{0.4pt}}
\vspace{0.9\baselineskip}
\noindent \begin{minipage}[t]{0.25\textwidth}
\vspace{0pt}
\centering
\includegraphics[page=363,width=\linewidth]{knotoids.pdf}
\end{minipage}
\hfill
\begin{minipage}[t]{0.73\textwidth}
\vspace{0pt}
\raggedright
\textbf{Name:} {\large{$\mathbf{K7_{244}}$}} (chiral, non-rotatable$^{*}$) \\ \textbf{PD:} {\scriptsize\texttt{[0],[0,1,2,3],[1,4,5,6],[2,6,7,8],[9,10,4,3],[11,12,13,5],[12,14,8,7],[14,11,10,9],[13]}} \\ \textbf{EM:} {\scriptsize\texttt{(B0, A0C0D0E3, B1E2F3D1, B2C3G3G2, H3H2C1B3, H1G0I0C2, F1H0D3D2, G1F0E1E0, F2)}} \\ \textbf{Kauffman bracket:} {\scriptsize $A^{22} + A^{20} - 2A^{18} - 3A^{16} + 3A^{14} + 5A^{12} - 4A^{8} - A^{6} + 2A^{4} - 1$} \\ \textbf{Arrow:} {\scriptsize $A^{-8} + L_1/A^{10} - 2/A^{12} - 3L_1/A^{14} + 3/A^{16} + 5L_1/A^{18} + L_2/A^{20} - 1/A^{20} - 4L_1/A^{22} - L_2/A^{24} + 2L_1/A^{26} - L_1/A^{30}$} \\ \textbf{Mock:} {\scriptsize $3w^{4} + 3w^{3} - 3w^{2} - 6w + 5/w + w^{-2} - 2/w^{3}$} \\ \textbf{Affine:} {\scriptsize $t^{2} + 2t - 6 + 2/t + t^{-2}$} \\ \textbf{Yamada:} {\scriptsize $-2A^{24} + A^{23} + 2A^{22} - 6A^{21} - 2A^{20} + 6A^{19} - 6A^{18} - 5A^{17} + 8A^{16} - 2A^{15} - 4A^{14} + 4A^{13} - A^{11} - 5A^{10} + 5A^{9} + A^{8} - 8A^{7} + 7A^{6} + 4A^{5} - 6A^{4} + 3A^{3} + 4A^{2} - 3A - 1$}
\end{minipage}

\noindent{\color{gray!40}\rule{\textwidth}{0.4pt}}
\vspace{0.9\baselineskip}
\noindent \begin{minipage}[t]{0.25\textwidth}
\vspace{0pt}
\centering
\includegraphics[page=364,width=\linewidth]{knotoids.pdf}
\end{minipage}
\hfill
\begin{minipage}[t]{0.73\textwidth}
\vspace{0pt}
\raggedright
\textbf{Name:} {\large{$\mathbf{K7_{245}}$}} (chiral, non-rotatable$^{*}$) \\ \textbf{PD:} {\scriptsize\texttt{[0],[0,1,2,3],[1,4,5,2],[3,6,7,4],[5,8,9,6],[7,10,11,8],[12,13,10,9],[13,12,14,11],[14]}} \\ \textbf{EM:} {\scriptsize\texttt{(B0, A0C0C3D0, B1D3E0B2, B3E3F0C1, C2F3G3D1, D2G2H3E1, H1H0F1E2, G1G0I0F2, H2)}} \\ \textbf{Kauffman bracket:} {\scriptsize $A^{21} + A^{19} - 2A^{17} - 3A^{15} + A^{13} + 4A^{11} - 3A^{7} - 2A^{5} + A^{3} + A$} \\ \textbf{Arrow:} {\scriptsize $-A^{12} - A^{10}L_1 + 2A^{8} + 3A^{6}L_1 + A^{4}L_2 - 2A^{4} - 4A^{2}L_1 - 2L_2 + 2 + 3L_1/A^{2} + 2L_2/A^{4} - L_1/A^{6} - L_2/A^{8}$} \\ \textbf{Mock:} {\scriptsize $-w^{5} + 3w^{3} - 4w + 1 + 3/w - 1/w^{3}$} \\ \textbf{Affine:} {\scriptsize $t^{2} - 2 + t^{-2}$} \\ \textbf{Yamada:} {\scriptsize $A^{26} - A^{25} - A^{24} + 3A^{23} - 4A^{22} - A^{21} + 7A^{20} - 6A^{19} - A^{18} + 6A^{17} - 3A^{16} - A^{15} + 2A^{13} - 3A^{12} - 4A^{11} + 4A^{10} - A^{9} - 6A^{8} + 4A^{7} + A^{6} - 6A^{5} + 2A^{4} + 3A^{3} - 3A^{2} + 2$}
\end{minipage}

\noindent{\color{gray!40}\rule{\textwidth}{0.4pt}}
\vspace{0.9\baselineskip}
\noindent \begin{minipage}[t]{0.25\textwidth}
\vspace{0pt}
\centering
\includegraphics[page=365,width=\linewidth]{knotoids.pdf}
\end{minipage}
\hfill
\begin{minipage}[t]{0.73\textwidth}
\vspace{0pt}
\raggedright
\textbf{Name:} {\large{$\mathbf{K7_{246}}$}} (chiral, non-rotatable$^{*}$) \\ \textbf{PD:} {\scriptsize\texttt{[0],[0,1,2,3],[1,4,5,2],[3,6,7,4],[5,8,9,6],[7,10,11,12],[8,13,10,9],[13,12,14,11],[14]}} \\ \textbf{EM:} {\scriptsize\texttt{(B0, A0C0C3D0, B1D3E0B2, B3E3F0C1, C2G0G3D1, D2G2H3H1, E1H0F1E2, G1F3I0F2, H2)}} \\ \textbf{Kauffman bracket:} {\scriptsize $A^{22} + 2A^{20} - 2A^{18} - 4A^{16} + A^{14} + 6A^{12} + 2A^{10} - 4A^{8} - 3A^{6} + A^{4} + A^{2}$} \\ \textbf{Arrow:} {\scriptsize $L_1/A^{2} + 2/A^{4} - 2L_1/A^{6} - 4/A^{8} + L_1/A^{10} + 6/A^{12} + 2L_1/A^{14} - 4/A^{16} - 3L_1/A^{18} + A^{-20} + L_1/A^{22}$} \\ \textbf{Mock:} {\scriptsize $-w^{4} + 5w^{2} + w - 6 - 1/w + 4/w^{2} - 1/w^{4}$} \\ \textbf{Affine:} {\scriptsize $t - 2 + 1/t$} \\ \textbf{Yamada:} {\scriptsize $-A^{26} + 2A^{25} + 3A^{24} - 6A^{23} + A^{22} + 11A^{21} - 9A^{20} - 3A^{19} + 13A^{18} - 6A^{17} - 4A^{16} + 8A^{15} - 2A^{13} - 3A^{12} + 7A^{11} - A^{10} - 11A^{9} + 10A^{8} + 2A^{7} - 12A^{6} + 6A^{5} + 5A^{4} - 6A^{3} + A^{2} + 2A - 1$}
\end{minipage}

\noindent{\color{gray!40}\rule{\textwidth}{0.4pt}}
\vspace{0.9\baselineskip}
\noindent \begin{minipage}[t]{0.25\textwidth}
\vspace{0pt}
\centering
\includegraphics[page=366,width=\linewidth]{knotoids.pdf}
\end{minipage}
\hfill
\begin{minipage}[t]{0.73\textwidth}
\vspace{0pt}
\raggedright
\textbf{Name:} {\large{$\mathbf{K7_{247}}$}} (chiral, non-rotatable$^{*}$) \\ \textbf{PD:} {\scriptsize\texttt{[0],[0,1,2,3],[1,4,5,2],[3,6,7,4],[5,8,9,6],[10,11,8,7],[9,12,13,10],[11,13,12,14],[14]}} \\ \textbf{EM:} {\scriptsize\texttt{(B0, A0C0C3D0, B1D3E0B2, B3E3F3C1, C2F2G0D1, G3H0E1D2, E2H2H1F0, F1G2G1I0, H3)}} \\ \textbf{Kauffman bracket:} {\scriptsize $-A^{18} + A^{16} + 3A^{14} + A^{12} - 4A^{10} - 2A^{8} + 3A^{6} + 3A^{4} - A^{2} - 2$} \\ \textbf{Arrow:} {\scriptsize $-L_1/A^{6} + A^{-8} + 3L_1/A^{10} + L_2/A^{12} - 4L_1/A^{14} - 2L_2/A^{16} + 3L_1/A^{18} + 2L_2/A^{20} + A^{-20} - L_1/A^{22} - L_2/A^{24} - 1/A^{24}$} \\ \textbf{Mock:} {\scriptsize $-w^{5} + 3w^{3} + 2w^{2} - 4w - 3 + 3/w + 2/w^{2} - 1/w^{3}$} \\ \textbf{Affine:} {\scriptsize $t^{2} - 2 + t^{-2}$} \\ \textbf{Yamada:} {\scriptsize $A^{27} - A^{25} + 3A^{24} + A^{23} - 5A^{22} + 5A^{21} + 5A^{20} - 6A^{19} + 3A^{18} + 5A^{17} - 4A^{16} + A^{14} + 3A^{13} - 4A^{12} - A^{11} + 6A^{10} - 6A^{9} - 3A^{8} + 6A^{7} - 2A^{6} - 4A^{5} + 3A^{4} + A^{3} - 2A^{2} + 1$}
\end{minipage}

\noindent{\color{gray!40}\rule{\textwidth}{0.4pt}}
\vspace{0.9\baselineskip}
\noindent \begin{minipage}[t]{0.25\textwidth}
\vspace{0pt}
\centering
\includegraphics[page=367,width=\linewidth]{knotoids.pdf}
\end{minipage}
\hfill
\begin{minipage}[t]{0.73\textwidth}
\vspace{0pt}
\raggedright
\textbf{Name:} {\large{$\mathbf{K7_{248}}$}} (chiral, non-rotatable$^{*}$) \\ \textbf{PD:} {\scriptsize\texttt{[0],[0,1,2,3],[1,4,5,2],[3,6,7,4],[8,9,6,5],[7,10,11,8],[12,13,10,9],[13,12,14,11],[14]}} \\ \textbf{EM:} {\scriptsize\texttt{(B0, A0C0C3D0, B1D3E3B2, B3E2F0C1, F3G3D1C2, D2G2H3E0, H1H0F1E1, G1G0I0F2, H2)}} \\ \textbf{Kauffman bracket:} {\scriptsize $-A^{20} - A^{18} + 2A^{16} + 3A^{14} - A^{12} - 5A^{10} - A^{8} + 4A^{6} + 3A^{4} - A^{2} - 1$} \\ \textbf{Arrow:} {\scriptsize $-A^{20} - A^{18}L_1 + 2A^{16} + 3A^{14}L_1 + A^{12}L_2 - 2A^{12} - 5A^{10}L_1 - 2A^{8}L_2 + A^{8} + 4A^{6}L_1 + 2A^{4}L_2 + A^{4} - A^{2}L_1 - L_2$} \\ \textbf{Mock:} {\scriptsize $-w^{5} + 3w^{3} - 5w + 4/w + w^{-2} - 1/w^{3}$} \\ \textbf{Affine:} {\scriptsize $t^{2} - t - 1/t + t^{-2}$} \\ \textbf{Yamada:} {\scriptsize $A^{28} - A^{27} - 2A^{26} + 4A^{25} - 7A^{23} + 4A^{22} + 5A^{21} - 8A^{20} + 9A^{18} - 4A^{17} - 2A^{16} + 5A^{15} - 3A^{13} - 3A^{12} + 7A^{11} - 5A^{10} - 8A^{9} + 8A^{8} - 3A^{7} - 8A^{6} + 4A^{5} + A^{4} - 2A^{3} + A + 1$}
\end{minipage}

\noindent{\color{gray!40}\rule{\textwidth}{0.4pt}}
\vspace{0.9\baselineskip}
\noindent \begin{minipage}[t]{0.25\textwidth}
\vspace{0pt}
\centering
\includegraphics[page=368,width=\linewidth]{knotoids.pdf}
\end{minipage}
\hfill
\begin{minipage}[t]{0.73\textwidth}
\vspace{0pt}
\raggedright
\textbf{Name:} {\large{$\mathbf{K7_{249}}$}} (chiral, non-rotatable$^{*}$) \\ \textbf{PD:} {\scriptsize\texttt{[0],[0,1,2,3],[1,4,5,2],[3,6,7,4],[8,9,6,5],[7,10,11,12],[13,10,9,8],[13,12,14,11],[14]}} \\ \textbf{EM:} {\scriptsize\texttt{(B0, A0C0C3D0, B1D3E3B2, B3E2F0C1, G3G2D1C2, D2G1H3H1, H0F1E1E0, G0F3I0F2, H2)}} \\ \textbf{Kauffman bracket:} {\scriptsize $-A^{20} - A^{18} + 2A^{16} + 4A^{14} - A^{12} - 5A^{10} - 2A^{8} + 4A^{6} + 3A^{4} - A^{2} - 1$} \\ \textbf{Arrow:} {\scriptsize $-A^{14}L_1 - A^{12} + 2A^{10}L_1 + 4A^{8} - A^{6}L_1 - 5A^{4} - 2A^{2}L_1 + 4 + 3L_1/A^{2} - 1/A^{4} - L_1/A^{6}$} \\ \textbf{Mock:} {\scriptsize $w^{4} - 4w^{2} - w + 6 + 1/w - 3/w^{2} + w^{-4}$} \\ \textbf{Affine:} {\scriptsize $-t + 2 - 1/t$} \\ \textbf{Yamada:} {\scriptsize $A^{26} - A^{25} - 3A^{24} + 4A^{23} + 3A^{22} - 8A^{21} + 3A^{20} + 9A^{19} - 9A^{18} - 2A^{17} + 8A^{16} - 6A^{15} - 4A^{14} + 4A^{13} - 3A^{11} - 4A^{10} + 8A^{9} - 2A^{8} - 10A^{7} + 10A^{6} - 9A^{4} + 5A^{3} + A^{2} - 2A + 1$}
\end{minipage}

\noindent{\color{gray!40}\rule{\textwidth}{0.4pt}}
\vspace{0.9\baselineskip}
\noindent \begin{minipage}[t]{0.25\textwidth}
\vspace{0pt}
\centering
\includegraphics[page=369,width=\linewidth]{knotoids.pdf}
\end{minipage}
\hfill
\begin{minipage}[t]{0.73\textwidth}
\vspace{0pt}
\raggedright
\textbf{Name:} {\large{$\mathbf{K7_{250}}$}} (chiral, non-rotatable$^{*}$) \\ \textbf{PD:} {\scriptsize\texttt{[0],[0,1,2,3],[1,4,5,2],[3,6,7,4],[8,9,6,5],[10,11,8,7],[9,12,13,10],[11,13,12,14],[14]}} \\ \textbf{EM:} {\scriptsize\texttt{(B0, A0C0C3D0, B1D3E3B2, B3E2F3C1, F2G0D1C2, G3H0E0D2, E1H2H1F0, F1G2G1I0, H3)}} \\ \textbf{Kauffman bracket:} {\scriptsize $-A^{26} + 2A^{22} + 3A^{20} - 2A^{18} - 5A^{16} + 6A^{12} + 3A^{10} - 3A^{8} - 3A^{6} + A^{2}$} \\ \textbf{Arrow:} {\scriptsize $-A^{2}L_1 + 2L_1/A^{2} + L_2/A^{4} + 2/A^{4} - 2L_1/A^{6} - 2L_2/A^{8} - 3/A^{8} + 2L_2/A^{12} + 4/A^{12} + 3L_1/A^{14} - L_2/A^{16} - 2/A^{16} - 3L_1/A^{18} + L_1/A^{22}$} \\ \textbf{Mock:} {\scriptsize $-w^{5} - 2w^{4} + w^{3} + 6w^{2} + w - 6 - 2/w + 3/w^{2} + w^{-3}$} \\ \textbf{Affine:} {\scriptsize $t^{2} - t - 1/t + t^{-2}$} \\ \textbf{Yamada:} {\scriptsize $-A^{27} + 4A^{25} - 3A^{24} - 7A^{23} + 10A^{22} + 3A^{21} - 16A^{20} + 7A^{19} + 9A^{18} - 15A^{17} + 8A^{15} - 6A^{14} - 4A^{13} + 2A^{12} + 8A^{11} - 10A^{10} - 4A^{9} + 17A^{8} - 9A^{7} - 8A^{6} + 14A^{5} - A^{4} - 7A^{3} + 3A^{2} + A - 1$}
\end{minipage}

\noindent{\color{gray!40}\rule{\textwidth}{0.4pt}}
\vspace{0.9\baselineskip}
\noindent \begin{minipage}[t]{0.25\textwidth}
\vspace{0pt}
\centering
\includegraphics[page=370,width=\linewidth]{knotoids.pdf}
\end{minipage}
\hfill
\begin{minipage}[t]{0.73\textwidth}
\vspace{0pt}
\raggedright
\textbf{Name:} {\large{$\mathbf{K7_{251}}$}} (chiral, non-rotatable$^{*}$) \\ \textbf{PD:} {\scriptsize\texttt{[0],[0,1,2,3],[1,4,5,2],[3,5,6,7],[4,7,8,9],[9,10,11,6],[8,11,12,13],[13,12,14,10],[14]}} \\ \textbf{EM:} {\scriptsize\texttt{(B0, A0C0C3D0, B1E0D1B2, B3C2F3E1, C1D3G0F0, E3H3G1D2, E2F2H1H0, G3G2I0F1, H2)}} \\ \textbf{Kauffman bracket:} {\scriptsize $A^{21} + A^{19} - 3A^{17} - 3A^{15} + 2A^{13} + 5A^{11} - 4A^{7} - 2A^{5} + A^{3} + A$} \\ \textbf{Arrow:} {\scriptsize $-A^{12}L_2 - A^{10}L_1 + 2A^{8}L_2 + A^{8} + 3A^{6}L_1 - 2A^{4}L_2 - 5A^{2}L_1 + L_2 - 1 + 4L_1/A^{2} + 2/A^{4} - L_1/A^{6} - 1/A^{8}$} \\ \textbf{Mock:} {\scriptsize $w^{4} + w^{3} - 3w^{2} - 4w + 3 + 5/w + w^{-2} - 3/w^{3} - 2/w^{4} + w^{-5} + w^{-6}$} \\ \textbf{Affine:} {\scriptsize $-t^{2} - t + 4 - 1/t - 1/t^{2}$} \\ \textbf{Yamada:} {\scriptsize $A^{28} + A^{27} - A^{26} - 2A^{25} + 3A^{24} + 2A^{23} - 8A^{22} + 2A^{21} + 8A^{20} - 9A^{19} - 2A^{18} + 7A^{17} - 4A^{16} - 3A^{15} + A^{14} + 4A^{13} - 5A^{12} - 3A^{11} + 9A^{10} - 4A^{9} - 7A^{8} + 7A^{7} + A^{6} - 7A^{5} + A^{4} + 4A^{3} - 2A^{2} - A + 1$}
\end{minipage}

\noindent{\color{gray!40}\rule{\textwidth}{0.4pt}}
\vspace{0.9\baselineskip}
\noindent \begin{minipage}[t]{0.25\textwidth}
\vspace{0pt}
\centering
\includegraphics[page=371,width=\linewidth]{knotoids.pdf}
\end{minipage}
\hfill
\begin{minipage}[t]{0.73\textwidth}
\vspace{0pt}
\raggedright
\textbf{Name:} {\large{$\mathbf{K7_{252}}$}} (chiral, non-rotatable$^{*}$) \\ \textbf{PD:} {\scriptsize\texttt{[0],[0,1,2,3],[1,4,5,2],[3,6,7,8],[4,8,9,10],[11,12,6,5],[13,14,9,7],[10,13,12,11],[14]}} \\ \textbf{EM:} {\scriptsize\texttt{(B0, A0C0C3D0, B1E0F3B2, B3F2G3E1, C1D3G2H0, H3H2D1C2, H1I0E2D2, E3G0F1F0, G1)}} \\ \textbf{Kauffman bracket:} {\scriptsize $-A^{24} + 3A^{20} + A^{18} - 5A^{16} - 4A^{14} + 3A^{12} + 5A^{10} - 2A^{6} + A^{2}$} \\ \textbf{Arrow:} {\scriptsize $-A^{30}L_1 + 3A^{26}L_1 + A^{24} - 5A^{22}L_1 - 4A^{20} + 3A^{18}L_1 + 5A^{16} - 2A^{12} + A^{8}$} \\ \textbf{Mock:} {\scriptsize $w^{4} - w^{3} - 2w^{2} + 3w + 3 - 5/w - 4/w^{2} + 3/w^{3} + 3/w^{4}$} \\ \textbf{Affine:} {\scriptsize $-2t + 4 - 2/t$} \\ \textbf{Yamada:} {\scriptsize $A^{26} - 3A^{24} + 2A^{23} + 4A^{22} - 7A^{21} - 3A^{20} + 10A^{19} - 3A^{18} - 9A^{17} + 10A^{16} + A^{15} - 6A^{14} + 4A^{13} + 4A^{12} - 6A^{10} + 7A^{9} + 3A^{8} - 12A^{7} + 4A^{6} + 6A^{5} - 10A^{4} - 2A^{3} + 4A^{2} - 2A - 3$}
\end{minipage}

\noindent{\color{gray!40}\rule{\textwidth}{0.4pt}}
\vspace{0.9\baselineskip}
\noindent \begin{minipage}[t]{0.25\textwidth}
\vspace{0pt}
\centering
\includegraphics[page=372,width=\linewidth]{knotoids.pdf}
\end{minipage}
\hfill
\begin{minipage}[t]{0.73\textwidth}
\vspace{0pt}
\raggedright
\textbf{Name:} {\large{$\mathbf{K7_{253}}$}} (chiral, non-rotatable$^{*}$) \\ \textbf{PD:} {\scriptsize\texttt{[0],[0,1,2,3],[1,4,5,2],[3,6,7,8],[4,8,9,10],[11,12,6,5],[12,11,13,7],[9,13,10,14],[14]}} \\ \textbf{EM:} {\scriptsize\texttt{(B0, A0C0C3D0, B1E0F3B2, B3F2G3E1, C1D3H0H2, G1G0D1C2, F1F0H1D2, E2G2E3I0, H3)}} \\ \textbf{Kauffman bracket:} {\scriptsize $-A^{24} + 3A^{20} + A^{18} - 4A^{16} - 3A^{14} + 4A^{12} + 4A^{10} - A^{8} - 3A^{6} + A^{2}$} \\ \textbf{Arrow:} {\scriptsize $-A^{12} + A^{8}L_2 + 2A^{8} + A^{6}L_1 - 2A^{4}L_2 - 2A^{4} - 3A^{2}L_1 + 2L_2 + 2 + 4L_1/A^{2} - L_2/A^{4} - 3L_1/A^{6} + L_1/A^{10}$} \\ \textbf{Mock:} {\scriptsize $w^{4} + w^{3} - 3w^{2} - 3w + 5 + 4/w - 3/w^{2} - 3/w^{3} + w^{-4} + w^{-5}$} \\ \textbf{Affine:} {\scriptsize $t^{2} - 2 + t^{-2}$} \\ \textbf{Yamada:} {\scriptsize $-A^{27} + 3A^{25} - 2A^{24} - 3A^{23} + 7A^{22} + A^{21} - 9A^{20} + 5A^{19} + 5A^{18} - 10A^{17} + 2A^{16} + 5A^{15} - 4A^{14} - A^{13} + 2A^{12} + 6A^{11} - 6A^{10} + 11A^{8} - 6A^{7} - 3A^{6} + 9A^{5} - 2A^{4} - 4A^{3} + 2A^{2} - 1$}
\end{minipage}

\noindent{\color{gray!40}\rule{\textwidth}{0.4pt}}
\vspace{0.9\baselineskip}
\noindent \begin{minipage}[t]{0.25\textwidth}
\vspace{0pt}
\centering
\includegraphics[page=373,width=\linewidth]{knotoids.pdf}
\end{minipage}
\hfill
\begin{minipage}[t]{0.73\textwidth}
\vspace{0pt}
\raggedright
\textbf{Name:} {\large{$\mathbf{K7_{254}}$}} (chiral, non-rotatable$^{*}$) \\ \textbf{PD:} {\scriptsize\texttt{[0],[0,1,2,3],[1,4,5,2],[3,6,7,8],[4,8,9,10],[11,12,6,5],[12,13,14,7],[14,11,10,9],[13]}} \\ \textbf{EM:} {\scriptsize\texttt{(B0, A0C0C3D0, B1E0F3B2, B3F2G3E1, C1D3H3H2, H1G0D1C2, F1I0H0D2, G2F0E3E2, G1)}} \\ \textbf{Kauffman bracket:} {\scriptsize $-A^{26} - A^{24} + A^{22} + 4A^{20} + A^{18} - 6A^{16} - 4A^{14} + 3A^{12} + 5A^{10} - 2A^{6} + A^{2}$} \\ \textbf{Arrow:} {\scriptsize $-A^{32} - A^{30}L_1 + A^{28} + 4A^{26}L_1 + A^{24} - 6A^{22}L_1 - 4A^{20} + 3A^{18}L_1 + 5A^{16} - 2A^{12} + A^{8}$} \\ \textbf{Mock:} {\scriptsize $w^{4} - w^{3} - w^{2} + 4w + 2 - 6/w - 4/w^{2} + 3/w^{3} + 3/w^{4}$} \\ \textbf{Affine:} {\scriptsize $-t + 2 - 1/t$} \\ \textbf{Yamada:} {\scriptsize $A^{26} - A^{25} - 4A^{24} + 4A^{23} + 6A^{22} - 10A^{21} - 2A^{20} + 14A^{19} - 7A^{18} - 11A^{17} + 13A^{16} - 8A^{14} + 6A^{13} + 6A^{12} - A^{11} - 7A^{10} + 10A^{9} + 2A^{8} - 17A^{7} + 7A^{6} + 7A^{5} - 14A^{4} - A^{3} + 6A^{2} - 2A - 3$}
\end{minipage}

\noindent{\color{gray!40}\rule{\textwidth}{0.4pt}}
\vspace{0.9\baselineskip}
\noindent \begin{minipage}[t]{0.25\textwidth}
\vspace{0pt}
\centering
\includegraphics[page=374,width=\linewidth]{knotoids.pdf}
\end{minipage}
\hfill
\begin{minipage}[t]{0.73\textwidth}
\vspace{0pt}
\raggedright
\textbf{Name:} {\large{$\mathbf{K7_{255}}$}} (chiral, non-rotatable$^{*}$) \\ \textbf{PD:} {\scriptsize\texttt{[0],[0,1,2,3],[1,4,5,2],[3,5,6,7],[7,8,9,4],[6,9,10,11],[11,12,13,8],[10,13,12,14],[14]}} \\ \textbf{EM:} {\scriptsize\texttt{(B0, A0C0C3D0, B1E3D1B2, B3C2F0E0, D3G3F1C1, D2E2H0G0, F3H2H1E1, F2G2G1I0, H3)}} \\ \textbf{Kauffman bracket:} {\scriptsize $-A^{21} - A^{19} - A^{17} + A^{15} + 2A^{13} + A^{11} - 2A^{9} - 2A^{7} + A^{3} + A$} \\ \textbf{Arrow:} {\scriptsize $A^{-12} + L_1/A^{14} + L_2/A^{16} - L_1/A^{18} - 2L_2/A^{20} - L_1/A^{22} + 2L_2/A^{24} + 2L_1/A^{26} - L_2/A^{28} + A^{-28} - L_1/A^{30} - 1/A^{32}$} \\ \textbf{Mock:} {\scriptsize $w^{6} + w^{5} - w^{3} - w^{2} - w + 1 + 2/w - 1/w^{2} - 1/w^{3} + w^{-4}$} \\ \textbf{Affine:} {\scriptsize $t^{2} + t - 4 + 1/t + t^{-2}$} \\ \textbf{Yamada:} {\scriptsize $A^{29} + A^{28} + 2A^{25} + A^{24} + 2A^{21} + A^{20} - 2A^{19} + 2A^{18} + 2A^{17} - 2A^{16} + A^{15} - 2A^{12} + A^{11} + A^{10} - 3A^{9} + A^{7} - A^{6} - 2A^{5} + A^{4} + A^{3} - A^{2} + 1$}
\end{minipage}

\noindent{\color{gray!40}\rule{\textwidth}{0.4pt}}
\vspace{0.9\baselineskip}
\noindent \begin{minipage}[t]{0.25\textwidth}
\vspace{0pt}
\centering
\includegraphics[page=375,width=\linewidth]{knotoids.pdf}
\end{minipage}
\hfill
\begin{minipage}[t]{0.73\textwidth}
\vspace{0pt}
\raggedright
\textbf{Name:} {\large{$\mathbf{K7_{256}}$}} (chiral, non-rotatable$^{*}$) \\ \textbf{PD:} {\scriptsize\texttt{[0],[0,1,2,3],[1,4,5,2],[3,6,7,8],[8,7,9,4],[5,10,11,6],[9,11,12,13],[10,13,14,12],[14]}} \\ \textbf{EM:} {\scriptsize\texttt{(B0, A0C0C3D0, B1E3F0B2, B3F3E1E0, D3D2G0C1, C2H0G1D1, E2F2H3H1, F1G3I0G2, H2)}} \\ \textbf{Kauffman bracket:} {\scriptsize $-2A^{17} - 2A^{15} + 3A^{13} + 4A^{11} - 4A^{7} - 2A^{5} + A^{3} + A$} \\ \textbf{Arrow:} {\scriptsize $2A^{2}L_1 + 2 - 3L_1/A^{2} - 4/A^{4} + 4/A^{8} + 2L_1/A^{10} - 1/A^{12} - L_1/A^{14}$} \\ \textbf{Mock:} {\scriptsize $w^{4} + w^{3} - 2w^{2} + 4 - 2/w - 3/w^{2} + w^{-3} + w^{-4}$} \\ \textbf{Affine:} {\scriptsize $t - 2 + 1/t$} \\ \textbf{Yamada:} {\scriptsize $-A^{26} + 3A^{24} - 5A^{22} + 3A^{21} + 3A^{20} - 8A^{19} + 4A^{17} - 4A^{16} - A^{15} + 2A^{14} + 3A^{13} - 4A^{12} - A^{11} + 5A^{10} - 6A^{9} - 3A^{8} + 6A^{7} - A^{6} - 4A^{5} + 2A^{4} + 3A^{3} - 2A^{2} - A + 1$}
\end{minipage}

\noindent{\color{gray!40}\rule{\textwidth}{0.4pt}}
\vspace{0.9\baselineskip}
\noindent \begin{minipage}[t]{0.25\textwidth}
\vspace{0pt}
\centering
\includegraphics[page=376,width=\linewidth]{knotoids.pdf}
\end{minipage}
\hfill
\begin{minipage}[t]{0.73\textwidth}
\vspace{0pt}
\raggedright
\textbf{Name:} {\large{$\mathbf{K7_{257}}$}} (chiral, non-rotatable$^{*}$) \\ \textbf{PD:} {\scriptsize\texttt{[0],[0,1,2,3],[1,4,5,2],[3,6,7,8],[8,9,10,4],[11,12,6,5],[13,14,9,7],[10,13,12,11],[14]}} \\ \textbf{EM:} {\scriptsize\texttt{(B0, A0C0C3D0, B1E3F3B2, B3F2G3E0, D3G2H0C1, H3H2D1C2, H1I0E1D2, E2G0F1F0, G1)}} \\ \textbf{Kauffman bracket:} {\scriptsize $A^{22} - 3A^{18} - A^{16} + 4A^{14} + 3A^{12} - 3A^{10} - 5A^{8} + A^{6} + 3A^{4} + A^{2}$} \\ \textbf{Arrow:} {\scriptsize $A^{28} - 3A^{24} - A^{22}L_1 + 4A^{20} + 3A^{18}L_1 - 3A^{16} - 5A^{14}L_1 + A^{12} + 3A^{10}L_1 + A^{8}$} \\ \textbf{Mock:} {\scriptsize $-w^{4} - w^{3} + 4w^{2} + 3w - 5 - 5/w + 2/w^{2} + 3/w^{3} + w^{-4}$} \\ \textbf{Affine:} {\scriptsize $-2t + 4 - 2/t$} \\ \textbf{Yamada:} {\scriptsize $-A^{26} + A^{25} + 2A^{24} - 5A^{23} - 2A^{22} + 10A^{21} - 2A^{20} - 8A^{19} + 11A^{18} - 8A^{16} + 2A^{15} - 3A^{13} - 6A^{12} + 6A^{11} + 4A^{10} - 8A^{9} + 7A^{8} + 9A^{7} - 8A^{6} + 7A^{4} - 2A^{3} - 3A^{2} + 2A + 1$}
\end{minipage}

\noindent{\color{gray!40}\rule{\textwidth}{0.4pt}}
\vspace{0.9\baselineskip}
\noindent \begin{minipage}[t]{0.25\textwidth}
\vspace{0pt}
\centering
\includegraphics[page=377,width=\linewidth]{knotoids.pdf}
\end{minipage}
\hfill
\begin{minipage}[t]{0.73\textwidth}
\vspace{0pt}
\raggedright
\textbf{Name:} {\large{$\mathbf{K7_{258}}$}} (chiral, non-rotatable$^{*}$) \\ \textbf{PD:} {\scriptsize\texttt{[0],[0,1,2,3],[1,4,5,2],[3,6,7,8],[8,9,10,4],[11,12,6,5],[12,13,14,7],[9,14,11,10],[13]}} \\ \textbf{EM:} {\scriptsize\texttt{(B0, A0C0C3D0, B1E3F3B2, B3F2G3E0, D3H0H3C1, H2G0D1C2, F1I0H1D2, E1G2F0E2, G1)}} \\ \textbf{Kauffman bracket:} {\scriptsize $-A^{26} + 2A^{22} + A^{20} - 3A^{18} - 4A^{16} + 2A^{14} + 6A^{12} + 2A^{10} - 3A^{8} - 2A^{6} + A^{2}$} \\ \textbf{Arrow:} {\scriptsize $-A^{20} + 2A^{16} + A^{14}L_1 - 3A^{12} - 4A^{10}L_1 + 2A^{8} + 6A^{6}L_1 + 2A^{4} - 3A^{2}L_1 - 2 + A^{-4}$} \\ \textbf{Mock:} {\scriptsize $w^{4} + w^{3} - 2w^{2} - 4w + 2 + 6/w + w^{-2} - 3/w^{3} - 1/w^{4}$} \\ \textbf{Affine:} {\scriptsize $t - 2 + 1/t$} \\ \textbf{Yamada:} {\scriptsize $A^{26} - 4A^{24} - A^{23} + 8A^{22} - 3A^{21} - 9A^{20} + 12A^{19} + 4A^{18} - 12A^{17} + 7A^{16} + 5A^{15} - 8A^{14} - A^{13} + 2A^{12} + A^{11} - 10A^{10} + 2A^{9} + 9A^{8} - 13A^{7} - 2A^{6} + 12A^{5} - 6A^{4} - 5A^{3} + 5A^{2} + A - 1$}
\end{minipage}

\noindent{\color{gray!40}\rule{\textwidth}{0.4pt}}
\vspace{0.9\baselineskip}
\noindent \begin{minipage}[t]{0.25\textwidth}
\vspace{0pt}
\centering
\includegraphics[page=378,width=\linewidth]{knotoids.pdf}
\end{minipage}
\hfill
\begin{minipage}[t]{0.73\textwidth}
\vspace{0pt}
\raggedright
\textbf{Name:} {\large{$\mathbf{K7_{259}}$}} (chiral, non-rotatable$^{*}$) \\ \textbf{PD:} {\scriptsize\texttt{[0],[0,1,2,3],[1,4,5,2],[6,7,4,3],[5,8,9,6],[7,10,11,12],[8,13,10,9],[11,13,12,14],[14]}} \\ \textbf{EM:} {\scriptsize\texttt{(B0, A0C0C3D3, B1D2E0B2, E3F0C1B3, C2G0G3D0, D1G2H0H2, E1H1F1E2, F2G1F3I0, H3)}} \\ \textbf{Kauffman bracket:} {\scriptsize $A^{18} + 2A^{16} + A^{14} - 3A^{12} - 3A^{10} + A^{8} + 4A^{6} + A^{4} - 2A^{2} - 1$} \\ \textbf{Arrow:} {\scriptsize $A^{-12} + 2L_1/A^{14} + A^{-16} - 3L_1/A^{18} - 3/A^{20} + L_1/A^{22} + 4/A^{24} + L_1/A^{26} - 2/A^{28} - L_1/A^{30}$} \\ \textbf{Mock:} {\scriptsize $w^{6} + 2w^{5} + w^{4} - 3w^{3} - 3w^{2} + w + 3 + 1/w - 2/w^{2} - 1/w^{3} + w^{-4}$} \\ \textbf{Affine:} {\scriptsize $2t - 4 + 2/t$} \\ \textbf{Yamada:} {\scriptsize $-A^{26} - A^{25} + A^{24} + A^{23} - 4A^{22} - A^{21} + 3A^{20} - 4A^{19} - 5A^{18} + 4A^{17} - A^{16} - 6A^{15} + 3A^{14} + 2A^{13} - 3A^{12} + 2A^{11} + 3A^{10} + 3A^{9} - 4A^{8} + 3A^{7} + 5A^{6} - 7A^{5} + A^{4} + 3A^{3} - 3A^{2} - A + 1$}
\end{minipage}

\noindent{\color{gray!40}\rule{\textwidth}{0.4pt}}
\vspace{0.9\baselineskip}
\noindent \begin{minipage}[t]{0.25\textwidth}
\vspace{0pt}
\centering
\includegraphics[page=379,width=\linewidth]{knotoids.pdf}
\end{minipage}
\hfill
\begin{minipage}[t]{0.73\textwidth}
\vspace{0pt}
\raggedright
\textbf{Name:} {\large{$\mathbf{K7_{260}}$}} (chiral, non-rotatable$^{*}$) \\ \textbf{PD:} {\scriptsize\texttt{[0],[0,1,2,3],[1,4,5,2],[6,7,8,3],[4,8,9,10],[5,11,12,6],[7,13,14,9],[13,12,11,10],[14]}} \\ \textbf{EM:} {\scriptsize\texttt{(B0, A0C0C3D3, B1E0F0B2, F3G0E1B3, C1D2G3H3, C2H2H1D0, D1H0I0E2, G1F2F1E3, G2)}} \\ \textbf{Kauffman bracket:} {\scriptsize $A^{24} + A^{22} - A^{20} - 2A^{18} - A^{16} + 2A^{14} + 3A^{12} - 2A^{8} - A^{6} + A^{2}$} \\ \textbf{Arrow:} {\scriptsize $A^{6}L_1 + A^{4} - A^{2}L_1 - 2 - L_1/A^{2} + 2/A^{4} + 3L_1/A^{6} - 2L_1/A^{10} - 1/A^{12} + A^{-16}$} \\ \textbf{Mock:} {\scriptsize $-w^{4} - w^{3} + 2w^{2} + 3w + 1 - 1/w - 1/w^{3} - 1/w^{4}$} \\ \textbf{Affine:} {\scriptsize $2t - 4 + 2/t$} \\ \textbf{Yamada:} {\scriptsize $A^{25} - 2A^{23} - A^{22} + A^{21} - 3A^{19} + 5A^{17} - 2A^{16} - 2A^{15} + 4A^{14} - 2A^{13} - 3A^{12} + A^{11} - 2A^{10} - 2A^{9} - 2A^{8} + 3A^{7} + A^{6} - 4A^{5} + 4A^{4} + 2A^{3} - 4A^{2} + 1$}
\end{minipage}

\noindent{\color{gray!40}\rule{\textwidth}{0.4pt}}
\vspace{0.9\baselineskip}
\noindent \begin{minipage}[t]{0.25\textwidth}
\vspace{0pt}
\centering
\includegraphics[page=380,width=\linewidth]{knotoids.pdf}
\end{minipage}
\hfill
\begin{minipage}[t]{0.73\textwidth}
\vspace{0pt}
\raggedright
\textbf{Name:} {\large{$\mathbf{K7_{261}}$}} (chiral, rotatable) \\ \textbf{PD:} {\scriptsize\texttt{[0],[0,1,2,3],[1,4,5,2],[6,7,8,3],[4,8,9,10],[11,12,6,5],[7,13,14,9],[10,13,12,11],[14]}} \\ \textbf{EM:} {\scriptsize\texttt{(B0, A0C0C3D3, B1E0F3B2, F2G0E1B3, C1D2G3H0, H3H2D0C2, D1H1I0E2, E3G1F1F0, G2)}} \\ \textbf{Kauffman bracket:} {\scriptsize $2A^{20} + A^{18} - 4A^{16} - 3A^{14} + 2A^{12} + 4A^{10} - 2A^{6} + A^{2}$} \\ \textbf{Arrow:} {\scriptsize $2A^{14}L_1 + A^{12} - 4A^{10}L_1 - 3A^{8} + 2A^{6}L_1 + 4A^{4} - 2 + A^{-4}$} \\ \textbf{Mock:} {\scriptsize $2w^{3} + 2w^{2} - 4w - 4 + 2/w + 4/w^{2} - 1/w^{4}$} \\ \textbf{Affine:} {\scriptsize $0$} \\ \textbf{Yamada:} {\scriptsize $-A^{25} - A^{24} + 3A^{23} + 2A^{22} - 3A^{21} - A^{20} + 5A^{19} - 3A^{18} - 7A^{17} + 6A^{16} - 4A^{14} + 3A^{13} + 2A^{12} - 3A^{10} + 4A^{9} - 8A^{7} + 2A^{6} + 3A^{5} - 6A^{4} - A^{3} + 4A^{2} - 2$}
\end{minipage}

\noindent{\color{gray!40}\rule{\textwidth}{0.4pt}}
\vspace{0.9\baselineskip}
\noindent \begin{minipage}[t]{0.25\textwidth}
\vspace{0pt}
\centering
\includegraphics[page=381,width=\linewidth]{knotoids.pdf}
\end{minipage}
\hfill
\begin{minipage}[t]{0.73\textwidth}
\vspace{0pt}
\raggedright
\textbf{Name:} {\large{$\mathbf{K7_{262}}$}} (chiral, non-rotatable$^{*}$) \\ \textbf{PD:} {\scriptsize\texttt{[0],[0,1,2,3],[1,4,5,2],[6,7,8,3],[8,9,10,4],[5,11,12,6],[7,12,13,14],[9,14,11,10],[13]}} \\ \textbf{EM:} {\scriptsize\texttt{(B0, A0C0C3D3, B1E3F0B2, F3G0E0B3, D2H0H3C1, C2H2G1D0, D1F2I0H1, E1G3F1E2, G2)}} \\ \textbf{Kauffman bracket:} {\scriptsize $A^{16} + 2A^{14} - 2A^{12} - 4A^{10} + 5A^{6} + 3A^{4} - 2A^{2} - 2$} \\ \textbf{Arrow:} {\scriptsize $L_1/A^{2} + 2/A^{4} - 2L_1/A^{6} - 4/A^{8} + 5/A^{12} + 3L_1/A^{14} - 2/A^{16} - 2L_1/A^{18}$} \\ \textbf{Mock:} {\scriptsize $-w^{4} - w^{3} + 4w^{2} + 2w - 4 + 3/w^{2} - 1/w^{3} - 1/w^{4}$} \\ \textbf{Affine:} {\scriptsize $t - 2 + 1/t$} \\ \textbf{Yamada:} {\scriptsize $-A^{26} + 2A^{25} + 3A^{24} - 5A^{23} - A^{22} + 8A^{21} - 4A^{20} - 5A^{19} + 8A^{18} + A^{17} - 3A^{16} + 5A^{15} + 3A^{14} - 5A^{12} + 4A^{11} + A^{10} - 10A^{9} + 5A^{8} + 4A^{7} - 7A^{6} + A^{5} + 4A^{4} - A^{3} - 2A^{2} + 1$}
\end{minipage}

\noindent{\color{gray!40}\rule{\textwidth}{0.4pt}}
\vspace{0.9\baselineskip}
\noindent \begin{minipage}[t]{0.25\textwidth}
\vspace{0pt}
\centering
\includegraphics[page=382,width=\linewidth]{knotoids.pdf}
\end{minipage}
\hfill
\begin{minipage}[t]{0.73\textwidth}
\vspace{0pt}
\raggedright
\textbf{Name:} {\large{$\mathbf{K7_{263}}$}} (chiral, non-rotatable$^{*}$) \\ \textbf{PD:} {\scriptsize\texttt{[0],[0,1,2,3],[1,4,5,2],[3,6,7,8],[4,8,9,10],[11,12,6,5],[12,13,9,7],[10,13,14,11],[14]}} \\ \textbf{EM:} {\scriptsize\texttt{(B0, A0C0C3D0, B1E0F3B2, B3F2G3E1, C1D3G2H0, H3G0D1C2, F1H1E2D2, E3G1I0F0, H2)}} \\ \textbf{Kauffman bracket:} {\scriptsize $-A^{24} + 2A^{20} + A^{18} - 4A^{16} - 3A^{14} + 3A^{12} + 4A^{10} - 2A^{6} + A^{2}$} \\ \textbf{Arrow:} {\scriptsize $-A^{30}L_1 + 2A^{26}L_1 + A^{24} - 4A^{22}L_1 - 3A^{20} + 3A^{18}L_1 + 4A^{16} - 2A^{12} + A^{8}$} \\ \textbf{Mock:} {\scriptsize $w^{4} - w^{3} - 2w^{2} + 2w + 2 - 4/w - 3/w^{2} + 3/w^{3} + 3/w^{4}$} \\ \textbf{Affine:} {\scriptsize $-3t + 6 - 3/t$} \\ \textbf{Yamada:} {\scriptsize $-A^{25} - A^{24} + 2A^{23} + 2A^{22} - 4A^{21} - A^{20} + 7A^{19} - 2A^{18} - 6A^{17} + 8A^{16} - 5A^{14} + 3A^{13} + A^{12} - A^{11} - 4A^{10} + 5A^{9} + 2A^{8} - 7A^{7} + 4A^{6} + 4A^{5} - 8A^{4} - 2A^{3} + 3A^{2} - 2A - 3$}
\end{minipage}

\noindent{\color{gray!40}\rule{\textwidth}{0.4pt}}
\vspace{0.9\baselineskip}
\noindent \begin{minipage}[t]{0.25\textwidth}
\vspace{0pt}
\centering
\includegraphics[page=383,width=\linewidth]{knotoids.pdf}
\end{minipage}
\hfill
\begin{minipage}[t]{0.73\textwidth}
\vspace{0pt}
\raggedright
\textbf{Name:} {\large{$\mathbf{K7_{264}}$}} (chiral, non-rotatable$^{*}$) \\ \textbf{PD:} {\scriptsize\texttt{[0],[0,1,2,3],[1,4,5,2],[3,6,7,8],[4,8,9,10],[11,12,6,5],[12,13,9,7],[10,14,13,11],[14]}} \\ \textbf{EM:} {\scriptsize\texttt{(B0, A0C0C3D0, B1E0F3B2, B3F2G3E1, C1D3G2H0, H3G0D1C2, F1H2E2D2, E3I0G1F0, H1)}} \\ \textbf{Kauffman bracket:} {\scriptsize $A^{22} + 2A^{20} - A^{18} - 5A^{16} - 2A^{14} + 5A^{12} + 4A^{10} - A^{8} - 3A^{6} + A^{2}$} \\ \textbf{Arrow:} {\scriptsize $A^{10}L_1 + A^{8}L_2 + A^{8} - A^{6}L_1 - 3A^{4}L_2 - 2A^{4} - 2A^{2}L_1 + 3L_2 + 2 + 4L_1/A^{2} - L_2/A^{4} - 3L_1/A^{6} + L_1/A^{10}$} \\ \textbf{Mock:} {\scriptsize $-2w^{2} - w + 6 + 3/w - 4/w^{2} - 3/w^{3} + w^{-4} + w^{-5}$} \\ \textbf{Affine:} {\scriptsize $0$} \\ \textbf{Yamada:} {\scriptsize $A^{26} + A^{25} - 4A^{24} + 7A^{22} - 5A^{21} - 7A^{20} + 10A^{19} - 10A^{17} + 8A^{16} + 4A^{15} - 5A^{14} + 2A^{13} + 4A^{12} + 2A^{11} - 8A^{10} + 4A^{9} + 7A^{8} - 12A^{7} + 2A^{6} + 10A^{5} - 6A^{4} - 2A^{3} + 4A^{2} - 1$}
\end{minipage}

\noindent{\color{gray!40}\rule{\textwidth}{0.4pt}}
\vspace{0.9\baselineskip}
\noindent \begin{minipage}[t]{0.25\textwidth}
\vspace{0pt}
\centering
\includegraphics[page=384,width=\linewidth]{knotoids.pdf}
\end{minipage}
\hfill
\begin{minipage}[t]{0.73\textwidth}
\vspace{0pt}
\raggedright
\textbf{Name:} {\large{$\mathbf{K7_{265}}$}} (chiral, non-rotatable$^{*}$) \\ \textbf{PD:} {\scriptsize\texttt{[0],[0,1,2,3],[1,4,5,2],[3,6,7,8],[8,9,10,4],[11,12,6,5],[12,13,9,7],[10,14,13,11],[14]}} \\ \textbf{EM:} {\scriptsize\texttt{(B0, A0C0C3D0, B1E3F3B2, B3F2G3E0, D3G2H0C1, H3G0D1C2, F1H2E1D2, E2I0G1F0, H1)}} \\ \textbf{Kauffman bracket:} {\scriptsize $A^{22} - 3A^{18} + 5A^{14} + 2A^{12} - 6A^{10} - 5A^{8} + 3A^{6} + 5A^{4} - 1$} \\ \textbf{Arrow:} {\scriptsize $A^{22}L_1 - 3A^{18}L_1 + 5A^{14}L_1 + A^{12}L_2 + A^{12} - 6A^{10}L_1 - 3A^{8}L_2 - 2A^{8} + 3A^{6}L_1 + 3A^{4}L_2 + 2A^{4} - L_2$} \\ \textbf{Mock:} {\scriptsize $-w^{5} - w^{4} + 4w^{3} + 3w^{2} - 7w - 4 + 5/w + 4/w^{2} - 1/w^{3} - 1/w^{4}$} \\ \textbf{Affine:} {\scriptsize $-t + 2 - 1/t$} \\ \textbf{Yamada:} {\scriptsize $-3A^{26} + A^{25} + 8A^{24} - 5A^{23} - 8A^{22} + 14A^{21} + A^{20} - 17A^{19} + 9A^{18} + 6A^{17} - 11A^{16} + 6A^{14} + A^{13} - 10A^{12} + 6A^{11} + 10A^{10} - 16A^{9} + 14A^{7} - 10A^{6} - 7A^{5} + 8A^{4} + A^{3} - 6A^{2} + 2$}
\end{minipage}

\noindent{\color{gray!40}\rule{\textwidth}{0.4pt}}
\vspace{0.9\baselineskip}
\noindent \begin{minipage}[t]{0.25\textwidth}
\vspace{0pt}
\centering
\includegraphics[page=385,width=\linewidth]{knotoids.pdf}
\end{minipage}
\hfill
\begin{minipage}[t]{0.73\textwidth}
\vspace{0pt}
\raggedright
\textbf{Name:} {\large{$\mathbf{K7_{266}}$}} (chiral, non-rotatable$^{*}$) \\ \textbf{PD:} {\scriptsize\texttt{[0],[0,1,2,3],[1,4,5,2],[3,6,7,8],[8,9,10,4],[11,12,6,5],[12,13,9,7],[13,14,11,10],[14]}} \\ \textbf{EM:} {\scriptsize\texttt{(B0, A0C0C3D0, B1E3F3B2, B3F2G3E0, D3G2H3C1, H2G0D1C2, F1H0E1D2, G1I0F0E2, H1)}} \\ \textbf{Kauffman bracket:} {\scriptsize $-A^{26} - A^{24} + 2A^{22} + 3A^{20} - A^{18} - 6A^{16} - A^{14} + 6A^{12} + 4A^{10} - 2A^{8} - 3A^{6} + A^{2}$} \\ \textbf{Arrow:} {\scriptsize $-A^{20} - A^{18}L_1 + 2A^{16} + 3A^{14}L_1 - A^{12} - 6A^{10}L_1 - A^{8} + 6A^{6}L_1 + 4A^{4} - 2A^{2}L_1 - 3 + A^{-4}$} \\ \textbf{Mock:} {\scriptsize $w^{4} + 2w^{3} - w^{2} - 7w - 2 + 7/w + 4/w^{2} - 2/w^{3} - 1/w^{4}$} \\ \textbf{Affine:} {\scriptsize $-t + 2 - 1/t$} \\ \textbf{Yamada:} {\scriptsize $2A^{26} - 5A^{24} + 2A^{23} + 8A^{22} - 9A^{21} - 9A^{20} + 15A^{19} - 2A^{18} - 16A^{17} + 13A^{16} + 6A^{15} - 9A^{14} + 4A^{13} + 7A^{12} - 11A^{10} + 8A^{9} + 8A^{8} - 17A^{7} + 4A^{6} + 15A^{5} - 10A^{4} - 4A^{3} + 9A^{2} - 3$}
\end{minipage}

\noindent{\color{gray!40}\rule{\textwidth}{0.4pt}}
\vspace{0.9\baselineskip}
\noindent \begin{minipage}[t]{0.25\textwidth}
\vspace{0pt}
\centering
\includegraphics[page=386,width=\linewidth]{knotoids.pdf}
\end{minipage}
\hfill
\begin{minipage}[t]{0.73\textwidth}
\vspace{0pt}
\raggedright
\textbf{Name:} {\large{$\mathbf{K7_{267}}$}} (chiral, non-rotatable$^{*}$) \\ \textbf{PD:} {\scriptsize\texttt{[0],[0,1,2,3],[1,4,5,2],[6,7,8,3],[4,8,9,10],[5,11,12,6],[7,12,13,9],[10,13,14,11],[14]}} \\ \textbf{EM:} {\scriptsize\texttt{(B0, A0C0C3D3, B1E0F0B2, F3G0E1B3, C1D2G3H0, C2H3G1D0, D1F2H1E2, E3G2I0F1, H2)}} \\ \textbf{Kauffman bracket:} {\scriptsize $-A^{20} - A^{18} + 2A^{16} + 4A^{14} - A^{12} - 5A^{10} - 2A^{8} + 4A^{6} + 3A^{4} - A^{2} - 1$} \\ \textbf{Arrow:} {\scriptsize $-A^{14}L_1 - A^{12} + 2A^{10}L_1 + 4A^{8} - A^{6}L_1 - 5A^{4} - 2A^{2}L_1 + 4 + 3L_1/A^{2} - 1/A^{4} - L_1/A^{6}$} \\ \textbf{Mock:} {\scriptsize $w^{4} - 4w^{2} - w + 6 + 1/w - 3/w^{2} + w^{-4}$} \\ \textbf{Affine:} {\scriptsize $-t + 2 - 1/t$} \\ \textbf{Yamada:} {\scriptsize $-2A^{27} + 5A^{25} - 6A^{23} + 4A^{22} + 8A^{21} - 9A^{20} - 4A^{19} + 9A^{18} - 6A^{17} - 6A^{16} + 5A^{15} - 4A^{13} - 2A^{12} + 8A^{11} - 3A^{10} - 8A^{9} + 10A^{8} + A^{7} - 10A^{6} + 3A^{5} + 4A^{4} - 3A^{3} - 2A^{2} + A + 1$}
\end{minipage}

\noindent{\color{gray!40}\rule{\textwidth}{0.4pt}}
\vspace{0.9\baselineskip}
\noindent \begin{minipage}[t]{0.25\textwidth}
\vspace{0pt}
\centering
\includegraphics[page=387,width=\linewidth]{knotoids.pdf}
\end{minipage}
\hfill
\begin{minipage}[t]{0.73\textwidth}
\vspace{0pt}
\raggedright
\textbf{Name:} {\large{$\mathbf{K7_{268}}$}} (chiral, non-rotatable$^{*}$) \\ \textbf{PD:} {\scriptsize\texttt{[0],[0,1,2,3],[1,4,5,2],[6,7,8,3],[4,8,9,10],[5,11,12,6],[7,12,13,9],[14,13,11,10],[14]}} \\ \textbf{EM:} {\scriptsize\texttt{(B0, A0C0C3D3, B1E0F0B2, F3G0E1B3, C1D2G3H3, C2H2G1D0, D1F2H1E2, I0G2F1E3, H0)}} \\ \textbf{Kauffman bracket:} {\scriptsize $A^{24} + A^{22} - 2A^{18} - 3A^{16} + 4A^{12} + 2A^{10} - A^{8} - 2A^{6} + A^{2}$} \\ \textbf{Arrow:} {\scriptsize $A^{12} + A^{10}L_1 + A^{8}L_2 - A^{8} - 2A^{6}L_1 - 3A^{4}L_2 + 3L_2 + 1 + 2L_1/A^{2} - L_2/A^{4} - 2L_1/A^{6} + L_1/A^{10}$} \\ \textbf{Mock:} {\scriptsize $-w^{4} - w^{3} + w + 4 + 1/w - 3/w^{2} - 2/w^{3} + w^{-4} + w^{-5}$} \\ \textbf{Affine:} {\scriptsize $-t + 2 - 1/t$} \\ \textbf{Yamada:} {\scriptsize $A^{27} - 2A^{25} + 3A^{23} - 4A^{21} + 4A^{19} - 4A^{18} - 3A^{17} + 7A^{16} - A^{14} + 5A^{13} + A^{12} - A^{11} - 2A^{10} + 3A^{9} - 2A^{8} - 6A^{7} + 5A^{6} + 3A^{5} - 4A^{4} + A^{3} + 3A^{2} - 1$}
\end{minipage}

\noindent{\color{gray!40}\rule{\textwidth}{0.4pt}}
\vspace{0.9\baselineskip}
\noindent \begin{minipage}[t]{0.25\textwidth}
\vspace{0pt}
\centering
\includegraphics[page=388,width=\linewidth]{knotoids.pdf}
\end{minipage}
\hfill
\begin{minipage}[t]{0.73\textwidth}
\vspace{0pt}
\raggedright
\textbf{Name:} {\large{$\mathbf{K7_{269}}$}} (chiral, non-rotatable$^{*}$) \\ \textbf{PD:} {\scriptsize\texttt{[0],[0,1,2,3],[1,4,5,6],[6,7,8,2],[9,10,4,3],[10,11,12,5],[7,12,11,13],[13,14,9,8],[14]}} \\ \textbf{EM:} {\scriptsize\texttt{(B0, A0C0D3E3, B1E2F3D0, C3G0H3B2, H2F0C1B3, E1G2G1C2, D1F2F1H0, G3I0E0D2, H1)}} \\ \textbf{Kauffman bracket:} {\scriptsize $2A^{22} + A^{20} - 4A^{18} - 4A^{16} + 4A^{14} + 6A^{12} - A^{10} - 5A^{8} + 3A^{4} - 1$} \\ \textbf{Arrow:} {\scriptsize $A^{4}L_2 + A^{4} + A^{2}L_1 - 2L_2 - 2 - 4L_1/A^{2} + 2L_2/A^{4} + 2/A^{4} + 6L_1/A^{6} - L_2/A^{8} - 5L_1/A^{10} + 3L_1/A^{14} - L_1/A^{18}$} \\ \textbf{Mock:} {\scriptsize $-w^{4} - 2w^{3} + 3w^{2} + 6w - 3 - 7/w + 3/w^{2} + 4/w^{3} - 1/w^{4} - 1/w^{5}$} \\ \textbf{Affine:} {\scriptsize $t^{2} - 2 + t^{-2}$} \\ \textbf{Yamada:} {\scriptsize $A^{27} - A^{26} - 3A^{25} + 3A^{24} + 5A^{23} - 6A^{22} - 3A^{21} + 12A^{20} - 2A^{19} - 13A^{18} + 10A^{17} + 5A^{16} - 13A^{15} + 5A^{14} + 8A^{13} - 6A^{12} - 2A^{11} + 6A^{10} + 6A^{9} - 11A^{8} + 2A^{7} + 14A^{6} - 12A^{5} - 4A^{4} + 11A^{3} - 4A^{2} - 4A + 2$}
\end{minipage}

\noindent{\color{gray!40}\rule{\textwidth}{0.4pt}}
\vspace{0.9\baselineskip}
\noindent \begin{minipage}[t]{0.25\textwidth}
\vspace{0pt}
\centering
\includegraphics[page=389,width=\linewidth]{knotoids.pdf}
\end{minipage}
\hfill
\begin{minipage}[t]{0.73\textwidth}
\vspace{0pt}
\raggedright
\textbf{Name:} {\large{$\mathbf{K7_{270}}$}} (chiral, non-rotatable$^{*}$) \\ \textbf{PD:} {\scriptsize\texttt{[0],[0,1,2,3],[1,4,5,6],[2,6,7,8],[9,10,4,3],[10,11,12,5],[13,12,14,7],[8,14,11,9],[13]}} \\ \textbf{EM:} {\scriptsize\texttt{(B0, A0C0D0E3, B1E2F3D1, B2C3G3H0, H3F0C1B3, E1H2G1C2, I0F2H1D2, D3G2F1E0, G0)}} \\ \textbf{Kauffman bracket:} {\scriptsize $A^{23} + A^{21} - 3A^{19} - 4A^{17} + A^{15} + 5A^{13} + A^{11} - 3A^{9} - A^{7} + 2A^{5} - A$} \\ \textbf{Arrow:} {\scriptsize $-A^{8} - A^{6}L_1 + 3A^{4} + 4A^{2}L_1 + L_2 - 2 - 5L_1/A^{2} - 2L_2/A^{4} + A^{-4} + 3L_1/A^{6} + L_2/A^{8} - 2L_1/A^{10} + L_1/A^{14}$} \\ \textbf{Mock:} {\scriptsize $w^{6} + w^{5} - 3w^{4} - 3w^{3} + 3w^{2} + 5w + 1 - 5/w - 2/w^{2} + 2/w^{3} + w^{-4}$} \\ \textbf{Affine:} {\scriptsize $0$} \\ \textbf{Yamada:} {\scriptsize $-A^{25} + 2A^{24} + 2A^{23} - 5A^{22} + A^{21} + 5A^{20} - 5A^{19} - 3A^{18} + 7A^{17} - A^{16} - 4A^{15} + 6A^{14} + 2A^{13} - 3A^{12} + A^{11} + 3A^{10} + A^{9} - 6A^{8} + 5A^{7} + 4A^{6} - 9A^{5} + 3A^{4} + 3A^{3} - 4A^{2} + A + 1$}
\end{minipage}

\noindent{\color{gray!40}\rule{\textwidth}{0.4pt}}
\vspace{0.9\baselineskip}
\noindent \begin{minipage}[t]{0.25\textwidth}
\vspace{0pt}
\centering
\includegraphics[page=390,width=\linewidth]{knotoids.pdf}
\end{minipage}
\hfill
\begin{minipage}[t]{0.73\textwidth}
\vspace{0pt}
\raggedright
\textbf{Name:} {\large{$\mathbf{K7_{271}}$}} (chiral, non-rotatable$^{*}$) \\ \textbf{PD:} {\scriptsize\texttt{[0],[0,1,2,3],[1,4,5,6],[6,7,8,2],[3,8,9,10],[4,10,11,12],[12,13,7,5],[13,14,11,9],[14]}} \\ \textbf{EM:} {\scriptsize\texttt{(B0, A0C0D3E0, B1F0G3D0, C3G2E1B2, B3D2H3F1, C1E3H2G0, F3H0D1C2, G1I0F2E2, H1)}} \\ \textbf{Kauffman bracket:} {\scriptsize $-A^{24} + 4A^{20} + 2A^{18} - 5A^{16} - 5A^{14} + 3A^{12} + 5A^{10} - 3A^{6} + A^{2}$} \\ \textbf{Arrow:} {\scriptsize $-A^{12} + A^{8}L_2 + 3A^{8} + 2A^{6}L_1 - 2A^{4}L_2 - 3A^{4} - 5A^{2}L_1 + L_2 + 2 + 5L_1/A^{2} - 3L_1/A^{6} + L_1/A^{10}$} \\ \textbf{Mock:} {\scriptsize $w^{4} + 2w^{3} - 3w^{2} - 5w + 5 + 5/w - 3/w^{2} - 3/w^{3} + w^{-4} + w^{-5}$} \\ \textbf{Affine:} {\scriptsize $0$} \\ \textbf{Yamada:} {\scriptsize $A^{26} - A^{25} - 4A^{24} + 3A^{23} + 6A^{22} - 7A^{21} - 5A^{20} + 14A^{19} - 2A^{18} - 13A^{17} + 11A^{16} + 2A^{15} - 10A^{14} + 2A^{13} + 4A^{12} - 2A^{11} - 9A^{10} + 7A^{9} + 5A^{8} - 15A^{7} + 4A^{6} + 11A^{5} - 10A^{4} - 3A^{3} + 8A^{2} - 3$}
\end{minipage}

\noindent{\color{gray!40}\rule{\textwidth}{0.4pt}}
\vspace{0.9\baselineskip}
\noindent \begin{minipage}[t]{0.25\textwidth}
\vspace{0pt}
\centering
\includegraphics[page=391,width=\linewidth]{knotoids.pdf}
\end{minipage}
\hfill
\begin{minipage}[t]{0.73\textwidth}
\vspace{0pt}
\raggedright
\textbf{Name:} {\large{$\mathbf{K7_{272}}$}} (chiral, non-rotatable$^{*}$) \\ \textbf{PD:} {\scriptsize\texttt{[0],[0,1,2,3],[1,4,5,2],[3,6,7,4],[5,8,9,6],[7,10,11,12],[8,12,13,9],[10,13,14,11],[14]}} \\ \textbf{EM:} {\scriptsize\texttt{(B0, A0C0C3D0, B1D3E0B2, B3E3F0C1, C2G0G3D1, D2H0H3G1, E1F3H1E2, F1G2I0F2, H2)}} \\ \textbf{Kauffman bracket:} {\scriptsize $A^{22} + A^{20} - 2A^{18} - 4A^{16} + 2A^{14} + 6A^{12} + 2A^{10} - 4A^{8} - 3A^{6} + A^{4} + A^{2}$} \\ \textbf{Arrow:} {\scriptsize $A^{4} + A^{2}L_1 - 2 - 4L_1/A^{2} - L_2/A^{4} + 3/A^{4} + 6L_1/A^{6} + 3L_2/A^{8} - 1/A^{8} - 4L_1/A^{10} - 3L_2/A^{12} + L_1/A^{14} + L_2/A^{16}$} \\ \textbf{Mock:} {\scriptsize $w^{5} - 4w^{3} + 6w + 1 - 4/w + w^{-3}$} \\ \textbf{Affine:} {\scriptsize $0$} \\ \textbf{Yamada:} {\scriptsize $-A^{26} + A^{25} + 2A^{24} - 5A^{23} + 2A^{22} + 9A^{21} - 10A^{20} - A^{19} + 13A^{18} - 8A^{17} - 3A^{16} + 8A^{15} - 2A^{14} - 2A^{13} - A^{12} + 8A^{11} - A^{10} - 8A^{9} + 11A^{8} - 12A^{6} + 7A^{5} + 3A^{4} - 8A^{3} + 2A^{2} + 3A - 1$}
\end{minipage}

\noindent{\color{gray!40}\rule{\textwidth}{0.4pt}}
\vspace{0.9\baselineskip}
\noindent \begin{minipage}[t]{0.25\textwidth}
\vspace{0pt}
\centering
\includegraphics[page=392,width=\linewidth]{knotoids.pdf}
\end{minipage}
\hfill
\begin{minipage}[t]{0.73\textwidth}
\vspace{0pt}
\raggedright
\textbf{Name:} {\large{$\mathbf{K7_{273}}$}} (chiral, non-rotatable$^{*}$) \\ \textbf{PD:} {\scriptsize\texttt{[0],[0,1,2,3],[1,4,5,2],[3,6,7,4],[5,8,9,6],[10,11,12,7],[8,12,13,9],[13,14,11,10],[14]}} \\ \textbf{EM:} {\scriptsize\texttt{(B0, A0C0C3D0, B1D3E0B2, B3E3F3C1, C2G0G3D1, H3H2G1D2, E1F2H0E2, G2I0F1F0, H1)}} \\ \textbf{Kauffman bracket:} {\scriptsize $A^{21} - A^{19} - 4A^{17} - A^{15} + 5A^{13} + 4A^{11} - 3A^{9} - 5A^{7} + 2A^{3} + A$} \\ \textbf{Arrow:} {\scriptsize $-L_1/A^{6} + A^{-8} + 4L_1/A^{10} + L_2/A^{12} - 5L_1/A^{14} - 3L_2/A^{16} - 1/A^{16} + 3L_1/A^{18} + 3L_2/A^{20} + 2/A^{20} - L_2/A^{24} - 1/A^{24} - L_1/A^{26}$} \\ \textbf{Mock:} {\scriptsize $-w^{5} + 4w^{3} + 4w^{2} - 4w - 5 + 2/w + 2/w^{2} - 1/w^{3}$} \\ \textbf{Affine:} {\scriptsize $2t - 4 + 2/t$} \\ \textbf{Yamada:} {\scriptsize $A^{27} - 3A^{25} + 4A^{24} + 5A^{23} - 11A^{22} + 4A^{21} + 13A^{20} - 10A^{19} + A^{18} + 11A^{17} - 4A^{16} - 3A^{15} + A^{14} + 5A^{13} - 7A^{12} - 6A^{11} + 11A^{10} - 5A^{9} - 10A^{8} + 11A^{7} + 2A^{6} - 9A^{5} + 4A^{4} + 4A^{3} - 3A^{2} - A + 1$}
\end{minipage}

\noindent{\color{gray!40}\rule{\textwidth}{0.4pt}}
\vspace{0.9\baselineskip}
\noindent \begin{minipage}[t]{0.25\textwidth}
\vspace{0pt}
\centering
\includegraphics[page=393,width=\linewidth]{knotoids.pdf}
\end{minipage}
\hfill
\begin{minipage}[t]{0.73\textwidth}
\vspace{0pt}
\raggedright
\textbf{Name:} {\large{$\mathbf{K7_{274}}$}} (chiral, non-rotatable$^{*}$) \\ \textbf{PD:} {\scriptsize\texttt{[0],[0,1,2,3],[1,4,5,2],[3,6,7,4],[5,8,9,10],[6,10,11,7],[8,11,12,13],[9,14,13,12],[14]}} \\ \textbf{EM:} {\scriptsize\texttt{(B0, A0C0C3D0, B1D3E0B2, B3F0F3C1, C2G0H0F1, D1E3G1D2, E1F2H3H2, E2I0G3G2, H1)}} \\ \textbf{Kauffman bracket:} {\scriptsize $-A^{18} + A^{16} + 3A^{14} - 5A^{10} - 3A^{8} + 3A^{6} + 4A^{4} - 1$} \\ \textbf{Arrow:} {\scriptsize $-A^{18}L_1 + A^{16} + 3A^{14}L_1 + A^{12}L_2 - A^{12} - 5A^{10}L_1 - 3A^{8}L_2 + 3A^{6}L_1 + 3A^{4}L_2 + A^{4} - L_2$} \\ \textbf{Mock:} {\scriptsize $-w^{5} - w^{4} + 3w^{3} + 2w^{2} - 5w - 2 + 3/w + 2/w^{2}$} \\ \textbf{Affine:} {\scriptsize $-2t + 4 - 2/t$} \\ \textbf{Yamada:} {\scriptsize $-A^{26} + A^{25} + 3A^{24} - 4A^{23} - A^{22} + 8A^{21} - 5A^{20} - 6A^{19} + 8A^{18} - 2A^{17} - 5A^{16} + 3A^{15} + A^{14} - A^{13} - 4A^{12} + 5A^{11} + 2A^{10} - 9A^{9} + 5A^{8} + 4A^{7} - 9A^{6} + A^{5} + 3A^{4} - 3A^{3} - 2A^{2} + A + 1$}
\end{minipage}

\noindent{\color{gray!40}\rule{\textwidth}{0.4pt}}
\vspace{0.9\baselineskip}
\noindent \begin{minipage}[t]{0.25\textwidth}
\vspace{0pt}
\centering
\includegraphics[page=394,width=\linewidth]{knotoids.pdf}
\end{minipage}
\hfill
\begin{minipage}[t]{0.73\textwidth}
\vspace{0pt}
\raggedright
\textbf{Name:} {\large{$\mathbf{K7_{275}}$}} (chiral, non-rotatable$^{*}$) \\ \textbf{PD:} {\scriptsize\texttt{[0],[0,1,2,3],[1,4,5,2],[3,6,7,4],[8,9,6,5],[7,10,11,12],[12,13,9,8],[10,13,14,11],[14]}} \\ \textbf{EM:} {\scriptsize\texttt{(B0, A0C0C3D0, B1D3E3B2, B3E2F0C1, G3G2D1C2, D2H0H3G0, F3H1E1E0, F1G1I0F2, H2)}} \\ \textbf{Kauffman bracket:} {\scriptsize $-A^{20} - A^{18} + 2A^{16} + 4A^{14} - A^{12} - 6A^{10} - 2A^{8} + 4A^{6} + 4A^{4} - A^{2} - 1$} \\ \textbf{Arrow:} {\scriptsize $-A^{20} - A^{18}L_1 + 2A^{16} + 4A^{14}L_1 + A^{12}L_2 - 2A^{12} - 6A^{10}L_1 - 3A^{8}L_2 + A^{8} + 4A^{6}L_1 + 3A^{4}L_2 + A^{4} - A^{2}L_1 - L_2$} \\ \textbf{Mock:} {\scriptsize $-w^{5} + 4w^{3} + w^{2} - 6w - 1 + 4/w + w^{-2} - 1/w^{3}$} \\ \textbf{Affine:} {\scriptsize $0$} \\ \textbf{Yamada:} {\scriptsize $A^{28} - A^{27} - 3A^{26} + 4A^{25} + 2A^{24} - 9A^{23} + 4A^{22} + 10A^{21} - 10A^{20} - 2A^{19} + 12A^{18} - 5A^{17} - 5A^{16} + 5A^{15} - 5A^{13} - 5A^{12} + 9A^{11} - 3A^{10} - 11A^{9} + 12A^{8} + A^{7} - 11A^{6} + 5A^{5} + 2A^{4} - 4A^{3} - A^{2} + A + 1$}
\end{minipage}

\noindent{\color{gray!40}\rule{\textwidth}{0.4pt}}
\vspace{0.9\baselineskip}
\noindent \begin{minipage}[t]{0.25\textwidth}
\vspace{0pt}
\centering
\includegraphics[page=395,width=\linewidth]{knotoids.pdf}
\end{minipage}
\hfill
\begin{minipage}[t]{0.73\textwidth}
\vspace{0pt}
\raggedright
\textbf{Name:} {\large{$\mathbf{K7_{276}}$}} (chiral, non-rotatable$^{*}$) \\ \textbf{PD:} {\scriptsize\texttt{[0],[0,1,2,3],[1,4,5,2],[3,6,7,4],[8,9,6,5],[10,11,12,7],[12,13,9,8],[13,14,11,10],[14]}} \\ \textbf{EM:} {\scriptsize\texttt{(B0, A0C0C3D0, B1D3E3B2, B3E2F3C1, G3G2D1C2, H3H2G0D2, F2H0E1E0, G1I0F1F0, H1)}} \\ \textbf{Kauffman bracket:} {\scriptsize $-A^{25} + 3A^{21} + 2A^{19} - 4A^{17} - 6A^{15} + A^{13} + 6A^{11} + 2A^{9} - 3A^{7} - 2A^{5} + A$} \\ \textbf{Arrow:} {\scriptsize $A^{10}L_1 - 3A^{6}L_1 - A^{4}L_2 - A^{4} + 4A^{2}L_1 + 3L_2 + 3 - L_1/A^{2} - 3L_2/A^{4} - 3/A^{4} - 2L_1/A^{6} + L_2/A^{8} + 2/A^{8} + 2L_1/A^{10} - L_1/A^{14}$} \\ \textbf{Mock:} {\scriptsize $w^{5} + 2w^{4} - 2w^{3} - 5w^{2} + 2w + 7 - 3/w^{2} - 1/w^{3}$} \\ \textbf{Affine:} {\scriptsize $2t - 4 + 2/t$} \\ \textbf{Yamada:} {\scriptsize $A^{27} - 2A^{26} - 3A^{25} + 8A^{24} - A^{23} - 13A^{22} + 14A^{21} + 7A^{20} - 18A^{19} + 8A^{18} + 9A^{17} - 10A^{16} - A^{15} + 5A^{14} + 4A^{13} - 10A^{12} + 2A^{11} + 15A^{10} - 13A^{9} - 4A^{8} + 18A^{7} - 7A^{6} - 9A^{5} + 9A^{4} + A^{3} - 5A^{2} + 1$}
\end{minipage}

\noindent{\color{gray!40}\rule{\textwidth}{0.4pt}}
\vspace{0.9\baselineskip}
\noindent \begin{minipage}[t]{0.25\textwidth}
\vspace{0pt}
\centering
\includegraphics[page=396,width=\linewidth]{knotoids.pdf}
\end{minipage}
\hfill
\begin{minipage}[t]{0.73\textwidth}
\vspace{0pt}
\raggedright
\textbf{Name:} {\large{$\mathbf{K7_{277}}$}} (chiral, non-rotatable$^{*}$) \\ \textbf{PD:} {\scriptsize\texttt{[0],[0,1,2,3],[1,4,5,2],[3,6,7,4],[8,9,10,5],[6,10,11,7],[11,12,13,8],[14,13,12,9],[14]}} \\ \textbf{EM:} {\scriptsize\texttt{(B0, A0C0C3D0, B1D3E3B2, B3F0F3C1, G3H3F1C2, D1E2G0D2, F2H2H1E0, I0G2G1E1, H0)}} \\ \textbf{Kauffman bracket:} {\scriptsize $A^{24} + A^{22} - A^{20} - 3A^{18} + 4A^{14} + 4A^{12} - 2A^{10} - 4A^{8} - A^{6} + A^{4} + A^{2}$} \\ \textbf{Arrow:} {\scriptsize $1 + L_1/A^{2} - 1/A^{4} - 3L_1/A^{6} - L_2/A^{8} + A^{-8} + 4L_1/A^{10} + 3L_2/A^{12} + A^{-12} - 2L_1/A^{14} - 3L_2/A^{16} - 1/A^{16} - L_1/A^{18} + L_2/A^{20} + L_1/A^{22}$} \\ \textbf{Mock:} {\scriptsize $-w^{6} - w^{5} + 2w^{4} + 3w^{3} + w^{2} - 3w - 4 + 1/w + 4/w^{2} - 1/w^{4}$} \\ \textbf{Affine:} {\scriptsize $0$} \\ \textbf{Yamada:} {\scriptsize $-A^{28} + 3A^{26} - 4A^{24} + 3A^{23} + 6A^{22} - 4A^{21} - 4A^{20} + 9A^{19} - 8A^{17} + 6A^{16} + 3A^{15} - 6A^{14} + 2A^{13} + 4A^{12} - 6A^{10} + 5A^{9} + 5A^{8} - 11A^{7} + A^{6} + 6A^{5} - 4A^{4} - 2A^{3} + 3A^{2} + A - 1$}
\end{minipage}

\noindent{\color{gray!40}\rule{\textwidth}{0.4pt}}
\vspace{0.9\baselineskip}
\noindent \begin{minipage}[t]{0.25\textwidth}
\vspace{0pt}
\centering
\includegraphics[page=397,width=\linewidth]{knotoids.pdf}
\end{minipage}
\hfill
\begin{minipage}[t]{0.73\textwidth}
\vspace{0pt}
\raggedright
\textbf{Name:} {\large{$\mathbf{K7_{278}}$}} (chiral, non-rotatable$^{*}$) \\ \textbf{PD:} {\scriptsize\texttt{[0],[0,1,2,3],[1,4,5,2],[6,7,4,3],[5,8,9,6],[10,11,12,7],[8,12,13,9],[13,14,11,10],[14]}} \\ \textbf{EM:} {\scriptsize\texttt{(B0, A0C0C3D3, B1D2E0B2, E3F3C1B3, C2G0G3D0, H3H2G1D1, E1F2H0E2, G2I0F1F0, H1)}} \\ \textbf{Kauffman bracket:} {\scriptsize $-A^{21} - A^{19} - A^{17} + A^{15} + 3A^{13} + A^{11} - 3A^{9} - 3A^{7} + 2A^{3} + A$} \\ \textbf{Arrow:} {\scriptsize $A^{-12} + L_1/A^{14} + L_2/A^{16} - L_1/A^{18} - 3L_2/A^{20} - L_1/A^{22} + 3L_2/A^{24} + 3L_1/A^{26} - L_2/A^{28} + A^{-28} - 2L_1/A^{30} - 1/A^{32}$} \\ \textbf{Mock:} {\scriptsize $w^{6} + w^{5} - w^{3} - w^{2} - w + 2 + 3/w - 2/w^{2} - 2/w^{3} + w^{-4}$} \\ \textbf{Affine:} {\scriptsize $2t - 4 + 2/t$} \\ \textbf{Yamada:} {\scriptsize $A^{29} + A^{28} - A^{27} - A^{26} + 2A^{25} + A^{24} + 5A^{21} + 3A^{20} - 4A^{19} + 4A^{18} + 3A^{17} - 4A^{16} + A^{15} - A^{13} - 4A^{12} + A^{11} + A^{10} - 5A^{9} + 4A^{7} - A^{6} - 3A^{5} + 3A^{4} + A^{3} - 2A^{2} + 1$}
\end{minipage}

\noindent{\color{gray!40}\rule{\textwidth}{0.4pt}}
\vspace{0.9\baselineskip}
\noindent \begin{minipage}[t]{0.25\textwidth}
\vspace{0pt}
\centering
\includegraphics[page=398,width=\linewidth]{knotoids.pdf}
\end{minipage}
\hfill
\begin{minipage}[t]{0.73\textwidth}
\vspace{0pt}
\raggedright
\textbf{Name:} {\large{$\mathbf{K7_{279}}$}} (chiral, non-rotatable$^{*}$) \\ \textbf{PD:} {\scriptsize\texttt{[0],[0,1,2,3],[1,4,5,2],[3,6,7,4],[5,8,9,6],[7,9,10,11],[8,11,12,13],[13,12,14,10],[14]}} \\ \textbf{EM:} {\scriptsize\texttt{(B0, A0C0C3D0, B1D3E0B2, B3E3F0C1, C2G0F1D1, D2E2H3G1, E1F3H1H0, G3G2I0F2, H2)}} \\ \textbf{Kauffman bracket:} {\scriptsize $A^{21} + A^{19} - 2A^{17} - 4A^{15} + A^{13} + 4A^{11} + A^{9} - 3A^{7} - 2A^{5} + A^{3} + A$} \\ \textbf{Arrow:} {\scriptsize $-A^{6}L_1 - A^{4} + 2A^{2}L_1 + 4 - L_1/A^{2} - 4/A^{4} - L_1/A^{6} + 3/A^{8} + 2L_1/A^{10} - 1/A^{12} - L_1/A^{14}$} \\ \textbf{Mock:} {\scriptsize $w^{4} - 3w^{2} + 5 - 3/w^{2} + w^{-4}$} \\ \textbf{Affine:} {\scriptsize $0$} \\ \textbf{Yamada:} {\scriptsize $A^{26} - A^{25} + 4A^{23} - 5A^{22} - 2A^{21} + 8A^{20} - 6A^{19} - 4A^{18} + 5A^{17} - 4A^{16} - 3A^{15} - A^{14} + 4A^{13} - A^{12} - 4A^{11} + 6A^{10} - 7A^{8} + 4A^{7} + 3A^{6} - 6A^{5} + A^{4} + 4A^{3} - 2A^{2} - A + 1$}
\end{minipage}

\noindent{\color{gray!40}\rule{\textwidth}{0.4pt}}
\vspace{0.9\baselineskip}
\noindent \begin{minipage}[t]{0.25\textwidth}
\vspace{0pt}
\centering
\includegraphics[page=399,width=\linewidth]{knotoids.pdf}
\end{minipage}
\hfill
\begin{minipage}[t]{0.73\textwidth}
\vspace{0pt}
\raggedright
\textbf{Name:} {\large{$\mathbf{K7_{280}}$}} (chiral, non-rotatable$^{*}$) \\ \textbf{PD:} {\scriptsize\texttt{[0],[0,1,2,3],[1,4,5,2],[3,6,7,4],[5,8,9,10],[10,11,12,6],[7,13,14,8],[13,12,11,9],[14]}} \\ \textbf{EM:} {\scriptsize\texttt{(B0, A0C0C3D0, B1D3E0B2, B3F3G0C1, C2G3H3F0, E3H2H1D1, D2H0I0E1, G1F2F1E2, G2)}} \\ \textbf{Kauffman bracket:} {\scriptsize $-A^{17} - A^{15} + 2A^{13} + 2A^{11} - A^{9} - 3A^{7} - A^{5} + A^{3} + A$} \\ \textbf{Arrow:} {\scriptsize $A^{2}L_1 + 1 - 2L_1/A^{2} - 2/A^{4} + L_1/A^{6} + 3/A^{8} + L_1/A^{10} - 1/A^{12} - L_1/A^{14}$} \\ \textbf{Mock:} {\scriptsize $-w^{3} - w^{2} + 2w + 4 - 2/w^{2} - 1/w^{3}$} \\ \textbf{Affine:} {\scriptsize $t - 2 + 1/t$} \\ \textbf{Yamada:} {\scriptsize $A^{23} - A^{22} - A^{21} + 3A^{20} - A^{19} - A^{18} + 4A^{17} + 2A^{13} - A^{11} + 3A^{10} - 2A^{8} + A^{7} + A^{6} - 2A^{5} - A^{4} + 2A^{3} - A^{2} - A + 1$}
\end{minipage}

\noindent{\color{gray!40}\rule{\textwidth}{0.4pt}}
\vspace{0.9\baselineskip}
\noindent \begin{minipage}[t]{0.25\textwidth}
\vspace{0pt}
\centering
\includegraphics[page=400,width=\linewidth]{knotoids.pdf}
\end{minipage}
\hfill
\begin{minipage}[t]{0.73\textwidth}
\vspace{0pt}
\raggedright
\textbf{Name:} {\large{$\mathbf{K7_{281}}$}} (chiral, non-rotatable$^{*}$) \\ \textbf{PD:} {\scriptsize\texttt{[0],[0,1,2,3],[1,4,5,2],[3,6,7,4],[5,8,9,10],[10,9,11,6],[7,11,12,13],[8,13,14,12],[14]}} \\ \textbf{EM:} {\scriptsize\texttt{(B0, A0C0C3D0, B1D3E0B2, B3F3G0C1, C2H0F1F0, E3E2G1D1, D2F2H3H1, E1G3I0G2, H2)}} \\ \textbf{Kauffman bracket:} {\scriptsize $A^{19} - A^{17} - 2A^{15} + A^{13} + 3A^{11} - 3A^{7} - 2A^{5} + A^{3} + A$} \\ \textbf{Arrow:} {\scriptsize $-A^{10}L_1 + A^{8} + 2A^{6}L_1 - A^{4} - 3A^{2}L_1 - L_2 + 1 + 3L_1/A^{2} + 2L_2/A^{4} - L_1/A^{6} - L_2/A^{8}$} \\ \textbf{Mock:} {\scriptsize $-w^{5} - w^{4} + 2w^{3} + w^{2} - 3w + 1 + 3/w - 1/w^{3}$} \\ \textbf{Affine:} {\scriptsize $-t + 2 - 1/t$} \\ \textbf{Yamada:} {\scriptsize $A^{24} + A^{23} - 2A^{22} + 3A^{20} - 4A^{19} - 2A^{18} + 3A^{17} - 2A^{16} - A^{15} + 2A^{13} - 2A^{12} - 2A^{11} + 3A^{10} - 2A^{9} - 3A^{8} + 2A^{7} - 3A^{5} + 3A^{3} - A^{2} - A + 1$}
\end{minipage}

\noindent{\color{gray!40}\rule{\textwidth}{0.4pt}}
\vspace{0.9\baselineskip}
\noindent \begin{minipage}[t]{0.25\textwidth}
\vspace{0pt}
\centering
\includegraphics[page=401,width=\linewidth]{knotoids.pdf}
\end{minipage}
\hfill
\begin{minipage}[t]{0.73\textwidth}
\vspace{0pt}
\raggedright
\textbf{Name:} {\large{$\mathbf{K7_{282}}$}} (chiral, non-rotatable$^{*}$) \\ \textbf{PD:} {\scriptsize\texttt{[0],[0,1,2,3],[1,4,5,2],[3,6,7,4],[5,8,9,10],[10,11,12,6],[13,14,8,7],[13,12,11,9],[14]}} \\ \textbf{EM:} {\scriptsize\texttt{(B0, A0C0C3D0, B1D3E0B2, B3F3G3C1, C2G2H3F0, E3H2H1D1, H0I0E1D2, G0F2F1E2, G1)}} \\ \textbf{Kauffman bracket:} {\scriptsize $A^{16} + 2A^{14} - A^{12} - 4A^{10} - A^{8} + 4A^{6} + 3A^{4} - A^{2} - 2$} \\ \textbf{Arrow:} {\scriptsize $L_1/A^{2} + 2/A^{4} - L_1/A^{6} - 4/A^{8} - L_1/A^{10} + 4/A^{12} + 3L_1/A^{14} - 1/A^{16} - 2L_1/A^{18}$} \\ \textbf{Mock:} {\scriptsize $-w^{4} - w^{3} + 4w^{2} + 3w - 3 - 1/w + 2/w^{2} - 1/w^{3} - 1/w^{4}$} \\ \textbf{Affine:} {\scriptsize $2t - 4 + 2/t$} \\ \textbf{Yamada:} {\scriptsize $-A^{26} + A^{25} + 3A^{24} - 2A^{23} - 2A^{22} + 7A^{21} - 6A^{19} + 5A^{18} + 2A^{17} - 4A^{16} + 2A^{15} + 2A^{14} + A^{13} - 4A^{12} + 3A^{11} + 4A^{10} - 7A^{9} + 2A^{8} + 5A^{7} - 5A^{6} - 2A^{5} + 3A^{4} - 2A^{2} + 1$}
\end{minipage}

\noindent{\color{gray!40}\rule{\textwidth}{0.4pt}}
\vspace{0.9\baselineskip}
\noindent \begin{minipage}[t]{0.25\textwidth}
\vspace{0pt}
\centering
\includegraphics[page=402,width=\linewidth]{knotoids.pdf}
\end{minipage}
\hfill
\begin{minipage}[t]{0.73\textwidth}
\vspace{0pt}
\raggedright
\textbf{Name:} {\large{$\mathbf{K7_{283}}$}} (chiral, non-rotatable$^{*}$) \\ \textbf{PD:} {\scriptsize\texttt{[0],[0,1,2,3],[1,4,5,2],[3,6,7,4],[8,9,6,5],[9,10,11,7],[8,11,12,13],[13,12,14,10],[14]}} \\ \textbf{EM:} {\scriptsize\texttt{(B0, A0C0C3D0, B1D3E3B2, B3E2F3C1, G0F0D1C2, E1H3G1D2, E0F2H1H0, G3G2I0F1, H2)}} \\ \textbf{Kauffman bracket:} {\scriptsize $A^{23} + A^{21} - A^{19} - 4A^{17} - 2A^{15} + 5A^{13} + 5A^{11} - A^{9} - 5A^{7} - 2A^{5} + A^{3} + A$} \\ \textbf{Arrow:} {\scriptsize $-A^{20} - A^{18}L_1 + A^{16} + 4A^{14}L_1 + 2A^{12} - 5A^{10}L_1 - 5A^{8} + A^{6}L_1 + 5A^{4} + 2A^{2}L_1 - 1 - L_1/A^{2}$} \\ \textbf{Mock:} {\scriptsize $-w^{5} + 4w^{3} + 2w^{2} - 5w - 4 + 1/w + 4/w^{2} + 2/w^{3} - 1/w^{4} - 1/w^{5}$} \\ \textbf{Affine:} {\scriptsize $0$} \\ \textbf{Yamada:} {\scriptsize $A^{28} - 3A^{26} + 6A^{24} - 3A^{23} - 8A^{22} + 11A^{21} + 3A^{20} - 15A^{19} + 7A^{18} + 7A^{17} - 9A^{16} + 6A^{14} + 2A^{13} - 9A^{12} + 3A^{11} + 8A^{10} - 14A^{9} - 3A^{8} + 12A^{7} - 7A^{6} - 7A^{5} + 7A^{4} + 2A^{3} - 5A^{2} + 2$}
\end{minipage}

\noindent{\color{gray!40}\rule{\textwidth}{0.4pt}}
\vspace{0.9\baselineskip}
\noindent \begin{minipage}[t]{0.25\textwidth}
\vspace{0pt}
\centering
\includegraphics[page=403,width=\linewidth]{knotoids.pdf}
\end{minipage}
\hfill
\begin{minipage}[t]{0.73\textwidth}
\vspace{0pt}
\raggedright
\textbf{Name:} {\large{$\mathbf{K7_{284}}$}} (chiral, rotatable) \\ \textbf{PD:} {\scriptsize\texttt{[0],[0,1,2,3],[1,4,5,2],[3,6,7,4],[8,9,6,5],[9,10,11,7],[11,12,13,8],[10,13,12,14],[14]}} \\ \textbf{EM:} {\scriptsize\texttt{(B0, A0C0C3D0, B1D3E3B2, B3E2F3C1, G3F0D1C2, E1H0G0D2, F2H2H1E0, F1G2G1I0, H3)}} \\ \textbf{Kauffman bracket:} {\scriptsize $-A^{26} + 3A^{22} + 2A^{20} - 4A^{18} - 6A^{16} + 2A^{14} + 8A^{12} + 3A^{10} - 4A^{8} - 3A^{6} + A^{2}$} \\ \textbf{Arrow:} {\scriptsize $-A^{8} + 3A^{4} + 2A^{2}L_1 - 4 - 6L_1/A^{2} + 2/A^{4} + 8L_1/A^{6} + 3/A^{8} - 4L_1/A^{10} - 3/A^{12} + A^{-16}$} \\ \textbf{Mock:} {\scriptsize $w^{6} + 2w^{5} - 2w^{4} - 6w^{3} + 2w^{2} + 8w + 1 - 4/w - 1/w^{2}$} \\ \textbf{Affine:} {\scriptsize $2t - 4 + 2/t$} \\ \textbf{Yamada:} {\scriptsize $-A^{27} + A^{26} + 4A^{25} - 6A^{24} - 5A^{23} + 15A^{22} - 4A^{21} - 19A^{20} + 19A^{19} + 7A^{18} - 23A^{17} + 9A^{16} + 10A^{15} - 13A^{14} - 4A^{13} + 6A^{12} + 5A^{11} - 17A^{10} + 3A^{9} + 20A^{8} - 20A^{7} - 5A^{6} + 22A^{5} - 9A^{4} - 9A^{3} + 10A^{2} + A - 3$}
\end{minipage}

\noindent{\color{gray!40}\rule{\textwidth}{0.4pt}}
\vspace{0.9\baselineskip}
\noindent \begin{minipage}[t]{0.25\textwidth}
\vspace{0pt}
\centering
\includegraphics[page=404,width=\linewidth]{knotoids.pdf}
\end{minipage}
\hfill
\begin{minipage}[t]{0.73\textwidth}
\vspace{0pt}
\raggedright
\textbf{Name:} {\large{$\mathbf{K7_{285}}$}} (chiral, non-rotatable$^{*}$) \\ \textbf{PD:} {\scriptsize\texttt{[0],[0,1,2,3],[1,4,5,2],[3,6,7,4],[8,9,10,5],[6,10,11,12],[7,13,14,8],[9,13,12,11],[14]}} \\ \textbf{EM:} {\scriptsize\texttt{(B0, A0C0C3D0, B1D3E3B2, B3F0G0C1, G3H0F1C2, D1E2H3H2, D2H1I0E0, E1G1F3F2, G2)}} \\ \textbf{Kauffman bracket:} {\scriptsize $A^{23} + A^{21} - 2A^{19} - 3A^{17} + A^{15} + 4A^{13} + A^{11} - 3A^{9} - A^{7} + A^{5} - A$} \\ \textbf{Arrow:} {\scriptsize $-A^{26}L_1 - A^{24} + 2A^{22}L_1 + 3A^{20} - A^{18}L_1 - 4A^{16} - A^{14}L_1 + 3A^{12} + A^{10}L_1 - A^{8} + A^{4}$} \\ \textbf{Mock:} {\scriptsize $w^{3} + 3w^{2} - 2w - 6 + 4/w^{2} + w^{-3}$} \\ \textbf{Affine:} {\scriptsize $-t + 2 - 1/t$} \\ \textbf{Yamada:} {\scriptsize $-A^{25} + A^{24} + 2A^{23} - 3A^{22} + 4A^{20} - 3A^{19} - 2A^{18} + 4A^{17} - A^{16} - 4A^{15} + 3A^{14} - 3A^{12} + A^{10} + 2A^{9} - 4A^{8} + 4A^{7} + 4A^{6} - 5A^{5} + 3A^{4} + 3A^{3} - 2A^{2} + 2A + 1$}
\end{minipage}

\noindent{\color{gray!40}\rule{\textwidth}{0.4pt}}
\vspace{0.9\baselineskip}
\noindent \begin{minipage}[t]{0.25\textwidth}
\vspace{0pt}
\centering
\includegraphics[page=405,width=\linewidth]{knotoids.pdf}
\end{minipage}
\hfill
\begin{minipage}[t]{0.73\textwidth}
\vspace{0pt}
\raggedright
\textbf{Name:} {\large{$\mathbf{K7_{286}}$}} (chiral, non-rotatable$^{*}$) \\ \textbf{PD:} {\scriptsize\texttt{[0],[0,1,2,3],[1,4,5,2],[3,6,7,4],[8,9,10,5],[6,10,9,11],[7,11,12,13],[8,13,14,12],[14]}} \\ \textbf{EM:} {\scriptsize\texttt{(B0, A0C0C3D0, B1D3E3B2, B3F0G0C1, H0F2F1C2, D1E2E1G1, D2F3H3H1, E0G3I0G2, H2)}} \\ \textbf{Kauffman bracket:} {\scriptsize $A^{18} + A^{16} - A^{14} - 2A^{12} + 2A^{8} + 2A^{6} - A^{4} - A^{2}$} \\ \textbf{Arrow:} {\scriptsize $1 + L_1/A^{2} - 1/A^{4} - 2L_1/A^{6} - L_2/A^{8} + A^{-8} + 2L_1/A^{10} + 2L_2/A^{12} - L_1/A^{14} - L_2/A^{16}$} \\ \textbf{Mock:} {\scriptsize $-w^{6} - w^{5} + 2w^{4} + 2w^{3} - w - 1 + 1/w + 2/w^{2} - 1/w^{3} - 1/w^{4}$} \\ \textbf{Affine:} {\scriptsize $t - 2 + 1/t$} \\ \textbf{Yamada:} {\scriptsize $A^{23} - 2A^{21} + A^{19} - 2A^{18} - 2A^{17} + A^{15} - 2A^{14} + 2A^{12} - 2A^{11} - A^{10} + A^{9} - 2A^{8} - A^{7} + 2A^{5} - A^{4} - A^{3} + 3A^{2} - 1$}
\end{minipage}

\noindent{\color{gray!40}\rule{\textwidth}{0.4pt}}
\vspace{0.9\baselineskip}
\noindent \begin{minipage}[t]{0.25\textwidth}
\vspace{0pt}
\centering
\includegraphics[page=406,width=\linewidth]{knotoids.pdf}
\end{minipage}
\hfill
\begin{minipage}[t]{0.73\textwidth}
\vspace{0pt}
\raggedright
\textbf{Name:} {\large{$\mathbf{K7_{287}}$}} (chiral, non-rotatable$^{*}$) \\ \textbf{PD:} {\scriptsize\texttt{[0],[0,1,2,3],[1,4,5,2],[3,6,7,4],[8,9,10,5],[6,10,11,12],[13,14,8,7],[9,13,12,11],[14]}} \\ \textbf{EM:} {\scriptsize\texttt{(B0, A0C0C3D0, B1D3E3B2, B3F0G3C1, G2H0F1C2, D1E2H3H2, H1I0E0D2, E1G0F3F2, G1)}} \\ \textbf{Kauffman bracket:} {\scriptsize $-A^{26} - A^{24} + 2A^{22} + 4A^{20} - 6A^{16} - 3A^{14} + 4A^{12} + 4A^{10} - A^{8} - 2A^{6} + A^{2}$} \\ \textbf{Arrow:} {\scriptsize $-A^{20} - A^{18}L_1 + 2A^{16} + 4A^{14}L_1 - 6A^{10}L_1 - 3A^{8} + 4A^{6}L_1 + 4A^{4} - A^{2}L_1 - 2 + A^{-4}$} \\ \textbf{Mock:} {\scriptsize $w^{4} + 3w^{3} - 7w - 3 + 5/w + 4/w^{2} - 1/w^{3} - 1/w^{4}$} \\ \textbf{Affine:} {\scriptsize $0$} \\ \textbf{Yamada:} {\scriptsize $-A^{27} + A^{26} + 3A^{25} - 4A^{24} - 3A^{23} + 9A^{22} - 2A^{21} - 12A^{20} + 11A^{19} + 5A^{18} - 15A^{17} + 6A^{16} + 7A^{15} - 7A^{14} - A^{13} + 5A^{12} + 5A^{11} - 10A^{10} + 2A^{9} + 13A^{8} - 12A^{7} - 3A^{6} + 14A^{5} - 5A^{4} - 5A^{3} + 6A^{2} + A - 2$}
\end{minipage}

\noindent{\color{gray!40}\rule{\textwidth}{0.4pt}}
\vspace{0.9\baselineskip}
\noindent \begin{minipage}[t]{0.25\textwidth}
\vspace{0pt}
\centering
\includegraphics[page=407,width=\linewidth]{knotoids.pdf}
\end{minipage}
\hfill
\begin{minipage}[t]{0.73\textwidth}
\vspace{0pt}
\raggedright
\textbf{Name:} {\large{$\mathbf{K7_{288}}$}} (chiral, rotatable) \\ \textbf{PD:} {\scriptsize\texttt{[0],[0,1,2,3],[1,4,5,2],[6,7,4,3],[5,8,9,6],[7,9,10,11],[11,12,13,8],[10,13,12,14],[14]}} \\ \textbf{EM:} {\scriptsize\texttt{(B0, A0C0C3D3, B1D2E0B2, E3F0C1B3, C2G3F1D0, D1E2H0G0, F3H2H1E1, F2G2G1I0, H3)}} \\ \textbf{Kauffman bracket:} {\scriptsize $A^{18} + 2A^{16} + A^{14} - 4A^{12} - 4A^{10} + 2A^{8} + 6A^{6} + 2A^{4} - 3A^{2} - 2$} \\ \textbf{Arrow:} {\scriptsize $A^{-12} + 2L_1/A^{14} + A^{-16} - 4L_1/A^{18} - 4/A^{20} + 2L_1/A^{22} + 6/A^{24} + 2L_1/A^{26} - 3/A^{28} - 2L_1/A^{30}$} \\ \textbf{Mock:} {\scriptsize $w^{6} + 2w^{5} + w^{4} - 4w^{3} - 4w^{2} + 2w + 5 + 2/w - 3/w^{2} - 2/w^{3} + w^{-4}$} \\ \textbf{Affine:} {\scriptsize $2t - 4 + 2/t$} \\ \textbf{Yamada:} {\scriptsize $A^{29} + A^{28} - 2A^{27} - A^{26} + 7A^{25} - 7A^{23} + 7A^{22} + 7A^{21} - 10A^{20} + 11A^{18} - 5A^{17} - 3A^{16} + 8A^{15} - 5A^{13} - 3A^{12} + 8A^{11} - 7A^{10} - 9A^{9} + 12A^{8} - 2A^{7} - 8A^{6} + 5A^{5} + 3A^{4} - 3A^{3} - A^{2} + A + 1$}
\end{minipage}

\noindent{\color{gray!40}\rule{\textwidth}{0.4pt}}
\vspace{0.9\baselineskip}
\noindent \begin{minipage}[t]{0.25\textwidth}
\vspace{0pt}
\centering
\includegraphics[page=408,width=\linewidth]{knotoids.pdf}
\end{minipage}
\hfill
\begin{minipage}[t]{0.73\textwidth}
\vspace{0pt}
\raggedright
\textbf{Name:} {\large{$\mathbf{K7_{289}}$}} (chiral, non-rotatable$^{*}$) \\ \textbf{PD:} {\scriptsize\texttt{[0],[0,1,2,3],[1,4,5,2],[6,7,4,3],[5,8,9,10],[10,11,12,6],[7,13,14,8],[13,12,11,9],[14]}} \\ \textbf{EM:} {\scriptsize\texttt{(B0, A0C0C3D3, B1D2E0B2, F3G0C1B3, C2G3H3F0, E3H2H1D0, D1H0I0E1, G1F2F1E2, G2)}} \\ \textbf{Kauffman bracket:} {\scriptsize $-A^{20} - A^{18} + 2A^{16} + 3A^{14} - 2A^{12} - 5A^{10} + 6A^{6} + 3A^{4} - 2A^{2} - 2$} \\ \textbf{Arrow:} {\scriptsize $-A^{14}L_1 - A^{12} + 2A^{10}L_1 + 3A^{8} - 2A^{6}L_1 - 5A^{4} + 6 + 3L_1/A^{2} - 2/A^{4} - 2L_1/A^{6}$} \\ \textbf{Mock:} {\scriptsize $w^{4} - w^{3} - 4w^{2} + w + 7 + 1/w - 4/w^{2} - 1/w^{3} + w^{-4}$} \\ \textbf{Affine:} {\scriptsize $0$} \\ \textbf{Yamada:} {\scriptsize $2A^{28} - 5A^{26} + 2A^{25} + 6A^{24} - 9A^{23} - 4A^{22} + 12A^{21} - 4A^{20} - 9A^{19} + 10A^{18} + 2A^{17} - 7A^{16} + 4A^{15} + 6A^{14} + A^{13} - 7A^{12} + 9A^{11} + 5A^{10} - 14A^{9} + 6A^{8} + 8A^{7} - 9A^{6} - A^{5} + 4A^{4} - A^{3} - 2A^{2} + 1$}
\end{minipage}

\noindent{\color{gray!40}\rule{\textwidth}{0.4pt}}
\vspace{0.9\baselineskip}
\noindent \begin{minipage}[t]{0.25\textwidth}
\vspace{0pt}
\centering
\includegraphics[page=409,width=\linewidth]{knotoids.pdf}
\end{minipage}
\hfill
\begin{minipage}[t]{0.73\textwidth}
\vspace{0pt}
\raggedright
\textbf{Name:} {\large{$\mathbf{K7_{290}}$}} (chiral, non-rotatable$^{*}$) \\ \textbf{PD:} {\scriptsize\texttt{[0],[0,1,2,3],[1,3,4,5],[2,6,7,8],[4,9,10,11],[5,11,12,6],[13,12,14,7],[14,10,9,8],[13]}} \\ \textbf{EM:} {\scriptsize\texttt{(B0, A0C0D0C1, B1B3E0F0, B2F3G3H3, C2H2H1F1, C3E3G1D1, I0F2H0D2, G2E2E1D3, G0)}} \\ \textbf{Kauffman bracket:} {\scriptsize $2A^{14} + 2A^{12} - 2A^{10} - 3A^{8} + 2A^{4} + A^{2} - 1$} \\ \textbf{Arrow:} {\scriptsize $2/A^{4} + 2L_1/A^{6} - 2/A^{8} - 3L_1/A^{10} - L_2/A^{12} + A^{-12} + 2L_1/A^{14} + L_2/A^{16} - L_1/A^{18}$} \\ \textbf{Mock:} {\scriptsize $-w^{4} - w^{3} + 3w^{2} + 4w - 3/w - 1/w^{2}$} \\ \textbf{Affine:} {\scriptsize $-t^{2} + 2t - 2 + 2/t - 1/t^{2}$} \\ \textbf{Yamada:} {\scriptsize $A^{22} + A^{21} - 3A^{20} - 6A^{17} - A^{16} + A^{15} - 3A^{14} + A^{13} + A^{12} + 3A^{11} + 3A^{8} - 3A^{7} - 2A^{6} + 2A^{5} - A^{4} - 2A^{3} + 2A^{2} + A - 1$}
\end{minipage}

\noindent{\color{gray!40}\rule{\textwidth}{0.4pt}}
\vspace{0.9\baselineskip}
\noindent \begin{minipage}[t]{0.25\textwidth}
\vspace{0pt}
\centering
\includegraphics[page=410,width=\linewidth]{knotoids.pdf}
\end{minipage}
\hfill
\begin{minipage}[t]{0.73\textwidth}
\vspace{0pt}
\raggedright
\textbf{Name:} {\large{$\mathbf{K7_{291}}$}} (chiral, rotatable) \\ \textbf{PD:} {\scriptsize\texttt{[0],[0,1,2,3],[1,4,5,6],[7,8,9,2],[10,11,4,3],[11,12,13,5],[6,13,8,7],[12,14,10,9],[14]}} \\ \textbf{EM:} {\scriptsize\texttt{(B0, A0C0D3E3, B1E2F3G0, G3G2H3B2, H2F0C1B3, E1H0G1C2, C3F2D1D0, F1I0E0D2, H1)}} \\ \textbf{Kauffman bracket:} {\scriptsize $-2A^{24} - A^{22} + 4A^{20} + 4A^{18} - 4A^{16} - 5A^{14} + 2A^{12} + 5A^{10} - 3A^{6} + A^{2}$} \\ \textbf{Arrow:} {\scriptsize $-2A^{18}L_1 - A^{16} + 4A^{14}L_1 + 4A^{12} - 4A^{10}L_1 - 5A^{8} + 2A^{6}L_1 + 5A^{4} - 3 + A^{-4}$} \\ \textbf{Mock:} {\scriptsize $-w^{4} + 3w^{2} - 4w - 6 + 6/w + 7/w^{2} - 2/w^{3} - 2/w^{4}$} \\ \textbf{Affine:} {\scriptsize $-2t + 4 - 2/t$} \\ \textbf{Yamada:} {\scriptsize $-A^{27} + A^{26} + A^{25} - 4A^{24} + 8A^{22} - 2A^{21} - 9A^{20} + 11A^{19} + 4A^{18} - 15A^{17} + 7A^{16} + 6A^{15} - 9A^{14} - A^{13} + 4A^{12} + 2A^{11} - 10A^{10} + 3A^{9} + 9A^{8} - 13A^{7} - 3A^{6} + 12A^{5} - 6A^{4} - 6A^{3} + 7A^{2} + A - 3$}
\end{minipage}

\noindent{\color{gray!40}\rule{\textwidth}{0.4pt}}
\vspace{0.9\baselineskip}
\noindent \begin{minipage}[t]{0.25\textwidth}
\vspace{0pt}
\centering
\includegraphics[page=411,width=\linewidth]{knotoids.pdf}
\end{minipage}
\hfill
\begin{minipage}[t]{0.73\textwidth}
\vspace{0pt}
\raggedright
\textbf{Name:} {\large{$\mathbf{K7_{292}}$}} (chiral, non-rotatable$^{*}$) \\ \textbf{PD:} {\scriptsize\texttt{[0],[0,1,2,3],[1,4,5,6],[6,7,8,2],[3,9,10,4],[11,12,7,5],[12,13,9,8],[13,14,11,10],[14]}} \\ \textbf{EM:} {\scriptsize\texttt{(B0, A0C0D3E0, B1E3F3D0, C3F2G3B2, B3G2H3C1, H2G0D1C2, F1H0E1D2, G1I0F0E2, H1)}} \\ \textbf{Kauffman bracket:} {\scriptsize $-A^{26} - A^{24} + 2A^{22} + 5A^{20} - 7A^{16} - 4A^{14} + 5A^{12} + 5A^{10} - A^{8} - 3A^{6} + A^{2}$} \\ \textbf{Arrow:} {\scriptsize $-A^{14}L_1 - A^{12} + 2A^{10}L_1 + A^{8}L_2 + 4A^{8} - 2A^{4}L_2 - 5A^{4} - 4A^{2}L_1 + L_2 + 4 + 5L_1/A^{2} - 1/A^{4} - 3L_1/A^{6} + L_1/A^{10}$} \\ \textbf{Mock:} {\scriptsize $w^{4} + w^{3} - 4w^{2} - 3w + 7 + 3/w - 5/w^{2} - 2/w^{3} + 2/w^{4} + w^{-5}$} \\ \textbf{Affine:} {\scriptsize $-t + 2 - 1/t$} \\ \textbf{Yamada:} {\scriptsize $-A^{27} + 2A^{26} + 4A^{25} - 7A^{24} - 5A^{23} + 14A^{22} - 3A^{21} - 19A^{20} + 16A^{19} + 9A^{18} - 20A^{17} + 8A^{16} + 12A^{15} - 9A^{14} - 3A^{13} + 8A^{12} + 7A^{11} - 15A^{10} + 2A^{9} + 19A^{8} - 18A^{7} - 7A^{6} + 20A^{5} - 6A^{4} - 9A^{3} + 7A^{2} + 2A - 2$}
\end{minipage}

\noindent{\color{gray!40}\rule{\textwidth}{0.4pt}}
\vspace{0.9\baselineskip}
\noindent \begin{minipage}[t]{0.25\textwidth}
\vspace{0pt}
\centering
\includegraphics[page=412,width=\linewidth]{knotoids.pdf}
\end{minipage}
\hfill
\begin{minipage}[t]{0.73\textwidth}
\vspace{0pt}
\raggedright
\textbf{Name:} {\large{$\mathbf{K7_{293}}$}} (chiral, non-rotatable$^{*}$) \\ \textbf{PD:} {\scriptsize\texttt{[0],[0,1,2,3],[1,4,5,2],[3,6,7,4],[5,8,9,6],[7,10,11,8],[9,11,12,13],[13,12,14,10],[14]}} \\ \textbf{EM:} {\scriptsize\texttt{(B0, A0C0C3D0, B1D3E0B2, B3E3F0C1, C2F3G0D1, D2H3G1E1, E2F2H1H0, G3G2I0F1, H2)}} \\ \textbf{Kauffman bracket:} {\scriptsize $A^{14} + 2A^{12} - 2A^{8} - 2A^{6} + A^{2} + 1$} \\ \textbf{Arrow:} {\scriptsize $A^{2}L_1 + 2 - 2/A^{4} - 2L_1/A^{6} - L_2/A^{8} + A^{-8} + L_1/A^{10} + L_2/A^{12}$} \\ \textbf{Mock:} {\scriptsize $w^{3} - w + 2 - 1/w^{2}$} \\ \textbf{Affine:} {\scriptsize $-t^{2} + t + 1/t - 1/t^{2}$} \\ \textbf{Yamada:} {\scriptsize $-A^{26} + 2A^{22} + A^{20} + 3A^{19} + A^{17} + A^{15} - A^{14} - A^{13} + A^{12} - A^{11} - A^{10} - A^{7} + A^{6} + 2A^{5} + A^{2} - 1$}
\end{minipage}

\noindent{\color{gray!40}\rule{\textwidth}{0.4pt}}
\vspace{0.9\baselineskip}
\noindent \begin{minipage}[t]{0.25\textwidth}
\vspace{0pt}
\centering
\includegraphics[page=413,width=\linewidth]{knotoids.pdf}
\end{minipage}
\hfill
\begin{minipage}[t]{0.73\textwidth}
\vspace{0pt}
\raggedright
\textbf{Name:} {\large{$\mathbf{K7_{294}}$}} (chiral, non-rotatable$^{*}$) \\ \textbf{PD:} {\scriptsize\texttt{[0],[0,1,2,3],[1,4,5,2],[3,6,7,4],[5,8,9,6],[7,10,11,8],[11,12,13,9],[10,13,12,14],[14]}} \\ \textbf{EM:} {\scriptsize\texttt{(B0, A0C0C3D0, B1D3E0B2, B3E3F0C1, C2F3G3D1, D2H0G0E1, F2H2H1E2, F1G2G1I0, H3)}} \\ \textbf{Kauffman bracket:} {\scriptsize $-A^{13} - A^{11} + A^{9} + 3A^{7} + A^{5} - 2A^{3} - 2A$} \\ \textbf{Arrow:} {\scriptsize $A^{16} + A^{14}L_1 - A^{12} - 3A^{10}L_1 - A^{8}L_2 + 2A^{6}L_1 + A^{4}L_2 + A^{4}$} \\ \textbf{Mock:} {\scriptsize $w^{3} - 3w + 2/w + w^{-2}$} \\ \textbf{Affine:} {\scriptsize $-t^{2} - t + 4 - 1/t - 1/t^{2}$} \\ \textbf{Yamada:} {\scriptsize $-A^{22} + A^{21} - 3A^{19} + 2A^{18} + A^{17} - A^{16} + 2A^{15} + A^{14} + A^{13} - A^{12} + A^{11} + 2A^{10} - A^{9} + A^{8} + 2A^{7} - 2A^{6} - A^{5} + A^{4} - A^{2} + A + 1$}
\end{minipage}

\noindent{\color{gray!40}\rule{\textwidth}{0.4pt}}
\vspace{0.9\baselineskip}
\noindent \begin{minipage}[t]{0.25\textwidth}
\vspace{0pt}
\centering
\includegraphics[page=414,width=\linewidth]{knotoids.pdf}
\end{minipage}
\hfill
\begin{minipage}[t]{0.73\textwidth}
\vspace{0pt}
\raggedright
\textbf{Name:} {\large{$\mathbf{K7_{295}}$}} (chiral, non-rotatable$^{*}$) \\ \textbf{PD:} {\scriptsize\texttt{[0],[0,1,2,3],[1,4,5,2],[3,6,7,4],[5,8,9,6],[10,11,8,7],[9,11,12,13],[13,12,14,10],[14]}} \\ \textbf{EM:} {\scriptsize\texttt{(B0, A0C0C3D0, B1D3E0B2, B3E3F3C1, C2F2G0D1, H3G1E1D2, E2F1H1H0, G3G2I0F0, H2)}} \\ \textbf{Kauffman bracket:} {\scriptsize $-A^{18} + A^{16} + 4A^{14} + A^{12} - 5A^{10} - 5A^{8} + 2A^{6} + 4A^{4} + A^{2} - 1$} \\ \textbf{Arrow:} {\scriptsize $-A^{12} + A^{10}L_1 + 4A^{8} + A^{6}L_1 - 5A^{4} - 5A^{2}L_1 - L_2 + 3 + 4L_1/A^{2} + L_2/A^{4} - L_1/A^{6}$} \\ \textbf{Mock:} {\scriptsize $w^{4} + 2w^{3} - 3w^{2} - 4w + 5 + 3/w - 2/w^{2} - 1/w^{3}$} \\ \textbf{Affine:} {\scriptsize $-t^{2} + t + 1/t - 1/t^{2}$} \\ \textbf{Yamada:} {\scriptsize $-A^{26} + A^{25} + A^{24} - 5A^{23} + 3A^{22} + 7A^{21} - 10A^{20} + A^{19} + 11A^{18} - 8A^{17} - A^{16} + 7A^{15} - 2A^{14} - 2A^{13} - A^{12} + 6A^{11} - 3A^{10} - 7A^{9} + 11A^{8} - A^{7} - 9A^{6} + 9A^{5} + 2A^{4} - 6A^{3} + 2A^{2} + 2A - 1$}
\end{minipage}

\noindent{\color{gray!40}\rule{\textwidth}{0.4pt}}
\vspace{0.9\baselineskip}
\noindent \begin{minipage}[t]{0.25\textwidth}
\vspace{0pt}
\centering
\includegraphics[page=415,width=\linewidth]{knotoids.pdf}
\end{minipage}
\hfill
\begin{minipage}[t]{0.73\textwidth}
\vspace{0pt}
\raggedright
\textbf{Name:} {\large{$\mathbf{K7_{296}}$}} (chiral, non-rotatable$^{*}$) \\ \textbf{PD:} {\scriptsize\texttt{[0],[0,1,2,3],[1,4,5,2],[3,6,7,4],[5,8,9,6],[10,11,8,7],[11,12,13,9],[10,13,12,14],[14]}} \\ \textbf{EM:} {\scriptsize\texttt{(B0, A0C0C3D0, B1D3E0B2, B3E3F3C1, C2F2G3D1, H0G0E1D2, F1H2H1E2, F0G2G1I0, H3)}} \\ \textbf{Kauffman bracket:} {\scriptsize $A^{22} + A^{20} - 3A^{18} - 5A^{16} + A^{14} + 7A^{12} + 3A^{10} - 3A^{8} - 3A^{6} + A^{4} + A^{2}$} \\ \textbf{Arrow:} {\scriptsize $A^{22}L_1 + A^{20} - 3A^{18}L_1 - A^{16}L_2 - 4A^{16} + A^{14}L_1 + A^{12}L_2 + 6A^{12} + 3A^{10}L_1 - 3A^{8} - 3A^{6}L_1 + A^{4} + A^{2}L_1$} \\ \textbf{Mock:} {\scriptsize $-w^{4} + 5w^{2} - 7 - 1/w + 4/w^{2} + w^{-3}$} \\ \textbf{Affine:} {\scriptsize $-t^{2} - t + 4 - 1/t - 1/t^{2}$} \\ \textbf{Yamada:} {\scriptsize $-A^{27} + 2A^{25} - 3A^{24} - 2A^{23} + 9A^{22} - A^{21} - 10A^{20} + 11A^{19} + 5A^{18} - 15A^{17} + 6A^{16} + 6A^{15} - 8A^{14} - A^{13} + 3A^{12} + 3A^{11} - 11A^{10} + A^{9} + 11A^{8} - 12A^{7} - 3A^{6} + 13A^{5} - 6A^{4} - 6A^{3} + 6A^{2} - 3$}
\end{minipage}

\noindent{\color{gray!40}\rule{\textwidth}{0.4pt}}
\vspace{0.9\baselineskip}
\noindent \begin{minipage}[t]{0.25\textwidth}
\vspace{0pt}
\centering
\includegraphics[page=416,width=\linewidth]{knotoids.pdf}
\end{minipage}
\hfill
\begin{minipage}[t]{0.73\textwidth}
\vspace{0pt}
\raggedright
\textbf{Name:} {\large{$\mathbf{K7_{297}}$}} (chiral, non-rotatable$^{*}$) \\ \textbf{PD:} {\scriptsize\texttt{[0],[0,1,2,3],[1,4,5,2],[3,6,7,4],[5,7,8,9],[6,9,10,11],[11,12,13,8],[10,13,12,14],[14]}} \\ \textbf{EM:} {\scriptsize\texttt{(B0, A0C0C3D0, B1D3E0B2, B3F0E1C1, C2D2G3F1, D1E3H0G0, F3H2H1E2, F2G2G1I0, H3)}} \\ \textbf{Kauffman bracket:} {\scriptsize $A^{22} + A^{20} - 2A^{18} - 5A^{16} + 6A^{12} + 4A^{10} - 2A^{8} - 3A^{6} + A^{2}$} \\ \textbf{Arrow:} {\scriptsize $A^{16} + A^{14}L_1 - A^{12}L_2 - A^{12} - 5A^{10}L_1 + A^{8}L_2 - A^{8} + 6A^{6}L_1 + 4A^{4} - 2A^{2}L_1 - 3 + A^{-4}$} \\ \textbf{Mock:} {\scriptsize $w^{3} - 5w - 1 + 6/w + 3/w^{2} - 2/w^{3} - 1/w^{4}$} \\ \textbf{Affine:} {\scriptsize $-t^{2} - t + 4 - 1/t - 1/t^{2}$} \\ \textbf{Yamada:} {\scriptsize $-A^{27} + 3A^{25} - 2A^{24} - 4A^{23} + 8A^{22} + 3A^{21} - 11A^{20} + 5A^{19} + 8A^{18} - 12A^{17} + 6A^{15} - 5A^{14} - 3A^{13} + A^{12} + 6A^{11} - 8A^{10} - 3A^{9} + 12A^{8} - 7A^{7} - 7A^{6} + 10A^{5} - A^{4} - 6A^{3} + 2A^{2} + A - 1$}
\end{minipage}

\noindent{\color{gray!40}\rule{\textwidth}{0.4pt}}
\vspace{0.9\baselineskip}
\noindent \begin{minipage}[t]{0.25\textwidth}
\vspace{0pt}
\centering
\includegraphics[page=417,width=\linewidth]{knotoids.pdf}
\end{minipage}
\hfill
\begin{minipage}[t]{0.73\textwidth}
\vspace{0pt}
\raggedright
\textbf{Name:} {\large{$\mathbf{K7_{298}}$}} (chiral, non-rotatable$^{*}$) \\ \textbf{PD:} {\scriptsize\texttt{[0],[0,1,2,3],[1,4,5,2],[3,6,7,4],[5,7,8,9],[9,10,11,6],[8,11,12,13],[13,12,14,10],[14]}} \\ \textbf{EM:} {\scriptsize\texttt{(B0, A0C0C3D0, B1D3E0B2, B3F3E1C1, C2D2G0F0, E3H3G1D1, E2F2H1H0, G3G2I0F1, H2)}} \\ \textbf{Kauffman bracket:} {\scriptsize $2A^{14} + 2A^{12} - 3A^{8} - 2A^{6} + A^{2} + 1$} \\ \textbf{Arrow:} {\scriptsize $L_2/A^{4} + A^{-4} + 2L_1/A^{6} - L_2/A^{8} + A^{-8} - 3L_1/A^{10} - 2/A^{12} + A^{-16} + L_1/A^{18}$} \\ \textbf{Mock:} {\scriptsize $w + 1 - 2/w + w^{-3}$} \\ \textbf{Affine:} {\scriptsize $t^{2} - 2 + t^{-2}$} \\ \textbf{Yamada:} {\scriptsize $A^{27} + A^{23} + 2A^{22} - A^{21} + A^{20} + 3A^{19} - A^{18} + A^{15} - 2A^{14} - A^{13} + A^{12} - 2A^{11} + 2A^{9} + A^{8} - A^{7} + A^{6} + 2A^{5} - A^{4} - A^{3} + A^{2} - 1$}
\end{minipage}

\noindent{\color{gray!40}\rule{\textwidth}{0.4pt}}
\vspace{0.9\baselineskip}
\noindent \begin{minipage}[t]{0.25\textwidth}
\vspace{0pt}
\centering
\includegraphics[page=418,width=\linewidth]{knotoids.pdf}
\end{minipage}
\hfill
\begin{minipage}[t]{0.73\textwidth}
\vspace{0pt}
\raggedright
\textbf{Name:} {\large{$\mathbf{K7_{299}}$}} (chiral, non-rotatable$^{*}$) \\ \textbf{PD:} {\scriptsize\texttt{[0],[0,1,2,3],[1,4,5,2],[3,6,7,4],[5,7,8,9],[9,10,11,6],[11,12,13,8],[10,13,12,14],[14]}} \\ \textbf{EM:} {\scriptsize\texttt{(B0, A0C0C3D0, B1D3E0B2, B3F3E1C1, C2D2G3F0, E3H0G0D1, F2H2H1E2, F1G2G1I0, H3)}} \\ \textbf{Kauffman bracket:} {\scriptsize $A^{16} + 2A^{14} - A^{12} - 4A^{10} - 2A^{8} + 3A^{6} + 3A^{4} - 1$} \\ \textbf{Arrow:} {\scriptsize $A^{10}L_1 + 2A^{8} - A^{6}L_1 - A^{4}L_2 - 3A^{4} - 2A^{2}L_1 + L_2 + 2 + 3L_1/A^{2} - L_1/A^{6}$} \\ \textbf{Mock:} {\scriptsize $-2w^{2} - w + 5 + 2/w - 2/w^{2} - 1/w^{3}$} \\ \textbf{Affine:} {\scriptsize $-t^{2} + 2 - 1/t^{2}$} \\ \textbf{Yamada:} {\scriptsize $-A^{26} + 2A^{24} - A^{23} - 3A^{22} + 4A^{21} + 2A^{20} - 5A^{19} + 3A^{18} + 3A^{17} - 3A^{16} + A^{15} + 2A^{14} + 2A^{13} - 3A^{12} + A^{11} + 4A^{10} - 5A^{9} + 5A^{7} - 2A^{6} - 2A^{5} + 2A^{4} + A^{3} - 2A^{2} + 1$}
\end{minipage}

\noindent{\color{gray!40}\rule{\textwidth}{0.4pt}}
\vspace{0.9\baselineskip}
\noindent \begin{minipage}[t]{0.25\textwidth}
\vspace{0pt}
\centering
\includegraphics[page=419,width=\linewidth]{knotoids.pdf}
\end{minipage}
\hfill
\begin{minipage}[t]{0.73\textwidth}
\vspace{0pt}
\raggedright
\textbf{Name:} {\large{$\mathbf{K7_{300}}$}} (chiral, non-rotatable$^{*}$) \\ \textbf{PD:} {\scriptsize\texttt{[0],[0,1,2,3],[1,4,5,2],[3,6,7,4],[8,9,6,5],[7,10,11,8],[9,11,12,13],[13,12,14,10],[14]}} \\ \textbf{EM:} {\scriptsize\texttt{(B0, A0C0C3D0, B1D3E3B2, B3E2F0C1, F3G0D1C2, D2H3G1E0, E1F2H1H0, G3G2I0F1, H2)}} \\ \textbf{Kauffman bracket:} {\scriptsize $A^{22} - 3A^{18} - A^{16} + 5A^{14} + 5A^{12} - 3A^{10} - 6A^{8} + 3A^{4} + A^{2} - 1$} \\ \textbf{Arrow:} {\scriptsize $A^{4} - 3 - L_1/A^{2} + 5/A^{4} + 5L_1/A^{6} - 3/A^{8} - 6L_1/A^{10} - L_2/A^{12} + A^{-12} + 3L_1/A^{14} + L_2/A^{16} - L_1/A^{18}$} \\ \textbf{Mock:} {\scriptsize $-2w^{4} - 2w^{3} + 6w^{2} + 7w - 3 - 6/w + w^{-3}$} \\ \textbf{Affine:} {\scriptsize $-t^{2} + 2t - 2 + 2/t - 1/t^{2}$} \\ \textbf{Yamada:} {\scriptsize $-2A^{28} + A^{27} + 5A^{26} - 6A^{25} - 4A^{24} + 14A^{23} - 3A^{22} - 11A^{21} + 14A^{20} + 3A^{19} - 10A^{18} + 5A^{17} + 5A^{16} - A^{15} - 8A^{14} + 7A^{13} + 5A^{12} - 16A^{11} + 5A^{10} + 10A^{9} - 12A^{8} - 2A^{7} + 10A^{6} - 2A^{5} - 4A^{4} + 3A^{3} + 2A^{2} - A - 1$}
\end{minipage}

\noindent{\color{gray!40}\rule{\textwidth}{0.4pt}}
\vspace{0.9\baselineskip}
\noindent \begin{minipage}[t]{0.25\textwidth}
\vspace{0pt}
\centering
\includegraphics[page=420,width=\linewidth]{knotoids.pdf}
\end{minipage}
\hfill
\begin{minipage}[t]{0.73\textwidth}
\vspace{0pt}
\raggedright
\textbf{Name:} {\large{$\mathbf{K7_{301}}$}} (chiral, non-rotatable$^{*}$) \\ \textbf{PD:} {\scriptsize\texttt{[0],[0,1,2,3],[1,4,5,2],[3,6,7,4],[8,9,6,5],[10,11,8,7],[11,12,13,9],[10,13,12,14],[14]}} \\ \textbf{EM:} {\scriptsize\texttt{(B0, A0C0C3D0, B1D3E3B2, B3E2F3C1, F2G3D1C2, H0G0E0D2, F1H2H1E1, F0G2G1I0, H3)}} \\ \textbf{Kauffman bracket:} {\scriptsize $A^{18} + A^{16} - A^{14} - 3A^{12} - 2A^{10} + 2A^{8} + 3A^{6} + A^{4} - A^{2}$} \\ \textbf{Arrow:} {\scriptsize $A^{18}L_1 + A^{16} - A^{14}L_1 - A^{12}L_2 - 2A^{12} - 2A^{10}L_1 + A^{8}L_2 + A^{8} + 3A^{6}L_1 + A^{4} - A^{2}L_1$} \\ \textbf{Mock:} {\scriptsize $-w^{2} - w + 3 + 2/w - 1/w^{2} - 1/w^{3}$} \\ \textbf{Affine:} {\scriptsize $-t^{2} + 2 - 1/t^{2}$} \\ \textbf{Yamada:} {\scriptsize $A^{23} - A^{22} - 2A^{21} + 2A^{20} - 3A^{18} + 2A^{17} + 2A^{16} - 3A^{15} + 2A^{14} + 2A^{13} - A^{12} + A^{11} + 3A^{9} - 2A^{8} + 4A^{6} - 3A^{5} - A^{4} + 3A^{3} - A^{2} + 1$}
\end{minipage}

\noindent{\color{gray!40}\rule{\textwidth}{0.4pt}}
\vspace{0.9\baselineskip}
\noindent \begin{minipage}[t]{0.25\textwidth}
\vspace{0pt}
\centering
\includegraphics[page=421,width=\linewidth]{knotoids.pdf}
\end{minipage}
\hfill
\begin{minipage}[t]{0.73\textwidth}
\vspace{0pt}
\raggedright
\textbf{Name:} {\large{$\mathbf{K7_{302}}$}} (chiral, non-rotatable$^{*}$) \\ \textbf{PD:} {\scriptsize\texttt{[0],[0,1,2,3],[1,4,5,2],[3,6,7,4],[7,8,9,5],[6,9,10,11],[11,12,13,8],[10,13,12,14],[14]}} \\ \textbf{EM:} {\scriptsize\texttt{(B0, A0C0C3D0, B1D3E3B2, B3F0E0C1, D2G3F1C2, D1E2H0G0, F3H2H1E1, F2G2G1I0, H3)}} \\ \textbf{Kauffman bracket:} {\scriptsize $-A^{21} - A^{19} + A^{17} + 4A^{15} + 2A^{13} - 3A^{11} - 4A^{9} + 2A^{5} - A$} \\ \textbf{Arrow:} {\scriptsize $A^{24} + A^{22}L_1 - A^{20} - 4A^{18}L_1 - A^{16}L_2 - A^{16} + 3A^{14}L_1 + A^{12}L_2 + 3A^{12} - 2A^{8} + A^{4}$} \\ \textbf{Mock:} {\scriptsize $w^{3} + w^{2} - 4w - 3 + 3/w + 3/w^{2}$} \\ \textbf{Affine:} {\scriptsize $-t^{2} - 2t + 6 - 2/t - 1/t^{2}$} \\ \textbf{Yamada:} {\scriptsize $-A^{25} - A^{24} + 3A^{23} + A^{22} - 4A^{21} + 3A^{20} + 5A^{19} - 5A^{18} - 2A^{17} + 6A^{16} - 4A^{15} - 3A^{14} + 3A^{13} - A^{12} - 2A^{11} - 3A^{10} + 5A^{9} - 2A^{8} - 5A^{7} + 7A^{6} - 6A^{4} + 2A^{3} + A^{2} - 2A - 1$}
\end{minipage}

\noindent{\color{gray!40}\rule{\textwidth}{0.4pt}}
\vspace{0.9\baselineskip}
\noindent \begin{minipage}[t]{0.25\textwidth}
\vspace{0pt}
\centering
\includegraphics[page=422,width=\linewidth]{knotoids.pdf}
\end{minipage}
\hfill
\begin{minipage}[t]{0.73\textwidth}
\vspace{0pt}
\raggedright
\textbf{Name:} {\large{$\mathbf{K7_{303}}$}} (chiral, non-rotatable$^{*}$) \\ \textbf{PD:} {\scriptsize\texttt{[0],[0,1,2,3],[1,4,5,2],[3,6,7,4],[7,8,9,5],[9,10,11,6],[8,11,12,13],[13,12,14,10],[14]}} \\ \textbf{EM:} {\scriptsize\texttt{(B0, A0C0C3D0, B1D3E3B2, B3F3E0C1, D2G0F0C2, E2H3G1D1, E1F2H1H0, G3G2I0F1, H2)}} \\ \textbf{Kauffman bracket:} {\scriptsize $-A^{18} + 2A^{14} + 3A^{12} - A^{10} - 3A^{8} - A^{6} + A^{4} + A^{2}$} \\ \textbf{Arrow:} {\scriptsize $-A^{6}L_1 + 2A^{2}L_1 + L_2 + 2 - L_1/A^{2} - L_2/A^{4} - 2/A^{4} - L_1/A^{6} + A^{-8} + L_1/A^{10}$} \\ \textbf{Mock:} {\scriptsize $-w^{3} - 2w^{2} + w + 4 - 1/w^{2}$} \\ \textbf{Affine:} {\scriptsize $t^{2} - t - 1/t + t^{-2}$} \\ \textbf{Yamada:} {\scriptsize $-A^{24} + 2A^{23} + A^{22} - 3A^{21} + 3A^{20} + 2A^{19} - 2A^{18} + 2A^{17} + A^{16} + A^{15} - 2A^{14} + A^{13} + 2A^{12} - 3A^{11} + 2A^{9} - 2A^{8} - A^{7} + 3A^{6} + A^{5} - A^{4} + A^{3} + A^{2} - A - 1$}
\end{minipage}

\noindent{\color{gray!40}\rule{\textwidth}{0.4pt}}
\vspace{0.9\baselineskip}
\noindent \begin{minipage}[t]{0.25\textwidth}
\vspace{0pt}
\centering
\includegraphics[page=423,width=\linewidth]{knotoids.pdf}
\end{minipage}
\hfill
\begin{minipage}[t]{0.73\textwidth}
\vspace{0pt}
\raggedright
\textbf{Name:} {\large{$\mathbf{K7_{304}}$}} (chiral, non-rotatable$^{*}$) \\ \textbf{PD:} {\scriptsize\texttt{[0],[0,1,2,3],[1,4,5,2],[6,7,4,3],[5,8,9,6],[7,10,11,8],[9,11,12,13],[13,12,14,10],[14]}} \\ \textbf{EM:} {\scriptsize\texttt{(B0, A0C0C3D3, B1D2E0B2, E3F0C1B3, C2F3G0D0, D1H3G1E1, E2F2H1H0, G3G2I0F1, H2)}} \\ \textbf{Kauffman bracket:} {\scriptsize $A^{16} + 2A^{14} - 3A^{10} - 3A^{8} + A^{6} + 2A^{4} + A^{2}$} \\ \textbf{Arrow:} {\scriptsize $A^{10}L_1 + 2A^{8} - 3A^{4} - 3A^{2}L_1 - L_2 + 2 + 2L_1/A^{2} + L_2/A^{4}$} \\ \textbf{Mock:} {\scriptsize $w^{3} - w^{2} - 2w + 3 + 1/w - 1/w^{2}$} \\ \textbf{Affine:} {\scriptsize $-t^{2} + 2 - 1/t^{2}$} \\ \textbf{Yamada:} {\scriptsize $-A^{26} - A^{23} + A^{22} + 2A^{21} - A^{20} + A^{19} + 4A^{18} - 2A^{17} + 2A^{15} - 2A^{14} - A^{12} + 2A^{11} - A^{10} - 2A^{9} + 3A^{8} - A^{7} - A^{6} + 3A^{5} + A^{4} - A^{3} + A^{2} + A - 1$}
\end{minipage}

\noindent{\color{gray!40}\rule{\textwidth}{0.4pt}}
\vspace{0.9\baselineskip}
\noindent \begin{minipage}[t]{0.25\textwidth}
\vspace{0pt}
\centering
\includegraphics[page=424,width=\linewidth]{knotoids.pdf}
\end{minipage}
\hfill
\begin{minipage}[t]{0.73\textwidth}
\vspace{0pt}
\raggedright
\textbf{Name:} {\large{$\mathbf{K7_{305}}$}} (chiral, non-rotatable$^{*}$) \\ \textbf{PD:} {\scriptsize\texttt{[0],[0,1,2,3],[1,4,5,2],[6,7,4,3],[5,8,9,6],[10,11,8,7],[11,12,13,9],[10,13,12,14],[14]}} \\ \textbf{EM:} {\scriptsize\texttt{(B0, A0C0C3D3, B1D2E0B2, E3F3C1B3, C2F2G3D0, H0G0E1D1, F1H2H1E2, F0G2G1I0, H3)}} \\ \textbf{Kauffman bracket:} {\scriptsize $-A^{20} - 2A^{18} - A^{16} + 3A^{14} + 4A^{12} - 3A^{8} - A^{6} + A^{4} + A^{2}$} \\ \textbf{Arrow:} {\scriptsize $-A^{26}L_1 - A^{24}L_2 - A^{24} - A^{22}L_1 + A^{20}L_2 + 2A^{20} + 4A^{18}L_1 - 3A^{14}L_1 - A^{12} + A^{10}L_1 + A^{8}$} \\ \textbf{Mock:} {\scriptsize $-w^{3} + 2w^{2} + 2w - 3 - 3/w + 2/w^{3} + 2/w^{4}$} \\ \textbf{Affine:} {\scriptsize $-t^{2} - 2t + 6 - 2/t - 1/t^{2}$} \\ \textbf{Yamada:} {\scriptsize $-A^{29} - A^{28} + A^{27} + A^{26} - 2A^{25} + 4A^{23} + A^{22} - 3A^{21} + 2A^{20} + 4A^{19} - 3A^{18} - 3A^{17} + 5A^{16} - 2A^{15} - 4A^{14} + 3A^{13} - A^{12} - 2A^{11} - 2A^{10} + 5A^{9} - A^{8} - 5A^{7} + 4A^{6} + A^{5} - 4A^{4} - A^{3} + A^{2} - A - 2$}
\end{minipage}

\noindent{\color{gray!40}\rule{\textwidth}{0.4pt}}
\vspace{0.9\baselineskip}
\noindent \begin{minipage}[t]{0.25\textwidth}
\vspace{0pt}
\centering
\includegraphics[page=425,width=\linewidth]{knotoids.pdf}
\end{minipage}
\hfill
\begin{minipage}[t]{0.73\textwidth}
\vspace{0pt}
\raggedright
\textbf{Name:} {\large{$\mathbf{K7_{306}}$}} (chiral, non-rotatable$^{*}$) \\ \textbf{PD:} {\scriptsize\texttt{[0],[0,1,2,3],[1,4,5,2],[6,7,4,3],[5,7,8,9],[6,9,10,11],[11,12,13,8],[10,13,12,14],[14]}} \\ \textbf{EM:} {\scriptsize\texttt{(B0, A0C0C3D3, B1D2E0B2, F0E1C1B3, C2D1G3F1, D0E3H0G0, F3H2H1E2, F2G2G1I0, H3)}} \\ \textbf{Kauffman bracket:} {\scriptsize $A^{17} + 2A^{15} - 3A^{11} - 3A^{9} + A^{7} + 2A^{5} - A$} \\ \textbf{Arrow:} {\scriptsize $-A^{20}L_2 - 2A^{18}L_1 + A^{16}L_2 - A^{16} + 3A^{14}L_1 + 3A^{12} - A^{10}L_1 - 2A^{8} + A^{4}$} \\ \textbf{Mock:} {\scriptsize $w^{2} - 2w - 2 + 3/w + 2/w^{2} - 1/w^{3}$} \\ \textbf{Affine:} {\scriptsize $-t^{2} - t + 4 - 1/t - 1/t^{2}$} \\ \textbf{Yamada:} {\scriptsize $-A^{27} - A^{26} + A^{25} + A^{24} - A^{23} + A^{22} + 3A^{21} - 2A^{19} + 2A^{18} + A^{17} - 3A^{16} + 2A^{14} - 3A^{13} - A^{12} + 2A^{11} - 2A^{10} - A^{9} - A^{8} + 3A^{7} - 2A^{6} - 2A^{5} + 3A^{4} - 2A^{3} - 2A^{2} - 1$}
\end{minipage}

\noindent{\color{gray!40}\rule{\textwidth}{0.4pt}}
\vspace{0.9\baselineskip}
\noindent \begin{minipage}[t]{0.25\textwidth}
\vspace{0pt}
\centering
\includegraphics[page=426,width=\linewidth]{knotoids.pdf}
\end{minipage}
\hfill
\begin{minipage}[t]{0.73\textwidth}
\vspace{0pt}
\raggedright
\textbf{Name:} {\large{$\mathbf{K7_{307}}$}} (chiral, non-rotatable$^{*}$) \\ \textbf{PD:} {\scriptsize\texttt{[0],[0,1,2,3],[1,4,5,2],[6,7,4,3],[8,9,6,5],[7,10,11,8],[9,11,12,13],[13,12,14,10],[14]}} \\ \textbf{EM:} {\scriptsize\texttt{(B0, A0C0C3D3, B1D2E3B2, E2F0C1B3, F3G0D0C2, D1H3G1E0, E1F2H1H0, G3G2I0F1, H2)}} \\ \textbf{Kauffman bracket:} {\scriptsize $A^{21} - 4A^{17} - 3A^{15} + 4A^{13} + 6A^{11} - 5A^{7} - 2A^{5} + A^{3} + A$} \\ \textbf{Arrow:} {\scriptsize $-A^{6}L_1 + 4A^{2}L_1 + 3 - 4L_1/A^{2} - L_2/A^{4} - 5/A^{4} + L_2/A^{8} + 4/A^{8} + 2L_1/A^{10} - 1/A^{12} - L_1/A^{14}$} \\ \textbf{Mock:} {\scriptsize $-w^{3} - 2w^{2} + 3w + 7 - 2/w - 5/w^{2} + w^{-4}$} \\ \textbf{Affine:} {\scriptsize $-t^{2} + t + 1/t - 1/t^{2}$} \\ \textbf{Yamada:} {\scriptsize $A^{28} + 2A^{27} - 2A^{26} - 3A^{25} + 6A^{24} - 13A^{22} + 6A^{21} + 8A^{20} - 14A^{19} + 2A^{18} + 9A^{17} - 7A^{16} - 3A^{15} + 3A^{14} + 5A^{13} - 8A^{12} - A^{11} + 13A^{10} - 9A^{9} - 7A^{8} + 12A^{7} - 2A^{6} - 9A^{5} + 4A^{4} + 3A^{3} - 3A^{2} + 1$}
\end{minipage}

\noindent{\color{gray!40}\rule{\textwidth}{0.4pt}}
\vspace{0.9\baselineskip}
\noindent \begin{minipage}[t]{0.25\textwidth}
\vspace{0pt}
\centering
\includegraphics[page=427,width=\linewidth]{knotoids.pdf}
\end{minipage}
\hfill
\begin{minipage}[t]{0.73\textwidth}
\vspace{0pt}
\raggedright
\textbf{Name:} {\large{$\mathbf{K7_{308}}$}} (chiral, non-rotatable$^{*}$) \\ \textbf{PD:} {\scriptsize\texttt{[0],[0,1,2,3],[1,4,5,2],[3,6,7,4],[5,7,8,9],[6,9,10,11],[12,13,10,8],[11,13,14,12],[14]}} \\ \textbf{EM:} {\scriptsize\texttt{(B0, A0C0C3D0, B1D3E0B2, B3F0E1C1, C2D2G3F1, D1E3G2H0, H3H1F2E2, F3G1I0G0, H2)}} \\ \textbf{Kauffman bracket:} {\scriptsize $-2A^{17} - 2A^{15} + 2A^{13} + 4A^{11} + A^{9} - 3A^{7} - 2A^{5} + A$} \\ \textbf{Arrow:} {\scriptsize $2A^{2}L_1 + 2 - 2L_1/A^{2} - L_2/A^{4} - 3/A^{4} - L_1/A^{6} + L_2/A^{8} + 2/A^{8} + 2L_1/A^{10} - L_1/A^{14}$} \\ \textbf{Mock:} {\scriptsize $-w^{2} + w + 5 - 3/w^{2} - 1/w^{3}$} \\ \textbf{Affine:} {\scriptsize $-t^{2} + 2t - 2 + 2/t - 1/t^{2}$} \\ \textbf{Yamada:} {\scriptsize $-A^{24} + 2A^{22} - 2A^{21} - 2A^{20} + 6A^{19} - 4A^{17} + 5A^{16} + A^{15} - 3A^{14} + 2A^{13} + A^{12} + A^{11} - 3A^{10} + 2A^{9} + 3A^{8} - 6A^{7} + 2A^{6} + 4A^{5} - 4A^{4} + 3A^{2} - 1$}
\end{minipage}

\noindent{\color{gray!40}\rule{\textwidth}{0.4pt}}
\vspace{0.9\baselineskip}